\documentclass[11pt]{article}
\usepackage{amssymb,latexsym}
\usepackage{epsfig}
\usepackage{eufrak}
\usepackage{amsmath}
\usepackage{mathrsfs}
\usepackage{color}

\newtheorem{theorem}{Theorem}
\newtheorem{lemma}{Lemma}
\newtheorem{corollary}{Corollary}

\newcommand{\be}{\begin{equation}}
\newcommand{\ee}{\end{equation}}
\newcommand{\bee}{\begin{eqnarray*}}
\newcommand{\eee}{\end{eqnarray*}}
\newcommand{\bel}{\begin{eqnarray}}
\newcommand{\eel}{\end{eqnarray}}
\newcommand{\bec}{\begin{cases}}
\newcommand{\eec}{\end{cases}}
\newcommand{\bem}{\begin{bmatrix}}
\newcommand{\eem}{\end{bmatrix}}

\newcommand{\la}{\label}
\newcommand{\li}{\left}
\newcommand{\ri}{\right}

\newcommand{\DEF}{\stackrel{\mathrm{def}}{=}}

\newcommand{\ovl}{\overline}
\newcommand{\udl}{\underline}

\newcommand{\lc}{\lceil}
\newcommand{\rc}{\rceil}
\newcommand{\lf}{\lfloor}
\newcommand{\rf}{\rfloor}

\newcommand{\ep}{\epsilon}
\newcommand{\vep}{\varepsilon}
\newcommand{\lm}{\lambda}

\newcommand{\Up}{\Upsilon}
\newcommand{\up}{\upsilon}
\newcommand{\si}{\sigma}

\newcommand{\de}{\delta}

\newcommand{\vDe}{\varDelta}

\newcommand{\ga}{\gamma}
\newcommand{\Ga}{\Gamma}
\newcommand{\vse}{\vartheta}
\newcommand{\se}{\theta}
\newcommand{\Se}{\Theta}

\newcommand{\ze}{\zeta}
\newcommand{\al}{\alpha}
\newcommand{\ba}{\beta}

\newcommand{\vro}{\varrho}
\newcommand{\ro}{\rho}
\newcommand{\ka}{\kappa}
\newcommand{\om}{\omega}
\newcommand{\Om}{\Omega}

\newcommand{\f}{\frac}
\newcommand{\sq}{\sqrt}
\newcommand{\cd}{\cdots}

\newcommand{\qu}{\quad}
\newcommand{\qqu}{\qquad}
\newcommand{\fa}{\forall}

\newcommand{\mscr}{\mathscr}
\newcommand{\mcal}{\mathcal}
\newcommand{\mbf}{\mathbf}
\newcommand{\bb}{\mathbb}
\newcommand{\fra}{\mathfrak}

\newcommand{\wh}{\widehat}
\newcommand{\wt}{\widetilde}
\newcommand{\mrm}{\mathrm}
\newcommand{\bs}{\boldsymbol}

\newcommand{\ap}{\approx}

\newcommand{\LRA}{\Longleftrightarrow}
\newcommand{\sh}{\slash}

\newcommand{\tx}{\text}

\newcommand{\iy}{\infty}

\newcommand{\leu}{\subseteq}

\newcommand{\pa}{\partial}

\newcommand{\bed}{\begin{description}}
\newcommand{\eed}{\end{description}}
\newcommand{\bei}{\begin{itemize}}
\newcommand{\eei}{\end{itemize}}
\newcommand{\ben}{\begin{enumerate}}
\newcommand{\een}{\end{enumerate}}
\newcommand{\bib}{\bibitem}
\newcommand{\beL}{\begin{lemma}}
\newcommand{\eeL}{\end{lemma}}
\newcommand{\beT}{\begin{theorem}}
\newcommand{\eeT}{\end{theorem}}
\newcommand{\beC}{\begin{corollary}}
\newcommand{\eeC}{\end{corollary}}
\newcommand{\sect}{\section}

\newcommand{\bpf}{\begin{pf}}
\newcommand{\epf}{\end{pf}}
\newcommand{\bsk}{\bigskip}
\newcommand{\bi}{\binom}

\newcommand{\pfbox}{\hfill\mbox{$\Box$}}

\newenvironment{pf}{\paragraph*{Proof{\rm.}}}{\pfbox\bigskip}

\begin{document}

\title{{\bf A New Framework of Multistage Estimation}
\thanks{The author had been previously working with Louisiana
State University at Baton Rouge, LA 70803, USA, and is now with Department of Electrical Engineering, Southern University and A\&M College,
Baton Rouge, LA 70813, USA; Email: chenxinjia@gmail.com. The main results of this paper have been presented in Proceedings of SPIE Conferences,
Orlando, April 5-9, 2010 and April 25-29, 2011.  The statistical methodology proposed in this paper has been applied to electrical engineering
and computer science, see recent literature \cite{Chen_SPIE8387, Chen_SPIE, Chen_SPIE11a, Chen_SPIE11b, Chen_SPIE11c, Chencontrol, Chen_J} and
the references therein.}}

\author{Xinjia Chen}

\date{First submitted in September 2008}

\maketitle

$\qqu \qqu  \qqu $ {\it In Memory of My Dear Father Hualong Chen
(1933--1990)}

\begin{abstract}

In this paper, we have established a unified framework of multistage parameter estimation.  We demonstrate that a wide variety of statistical
problems such as fixed-sample-size interval estimation,  point estimation with error control, bounded-width confidence intervals, interval
estimation following hypothesis testing, construction of confidence sequences, can be cast into the general framework of constructing sequential
random intervals with prescribed coverage probabilities. We have developed exact methods for the construction of such sequential random
intervals in the context of multistage sampling.  In particular, we have established inclusion principle and coverage tuning techniques to
control and adjust the coverage probabilities of sequential random intervals. We have obtained concrete sampling schemes which are
unprecedentedly efficient in terms of sampling effort as compared to existing procedures.

\end{abstract}

\tableofcontents

\sect{Introduction}

Parameter estimation is a fundamental area of statistical inference,
which enjoys numerous applications in various fields of sciences and
engineering.  Specially, it is of ubiquitous significance to
estimate, via sampling,  the parameters of binomial, Poisson,
hypergeometrical, and normal distributions.   In general, a
parameter estimation problem can be formulated as follows. Let $X$
be a random variable defined in a probability space $(\Om, \mscr{F},
\Pr )$. Suppose the distribution of $X$ is determined by an unknown
parameter $\se$ in a parameter space $\Se$. In many applications, it
is desirable to construct a random interval which includes $\se$
with a prescribed level of confidence from random samples $X_1, X_2,
\cd$ of $X$. This problem is so fundamental that it has been
persistent issues of research in probability, statistics and other
relevant fields (see, e.g., \cite{Desu, Gosh, Govindarajulu, MukSi,
GR1986, Wald} and the references therein). Despite the richness of
literature devoted to such issues, existing approaches may suffer
from the drawbacks of lacking either efficiency or rigorousness.
Such drawbacks are mainly due to two frequently-used methods of
designing sampling schemes. The first method is to seek a worst-case
solution based on the assumption that the true parameter $\se$ is
included in an interval $[a, b] \subseteq \Se$.  Since it is
difficult to have tight bounds for the unknown parameter $\se$, such
a worst-case method can lead to overly wasteful sample size if the
interval $[a, b]$ is too wide. Moreover, if the true value of $\se$
is not included in $[a, b]$, the resultant sample size can be
misleading. The second method is to employ asymptotic theories such
as large deviations theory, Brownian motion theory, diffusion theory
and nonlinear renewal theory  in the design and analysis of sampling
schemes (see, \cite{Dembo, Lai, Siegmund, Vara, Woodroofe} and the
references therein). Undoubtedly, asymptotic techniques may offer
approximate solutions and important insight for the relevant
problems. Since any asymptotic theory holds only if the sample size
tends to infinity and, unfortunately, any practical sampling scheme
must be of a finite sample size, it is inevitable for an asymptotic
method to introduce unknown error in the resultant approximate
solution.

Motivated by the limitations of existing approaches of parameter estimation, we have established a new framework of multistage estimation, which
accommodates fixed-sample-size estimation and fully sequential estimation as special cases. The main characteristics of our multistage
estimation methods is as follows: i) Prior information of the parameter $\se$ is not necessary but can be used to make the estimation more
efficient; ii) The sampling schemes are asymptotically optimal in the sense that as the required precision gets high, the average sample number
is almost the same as the exact sample size computed as the true value of $\se$ were available; iii) The prescribed level of confidence is
rigorously guaranteed. Our new estimation techniques are developed under the spirit that parameter estimation, as an important branch of
statistical inference, should be accomplished with minimum cost in sampling and absolute rigorousness in quantifying uncertainty. In other
words, as many other researchers advocated, we propose to offer statistical inferential statements which guarantee prescribed level of
credibility and minimize conservatism as well. For example, we seek to provide statistical statements like ``with confidence level at least $100
(1 - \de)\%$, an estimator differs from its true value less than an {\it a priori} number $\vep$.'' In addition to guaranteeing  the desired
confidence level $100 (1 - \de)\%$, we try to make the true confidence level for each parametric value as close as possible to $100 (1 -
\de)\%$. Some aspects of our general framework are outlined as follows.

\bed

\item [(I):] We unify classical problems such as, point estimation
with precision requirements, construction of fixed-width confidence intervals, interval estimation based on a given sampling scheme, as a much
more general problem of constructing {\it sequential random intervals} with prescribed coverage probabilities. For example, the point estimation
problem of obtaining a point estimator $\wh{\bs{\se}}$ for $\se$ such that $\Pr \{ | \wh{\bs{\se}} - \se | < \vep \} > 1 - \de$ based on
multistage estimation can be considered as the problem of constructing sequential random interval $(\wh{\bs{\se}} - \vep, \wh{\bs{\se}} + \vep)$
with coverage probabilities greater than $1 - \de$ for all $\se \in \Se$.  The word ``sequential'' is used to indicate the fact that the random
interval is constructed from samples of random size.

\item [(II):]  We propose to construct stopping rules which are parameterized by a
number $\ze > 0$, referred to as {\it coverage tuning parameter}, such that the coverage probabilities of the associated sequential random
intervals can be controlled by $\ze$.  Here, by ``controlled'', we mean that the coverage probabilities can be adjusted to be above any
desirable level by making $\ze > 0$ sufficiently small.  To make the coverage probability of a sequential random interval controllable by $\ze$,
we propose to use a sequence of confidence intervals whose coverage probabilities can be controlled by $\ze$, referred to as {\it controlling
confidence sequence}, to determine a stopping rule such that {\it the sequential random interval must include the controlling confidence
sequence at the termination of sampling process} (see, e.g., Section 3 of the fifth version of our paper \cite{Chen_rule} published in arXiv on
April 7, 2009, our SPIE paper \cite{Chen_SPIE} published in April 2010, and our earlier versions of this paper from September 2008 to present).
 We call such a methodology of using confidence sequences to define stopping rules to control the coverage probabilities of the associated sequential
random intervals as {\it inclusion principle}. We have shown that if the coverage probability of the controlling confidence sequence can be
controlled by $\ze$, then the coverage probability of the sequential random interval can also be controlled by $\ze$. To make the coverage
probability of the controlling confidence sequence controllable by $\ze$, we propose to use lower and upper confidence limits $\mcal{L}_\ell, \;
\mcal{U}_\ell$ for the $\ell$-th stage such that the probability of $\{ \se \leq \mcal{L}_\ell \}$ is no greater than $\ze \de_\ell$ and that
the probability of $\{ \se \geq \mcal{U}_\ell \}$ is no greater than $\ze \de_\ell$, where $\de_\ell \in (0, 1)$ is independent of $\ze$. Of
course, conservative bounds or approximations of exact confidence limits may be used to construct stopping rules by the inclusion principle so
that the coverage probability of the desired sequential random interval can be controlled by $\ze$. Since the calculation of confidence limits
can be cumbersome and may involve solving complicated equations, we have managed to avoid such computation to make stopping rules as simple as
possible.

\item [(III):]  Once we have constructed stopping rules such that
the coverage probabilities of the associated sequential random interval is controllable by $\ze$.  Our next task is to seek the largest value of
the coverage tuning parameter $\ze$ such that the coverage probabilities of the sequential random interval is above the desired level. The
purpose of making $\ze$ as large as possible is to avoid unnecessary sampling effort.  The desired value of $\ze$ can be obtained by a method we
called {\it bisection coverage tuning}.  To achieve higher computational accuracy, we propose to evaluate the complementary coverage
probabilities. This is increasingly important as the desired level of coverage probabilities becomes higher, e.g., $0.9999$. A critical
subroutine of bisection coverage tuning is to determine whether the complementary coverage probabilities of the sequential random interval
corresponding to a fixed value of $\ze$ are no greater than the desired level for all parametric values of $\se \in \Se$. The major difficulty
of this subroutine is the computational complexity. First, for each parametric value, the evaluation of the complementary coverage probability
of the sequential random interval can be time-consuming. Second, the number of parametric values can be infinity or extremely large. Therefore,
we must avoid the exhaustive method of computing complementary coverage probabilities of the sequential random interval for all parametric
values.  In this direction, we have developed two algorithms to overcome the difficulty. The first algorithm is adapted from Branch and Bound
method in global optimization (see our earlier versions of this paper published on arXiv before July 2009). The second algorithm is called
Adaptive Maximum Checking Algorithm (AMCA). An indispensable technique for these two algorithms is the method of {\it interval bounding}.  That
is, how to bound the complementary coverage probabilities of the sequential random interval for parameter $\se \in [a, b]$. The tightness of
such bounds is extremely important for the efficiency of bisection coverage tuning. A simple idea of interval bounding is to express the
complementary coverage probability as a number of polynomial functions of $\se$, bound each function for $\se \in [a, b]$ by virtue of
monotonicity, and obtain bounds for the complementary coverage probability for $\se \in [a, b]$ using the relationship $\udl{q}_i < q_i <
\ovl{q}_i, \; i = 1, \cd, m \Rightarrow \sum_{i=1}^m \udl{q}_i < \sum_{i=1}^m q_i < \sum_{i=1}^m\ovl{q}_i$. We call this as {\it over-bounding}
method. Clearly, for a large $m$, the bounds derived from this method can be very conservative. In contrast to the over-bounding method, we have
obtained very tight bounds for the complementary coverage probabilities by exploiting the statistical properties of the sequential random
interval and the estimator of $\se$. In this regard, we have introduced the concept of unimodal-likelihood estimator (ULE).

\item [(IV):]  To start the bisection coverage tuning, we need to
find an initial interval of $\ze$.  For this purpose, we first use
results from asymptotic analysis of the coverage probabilities to
find a value $\ze_0$ for $\ze$ such that the corresponding coverage
probabilities are close to the desired level. Afterward, we use the
subroutine described above to find non-negative integers $i$ and $j$
as small as possible such that the complementary coverage
probabilities satisfy the requirement for $\ze = \ze_0
 2^{-i}$, but violate the requirement for $\ze = \ze_0
 2^j$.  Using $[\ze_0 2^{-i}, \ze_0 2^j ]$ as the starting interval,
 we can apply a bisection search to find a value of $\ze$ as large as
 possible such that the complementary coverage probability of the sequential
 random interval is not exceeding the pre-specified level for any
 parametric value.

\eed

The remainder of the paper is organized as follows. In Section 2, we present our general theory for the design and analysis of multistage
sampling schemes. Especially, we propose inclusion principle for construction of sampling schemes. We establish a general theory on coverage
probability of sequential random intervals which eliminates the necessity of exhaustive computation of coverage probability for designing
sampling schemes. In Section 3, we introduce powerful techniques such as bisection coverage tuning, consecutive-decision-variable bounding,
recursive computation, adaptive maximum checking, domain truncation and triangular partition that are crucial for a successful design of a
multistage sampling scheme. In Section 4, we present sampling schemes for estimation of binomial parameters. In Section 5, we discuss the
estimation of functions of two binomial proportions. Section 6 is devoted to the estimation of multinomial proportions. In Section 7, we
consider the estimation of means of bounded variables. In Section 8, we discuss the multistage estimation of Poisson parameters.   In Section 9,
we address the problem of estimating the proportion of a finite population.  In Section 10, we propose a general method for constructing
sampling schemes for parametric estimation associated with the absolute and relative error criteria based on prior information about parameters.
We consider the estimation of normal mean with unknown variance in Section 11. In Section 12, we discuss the estimation of the scale parameter
of a Gamma distribution.  In Section 13, we propose our exact methods for the construction of bounded-width confidence intervals. In Section 14,
we discuss the interval estimation based on a given sample scheme. In Section 15, we consider the exact construction of confidence sequences. In
Section 16, we address the problem of multistage linear regression.  In Section 17, we investigate the multistage estimation of quantile.
Section 18 is the conclusion. The proofs of all theorems are given in Appendices. We discuss various branch and bound methods in Appendix T.

Throughout this paper, we shall use the following notations.  The set of real numbers is denoted by $\bb{R}$.  The set of $n$-dimensional
vectors of real numbers is denoted by $\bb{R}^n$. The set of integers is denoted by $\bb{Z}$. The set of positive integers is denoted by
$\bb{N}$. The element of matrix $A$ in the $i$-th row and $j$-th column is denoted by $[A]_{i, j}$. The ceiling function and floor function are
denoted respectively by $\lc . \rc$ and $\lf . \rf$ (i.e., $\lc x \rc$ represents the smallest integer no less than $x$; $\lf x \rf$ represents
the largest integer no greater than $x$). The notation $\mrm{sgn}(x)$ denotes the sign function which assumes value $1$ for $x
> 0$, value $0$ for $x = 0$, and value $-1$ for $x < 0$. The gamma function is denoted by $\Ga(.)$.  For any integer $m$, the combinatoric
function $\bi{m}{z}$ with respect to integer $z$ takes value {\small $\f{ \Ga( m + 1) } { \Ga( z + 1) \Ga (m- z + 1) }$} for $z \leq m$ and
value $0$ otherwise.  The left limit as $\ep$ tends to $0$ is denoted as $\lim_{\ep \downarrow 0}$. The notation ``$\LRA$'' means ``if and only
if''.  The expectation of a random variable is denoted by $\bb{E}[.]$.  The notation $\Pr \{ . \mid \se \}$ denotes the probability of an event
associated with random samples $X_1, X_2, \cd$ parameterized by $\se \in \Se$, where $\se$ may be dropped if it can be done without introducing
confusion.  The parameter $\se$ in $\Pr \{ . \mid \se \}$  may be dropped whenever this can be done without introducing confusion. The
cumulative distribution function of a Gaussian random variable is denoted by $\Phi(.)$.  For $\al \in (0, 1)$, let $\mcal{Z}_\al$ and $t_{n,
\al}$ denote, respectively, the $100(1 - \al) \%$ percentiles of a standard normal distribution and a Student $t$-distribution of $n$ degrees of
freedom. For $\al \in (0, 1)$, let $\chi_{n, \al}^2$ denote the $100 \al \%$ percentile of a chi-square distribution of $n$ degrees of freedom.
In the presentation of our sampling schemes, we need to use the following functions: {\small \bee &  & S_{\mrm{B}} (k, n, \se) = \bec \sum_{i =
0}^k \bi{n}{i} \se^i (1 -
\se)^{n - i} & \tx{for} \; \se \in [0, 1],\\
1 & \tx{for} \; \se < 0,\\
0 & \tx{for} \; \se > 1   \eec \\
&  & S_N (k, n, \se) = \bec \sum_{i= 0 }^k {\se N \choose i} {N -
\se N
\choose n - i  }  \slash {N \choose n}  & \tx{for} \; \se \in \{ \f{m}{N}: m = 0, 1, \cd, N \},\\
1 & \tx{for} \; \se  < 0,\\
0 & \tx{for} \; \se  > 1  \eec \\
&  & S_{\mrm{P}} (k, \se ) = \bec \sum_{i = 0}^k \f{ \se^i e^{-  \se}} { i! }  & \tx{for} \; \se \geq 0,\\
0 & \tx{for} \; \se < 0  \eec \\
&  & \mscr{M}(z, \se) = \bec \f{ 9 (z - \se)^2 } {2 \li ( z + 2 \se
\ri ) \li ( z + 2 \se - 3 \ri ) } & \tx{for}
\; 0 \leq z \leq 1 \; \tx{and} \;  \se \in (0, 1),\\
- \iy & \tx{for} \; 0 \leq z \leq 1 \; \tx{and} \; \se \notin (0, 1)
\eec \\
&   &  \mscr{M}_{\mrm{B}} (z,\se) = \bec z \ln \f{\se}{z} + (1 - z)
\ln \f{1 - \se}{1 - z} &
\tx{for} \; z \in (0,1) \; \tx{and} \; \se \in (0, 1),\\
\ln(1-\se) & \tx{for} \; z = 0 \; \tx{and} \; \se \in (0, 1),\\
\ln \se &  \tx{for} \; z = 1 \; \tx{and} \; \se \in (0, 1),\\
- \iy &  \tx{for} \; z \in [0, 1] \; \tx{and} \; \se \notin (0, 1)
\eec\\
&  & \mscr{M}_{\mrm{I}} (z,\se) = \bec  \ln \f{\se}{z} + \li (
\f{1}{z} - 1 \ri ) \ln \f{1 - \se}{1 - z} &
\tx{for} \; z \in (0,1) \; \tx{and} \; \se \in (0, 1),\\
\ln \se &  \tx{for} \; z = 1 \; \tx{and} \; \se \in (0, 1),\\
- \iy & \tx{for} \; z = 0 \; \tx{and} \; \se \in (0, 1),\\
- \iy &  \tx{for} \; z \in [0, 1] \; \tx{and} \; \se \notin (0, 1) \eec\\
&  & \mscr{M}_{\mrm{P}} (z, \se) =
\bec z - \se + z \ln \li ( \f{\se}{z} \ri ) & \tx{for} \; z > 0 \; \tx{and} \;  \se > 0,\\
- \se & \tx{for} \; z = 0 \; \tx{and} \; \se > 0,\\
- \iy & \tx{for} \; z \geq 0 \; \tx{and} \; \se \leq 0.  \eec \eee}
In the design of multistage sampling schemes, one of our methods for
defining sample sizes is to use a descending sequence $C_\ell, \;
\ell \in \bb{Z}$ such that $C_0 = 1$ and $1 < \inf_{\ell \in \bb{Z}
} \f{C_\ell}{C_{\ell + 1}} \leq \sup_{\ell \in \bb{Z} }
\f{C_\ell}{C_{\ell + 1}} < \iy$.  Throughout the remainder of this
paper, $\de$ and $\ze$ are reserved, respectively, for the
``confidence parameter'' and the ``coverage tuning parameter'',
where these concepts will be illustrated later. It is assumed that
$0 < \de < 1$ and $0 < \ze < \f{1}{\de}$.   The other notations will
be made clear as we proceed.

\sect{General Theory}

In this section, we shall discuss the general theory of multistage
estimation.  A central theme of our theory is on the reduction of
the computational complexity associated with the design and analysis
of multistage sampling schemes.

\subsection{Basic Structure of Multistage Estimation} \la{gen_structure}

In our proposed framework of multistage estimation, a sampling
process consists of $s$ stages, where $s$ can be a finite number or
infinity. The continuation or termination of sampling is determined
by decision variables.  For the $\ell$-th stage, a decision variable
$\bs{D}_\ell = \mscr{D}_\ell (X_1, \cd, X_{\mbf{n}_\ell})$ is
defined in terms of samples $X_1, \cd, X_{\mbf{n}_\ell}$, where
$\mbf{n}_\ell$ is the number of samples available at the $\ell$-th
stage.  It should be noted that $\mbf{n}_\ell$ can be a random
number, depending on specific sampling schemes. The decision
variable $\bs{D}_\ell$ assumes only two possible values $0, \; 1$
with the notion that the sampling process is continued until
$\bs{D}_\ell = 1$ for some $\ell \in \{1, \cd, s\}$. Since the
sampling must be terminated at or before the $s$-th stage, it is
required that $\bs{D}_s = 1$. For simplicity of notations, we also
define $\bs{D}_\ell = 0$ for $\ell < 1$ and $\bs{D}_\ell = 1$ for
$\ell > s$ throughout the remainder of the paper. Let $\bs{l}$
denote the index of stage when the sampling is terminated.  Then,
the sample number when the sampling is terminated, denoted by
$\mbf{n}$, is equal to $\mbf{n}_{\bs{l}}$. Since  a sampling scheme
with the above structure is like a multistage version of the
conventional fixed-size sampling procedure, we call it {\it
multistage sampling} in this paper.

As mentioned earlier, the number of available samples,
$\mbf{n}_\ell$, for the $\ell$-th stage can be a random number.  An
important case can be made in the estimation of the parameter of a
Bernoulli random variable $X$ with distribution $\Pr \{ X = 1 \} = 1
- \Pr \{ X = 0 \} = p \in (0, 1)$. To estimate $p$, we can choose a
sequence of positive integers $\ga_1 < \ga_2 < \cd < \ga_s$ and
define decision variables such that $\bs{D}_\ell$ is expressed in
terms of i.i.d. samples $X_1, \cd, X_{\mbf{n}_\ell}$ of Bernoulli
random variable $X$,  where $\mathbf{n}_\ell$ is the minimum integer
such that $\sum_{i = 1}^{\mathbf{n}_\ell} X_i = \ga_\ell$ for $\ell
= 1, \cd, s$.  A sampling scheme with such a structure is called a
{\it multistage inverse binomial sampling}, which is a special class
of multistage sampling schemes and is a multistage version of the
inverse binomial sampling (see, e.g., \cite{H, H2} and the
references therein).

If the sample sizes of a multistage sampling scheme is desired to be
deterministic, the following criteria can be applied to determine
the minimum and maximum sample sizes:

(I) The minimum sample size $n_1$ guarantees that $\{ \bs{D}_1 = 1
\}$ is not an impossible event.

(II) The maximum sample size $n_s$ guarantees that $\{ \bs{D}_s = 1
\}$ is a sure event.

For the purpose of reducing sample number, the minimum and maximum
sample sizes should be as small as possible, while satisfying
criteria (I) and (II).  Once the minimum and maximum sample sizes
are fixed, the sample sizes for other stages can be determined, for
example,  as an arithmetic or geometric progression.

\subsection{Truncated Inverse Sampling}

It should be noted that the conventional single stage sampling
procedures can be accommodated in the general framework of
multistage sampling.  A common stopping rule for single stage
sampling procedures is that ``the sampling is continued until the
sample sum reach a prescribed threshold $\ga$ or the number of
samples reach a pre-specified integer $m$''.  Such a sampling scheme
is referred to as {\it truncated inverse sampling}, for which we
have derived the following results.

 \beT
\la{Truncated_sampling}
 Let $\ga > 1, \; 0 < \vep_a < \vep_r < 1$ and $p^\star =
\f{\vep_a}{\vep_r}$. Let $X_1, \; X_2, \; \cd$ be a sequence of
i.i.d. random variables such that $0 \leq X_i \leq 1$ and
$\bb{E}[X_i] = \mu \in (0, 1)$ for any positive integer $i$.  Let
$\mbf{n}$ be a random variable such that {\small $\li \{ \sum_{i
=1}^{\mbf{n} - 1} X_i < \ga \leq \sum_{i =1}^{\mbf{n}} X_i \ri \}$}
is a sure event. Let {\small $\mbf{m} = \min \{ \mbf{n}, m \}$},
where $m$ is a positive integer. The following statements hold true.

(I) {\small $\Pr \{  | \f{\ga}{ \mbf{n} } - \mu | < \vep \mu \}
> 1 - \de$} and {\small $\Pr \{ | \f{\ga -
1}{ \mbf{n} - 1 } - \mu | < \vep \mu \} > 1 - \de$} provided that
{\small $\ga > \f{ (1 + \vep) \ln (2 \sh \de) } { (1 + \vep) \ln (1
+ \vep) - \vep }$}.

(II) {\small $\Pr \li \{ | \f{\ga}{ \mbf{m} } - \mu | < \vep_a \;
\mrm{or} \; | \f{\ga}{ \mbf{m} } - \mu | < \vep_r \mu \ri \}
> 1 - \de$} provided that {\small $p^\star + \vep_a \leq \f{1}{2},
\; \ga > \f{1 - \vep_r}{\vep_r}, \; \ga
> \f{ \ln (\de \sh 2) } { \mscr{M}_{\mrm{I}}  \li ( \f{\ga (p^\star -
\vep_a)}{\ga - 1 + \vep_r}, p^\star \ri )  }, \; \ga > \f{ \ln (\de
\sh 2) } { \mscr{M}_{\mrm{I}}  ( p^\star + \vep_a, p^\star ) }$} and
{\small $m > \f{ \ln (\de \sh 2) } { \mscr{M}_{\mrm{B}} ( p^\star +
\vep_a, p^\star )  }$}.

(III) If $X_1, X_2, \cd$ are i.i.d. Bernoulli variables,  then
{\small $\Pr \li \{ | \f{\ga}{ \mbf{m} } - \mu | < \vep_a \;
\mrm{or} \; | \f{\ga}{ \mbf{m} } - \mu | < \vep_r \mu \ri \} > 1 -
\de$} provided that {\small $p^\star + \vep_a \leq \f{1}{2}, \; \ga
> \f{ \ln (\de \sh 2) } { \mscr{M}_{\mrm{I}} ( p^\star + \vep_a,
p^\star ) }$} and {\small $m > \f{ \ln (\de \sh 2) } {
\mscr{M}_{\mrm{B}} ( p^\star + \vep_a, p^\star )  }$}.  \eeT

The proof of Theorem \ref{Truncated_sampling} can be found in
\cite{Chen_TIV, Chen_IV}.

\subsection{Sequential Random Intervals}

A primary goal of multistage sampling is to construct, based on samples of $X$,  a sequential random interval with lower limit {\small $\mscr{L}
(X_1, \cd, X_{\mbf{n}})$} and upper limit {\small $\mscr{U} (X_1, \cd, X_{\mbf{n}})$} such that, for {\it a priori } specified confidence
parameter $\de$,
\[
\Pr \{ \mscr{L} (X_1, \cd, X_{\mbf{n}}) < \se < \mscr{U} (X_1, \cd,
X_{\mbf{n}}) \mid \se \} \geq 1 - \de \] for any $\se \in \Se$.  For
the $\ell$-th stage, an estimator $\wh{\bs{\se}}_\ell$ for $\se$ can
be defined in terms of samples $X_1, \cd, X_{\mbf{n}_\ell}$.
Consequently, the overall estimator for $\se$, denoted by
$\wh{\bs{\se}}$, is equal to $\wh{\bs{\se}}_{\bs{l}}$.  In many
cases, $\mscr{L} (X_1, \cd, X_{\mbf{n}_\ell})$ and $\mscr{U} (X_1,
\cd, X_{\mbf{n}_\ell})$ can be expressed as a function of
$\wh{\bs{\se}}_\ell$ and $\mbf{n}_\ell$. For simplicity of
notations, we abbreviate $\mscr{L} (X_1, \cd, X_{\mbf{n}_\ell})$ and
$\mscr{U} (X_1, \cd, X_{\mbf{n}_\ell})$ as $\mscr{L}
(\wh{\bs{\se}}_\ell, \mbf{n}_\ell)$ and $\mscr{U}
(\wh{\bs{\se}}_\ell, \mbf{n}_\ell)$ respectively. Accordingly,
$\mscr{L} (X_1, \cd, X_{\mbf{n}})$ and $\mscr{U} (X_1, \cd,
X_{\mbf{n}})$ are abbreviated as $\mscr{L} (\wh{\bs{\se}}, \mbf{n})$
and $\mscr{U} (\wh{\bs{\se}}, \mbf{n})$.  In the special case that
the lower and upper limits are independent of $\mbf{n}$, we will
drop the argument $\mbf{n}$ for further simplification of notations.

In the sequel, we shall focus on the construction of sequential random intervals of lower limit $\mscr{L} (\wh{\bs{\se}}, \mbf{n})$ and upper
limit $\mscr{U} (\wh{\bs{\se}}, \mbf{n})$ such that $\Pr \{ \mscr{L} (\wh{\bs{\se}}, \mbf{n})  < \se < \mscr{U} (\wh{\bs{\se}}, \mbf{n}) \mid
\se \} \geq 1 - \de$ for any $\se \in \Se$. Such a framework is general enough to address a wide spectrum of traditional problems in parametric
estimation.  First, it is obvious that the problem of interval estimation based on a given sampling scheme can be cast in this framework.
Second, the issue of error control in the point estimation of parameter $\se$ can be addressed in the framework of sequential random intervals.
Let $\udl{\se}$ and $\ovl{\se}$ be two numbers such that $\inf \Se \leq \udl{\se} \leq \ovl{\se} \leq \sup \Se$.  Let $\varTheta = \{ \se \in
\Se: \udl{\se} \leq \se \leq \ovl{\se} \}$. Based on different error criteria, the point estimation problems are typically posed in the
following ways:

(i) Given {\it a priori} margin of absolute error $\vep > 0$,
construct an estimator $\wh{\bs{\se}}$ for $\se$ such that \be
\la{abs888} \Pr \{ | \wh{\bs{\se}} - \se | < \vep \mid \se \}
> 1 - \de \qqu  \tx{for any $\se \in \varTheta$}. \ee

(ii) Given {\it a priori} margin of relative error $\vep \in (0,
1)$, construct an estimator $\wh{\bs{\se}}$ for $\se$ such that \be
\la{rev888} \Pr \{ | \wh{\bs{\se}} - \se | < \vep |\se| \} > 1 - \de
\qqu \tx{ for any $\se \in \varTheta$}. \ee

(iii) Given {\it a priori} margin of absolute error $\vep_a \geq 0$
and margin of relative error $\vep_r \in [0, 1)$, construct an
estimator $\wh{\bs{\se}}$ for $\se$ such that \be \la{mix999}
 \Pr \{ | \wh{\bs{\se}} - \se | < \vep_a \; \tx{or} \;  |
\wh{\bs{\se}} - \se | < \vep_r |\se| \mid \se \} > 1 - \de \qqu
\tx{for any $\se \in \varTheta$}. \ee

Clearly, problem (iii) can be reduced to problems (i) and (ii) by,
respectively,  setting $\vep_r = 0$ and $\vep_a = 0$.  As can be
seen from Appendix \ref{App_RI_Indentity}, putting \bee &  &
\mscr{L} (\wh{\bs{\se}}) = \bec \min \li \{ \wh{\bs{\se}} - \vep_a,
\; \f{ \wh{\bs{\se}} } { 1 + \vep_r } \ri \}
& \tx{if $\udl{\se} > 0$},\\
\min \li \{ \wh{\bs{\se}} - \vep_a, \; \f{ \wh{\bs{\se}} } { 1 -
\vep_r } \ri \}
& \tx{if $\ovl{\se} < 0$},\\
\min \li \{ \wh{\bs{\se}} - \vep_a, \; \f{ \wh{\bs{\se}} } { 1 +
\mrm{sgn} (\wh{\bs{\se}}) \; \vep_r } \ri \} & \tx{if $0 \in [
\udl{\se}, \ovl{\se}]$ } \eec\\
&  & \mscr{U} (\wh{\bs{\se}}) = \bec \max \li \{ \wh{\bs{\se}} +
\vep_a, \; \f{ \wh{\bs{\se}} } { 1 - \vep_r } \ri \}
& \tx{if $\udl{\se} > 0$},\\
\max \li \{ \wh{\bs{\se}} + \vep_a, \; \f{ \wh{\bs{\se}} } { 1 +
\vep_r } \ri \}
& \tx{if $\ovl{\se} < 0$},\\
\max \li \{ \wh{\bs{\se}} + \vep_a, \; \f{ \wh{\bs{\se}} } { 1 -
\mrm{sgn} (\wh{\bs{\se}}) \; \vep_r } \ri \} & \tx{if $0 \in [
\udl{\se}, \ovl{\se}]$ } \eec
 \eee
we can show that \be \la{RI_Indentity} \{ | \wh{\bs{\se}} - \se | < \vep_a \; \tx{or} \; | \wh{\bs{\se}} - \se | < \vep_r |\se| \}  = \{
\mscr{L} (\wh{\bs{\se}})  < \se < \mscr{U} (\wh{\bs{\se}}) \}. \ee This implies that problems (i)-(iii) can be accommodated in the general
framework of sequential random intervals.

Third, the framework of sequential random intervals accommodates an important class of problems concerned with the construction of bounded-width
confidence intervals. The objective is to construct lower confidence limit $\mscr{L} (\wh{\bs{\se}}, \mbf{n})$ and upper confidence limit
$\mscr{U} (\wh{\bs{\se}}, \mbf{n})$ such that $|\mscr{U} (\wh{\bs{\se}}, \mbf{n}) - \mscr{L} (\wh{\bs{\se}}, \mbf{n}) | \leq 2 \vep$ for some
prescribed number $\vep > 0$ and that $\Pr \{ \mscr{L} (\wh{\bs{\se}}, \mbf{n}) < \se < \mscr{U} (\wh{\bs{\se}}, \mbf{n}) \mid \se \} \geq 1 -
\de$ for any $\se \in \Se$. Obviously, this class of problems can be cast into the framework of sequential random intervals.

In order to construct a sequential random interval of desired level of confidence, our global strategy is to construct a sampling scheme such
that the coverage probability $\Pr \{ \mscr{L} (\wh{\bs{\se}}, \mbf{n}) < \se < \mscr{U} (\wh{\bs{\se}}, \mbf{n}) \mid \se \}$ can be adjusted
by some parameter $\ze$.  This parameter $\ze$ is referred to as ``coverage tuning parameter''.   Obviously, the coverage probability is a
function of the unknown parameter $\se$. In practice, it is impossible or extremely difficult to evaluate the coverage probability for every
value of $\se$ in the parameter space.  Such an issue presents in the estimation of binomial parameters, Poisson parameters and the proportion
of a finite population.  For the cases of estimating binomial and Poisson parameters, the parameter spaces are continuous and thus the number of
parametric values is infinity. For the case of estimating the proportion of a finite population, the number of parametric values can be as large
as the population size.  To overcome the difficulty associated with the number of parametric values, we have developed a general theory of
coverage probability of sequential random intervals which eliminates the need of exhaustive evaluation of coverage probabilities to determine
whether the minimum coverage probability achieves the desired level of confidence.  In this direction, the concept of {\it Unimodal-Likelihood
Estimator}, to be discussed in the following subsection,  play a crucial role.

\subsection{Unimodal-Likelihood Estimator}

The concept of maximum-likelihood estimator (MLE) is classical and widely used in numerous areas. However,  a MLE may not be unbiased and its
associated likelihood function need not be monotone.  For the purpose of developing a rigorous theory on coverage probability of random
intervals, we shall introduce the concept of {\it unimodal-likelihood estimator} (ULE) in this paper.  Let $X_1, X_2, \cd$ be a sequence of
random samples parameterized by $\se \in \Se$. Let $\mbf{m}$ be a positive integer-valued random variable such that for any positive integer
$m$, event $\{\mbf{m} = m \}$ depends only on $X_1, \cd, X_m$.  In a rigorous probabilistic terminology, $\mbf{m}$ is a stopping time.  For
samples $X_1, \cd, X_{\mbf{m}}$ of random length $\mbf{m}$, we say that the estimator $\varphi (X_1, \cd, X_{\mbf{m}})$ is a ULE of $\se$ if
$\varphi$ is a multivariate function such that, for any observation $(x_1, \cd, x_{m})$ of $(X_1, \cd, X_{\mbf{m}})$, the likelihood function is
non-decreasing with respect to $\se$ no greater than $\varphi (x_1, \cd, x_{m})$ and is non-increasing with respect to $\se$ no less than
$\varphi (x_1, \cd, x_{m})$. For discrete random samples $X_1, \cd, X_{m}$, the associated likelihood function is $\Pr \{ X_i = x_i, \; i = 1,
\cd, m \mid \se \}$. For continuous random samples $X_1, \cd, X_{m}$, the corresponding likelihood function is, $f_{X_1, \cd, X_m} (x_1, \cd,
x_m, \se)$, the joint probability density function of random samples $X_1, \cd, X_m$.  We emphasize that a MLE may not be a ULE and that a ULE
may not be a MLE.  In contrast to a MLE,  a ULE can assume values not contained in the parameter space.

Clearly,  for the cases that $X$ is a Bernoulli or Poisson variable,
{\small $\varphi (X_1, \cd, X_{\mbf{n}_\ell}) = \f{ \sum_{i =
1}^{\mbf{n}_\ell} X_i } {\mbf{n}_\ell}$} is a ULE of $\se$ at the
$\ell$-th stage.  As another illustration of ULE, consider the
multistage inverse binomial sampling scheme described in Section
\ref{gen_structure}. For $\ell = 1, \cd s$, a ULE of $p$ can be
defined as $\wh{\bs{p}}_\ell = \f{\ga_\ell} { \mathbf{n}_\ell }$. At
the termination of sampling, the estimator, $\wh{\bs{p}} =
\wh{\bs{p}}_{\bs{l}}$, of $p$ is also a ULE.  More generally, if the
distribution of $X$ belongs to the exponential family, then {\small
$\varphi (X_1, \cd, X_{\mbf{n}_\ell}) = \f{ \sum_{i =
1}^{\mbf{n}_\ell} X_i } {\mbf{n}_\ell}$} is a ULE of $\se$ at the
$\ell$-th stage.

\subsection{Inclusion Principle for Construction of Sampling Schemes}

In this subsection, we shall discuss the fundamental principle for
the design of multistage sampling schemes.  We shall address two
critical problems:

(I) Determine sufficient conditions for a multistage sampling scheme
such that the coverage probability $\Pr \{ \mscr{L} (\wh{\bs{\se}},
\mbf{n})  < \se < \mscr{U} (\wh{\bs{\se}}, \mbf{n}) \mid \se \}$ can
be adjusted by  a positive number $\ze$.

(II)  For a given sampling scheme, determine whether the coverage
probability $\Pr \{ \mscr{L} (\wh{\bs{\se}}, \mbf{n}) < \se <
\mscr{U} (\wh{\bs{\se}}, \mbf{n}) \mid \se \}$ is no less than $1 -
\de$ for any $\se \in \Se$.

To describe our sampling schemes, define cumulative distribution
function (CDF) and complementary cumulative distribution function
(CCDF) respectively as {\small \[ \la{defFG}
 F_{\wh{\bs{\se}}_\ell} (z, \se) = \bec \Pr \{
\wh{\bs{\se}}_\ell
\leq z \mid \se \} & \tx{for} \; \se \in \Se,\\
1  & \tx{for} \; \se < \inf \Se,\\
0  & \tx{for} \; \se > \sup \Se \eec  \qqu \qqu
G_{\wh{\bs{\se}}_\ell} (z, \se) = \bec \Pr \{ \wh{\bs{\se}}_\ell
\geq z \mid \se \} & \tx{for} \; \se \in \Se,\\
0  & \tx{for} \; \se < \inf \Se,\\
1  & \tx{for} \; \se > \sup \Se \eec \]} for $\ell = 1, \cd, s$,
where $z$ assumes values in the support of $\wh{\bs{\se}}_\ell$.

Let $\de_\ell \in (0, \f{1}{\ze}), \; \ell = 1, \cd, s$.  For
sampling schemes of structure described in Section
\ref{gen_structure}, we have the following results on the coverage
probability of random intervals.

\beT \la{Monotone_second}

Suppose that a multistage sampling scheme satisfies the following
requirements:

(i) For $\ell = 1, \cd, s$,  $\wh{\bs{\se}}_\ell$ is a ULE of $\se$.

(ii) For $\ell = 1, \cd, s$, $\{ \mscr{L}  (\wh{\bs{\se}}_\ell,
\mbf{n}_\ell) \leq \wh{\bs{\se}}_\ell \leq \mscr{U}
(\wh{\bs{\se}}_\ell, \mbf{n}_\ell) \}$ is a sure event.

(iii) {\small $\li \{ \bs{D}_\ell = 1 \ri \} \subseteq \li \{
F_{\wh{\bs{\se}}_\ell} \li ( \wh{\bs{\se}}_\ell, \mscr{U}
(\wh{\bs{\se}}_\ell, \mbf{n}_\ell) \ri ) \leq \ze \de_\ell , \;
G_{\wh{\bs{\se}}_\ell} \li ( \wh{\bs{\se}}_\ell, \mscr{L}
(\wh{\bs{\se}}_\ell, \mbf{n}_\ell) \ri ) \leq \ze \de_\ell \ri \}$}
for $\ell = 1, \cd, s$.

(iv) $\{ \bs{D}_s = 1 \}$ is a sure event.

Then, \bee &   & \Pr \{ \mscr{L} (\wh{\bs{\se}}, \mbf{n})  \geq \se
\mid \se \} \leq \sum_{\ell = 1}^s \Pr \{ \mscr{L}
(\wh{\bs{\se}}_\ell, \mbf{n}_\ell) \geq \se, \; \bs{D}_\ell = 1 \mid
\se
\} \leq \ze \sum_{\ell = 1}^s \de_\ell,\\
&   & \Pr \{ \mscr{U} (\wh{\bs{\se}}, \mbf{n})  \leq \se \mid \se \}
\leq \sum_{\ell = 1}^s \Pr \{ \mscr{U} (\wh{\bs{\se}}_\ell,
\mbf{n}_\ell) \leq \se, \; \bs{D}_\ell = 1 \mid \se \} \leq \ze
\sum_{\ell = 1}^s
\de_\ell,\\
&  & \Pr \{ \mscr{L} (\wh{\bs{\se}}, \mbf{n})  < \se < \mscr{U}
(\wh{\bs{\se}}, \mbf{n}) \mid \se \} \geq 1 - 2 \ze \sum_{\ell =
1}^s \de_\ell \eee for any $\se \in \Se$.

 \eeT

See Appendix \ref{App_Monotone_second} for a proof.   Theorem \ref{Monotone_second} addresses the first problem posed at the beginning of this
subsection. It tells how to define a stopping rule such that the coverage probability of the sequential random interval can be bounded by a
function of $\ze$ and $\sum_{\ell = 1}^s \de_\ell$. If $\sum_{\ell = 1}^s \de_\ell$ is bounded with respect to $\ze$, then the coverage
probability can be ``tuned'' to be no less than the prescribed level $1 - \de$. This process is referred to as ``coverage tuning'', which will
be illustrated in details in the sequel.  The intuition behind the definition of the stopping rule in Theorem \ref{Monotone_second} is as
follows.

 At the $\ell$-th stage, in order to determine whether the sampling should stop,
 two tests are performed based on the observations of
$\wh{\bs{\se}}_\ell, \; \mscr{L} (\wh{\bs{\se}}_\ell, \mbf{n}_\ell)$
and $\mscr{U} (\wh{\bs{\se}}_\ell, \mbf{n}_\ell)$, which are denoted
by $\vse_\ell, \; L_\ell$ and $U_\ell$ respectively.  The first test
is  $\mscr{H}_0: \se < U_\ell$ versus $\mscr{H}_1: \se \geq U_\ell$,
and the second test is
 $\mscr{H}_0^\prime: \se \leq L_\ell$ versus $\mscr{H}_1^\prime: \se > L_\ell$.
 Hypothesis $\mscr{H}_0$ is accepted if $F_{\wh{\bs{\se}}_\ell} (\vse_\ell, U_\ell) \leq \ze
\de_\ell$, and is rejected otherwise.  On the other side, hypothesis
$\mscr{H}_0^\prime$ is rejected if $G_{\wh{\bs{\se}}_\ell} (
\vse_\ell,  L_\ell ) \leq \ze \de_\ell$, and is accepted otherwise.
If $\mscr{H}_0$ is accepted and $\mscr{H}_0^\prime$ is rejected,
then, the decision variable $\bs{D}_\ell$ assumes value $1$ and
accordingly the sampling is terminated. Otherwise, $\bs{D}_\ell$
assumes value $0$ and the sampling is continued.  It can be seen
that, if $\ze \de_\ell$ is small, then $\mscr{H}_0$ and
$\mscr{H}_1^\prime$ are accepted with high credibility and
consequently, $L_\ell < \se < U_\ell$ is highly likely to be true.
Therefore, by making $\ze \sum_{\ell = 1}^s \de_\ell$ sufficiently
small, it is possible to ensure that the coverage probability of the
random interval is above the desired level.

Since there is a close relationship between hypothesis testing and confidence intervals, it is natural to imagine that the method described by
Theorem \ref{Monotone_second} for defining stopping rules to control the coverage probabilities of sequential random intervals  can be
interpreted with the concept of confidence intervals.  Since $\wh{\bs{\se}}_\ell$ is a ULE of $\se$ for $\ell = 1, \cd, s$, it follows from
Lemma \ref{ULE_Basic} in Appendix \ref{App_Basic_ULE} that $F_{\wh{\bs{\se}}_\ell} (z, \se)$ is non-increasing with respect to $\se \in \Se$ no
less than $z \in I_{\wh{\bs{\se}}_\ell }$ and that $G_{\wh{\bs{\se}}_\ell} (z, \se)$ is non-decreasing with respect to $\se \in \Se$ no greater
than $z \in I_{\wh{\bs{\se}}_\ell }$. Therefore, for the $\ell$-th stage, we can construct lower confidence limit $\mcal{L}_\ell (
\wh{\bs{\se}}_\ell, \mbf{n}_\ell, \ze \de_\ell)$ and upper confidence limit $\mcal{U}_\ell (\wh{\bs{\se}}_\ell, \mbf{n}_\ell, \ze \de_\ell)$
such that \bel & & \mcal{L}_\ell (\wh{\bs{\se}}_\ell, \mbf{n}_\ell, \ze \de_\ell) = \sup \li \{ \se \in \Se: G_{\wh{\bs{\se}}_\ell}
(\wh{\bs{\se}}_\ell, \se) \leq \ze \de_\ell, \; \se \leq
\wh{\bs{\se}}_\ell \ri \}, \la{lowCI} \\
&  & \mcal{U}_\ell (\wh{\bs{\se}}_\ell, \mbf{n}_\ell, \ze \de_\ell) = \inf \li \{ \se \in \Se: F_{\wh{\bs{\se}}_\ell} (\wh{\bs{\se}}_\ell, \se)
\leq \ze \de_\ell, \; \se \geq \wh{\bs{\se}}_\ell \ri \}. \la{uppCI} \eel As a consequence of (\ref{lowCI}) and (\ref{uppCI}), we have that
\[
\Pr \{ \se \leq \mcal{L}_\ell (\wh{\bs{\se}}_\ell, \mbf{n}_\ell, \ze
\de_\ell) \mid \se \} \leq \ze \de_\ell, \qu \qqu \Pr \{ \se \geq
\mcal{U}_\ell (\wh{\bs{\se}}_\ell, \mbf{n}_\ell, \ze \de_\ell) \mid
\se \} \leq \ze \de_\ell, \] \[ \Pr \{ \mcal{L}_\ell
(\wh{\bs{\se}}_\ell, \mbf{n}_\ell, \ze \de_\ell) < \se <
\mcal{U}_\ell (\wh{\bs{\se}}_\ell, \mbf{n}_\ell, \ze \de_\ell) \mid
\se \} \geq 1 - 2 \ze \de_\ell,
\]
which implies that $\mcal{L}_\ell (\wh{\bs{\se}}_\ell, \mbf{n}_\ell, \ze \de_\ell)$ and $\mcal{U}_\ell (\wh{\bs{\se}}_\ell, \mbf{n}_\ell, \ze
\de_\ell)$ are confidence limits with coverage probabilities controllable by $\ze$. It should be noted that such confidence limits are not
necessarily fixed-sample-size confidence limits, since the sample size $\mbf{n}_\ell$ can be a random number. Due to the monotonicity of
functions $F_{\wh{\bs{\se}}_\ell}(.,.)$ and $G_{\wh{\bs{\se}}_\ell} (.,.)$, we have that the requirement (iii) of Theorem \ref{Monotone_second}
can be restated as {\small \bel &   &  \li \{ \bs{D}_\ell = 1 \ri \} \subseteq \li \{ F_{\wh{\bs{\se}}_\ell} \li ( \wh{\bs{\se}}_\ell, \mscr{U}
(\wh{\bs{\se}}_\ell, \mbf{n}_\ell) \ri ) \leq \ze \de_\ell , \; G_{\wh{\bs{\se}}_\ell} \li ( \wh{\bs{\se}}_\ell, \mscr{L}
(\wh{\bs{\se}}_\ell, \mbf{n}_\ell) \ri ) \leq \ze \de_\ell \ri \} \nonumber\\
&  &  = \{  \mscr{L} (\wh{\bs{\se}}_\ell, \mbf{n}_\ell) \leq \mcal{L}_\ell (\wh{\bs{\se}}_\ell, \mbf{n}_\ell, \ze \de_\ell) \leq \mcal{U}_\ell
(\wh{\bs{\se}}_\ell, \mbf{n}_\ell, \ze \de_\ell) \leq \mscr{U} (\wh{\bs{\se}}_\ell, \mbf{n}_\ell) \} \la{control} \eel} for $\ell = 1, \cd, s$.
Clearly, such a requirement is set for controlling the coverage probability of the sequential random interval.  Since $[ \mcal{L}_\ell
(\wh{\bs{\se}}_\ell, \mbf{n}_\ell, \ze \de_\ell), \; \mcal{U}_\ell (\wh{\bs{\se}}_\ell, \mbf{n}_\ell, \ze \de_\ell) ], \; \ell = 1, \cd, s$
constitutes a sequence of confidence intervals,  it is referred to as a {\it confidence sequence} by the convention of statistical terminology
\cite{darling2}.  As our purpose of using such a confidence sequence is to control the coverage probability of the desired sequential random
interval, we call it {\it controlling confidence sequence}.  We also hope to use this term to avoid the confusion between the sequence of
confidence intervals and the desired sequential random interval.    By virtue of (\ref{control}) and the concepts of sequential random interval
and controlling confidence sequence,  the method described by Theorem \ref{Monotone_second} for constructing sampling schemes can be interpreted
as follows:

\bsk

\begin{tabular} {|l |}
\hline $\tx{ {\bf A necessary condition for the sampling process to be terminated is that}}$\\
$\; \tx{{\bf the sequential random interval includes the controlling confidence sequence.} }$
\\ \hline
\end{tabular}

\bsk

In view of the inclusion relationship imposed by the stopping condition and the versatility of the approach, we call such a methodology of using
a confidence sequence to define a stopping rule to control the coverage probability of the desired sequential random interval as {\it Inclusion
Principle}.

In the above statement of the inclusion principle, no sufficient condition is posed for the termination of the sampling process.  Our
consideration is that the inclusion relationship is not sufficient for the sampling process to stop if there are some other requirements imposed
on the sequential random interval.  A familiar example is the construction of bounded-length confidence interval. For this problem, the
sequential random interval coincides with the controlling confidence sequence. So, the inclusion relationship is something automatic. The only
extra requirement is that the sequential random interval has a pre-specified length.

In situations that no other requirement imposed on the sequential random interval except the specification of coverage probability,  we have
proposed a more specific version of inclusion principle for constructing sampling schemes as follows:

\bsk

\begin{tabular} {|l |}
\hline $\tx{ {\bf The sampling process is continued until the controlling confidence sequence is}}$\\
$\; \tx{{\bf included by the sequential random interval at some stage.} }$
\\ \hline
\end{tabular}

 \bsk

As a consequence of the inclusion principle, the coverage probabilities of the sequential random intervals can be controlled by $\ze
> 0$ provided that the coverage probability of the controlling confidence sequence can be controlled by $\ze$.
More precisely, this connection can be seen from the following theorem.

\beT \la{Monotone_second_CI_ST}

Suppose that a multistage sampling scheme satisfies the following requirements:

(i) For $\ell = 1, \cd, s$,  $\mcal{L}_\ell (\wh{\bs{\se}}_\ell, \mbf{n}_\ell, \ze \de_\ell)$ and $\mcal{U}_\ell (\wh{\bs{\se}}_\ell,
\mbf{n}_\ell, \ze \de_\ell)$ are lower and upper limits such that $\Pr \{ \mcal{L}_\ell (\wh{\bs{\se}}_\ell, \mbf{n}_\ell, \ze \de_\ell) < \se <
\mcal{U}_\ell (\wh{\bs{\se}}_\ell, \mbf{n}_\ell, \ze \de_\ell) \mid \se \} \geq 1 - \ze \de_\ell$ for any $\se \in \Se$.

(ii) For $\ell = 1, \cd, s$, {\small $\{ \tx{The sampling process is terminated at the $\ell$-th stage} \}$} implies {\small $\{  \mscr{L}
(\wh{\bs{\se}}_\ell, \mbf{n}_\ell) \leq \mcal{L}_\ell (\wh{\bs{\se}}_\ell, \mbf{n}_\ell, \ze \de_\ell) \leq \mcal{U}_\ell (\wh{\bs{\se}}_\ell,
\mbf{n}_\ell, \ze \de_\ell) \leq \mscr{U} (\wh{\bs{\se}}_\ell, \mbf{n}_\ell) \}$ }.

(iii) $\{ \tx{The sampling process is terminated at some finite stage}  \}$ is a sure event.

Then, $\Pr \{ \mscr{L} (\wh{\bs{\se}}, \mbf{n})  < \se < \mscr{U} (\wh{\bs{\se}}, \mbf{n}) \mid \se \} \geq \Pr \{ \mcal{L}_\ell
(\wh{\bs{\se}}_\ell, \mbf{n}_\ell, \ze \de_\ell) < \se < \mcal{U}_\ell (\wh{\bs{\se}}_\ell, \mbf{n}_\ell, \ze \de_\ell) \; \tx{for} \; \ell = 1,
\cd, s \mid \se \} \geq 1 - \ze \sum_{\ell = 1}^s \de_\ell$ for any $\se \in \Se$.

 \eeT

It should be noted that Theorem \ref{Monotone_second_CI_ST} remains valid if the controlling confidence sequence and sequential random interval
are defined as more complex functions of samples instead of point estimators of $\se$ and sample sizes.  Theorem \ref{Monotone_second_CI_ST} is
an immediate consequence of the following probabilistic result.

\beT \la{General Inclusion Principle}

Let $(\Om, \mscr{F}, \{ \mscr{F}_\ell \}, \Pr )$ be a filtered space.  Let $\bs{\tau}$ be a proper stopping time with support $I_{\bs{\tau}}$.
 For $\ell \in I_{\bs{\tau}}$, let $\bs{\mcal{A}}_\ell \subseteq \bb{R}^n$ and
$\bs{\mcal{B}}_\ell \subseteq \bb{R}^n$ be $n$-dimensional random regions defined by random variables measurable in $\mscr{F}_\ell$.   Let
$\bs{\se} \in \bb{R}^n$ be an $n$-dimensional vector.  Assume that {\small $\{ \bs{\tau} = \ell \} \subseteq \{ \bs{\mcal{A}}_\ell  \subseteq
\bs{\mcal{B}}_\ell \}$ } for $\ell \in I_{\bs{\tau}}$.  Then, $\Pr \{ \bs{\se} \in \bs{\mcal{B}}_{\bs{\tau}} \} \geq \Pr \{ \bs{\se} \in
\bs{\mcal{A}}_\ell   \; \tx{for} \; \ell \in I_{\bs{\tau}} \} \geq 1 - \sum_{\ell \in I_{\bs{\tau}} } \Pr \{ \bs{\se} \notin \bs{\mcal{A}}_\ell
\}$. \eeT

The proof of Theorem \ref{General Inclusion Principle} is given in Appendix \ref{General Inclusion Principle_App}.  This theorem implies that
the inclusion principle can be applied to construct sampling schemes for estimating multiple parameters. In the general case, the sequential
random interval is generalized as sequential random region and accordingly, the controlling confidence sequence is generalized as a sequence of
confidence regions.  This method has been used in Section 17 for multistage linear region.

Actually, we have already extensively used the inclusion principle to derive stopping rules in the first version of this paper published in
arXiv on September 8, 2008. However, due to the simplification of the stopping rules, the link between stopping rules and controlling confidence
sequences is not obvious at the first glance, though it can be seen by a careful reading of the relevant proofs. In the first version of our
paper \cite{Chen_BP} published in October 2, 2008, we have derived stopping rules from which the connection between stopping rules and
controlling confidence sequences can be readily identified (see Theorem 1 and its proof in subsequent versions). About six months later, we have
proposed a systematic method of using a sequence of confidence intervals to define stopping rules to control coverage probabilities of
sequential random intervals in Section 3 of the fifth version of our paper \cite{Chen_rule} published in arXiv on April 7, 2009.

A fundamental fact disclosed by Theorem \ref{Monotone_second_CI_ST} is that {\it the coverage probability of the sequential random interval is
bounded from below by that of the controlling confidence sequence}.  This implies that  in the design of stopping rules, the coverage
probabilities of the sequential random intervals may still be controllable by $\ze$ if the confidence limits of an exact controlling confidence
sequence are replaced by their approximations or conservative bounds.

In situations that the parameter $\se$ to be estimated is the expectation of $X$, we can apply normal approximation to obtain confidence limits
for $\se$ as follows. Assume that $X_1, X_2, \cd$ are identical samples of $X$ and that the variance of the sample mean $\ovl{X}_n \DEF
\f{\sum_{i=1}^n X_i}{n}$ is a bivariate function, denoted by $\mscr{V} (\se, n)$, of $\se$ and $n$. If the sample size $n$ is large, then the
central limit theorem  may be applied to establish the normal approximation \bel &  & F_{\ovl{X}_n} \li ( z, \se \ri ) \DEF \Pr \{ \ovl{X}_n
\leq z \mid \se \} \ap \Phi \li ( \f{ z - \se  } { \sqrt{ \mscr{V} ( \se, n )
} } \ri ), \la{normap1}\\
&  & G_{\ovl{X}_n} \li ( z, \se  \ri ) \DEF \Pr \{ \ovl{X}_n \geq z \mid \se  \} \ap \Phi \li ( \f{ \se - z } { \sqrt{ \mscr{V} ( \se, n) } }
\ri ). \la{normap2} \eel Let $\al \in (0, 1)$. If $\ovl{X}_n$ is a ULE of $\se$, then we can obtain lower and upper confidence limits
respectively as
\[
\mcal{L} (\ovl{X}_n, n, \al) = \inf \li \{  \se \in \Se : \Phi \li ( \f{ \se - \ovl{X}_n } { \sqrt{ \mscr{V} ( \se, n) } } \ri ) > \f{\al}{2}
\ri \}
\]
and
\[
\mcal{U} (\ovl{X}_n, n, \al) = \sup \li \{  \se \in \Se : \Phi \li ( \f{ \ovl{X}_n - \se  } { \sqrt{ \mscr{V} ( \se, n ) } } \ri ) > \f{\al}{2}
\ri \}
\]
such that $\Pr \{ \mcal{L} (\ovl{X}_n, n, \al) < \se < \mcal{U} (\ovl{X}_n, n, \al) \mid \se \} \ap 1 - \al$ for $\se \in \Se$.  To improve the
accuracy of normal approximation (\ref{normap1}) and (\ref{normap2}), we propose to replace $\se$ in $\mscr{V} ( \se, n )$ by $z + \ro (\se -
z)$ with $\ro \in [0, 1]$. That is, we suggest modifying (\ref{normap1}) and (\ref{normap2}) as follows: \bel & & F_{\ovl{X}_n} \li ( z, \se \ri
) \ap
\Phi \li ( \f{ z - \se  } { \sqrt{ \mscr{V} ( z + \ro (\se - z), n ) } } \ri ), \la{normap3}\\
&  & G_{\ovl{X}_n} \li ( z, \se  \ri ) \ap \Phi \li ( \f{ \se - z } { \sqrt{ \mscr{V} ( z + \ro (\se - z), n) } } \ri ). \la{normap4} \eel The
new parameter $\ro$ is introduced to improve the accuracy of approximation.   Based on the new approximation (\ref{normap3}) and
(\ref{normap4}), we propose to obtain lower and upper confidence limits for $\se$ respectively as \be \la{best_CI_A}
 \mcal{L} (\ovl{X}_n, n, \al) = \inf
\li \{ \se \in \Se : \Phi \li ( \f{ \se - \ovl{X}_n } { \sqrt{ \mscr{V} ( \ovl{X}_n + \ro (\se - \ovl{X}_n), n) } } \ri ) > \f{\al}{2} \ri \}
\ee and \be \la{best_CI_B} \mcal{U} (\ovl{X}_n, n, \al) = \sup \li \{  \se \in \Se : \Phi \li ( \f{ \ovl{X}_n - \se  } { \sqrt{ \mscr{V} (
\ovl{X}_n + \ro (\se - \ovl{X}_n), n ) } } \ri ) > \f{\al}{2} \ri \} \ee so that $\Pr \{ \mcal{L} (\ovl{X}_n, n, \al) < \se < \mcal{U}
(\ovl{X}_n, n, \al) \mid \se \} \ap 1 - \al$ for $\se \in \Se$. Our computational experiences indicate that the coverage performance of this
class of  confidence intervals can be very close to the nominal level $1- \al$ by choosing appropriate $\ro \in (0, 1]$.  It should be noted
that for estimating parameters of binomial, Poisson, negative binomial, geometric and hypergeometric distributions, explicit confidence
intervals can be obtained from (\ref{best_CI_A}) and (\ref{best_CI_B}) by solving quadratic equations.  To illustrate, let $\mcal{Z}$ denote the
critical value such that $\Phi( \mcal{Z} ) = 1 - \f{\al}{2}$.   In scenarios that $X_1, \cd, X_n$ are i.i.d. samples of Bernoulli random
variable $X$ such that $\Pr \{ X = 1 \} = 1 - \Pr \{ X = 0 \} = p \in (0, 1)$, the lower and upper limits of the approximate $100 (1 - \al) \%$
confidence interval for $p$ can be readily derived from (\ref{best_CI_A}) and (\ref{best_CI_B}) respectively as \be \la{CI381} \mcal{L}
(\ovl{X}_n, n, \al) = \f{ \ovl{X}_n + \f{ \ro \mcal{Z}^2 }{ 2 n} [ 1 - 2 (1 - \ro) \ovl{X}_n ]  - \mcal{Z} \sq{ \f{ \ovl{X}_n ( 1 - \ovl{X}_n )
}{n} + \li ( \f{ \ro \mcal{Z}  }{ 2 n } \ri )^2 } }{ 1 + \f{ ( \ro \mcal{Z} )^2 }{n}  } \ee and \be \la{CI382}
 \mcal{U} (\ovl{X}_n, n, \al) =  \f{ \ovl{X}_n  +  \f{ \ro \mcal{Z}^2 }{ 2 n} [ 1 - 2 (1 - \ro)
\ovl{X}_n ]  + \mcal{Z} \sq{ \f{ \ovl{X}_n ( 1 - \ovl{X}_n ) }{n} + \li ( \f{ \ro \mcal{Z}  }{ 2 n } \ri )^2 } }{ 1 + \f{ ( \ro \mcal{Z} )^2
}{n}  }. \ee  It can be checked that taking $\ro = 0, \; 1$ and $\f{2}{3}$ respectively leads to Wald's interval, Wilson's interval and the
confidence interval proposed by Chen et. al. \cite{Chen_CI}. By the coverage theory of \cite{Chen_CI}, it can be shown that for $\ro =
\f{2}{3}$,
\[
\Pr \{  \mcal{L} (\ovl{X}_n, n, \al) < p < \mcal{U} (\ovl{X}_n, n, \al) \mid p \} \geq 1 - 2 \exp \li ( - \f{ \mcal{Z}^2 }{2} \ri )
\]
for any $p \in (0, 1)$.  Using approximation $\ovl{X}_n +  \f{ \ro \mcal{Z}^2 }{ 2 n}  [ 1 - 2 (1 - \ro) \ovl{X}_n ]  \ap \ovl{X}_n, \; 1 +  \f{
( \ro \mcal{Z} )^2 }{n} \ap 1$ and introducing parameter $\vro = \f{ (\ro \mcal{Z})^2 }{4}$, we can simplify (\ref{CI381}) and (\ref{CI382}) as
\[
\mcal{L} (\ovl{X}_n, n, \al) \ap \ovl{X}_n - \mcal{Z} \sq{ \f{ \ovl{X}_n ( 1 - \ovl{X}_n ) }{n} + \f{ \vro }{ n^2 }  }, \qqu \mcal{U}
(\ovl{X}_n, n, \al) \ap \ovl{X}_n + \mcal{Z} \sq{ \f{ \ovl{X}_n ( 1 - \ovl{X}_n ) }{n} + \f{ \vro }{ n^2 }  }.
\]

In the context that $X_1, \cd, X_n$ are i.i.d. samples of Poisson random variable $X$ of mean $\lm > 0$, we can apply (\ref{best_CI_A}) and
(\ref{best_CI_B}) to derive explicit lower and upper confidence limits for $\lm$ as
\[
\mcal{L} (\ovl{X}_n, n, \al) = \ovl{X}_n + \f{ \ro \mcal{Z}^2 } { 2 n } - \mcal{Z} \sq{ \f{ \ovl{X}_n } {n} +  \li (  \f{ \ro \mcal{Z} } { 2 n }
\ri )^2 }, \qqu \mcal{U} (\ovl{X}_n, n, \al) = \ovl{X}_n + \f{ \ro \mcal{Z}^2 } { 2 n } + \mcal{Z} \sq{ \f{ \ovl{X}_n } {n} +  \li (  \f{ \ro
\mcal{Z} } { 2 n } \ri )^2 }
\]
so that $\Pr \{ \mcal{L} (\ovl{X}_n, n, \al) < \lm < \mcal{U} (\ovl{X}_n, n, \al) \mid \lm \} \ap 1 - \al$ for all $\lm > 0$. Using
approximation $\ovl{X}_n + \f{ \ro \mcal{Z}^2 } { 2 n } \ap \ovl{X}_n$ and introducing parameter $\vro = \f{ \ro^2 \mcal{Z}^2 } { 4 } $, we have
 simplified confidence limits for $\lm$ as \[
\mcal{L} (\ovl{X}_n, n, \al) \ap \ovl{X}_n - \mcal{Z} \sq{ \f{ \ovl{X}_n } {n} +  \f{ \vro  } { n^2 } }, \qqu \mcal{U} (\ovl{X}_n, n, \al) \ap
\ovl{X}_n + \mcal{Z} \sq{ \f{ \ovl{X}_n } {n} + \f{ \vro } { n^2 } }.
\]

In the context that $X_1, \cd, X_n$ are i.i.d. samples of random variable $X \DEF Y + r$ with mean $\mu \DEF \f{r}{p}$, where $r
> 0, \; p \in (0, 1)$ are real numbers and $Y$ is a random variable possessing a negative binomial distribution, we can apply (\ref{best_CI_A}) and
(\ref{best_CI_B}) to derive explicit lower and upper confidence limits for $\mu$ as
\[
\mcal{L} (\ovl{X}_n, n, \al) = \f{ \ovl{X}_n - \f{ \ro \mcal{Z}^2 }{ 2 n r} [ r - 2 (1 - \ro) \ovl{X}_n ]  - \mcal{Z} \sq{ \f{ \ovl{X}_n (
\ovl{X}_n - r ) }{n r} + \li ( \f{ \ro \mcal{Z}  }{ 2 n } \ri )^2 } }{ 1 - \f{ ( \ro \mcal{Z} )^2 }{n r}  }
\]
\[
\mcal{U} (\ovl{X}_n, n, \al) = \f{ \ovl{X}_n -  \f{ \ro \mcal{Z}^2 }{ 2 n r} [ r - 2 (1 - \ro) \ovl{X}_n ]  + \mcal{Z} \sq{ \f{ \ovl{X}_n (
\ovl{X}_n - r ) }{n r} + \li ( \f{ \ro \mcal{Z}  }{ 2 n } \ri )^2 } }{ 1 - \f{ ( \ro \mcal{Z} )^2 }{n r}  }
\]
so that $\Pr \{ \mcal{L} (\ovl{X}_n, n, \al) < \mu < \mcal{U} (\ovl{X}_n, n, \al) \mid \mu \} \ap 1 - \al$ for all $\mu \in (1, \iy)$. Using
approximation $\ovl{X}_n - \f{ \ro \mcal{Z}^2 }{ 2 n r} [ r - 2 (1 - \ro) \ovl{X}_n ] \ap \ovl{X}_n, \; \;  1 - \f{ ( \ro \mcal{Z} )^2 }{n r}
\ap 1$ and introducing parameter $\vro = \f{ \ro^2 \mcal{Z}^2  }{ 4}$,  we can simplify the confidence limits for $\mu$ as
\[
\mcal{L} (\ovl{X}_n, n, \al) \ap \ovl{X}_n  - \mcal{Z} \sq{ \f{ \ovl{X}_n ( \ovl{X}_n - r ) }{n r} + \f{ \vro }{ n^2} }, \qqu \mcal{U}
(\ovl{X}_n, n, \al) \ap \ovl{X}_n + \mcal{Z} \sq{ \f{ \ovl{X}_n ( \ovl{X}_n - r ) }{n r} + \f{ \vro }{ n^2}  }.
\]

As the last example of constructing approximate confidence interval based on normal approximation, consider the procedure of sampling without
replacement from a population of $N$ units, among which there are $p N$ units having a certain attribute.  Let $n$ be the sample size. Let $X_1,
\cd, X_n$ be random variables associated with the sampling procedure in a way as follows:  $X_i$ assumes value $1$ if the $i$-th sample have the
attribute and assumes value $0$ otherwise.  Then, $X_1, \cd, X_n$ are identical but dependent random variables with mean $p \in \{0, \f{1}{N},
\cd, \f{N-1}{N}, 1 \}$.  The sample mean $\ovl{X}_n = \f{ \sum_{i = 1}^n X_i }{n}$ is unbiased.  Introducing modified sample size $n^* = \f{n(N
- 1)}{N - n}$ and making use of (\ref{best_CI_A}) and (\ref{best_CI_B}), we obtain explicit lower and upper confidence limits for $p$ as
\[
\mcal{L} (\ovl{X}_n, n, \al) =  \f{ \ovl{X}_n +  \f{ \ro \mcal{Z}^2 }{ 2 n^*} [ 1 - 2 (1 - \ro) \ovl{X}_n ]  - \mcal{Z} \sq{ \f{ \ovl{X}_n ( 1 -
\ovl{X}_n ) }{n^*} + \li ( \f{ \ro \mcal{Z} }{ 2 n^* } \ri )^2 } }{ 1 +  \f{ ( \ro \mcal{Z} )^2 }{n^*}  }
\]
\[
\mcal{U} (\ovl{X}_n, n, \al) = \f{ \ovl{X}_n +  \f{ \ro \mcal{Z}^2 }{ 2 n^*} [ 1 - 2 (1 - \ro) \ovl{X}_n ]   + \mcal{Z} \sq{ \f{ \ovl{X}_n ( 1 -
\ovl{X}_n ) }{n^*} + \li ( \f{ \ro \mcal{Z}  }{ 2 n^* } \ri )^2 } }{ 1 +  \f{ ( \ro \mcal{Z} )^2 }{n^*}  }.
\]
so that $\Pr \{ \mcal{L} (\ovl{X}_n, n, \al) \leq p \leq \mcal{U} (\ovl{X}_n, n, \al) \mid p \} \ap 1 - \al$ for all $p \in \{0, \f{1}{N}, \cd,
\f{N-1}{N}, 1 \}$.  Using approximation $\ovl{X}_n +  \f{ \ro \mcal{Z}^2 }{ 2 n^*} [ 1 - 2 (1 - \ro) \ovl{X}_n ] \ap \ovl{X}_n, \;\;  1 +  \f{ (
\ro \mcal{Z} )^2 }{n^*}  \ap 1$ and introducing new parameter $\f{ \ro^2 \mcal{Z}^2 }{ 4 }$,  we have simple confidence limits for $p$ as
\[
\mcal{L} (\ovl{X}_n, n, \al) \ap  \ovl{X}_n    - \mcal{Z} \sq{ \f{ \ovl{X}_n ( 1 - \ovl{X}_n ) }{n^*} + \f{ \vro }{ (n^*)^2 }  }, \qqu \mcal{U}
(\ovl{X}_n, n, \al) \ap \ovl{X}_n +   \mcal{Z} \sq{ \f{ \ovl{X}_n ( 1 - \ovl{X}_n ) }{n^*} + \f{ \vro }{ (n^*)^2 } }.
\]
From above examples, it can be seen that the lower and upper confidence limits of confidence level $1 - \al$ for a mean can be generally
constructed as
\[
\mcal{L} (\ovl{X}_n, n, \al) \ap \ovl{X}_n - \mcal{Z} \sq{  \mscr{V} ( \ovl{X}_n, n) + \f{\vro}{n^2}  }, \qqu \mcal{U} (\ovl{X}_n, n, \al) \ap
\ovl{X}_n + \mcal{Z} \sq{  \mscr{V} ( \ovl{X}_n, n) + \f{\vro}{n^2} },
\]
where the term $\f{\vro}{n^2}$ is introduced to remedy the inaccuracy of normal approximation.

The approximate confidence limits described above can be used to derive simple stopping rules by virtue of the inclusion principle. Although the
confidence limits are not rigorously constructed, for each
 fixed value of $\ro$, we can apply bisection coverage tuning technique to determine appropriate value of $\ze$ to guarantee the required
 confidence level for the desired sequential random interval.
 Stopping rules of excellent performance can be obtained by trying various values of $\ro$.

Although the stopping rules can be expressed in the terms of confidence limits by the inclusion principle,  we propose to eliminate the need of
computing confidence limits in order to make stopping rules as simple as possible.  In scenarios that the parameter $\se$ to be estimated is the
expectation of $X$, we can apply normal approximation to simplify the general stopping rule proposed by Theorem \ref{Monotone_second}. The
stopping rule described by Theorem \ref{Monotone_second} can be interpreted as ``sampling is continued until {\small $F_{\wh{\bs{\se}}_\ell} \li
( \wh{\bs{\se}}_\ell, \mscr{U} (\wh{\bs{\se}}_\ell, \mbf{n}_\ell) \ri ) \leq \ze \de_\ell , \; G_{\wh{\bs{\se}}_\ell} \li ( \wh{\bs{\se}}_\ell,
\mscr{L} (\wh{\bs{\se}}_\ell, \mbf{n}_\ell) \ri ) \leq \ze \de_\ell$} for some $\ell \in \{ 1, \cd, s \}$''.   Since $\se = \bb{E}[X]$, the
estimators $\wh{\bs{\se}}_\ell$ are naturally defined as sample means such that $\wh{\bs{\se}}_\ell = \f{ \sum_{i=1}^{\mbf{n}_\ell} X_i
}{\mbf{n}_\ell}$.  As before, assume that $X_1, X_2, \cd$ are identical samples of $X$ and that the variance of $\f{\sum_{i=1}^n X_i}{n}$ is a
bivariate function, denoted by $\mscr{V} (\se, n)$, of $\se$ and $n$. If all sample sizes are large, then we can apply the normal approximation
(\ref{normap1}) and (\ref{normap2}) to simplify the stopping rule described by Theorem \ref{Monotone_second} as ``sampling is continued until
\be \la{normalgen} \Phi \li ( \f{ \wh{\bs{\se}}_\ell - \mscr{U} (\wh{\bs{\se}}_\ell, \mbf{n}_\ell) } { \sqrt{ \mscr{V} ( \mscr{U}
(\wh{\bs{\se}}_\ell, \mbf{n}_\ell), \; \mbf{n}_\ell ) } } \ri ) \leq \ze \de_\ell, \qqu  \Phi \li ( \f{  \mscr{L} (\wh{\bs{\se}}_\ell,
\mbf{n}_\ell) - \wh{\bs{\se}}_\ell } { \sqrt{ \mscr{V} ( \mscr{L} (\wh{\bs{\se}}_\ell, \mbf{n}_\ell), \; \mbf{n}_\ell ) } } \ri ) \leq \ze
\de_\ell \ee for some $\ell \in \{ 1, \cd, s \}$''.  Moreover, for better performance of coverage probability, we can apply the normal
approximation (\ref{normap3}) and (\ref{normap4}) to simplify the stopping rule described by Theorem \ref{Monotone_second} as ``sampling is
continued until \bel  &  &  \Phi \li ( \f{ \wh{\bs{\se}}_\ell - \mscr{U} (\wh{\bs{\se}}_\ell, \mbf{n}_\ell) } { \sqrt{ \mscr{V} \li (
\wh{\bs{\se}}_\ell + \ro [ \mscr{U} (\wh{\bs{\se}}_\ell, \mbf{n}_\ell)
- \wh{\bs{\se}}_\ell ], \; \mbf{n}_\ell \ri ) } } \ri ) \leq \ze \de_\ell, \la{st38a}\\
&  & \Phi \li ( \f{  \mscr{L} (\wh{\bs{\se}}_\ell, \mbf{n}_\ell) - \wh{\bs{\se}}_\ell } { \sqrt{ \mscr{V} \li ( \wh{\bs{\se}}_\ell + \ro [
\mscr{L} ( \wh{\bs{\se}}_\ell, \mbf{n}_\ell) - \wh{\bs{\se}}_\ell ], \; \mbf{n}_\ell \ri ) } } \ri ) \leq \ze \de_\ell \la{st38b} \eel for some
$\ell \in \{ 1, \cd, s \}$''.   Our computational experiments show that, for a given $\ro \in (0, 1]$, the coverage probability can be
guaranteed by choosing $\ze$ to be a sufficiently small number.  Optimization of sampling schemes with respect to $\ro \in (0, 1]$ leads to
better performance. Clearly, this approach of constructing simple stopping rules applies to the problems of estimating parameters of binomial,
Poisson, geometric and hypergeometric distributions, etc.  More details are presented in the sequel.

\bed

\item [Estimating a Proportion:] Let $X$ be a Bernoulli random variable such that $\Pr \{ X = 1 \} = 1 - \Pr \{ X = 0 \} = p \in (0, 1)$.
Let $X_1, X_2, \cd$ be i.i.d. samples of $X$.   We consider the construction of multistage estimation procedures with deterministic sample sizes
$n_1, n_2, \cd, n_s$, where the number of stages $s$ can be finite or infinite.  For $\ell = 1, \cd, s$, the estimator for $p$ is defined as
$\wh{\bs{p}}_\ell = \f{ \sum_{i = 1}^{n_\ell} X_i }{n_\ell}$.   The sequential estimator $\wh{\bs{p}}$ is defined as $\wh{\bs{p}}_{\bs{l}}$,
where $\bs{l}$ is the index of stage at the termination of the sampling process.   The prescribed confidence interval is $1 - \de$, where $\de
\in (0, 1)$.  For $\ell = 1, 2, \cd, s$, let $\mcal{Z}_{\ze \de_\ell}$ be the critical value such that $\Phi( \mcal{Z}_{\ze \de_\ell} ) = 1 -
\ze \de_\ell$.

Let $\vep \in (0, 1)$ be a margin of absolute error. To design a multistage procedure such that the sequential estimator $\wh{\bs{p}}$
guarantees that $\Pr \{ | \wh{\bs{p}} - p| < \vep \mid p \} \geq 1 - \de$ for all $p \in (0, 1)$, we can apply the general stopping rule
associated with (\ref{st38a}) and (\ref{st38b}) by identifying the lower and upper bounds of the random interval as $\mscr{L} (\wh{\bs{p}}_\ell)
= \wh{\bs{p}}_\ell - \vep$ and $\mscr{U} (\wh{\bs{p}}_\ell) = \wh{\bs{p}}_\ell + \vep$ respectively to derive a simple stopping rule as follows:
Continue sampling until $\li ( \li | \wh{\bs{p}}_\ell - \f{1}{2} \ri | - \ro \vep \ri )^2 \geq \f{1}{4} - n_\ell \li ( \f{\vep}{\mcal{Z}_{\ze
\de_\ell}} \ri )^2$ for some $\ell \in \{1, 2, \cd, s\}$.  Consequently, using bisection coverage tuning techniques, we can obtain an
appropriate value of $\ze$ such that the coverage probability $\Pr \{ | \wh{\bs{p}} - p| < \vep \mid p \}$ is no less than but close to the
prescribed level $1 - \de$ for all $p \in (0, 1)$.

Let $\vep \in (0, 1)$ be a margin of relative error. To design a multistage procedure such that the resultant estimator $\wh{\bs{p}}$ guarantees
that $\Pr \{ | \wh{\bs{p}} - p| < \vep p \mid p \} \geq 1 - \de$ for all $p \in (0, 1)$, we can apply the general stopping rule by identifying
the lower and upper bounds of the random interval as $\mscr{L} (\wh{\bs{p}}_\ell) = \f{\wh{\bs{p}}_\ell}{1 + \vep}$ and  $\mscr{U}
(\wh{\bs{p}}_\ell) = \f{\wh{\bs{p}}_\ell}{1 - \vep}$ respectively to derive the following stopping rule: Continue sampling until
$\wh{\bs{p}}_\ell \geq \f{ (1 + \vep) ( 1 + \vep - \ro \vep ) \mcal{Z}_{\ze \de_\ell}^2 } { n_\ell \vep^2 + ( 1 + \vep - \ro \vep  )^2
\mcal{Z}_{\ze \de_\ell}^2   }$ for some $\ell \in \{1, 2, 3, \cd, \iy \}$.  The parameters $\de_\ell, \; \ell = 1, 2, \cd$ are positive numbers
less than $1$ such that their sum $\sum_{\ell = 1}^\iy \de_\ell$ is bounded.  By virtue of bisection coverage tuning techniques, we can obtain
an appropriate value of $\ze$ such that the coverage probability $\Pr \{ | \wh{\bs{p}} - p| < \vep p \mid p \}$ is no less than but close to the
prescribed level $1 - \de$ for all $p \in (0, 1)$.

In many situations, it is desirable to estimate the binomial parameter $p$ with a mixed error criterion.  More formally, let $\vep_a \in (0, 1)$
and $\vep_r \in (0, 1)$, we wish to design a multistage procedure such that the resultant estimator $\wh{\bs{p}}$ guarantees that $\Pr \{ |
\wh{\bs{p}} - p| < \vep_a \; \tx{or} \;  | \wh{\bs{p}} - p| < \vep_r p \mid p \} \geq 1 - \de$ for all $p \in (0, 1)$.  For this purpose, we can
apply the general stopping rule by identifying the lower and upper bounds of the random interval as $\mscr{L} (\wh{\bs{p}}_\ell) = \min \li \{
\wh{\bs{p}}_\ell - \vep_a, \; \f{\wh{\bs{p}}_\ell}{1 + \vep_r} \ri \}$ and  $\mscr{U} (\wh{\bs{p}}_\ell) = \max \li \{ \wh{\bs{p}}_\ell +
\vep_a, \; \f{\wh{\bs{p}}_\ell}{1 - \vep_r} \ri \}$ respectively to derive the following stopping rule: Continue sampling until either
\[
\f{1}{2} + \ro \vep_a - \sq{ \f{1}{4} - n_\ell \li ( \f{\vep_a}{\mcal{Z}_{\ze \de_\ell}}  \ri )^2  } < \wh{\bs{p}}_\ell < \f{ (1 + \vep_r) ( 1 +
\vep_r - \ro \vep_r ) \mcal{Z}_{\ze \de_\ell}^2 } {  n_\ell \vep_r^2 + ( 1 + \vep_r - \ro \vep_r  )^2 \mcal{Z}_{\ze \de_\ell}^2   }
\]
or
\[
\f{1}{2} - \ro \vep_a - \sq{ \f{1}{4} - n_\ell \li ( \f{\vep_a}{\mcal{Z}_{\ze \de_\ell}}  \ri )^2  } < \wh{\bs{p}}_\ell < \f{ (1 - \vep_r) ( 1 -
\vep_r + \ro \vep_r ) \mcal{Z}_{\ze \de_\ell}^2 } {  n_\ell \vep_r^2 + ( 1 - \vep_r + \ro \vep_r  )^2 \mcal{Z}_{\ze \de_\ell}^2   }
\]
is violated for some $\ell \in \{1, 2, \cd, s \}$.   To ensure that the coverage probability $\Pr \{ | \wh{\bs{p}} - p| < \vep_a \; \tx{or} \; |
\wh{\bs{p}} - p| < \vep_r p \mid p \}$ is no less than but close to the prescribed level $1 - \de$ for all $p \in (0, 1)$, we can employ the
bisection coverage tuning techniques.

We would like to point out that the above stopping rules can be readily adapted for estimating the proportion $p$ of a finite population having
a certain attribute based on multistage sampling without replacement.   The adaptation is just to replace $n_\ell$ in the above stopping
conditions with $\f{ n_\ell (N - 1) }{N - n_\ell}$, while the estimator $\wh{\bs{p}}_\ell$ is interpreted as the sample mean $\f{ \sum_{i=1}^n
X_i }{n_\ell }$, where $X_1, X_2, \cd X_N$ are random variables associated with the scheme of sampling without replacement such that $X_i$
assumes value $1$ if the $i$-th drawn unit has the attribute, and assumes value $0$ otherwise.  It should be noted that the stopping rules from
such adaptation can be equivalently derived from the general stopping rule associated with (\ref{st38a}) and (\ref{st38b}).

With regard to the estimation of the binomial parameter $p$ to satisfy prescribed levels of relative precision and confidence, we can derive a
multistage version of inverse binomial sampling scheme from the general stopping rule associated with (\ref{st38a}) and (\ref{st38b}).  In terms
of i.i.d. Bernoulli random variables $X_1, X_2, \cd$ with common mean $p \in (0, 1)$, a random variable $Y$ possessing a geometric distribution
can be defined as the minimum integer such that $\sum_{i =1}^Y X_i = 1$.    Clearly, the mean of $Y$ is $\se = \f{1}{p}$.  So, to estimate $p$,
it suffices to estimate $\se$ based on the samples of $Y$.  For this purpose, note that the sequence of  i.i.d. Bernoulli samples defines a
sequence of i.i.d. samples $Y_1, Y_2, \cd$ of $Y$ such that $Y_j$ is the minimum integer satisfying $\sum_{i = 1 + Y_{j -1}}^{Y_j} X_i = 1$ for
$j = 1, 2, \cd$, where $Y_0 = 0$.  Consider a multistage sampling scheme for $Y$ with deterministic sample sizes $n_1, n_2, \cd, n_s$, where the
number of stages $s$ is finite. For $\ell = 1, \cd, s$, the estimator for $\se$ is defined as $\wh{\bs{\se}}_\ell = \f{ \sum_{i = 1}^{n_\ell}
Y_i }{n_\ell}$. The sequential estimator $\wh{\bs{\se}}$ is defined as $\wh{\bs{\se}}_{\bs{l}}$, where $\bs{l}$ is the index of stage at the
termination of the sampling process. Let $\mcal{Z}_{\ze \de}$ be the critical value such that $\Phi( \mcal{Z}_{\ze \de} ) = 1 - \ze \de$.  To
design a multistage procedure such that the resultant estimator $\wh{\bs{\se}}$ guarantees that $\Pr \{ (1 - \vep) \wh{\bs{\se}} < \se < (1 +
\vep) \wh{\bs{\se}} \mid \se \} \geq 1 - \de$ for all $\se \in (1, \iy)$, we can apply the general stopping rule by identifying the lower and
upper bounds of the random interval as $\mscr{L} (\wh{\bs{\se}}_\ell) = (1 - \vep) \wh{\bs{\se}}_\ell$ and  $\mscr{U} (\wh{\bs{\se}}_\ell) = (1
+ \vep) \wh{\bs{\se}}_\ell$ respectively to derive the following stopping rule: Continue sampling until $\f{1}{\wh{\bs{\se}}_\ell}  \geq 1 + \ro
\vep  - \f{ n_\ell  } { 1 + \ro \vep } \li (  \f{\vep}{ \mcal{Z}_{\ze \de} } \ri)^2$ for some $\ell \in \{1, 2, \cd, s\}$.   By virtue of
bisection coverage tuning techniques, we can obtain an appropriate value of $\ze$ such that the coverage probability $\Pr \{ (1 - \vep)
\wh{\bs{\se}} < \se < (1 + \vep) \wh{\bs{\se}}  \mid \se \}$ is no less than but close to the prescribed level $1 - \de$ for all $\se \in (1,
\iy)$. Observing that the reciprocals of $\se$ and $\wh{\bs{\se}}$ are $p$ and $\wh{\bs{p}}$ respectively,  we have that such sampling scheme is
actually a multistage version of inverse binomial sampling scheme which ensures that $\Pr \{ | \wh{\bs{p}} - p| < \vep p \mid p \}$ for all $p
\in (0, 1)$.

\item [Estimating Poisson Parameters:]  Let $X$ be a Poisson random variable with mean $\lm > 0$.
Let $X_1, X_2, \cd$ be i.i.d. samples of $X$.   Similar to the estimation of a proportion, we can apply the general stopping rule to the
multistage estimation of $\lm$ with deterministic sample sizes $n_1, n_2, \cd, n_s$, where the number of stages $s$ can be finite or infinite.
For $\ell = 1, \cd, s$, the estimator for $\lm$ is defined as $\wh{\bs{\lm}}_\ell = \f{ \sum_{i = 1}^{n_\ell} X_i }{n_\ell}$.   The sequential
estimator $\wh{\bs{\lm}}$ is defined as $\wh{\bs{\lm}}_{\bs{l}}$, where $\bs{l}$ is the index of stage at the termination of the sampling
process. The prescribed confidence interval is $1 - \de$, where $\de \in (0, 1)$.  For $\ell = 1, 2, \cd, s$, let $\mcal{Z}_{\ze \de_\ell}$ be
the critical value such that $\Phi( \mcal{Z}_{\ze \de_\ell} ) = 1 - \ze \de_\ell$.

Let $\vep > 0$ be a margin of absolute error. To design a multistage procedure such that the sequential estimator $\wh{\bs{\lm}}$ guarantees
that $\Pr \{ | \wh{\bs{\lm}} - \lm| < \vep \mid \lm \} \geq 1 - \de$ for all $\lm \in (0, \iy)$, we can apply the general stopping rule to
derive a stopping rule as follows: Continue sampling until $\wh{\bs{\lm}}_\ell \leq \li ( \f{\vep}{\mcal{Z}_{\ze \de_\ell}} \ri )^2 n_\ell - \ro
\vep$ for some $\ell \in \{1, 2, \cd, \iy \}$.  For the controllability of the coverage probability, the parameters $\de_\ell, \; \ell = 1, 2,
\cd$ are chosen to be positive numbers less than $1$ such that their sum $\sum_{\ell = 1}^\iy \de_\ell$ is bounded.

Let $\vep > 0$ be a margin of relative error.  To design a multistage procedure such that the resultant estimator $\wh{\bs{\lm}}$ guarantees
that $\Pr \{ | \wh{\bs{\lm}} - \lm| < \vep \lm \mid \lm \} \geq 1 - \de$ for all $\lm \in (0, \iy)$, we can apply the general stopping rule to
derive the following stopping rule: Continue sampling until $\wh{\bs{\lm}}_\ell \geq \f{1}{n_\ell} \li ( \f{\mcal{Z}_{\ze \de_\ell}}{\vep} \ri
)^2 ( 1 + \vep) ( 1 + \vep - \ro \vep )$ for some $\ell \in \{1, 2, 3, \cd, \iy \}$.  The parameters $\de_\ell, \; \ell = 1, 2, \cd$ are
positive numbers less than $1$ such that their sum $\sum_{\ell = 1}^\iy \de_\ell$ is bounded.

Let $\vep_a \in (0, \iy)$ and $\vep_r \in (0, 1)$ be respectively the margins of absolute and relative errors. To design a multistage procedure
such that the resultant estimator $\wh{\bs{\lm}}$ guarantees that $\Pr \{ | \wh{\bs{\lm}} - \lm| < \vep_a \; \tx{or} \; | \wh{\bs{\lm}} - \lm| <
\vep_r \lm \mid \lm \} \geq 1 - \de$ for all $\lm \in (0, \iy)$, we can apply the general stopping rule  to derive the following stopping rule:
Continue sampling until either
\[
\li ( \f{\vep_a}{\mcal{Z}_{\ze \de_\ell}} \ri )^2 n_\ell + \ro \vep_a < \wh{\bs{\lm}}_\ell < \f{1}{n_\ell} \li ( \f{\mcal{Z}_{\ze
\de_\ell}}{\vep_r} \ri )^2 ( 1 + \vep_r) ( 1 + \vep_r - \ro \vep_r )
\]
or
\[
\li ( \f{\vep_a}{\mcal{Z}_{\ze \de_\ell}} \ri )^2 n_\ell - \ro \vep_a < \wh{\bs{\lm}}_\ell < \f{1}{n_\ell} \li ( \f{\mcal{Z}_{\ze
\de_\ell}}{\vep_r} \ri )^2 ( 1 - \vep_r) ( 1 - \vep_r + \ro \vep_r )
\]
is violated for some $\ell \in \{1, 2, \cd, s\}$.

For the above estimation problems with various error criteria, to ensure that the coverage probability is no less than but close to the
prescribed level $1 - \de$ for all $\lm \in (0, \iy)$, we can employ the bisection coverage tuning techniques.

\eed

In the above presentation, we have extensively discussed approximation methods for simplifying stopping rules.  In addition to the normal
approximation, bounds of the CDF $\&$ CCDF of $\wh{\bs{\se}}_\ell$ can be used to simplify stopping rules. Specially, in situations that sample
sizes are deterministic numbers $n_1 < n_2 < \cd < n_s$ and that $\wh{\bs{\se}}_\ell = \f{\sum_{i = 1}^{n_\ell} X_i }{n_\ell}$ for $\ell = 1,
\cd, s$, we have established multistage sampling schemes by virtue of Theorem \ref{Monotone_second} and Chernoff bounds as follows.

\beC

\la{Monotone_third}

Suppose that a multistage sampling scheme satisfies the following
requirements:

(i) For $\ell = 1, \cd, s$, $\wh{\bs{\se}}_\ell= \f{\sum_{i =
1}^{n_\ell} X_i }{n_\ell}$ is a ULE of $\se$, where $X_1, X_2, \cd$
are i.i.d. samples of $X$.

(ii) The moment generating function $\bb{E} [e^{ t X}  ]$ exists for
any real number $t$.

(iii) For $\ell = 1, \cd, s$, $\{ \mscr{L} (\wh{\bs{\se}}_\ell,
n_\ell) \leq \wh{\bs{\se}}_\ell \leq \mscr{U} (\wh{\bs{\se}}_\ell,
n_\ell) \}$ is a sure event.

(iv) {\small $\li \{ \bs{D}_\ell = 1 \ri \} \subseteq \li \{ \li [
\mcal{F} \li ( \wh{\bs{\se}}_\ell, \mscr{U} (\wh{\bs{\se}}_\ell,
n_\ell) \ri ) \ri ]^{n_\ell} \leq \ze \de_\ell, \; \li [ \mcal{G}
\li ( \wh{\bs{\se}}_\ell, \mscr{L} (\wh{\bs{\se}}_\ell, n_\ell) \ri
) \ri ]^{n_\ell} \leq \ze \de_\ell \ri \}$} for $\ell = 1, \cd, s$,
where $\mcal{F} (.,.)$
 and $\mcal{G} (.,.)$ are functions such that
{\small \[ \mcal{F} (z, \se) = \bec  \inf_{t <
0} \bb{E} [e^{t ( X - z)} ]  & \tx{for} \; \se \in \Se,\\
1  & \tx{for} \; \se < \inf \Se,\\
0  & \tx{for} \; \se > \sup \Se \eec \qqu \mcal{G} (z, \se) = \bec
\inf_{t >
0} \bb{E} [e^{t ( X - z)} ]  & \tx{for} \; \se \in \Se,\\
0  & \tx{for} \; \se < \inf \Se,\\
1  & \tx{for} \; \se > \sup \Se \eec
 \]}

(v) $\{ \bs{D}_s = 1 \}$ is a sure event.

Then, \bee &   & \Pr \{ \mscr{L} (\wh{\bs{\se}}, \mbf{n})  \geq \se
\mid \se \} \leq \sum_{\ell = 1}^s \Pr \{ \mscr{L}
(\wh{\bs{\se}}_\ell, n_\ell) \geq \se, \; \bs{D}_\ell = 1 \mid \se
\} \leq \ze \sum_{\ell = 1}^s \de_\ell,\\
&   & \Pr \{ \mscr{U} (\wh{\bs{\se}}, \mbf{n})  \leq \se \mid \se \}
\leq \sum_{\ell = 1}^s \Pr \{ \mscr{U} (\wh{\bs{\se}}_\ell, n_\ell)
\leq \se, \; \bs{D}_\ell = 1 \mid \se \} \leq \ze \sum_{\ell = 1}^s
\de_\ell,\\
&  & \Pr \{ \mscr{L} (\wh{\bs{\se}}, \mbf{n})  < \se < \mscr{U}
(\wh{\bs{\se}}, \mbf{n}) \mid \se \} \geq 1 - 2 \ze \sum_{\ell =
1}^s \de_\ell \eee for any $\se \in \Se$.

\eeC

To establish Corollary \ref{Monotone_third}, it suffices to show that the assumption (iv) of Corollary \ref{Monotone_third} implies the
assumption (iii) of Theorem \ref{Monotone_second}, which can be seen from Chernoff bounds \bee &  & F_{\wh{\bs{\se}}_\ell} (z, \se) \leq [
\mcal{F} (z, \se) ]^{n_\ell}, \qqu G_{\wh{\bs{\se}}_\ell} (z, \se) \leq [ \mcal{G} (z, \se) ]^{n_\ell} \eee for $\se \in \Se$ and $z$ assuming
values from the support of $\wh{\bs{\se}}_\ell$. It can seen that the method of defining stopping rules proposed in Corollary
\ref{Monotone_third} is in the same spirit of the inclusion principle, except that the confidence limits are more conservative since the bounds
of CDF $\&$ CDF are used.  As will be seen in the sequel, the conservativeness can be significantly reduced by virtue of coverage tuning.

It should be noted that explicit forms for functions $\mcal{F} (z,
\se)$ and $\mcal{G} (z, \se)$ in Corollary \ref{Monotone_third} can
be derived for the exponential family.  A single-parameter
exponential family is a set of probability distributions whose
probability density function (or probability mass function, for the
case of a discrete distribution) can be expressed in the form \be
\la{expdef}
 f_X (x, \se) = c(x) \exp ( \eta (\se) x - \psi (\se) ), \qqu \se
 \in \Se
\ee where $c(x), \eta(\se)$, and $\psi (\se)$ are known functions.
Regarding the sample mean of $X$, we have the following results.

\beT  \la{explicit_Chernoff}   Let $\ovl{X}_n = \f{ \sum_{i=1}^n X_i
}{n}$, where $X_1, \cd, X_n$ are i.i.d. samples of random variable
$X$ possessing a probability density function or probability mass
function defined by (\ref{expdef}).  Suppose that $\f{d \eta
(\se)}{d \se}$ is positive and that $\f{ d \psi ( \se ) }{ d \se } =
\se \f{ d \eta ( \se ) }{d \se}$ for $\se \in \Se$.  Then,
$\ovl{X}_n$ is a ULE and an unbiased estimator of $\se$. Moreover,
\bee & & \Pr \{ \ovl{X}_n \leq z \mid \se \} \leq \li ( \inf_{t < 0
} \bb{E} \li [ e^{ t ( X - z) } \ri ] \ri )^n = [w (z, \se )]^n \qqu
\tx{for $z \leq
\se$},\\
&  & \Pr \{ \ovl{X}_n \geq z \mid \se  \} \leq \li ( \inf_{t > 0 }
\bb{E} \li [ e^{ t  ( X - z) } \ri ] \ri )^n =   [ w (z, \se )]^n
\qqu \tx{for $z \geq \se$} \eee where $w (z, \se ) = \f{ \exp \li (
\eta (\se ) z -  \psi (\se) \ri ) } { \exp \li ( \eta (z ) z - \psi
(z) \ri ) }$.  \eeT

See Appendix \ref{explicit_Chernoff_app} for a proof.  Applying
Theorem \ref{explicit_Chernoff} and Corollary \ref{Monotone_third}
to the estimation of the parameter of the exponential family, we
have $\mcal{F} (z, \se) = \mcal{G} (z, \se) = w(z, \se)$ for $\se
\in \Se$ and consequently the sampling scheme can be simplified.

\bsk

From Theorems  \ref{Monotone_second}, \ref{Monotone_second_CI_ST} and Corollary \ref{Monotone_third}, it can be seen that, if the number of
stage $s$ is independent of the coverage tuning parameter $\ze$, then the coverage probability of the sequential random interval $(\mscr{L}
(\wh{\bs{\se}}, \mbf{n}), \mscr{U} (\wh{\bs{\se}}, \mbf{n}))$ can be adjusted to be above $1 - \de$ if $\ze$ is sufficiently small.  In the
design of sampling scheme, the number of stages and the sample sizes at all stages can be dependent on the coverage tuning parameter $\ze$.  To
satisfy the coverage requirement, we hope that the coverage probability of the sequential random interval can still be controlled by $\ze$. Such
controllability can be established under mild conditions. Specially, to construct sequential random interval $(\mscr{L} (\wh{\bs{\se}}),
\mscr{U} (\wh{\bs{\se}}))$, where $\mscr{L} (.)$ and $\mscr{U} (. )$ are univariate functions of $ \wh{\bs{\se}}$, we have the following result
regarding the sampling schemes described by Theorem \ref{Monotone_second} and Corollary \ref{Monotone_third}, where the stage number $s =
s(\ze)$ and the sample sizes $n_\ell = n_\ell (\ze), \; \ell = 1, \cd, s$ are functions of $\ze$.

\beT \la{fullcontrol} Let $X_1, X_2, \cd$ be i.i.d. samples of
random variable $X$ possessing a probability density function or
probability mass function defined by (\ref{expdef}). Suppose that
$\f{d \eta (\se)}{d \se}$ is positive and that $\f{ d \psi ( \se )
}{ d \se } = \se \f{ d \eta ( \se ) }{d \se}$ for $\se \in \Se$.
Suppose that $\mscr{L}( \se ) < \se < \mscr{U}( \se )$ for all $\se
\in \Se$. Let $\de_\ell, \; \ell = 1, 2, \cd$ be given such that
$\sum_{\ell = 1}^\iy \de_\ell < \iy$ or $0 < \al < \de_\ell < \ba$
for all $\ell$, where $\al $ and $\ba$ are some positive numbers.
Then, the coverage probability $\Pr \{ \mscr{L} (\wh{\bs{\se}}) <
\se < \mscr{U} (\wh{\bs{\se}}) \mid \se \}$ tends to $1$ as $\ze$
tends to $0$.

\eeT

See Appendix \ref{fullcontrol_app} for a proof.

Now, we turn to consider the second problem posed at the beginning of this subsection.  For the sampling schemes of structure described in
Section \ref{gen_structure}, we have the following results regarding the coverage probability of sequential random intervals.

\beT \la{Main_Bound_Gen}

Let $X_1, X_2, \cd$ be a sequence of identical samples of discrete
random variable $X$ parameterized by $\se \in \Se$. For $\ell = 1,
\cd, s$, let $\wh{\bs{\se}}_\ell =\varphi (X_1, \cd,
X_{\mbf{n}_\ell})$ be a ULE of $\se$. Define estimator
$\wh{\bs{\se}} = \wh{\bs{\se}}_{\bs{l}}$, where $\bs{l}$ is the
index of stage when the sampling is terminated.  Let $\mscr{L}(.,.)$
and $\mscr{U}(.,.)$  be bivariate functions such that $\{
\mscr{L}(\wh{\bs{\se}}, \mbf{n}) \leq \wh{\bs{\se}} \leq \mscr{U}
(\wh{\bs{\se}}, \mbf{n}) \}$ is a sure event.   Let $[a, b]$ be a
subset of $\Se$.  Let $I_\mscr{L}$ denote the intersection of
interval $(a, b)$ and the support of $\mscr{L}(\wh{\bs{\se}},
\mbf{n} )$. Let $I_\mscr{U}$ denote the intersection of interval
$(a, b)$ and the support of $\mscr{U}(\wh{\bs{\se}}, \mbf{n} )$. Let
$\mscr{E}$ be an event dependent only on the random tuple $(X_1,
\cd, X_{\mbf{n}})$. The following statements hold true:

(I) Both $\Pr \{ \mscr{L}( \wh{\bs{\se}}, \mbf{n} ) \geq \se \;
\tx{and} \; \mscr{E} \; \tx{occurs} \mid \se \}$ and $\Pr \{
\mscr{L}( \wh{\bs{\se}}, \mbf{n} )
> \se \; \tx{and} \; \mscr{E} \; \tx{occurs} \mid \se
\}$ are no-decreasing with respect to $\se$ in any open interval
with endpoints being consecutive distinct elements of $I_\mscr{L}
\cup \{ a, b\}$. Moreover, both the maximum of $\Pr \{ \mscr{L}(
\wh{\bs{\se}}, \mbf{n} ) \geq \se \; \tx{and} \; \mscr{E} \;
\tx{occurs} \mid \se \}$  and the supremum of $\Pr \{ \mscr{L}(
\wh{\bs{\se}}, \mbf{n} ) > \se \; \tx{and} \; \mscr{E} \;
\tx{occurs} \mid \se \}$ with respect to $\se \in [a, b]$ are equal
to the maximum of $\Pr \{ \mscr{L}( \wh{\bs{\se}}, \mbf{n} ) \geq
\se \; \tx{and} \; \mscr{E} \; \tx{occurs} \mid \se \}$ for $\se \in
I_\mscr{L} \cup \{ a, b\}$.

(II)  Both $\Pr \{ \mscr{U}( \wh{\bs{\se}}, \mbf{n} ) \leq \se \;
\tx{and} \; \mscr{E} \; \tx{occurs} \mid \se \}$ and $\Pr \{
\mscr{U}( \wh{\bs{\se}}, \mbf{n} ) < \se \; \tx{and} \; \mscr{E} \;
\tx{occurs} \mid \se \}$ are non-increasing with respect to $\se$ in
any open interval with endpoints being consecutive distinct elements
of $I_\mscr{U} \cup \{ a, b\}$.  Moreover, both the maximum of $\Pr
\{ \mscr{U}( \wh{\bs{\se}}, \mbf{n} ) \leq \se \; \tx{and} \;
\mscr{E} \; \tx{occurs} \mid \se \}$  and the supremum of $\Pr \{
\mscr{U}( \wh{\bs{\se}}, \mbf{n} ) < \se \; \tx{and} \; \mscr{E} \;
\tx{occurs} \mid \se \}$ with respect to $\se \in [a, b]$ are equal
to the maximum of $\Pr \{ \mscr{U}( \wh{\bs{\se}}, \mbf{n} ) \leq
\se \; \tx{and} \; \mscr{E} \; \tx{occurs} \mid \se \}$ for $\se \in
I_\mscr{U} \cup \{ a, b\}$.

(III) If $\{ \mscr{L}( \wh{\bs{\se}}, \mbf{n} ) \geq a \} \subseteq
\{ \wh{\bs{\se}} \geq b   \}$,  then {\small $\Pr \{ \mscr{L}(
\wh{\bs{\se}}, \mbf{n} ) \geq b \; \tx{and} \; \mscr{E} \;
\tx{occurs} \mid a \} \leq \Pr \{ \mscr{L}( \wh{\bs{\se}}, \mbf{n} )
\geq \se \; \tx{and} \; \mscr{E} \; \tx{occurs} \mid \se \} \leq \Pr
\{ \mscr{L}( \wh{\bs{\se}}, \mbf{n} ) \geq a \; \tx{and} \; \mscr{E}
\; \tx{occurs} \mid b \}$} and {\small $\Pr \{ \mscr{L}(
\wh{\bs{\se}}, \mbf{n} ) > b \; \tx{and} \; \mscr{E} \; \tx{occurs}
\mid a \} \leq \Pr \{ \mscr{L}( \wh{\bs{\se}}, \mbf{n} ) > \se \;
\tx{and} \; \mscr{E} \; \tx{occurs} \mid \se \} \leq \Pr \{
\mscr{L}( \wh{\bs{\se}}, \mbf{n} ) > a \; \tx{and} \; \mscr{E} \;
\tx{occurs} \mid b \}$} for any $\se \in [a, b]$. Similarly,  if $\{
\mscr{U}( \wh{\bs{\se}}, \mbf{n} ) \leq b \} \subseteq \{
\wh{\bs{\se}} \leq a \}$,  then {\small $\Pr \{ \mscr{U}(
\wh{\bs{\se}}, \mbf{n} ) \leq a \; \tx{and} \; \mscr{E} \;
\tx{occurs} \mid b \} \leq \Pr \{ \mscr{U}( \wh{\bs{\se}}, \mbf{n} )
\leq \se \; \tx{and} \; \mscr{E} \; \tx{occurs} \mid \se \} \leq \Pr
\{ \mscr{U}( \wh{\bs{\se}}, \mbf{n} ) \leq b \; \tx{and} \; \mscr{E}
\; \tx{occurs} \mid a \}$} and {\small $\Pr \{ \mscr{U}(
\wh{\bs{\se}}, \mbf{n} ) < a \; \tx{and} \; \mscr{E} \; \tx{occurs}
\mid b \} \leq \Pr \{ \mscr{U}( \wh{\bs{\se}}, \mbf{n} ) < \se \;
\tx{and} \; \mscr{E} \; \tx{occurs} \mid \se \} \leq \Pr \{
\mscr{U}( \wh{\bs{\se}}, \mbf{n} ) < b \; \tx{and} \; \mscr{E} \;
\tx{occurs} \mid a \}$} for any $\se \in [a, b]$.

(IV) If $\mscr{L}( \wh{\bs{\se}}, \mbf{n} )$ and $\mscr{U}(
\wh{\bs{\se}}, \mbf{n} )$ can be expressed as non-decreasing
univariate functions $\mscr{L}( \wh{\bs{\se}} )$ and $\mscr{U}(
\wh{\bs{\se}})$ of $\wh{\bs{\se}}$ respectively, then, without the
assumption that $\{ \mscr{L}(\wh{\bs{\se}}, \mbf{n}) \leq
\wh{\bs{\se}} \leq \mscr{U} (\wh{\bs{\se}}, \mbf{n}) \}$ is a sure
event, {\small \bee & & \Pr \{ \mscr{L}( \wh{\bs{\se}} ) \geq b \mid
a \} \leq \Pr \{ \mscr{L}( \wh{\bs{\se}} ) \geq \se \mid \se \} \leq
\Pr \{ \mscr{L}(
\wh{\bs{\se}} ) \geq a \mid b \}, \\
&  & \Pr \{ \mscr{L}( \wh{\bs{\se}} )
> b \mid a \} \leq \Pr \{ \mscr{L}( \wh{\bs{\se}} )
> \se \mid \se \} \leq \Pr \{ \mscr{L}( \wh{\bs{\se}} ) > a \mid b \},\\
&  & \Pr \{ \mscr{U}( \wh{\bs{\se}} ) \leq a \mid b \} \leq \Pr \{
\mscr{U}( \wh{\bs{\se}} ) \leq \se  \mid \se \} \leq \Pr \{
\mscr{U}( \wh{\bs{\se}} ) \leq b \mid a \},\\
&  & \Pr \{ \mscr{U}( \wh{\bs{\se}} ) < a \mid b \} \leq \Pr \{
\mscr{U}( \wh{\bs{\se}} ) < \se \mid \se \} \leq \Pr \{ \mscr{U}(
\wh{\bs{\se}} ) < b \mid a \} \eee} for any $\se \in [a, b]$.

\eeT

\bsk

See Appendix \ref{App_Main_Bound_Gen} for a proof.   Actually, as special results of Theorem \ref{Main_Bound_Gen}, we have established ``Theorem
8''  and other similar theorems in the $12$th version of this paper published in arXiv on April 27, 2009.   In Theorem \ref{Main_Bound_Gen}, we
have used the concept of support in probability theory. The support of a random variable $Z$ refers to $\{Z (\om): \om \in \Om \}$, which is the
set of all possible values of $Z$.  We say that ``an event $\mscr{E}$ is dependent only on the random tuple $(X_1, \cd, X_{\mbf{n}})$'' if, for
any $n$ in the support of $\mbf{n}$, the event $\{ \tx{$\mscr{E}$  occurs and} \; \mbf{n} = n \}$ can be expressed in terms of random variables
$X_1, \cd, X_n$.  In a more rigorous probabilistic terminology, $\mbf{n}$ is a {\it stopping time} and $\mscr{E}$ is an event of the
$\si$-subalgebra generated by $X_1, \cd, X_{\mbf{n}}$.

Based on Theorem \ref{Main_Bound_Gen} in the special case that
$\mscr{E}$ is a sure event, two different approaches can be
developed to address the second problem proposed at the beginning of
this subsection.

First, as a consequence of statements (I) and (II) of Theorem
\ref{Main_Bound_Gen}, it is true that $\Pr \{ \mscr{L}(
\wh{\bs{\se}}, \mbf{n} ) < \se < \mscr{U}( \wh{\bs{\se}}, \mbf{n} )
\mid \se \} \geq 1 - \de$ for any $\se \in [a, b]$ provided that
\[
\Pr \{ \se \leq \mscr{L}( \wh{\bs{\se}}, \mbf{n} ) \mid \se \} \leq
\f{\de}{2}, \qqu \fa \se \in I_\mscr{L} \cup \{ a, b\},
\]
\[
\Pr \{ \se \geq \mscr{U}( \wh{\bs{\se}}, \mbf{n} ) \mid \se \} \leq
\f{\de}{2}, \qqu \fa \se \in I_\mscr{U} \cup \{ a, b\}.
\]
As can be seen from the proofs of Theorems 1 and 2, under certain
conditions, the probabilities $\Pr \{ \se \leq \mscr{L}(
\wh{\bs{\se}}, \mbf{n} ) \mid \se \}$ and $\Pr \{ \se \geq \mscr{U}(
\wh{\bs{\se}}, \mbf{n} ) \mid \se \}$ can be adjusted by $\ze$.
Hence, it is possible to obtain appropriate value of $\ze$, without
exhaustive evaluation of probabilities, such that $\Pr \{ \mscr{L}(
\wh{\bs{\se}}, \mbf{n} ) < \se < \mscr{U}( \wh{\bs{\se}}, \mbf{n} )
\mid \se \} \geq 1 - \de$ for any $\se \in [a, b]$.

Second, statements (III) and (IV) of Theorem \ref{Main_Bound_Gen}
will be used to develop Adaptive Maximum Checking Algorithm in
Section \ref{AMCA} to determine an appropriate value of coverage
tuning parameter $\ze$.

In the special case that the number of stages $s$ is equal to $1$
and that the sample number is a deterministic integer $n$, we have
the following results.

\beT \la{Fundamental} Let $X_1, X_2, \cd, X_n$ be a sequence of
discrete random variables parameterized by $\se \in \Se$. Let
$\wh{\bs{\se}} =\varphi (X_1, \cd, X_n)$ be an estimator of $\se$.
Let $\mscr{L}(.)$ and $\mscr{U}(.)$ be functions such that, for any
$\vartheta \in \Se$, $\Pr \{ \mscr{L}(\wh{\bs{\se}}) \leq \vartheta
\leq \mscr{U} (\wh{\bs{\se}}) \mid \se \}$ is a continuous and
unimodal function of $\se \in \Se$.  Let $[a, b]$ be an interval
contained in $\Se$.  Let $I_\mscr{L}$ denote the intersection of the
interval $(a, b)$ and the support of $\mscr{L}(\wh{\bs{\se}})$. Let
$I_\mscr{U}$ denote the intersection of the interval $(a, b)$ and
the support of $\mscr{U}(\wh{\bs{\se}})$.  Then, the minimum of $\Pr
\{ \mscr{L}(\wh{\bs{\se}}) < \se < \mscr{U}(\wh{\bs{\se}}) \mid \se
\}$ with respect to $\se \in [a, b]$ is attained at the set
$I_{\mscr{L}} \cup I_{\mscr{U}} \cup \{a, b \} $ and the infimum of
$\Pr \{ \mscr{L}(\wh{\bs{\se}}) \leq \se \leq
\mscr{U}(\wh{\bs{\se}}) \mid \se \}$ with respect to $\se \in [a,
b]$ is equal to the minimum of the set $\li \{ C_L (\se):  \se \in
I_{\mscr{L}} \ri \} \cup \li \{ C_U (\se): \se \in I_{\mscr{U}} \ri
\} \cup \{C(a), \; C_U(a), \;  C(b), \; C_L (b) \} $, where $ C_L
(\se) = \Pr \{ \mscr{L}(\wh{\bs{\se}}) < \se \leq
\mscr{U}(\wh{\bs{\se}}) \mid \se \}, \; C_U (\se) = \Pr \{
\mscr{L}(\wh{\bs{\se}}) \leq \se < \mscr{U}(\wh{\bs{\se}}) \mid \se
\}$ and $C(\se) = \Pr \{ \mscr{L}(\wh{\bs{\se}}) \leq \se \leq
\mscr{U}(\wh{\bs{\se}}) \mid \se \}$. Moreover, for both open random
interval $((\mscr{L}(\wh{\bs{\se}}), \mscr{U}(\wh{\bs{\se}}))$ and
closed random interval $[\mscr{L}(\wh{\bs{\se}}),
\mscr{U}(\wh{\bs{\se}})]$, the coverage probability is continuous
and unimodal for $\se \in (\se^\prime, \se^{\prime \prime})$, where
$\se^\prime$ and $\se^{\prime \prime}$ are arbitrary consecutive
distinct elements of $I_{\mscr{L}} \cup I_{\mscr{U}} \cup \{a, b
\}$.

\eeT

The proof of Theorem \ref{Fundamental} can be found in
\cite{Chen_RI}.

\subsection{Multistage Sampling without Replacement}
\la{finite_proportion}

It should be noted that the theories in preceding discussion can be
applied to the multistage estimation of the proportion of a finite
population, where the random samples are dependent if a sampling
without replacement is used.   Consider a population of $N$ units,
among which there are $p N$ units having a certain attribute, where
$p \in \Se = \{\f{M}{N}: M = 0, 1, \cd, N \}$. In many situations,
it is desirable to estimate the population proportion $p$ by
sampling without replacement. The procedure of sampling without
replacement can be precisely described as follows:

 Each time a single unit is drawn without replacement from the remaining population so
that every unit of the remaining population has equal chance of
being selected.

Such a sampling process can be exactly characterized by random
variables $X_1, \cd, X_N$ defined in a probability space $(\Omega,
\mscr{F}, \Pr)$ such that $X_i$ assumes value $1$ if the $i$-th
sample has the attribute and assumes value $0$ otherwise. By the
nature of the sampling procedure, it can be shown that
\[ \Pr \{ X_i = x_i, \; i = 1, \cd, n \} =
\bi{p N}{\sum_{i = 1}^n x_i} \bi{N - p N}{n - \sum_{i = 1}^n x_i}
 \li \slash \li [ \bi{n}{\sum_{i = 1}^n x_i} \bi{N}{n} \ri. \ri ]
\]
for any $n \in \{1, \cd, N\}$ and any $x_i \in \{0, 1\}, \; i = 1,
\cd, n$.   Clearly, for any $n \in \{1, \cd, N\}$, the sample mean
{\small $\f{ \sum_{i=1}^n X_i }{n}$} is unbiased but is not a MLE
for $p \in \Se$. However, we have shown in Appendix
\ref{FiniteULE_Ap} the following result:

\beT \la{FiniteULE} For any $n \in \{1, \cd, N\}$, {\small $\f{
\sum_{i=1}^n X_i }{n}$} is a ULE for $p \in \Se$. \eeT

Based on random variables $X_1, \cd, X_N$, we can define a
multistage sampling scheme in the same way as that of the multistage
sampling described in Section \ref{gen_structure}.  More specially,
we can define decision variables such that, for the $\ell$-th stage,
$\bs{D}_\ell$ is a function of $X_1, \cd, X_{\mbf{n}_\ell}$, where
the random variable $\mbf{n}_\ell$ is the number of samples
available at the $\ell$-th stage. For $\ell = 1, \cd, s$, an
estimator of $p$ at the $\ell$-stage can be defined as {\small
$\wh{\bs{p}}_\ell = \f{\sum_{i = 1}^{\mbf{n}_\ell} X_i}{
\mbf{n}_\ell }$}. Letting $\bs{l}$ be the index of stage when the
sampling is terminated, we can define an estimator for  $p$ as
$\wh{\bs{p}} = \wh{\bs{p}}_{\bs{l}} = \f{\sum_{i = 1}^{\mathbf{n}}
X_i}{\mbf{n}}$, where $\mathbf{n} = \mbf{n}_{\bs{l}}$ is the sample
size at the termination of sampling. A sampling scheme described in
this setting is referred to as a {\it multistage sampling without
replacement} in this paper. Regarding the coverage probability of
random intervals, we have the following results which are direct
consequence of  Theorems \ref{Main_Bound_Gen} and \ref{FiniteULE}.

\beC \la{Finite_Main_Bound_Gen}  Let $\mscr{L}(.,.)$ and
$\mscr{U}(.,.)$  be bivariate functions such that $\{ \mscr{L}
(\wh{\bs{p}}, \mbf{n}) \leq \wh{\bs{p}} \leq \mscr{U} (\wh{\bs{p}},
\mbf{n}) \}$ is a sure event and that both $N \mscr{L} (\wh{\bs{p}},
\mbf{n})$ and $N \mscr{U} (\wh{\bs{p}}, \mbf{n})$ are integer-valued
random variables. Let $a \leq b$ be two parametric values in $\Se$.
Let $I_\mscr{L}$ denote the intersection of interval $(a, b)$ and
the support of $\mscr{L}(\wh{\bs{p}}, \mbf{n} )$. Let $I_\mscr{U}$
denote the intersection of interval $(a, b)$ and the support of
$\mscr{U}(\wh{\bs{p}}, \mbf{n} )$.  The following statements hold
true:

(I) $\Pr \{ \mscr{L}( \wh{\bs{p}}, \mbf{n} ) \geq p \mid p \}$ is
non-decreasing with respect to $p \in \Se$ in any interval with
endpoints being consecutive distinct elements of $I_\mscr{L} \cup \{
a, b\}$.  Moreover, the maximum of $\Pr \{ \mscr{L}( \wh{\bs{p}},
\mbf{n} ) \geq p \mid p \}$ with respect to $p \in [a, b] \cap \Se$
is achieved at $I_\mscr{L} \cup \{ a, b\}$.

(II) $\Pr \{ \mscr{U}( \wh{\bs{p}}, \mbf{n} ) \leq p \mid p \}$ is
non-increasing with respect to $p \in \Se$ in any interval with
endpoints being consecutive distinct elements of $I_\mscr{U}
 \cup \{ a, b\}$.  Moreover, the maximum of
 $\Pr \{ \mscr{U}( \wh{\bs{p}}, \mbf{n} ) \leq p \mid p \}$ with respect to $p
\in [a, b] \cap \Se$ is achieved at $I_\mscr{U} \cup \{ a, b\}$.

(III) If $\{ \mscr{L}( \wh{\bs{p}}, \mbf{n} ) \geq a \} \subseteq \{
\wh{\bs{p}} \geq b   \}$,  then $\Pr \{ \mscr{L}( \wh{\bs{p}},
\mbf{n} ) \geq b \mid a \} \leq \Pr \{ \mscr{L}( \wh{\bs{p}},
\mbf{n} ) \geq p \mid p \} \leq \Pr \{ \mscr{L}( \wh{\bs{p}},
\mbf{n} ) \geq a \mid b \}$ for any $p \in [a, b] \cap \Se$.
Similarly, if $\{ \mscr{U}( \wh{\bs{p}}, \mbf{n} ) \leq b \}
\subseteq \{ \wh{\bs{p}} \leq a \}$,  then $\Pr \{ \mscr{U}(
\wh{\bs{p}}, \mbf{n} ) \leq a \mid b \} \leq \Pr \{ \mscr{U}(
\wh{\bs{p}}, \mbf{n} ) \leq p \mid p \} \leq \Pr \{ \mscr{U}(
\wh{\bs{p}}, \mbf{n} ) \leq b \mid a \}$ for any $p \in [a, b] \cap
\Se$.

(IV) If $\mscr{L}( \wh{\bs{p}}, \mbf{n} )$ and $\mscr{U}(
\wh{\bs{p}}, \mbf{n} )$ can be expressed as non-decreasing
univariate functions $\mscr{L}( \wh{\bs{p}} )$ and $\mscr{U}(
\wh{\bs{p}})$ of $\wh{\bs{p}}$ respectively, then,  without the
assumption that $\{ \mscr{L} (\wh{\bs{p}}) \leq \wh{\bs{p}} \leq
\mscr{U} (\wh{\bs{p}}) \}$ is a sure event, {\small \bee & & \Pr \{
\mscr{L}( \wh{\bs{p}} ) \geq b  \mid a \} \leq \Pr \{ \mscr{L}(
\wh{\bs{p}} ) \geq p \mid p \} \leq \Pr \{ \mscr{L}(
\wh{\bs{p}} ) \geq a \mid b \}, \\
&  & \Pr \{ \mscr{U}( \wh{\bs{p}} ) \leq a \mid b \} \leq \Pr \{
\mscr{U}( \wh{\bs{p}} ) \leq p  \mid p \} \leq \Pr \{ \mscr{U}(
\wh{\bs{p}} ) \leq b \mid a \} \eee} for any $p \in [a, b] \cap
\Se$.

\eeC

In the special case that the number of stages $s$ is equal to $1$
and that the sample number is a deterministic integer $n$, we have
the following results.

\beT \la{Finite_RI_07} Let $a < b$ be two parametric values in
$\Se$. Suppose that $\mscr{L} (.)$ and $\mscr{U} (.)$ are
non-decreasing functions such that both $N \mscr{L} (\wh{\bs{p}})$
and $N \mscr{U} (\wh{\bs{p}})$ are integer-valued random variables.
Then, the minimum of $\Pr \{ \mscr{L}( \wh{\bs{p}}) < p < \mscr{U}(
\wh{\bs{p}}) \mid p \}$ with respect to $p \in [a, b] \cap \Se$ is
attained at a discrete set $I_{UL}$ which is the union of $\{a, b
\}$ and the supports of $\mscr{L}( \wh{\bs{p}})$ and $\mscr{U}(
\wh{\bs{p}})$. Moreover, $\Pr \{ \mscr{L}( \wh{\bs{p}}) < p <
\mscr{U}( \wh{\bs{p}}) \mid p \}$ is unimodal with respect to $p$ in
between consecutive distinct elements of $I_{UL}$. \eeT

The proof of Theorem \ref{Finite_RI_07} can be found in
\cite{Chen_RI}.

\subsection{Asymptotically Unbiased Estimators of Mean Values}

Some important distributions are determined by the mean values of
associated random variables. Familiar examples are binomial
distribution, Poisson distribution, normal distribution, and
exponential distribution. To estimate the expectation, $\mu$, of a
random variable $X$ based on i.i.d. samples $X_1, X_2, \cd$, we can
use a multistage sampling scheme with a structure described in
Section \ref{gen_structure}. Specially, an estimator of $\mu$ can be
defined as the sample mean $\wh{\bs{\mu}} = \f{ \sum_{i =
1}^{\mbf{n} } X_i } { \mbf{n} }$, where $\mbf{n}$ is the sample
number at the termination of sampling. To justify that the estimator
$\wh{\bs{\mu}}$ is superior than other estimators, we shall show its
asymptotic unbiasedness and relevant properties.  For a multistage
sampling scheme with deterministic sample sizes $n_1 < n_2 < \cd <
n_s$,  we have established the following general results.

\beT \la{Unbiased_Gen}  Suppose that the sampling process is sure to be
terminated at a finite stage and that $\inf_{\ell > 0} \f{n_{\ell +
1}}{n_\ell}$ is greater than $1$. The following statements hold true.

(I) If $X$ has a finite variance, then $\bb{E} [\wh{\bs{\mu}} -
\mu], \; \bb{E} |\wh{\bs{\mu}} - \mu|$ and $\bb{E} |\wh{\bs{\mu}} -
\mu|^2$ tend to $0$  as the minimum sample size tends to infinity.

(II) If $X$ is a bounded random variable, then $\bb{E}
[\wh{\bs{\mu}} - \mu]$ and $\bb{E} |\wh{\bs{\mu}} - \mu|^k, \; k =
1, 2, \cd$ tend to $0$ as the minimum sample size tends to infinity.

\eeT

See Appendix \ref{App_Unbiased_Gen} for a proof.

\sect{Computational Machinery}  \la{Computing}

\subsection{Bisection Coverage Tuning}

To avoid prohibitive burden of computational complexity in the design process, we shall focus on a class of multistage sampling schemes for
which the coverage probability can be adjusted by a single parameter $\ze$.  Such a parameter $\ze$ is referred to as the {\it coverage tuning
parameter} in this paper to convey the idea that $\ze$ is used to ``tune'' the coverage probability to meet the desired confidence level. As
will be seen in the sequel, we are able to construct a class of multistage sampling schemes such that the coverage probability can be ``tuned''
to ensure prescribed level of confidence by making the coverage tuning parameter sufficiently small. One great advantage of our sampling schemes
is that the tuning can be accomplished by a bisection search method. To apply a bisection method, it is required to determine whether the
coverage probability for a given $\ze$ is exceeding the prescribed level of confidence. Such a task is addressed in the following subsection.

\subsection{Adaptive Maximum Checking} \la{AMCA}

A wide class of computational problems depends on the following
critical subroutine:

Determine whether a function $C(\se)$ is smaller than a prescribed
number $\de$ for every value of $\se$ in interval $[ \udl{\se},
\ovl{\se}]$.

In many situations, it is impossible or very difficult to evaluate $C(\se)$ for every value of $\se$ in interval $[ \udl{\se}, \ovl{\se}]$,
since the interval may contain infinitely many or an extremely large number of values. To overcome such an issue of computational complexity, we
have developed Adapted Branch and Bound (ABB) algorithms in Appendix \ref{BBA}, which is a generalization of our previous ABB algorithm to the
multidimensional parameter space. To further reduce computational complexity, we propose an {\it Adaptive Maximum Checking Algorithm},
abbreviated as AMCA, to determine whether the maximum of $C(\se)$ over $[ \udl{\se}, \ovl{\se}]$ is less than $\de$. The only assumption
required for our AMCA is that, for any interval $[a, b] \subseteq [ \udl{\se}, \ovl{\se}]$, it is possible to compute an upper bound $\ovl{C}
(a, b)$ such that $C(\se) \leq \ovl{C} (a, b)$ for any $\se \in [a, b]$ and that the upper bound converges to $C(\se)$ as the interval width $b
- a$ tends to $0$.

Our backward AMCA proceeds as follows:

\bei

\item Choose initial step size $d > \eta$.

\item Let $F \leftarrow 0, \; T \leftarrow 0$ and $b \leftarrow \ovl{\se}$.

\item While $F = T = 0$, do the following:

    \bei

       \item Let $\tx{st} \leftarrow 0$ and $\ell \leftarrow 2$;

        \item While $\tx{st} = 0$, do the following:

             \bei

               \item Let $\ell \leftarrow \ell - 1$ and $d \leftarrow d  2^\ell$.
               \item If $b - d > \udl{\se}$, let $a \leftarrow b - d$ and $T
               \leftarrow 0$.   Otherwise, let $a \leftarrow \udl{\se}$ and $T \leftarrow 1$.

               \item If $\ovl{C} (a, b) < \de$, let $\tx{st} \leftarrow 1$ and $b \leftarrow a$.

                \item If $d < \eta$, let $\tx{st} \leftarrow 1$ and $F \leftarrow 1$.
            \eei

     \eei

\item Return $F$.
\eei

The output of our backward AMCA is a binary variable $F$ such that
``$F = 0$'' means ``$C(\se) < \de$'' and ``$F = 1$'' means ``$C(\se)
\geq \de$''. An intermediate variable $T$ is introduced in the
description of AMCA such that ``$T = 1$'' means that the left
endpoint of the interval is reached. The backward AMCA starts from
the right endpoint of the interval (i.e., $b = \ovl{\se}$) and
attempts to find an interval $[a, b]$ such that $\ovl{C} (a, b) <
\de$. If such an interval is available, then, attempt to go backward
to find the next consecutive interval with twice width. If doubling
the interval width fails to guarantee $\ovl{C} (a, b) < \de$, then
try to repeatedly cut the interval width in half to ensure that
$\ovl{C} (a, b) < \de$. If the interval width becomes smaller than a
prescribed tolerance $\eta$, then AMCA declares that ``$F = 1$''.
For our relevant statistical problems, if $C(\se) \geq \de$ for some
$\se \in [\udl{\se}, \ovl{\se}]$, it is sure that ``$F = 1$'' will
be declared.  On the other hand, it is possible that ``$F = 1$'' is
declared even though $C(\se) < \de$ for any $\se \in [\udl{\se},
\ovl{\se}]$. However, such situation can be made extremely rare and
immaterial if we choose $\eta$ to be a very small number.  Moreover,
this will only introduce negligible conservativeness in the
evaluation of coverage probabilities of random intervals if we
choose $\eta$ to be sufficiently small (e.g., $\eta = 10^{-15}$).

To see the practical importance of AMCA in our statistical problems, consider the construction of a sequential random interval with lower limit
$\mscr{L}( \wh{\bs{\se}}, \mbf{n} )$ and upper limit $\mscr{U}( \wh{\bs{\se}}, \mbf{n} )$ such that $\Pr \{ \mscr{L}( \wh{\bs{\se}}, \mbf{n} ) <
\se  < \mscr{U}( \wh{\bs{\se}}, \mbf{n} ) \mid \se \}
> 1 - \de$, or equivalently, $C (\se)  < \de$ for any $\se \in [ \udl{\se},
\ovl{\se}]$,  where $C(\se) = \Pr \{ \mscr{L}( \wh{\bs{\se}}, \mbf{n} ) \geq \se \mid \se \} + \Pr \{ \mscr{U}( \wh{\bs{\se}}, \mbf{n} ) \leq
\se \mid \se \}$ and $[ \udl{\se}, \ovl{\se}]$ is a subset of $\Se$.  For our statistical problems, $C(\se)$ is dependent on the coverage tuning
parameter $\ze$.   By choosing small enough $\ze$, it is possible to ensure $C(\se) < \de$ for any $\se \in [ \udl{\se}, \ovl{\se}]$. To avoid
unnecessary conservativeness, it is desirable to obtain $\ze$ as large as possible such that $C(\se) < \de$ for any $\se \in [ \udl{\se},
\ovl{\se}]$. This can be accomplished by a computational approach. Clearly, an essential step is to determine, for a given value of $\ze$,
whether $C(\se) < \de$ holds for any $\se \in [ \udl{\se}, \ovl{\se}]$.  Here, $C(\se)$ is defined as the complementary probability of coverage.

In the case that $\Se$ is a discrete set, special care needs for $d$
to ensure that $a$ and $b$ are numbers in $\Se$. The backward AMCA
can be easily modified as forward AMCA.

\subsection{Adapted Branch and Bound}

Actually, we first proposed ABB algorithm in the 6th version of this paper, published in arXiv on March 2, 2009, to determine whether a function
$C(\se)$ is smaller than a prescribed number $\de$ for every value of $\se$ in interval $[ \udl{\se}, \ovl{\se}]$.  It should be noted that the
ABB algorithm can be applied to determine whether a multivariate function $C(\se)$ is smaller than a prescribed number $\de$ for every value of
$\se$ in a multidimensional region, while AMCA limits its applications to univariate functions.  In applications, for purpose of obtaining the
exact values of minimum or maximum of $C(\se)$ for $\se$ in a subset of the parameter space, the B\&B algorithms described in Appendix \ref{BBA}
can be used. For example, it is of significant interests to compute the minimum coverage probability of a confidence interval. In this problem
area, we advocate the use of the B\&B algorithms, which are extremely powerful.

\subsection{Interval Bounding} \la{ITVB}

   Given that the levels of relative
precision of computation are equivalent for different methods and
that the complementary coverage probabilities are much smaller than
the coverage probabilities, the numerical error will be
significantly smaller if we choose to evaluate the complementary
coverage probabilities in the design of stopping rules. Therefore,
for computational accuracy, we propose to evaluate the complementary
coverage probabilities of the form $\Pr \{ \mscr{L}( \wh{\bs{\se}},
\mbf{n} ) \geq \se \; \tx{or} \; \mscr{U}( \wh{\bs{\se}}, \mbf{n} )
\leq \se \mid \se \}$.  By virtue of statement (III) of Theorem
\ref{Main_Bound_Gen}, we have \bel &  & \Pr \{ \mscr{L}(
\wh{\bs{\se}}, \mbf{n} ) \geq \se \; \tx{or} \; \mscr{U}(
\wh{\bs{\se}}, \mbf{n} ) \leq \se \mid \se \} \geq \Pr \{ b \leq
\mscr{L}( \wh{\bs{\se}}, \mbf{n} ) \mid a \} + \Pr \{ a \geq
\mscr{U}(
\wh{\bs{\se}}, \mbf{n} ) \mid b \}, \qqu \qqu \la{ineqa8} \\
&  & \Pr \{ \mscr{L}( \wh{\bs{\se}}, \mbf{n}) \geq \se \; \tx{or} \; \mscr{U}( \wh{\bs{\se}}, \mbf{n} ) \leq \se \mid \se \} \leq \Pr \{ a \leq
\mscr{L}( \wh{\bs{\se}}, \mbf{n} ) \mid b \} + \Pr \{ b \geq \mscr{U}( \wh{\bs{\se}}, \mbf{n} ) \mid a \} \la{ineqb8}\eel for any $\se \in [a,
b]$ provided that \be \la{narrow} \{ a \leq \mscr{L}( \wh{\bs{\se}}, \mbf{n} ) \} \subseteq \{ \wh{\bs{\se}} \geq b   \}, \qqu \{ b \geq
\mscr{U}( \wh{\bs{\se}}, \mbf{n} ) \} \subseteq \{ \wh{\bs{\se}} \leq a \}. \ee For many problems,  if interval $[a, b]$ is narrow enough, then,
condition (\ref{narrow}) can be satisfied and the upper and lower bounds of $\Pr \{ \mscr{L}( \wh{\bs{\se}}, \mbf{n} ) \geq \se \; \tx{or} \;
\mscr{U}( \wh{\bs{\se}}, \mbf{n} ) \leq \se \mid \se \}$ in (\ref{ineqa8}) and (\ref{ineqb8}) can be used to determine whether $\Pr \{ \mscr{L}(
\wh{\bs{\se}}, \mbf{n} ) \geq \se \; \tx{or} \; \mscr{U}( \wh{\bs{\se}}, \mbf{n} ) \leq \se \mid \se \} \leq \de$ for any $\se \in [a, b]$. This
suggests an alternative approach for constructing sequential
random intervals to guarantee prescribed confidence level for any $\se \in
[\udl{\se}, \ovl{\se}]$, where $[\udl{\se}, \ovl{\se}]$ is a subset of parameter space $\Se$.  The basis idea is as follows:

(i) Construct sampling scheme such that the probabilities $\Pr \{
\se \leq \mscr{L}( \wh{\bs{\se}}, \mbf{n} ) \mid \se \}$ and $\Pr \{
\se \geq \mscr{U}( \wh{\bs{\se}}, \mbf{n} ) \mid \se \}$ can be
adjusted by $\ze$.

(ii) Partition $[\udl{\se}, \ovl{\se}]$ as small subintervals $[a,
b]$ such that (\ref{ineqa8}) and (\ref{ineqb8}) can be used to
determine whether $\Pr \{ \mscr{L}( \wh{\bs{\se}}, \mbf{n} ) \geq
\se \; \tx{or} \; \mscr{U}( \wh{\bs{\se}}, \mbf{n} ) \leq \se \mid
\se \} \leq \de$ for any $\se \in [a, b]$.

It should be noted that, in some cases, especially for point
estimation with precision requirements, we can use statement (IV) of
Theorem \ref{Main_Bound_Gen} for the purpose of interval bounding as
above.

\subsection{Recursive Computation} \la{Recursive}

As will be seen in the sequel, for most multistage sampling plans
with deterministic sample sizes $n_1, n_2, \cd, n_s$ for estimating
parameters of discrete variables, the probabilistic terms involving
$\wh{\bs{\se}}, \; \mbf{n}$ or $\wh{\bs{\se}}_\ell, \; n_\ell$ can
usually be expressed as a summation of terms $\Pr \{ K_i \in
\mscr{K}_i, \; i = 1, \cd, \ell  \}, \; \ell = 1, \cd, s$,  where
$K_\ell = \sum_{i = 1}^{n_\ell} X_i$ and $\mscr{K}_i$ is a subset of
integers. The calculation of such terms can be performed by virtue
of the following recursive relationship: \bel & & \Pr \{ K_i \in
\mscr{K}_i, \; i
= 1, \cd, \ell; \; K_{\ell + 1} = k_{\ell + 1} \} \nonumber\\
&  = & \sum_{k_\ell \in \mscr{K}_\ell} [ \Pr \{ K_i \in \mscr{K}_i, \; i = 1, \cd, \ell - 1; \; K_\ell = k_\ell  \} \nonumber \\
&   & \times \Pr \{ K_{\ell + 1} - K_\ell = k_{\ell + 1} - k_\ell \mid K_\ell = k_\ell; \;  K_i \in \mscr{K}_i, \; i = 1, \cd, \ell - 1 \} ],
\la{recur1}\eel where the computation of the conditional probability $\Pr \{ K_{\ell + 1} - K_\ell = k_{\ell + 1} - k_\ell \mid K_\ell = k_\ell;
\; K_i \in \mscr{K}_i, \; i = 1, \cd, \ell - 1 \}$ depends on specific estimation problems.  For estimating a binomial parameter $p$ with
deterministic sample sizes $n_1, n_2, \cd, n_s$, we have \bel &  & \Pr \{ K_{\ell + 1} -
K_\ell = k_{\ell + 1} - k_\ell \mid K_\ell = k_\ell; \;  K_i \in \mscr{K}_i, \; i = 1, \cd, \ell - 1 \} \nonumber\\
&  &  = \Pr \{ K_{\ell + 1} - K_\ell = k_{\ell + 1} - k_\ell \} \nonumber\\
&  &  = \bi{ n_{\ell + 1} - n_\ell }{ k_{\ell + 1} - k_\ell } p^{ k_{\ell + 1} - k_\ell } ( 1 - p )^{ n_{\ell + 1} - n_\ell -  k_{\ell + 1} +
k_\ell }. \la{rec8} \eel As an immediate consequence of (\ref{recur1}) and (\ref{rec8}), we have \be \la{for888} \Pr \{ K_i \in \mscr{K}_i, \; i
= 1, \cd, \ell; \; K_{\ell + 1} = k_{\ell + 1} \} = \nu ( k_{\ell + 1}, \ell + 1) \; p^{ k_{\ell + 1} } ( 1 - p )^{n_{\ell + 1} - k_{\ell + 1}},
\ee where $\nu (k, 1) = \bi{n_1}{k}$ for $k \in \mscr{K}_1$, and  \be \la{fracite} \nu (k, i) = \sum_{ k_{i - 1} \in \mcal{K}_{i - 1} \atop{ k +
n_{i - 1} - n_i \leq k_{i - 1} \leq k  } } \nu ( k_{i - 1}, i - 1) \bi{n_i - n_{i - 1}} { k - k_{i - 1} } \qqu \tx{for} \qu k \in \mscr{K}_i,
\qu 2 \leq i \leq \ell + 1. \ee  It should be noted that, in the case that $X$ is a Bernoulli variable,  the recursive relationship had been
used in \cite{sch} for designing hypothesis tests for drug screening.  In the special case of fully sequential sampling (i.e., the increment of
sample sizes is unity), (\ref{fracite}) reduces to the recursive formula given at the bottom of page 49 of \cite{frazen} for computing $\nu
(s)$.

 For estimating a Poisson parameter $\lm$
with deterministic sample sizes $n_1, n_2, \cd, n_s$, we have \bee &  & \Pr \{ K_{\ell + 1} - K_\ell = k_{\ell + 1} - k_\ell \mid K_\ell =
k_\ell; \; K_i \in \mscr{K}_i, \; i = 1, \cd, \ell - 1 \}\\
&  & \Pr \{ K_{\ell + 1} - K_\ell = k_{\ell + 1} - k_\ell \}\\
&  & = \f{ [ (n_{\ell + 1} - n_\ell ) \lm ]^{k_{\ell + 1} - k_\ell} \exp ( - (n_{\ell + 1} - n_\ell) \lm  ) } { (k_{\ell + 1} - k_\ell)!  }.
\eee For estimating the proportion, $p$, of a finite population using multistage sampling schemes described in Section \ref{finite_proportion},
we have {\small \be \la{confinite} \Pr \{ K_{\ell + 1} - K_\ell = k_{\ell + 1} - k_\ell \mid K_\ell = k_\ell; \;  K_i \in \mscr{K}_i, \; i = 1,
\cd, \ell - 1 \} = \f{ \bi{p N - k_\ell}{k_{\ell + 1} - k_\ell} \bi{N - p N - n_\ell + k_\ell} { n_{\ell + 1} - n_\ell - k_{\ell + 1} + k_\ell }
} { \bi{N - n_\ell}{n_{\ell + 1} - n_\ell} }, \ee} where the sample sizes are deterministic numbers $n_1, n_2, \cd, n_s$.  The conditional
probability in (\ref{confinite}) can be viewed as the probability of seeing $k_{\ell + 1} - k_\ell$ units having a certain attribute in the
course of drawing $n_{\ell + 1} - n_\ell$ units, based on a simple sampling without replacement, from a population of $N - n_\ell$ units, among
which  $p N - k_\ell$ units having the attribute.

It should be noted that such idea of recursive computation can be applied to general multistage sampling plans with random sample sizes
$\mbf{n}_1, \mbf{n}_2, \cd, \mbf{n}_s$. Moreover, the domain truncation technique described in the next subsection can be used to significantly
reduce computation.

\subsection{Domain Truncation}

The bounding methods described in the previous subsection reduce the computational problem of designing a multistage sampling scheme to the
evaluation of low-dimensional summation or integration. Despite the reduction of dimensionality, the associated computational complexity is
still high because the domain of summation or integration is large.  The truncation techniques recently established in \cite{Chen1} have the
power to considerably simplify the computation by reducing the domain of summation or integration to a much smaller subset.  The following
result derived from a similar method as that of \cite{Chen1}, shows that the truncation can be done with controllable error.

\beT \la{Trun_THM5} Let $\eta \in (0, 1)$. Let $\udl{\se}_\ell, \; \ovl{\se}_\ell, \; \ell = 1, \cd, s$ be real numbers such that $\Pr \{
\udl{\se}_\ell \leq \wh{\bs{\se}}_\ell \leq \ovl{\se}_\ell \; \tx{for} \; \ell = 1, \cd, s \} \geq 1 - \eta$.  Assume that there exist subsets
of real numbers $\mscr{A}_\ell, \; \ell = 1, \cd, s$ such that $\{ \bs{l} = \ell \} = \{ \wh{\bs{\se}}_i \in \mscr{A}_i \; \tx{for} \; 1 \le i
\leq \ell \}$ for $\ell = 1, \cd, s$.  Then, {\small \[ \Pr \{ \mscr{W} (\wh{\bs{\se}}, \mbf{n} ) \in \mscr{R} \} - \eta \leq \sum_{\ell = 1}^s
\Pr \{ \mscr{W} (\wh{\bs{\se}}_\ell, \mbf{n}_\ell ) \in \mscr{R} \; \tx{and} \; \wh{\bs{\se}}_i \in \mscr{B}_i \; \tx{for} \; 1 \leq i \leq \ell
\} \leq \Pr \{ \mscr{W} (\wh{\bs{\se}}, \mbf{n} ) \in \mscr{R} \}, \]} where $\mscr{B}_\ell = \{ \vse \in \mscr{A}_\ell: \; \udl{\se}_\ell \leq
\vse \leq \ovl{\se}_\ell  \}$ for $\ell = 1, \cd, s$.  \eeT

To determine numbers $\udl{\se}_\ell, \; \ovl{\se}_\ell, \; \ell = 1, \cd, s$ such that  $\Pr \{ \udl{\se}_\ell \leq \wh{\bs{\se}}_\ell \leq
\ovl{\se}_\ell \; \tx{for} \; \ell = 1, \cd, s \} \geq 1 - \eta$, we can follow a similar method as that of \cite{Chen1}.

As an illustration the truncation technique, consider probabilistic terms like $\Pr \{ \mscr{W} (\wh{\bs{\se}}, \mbf{n} ) \in \mscr{R} \}$
involved in a multistage sampling scheme. If $\udl{\se}_\ell$ and $\ovl{\se}_\ell$ can be found such that $\Pr \{ \udl{\se}_\ell \leq
\wh{\bs{\se}}_\ell \leq \ovl{\se}_\ell \} \geq 1 - \f{\eta}{s}$ for $\ell = 1, \cd, s$, then, by Bonferroni's inequality, {\small \be \la{good8}
\Pr \{ \mscr{W} (\wh{\bs{\se}}, \mbf{n} ) \in \mscr{R} \} - \eta \leq \sum_{\ell = 1}^s \Pr \{ \mscr{W} (\wh{\bs{\se}}_\ell, \mbf{n}_\ell ) \in
\mscr{R}, \; \udl{\se}_\ell \leq \wh{\bs{\se}}_\ell \leq \ovl{\se}_\ell, \; \bs{l} = \ell  \} \leq \Pr \{ \mscr{W} (\wh{\bs{\se}}, \mbf{n} ) \in
\mscr{R} \}, \ee} where $\bs{l}$ denotes the index of stage at the termination of the sampling process as before. For most multistage sampling
plans for estimating parameters of discrete variables, the probabilities $\Pr \{ \mscr{W} (\wh{\bs{\se}}_\ell, \mbf{n}_\ell ) \in \mscr{R}, \;
\udl{\se}_\ell \leq \wh{\bs{\se}}_\ell \leq \ovl{\se}_\ell, \; \bs{l} = \ell \}, \; \ell = 1, \cd, s$ can be evaluated recursively as described
in Section \ref{Recursive}.
 Specially, we can apply (\ref{good8}) to multistage sampling plans
for estimating the parameter $p$ of a Bernoulli random variable $X$
such that $\Pr \{ X = 1 \} = 1 - \Pr \{X = 0 \} = p \in (0, 1)$. Let
$X_1, X_2, \cd$ be i.i.d. samples of $X$. Suppose that the sampling
plan has $s$ stages with deterministic sample sizes $n_1, \cd, n_s$.
Let $\mcal{K}_\ell \subseteq \{0, 1, \cd, n_\ell \}$ and
$\mcal{K}_\ell^c \subseteq \{0, 1, \cd, n_\ell \} \setminus
\mcal{K}_\ell$ for $\ell = 1, \cd, s$.   Suppose the decision
variables  are defined such that
\[
\{  \bs{l} = \ell \} = \li \{ \sum_{i = 1}^{n_j} X_i \in \mcal{K}_j,
\; 1 \leq j < \ell; \; \sum_{i = 1}^{n_\ell} X_i \in \mcal{K}_\ell^c
\ri \}
\]
for $\ell = 1, \cd, s$ .  Define $\wh{\bs{p}}_\ell = \f{ \sum_{i =
1}^{n_\ell} X_i } {n_\ell}$ as an estimator of $p$ for $\ell = 1,
\cd, s$. As before, define $\wh{\bs{p}} = \f{ \sum_{i = 1}^{\mbf{n}}
X_i } {\mbf{n}}$, where $\mbf{n}$ is the number of samples at the
termination of sampling.  Define \bee &  & \psi (k, 1) =
\bi{n_1}{k} \qqu \tx{for} \qu  k \in \mcal{K}_1 \cup \mcal{K}_1^c,\\
&  & \psi (k, \ell) = \sum_{ k_{\ell - 1} \in \mcal{K}_{\ell - 1}
\atop{ k + n_{\ell - 1} - n_\ell \leq k_{\ell - 1} \leq k  } } \psi(
k_{\ell - 1}, \ell - 1) \bi{n_\ell - n_{\ell - 1}} { k - k_{\ell -
1} } \qqu \tx{for} \qu k \in \mcal{K}_\ell \cup \mcal{K}_\ell^c, \qu
2 \leq \ell \leq s. \eee Let $\udl{k}_\ell \leq \ovl{k}_\ell$ be
integers from the set $\{0, 1, \cd, n_\ell \}$ such that {\small
$\Pr \li \{ \f{\udl{k}_\ell}{n_\ell} \leq \wh{\bs{p}}_\ell \leq \f{
\ovl{k}_\ell }{n_\ell} \mid p \ri \} \geq 1 - \f{\eta}{s}$} for
$\ell = 1, \cd, s$.   Let {\small $\mscr{V}_\ell = \li \{ k \in
\mcal{K}_\ell^c: \udl{k}_\ell \leq k \leq \ovl{k}_\ell, \; \mscr{W}
\li ( \f{k}{n_\ell},  n_\ell \ri ) \in \mscr{R} \ri \}$} for $\ell =
1, \cd, s$.   Then, by (\ref{for888}), (\ref{good8}) and the
definition of $\psi(.,.)$,
\[
\sum_{\ell = 1}^s \sum_{k \in \mscr{V}_\ell } \psi (k, \ell) \; p^k
(1 - p)^{n_\ell - k} \leq \Pr \{ \mscr{W} ( \wh{\bs{p}},  \mbf{n} )
\in \mscr{R} \mid p \} \leq \eta + \sum_{\ell = 1}^s \sum_{k \in
\mscr{V}_\ell } \psi (k, \ell) \; p^k (1 - p)^{n_\ell - k}.
\]

\subsection{Consecutive-Decision-Variable Bounding} \la{consec}

One major problem in the design and analysis of multistage sampling schemes is the high-dimensional summation or integration involved in the
evaluation of probabilities.   For instance, a basic problem is to evaluate the coverage probabilities involving $\wh{\bs{\se}}$ and $\mbf{n}$.
Another example is to evaluate the distribution or the expectation of sample number $\mbf{n}$. Clearly, $\wh{\bs{\se}}$ depends on random
samples $X_1, \cd, X_{\mbf{n}}$. Since the sample number $\mbf{n}$ can assume very large values, the computational complexity associated with
the high-dimensionality can be a prohibitive burden to modern computers. In order to break the curse of dimensionality, we propose to obtain
tight bounds for those types of probabilities. In this regard, we have \beT \la{Thm_CDV} Let $\mscr{W} (.,.)$ be a bivariate function. Let
$\mscr{R}$ be a subset of real numbers. Then, {\small \bee & & \Pr \li \{ \mscr{W} (\wh{\bs{\se}}, \mbf{n}) \in \mscr{R} \ri \} \leq \sum_{\ell
= 1}^s \Pr \li \{ \mscr{W} ( \wh{\bs{\se}}_\ell, \mbf{n}_\ell) \in \mscr{R}, \; \bs{D}_\ell = 1 \; \tx{and} \; \bs{D}_{j} = 0 \; \tx{for} \;
\max(1, \ell -
r) \leq j < \ell \ri \}, \\
&   & \Pr \li \{ \mscr{W} (\wh{\bs{\se}}, \mbf{n}) \in \mscr{R} \ri \} \geq 1 - \sum_{\ell = 1}^s \Pr \li \{ \mscr{W} ( \wh{\bs{\se}}_\ell,
\mbf{n}_\ell) \notin \mscr{R}, \; \bs{D}_\ell = 1 \; \tx{and} \; \bs{D}_{j} = 0 \; \tx{for} \; \max(1, \ell - r) \leq j < \ell \ri \} \eee} for
$0 \leq r < s$. Moreover, \bee & & \Pr \{ \bs{l} > \ell \} \leq \Pr \{ \bs{D}_\ell = 0, \; \bs{D}_{j} =
0 \; \tx{for} \; \max(1, \ell - r) \leq j < \ell \},\\
&   & \Pr \{ \bs{l} > \ell \} \geq  1 - \sum_{j = 1}^\ell \Pr \{ \bs{D}_j = 1, \; \; \bs{D}_{i} = 0 \; \tx{for} \; \max(1, j - r) \leq i < j \}
\eee for $1 \leq \ell \leq s$ and $0 \leq r < s$. Furthermore, if the number of available samples at the $\ell$-th stage is a deterministic
number $n_\ell$ for $1 \leq \ell \leq s$, then $\bb{E} [ \mathbf{n} ]  = n_1 + \sum_{\ell = 1}^{s - 1} \; (n_{\ell + 1} - n_{\ell}) \; \Pr \{
\bs{l} > \ell \}$. \eeT

\bsk

See Appendix \ref{App_Thm_CDV} for a proof.  As can be seen from Theorem \ref{Thm_CDV}, the bounds are constructed by summing up probabilistic
terms involving one or multiple consecutive decision variables (CDV). Such general technique is referred to as CDV bounding.  A particular
interesting special case of CDV method is to construct bounds with every probabilistic term involving consecutive decision variables (i.e., $r =
1$ in Theorem \ref{Thm_CDV}). Such method is referred to as {\it double-decision-variable}  or DDV bounding for brevity. Similarly, the bounds
with each probabilistic term involving a single decision variable are referred to as  {\it single-decision-variable} bounds or  SDV bounds
(i.e., $r = 0$ in Theorem \ref{Thm_CDV}). Our computational experiences indicate that the bounds in Theorem \ref{Thm_CDV} become very tight as
the spacing between sample sizes increases. As can be seen from Theorem \ref{Thm_CDV}, DDV bounds are tighter than SDV bounds.  Needless to say,
the tightness of bounds is achieved at the price of computational complexity.  The reason that such bounding methods allow for powerful
dimension reduction is that, for many important estimation problems, $\bs{D}_{\ell - 1}, \; \bs{D}_{\ell}$ and $\wh{\bs{\se}}_\ell$ can be
expressed in terms of two independent variables $U$ and $V$. For instance, for the estimation of a binomial parameter, it is possible to design
a multistage sampling scheme such that $\bs{D}_{\ell - 1}, \; \bs{D}_{\ell}$ and $\wh{\bs{\se}}_\ell$ can be expressed in terms of $U = \sum_{i
=1}^{\mbf{n}_{\ell - 1}} X_i$ and $V = \sum_{i = \mbf{n}_{\ell - 1} + 1}^{\mbf{n}_\ell} X_i$. For the double decision variable method, it is
evident that $U$ and $V$ are two independent binomial random variables and accordingly the computation of probabilities such as $\Pr \{ \mscr{W}
(\wh{\bs{\se}}, \mbf{n})  \in \mscr{R} \}$ and $\Pr \{ \bs{l} > \ell \}$ can be reduced to two-dimensional problems. Clearly, the dimension of
these computational problems can be reduced to one if the single-decision-variable method is employed. As will be seen in the sequel, DDV bounds
can be shown to be asymptotically tight for a large class of multistage sampling schemes. Moreover, our computational experiences indicate that
SDV bounds are not very conservative.

For computational simplicity, the DDV upper bound can be relaxed as follows: \bel \Pr \{ \mscr{W} (\wh{\bs{\se}}, \mbf{n})  \in \mscr{R}   \} &
\leq &
\sum_{ \ell = 1 }^s \min \li [ \Pr \{ \bs{D}_{\ell - 1} = 0 \}, \;
\Pr \{ \mscr{W} (\wh{\bs{\se}}_\ell, \mbf{n}_\ell)  \in \mscr{R}, \; \bs{D}_{\ell} = 1 \} \ri ] \la{relaxedA}\\
& \leq & \sum_{ \ell = 1 }^s \min \li [ \Pr \{ \bs{D}_{\ell - 1} = 0 \}, \; \Pr \{ \mscr{W} (\wh{\bs{\se}}_\ell, \mbf{n}_\ell)  \in \mscr{R} \},
\; \Pr \{ \bs{D}_{\ell} = 1 \} \ri ].  \qqu \la{relaxedB} \eel In situations that for $\ell = 1, \cd, s$, the decision variable $\bs{D}_{\ell}$
is defined in terms of $\wh{\bs{\se}}_\ell$ and $\mbf{n}_\ell$,  such relaxed bounds further reduce the dimensionality of computation, since the
bounds can be evaluated by summing up probabilistic terms, among which each involves $\wh{\bs{\se}}_\ell$ and $\mbf{n}_\ell$ for a single stage.
To illustrate, consider the bounding of the complementary coverage probability $\Pr \{ \mscr{L}( \wh{\bs{\se}}, \mbf{n}) \geq \se \; \tx{or} \;
\mscr{U}( \wh{\bs{\se}}, \mbf{n} ) \leq \se \mid \se \}$ in the context associated with (\ref{ineqb8}).  For the purpose of reducing the
computational complexity, we can apply either (\ref{relaxedA}) or (\ref{relaxedB}) to relax the upper bound in the right side of (\ref{ineqb8}).
We shall first discuss the application of (\ref{relaxedA}).  Making use of (\ref{relaxedA}) with the event $\{ \mscr{W} (\wh{\bs{\se}}, \mbf{n})
\in \mscr{R} \}$ identified as $\{ a \leq \mscr{L}( \wh{\bs{\se}}, \mbf{n} ) \}$, we have \be \la{boundreA} \Pr \{ a \leq \mscr{L}(
\wh{\bs{\se}}, \mbf{n} ) \mid b \} \leq \sum_{ \ell = 1 }^s \min \li [ \Pr \{ \bs{D}_{\ell - 1} = 0 \mid b \}, \; \Pr \{ \mscr{L}(
\wh{\bs{\se}}_\ell, \mbf{n}_\ell ) \geq a, \; \bs{D}_{\ell} = 1 \mid b \} \ri ], \ee where the probabilistic terms $\Pr \{ \bs{D}_{\ell - 1} = 0
\mid b \}$ and $\Pr \{ \mscr{L}( \wh{\bs{\se}}_\ell, \mbf{n}_\ell ) \geq a, \; \bs{D}_{\ell} = 1 \mid b \}$ can be computed for all meaningful
values of index $\ell$. Similarly, identifying the event $\{ \mscr{W} (\wh{\bs{\se}}, \mbf{n}) \in \mscr{R} \}$ as $\{ \mscr{U}( \wh{\bs{\se}},
\mbf{n} ) \leq b \}$, we have \be \la{boundreB} \Pr \{ b \geq \mscr{U}( \wh{\bs{\se}}, \mbf{n} ) \mid a \} \leq \sum_{ \ell = 1 }^s \min \li [
\Pr \{ \bs{D}_{\ell - 1} = 0 \mid a \}, \; \Pr \{ \mscr{U}( \wh{\bs{\se}}_\ell, \mbf{n}_\ell ) \leq b, \; \bs{D}_{\ell} = 1 \mid a \} \ri ]. \ee
Adding up the upper bounds in (\ref{boundreA}) and (\ref{boundreB}) gives an upper bound for the complementary coverage probability $\Pr \{
\mscr{L}( \wh{\bs{\se}}, \mbf{n}) \geq \se \; \tx{or} \; \mscr{U}( \wh{\bs{\se}}, \mbf{n} ) \leq \se \mid \se \}$ for $\se \in [a, b]$. Clearly,
the computational complexity of such upper bound is low, since it does not rely on the recursive method as described in Section \ref{Recursive}.
Of course, the reduction of computation comes at the expense of conservatism.

In a similar manner, applying  (\ref{relaxedB}), we have {\small \bel  &  & \Pr \{ a \leq \mscr{L}( \wh{\bs{\se}}, \mbf{n} ) \mid b \} \leq
\sum_{ \ell = 1 }^s \min \li [ \Pr \{ \bs{D}_{\ell - 1} = 0 \mid b
\}, \; \Pr \{ \mscr{L}( \wh{\bs{\se}}_\ell, \mbf{n}_\ell ) \geq a \mid b \}, \; \Pr \{ \bs{D}_{\ell} = 1 \mid b \} \ri ], \qqu \la{boundreAA}\\
&  &  \Pr \{ b \geq \mscr{U}( \wh{\bs{\se}}, \mbf{n} ) \mid a \} \leq \sum_{ \ell = 1 }^s \min \li [ \Pr \{ \bs{D}_{\ell - 1} = 0 \mid a \}, \;
\Pr \{ \mscr{U}( \wh{\bs{\se}}_\ell, \mbf{n}_\ell ) \leq b \mid a \}, \; \Pr \{ \bs{D}_{\ell} = 1 \mid a \} \ri ]. \qqu \la{boundreBB} \eel}
Summing up the upper bounds in (\ref{boundreAA}) and (\ref{boundreBB}) yields an upper bound for the complementary coverage probability $\Pr \{
\mscr{L}( \wh{\bs{\se}}, \mbf{n}) \geq \se \; \tx{or} \; \mscr{U}( \wh{\bs{\se}}, \mbf{n} ) \leq \se \mid \se \}$ for $\se \in [a, b]$.

If one is willing to tolerate more conservatism for the reward of lower computational complexity, one could replace the probabilistic terms in
(\ref{boundreA}),  (\ref{boundreB}), (\ref{boundreAA}) and (\ref{boundreBB}) by their upper bounds derived from probabilistic inequalities such
as Chernoff-Hoeffding bounds \cite{Chernoff, Hoeffding} or likelihood bounds \cite{ChenR}.  To illustrate, consider the case that the sample
sizes $\mbf{n}_1, \cd, \mbf{n}_s$ are deterministic numbers and that for $\ell = 1, \cd, s$, $\wh{\bs{\se}}_\ell$ is the sample mean $\f{
\sum_{i = 1}^{\mbf{n}_\ell} X_i }{ \mbf{n}_\ell }$, where $X_i, \; i = 1, \cd, \mbf{n}_\ell$ are i.i.d. samples of $X$ parameterized by $\se \in
\Se$ such that $\bb{E} [ X ] = \se$.  Assume that the moment generating function of $X$ exists.  Let $E$ be an event defined in terms of
$\wh{\bs{\se}}_\ell$. To bound the probability, $\Pr \{ E \mid \vse \}$, of event  $E$ with the associated parameter $\se$ assuming a value
$\vse \in \Se$,  one can seek a number $c$ such that $E \subseteq \{ \wh{\bs{\se}}_\ell < c < \vse \}$ or $E \subseteq \{ \wh{\bs{\se}}_\ell
> c
> \vse \}$.  If $E \subseteq \{ \wh{\bs{\se}}_\ell < c < \vse \}$ is true, then the probability $\Pr \{ E \mid \vse \}$ is bounded from above by
the Chernoff-Hoeffding bound of $\Pr \{ \wh{\bs{\se}}_\ell < c < \vse \mid \vse \}$.  Similarly, if $E \subseteq \{ \wh{\bs{\se}}_\ell > c >
\vse \}$ is true, then the probability $\Pr \{ E \mid \vse \}$ is bounded from above by the Chernoff-Hoeffding bound of $\Pr \{
\wh{\bs{\se}}_\ell > c > \vse \mid \vse \}$.

\subsection{Triangular Partition} \la{Tri_par}

As can be seen from the preceding discussion, by means of the
double-decision-variable method,  the design of multistage sampling
schemes may be reduced to the evaluation of probabilities of the
form $\Pr \{ (U,V) \in \mscr{G} \}$, where $U$ and $V$ are
independent random variables, and $\mscr{G} = \{ (u,v) : \; a \leq u
\leq b, \; c \leq v \leq d, \; e \leq u + v \leq f \}$ is a
two-dimensional domain. It should be noted that such a domain can be
fairly complicated. It can be an empty set or a polygon with $3$ to
$6$ sides. Therefore, it is important to develop a systematic method
for computing $\Pr \{ (U,V) \in \mscr{G} \}$. For this purpose, we
have

\beT \la{polygon_split} Let $a \leq b, \; c \leq d$ and $e \leq f$.
Let $\ovl{e} = \max \{ e, a + c \}, \; \udl{f} = \min \{ f, b + d
\}, \; \underline{u} = \max \{a,  \ovl{e} - d \}, \; \overline{u} =
\min \{b, \udl{f} - c \} , \; \underline{v} = \max \{c,  \ovl{e} - b
\}$ and $\overline{v} = \min \{d,  \udl{f} - a \}$. Then, for any
independent random variables $U$ and $V$, \bee \Pr \{ (U,V) \in
\mscr{G} \} & =  &
\Pr \{\underline{u} \leq U \leq \overline{u} \} \Pr \{\underline{v} \leq V \leq \overline{v} \}\\
&   & -  \Pr \{U \leq \overline{u}, \; V \leq \overline{v},
\;  U + V > \udl{f} \} - \Pr \{U \geq \underline{u}, \; V \geq
\underline{v}, \; U + V < \ovl{e} \}. \eee \eeT

\bsk

The goal of using Theorem \ref{polygon_split} is to separate
variables and thus reduce computation. As can be seen from Theorem
\ref{polygon_split}, random variables $U$ and $V$ have been
separated in the product and thus the dimension of the corresponding
computation is reduced to one.  The last two terms on the left side
of equality are probabilities that $(U, V)$ is included in
rectangled triangles. The idea of separating variables can be
repeatedly used by partitioning rectangled triangles as smaller
rectangles and rectangled triangles. Specifically, if $U$ and $V$
are discrete random variables assuming integer values, we have
{\small \bel &  & \Pr \{ U \geq i, \; V \geq j, \; U + V \leq k \}
\; = \; \Pr \li \{ i \leq U \leq \li \lf \f{k + i - j}{2} \ri \rf
\ri \} \; \Pr \li \{ j \leq V < \li \lc \f{ k -
i + j}{2} \ri \rc \ri \} \nonumber\\
&  & \qu + \Pr \li \{ U > \li \lf \f{k + i - j}{2} \ri \rf, \; V
\geq j, \; U + V \leq k \ri \}  +  \Pr \li \{ U \geq i, \; V \geq
\li \lc \f{ k - i + j}{2} \ri \rc, \; U + V \leq k \ri \} \qu \qu
\qu \qu \la{m1} \eel} for integers $i, \; j$ and $k$ such that $i +
j \leq k$;  and {\small \bel &  & \Pr \{ U \leq i, \; V \leq j, \; U
+ V \geq k  \} \; = \; \Pr \li \{ \li \lc \f{k + i - j}{2} \ri \rc
\leq U \leq i \ri \} \; \Pr \li \{ \li \lf \f{ k - i
+ j}{2} \ri \rf < V \leq j \ri \} \nonumber\\
&   & \qu + \Pr \li \{ U \leq i, \; V \leq \li \lf \f{k - i + j}{2}
\ri \rf, \; U + V \geq k \ri \}  + \Pr \li \{ U < \li \lc \f{ k + i
- j}{2} \ri \rc, \; V \leq j, \; U + V \geq k \ri \} \qu \qu \qu \qu
\la{m3} \eel} for integers $i, \; j$ and $k$ such that $i + j \geq
k$. It is seen that the terms in (\ref{m1}) and (\ref{m3})
correspond to probabilities that $(U, V)$ is included in rectangled
triangles. Hence, the above method  of triangular partition can be
repeatedly applied.  For the sake of efficiency, we can save the
probabilities that $U$ and $V$ are respectively included in the
intervals corresponding to the rectangular sides of a parent
triangle, then when partitioning this triangle, it suffices to
compute the probabilities that $U$ and $V$ are included in the
intervals corresponding to two orthogonal sides of the  smaller
rectangle. The probabilities that $U$ and $V$ are included in the
intervals corresponding to the rectangular sides of the smaller
triangles can be readily obtained from the results of the  smaller
rectangle and the record of the probabilities for the parent
triangle.  This trick can be repeatedly used to save computation.

Since a crucial step in designing a sampling scheme is to compare
the coverage probability with a prescribed level of confidence, it
is useful to compute upper and lower bounds of the probabilities
that $U$ and $V$ are covered by a triangular domain. As the
triangular partition goes on, the rectangled triangles become
smaller and smaller. Clearly, the upper bounds of the probabilities
that $(U, V)$ is included in rectangled triangles can be obtained by
inequalities {\small
\[ \Pr \{ U \geq i, \; V \geq j, \; U + V \leq k \} \leq \Pr \{ i
\leq U \leq k - j \} \Pr \{ j \leq V \leq k - i \}, \qqu \qqu \qqu
\qqu
\]
\[
\Pr \{  U \leq i, \; V \leq j, \; U + V \geq k \} \leq \Pr \{ k - j
\leq U \leq i \} \Pr \{ k - i \leq V \leq j \}. \; \qqu \qqu \qqu
\qqu
\]}
Of course, the lower bounds can be taken as $0$.  As the triangular
partition goes on, the rectangled triangles become smaller and
smaller and accordingly such bounds becomes tighter. To avoid the
exponential growth of the number of rectangled triangles, we can
split the rectangled triangle with the largest gap between upper and
lower bounds in every triangular partition.

\subsection{Interval Splitting}

In the design of sampling schemes and other applications, it is a
frequently-used routine to evaluate the probability that a random
variable is bounded in an interval.  Note that, for most basic
random variables,  the probability mass (or density) functions
$f(.)$ possess nice concavity or convexity properties.  In many
cases, we can readily compute inflexion points which can be used to
partition the interval as subintervals such that $f(.)$ is either
convex or concave in each subinterval. By virtue of concavity or
convexity, we can calculate the upper and lower bounds of the
probability that the random variable is included in a subinterval.
The overall upper and lower bounds of the probability that the
random variable is included in the initial interval can be obtained
by summing up the upper and lower bounds for all subintervals
respectively. The gap between the overall upper and lower bounds can
be reduced by repeatedly partitioning the subinterval with the
largest gap of upper and lower bounds. This strategy is referred to
as {\it interval splitting} in this paper.

For a discrete random variable with probability mass function
$f(k)$, we can apply the following result to compute upper and
 lower bounds of $\sum_{k = a}^b f(k)$ over subinterval $[a, b]$.

\beT  \la{discrete_sum} Let $a < b$ be two integers. Define {\small
$r_a = \f{ f(a + 1) }{ f(a) }, \; r_b = \f{ f(b - 1) } { f (b)}, \;
r_{a, b} = \f{ f(a) }{ f(b) }$} and {\small $j = a + \f{ b - a - (1
- r_{a,b}) (1 - r_b )^{-1} } { 1 + r_{a, b} ( 1 - r_a ) ( 1 - r_b
)^{-1} }$}. Define {\small $\al (i) = (i + 1 - a) \li [ 1 + \f{(i -
a) (  r_a  - 1 ) }{2}  \ri ]$} and {\small $\ba (i) = (b - i) \li [
1 + \f{( b - i - 1) ( r_b - 1 ) }{2}  \ri ]$}.  The following
statements hold true:

 (I): If $f(k+1) - f(k) \leq f(k) - f(k-1)$  for $a < k < b$, then
\be \la{concavity} \f{ (b-a+1) [f(a) + f(b)] }{2}  \leq \sum_{k=a}^b
f(k) \leq \alpha (i) f(a) + \beta (i)  f(b) \ee for $a < i < b$. The
minimum gap between the lower and upper bounds is achieved at $i$
such that $\lf j \rf \leq i \leq \lc j \rc$.

(II): If $f(k+1) - f(k) \geq f(k) - f(k-1)$  for $a < k < b$, then
\[
\f{ (b-a+1) [f(a) + f(b)] }{2}  \geq \sum_{k=a}^b f(k) \geq \alpha
(i) f(a) + \beta (i)  f(b)
\]
for $a < i < b$. The minimum gap between the lower and upper bounds
is achieved at $i$ such that $\lf j \rf \leq i \leq \lc j \rc$.

\eeT

See Appendix \ref{App_discrete_sum} for a proof.  For a continuous
random variable with probability density function $f(x)$, we can
apply the following result to compute upper and lower bounds of
$\int_{a}^b f(x) dx$ over subinterval $[a, b]$.

\beT \la{continuous_integration}  Suppose $f(x)$ is differentiable
over interval $[a, b]$.  The following statements hold true:

(I): If $f(x)$ is concave over $[a, b]$, then {\small $\f{[f(a) +
f(b)](b-a) } { 2 } \leq   \int_{a}^b f(x) dx \leq \f{[f(a) +
f(b)](b-a) } { 2 } + \vDe (t)$,} where {\small $\vDe (t) = \li [
f^\prime(a)  - \f{ f(b) - f(a)  } { b - a } \ri ] \f{(t-a)^2}{2} -
\li [ f^\prime(b)  - \f{ f(b) - f(a)  } { b - a } \ri ]
\f{(b-t)^2}{2}$}.

(II): If $f(x)$ is convex over $[a, b]$, then {\small $\f{[f(a) +
f(b)](b-a) } { 2 } - \vDe (t) \leq \int_{a}^b f(x) dx \leq \f{[f(a)
+ f(b)](b-a) } { 2 }$}.

 The minimum of $\vDe (t)$ is achieved at {\small $t = \f{
f(b) - f(a) + a f^\prime(a) - b  f^\prime(b) } { f^\prime(a) -
f^\prime(b) }$}.

\eeT

See Appendix \ref{App_continuous_integration} for a proof.

\subsection{Factorial Evaluation}

In the evaluation of the coverage probability of a sampling scheme,
a frequent routine is the computation of the logarithm of the
factorial of an integer. To reduce computational complexity, we can
develop a table of $\ln(n!)$ and store it in computer for repeated
use. Such a table can be readily made by the recursive relationship
$\ln( (n+1)! ) = \ln(n + 1) + \ln(n!)$.  Modern computers can easily
support a table of $\ln (n!)$ of size in the order of $10^7$ to
$10^8$, which suffices most needs of our computation.  Another
method to calculate $\ln(n!)$ is to use the following double-sized
bounds:
 \[
 \ln ( \sq{2 \pi n} \; n^n) -n + \f{1}{12n} - \f{1}{360 n^3}  < \ln(n!) <
 \ln ( \sq{2 \pi n} \; n^n ) -n + \f{1}{12n} - \f{1}{360 n^3} +  \f{1}{1260 n^5}
 \]
for all $n \geq 1$. A proof for such bounds can be available in
pages 481-482 of \cite{Concrete}.

\sect{Estimation of Binomial Parameters}

Let $X$ be a Bernoulli random variable with distribution $\Pr \{ X =
1 \} = 1 - \Pr \{ X = 0 \} = p \in (0, 1)$.  In this section, we
shall consider the multistage estimation of binomial parameter $p$,
in the general framework proposed in Section \ref{gen_structure},
based on i.i.d. random samples $X_1, X_2, \cd $ of $X$.

To describe our estimation methods, we shall introduce the following
notations, which will be used throughout this section.

 Define $K_\ell =
\sum_{i=1}^{\mbf{n}_\ell} X_i$ and $\wh{\bs{p}}_\ell =
\f{K_\ell}{\mbf{n}_\ell}$
 for $\ell = 1, \cd, s$, where $\mbf{n}_\ell$ is the number of
 samples available at the $\ell$-th stage.
Specially, if the sample sizes are deterministic numbers
 $n_1 < n_2 < \cd < n_s$, then $\mbf{n}_\ell = n_\ell$ for $\ell =
 1, \cd, s$.  As described in Section \ref{gen_structure},
 the stopping rule is that sampling is continued until
$\bs{D}_\ell = 1$ for some $\ell \in \{1, \cd, s\}$, where
$\bs{D}_\ell$ is the decision variable for the $\ell$-th stage. Let
$\bs{\wh{p}} = \f{\sum_{i=1}^{\mathbf{n}} X_i}{\mathbf{n}}$, where
$\mathbf{n}$ is the sample size when the sampling is terminated.
Clearly, $\bs{\wh{p}} = \bs{\wh{p}}_{\bs{l}}$ and $\mbf{n} =
\mbf{n}_{\bs{l}}$, where $\bs{l}$ is the index of stage when the
sampling is terminated.  As mentioned before, the number of stage,
$s$, can be a finite number or infinity.

In the development of our multistage sampling schemes, we need to
use the following probability inequalities related to bounded
variables.

\beL  \la{Hoe_Mas} Let $\ovl{X}_n = \f{\sum_{i=1}^n X_i}{n}$, where
$X_1, \; \cd,\; X_n$ are i.i.d.  random variables such that $0 \leq
X_i \leq 1$ and $\bb{E}[X_i] = \mu \in (0, 1)$ for $i = 1, \; \cd,
n$. Then, \bel \Pr \li \{ \ovl{X}_n \geq z \ri \} & \leq & \exp \li
(n \mscr{M}_{\mrm{B}} \li (z, \mu \ri ) \ri ) \la{Hineq1} \\
& < & \exp \li ( n \mscr{M} (z, \mu) \ri ) \la{Mineq1} \eel for any
$z \in (\mu, 1)$. Similarly, \bel \Pr \li \{ \ovl{X}_n \leq z \ri \}
& \leq & \exp \li (n \mscr{M}_{\mrm{B}} \li (z, \mu \ri ) \ri ) \la{Hineq2} \\
& < & \exp \li ( n \mscr{M} (z, \mu)  \ri ) \la{Mineq2} \eel for any
$z \in (0, \mu)$. \eeL

Inequalities  (\ref{Hineq1}) and (\ref{Hineq2}) are classical
results established by Hoeffding in 1963 (see, \cite{Hoeffding}).
Inequalities (\ref{Mineq1}) and (\ref{Mineq2}) are recent results
due to Massart \cite{Massart:90}.  In this paper, (\ref{Hineq1}) and
(\ref{Hineq2}) are referred to as Hoeffding's inequalities.
Similarly, (\ref{Mineq1}) and (\ref{Mineq2}) are referred to as
Massart's inequalities. If $X_1, \cd, X_n$ are i.i.d. samples of
Bernoulli random variable $X$, then it can be shown that \bee &  &
\exp ( \mscr{M}_{\mrm{B}} (z, \mu) ) = \inf_{t < 0} \bb{E} [e^{ t (X
- z) } ] = \mcal{F} (z, \mu) \qu \tx{for} \; z \leq \mu, \\
&  & \exp ( \mscr{M}_{\mrm{B}} (z, \mu) ) = \inf_{t > 0} \bb{E} [e^{
t (X - z) } ] = \mcal{G} (z, \mu) \qu \tx{for} \; z \geq \mu \eee
which implies that (\ref{Hineq1}) and (\ref{Hineq2}) are actually
Chernoff bounds in the special case.

\subsection{Control of Absolute Error}

In this subsection, we shall propose multistage sampling schemes for
estimating $p$ with an absolute error criterion.  Specifically, for
margin of absolute error $\vep \in (0, \f{1}{2})$, we want to design
a multistage sampling scheme such that the estimator $\wh{\bs{p}}$
satisfies the requirement that $\Pr \{ | \wh{\bs{p}} - p | < \vep
\mid p \} > 1 - \de$ for any $p \in (0, 1)$.

\subsubsection{Stopping Rules from CDF $\&$ CCDF, Chernoff Bounds and Massart's
Inequality} \la{Abs_Ch_ST}

To construct an estimator satisfying an absolute error criterion
with a prescribed confidence level, we propose three types of
multistage sampling schemes with different stopping rules as
follows.

\bed

\item [Stopping Rule (i):] For $\ell = 1, \cd, s$, decision variable $\bs{D}_\ell$  assumes
value $1$ if {\small $F_{\wh{\bs{p}}_\ell}  \li (\wh{\bs{p}}_\ell,
\wh{\bs{p}}_\ell + \vep \ri ) \leq \ze \de, \; G_{\wh{\bs{p}}_\ell}
\li (\wh{\bs{p}}_\ell, \wh{\bs{p}}_\ell - \vep \ri ) \leq \ze \de$};
and assumes value $0$ otherwise.

\item [Stopping Rule (ii):] For $\ell = 1, \cd, s$, decision variable
$\bs{D}_\ell$ assumes value $1$  if $\mscr{M}_{\mrm{B}} ( \f{1}{2} -
|\f{1}{2} - \wh{\bs{p}}_\ell | , \f{1}{2} - |\f{1}{2} -
\wh{\bs{p}}_\ell | + \vep) \leq \f{ \ln ( \ze \de  ) } { n_\ell }$;
and assumes value $0$ otherwise.

\item [Stopping Rule (iii):]  For $\ell = 1, \cd, s$, decision variable
$\bs{D}_\ell$ assumes value $1$  if \be \la{massart836} \li ( \li |
\wh{\bs{p}}_\ell - \f{1}{2} \ri | - \f{2 \vep }{3} \ri )^2 \geq
\f{1}{4} + \f{ \vep^2 n_\ell } {2 \ln (\ze \de) }; \ee and assumes
value $0$ otherwise.

\eed

Stopping rule (i) is derived by virtue of the CDF $\&$ CCDF of
$\wh{\bs{p}}_\ell$.  Stopping rule (ii) is derived by virtue of
Chernoff bounds of the CDF $\&$ CCDF of $\wh{\bs{p}}_\ell$. Stopping
rule (iii) is derived by virtue of Massart's inequality for the CDF
$\&$ CCDF of $\wh{\bs{p}}_\ell$.

For stopping rules (ii) and (iii), we have the following results.

\beT \la{Bino_ABS_CDF_CH_MA}

Suppose that the sample size at the $s$-th stage is no less than
{\small $\li \lc \f{ \ln \f{1}{\ze \de} } { 2 \vep^2 } \ri \rc $}.
Then, {\small \bee & & \Pr \{ p \leq \wh{\bs{p}} - \vep \mid p \}
\leq \sum_{\ell = 1}^s \Pr \{ p \leq \wh{\bs{p}}_\ell - \vep, \;
\bs{D}_\ell = 1 \mid p \}
\leq  s \ze \de, \\
&  & \Pr \{ p \geq \wh{\bs{p}} +  \vep \mid p \} \leq \sum_{\ell =
1}^s \Pr \{ p \geq \wh{\bs{p}}_\ell  + \vep, \; \bs{D}_\ell = 1 \mid
p \} \leq  s \ze \de
 \eee}
and $\Pr \{  | \wh{\bs{p}} - p | < \vep \mid p \} \geq 1 - 2 s \ze
\de$ for any $p \in (0, 1)$.

\eeT

\bsk

See Appendix \ref{App_Bino_ABS_CDF_CH_MA} for a proof.

For stopping rules derived from CDFs $\&$ CCDFs, we can choose the
smallest sample sizes and the largest sample sizes  based on the
criteria proposed in Section \ref{gen_structure} such that {\small
$n_1 \geq \f{ \ln (\ze \de) } { \ln (1- \vep) } $} and  $n_s$ is the
smallest integer which ensures that {\small $F_{\wh{\bs{p}}_s}  \li
(\wh{\bs{p}}_s, \wh{\bs{p}}_s + \vep \ri ) \leq \ze \de, \;
G_{\wh{\bs{p}}_s} \li (\wh{\bs{p}}_s, \wh{\bs{p}}_s - \vep \ri )
\leq \ze \de$} is a sure event.

For stopping rules derived from Chernoff bounds, we can choose the
smallest sample sizes and the largest sample sizes  based on the
criteria proposed in Section \ref{gen_structure} such that {\small
$n_1 \geq \f{ \ln (\ze \de) } { \ln (1- \vep) } $} and {\small $n_s
\geq \f{ \ln \f{1}{\ze \de} } { 2 \vep^2 } $}. Specifically, the
sample sizes $n_1 < n_2 < \cd < n_s$ can be be chosen as the
ascending arrangement of all distinct elements of {\small \be
\la{SSnums}
 \li \{ \li \lc \f{ C_{\tau - \ell} \ln \f{1}{\ze \de} }
{ 2 \vep^2 } \ri \rc : \ell = 1, \cd, \tau \ri \}, \ee } where
$\tau$ is the maximum integer such that $\f{ C_{\tau - 1} \ln
\f{1}{\ze \de} } { 2 \vep^2 } \geq \f{ \ln (\ze \de) } { \ln (1-
\vep) }$, i.e., $C_{\tau - 1} \geq \f{ 2 \vep^2 } { \ln \f{1}{1 -
\vep} }$.  In a similar manner, for stopping rules derived from
Massart's inequality, the sample sizes $n_1 < n_2 < \cd < n_s$ can
be defined as (\ref{SSnums}) with $\tau$ chosen as the maximum
integer such that $\f{ C_{\tau - 1} \ln \f{1}{\ze \de} } { 2 \vep^2
} \geq \li ( \f{24 \vep - 16 \vep^2}{9} \ri ) \f{\ln \f{1}{\ze
\de}}{2 \vep^2}$, i.e., $C_{\tau - 1} \geq \f{24 \vep - 16
\vep^2}{9}$.

For above sampling methods of choosing sample sizes, we have $\Pr \{
| \wh{\bs{p}} - p | < \vep \mid p \} > 1 - \de$ for any $p \in (0,
1)$ if $\ze < \f{1}{ 2 \tau}$, where $\tau$ is independent of $\de$.
Hence, we can determine a value of $\ze$ as large as possible such
that $\Pr \{  | \wh{\bs{p}} - p | < \vep \mid p \} > 1 - \de$ by
virtue of the computational machinery described in Section
\ref{Computing}.

To evaluate the coverage probability associated with the stopping
rule derived from Chernoff bounds with sample sizes defined by
(\ref{SSnums}), we need to express events $\{ \bs{D}_\ell = i \}, \;
i = 0, 1$ in terms of $K_\ell$. This can be accomplished by using
the following results.

\beT \la{Range_Bino_Chernoff} Let $z^*$ be the unique solution of
equation {\small $ \ln \f{(z + \vep) ( 1 - z) }{z (1 - z - \vep) } =
\f{ \vep } { (z + \vep) (1 - z - \vep) }$ } with respect to $z \in
(\f{1}{2} - \vep, \f{1}{2} )$. Let $n_\ell$ be a sample size smaller
than $\f{ \ln (\ze \de) } { \mscr{M}_{\mrm{B}} (z^*, z^* + \vep) }$.
Let $\udl{z}$ be the unique solution of equation $\mscr{M}_{\mrm{B}}
(z, z + \vep) = \f{ \ln ( \ze \de  ) } { n_\ell }$ with respect to
$z \in [0, z^*)$. Let $\ovl{z}$ be the unique solution of equation
$\mscr{M}_{\mrm{B}} (z, z + \vep) = \f{ \ln ( \ze \de  ) } { n_\ell
}$ with respect to $z \in (z^*, 1 - \vep)$.  Then, $\{ \bs{D}_\ell =
0 \} = \{ n_\ell \udl{z} < K_\ell < n_\ell \ovl{z} \} \cup \{ n_\ell
(1 - \ovl{z}) < K_\ell < n_\ell (1 - \udl{z}) \}$.

\eeT

\bsk

See Appendix \ref{App_Range_Bino_Chernoff} for a proof.

\subsubsection{Asymptotic Stopping Rules}

It should be noted that, for a small $\vep$, we can simplify, by
using Taylor's series expansion formula $\ln (1 + x) = x -
\f{x^2}{2} + o (x^2)$, the sampling schemes described in Section
\ref{Abs_Ch_ST} as follows:

(i) The sequence of sample sizes $n_1, \cd, n_s$ is defined as the
ascending arrangement of all distinct elements of {\small $ \li \{
\li \lc \f{ C_{\tau - \ell} \; \ln \f{1}{\ze \de} } { 2 \vep^2 } \ri
\rc : \ell = 1, \cd, \tau \ri \}$}, where $\tau$ is the maximum
integer such that $C_{\tau - 1} \geq 2 \vep$.

(ii) The decision variables are defined such that $\bs{D}_\ell = 1$
if \be \la{Wald89} n_\ell \geq \f{\wh{\bs{p}}_\ell (1 -
\wh{\bs{p}}_\ell) \; 2 \ln \f{1}{\ze \de}  }{\vep^2}; \ee and
$\bs{D}_\ell = 0$ otherwise.

\bsk

For such a simplified sampling scheme, we have \bel \sum_{\ell =
1}^s \Pr \li \{ | \wh{\bs{p}}_\ell - p | \geq \vep, \;  \bs{D}_\ell
= 1 \ri \} & \leq & \sum_{\ell = 1}^s \Pr \li \{ | \wh{\bs{p}}_\ell
- p | \geq \vep \ri \} \leq  \sum_{\ell = 1}^\tau \Pr \li \{ |
\wh{\bs{p}}_\ell - p |
\geq \vep \ri \} \nonumber\\
& \leq &  \sum_{\ell = 1}^\tau 2 e^{- 2 n_\ell \vep^2} \la{ineCH}\\
& < & 2 \tau e^{- 2 n_1 \vep^2} \leq 2 \tau \exp \li ( - 2 \vep \ln
\f{1}{\ze \de} \ri ), \la{LastB} \eel where (\ref{ineCH}) is due to
the Chernoff bound. As can be seen from (\ref{LastB}), the last
bound is independent of $p$ and can be made smaller than $\de$ if
$\ze$ is sufficiently small.  This establishes the claim and it
follows that $\Pr \li \{ \li | \bs{\wh{p}} - p \ri | < \vep \mid p
\ri \}
> 1 - \de$ for any $p \in (0, 1)$ if $\ze$ is sufficiently
small.

By virtue of the normal approximation method as used in the
derivation of stopping rule (\ref{normalgen}), we can simplify
Stopping Rule (i) described in section 4.1.1 as follows: For small
$\vep$, the sample sizes $n_1, \cd, n_s$ are large. Hence, by the
central limit theorem,
\[
F_{\wh{\bs{p}}_\ell} \li (\wh{\bs{p}}_\ell, \wh{\bs{p}}_\ell + \vep
\ri ) \ap \Phi \li (  \f{ - \vep } { \sqrt{ \f{ (\wh{\bs{p}}_\ell +
\vep)  (1 - \wh{\bs{p}}_\ell - \vep )  } { n_\ell }  }   } \ri ),
\qqu G_{\wh{\bs{p}}_\ell} \li (\wh{\bs{p}}_\ell, \wh{\bs{p}}_\ell -
\vep \ri ) \ap 1 - \Phi \li (  \f{  \vep } { \sqrt{ \f{
(\wh{\bs{p}}_\ell - \vep)  (1 - \wh{\bs{p}}_\ell + \vep )  } {
n_\ell }  }   } \ri )
\]
for $\ell = 1, \cd, s$.  Therefore, the stopping condition
$F_{\wh{\bs{p}}_\ell} \li (\wh{\bs{p}}_\ell, \wh{\bs{p}}_\ell + \vep
\ri ) \leq \ze \de, \; G_{\wh{\bs{p}}_\ell} \li (\wh{\bs{p}}_\ell,
\wh{\bs{p}}_\ell - \vep \ri ) \leq \ze \de$ is roughly equivalent to
\[
\Phi \li (  \f{ - \vep } { \sqrt{ \f{ (\wh{\bs{p}}_\ell + \vep)  (1
- \wh{\bs{p}}_\ell - \vep )  } { n_\ell }  }   } \ri ) \leq \ze \de,
\qqu 1 - \Phi \li (  \f{  \vep } { \sqrt{ \f{ (\wh{\bs{p}}_\ell -
\vep)  (1 - \wh{\bs{p}}_\ell + \vep )  } { n_\ell }  }   } \ri )
\leq \ze \de,
\]
which can be written as \be \la{normalrule8} \li ( \li |
\wh{\bs{p}}_\ell - \f{1}{2} \ri | - \vep  \ri )^2 \geq \f{1}{4} -
n_\ell \li ( \f{ \vep } {\mcal{Z}_{\ze \de} } \ri )^2 \ee after some
tedious algebraic manipulations. This implies that Stopping Rule (i)
can be simplified as: For $\ell = 1, \cd, s$, decision variable
$\bs{D}_\ell$ assumes value $1$ if (\ref{normalrule8}) is satisfied;
and assumes value $0$ otherwise.

Since for any $\ze \in (0, \f{1}{\de})$, there exists a unique
number $\ze^\prime \in(0, \f{1}{\de})$ such that $\mcal{Z}_{\ze \de}
= \sqrt{2 \ln \f{1}{\ze^\prime \de}}$, the above simplified stopping
rule is equivalent to the following stopping rule:  For $\ell = 1,
\cd, s$, decision variable $\bs{D}_\ell$ assumes value $1$ if \be
\la{normalmassart} \li ( \li | \wh{\bs{p}}_\ell - \f{1}{2} \ri | -
\vep  \ri )^2 \geq \f{1}{4} +  \f{ \vep^2 n_\ell } { 2 \ln (\ze \de)
} \ee is satisfied; and assumes value $0$ otherwise.

Comparing (\ref{Wald89}), (\ref{normalmassart}) and
(\ref{massart836}), we can see that the stopping conditions can be
put in a general form \be \la{simplegeneral} \li ( \li |
\wh{\bs{p}}_\ell - \f{1}{2} \ri | - w \vep \ri )^2 \geq \f{1}{4} +
\f{ \vep^2 n_\ell } { 2 \ln (\ze \de) }, \ee where $w \geq 0$ is a
parameter affecting the shape of the stopping boundary.  Taking $w
=0, \;  \f{2}{3}$ and $1$ leads to (\ref{Wald89}),
(\ref{normalmassart}) and (\ref{massart836}) respectively.
Therefore, a very general and simple stopping rule can be stated as
follows:

For $\ell = 1, \cd, s$, decision variable $\bs{D}_\ell$ assumes
value $1$ if (\ref{simplegeneral}) is satisfied; and assumes value
$0$ otherwise.

For multistage sampling schemes with such a stopping rule, we have
established the following general results:

\beT \la{them18889}
Suppose that  $0 < w \leq 1$ and $n_s \geq \f{ \ln \f{1}{\ze \de} }{2 \vep^2}$. Then, $\Pr \li \{ \li | \bs{\wh{p}} - p \ri |
< \vep \mid p \ri \}$ is greater than $1 - \de$ for any $p \in (0, 1)$ if $\ze > 0$ is sufficiently small.  Moreover, \be \la{good9} \Pr \li \{
\li | \bs{\wh{p}} - p \ri | < \vep \mid p \ri \} \geq 1 - 2 s (\ze \de)^{9 w^2 \sh 4} \ee for any $p \in (0, 1)$ provided that $0 < w \leq
\f{2}{3}$ and $n_s \geq \f{ \ln \f{1}{\ze \de} }{2 \vep^2}$. \eeT

See Appendix \ref{them18889_app} for a proof.

Under the restriction that $0 \leq w \vep \leq \f{1}{4}$, the sample
sizes $n_1, \cd, n_s$ for the above stopping rule can be chosen
based on the following analysis: As a consequence of $0 \leq w \vep
\leq \f{1}{4}$ and $\li | \wh{\bs{p}}_\ell - \f{1}{2} \ri | \leq
\f{1}{2}$, it must be true that $\li ( \li | \wh{\bs{p}}_\ell -
\f{1}{2} \ri | - w \vep \ri )^2 \leq \li ( \f{1}{2} - w \vep \ri
)^2$.  Thus, (\ref{simplegeneral}) will not be satisfied if $\li (
\f{1}{2} - w \vep \ri )^2 < \f{1}{4} + \f{ \vep^2 n_\ell } { 2 \ln
(\ze \de) }$, or equivalently, $n_\ell < \f{2 w ( 1 - w \vep) \ln
\f{1}{\ze \de} }{ \vep }$.  This implies that the minimum sample
size $n_1$, i.e.,  the sample size for the first stage, should be
chosen to be no less than $\f{2 w ( 1 - w \vep) \ln \f{1}{\ze \de}
}{ \vep }$.  To determine the maximum sample size $n_s$, i.e., the
sample size for the last stage, observe that $\li ( \li |
\wh{\bs{p}}_\ell - \f{1}{2} \ri | - w \vep \ri )^2 \geq 0$ and thus
(\ref{simplegeneral}) will always be satisfied if $\f{1}{4} + \f{
\vep^2 n_\ell } { 2 \ln (\ze \de) } \leq 0$, or equivalently,
$n_\ell \geq \f{ \ln \f{1}{\ze \de} }{ 2 \vep^2}$. Therefore, the
maximum sample size $n_s$ should be chosen to be the smallest
integer no less than $\f{ \ln \f{1}{\ze \de} }{ 2 \vep^2}$. For a
fixed value of $w$, the appropriate value of $\ze$ can be obtained
by bisection coverage tuning. By optimizing the performance of the
stopping rules over $w$, it is possible to obtain sampling schemes
better than those associated with special values $w = \f{2}{3}$ or
$1$.  Of course, the improvement of performance should be available
with the price of more computational effort, since an extra boundary
parameter $w$ is introduced.

Before concluding this subsection, we also want to point out that it
is possible to modify (\ref{simplegeneral}) to obtain the following
stopping rule: For $\ell = 1, \cd, s$, decision variable
$\bs{D}_\ell$ assumes value $1$ if
\[ \li ( \wh{\bs{p}}_\ell - \f{1}{2} \ri )^2 \geq \f{1}{4} +
\f{w}{n_\ell} + \f{ \vep^2 n_\ell } { 2 \ln (\ze \de) } \] is
satisfied; and assumes value $0$ otherwise, where $w \geq 0$ is a
parameter affecting the shape of the stopping boundary.  As
suggested in Section 2.1, the maximum sample size $n_s$, i.e., the
sample size of the last stage,  should be defined as the smallest
integer such that $\{ \bs{D}_s = 1 \}$ is a sure event. The minimum
sample size $n_1$, i.e., the sample size of the first stage, should
be defined as the smallest integer such that $\{ \bs{D}_1 = 1 \}$ is
an event of a positive probability.  We can show that the coverage
probability $\Pr \li \{ \li | \bs{\wh{p}} - p \ri | < \vep \mid p
\ri \}$ is greater than $1 - \de$ for any $p \in (0, 1)$ if $\ze$ is
sufficiently small.

\subsubsection{Asymptotic Analysis of Sampling Schemes}

In this subsection, we shall focus on the asymptotic analysis of
multistage sampling schemes.  Throughout this subsection, we assume
that the multistage sampling schemes follow stopping rules derived
from Chernoff bounds as described in Section \ref{Abs_Ch_ST}.
Moreover, we assume that the sample sizes $n_1, \cd, n_s$ are chosen
as the ascending arrangement of all distinct elements of  the set
defined by (\ref{SSnums}).

With regard to the tightness of the DDV bound, we have

\beT \la{Bino_DDV_Asp} Let $\mscr{R}$ be a subset of real numbers.
Define {\small \[ \ovl{P} = \sum_{\ell = 1}^s \Pr \{
\wh{\bs{p}}_\ell \in \mscr{R}, \; \bs{D}_{\ell - 1} = 0, \;
\bs{D}_\ell = 1 \}, \qqu \udl{P} = 1 - \sum_{\ell = 1}^s \Pr \{
\wh{\bs{p}}_\ell \notin \mscr{R}, \; \bs{D}_{\ell - 1} = 0, \;
\bs{D}_\ell = 1 \}.
\]}
Then, $\udl{P} \leq \Pr \{ \wh{\bs{p}} \in \mscr{R} \} \leq \ovl{P}$
and $\lim_{\vep \to 0} | \Pr \{ \wh{\bs{p}} \in \mscr{R} \} -
\ovl{P} | = \lim_{\vep \to 0}  | \Pr \{ \wh{\bs{p}} \in \mscr{R} \}
- \udl{P}
 | = 0$ for any $p \in (0, 1)$.  \eeT

\bsk

See Appendix \ref{App_Bino_DDV_Asp} for a proof.

For $\ro > 0, \; d > 0, \; 0 < \nu \leq 1$, define
\[
\Psi (\ro, \nu, d) = \f{1}{2 \pi} \li [ \int_{- \phi_L}^{\phi_U} \exp \li (  - \f{ \nu^2 d^2  } { 2 \cos^2 \phi } \ri ) d \phi  + \int_{\phi_U -
\phi_\ro}^{2 \pi - \phi_L - \phi_\ro} \exp \li (  - \f{(1 + \ro) d^2  } { 2 \cos^2 \phi } \ri ) d \phi  \ri ]
\]
with {\small $\phi_\ro = \arctan (\sq{\ro}), \; \phi_L = \arctan \li ( \f{1 + \ro}{\nu \sq{\ro} }  +  \f{1}{\sq{\ro}} \ri )$} and {\small
$\phi_U = \arctan \li ( \f{1 + \ro}{\nu \sq{\ro} } - \f{1}{\sq{\ro}} \ri )$}.   With regard to the asymptotic performance of the sampling
scheme, we have

\beT \la{Bino_Asp_Analysis} Let {\small $\mcal{N}_{\mrm{a}} (p,
\vep) = \f{ \ln ( \ze \de ) } { \mscr{M}_{\mrm{B}} ( \f{1}{2} -
|\f{1}{2} - p | , \f{1}{2} - |\f{1}{2} - p | + \vep) }$}. Let
$\mcal{N}_{\mrm{f}} (p, \vep)$ be the minimum sample number $n$ such
that {\small $\Pr  \{  | \f{\sum_{i = 1}^n X_i}{n} - p  | < \vep
\mid p \} > 1 - \ze \de$} for a fixed-size sampling procedure. Let
$j_p$ be the maximum integer $j$ such that $C_j \geq 4 p (1 - p)$.
Let {\small $\nu = \f{2}{3}, \; d = \sq{2 \ln \f{1}{\ze \de}}$} and
{\small $\ka_p = \f{C_{j_p}}{4 p ( 1 - p)} $}.  Let {\small $\ro_p =
\f{C_{j_p - 1}} { 4 p ( 1 - p) } - 1$} for $\ka_p = 1, \; j_p > 0$
and $\ro_p = \ka_p - 1$ otherwise. The following statements hold
true:

(I): {\small $\Pr \li \{  1 \leq \limsup_{\vep \to 0} \f{ \mbf{n} }
{ \mcal{N}_{\mrm{a}}  (p, \vep) } \leq 1 + \ro_p \ri \} = 1$}.
Specially, {\small $\Pr \li \{  \lim_{\vep \to 0} \f{ \mbf{n} } {
\mcal{N}_{\mrm{a}} (p, \vep) } = \ka_p \ri \} = 1$} if $\ka_p > 1$.

(II): $\lim_{\vep \to 0} \f{ \bb{E} [ \mbf{n} ] } {
\mcal{N}_{\mrm{f}} (p, \vep)} = \li ( \f{ d } { \mcal{Z}_{\ze \de} }
\ri )^2 \times \lim_{\vep \to 0} \f{ \bb{E} [ \mbf{n} ] } {
\mcal{N}_{\mrm{a}} (p, \vep) }$, where
\[ \lim_{\vep \to 0} \f{ \bb{E} [ \mbf{n} ] }  {
\mcal{N}_{\mrm{a}} (p, \vep) } = \bec \ka_p
& \tx{if} \; \ka_p > 1,\\
1  + \ro_p \Phi (\nu d ) &  \tx{otherwise} \eec
\]
and $1 \leq \lim_{\vep \to 0} \f{ \bb{E} [ \mbf{n} ] }  {
\mcal{N}_{\mrm{a}} (p, \vep) } \leq 1 + \ro_p$.

(III): If $\ka_p > 1$, then $\lim_{\vep \to 0} \Pr \{ | \wh{\bs{p}} - p | < \vep \} = 2 \Phi \li ( d \sq{\ka_p} \ri ) - 1 > 2 \Phi \li ( d \ri )
- 1 > 1 - 2 \ze \de$. Otherwise, $\Phi \li ( d  \ri )  + \Phi \li ( d \sq{1 + \ro_p}  \ri )  - 1
> \lim_{\vep \to 0} \Pr \{ | \wh{\bs{p}} - p | < \vep  \} = 1 +
\Phi(d) - \Phi(\nu d) - \Psi (\ro_p, \nu, d) >  \Phi \li ( d \ri ) + 2 \Phi \li ( d \sq{1 + \ro_p}  \ri )  - 2 > 1 - 3 \ze \de$.

 \eeT

See Appendix \ref{App_Bino_Asp_Analysis} for a proof.

\subsection{Control of Relative Error}

In this section, we shall focus on the design of multistage sampling
schemes for estimating the binomial parameter $p$ with a relative
error criterion.  Specifically, for $\vep \in (0, 1)$, we wish to
construct a multistage sampling scheme and its associated estimator
$\wh{\bs{p}}$ for $p$ such that $\Pr \{ | \wh{\bs{p}} - p | < \vep p
\mid p \} > 1 - \de$ for any $p \in (0, 1)$.

\subsubsection{Multistage Inverse Sampling} \la{multiinv}

In this subsection, we shall develop multistage sampling schemes, of
which the number of stages, $s$, is a finite number.   Let $\ga_1 <
\ga_2 < \cd < \ga_s$ be a sequence of positive integers. The number,
$\ga_\ell$, is referred to as the {\it threshold of sample sum} of
the $\ell$-th stage. For $\ell = 1, \cd, s$, let $\wh{\bs{p}}_\ell =
\f{\ga_\ell}{ \mathbf{n}_\ell }$, where $\mathbf{n}_\ell$ is the
minimum number of samples such that $\sum_{i = 1}^{\mathbf{n}_\ell}
X_i = \ga_\ell$. As described in Section \ref{gen_structure}, the
stopping rule is that sampling is continued until $\bs{D}_\ell = 1$
for some $\ell \in \{1, \cd, s\}$, where $\bs{D}_\ell$ is the
decision variable for the $\ell$-th stage.  Define estimator
$\wh{\bs{p}} = \f{\sum_{i=1}^{\mathbf{n}} X_i}{\mathbf{n}}$,  where
$\mathbf{n}$ is the sample size when the sampling is terminated.

The rationale for choosing $\wh{\bs{p}}$ as an estimator for $p$ can
be illustrated by the following results.

\beT

\la{Unbiased_Bino}  Suppose that $\ga_{\ell + 1} - \ga_\ell \geq 1$
for any $\ell > 0$. Then $\bb{E} [\wh{\bs{p}} - p]$ and $\bb{E} |
\wh{\bs{p}} - p |^k, \; k = 1, 2, \cd$ tend to $0$ as the minimum
threshold of sample sum tends to infinity.

\eeT

See Appendix \ref{App_Unbiased_Bino} for a proof.

It should be noted that there exists an inherent connection between
the multistage inverse sampling scheme for Bernoulli random variable
$X$ and a multistage sampling scheme of sample sizes $\ga_1 < \ga_2
< \cd < \ga_s$ for a random variable $Y$ possessing a geometric
distribution with parameter $\se = \f{1}{p} = \bb{E} [ Y ]$. To see
this, for $j = 1, \cd, \ga_s$, let $Y_j$ be a random variable such
that $\sum_{i = 1}^{Y_j} X_i = j
> \sum_{i = 1}^{Y_j - 1} X_i$.  Then, \[ \mbf{n}_\ell = \sum_{i =
1}^{\ga_\ell} Y_i, \qqu \ell = 1, \cd, s
\]
and $Y_i, \; i = 1, \cd, \ga_s$ are i.i.d. samples of the geometric
random variable $Y$.  Clearly, for $\ell = 1, \cd, s$,
$\wh{\bs{\se}}_\ell = \f{\sum_{i = 1}^{\ga_\ell} Y_i }{\ga_\ell} =
\f{1}{\wh{\bs{p}}_\ell}$ is a ULE for $\se$.  Let $\bs{l}$ be the
index stage at the termination of the multistage inverse sampling
process as before.  Then, a ULE for $\se$ can be defined as
$\wh{\bs{\se}} = \wh{\bs{\se}}_{\bs{l}} =
\f{1}{\wh{\bs{p}}_{\bs{l}}} = \f{1}{\wh{\bs{p}}} =
\f{\mathbf{n}}{\sum_{i=1}^{\mathbf{n}} X_i}$. It follows that the
problem of constructing a ULE $\wh{\bs{p}} =
\f{\sum_{i=1}^{\mathbf{n}} X_i}{\mathbf{n}}$ for $p$ to ensure $\Pr
\{  | \wh{\bs{p}} - p | < \vep p \mid p \} > 1 - \de$,  where $\vep,
\de \in (0, 1)$, is equivalent to the problem of constructing a ULE
$\wh{\bs{\se}} = \f{\mathbf{n}}{\sum_{i=1}^{\mathbf{n}} X_i}$ for
$\se = \f{1}{p}$ such that $\Pr  \{ (1 - \vep)  \wh{\bs{\se}}  < \se
< (1 + \vep) \wh{\bs{\se}} \mid \se \} > 1 - \de$.  Thus, the
general stopping rule proposed in Section 2.5 can be applied to
construct a random interval $( (1 - \vep)  \wh{\bs{\se}},  (1 +
\vep) \wh{\bs{\se}} )$ for $\se$, or equivalently, a random interval
$\li ( \f{ \wh{\bs{p}} }{1 + \vep},   \f{ \wh{\bs{p}} }{1 - \vep}
\ri )$ for $p$ to guarantee that the coverage probability is greater
than $1 - \de$.  In this direction, we can use CDF $\&$ CCDF
functions of $\wh{\bs{\se}}_\ell$ or $\wh{\bs{p}}_\ell$,  their
approximations and bounds to design stopping rules.

By virtue of the CDF $\&$ CCDF of $\wh{\bs{p}}_\ell$, we propose a
class of multistage sampling schemes as follows.

\beT \la{Bino_Rev_CDF}  Suppose that, for $\ell = 1, \cd, s$,
decision variable $\bs{D}_\ell$ assumes values $1$ if {\small
$F_{\wh{\bs{p}}_\ell} (\wh{\bs{p}}_\ell, \f{\wh{\bs{p}}_\ell}{1 -
\vep} ) \leq \ze \de, \; G_{\wh{\bs{p}}_\ell} (\wh{\bs{p}}_\ell,
\f{\wh{\bs{p}}_\ell}{1 + \vep} ) \leq \ze \de$}; and assumes $0$
otherwise. Suppose that the threshold of sample sum for the $s$-th
stage is equal to {\small $\li \lc \f{ (1 + \vep) \ln (\ze \de) } {
\vep - (1 + \vep) \ln (1 + \vep) } \ri \rc$}.   Then,  {\bel & & \Pr
\li \{ p \geq \f{ \wh{\bs{p}} }{1 - \vep} \mid  p \ri \} \leq
\sum_{\ell = 1}^s \Pr \{ \wh{\bs{p}}_\ell \leq (1 - \vep) p, \;
\bs{D}_\ell = 1 \mid p \} \leq s \ze \de, \la{rev1cdf}\\
&   & \Pr \li \{ p \leq \f{ \wh{\bs{p}} }{1 + \vep} \mid  p \ri \}
\leq \sum_{\ell = 1}^s \Pr \{ \wh{\bs{p}}_\ell \geq (1 + \vep) p, \;
\bs{D}_\ell = 1 \mid p \} \leq s \ze \de  \la{rev2cdf} \eel} for any
$p \in (0, 1)$.  Moreover, {\small $\Pr \li \{ \li | \f{ \wh{\bs{p}}
- p } { p } \ri | \leq \vep \mid p \ri \} \geq 1 - \de$} for any $p
\in (0, 1)$ provided that $\ze$ is sufficiently small to guarantee
$1 - S_{\mrm{P}} (  \ga_s - 1, \f{\ga_s}{ 1 + \vep}  ) + S_{\mrm{P}}
(  \ga_s - 1, \f{\ga_s}{ 1 - \vep}  )  < \de$ and
 {\small \bee &  & \ln ( \ze \de )
< \li [ \f{ \li ( 1 + \vep  +  \sq{ 1 + 4 \vep + \vep^2 }  \ri )^2}
{ 4 \vep^2 } + \f{1}{2} \ri ] \li [ \f{\vep}{ 1 + \vep}  - \ln (1 +
\vep) \ri ], \\
&  & \Pr \li \{ \li | \f{ \wh{\bs{p}} - p } { p }  \ri | \leq \vep
\mid p \ri \} \geq 1 - \de \eee} for any $p \in [p^*, 1)$,  where
$p^* \in (0, z_{s-1})$ denotes the unique number satisfying
\[
1 - S_{\mrm{P}} \li (  \ga_s - 1, \f{\ga_s} {1 + \vep} \ri ) +
S_{\mrm{P}} \li (  \ga_s - 1, \f{\ga_s} {1 - \vep} \ri ) +
\sum_{\ell = 1}^{s - 1} \exp ( \ga_\ell \mscr{M}_{\mrm{I}} ( z_\ell,
p^* ) ) = \de \] with $z_\ell = \min \{ z \in I_{\wh{\bs{p}}_\ell}:
F_{ \wh{\bs{p}}_\ell } (z, \f{z}{1 - \vep} ) > \ze \de \; \tx{or} \;
G_{ \wh{\bs{p}}_\ell } ( z, \f{z}{1 + \vep}  )
> \ze \de \}$, where $I_{\wh{\bs{p}}_\ell}$ represents the support
of $\wh{\bs{p}}_\ell$, for $\ell = 1, \cd, s$.

\eeT

See Appendix \ref{App_Bino_Rev_CDF} for a proof. Based on the
criteria proposed in Section \ref{gen_structure}, the thresholds of
sample sum $\ga_1 < \ga_2 < \cd < \ga_s$ can be chosen as the
ascending arrangement of all distinct elements of {\small \be
\la{THolds} \li \{ \li \lc  \f{ C_{\tau - \ell} \; (1 + \vep) \ln (
\ze \de )} { \vep - (1 + \vep) \ln (1 + \vep) }  \ri \rc: \ell = 1,
\cd, \tau \ri \}, \ee} where $\tau$ is the maximum integer such that
$\f{ C_{\tau - 1} \; (1 + \vep) \ln ( \ze \de )} { \vep - (1 + \vep)
\ln (1 + \vep) }  \geq \f{ \ln \f{1} {\ze \de}} { \ln (1 + \vep) }$,
i.e., $C_{\tau - 1} \geq 1 - \f{\vep} { (1 + \vep) \ln (1 + \vep)
}$.

By virtue of Chernoff bounds of the CDF $\&$ CCDF of
$\wh{\bs{p}}_\ell$, we propose a class of multistage sampling
schemes as follows.

\beT \la{Bino_Rev_Chernoff}  Suppose that, for $\ell = 1, \cd, s$,
decision variable $\bs{D}_\ell$ assumes values $1$ if {\small
$\mscr{M}_{\mrm{I}} (\wh{\bs{p}}_\ell, \f{\wh{\bs{p}}_\ell}{1 +
\vep}  )  \leq \f{ \ln ( \ze \de  ) } { \ga_\ell }$}; and assumes
$0$ otherwise.  Suppose that the threshold of sample sum for the
$s$-th stage is equal to {\small $\li \lc \f{ (1 + \vep) \ln ( \ze
\de ) } { \vep  - (1 + \vep) \ln (1 + \vep) } \ri \rc$}. Then, {\bel
& & \Pr \li \{ p \geq \f{ \wh{\bs{p}} }{1 - \vep} \mid p \ri \} \leq
\sum_{\ell = 1}^s \Pr \{ \wh{\bs{p}}_\ell \leq (1 - \vep) p, \;
\bs{D}_\ell = 1 \mid p \} \leq s \ze \de, \la{rev1ch}\\
&   & \Pr \li \{ p \leq \f{ \wh{\bs{p}} }{1 + \vep} \mid  p \ri \}
\leq \sum_{\ell = 1}^s \Pr \{ \wh{\bs{p}}_\ell \geq (1 + \vep) p, \;
\bs{D}_\ell = 1 \mid p \} \leq s \ze \de  \la{rev2ch} \eel} for any
$p \in (0, 1)$.  Moreover, {\small $\Pr \li \{ \li | \f{ \wh{\bs{p}}
- p } { p } \ri | \leq \vep \mid p \ri \} \geq 1 - \de$} for any $p
\in (0, 1)$ provided that $\ze$ is sufficiently small to guarantee
$1 - S_{\mrm{P}} (  \ga_s - 1, \f{\ga_s}{ 1 + \vep}  ) + S_{\mrm{P}}
(  \ga_s - 1, \f{\ga_s}{ 1 - \vep}  )  < \de$ and
 {\small \bel &  & \ln ( \ze \de )
< \li [ \f{ \li ( 1 + \vep  +  \sq{ 1 + 4 \vep + \vep^2 }  \ri )^2}
{ 4 \vep^2 } + \f{1}{2} \ri ] \li [ \f{\vep}{ 1 + \vep}  - \ln (1 +
\vep) \ri ], \la{cona}\\
&  & \Pr \li \{ \li | \f{ \wh{\bs{p}} - p } { p }  \ri | \leq \vep
\mid p \ri \} \geq 1 - \de \nonumber \eel} for any $p \in [p^*, 1)$,
where $p^* \in (0, z_{s-1})$ denotes the unique number satisfying
\[
1 - S_{\mrm{P}} \li (  \ga_s - 1, \f{\ga_s} {1 + \vep} \ri ) +
S_{\mrm{P}} \li (  \ga_s - 1, \f{\ga_s} {1 - \vep} \ri )  +
\sum_{\ell = 1}^{s - 1} \exp ( \ga_\ell \mscr{M}_{\mrm{I}} ( z_\ell,
p^* ) ) = \de \] where $z_\ell \in (0, 1)$ is the unique number such
that {\small $\mscr{M}_{\mrm{I}} \li ( z_\ell, \f{ z_\ell } { 1 +
\vep } \ri ) = \f{ \ln ( \ze \de  ) } { \ga_\ell }$} for $\ell = 1,
\cd, s - 1$.
 \eeT

 \bsk

See Appendix \ref{App_Bino_Rev_Chernoff} for a proof.  Based on the
criteria proposed in Section \ref{gen_structure}, the thresholds of
sample sum $\ga_1 < \ga_2 < \cd < \ga_s$ can be chosen as the
ascending arrangement of all distinct elements of the set defined by
(\ref{THolds}).

It should be noted that both $z_\ell$ and $p^*$ can be readily
 computed by a bisection search method due to the monotonicity of
 the function $\mscr{M}_{\mrm{I}} (., .)$.

\bsk

By virtue of Massart's inequality for the CDF $\&$ CCDF of
$\wh{\bs{p}}_\ell$, we propose a class of multistage sampling
schemes as follows.

\beT \la{Bino_Rev_Inverse_Massart} Suppose that, for $\ell = 1, \cd,
s$, decision variable $\bs{D}_\ell$ assumes values $1$ if {\small
$\wh{\bs{p}}_\ell \geq 1 + \f{ 2 \vep}{3 + \vep} + \f{ 9 \vep^2
\ga_\ell } { 2 (3 + \vep)^2 \ln (\ze \de)}$}; and assumes $0$
otherwise. Suppose the threshold of sample sum for the $s$-th stage
is equal to {\small $\li \lc  \f{ 2 (1 + \vep) (3 + \vep) } {3
\vep^2 } \ln \f{1} { \ze \de } \ri \rc$}.  Then, {\bee  & & \Pr \li
\{ p \geq \f{ \wh{\bs{p}} }{1 - \vep} \mid  p \ri \} \leq \sum_{\ell
= 1}^s \Pr \{ \wh{\bs{p}}_\ell \leq (1 - \vep) p, \;
\bs{D}_\ell = 1 \mid p \} \leq s \ze \de, \\
&   & \Pr \li \{ p \leq \f{ \wh{\bs{p}} }{1 + \vep} \mid  p \ri \}
\leq \sum_{\ell = 1}^s \Pr \{ \wh{\bs{p}}_\ell \geq (1 + \vep) p, \;
\bs{D}_\ell = 1 \mid p \} \leq s \ze \de \eee} for any $p \in (0,
1)$.  Moreover, {\small $\Pr \li \{ \li | \f{ \wh{\bs{p}} - p } { p
} \ri | \leq \vep \mid p \ri \} \geq 1 - \de$} for any $p \in (0,
1)$ provided that $\ze$ is sufficiently small to guarantee $1 -
S_{\mrm{P}} (  \ga_s - 1, \f{\ga_s}{ 1 + \vep}  ) + S_{\mrm{P}} (
\ga_s - 1, \f{\ga_s}{ 1 - \vep}  )  < \de$ and
 {\small \bee &  & \ln ( \ze \de )
< \li [ \f{ \li ( 1 + \vep  +  \sq{ 1 + 4 \vep + \vep^2 }  \ri )^2}
{ 4 \vep^2 } + \f{1}{2} \ri ] \li [ \f{\vep}{ 1 + \vep}  - \ln (1 +
\vep) \ri ], \\
&  & \Pr \li \{ \li | \f{ \wh{\bs{p}} - p } { p }  \ri | \leq \vep
\mid p \ri \} \geq 1 - \de \eee} for any $p \in [p^*, 1)$,  where
$p^* \in (0, z_{s-1})$ denotes the unique number satisfying \[ 1 -
S_{\mrm{P}} \li (  \ga_s - 1, \f{\ga_s} {1 + \vep} \ri ) +
S_{\mrm{P}} \li (  \ga_s - 1, \f{\ga_s} {1 - \vep} \ri )  +
\sum_{\ell = 1}^{s - 1} \exp \li ( \f{\ga_\ell}{z_\ell} \mscr{M} (
z_\ell, p^* )  \ri ) = \de
\] with {\small $z_\ell =  1 + \f{ 2 \vep}{3 + \vep} + \f{ 9 \vep^2
\ga_\ell } { 2 (3 + \vep)^2 \ln (\ze \de)}$} for $\ell = 1, \cd, s -
1$.
 \eeT

 \bsk

See Appendix \ref{App_Bino_Rev_Inverse_Massart} for a proof.   Based
on the criteria proposed in Section \ref{gen_structure}, the
thresholds of sample sum $\ga_1 < \ga_2 < \cd < \ga_s$ can be chosen
as the ascending arrangement of all distinct elements of {\small
\[
\li \{ \li \lc 2 C_{\tau - \ell}   \li ( \f{1}{\vep}  +  1 \ri ) \li
( \f{1}{\vep} + \f{1}{3} \ri ) \ln \f{1} { \ze \de } \ri \rc: \ell =
1, \cd, \tau \ri \},
\]} where
$\tau$ is the maximum integer such that $2 C_{\tau - \ell}   \li (
\f{1}{\vep}  +  1 \ri ) \li ( \f{1}{\vep} + \f{1}{3} \ri ) \ln \f{1}
{ \ze \de } \geq \f{ 4(3 + \vep) } {9 \vep } \ln \f{1} { \ze \de }$,
i.e., $C_{\tau - 1} \geq  \f{2 \vep}{ 3 ( 1 + \vep ) } $.

It should be noted that $\{ \bs{D}_\ell = i \}$ can be expressed
 in terms of $\mathbf{n}_\ell$.  Specially,
 we have $\bs{D}_0 = 0, \; \bs{D}_s = 1$ and {\small $\{ \bs{D}_\ell = 0 \}
= \{ \mathbf{n}_\ell > \f{\ga_\ell} {z_\ell} \}$} for $\ell = 1,
\cd, s - 1$.

To apply the truncation techniques of \cite{Chen1} to reduce
computation,  we can make use of the bounds in Lemma \ref{lem_chen}
and a bisection search to truncate the domains of $\mathbf{n}_{\ell
- 1}$ and $\mathbf{n}_\ell$ to much smaller sets.   Since
$\mathbf{n}_\ell - \mathbf{n}_{\ell - 1}$ can be viewed as the
number of binomial trials to come up with $\ga_\ell - \ga_{\ell -
1}$ occurrences of successes, we have that $\mathbf{n}_\ell -
\mathbf{n}_{\ell - 1}$ is independent of $\mathbf{n}_{\ell - 1}$.
Hence, the technique of triangular partition described in Section
\ref{Tri_par} can be used by identifying $\mbf{n}_{\ell -1}$ as $U$
and $\mbf{n}_{\ell} - \mbf{n}_{\ell - 1}$ as $V$ respectively.  The
computation can be reduced to computing the following types of
probabilities: \bee &   & \Pr \{ u \leq \mathbf{n}_{\ell - 1} \leq v
\mid p \} = \sum_{n = u}^v \bi{ n - 1 }{\ga_{\ell - 1} - 1} \li (
\f{p}{1 - p} \ri )^{\ga_{\ell - 1}} (1 - p)^n, \\
&   &\Pr \{ u \leq \mathbf{n}_\ell - \mathbf{n}_{\ell - 1} \leq v
\mid p \} = \sum_{n = u }^v \bi{ n - 1 }{\ga_\ell - \ga_{\ell - 1} -
1} \li ( \f{p}{1 - p} \ri )^{\ga_\ell - \ga_{\ell - 1}} (1 - p)^n
\eee where $u$ and $v$ are integers.

\bsk

From the definition of the sampling scheme, it can be seen that the
probabilities  that $\wh{\bs{p}}$ is greater or smaller than certain
values can be expressed in terms of probabilities of the form $\Pr
\{ \mbf{n}_i \in \bb{N}_i, \; i = 1, \cd, \ell \}, \; 1 \leq \ell
\leq s$,  where $\bb{N}_1, \cd, \bb{N}_s$ are subsets of natural
numbers. Such probabilities can be computed by using the recursive
relationship {\small \bee & & \Pr \{ \mbf{n}_i \in \bb{N}_i, \; i =
1, \cd, \ell; \; \mbf{n}_{\ell +
1} = n_{\ell + 1} \}\\
& = & \sum_{n_\ell \in \bb{N}_\ell} \Pr \{ \mbf{n}_i \in \bb{N}_i,
\; i = 1, \cd, \ell - 1; \; \mbf{n}_\ell = n_\ell \} \Pr \{
\mbf{n}_{\ell + 1} - \mbf{n}_\ell = n_{\ell + 1} - n_\ell \}\\
& = & \sum_{n_\ell \in \bb{N}_\ell} \Pr \{ \mbf{n}_i \in \bb{N}_i,
\; i = 1, \cd, \ell - 1; \; \mbf{n}_\ell = n_\ell \}  \times \bi{
n_{\ell + 1} - n_\ell - 1 }{\ga_\ell - \ga_{\ell - 1} - 1} \li (
\f{p}{1 - p} \ri )^{\ga_\ell - \ga_{\ell - 1}} (1 - p)^{n_{\ell + 1}
- n_\ell} \eee} for $\ell = 1, \cd, s - 1$.

\bsk

With regard to the average sample number, we have

\beT \la{ASN_Bino_Inverse} For any $p \in (0, 1]$,  $\bb{E}[
\mathbf{n} ] = \f{ \bb{E}[ \bs{\ga} ] } {p}$ with $\bb{E}[ \bs{\ga}
] = \ga_1 + \sum_{\ell = 1}^{s - 1} (\ga_{\ell + 1} - \ga_\ell ) \Pr
\{\bs{l}
> \ell \}$.  \eeT

See Appendix \ref{App_ASN_Bino_Inverse} for a proof.

\subsubsection{Asymptotic Stopping Rule}  \la{Bino_inv_asp}

 We would like to remark that, for a small $\vep$, we can simplify,
by using Taylor's series expansion formula $\ln (1 + x) = x -
\f{x^2}{2} + o (x^2)$, the multistage inverse sampling schemes
described in Section \ref{multiinv} as follows:

(i) The sequence of thresholds $\ga_1, \cd, \ga_s$ is defined as the
ascending arrangement of all distinct elements of {\small $ \li \{
\li \lc  \f{ 2 C_{\tau - \ell} \; \ln \f{1}{\ze \de} } { \vep^2 }
\ri \rc : \ell = 1, \cd, \tau \ri \}$}, where $\tau$ is the maximum
integer such that $C_{\tau - 1} \geq \f{\vep}{2}$.

(ii) The decision variables are defined such that $\bs{D}_\ell = 1$
if \be \la{invrev} \ga_\ell \geq \f{ (1 - \wh{\bs{p}}_\ell) \; 2 \ln
\f{1}{\ze \de} }{\vep^2}; \ee
 and $\bs{D}_\ell = 0$ otherwise.

\bsk

For such a simplified sampling scheme, we have \bel \sum_{\ell =
1}^s \Pr \li \{ | \wh{\bs{p}}_\ell - p | \geq \vep p, \; \bs{D}_\ell
= 1 \ri \} & \leq & \sum_{\ell = 1}^s \Pr \li \{ | \wh{\bs{p}}_\ell
- p | \geq \vep p \ri \} \leq  \sum_{\ell = 1}^\tau \Pr \li \{ |
\wh{\bs{p}}_\ell - p |
\geq \vep p \ri \} \nonumber\\
& \leq &  \sum_{\ell = 1}^\tau 2 \exp \li ( \ga_\ell
\li [ \f{\vep}{1 + \vep} - \ln (1 + \vep) \ri ] \ri ) \la{RineCH}\\
& < & 2 \tau \exp \li ( \ga_1 \li [ \f{\vep}{1 + \vep} - \ln (1 +
\vep) \ri ] \ri ), \la{RLastB} \eel where (\ref{RineCH}) is due to
Corollary of \cite{Chen3}. As can be seen from (\ref{RLastB}), the
last bound is independent of $p$ and can be made smaller than $\de$
if $\ze$ is sufficiently small.  This establishes the claim and it
follows that $\Pr \li \{ \li | \bs{\wh{p}} - p \ri | < \vep p \mid p
\ri \}
> 1 - \de$ for any $p \in (0, 1)$ if $\ze$ is sufficiently
small.

To improve the performance of coverage probability, we propose to
revise the stopping rule associated with (\ref{invrev}) as follows:
For $\ell = 1, \cd, s$, decision variable $\bs{D}_\ell$ assumes
value $1$ if  \be \la{invnora}
 \li ( \li | \wh{\bs{p}}_\ell -
\f{1}{2} \ri | - w \vep \wh{\bs{p}}_\ell \ri )^2 \geq \f{1}{4} + \f{
\vep^2 \ga_\ell \wh{\bs{p}}_\ell } { 2 \ln (\ze \de) } \ee is
satisfied; and assumes value $0$ otherwise, where $w \geq 0$ is a
parameter affecting the shape of the stopping boundary.

Before concluding this subsection, we would like to point out that
it is possible to modify (\ref{invnora}) to produce the following
stopping rule: For $\ell = 1, \cd, s$, decision variable
$\bs{D}_\ell$ assumes value $1$ if \be \la{invnorb} \li (
\wh{\bs{p}}_\ell - \f{1}{2} \ri )^2 \geq \f{1}{4} + \f{w
\wh{\bs{p}}_\ell}{\ga_\ell}  + \f{ \vep^2 \ga_\ell \wh{\bs{p}}_\ell
} { 2 \ln (\ze \de) } \ee is satisfied; and assumes value $0$
otherwise, where $w \geq 0$ is a parameter affecting the shape of
the stopping boundary.

The thresholds $\ga_1 < \cd < \ga_s$ for the two stopping rules
associated with (\ref{invnora}) and (\ref{invnorb}) can be chosen in
a similar spirit as suggested in Section 2.1.  Specifically,  the
maximum threshold of sample sum $\ga_s$, i.e., the threshold of
sample sum of the last stage, should be defined as the smallest
integer such that $\{ \bs{D}_s = 1 \}$ is a sure event. The minimum
threshold of sample sum $\ga_1$, i.e., the threshold of sample sum
of the first stage, should be defined as the smallest integer such
that $\{ \bs{D}_1 = 1 \}$ is an event of a positive probability. For
both stopping rules, we can show that $\Pr \li \{ \li | \bs{\wh{p}}
- p \ri | < \vep p \mid p \ri \}
> 1 - \de$ for any $p \in (0, 1)$ if $\ze$ is sufficiently small.

\subsubsection{Noninverse Multistage Sampling}

In Sections \ref{multiinv} and \ref{Bino_inv_asp}, we have proposed
a multistage inverse sampling plan for estimating a binomial
parameter, $p$, with relative precision. In some situations, the
cost of sampling operation may be high since samples are obtained
one by one when inverse sampling is involved. In view of this fact,
it is desirable to develop multistage estimation methods without
using inverse sampling.

In contrast to the multistage inverse sampling schemes described in
Sections \ref{multiinv} and \ref{Bino_inv_asp}, our noninverse
multistage sampling schemes have infinitely many stages and
deterministic sample sizes $n_1 < n_2 < n_3 < \cd$. Moreover, the
confidence parameter for the $\ell$-th stage, $\de_\ell$,  is
dependent on $\ell$ such that $\de_\ell = \de$ for $1 \leq \ell \leq
\tau$ and $\de_\ell = \de 2^{\tau - \ell}$ for $\ell > \tau$,  where
$\tau$ is a positive integer.

By virtue of the CDF $\&$ CCDF of $\wh{\bs{p}}_\ell$, we propose a
class of multistage sampling schemes as follows.

\beT \la{Bino_Rev_noninverse_CDF}   Suppose that, for $\ell = 1, 2,
\cd$, decision variable $\bs{D}_\ell$ assumes values $1$ if {\small
$F_{ \wh{\bs{p}}_\ell } ( \wh{\bs{p}}_\ell, \f{\wh{\bs{p}}_\ell}{1 -
\vep} ) \leq \ze \de_\ell, \; G_{ \wh{\bs{p}}_\ell } (
\wh{\bs{p}}_\ell, \f{\wh{\bs{p}}_\ell}{1 + \vep} ) \leq \ze
\de_\ell$}; and assumes $0$ otherwise. The following statements hold
true.

(I): $\Pr \{ \mbf{n} < \iy \} = 1$  provided that $\inf_{\ell > 0}
\f{n_{\ell + 1}}{n_\ell} > 1$.

(II): $\bb{E} [ \mbf{n} ] < \iy$ provided that $1 < \inf_{\ell > 0}
\f{n_{\ell + 1}}{n_\ell} \leq \sup_{\ell > 0} \f{n_{\ell +
1}}{n_\ell} < \iy$.

(III): {\small $\Pr \li \{ \li | \f{ \wh{\bs{p}} - p } { p } \ri | <
\vep \mid p \ri \} \geq 1 - \de$} for any $p \in (0, 1)$ provided
that $\ze \leq \f{1}{ 2 (\tau + 1 ) }$.

(IV):  Let $0 < \eta < \ze \de$ and {\small ${\ell^\star} = \tau + 1
+ \li \lc \f{ \ln ( \ze \de \sh \eta )  } { \ln 2 } \ri \rc$}. Then,
$\Pr \li \{ \li | \wh{\bs{p}} - p \ri | \geq \vep p \ri \} < \de$
for any $p \in (0, p^*)$, where $p^*$ is a number such that {\small
$0 < p^* < z_\ell, \; \; \ell = 1, \cd, {\ell^\star}$} and that
$\sum_{\ell = 1}^{\ell^\star} \exp ( n_\ell \mscr{M}_{\mrm{B}} (
z_\ell, p^*) ) < \de - \eta$ with $z_\ell = \min \{  z \in
I_{\wh{\bs{p}}_\ell}: F_{ \wh{\bs{p}}_\ell }  (z, \f{z}{1 - \vep} )
> \ze \de_\ell \; \tx{or} \; G_{ \wh{\bs{p}}_\ell } ( z, \f{z}{1
+ \vep} ) > \ze \de_\ell \}$, where $I_{\wh{\bs{p}}_\ell}$
represents the support of $\wh{\bs{p}}_\ell$, for $\ell = 1, 2,
\cd$. Moreover, {\small \bee &  & \Pr \li \{ b \leq \f{ \wh{\bs{p}}
}{ 1 + \vep}, \; \bs{l} \leq \ell^\star \mid a \ri \} \leq \Pr \li
\{ p \leq \f{ \wh{\bs{p}} }{ 1 + \vep} \mid p \ri \} \leq
\f{\eta}{2} + \Pr
\li \{ a \leq \f{ \wh{\bs{p}} }{ 1 + \vep}, \; \bs{l} \leq \ell^\star \mid b \ri \}, \\
&   & \Pr \li \{ a \geq \f{ \wh{\bs{p}} }{ 1 - \vep}, \; \bs{l} \leq
\ell^\star \mid b \ri \}  \leq \Pr \li \{ p \geq \f{ \wh{\bs{p}} }{
1 - \vep} \mid p   \ri \} \leq \f{\eta}{2} + \Pr \li \{ b \geq \f{
\wh{\bs{p}} }{ 1 - \vep}, \; \bs{l} \leq \ell^\star \mid a \ri \}
\eee} for any $p \in [a, b]$, where $a$ and $b$ are numbers such
that $0 < b < ( 1 + \vep ) a < 1$.

(V): Let the sample sizes of the multistage sampling scheme be a
sequence $n_\ell = \li \lc m \ga^{\ell - 1} \ri \rc, \; \ell = 1, 2,
\cd$, where $\ga \geq 1 + \f{1}{m} > 1$. Let $0 < \ep < \f{1}{2}, \;
0 < \eta < 1$ and $c = \f{p (1 - \eta)^2}{2}$. Let $\ka$ be an
integer such that {\small $\ka > \max \li \{ \tau, \;  \f{1}{\ln
\ga} \ln \li ( \f{1}{ c m } \ln \f{\ga}{c \ep} \ri ) + 1, \; \tau +
\f{1}{\ga - 1} + \f{ \ln (\ze \de) } { \ln 2 }  \ri \} $} and
{\small $\mscr{M}_{\mrm{B}} ( \eta p, \f{\eta p}{1 + \vep} ) < \f{
\ln (\ze \de_\ka )}{n_\ka}$}.  Then, $\bb{E} [\mbf{n}] < \ep + n_1 +
\sum_{\ell = 1}^{\ka} (n_{\ell + 1} - n_\ell) \Pr \{ \bs{l}
> \ell \}$.

\eeT

The proof of Theorem \ref{Bino_Rev_noninverse_CDF} is similar to that of Theorem \ref {Bino_Rev_noninverse_Chernoff}, which is given at Appendix
\ref{App_Bino_Rev_noninverse_Chernoff}.

By virtue of Chernoff bounds of the CDF $\&$ CCDF of
$\wh{\bs{p}}_\ell$, we propose a class of multistage sampling
schemes as follows.

\beT \la{Bino_Rev_noninverse_Chernoff}   Suppose that, for $\ell =
1, 2, \cd$, decision variable $\bs{D}_\ell$ assumes values $1$ if
{\small $\mscr{M}_{\mrm{B}} (\wh{\bs{p}}_\ell,
\f{\wh{\bs{p}}_\ell}{1 + \vep} )  \leq \f{ \ln ( \ze \de_\ell ) } {
n_\ell }$}; and assumes $0$ otherwise. The following statements hold
true.

(I): $\Pr \{ \mbf{n} < \iy \} = 1$  provided that $\inf_{\ell > 0}
\f{n_{\ell + 1}}{n_\ell} > 1$.

(II): $\bb{E} [ \mbf{n} ] < \iy$ provided that $1 < \inf_{\ell > 0}
\f{n_{\ell + 1}}{n_\ell} \leq \sup_{\ell > 0} \f{n_{\ell +
1}}{n_\ell} < \iy$.

(III): {\small $\Pr \li \{ \li | \f{ \wh{\bs{p}} - p } { p } \ri | <
\vep \mid p \ri \} \geq 1 - \de$} for any $p \in (0, 1)$ provided
that $ \ze \leq \f{1}{ 2 (\tau + 1 ) }$.

(IV):  Let $0 < \eta < \ze \de$ and {\small ${\ell^\star} = \tau + 1
+ \li \lc \f{ \ln ( \ze \de \sh \eta )  } { \ln 2 } \ri \rc$}. Then,
$\Pr \li \{ \li | \wh{\bs{p}} - p \ri | \geq \vep p \ri \} < \de$
for any $p \in (0, p^*)$, where $p^*$ is a number such that {\small
$0 < p^* < z_\ell, \; \; \ell = \tau, \cd, {\ell^\star}$} and that
$\sum_{\ell = 1}^{\ell^\star} \exp ( n_\ell \mscr{M}_{\mrm{B}} (
z_\ell, p^*) ) < \de - \eta$ with $z_\ell$ satisfying {\small
$\mscr{M}_{\mrm{B}} \li (z_\ell, \f{z_\ell}{1 + \vep} \ri  )  = \f{
\ln ( \ze \de_\ell ) } { n_\ell }$} for $\ell = 1, 2, \cd$.
Moreover, {\small \bee & & \Pr \li \{ b \leq \f{ \wh{\bs{p}} }{ 1 +
\vep}, \; \bs{l} \leq \ell^\star \mid a \ri \} \leq \Pr \li \{ p
\leq \f{ \wh{\bs{p}} }{ 1 + \vep} \mid p \ri \} \leq \f{\eta}{2} +
\Pr
\li \{ a \leq \f{ \wh{\bs{p}} }{ 1 + \vep}, \; \bs{l} \leq \ell^\star \mid b \ri \}, \\
&   & \Pr \li \{ a \geq \f{ \wh{\bs{p}} }{ 1 - \vep}, \; \bs{l} \leq
\ell^\star \mid b \ri \}  \leq \Pr \li \{ p \geq \f{ \wh{\bs{p}} }{
1 - \vep} \mid p   \ri \} \leq \f{\eta}{2} + \Pr \li \{ b \geq \f{
\wh{\bs{p}} }{ 1 - \vep}, \; \bs{l} \leq \ell^\star \mid a \ri \}
\eee} for any $p \in [a, b]$, where $a$ and $b$ are numbers such
that $0 < b < ( 1 + \vep ) a < 1$.

(V): Let the sample sizes of the multistage sampling scheme be a
sequence $n_\ell = \li \lc m \ga^{\ell - 1} \ri \rc, \; \ell = 1, 2,
\cd$, where $\ga \geq 1 + \f{1}{m} > 1$. Let $0 < \ep < \f{1}{2}, \;
0 < \eta < 1$ and $c = \f{p (1 - \eta)^2}{2}$. Let $\ka$ be an
integer such that {\small $\ka > \max \li \{ \tau, \;  \f{1}{\ln
\ga} \ln \li ( \f{1}{ c m } \ln \f{\ga}{c \ep} \ri ) + 1, \; \tau +
\f{1}{\ga - 1} + \f{ \ln (\ze \de) } { \ln 2 }  \ri \} $} and
{\small $\mscr{M}_{\mrm{B}} ( \eta p, \f{\eta p}{1 + \vep} ) < \f{
\ln (\ze \de_\ka )}{n_\ka}$}.  Then, $\bb{E} [\mbf{n}] < \ep + n_1 +
\sum_{\ell = 1}^{\ka} (n_{\ell + 1} - n_\ell) \Pr \{ \bs{l}
> \ell \}$.

\eeT

See Appendix \ref{App_Bino_Rev_noninverse_Chernoff} for a proof.

By virtue of Massart's inequality for the CDF $\&$ CCDF of
$\wh{\bs{p}}_\ell$, we propose a class of multistage sampling
schemes as follows.

\beT \la{Bino_Rev_noninverse_Massart}   Suppose that, for $\ell = 1,
2, \cd$, decision variable $\bs{D}_\ell$ assumes values $1$ if
{\small $\wh{\bs{p}}_\ell \geq \f{ 6(1 + \vep) (3 + \vep) \ln (\ze
\de_\ell) } { 2 (3 + \vep)^2 \ln (\ze \de_\ell) - 9 n_\ell
\vep^2}$}; and assumes $0$ otherwise. The following statements hold
true.

(I): $\Pr \{ \mbf{n} < \iy \} = 1$  provided that $\inf_{\ell > 0}
\f{n_{\ell + 1}}{n_\ell} > 1$.

(II): $\bb{E} [ \mbf{n} ] < \iy$ provided that $1 < \inf_{\ell > 0}
\f{n_{\ell + 1}}{n_\ell} \leq \sup_{\ell > 0} \f{n_{\ell +
1}}{n_\ell} < \iy$.

(III): {\small $\Pr \li \{ \li | \f{ \wh{\bs{p}} - p } { p } \ri | <
\vep \mid p \ri \} \geq 1 - \de$} for any $p \in (0, 1)$ provided
that $ \ze \leq \f{1}{ 2 (\tau + 1 ) }$.

(IV):  Let $0 < \eta < \ze \de$ and {\small ${\ell^\star} = \tau + 1
+ \li \lc \f{ \ln ( \ze \de \sh \eta )  } { \ln 2 } \ri \rc$}. Then,
$\Pr \li \{ \li | \wh{\bs{p}} - p \ri | \geq \vep p \ri \} < \de$
for any $p \in (0, p^*)$, where $p^*$ is a number such that {\small
$0 < p^* < z_\ell, \; \; \ell = \tau, \cd, {\ell^\star}$} and that
$\sum_{\ell = 1}^{\ell^\star} \exp ( n_\ell \mscr{M} ( z_\ell, p^*)
) < \de - \eta$ with {\small $z_\ell = \f{ 6 (1 + \vep) (3 + \vep)
\ln (\ze \de_\ell) } { 2 (3 + \vep)^2 \ln (\ze \de_\ell) -  9 \vep^2
n_\ell }$} for $\ell = 1, 2, \cd$.  Moreover, {\small \bee & & \Pr
\li \{ b \leq \f{ \wh{\bs{p}} }{ 1 + \vep}, \; \bs{l} \leq
\ell^\star \mid a \ri \} \leq \Pr \li \{ p \leq \f{ \wh{\bs{p}} }{ 1
+ \vep} \mid p \ri \} \leq \f{\eta}{2} + \Pr
\li \{ a \leq \f{ \wh{\bs{p}} }{ 1 + \vep}, \; \bs{l} \leq \ell^\star \mid b \ri \}, \\
&   & \Pr \li \{ a \geq \f{ \wh{\bs{p}} }{ 1 - \vep}, \; \bs{l} \leq
\ell^\star \mid b \ri \}  \leq \Pr \li \{ p \geq \f{ \wh{\bs{p}} }{
1 - \vep} \mid p   \ri \} \leq \f{\eta}{2} + \Pr \li \{ b \geq \f{
\wh{\bs{p}} }{ 1 - \vep}, \; \bs{l} \leq \ell^\star \mid a \ri \}
\eee} for any $p \in [a, b]$, where $a$ and $b$ are numbers such
that $0 < b < ( 1 + \vep ) a < 1$.

(V): Let the sample sizes of the multistage sampling scheme be a
sequence $n_\ell = \li \lc m \ga^{\ell - 1} \ri \rc, \; \ell = 1, 2,
\cd$, where $\ga \geq 1 + \f{1}{m} > 1$. Let $0 < \ep < \f{1}{2}, \;
0 < \eta < 1$ and $c = \f{p (1 - \eta)^2}{2}$. Let $\ka$ be an
integer such that {\small $\ka > \max \li \{ \tau, \;  \f{1}{\ln
\ga} \ln \li ( \f{1}{ c m } \ln \f{\ga}{c \ep} \ri ) + 1, \; \tau +
\f{1}{\ga - 1} + \f{ \ln (\ze \de) } { \ln 2 }  \ri \} $} and
{\small $\mscr{M} ( \eta p, \f{\eta p}{1 + \vep} ) < \f{ \ln (\ze
\de_\ka )}{n_\ka}$}.  Then, $\bb{E} [\mbf{n}] < \ep + n_1 +
\sum_{\ell = 1}^{\ka} (n_{\ell + 1} - n_\ell) \Pr \{ \bs{l}
> \ell \}$.

\eeT

The proof of Theorem \ref{Bino_Rev_noninverse_Massart} is similar to
that of Theorem \ref {Bino_Rev_noninverse_Chernoff}, which is given at
Appendix \ref{App_Bino_Rev_noninverse_Chernoff}.

\subsubsection{Asymptotic Analysis of Multistage Inverse  Sampling Schemes}

In this subsection, we shall focus on the asymptotic analysis of
multistage inverse sampling schemes.  Throughout this subsection, we
assume that the multistage inverse sampling schemes follow stopping
rules derived from Chernoff bounds as described in Section
\ref{multiinv}. Moreover, we assume that the thresholds of sample
sum $\ga_1, \cd, \ga_s$ are chosen as the ascending arrangement of
all distinct elements of  the set defined by (\ref{THolds}).

With regard to the tightness of the double-decision-variable method,
we have

\beT \la{Bino_Inverse_DDV_Asp} Let $\mscr{R}$ be a subset of real
numbers. Define {\small
\[ \ovl{P} = \sum_{\ell = 1}^s \Pr \{ \wh{\bs{p}}_\ell \in \mscr{R},
\; \bs{D}_{\ell - 1} = 0, \; \bs{D}_\ell = 1 \}, \qqu \udl{P} = 1 -
\sum_{\ell = 1}^s \Pr \{ \wh{\bs{p}}_\ell \notin \mscr{R}, \;
\bs{D}_{\ell - 1} = 0, \; \bs{D}_\ell = 1 \}.
\]}
Then, $\udl{P} \leq \Pr \{ \wh{\bs{p}} \in \mscr{R} \} \leq \ovl{P}$
and $\lim_{\vep \to 0} | \Pr \{ \wh{\bs{p}} \in \mscr{R} \} -
\ovl{P} | = \lim_{\vep \to 0} | \Pr \{ \wh{\bs{p}} \in \mscr{R} \} -
\udl{P} | = 0$ for any $p \in (0, 1)$.  \eeT

\bsk

See Appendix \ref{App_Bino_Inverse_DDV_Asp} for a proof.

Recall that $\bs{l}$ is the index of stage when the sampling is
terminated. Define $\bs{\ga} = \ga_{\bs{l}}$.  Then, $\bs{\ga} =
\sum_{i=1}^{ \mathbf{n} } X_i$.   With regard to the asymptotic
performance of the sampling scheme, we have

\beT \la{Bino_Inverse_Asp_Analysis}  Let {\small $\ga (p, \vep) =
\f{ \ln ( \ze \de ) } { \mscr{M}_{\mrm{I}} \li ( p, \f{ p } { 1 +
\vep } \ri ) } $}.  Let $\mcal{N}_{\mrm{f}} (p, \vep)$ be the
minimum sample number $n$ such that {\small $\Pr \{ | \f{\sum_{i =
1}^n X_i}{n} - p | < \vep p \mid p \} > 1 - \ze \de$} for a
fixed-size sampling procedure.  Let $j_p$ be the maximum integer $j$
such that $C_j \geq 1 - p$.  Let {\small $\nu = \f{2}{3}, \; d =
\sq{ 2 \ln \f{1}{\ze \de} }$} and {\small $\ka_p = \f{C_{j_p} }{1 -
p}$}. Let {\small $\ro_p = \f{C_{j_p - 1} }{1 - p} - 1$} if $\ka_p =
1$ and $\ro_p = \ka_p - 1$ otherwise.   The following statements
hold true:

(I): {\small $\Pr \li \{  1 \leq \limsup_{\vep \to 0} \f{ \bs{\ga} }
{ \ga (p, \vep) } \leq 1 + \ro_p \ri \} = 1$}.  Specially, {\small
$\Pr \li \{  \lim_{\vep \to 0} \f{ \bs{\ga} } { \ga (p, \vep) } =
\ka_p \ri \} = 1$} if $\ka_p > 1$.

(II): $\lim_{\vep \to 0} \f{ \bb{E} [ \mbf{n} ] } {
\mcal{N}_{\mrm{f}} (p, \vep)} = \li ( \f{ d } { \mcal{Z}_{\ze \de} }
\ri )^2 \times \lim_{\vep \to 0} \f{ \bb{E} [ \bs{\ga} ] } { \ga (p,
\vep) }$, where
\[
\lim_{\vep \to 0} \f{ \bb{E} [ \bs{\ga} ] } { \ga (p, \vep) }   =
\bec \ka_p &
\tx{if} \; \ka_p > 1,\\
1 + \ro_p \Phi (\nu d ) &  \tx{otherwise} \eec
\]
and $1 \leq \lim_{\vep \to 0} \f{ \bb{E} [ \bs{\ga} ] } { \ga (p,
\vep) } \leq 1 + \ro_p$.

(III): If $\ka_p > 1$, then $\lim_{\vep \to 0} \Pr \{ | \wh{\bs{p}} - p | < \vep p \} = 2 \Phi \li ( d \sq{\ka_p} \ri ) - 1 > 2 \Phi \li ( d \ri
) - 1 > 1 - 2 \ze \de$. Otherwise, $\Phi \li ( d  \ri )  + \Phi \li ( d \sq{1 + \ro_p}  \ri )  - 1  > \lim_{\vep \to 0} \Pr \{ | \wh{\bs{p}} - p
| < \vep p \} = 1 + \Phi(d) - \Phi(\nu d) - \Psi (\ro_p, \nu, d) > \Phi \li ( d \ri ) + 2 \Phi \li ( d \sq{1 + \ro_p}  \ri )  - 2 > 1 - 3 \ze
\de$.

 \eeT

See Appendix \ref{App_Bino_Inverse_Asp_Analysis}.

\subsubsection{Asymptotic Analysis of Noninverse Multistage Sampling Schemes}

In this subsection, we shall focus on the asymptotic analysis of the
noninverse multistage sampling schemes which follow stopping rules
derived from Chernoff bounds of CDF $\&$ CCDF of $\wh{\bs{p}}_\ell$
as described in Theorem \ref{Bino_Rev_noninverse_Chernoff}.

We assume that the sample sizes $n_1, n_2, \cd$ are chosen as the
ascending arrangement of all distinct elements of the set {\small
\be \la{defss} \li \{ \li \lc   \f{ C_{\tau - \ell} \; \ln (\ze \de)
}{ \mscr{M}_{\mrm{B}} (p^*, \f{p^*}{1 + \vep} ) } \ri \rc : \ell =
1, 2, \cd \ri \} \ee} with $p^* \in (0, 1)$, where $\tau$ is the
maximum integer such that $\f{ C_{\tau - 1} \ln (\ze \de) } {
\mscr{M}_{\mrm{B}} (p^*, \f{p^*}{1 + \vep} ) } \geq \f{ \ln
\f{1}{\ze \de} } { \ln (1 + \vep) }$, i.e., $C_{\tau - 1} \geq -
\f{\mscr{M}_{\mrm{B}} (p^*, \f{p^*}{1 + \vep} )  } { \ln (1 + \vep)
}$.

With regard to the asymptotic performance of the sampling scheme, we
have

\beT \la{Bino_Asp_Analysis_Noninverse}  Let {\small $\mcal{N}_{\mrm{r}} (p, \vep) = \f{ \ln ( \ze \de  ) } { \mscr{M}_{\mrm{B}} ( p , \f{p}{1 +
\vep}) } $}. Let $\mcal{N}_{\mrm{f}} (p, \vep)$ be the minimum sample number $n$ such that {\small $\Pr \{ | \f{\sum_{i = 1}^n X_i}{n} - p  | <
\vep p \mid p \} > 1 - \ze \de$} for a fixed-size sampling procedure.  Let $j_p$ be the maximum integer $j$ such that $C_j \geq r (p)$, where
{\small $r (p) = \f{p^* (1 - p)}{p (1 - p^*)}$}. Let {\small $\nu = \f{2}{3} \f{p - p^*}{1 - p^*}, \; d = \sq{ 2 \ln \f{1}{\ze \de} }$} and
{\small $\ka_p = \f{C_{j_p}}{r (p)} $}.  Let {\small $\ro_p = \f{C_{j_p - 1}} { r (p) } - 1$} if $\ka_p = 1$ and $\ro_p  = \ka_p - 1$ otherwise.
For $p \in (p^*, 1)$, the following statements hold true:

(I): {\small $\Pr \li \{  1 \leq \limsup_{\vep \to 0} \f{ \mbf{n} }
{ \mcal{N}_{\mrm{r}} (p, \vep) } \leq 1 + \ro_p  \ri \} = 1$}.
Specially, {\small $\Pr \li \{  \lim_{\vep \to 0} \f{ \mbf{n} } {
\mcal{N}_{\mrm{r}} (p, \vep) } = \ka_p  \ri \} = 1$} if $\ka_p > 1$.

(II): {\small $\lim_{\vep \to 0} \f{ \bb{E} [ \mbf{n} ] } {
\mcal{N}_{\mrm{f}} (p, \vep)} = \li ( \f{ d } { \mcal{Z}_{\ze \de} }
\ri )^2 \times \lim_{\vep \to 0} \f{ \bb{E} [ \mbf{n} ] } {
\mcal{N}_{\mrm{r}} (p, \vep) }$}, where
\[ \lim_{\vep \to 0} \f{ \bb{E} [ \mbf{n} ] }  {
\mcal{N}_{\mrm{r}} (p, \vep) } = \bec \ka_p
& \tx{if} \; \ka_p > 1,\\
1  + \ro_p \Phi (\nu d ) &  \tx{otherwise} \eec
\]
and $1 \leq \lim_{\vep \to 0} \f{ \bb{E} [ \mbf{n} ] }  {
\mcal{N}_{\mrm{r}} (p, \vep) } \leq 1 + \ro_p$.

(III): If $\ka_p > 1$, then $\lim_{\vep \to 0} \Pr \{ | \wh{\bs{p}} - p | < \vep p \} = 2 \Phi \li ( d \sq{\ka_p} \ri ) - 1 > 2 \Phi \li ( d \ri
) - 1 > 1 - 2 \ze \de$. Otherwise, $\Phi \li ( d  \ri )  + \Phi \li ( d \sq{1 + \ro_p}  \ri )  - 1  > \lim_{\vep \to 0} \Pr \{ | \wh{\bs{p}} - p
| < \vep p \} = 1 + \Phi(d) - \Phi(\nu d) - \Psi (\ro_p, \nu, d) > \Phi \li ( d \ri ) + 2 \Phi \li ( d \sq{1 + \ro_p}  \ri )  - 2  > 1 - 3 \ze
\de$.

 \eeT

\bsk

See Appendix \ref{App_Bino_Asp_Analysis_Noninverse} for a proof.

\subsection{Control of Absolute and Relative Errors}

In this section, we shall focus on the design of multistage sampling
schemes for estimating the binomial parameter $p$ with a mixed error
criterion.  Specifically, for $0 < \vep_a < 1$ and $0 < \vep_r < 1$,
we wish to construct a multistage sampling scheme and its associated
estimator $\wh{\bs{p}}$ for $p$ such that $\Pr \{ | \wh{\bs{p}} - p
| < \vep_a, \;  | \wh{\bs{p}} - p | < \vep_r p \mid p \} > 1 - \de$
for any $p \in (0, 1)$.  This is equivalent to the construction of a
random interval with lower limit $\mscr{L} ( \wh{\bs{p}} )$ and
upper limit $\mscr{U} ( \wh{\bs{p}} )$ such that $\Pr \{  \mscr{L} (
\wh{\bs{p}} ) < p <  \mscr{U} ( \wh{\bs{p}} ) \mid p \} > 1 - \de$
for any $p \in (0, 1)$, where $\mscr{L} ( . )$ and $\mscr{U} ( . )$
are functions such that $\mscr{L} ( z ) = \min \{ z - \vep_a, \;
\f{z}{1 + \vep_r} \}$ and $\mscr{U} ( z ) = \max \{ z + \vep_a, \;
\f{z}{1 - \vep_r} \}$ for $z \in [0, 1]$.  In the sequel, we shall
propose multistage sampling schemes such that the number of stages,
$s$, is finite and that the sample sizes are deterministic numbers
$n_1 < n_2 < \cd < n_s$.

\subsubsection{Stopping Rules from CDF $\&$ CCDF and Chernoff Bounds} \la{multiinv_mix}

To construct an estimator satisfying a mixed criterion in terms of
absolute and relative errors with a prescribed confidence level, we
have developed two types of multistage sampling schemes with
different stopping rules as follows.

\bed

\item [Stopping Rule (i):] For $\ell = 1, \cd, s$, decision variable
$\bs{D}_\ell$  assumes value $1$ if {\small $F_{\wh{\bs{p}}_\ell}
(\wh{\bs{p}}_\ell, \mscr{U} ( \wh{\bs{p}}_\ell ) ) \leq \ze \de, \;
G_{\wh{\bs{p}}_\ell} (\wh{\bs{p}}_\ell, \mscr{L} ( \wh{\bs{p}}_\ell
) ) \leq \ze \de$}; and assumes value $0$ otherwise.

\item [Stopping Rule (ii):] For $\ell = 1, \cd, s$, decision variable
$\bs{D}_\ell$ assumes value $1$  if

{\small $\max \{ \mscr{M}_{\mrm{B}} (\wh{\bs{p}}_\ell, \mscr{L} (
\wh{\bs{p}}_\ell ) ), \; \mscr{M}_{\mrm{B}} (\wh{\bs{p}}_\ell,
\mscr{U} ( \wh{\bs{p}}_\ell ) ) \} \leq \f{ \ln ( \ze \de  ) } {
n_\ell }$}; and assumes value $0$ otherwise.

\eed

Stopping rule (i) is derived by virtue of the CDF $\&$ CCDF of
$\wh{\bs{p}}_\ell$.  Stopping rule (ii) is derived by virtue of
Chernoff bounds of the CDF $\&$ CCDF of $\wh{\bs{p}}_\ell$.   For
both types of multistage sampling schemes described above, we have
the following results.

\beT \la{Bino_mix_CDF_CH}   Let $\vep_a$ and $\vep_r$ be positive
numbers such that $0 < \vep_a < \f{35}{94}$ and {\small $\f{70
\vep_a}{35 - 24 \vep_a } < \vep_r < 1$}.  Suppose that the sample
size for the $s$-th stage is no less than {\small $\li \lc \f{ \ln
(\ze \de) } { \mscr{M}_{\mrm{B}} (\f{\vep_a}{\vep_r} + \vep_a,
\f{\vep_a}{\vep_r} ) } \ri \rc$}. Then, \bee & & \Pr \{ p \leq
\mscr{L} ( \wh{\bs{p}} ) \mid p \} \leq \sum_{\ell = 1}^s \Pr \{  p
\leq \mscr{L} ( \wh{\bs{p}}_\ell ), \;
\bs{D}_\ell = 1 \mid p \} \leq s \ze \de,\\
&    & \Pr \{  p \geq \mscr{U} ( \wh{\bs{p}} ) \mid p \} \leq
\sum_{\ell = 1}^s  \Pr \{  p \geq \mscr{U} ( \wh{\bs{p}}_\ell ), \;
\bs{D}_\ell = 1 \mid p \} \leq s \ze \de \eee and $\Pr \{  |
\wh{\bs{p}} - p | < \vep_a \; \tx{or} \;  | \wh{\bs{p}} - p | <
\vep_r p \mid p \} \geq 1 - 2 s \ze \de $ for any $p \in (0, 1)$.
\eeT

\bsk

See Appendix \ref{App_Bino_mix_CDF_CH} for a proof.  Based on the
criteria proposed in Section \ref{gen_structure}, the sample sizes
$n_1 < n_2 < \cd < n_s$ can be chosen as the ascending arrangement
of all distinct elements of the set {\small \be \la{Bino_mix_ss} \li
\{ \li \lc \f{ C_{\tau - \ell} \; \ln (\ze \de) }{
\mscr{M}_{\mrm{B}} (\f{\vep_a}{\vep_r} + \vep_a, \f{\vep_a}{\vep_r}
) } \ri \rc : \ell = 1, \cd, \tau \ri \}, \ee} where $\tau$ is the
maximum integer such that $\f{ C_{\tau - \ell} \; \ln (\ze \de) }{
\mscr{M}_{\mrm{B}} (\f{\vep_a}{\vep_r} + \vep_a, \f{\vep_a}{\vep_r}
) } \geq \f{ \ln \f{1}{\ze \de} }{ \ln (1 + \vep_r ) }$, i.e.,
$C_{\tau - 1} \geq - \f{ \mscr{M}_{\mrm{B}} (\f{\vep_a}{\vep_r} +
\vep_a, \f{\vep_a}{\vep_r} ) } { \ln (1 + \vep_r) }$.

For such a choice of sample sizes, as a result of Theorem
\ref{Bino_mix_CDF_CH}, we have that $\Pr \{  | \wh{\bs{p}} - p | <
\vep_a \; \tx{or} \;  | \wh{\bs{p}} - p | < \vep_r p \mid p \}
> 1 - \de $ for any $p \in (0, 1)$ provided that $\ze < \f{1}{2
\tau}$.

For computing the coverage probability  associated with a multistage
sampling scheme following a stopping rule derived from Chernoff
bounds, events $\{ \bs{D}_\ell = i \}, \; i = 0, 1$ need to be
expressed as events involving only $K_\ell$. This can be
accomplished by using the following results.

\beT \la{Bino_range_mix}

Let $p^\star = \f{\vep_a}{\vep_r}$.  For $\ell = 1, \cd, s - 1$,
{\small $\{ \bs{D}_\ell = 0 \}
 =  \{ \mscr{M}_{\mrm{B}} (\wh{\bs{p}}_\ell, \mscr{L} ( \wh{\bs{p}}_\ell ) )  >
\f{ \ln ( \ze \de  ) } { n_\ell }  \} \cup  \{ \mscr{M}_{\mrm{B}}
(\wh{\bs{p}}_\ell, \mscr{U} ( \wh{\bs{p}}_\ell ) )
> \f{ \ln ( \ze \de  ) } { n_\ell } \}$} and the following
statements hold true:

(I) {\small $\{  \mscr{M}_{\mrm{B}} (\wh{\bs{p}}_\ell, \mscr{L} (
\wh{\bs{p}}_\ell ) ) > \f{ \ln ( \ze \de ) } { n_\ell }  \} =  \{
n_\ell \; z_a^- < K_\ell < n_\ell \; z_r^+ \}$} where $z_r^+$ is the
unique solution of equation {\small $\mscr{M}_{\mrm{B}} (z, \f{z}{1
+ \vep_r} ) = \f{ \ln ( \ze \de  ) } { n_\ell }$} with respect to $z
\in (p^\star + \vep_a, 1]$, and $z_a^-$ is the unique solution of
equation {\small $\mscr{M}_{\mrm{B}} (z, z - \vep_a ) = \f{ \ln (
\ze \de  ) } { n_\ell }$} with respect to $z \in (\vep_a, p^\star +
\vep_a)$.

(II) {\small \[ \li \{ \mscr{M}_{\mrm{B}} (\wh{\bs{p}}_\ell,
\mscr{U} ( \wh{\bs{p}}_\ell ) )  > \f{ \ln ( \ze \de ) } { n_\ell }
\ri \} = \bec  \{ 0 \leq K_\ell < n_\ell \; z_r^- \} &
\tx{for} \;  n_\ell < \f{ \ln (\ze \de) } { \ln (1 - \vep_a) },\\
\{ n_\ell \; z_a^+ < K_\ell < n_\ell \; z_r^- \} & \tx{for} \; \f{
\ln (\ze \de) } { \ln (1 - \vep_a) } \leq n_\ell < \f{ \ln (\ze \de)
}
{ \mscr{M}_{\mrm{B}} (p^\star - \vep_a, p^\star) },\\
\emptyset & \tx{for} \; n_\ell \geq \f{ \ln (\ze \de) } {
\mscr{M}_{\mrm{B}} (p^\star - \vep_a, p^\star) }
 \eec
\]}
where $z_r^-$ is the unique solution of equation {\small
$\mscr{M}_{\mrm{B}} (z, \f{z}{1 - \vep_r} ) = \f{ \ln ( \ze \de  ) }
{ n_\ell }$} with respect to $z \in (p^\star - \vep_a, 1 - \vep_r)$,
and $z_a^+$ is the unique solution of equation {\small
$\mscr{M}_{\mrm{B}} (z, z + \vep_a ) = \f{ \ln ( \ze \de  ) } {
n_\ell }$} with respect to $z \in [0, p^\star - \vep_a)$.

\eeT

\bsk

See Appendix \ref{App_Bino_range_mix} for a proof.

\subsubsection{Stopping Rule from Massart's Inequality}

 By virtue of Massart's inequality of the CDF $\&$ CCDF of $\wh{\bs{p}}_\ell$, we can construct
 a multistage sampling scheme such that its associated estimator for $p$ satisfies the mixed
 criterion.  Such a sampling scheme and its properties are described
 by the following theorem.

\beT \la{Bino_mix_Massart}  Let $\vep_a$ and $\vep_r$ be positive
numbers such that $0 < \vep_a < \f{3}{8}$ and {\small $\f{6
\vep_a}{3 - 2 \vep_a } < \vep_r < 1$}.   Suppose the sample size for
the $s$-th stage is no less than {\small $\li \lc  \f{ \ln (\ze \de)
} { \mscr{M} (\f{\vep_a}{\vep_r} + \vep_a, \f{\vep_a}{\vep_r} ) }
\ri \rc$}. Define {\small
\[ \bs{D}_\ell = \bec 0 & \mrm{for} \; \f{1}{2} - \f{2}{3} \vep_a -
\sq{ \f{1}{4} + \f{ n_\ell \vep_a^2 } {2 \ln (\ze \de) } } <
\wh{\bs{p}}_\ell < \f{ 6(1 - \vep_r) (3 - \vep_r) \ln (\ze
\de) } { 2 (3 - \vep_r)^2 \ln (\ze \de) - 9 n_\ell \vep_r^2} \; \mrm{or}\\
  &  \qu \;\; \f{1}{2} + \f{2}{3} \vep_a - \sq{ \f{1}{4} + \f{
n_\ell \vep_a^2 } {2 \ln (\ze \de) } } < \wh{\bs{p}}_\ell < \f{ 6(1
+ \vep_r) (3 + \vep_r) \ln (\ze \de) } { 2 (3 +
\vep_r)^2 \ln (\ze \de) - 9 n_\ell  \vep_r^2},\\
1 & \mrm{else}
 \eec
\]}
for $\ell = 1, \cd, s$.  Then, \bee &    & \Pr \{  p \leq \mscr{L} (
\wh{\bs{p}} ) \mid p \} \leq \sum_{\ell = 1}^s  \Pr \{  p \leq
\mscr{L} ( \wh{\bs{p}}_\ell ), \;
\bs{D}_\ell = 1 \mid p \} \leq s \ze \de,\\
&    & \Pr \{  p \geq \mscr{U} ( \wh{\bs{p}} ) \mid p \} \leq
\sum_{\ell = 1}^s  \Pr \{  p \geq \mscr{U} ( \wh{\bs{p}}_\ell ), \;
\bs{D}_\ell = 1 \mid p \} \leq s \ze \de \eee and $\Pr \{  |
\wh{\bs{p}} - p | < \vep_a \; \tx{or} \;  | \wh{\bs{p}} - p | <
\vep_r p \mid p \} \geq 1 - 2 s \ze \de $ for any $p \in (0, 1)$.
\eeT

See Appendix \ref{App_Bino_mix_Massart} for a proof.  Based on the
criteria proposed in Section \ref{gen_structure}, the sample sizes
$n_1 < n_2 < \cd < n_s$ can be chosen as the ascending arrangement
of all distinct elements of {\small
\[
\li \{ \li \lc 2 C_{\tau - \ell} \; \li ( \f{1}{\vep_a} -
\f{1}{\vep_r}  - \f{1}{3}  \ri ) \li ( \f{1}{\vep_r} + \f{1}{3} \ri
) \ln \f{1} { \ze \de } \ri \rc: \ell = 1, \cd, \tau \ri \},
 \]}
where $\tau$ is the maximum integer such that $2 C_{\tau - \ell} \;
\li ( \f{1}{\vep_a} - \f{1}{\vep_r}  - \f{1}{3}  \ri ) \li (
\f{1}{\vep_r} + \f{1}{3} \ri ) \ln \f{1} { \ze \de }  \geq \f{ 4(3 +
\vep_r) } {9 \vep_r } \ln \f{1} { \ze \de }$, i.e., $C_{\tau - 1}
\geq \f{2}{3} \li ( \f{1}{\vep_a} - \f{1}{\vep_r}  - \f{1}{3}  \ri
)^{- 1}$.  For such a choice of sample sizes, as a result of Theorem
\ref{Bino_mix_Massart}, we have that $\Pr \{ | \wh{\bs{p}} - p | <
\vep_a \; \tx{or} \;  | \wh{\bs{p}} - p | < \vep_r p \mid p \} > 1 -
\de $ for any $p \in (0, 1)$ provided that $\ze < \f{1}{2 \tau}$.

\subsubsection{Asymptotic Stopping Rule}

It should be noted that, for small $\vep_a$ and $\vep_r$, we can
simplify, by using Taylor's series expansion formula $\ln (1 + x) =
x - \f{x^2}{2} + o (x^2)$, the sampling schemes described in Section
\ref{multiinv_mix} as follows:

(i) The sequence of sample sizes $n_1, \cd, n_s$ is defined as the
ascending arrangement of all distinct elements of {\small $ \li \{
\li \lc 2 C_{\tau - \ell} \li ( \f{1}{\vep_a} - \f{1}{\vep_r} \ri )
\f{ \ln \f{1}{\ze \de} } { \vep_r } \ri \rc : \ell = 1, \cd, \tau
\ri \}$} with $\vep_a < \f{\vep_r}{2}$, where $\tau$ is the maximum
integer such that $C_{\tau - 1} \geq \li ( \f{2}{\vep_a} -
\f{2}{\vep_r} \ri )^{-1}$.

(ii) The decision variables are defined such that $\bs{D}_\ell = 1$
if \be \la{mixnor} n_\ell \geq \f{ \wh{\bs{p}}_\ell (1 -
\wh{\bs{p}}_\ell) \; 2 \ln \f{1}{\ze \de}  }{\max \{ \vep_a^2, \; (
\vep_r \wh{\bs{p}}_\ell )^2 \} }; \ee and $\bs{D}_\ell = 0$
otherwise.

\bsk

For such a simplified sampling scheme, we have \bel \sum_{\ell =
1}^s \Pr \li \{ | \wh{\bs{p}}_\ell - p | \geq \max \{\vep_a, \vep_r
p \}, \;  \bs{D}_\ell = 1 \ri \} & \leq & \sum_{\ell = 1}^s \Pr \li
\{ |
\wh{\bs{p}}_\ell - p | \geq \max \{\vep_a, \vep_r p \} \ri \} \nonumber\\
& \leq & \sum_{\ell = 1}^\tau \Pr \li \{ | \wh{\bs{p}}_\ell - p |
\geq \max \{\vep_a, \vep_r p \} \ri \} \nonumber\\
& \leq &  \sum_{\ell = 1}^\tau 2 \exp \li ( n_\ell \mscr{M}_{\mrm{B}}
\li (\f{\vep_a}{\vep_r} + \vep_a, \f{\vep_a}{\vep_r} \ri )  \ri ) \la{MineCH}\\
& < & 2 \tau \exp \li ( n_1 \mscr{M}_{\mrm{B}} \li
(\f{\vep_a}{\vep_r} + \vep_a, \f{\vep_a}{\vep_r} \ri ) \ri ),
\la{MLastB} \eel where (\ref{MineCH}) is due to Theorem 1 of
\cite{Chen0}. As can be seen from (\ref{MLastB}), the last bound is
independent of $p$ and can be made smaller than $\de$ if $\ze$ is
sufficiently small.  This establishes the claim and it follows that
{\small $\Pr \li \{ \li | \wh{\bs{p}} - p
 \ri | < \vep_a  \; \mrm{or} \; \li | \f{\wh{\bs{p}} - p } {p }
 \ri | < \vep_r \mid p \ri \} > 1 - \de$} for any $p \in (0, 1)$ if $\ze$ is sufficiently
small.

To improve the performance of coverage probability, we propose to
modify (\ref{mixnor}) to produce a new stopping rule as follows:

For $\ell = 1, \cd, s$, decision variable $\bs{D}_\ell$ assumes
value $1$ if \be \la{mixgood}
 \li ( \li | \wh{\bs{p}}_\ell -
\f{1}{2} \ri | - w \bs{\ep}_\ell \ri )^2 \geq \f{1}{4} + \f{
\bs{\ep}_\ell^2 n_\ell } { 2 \ln (\ze \de) } \ee  is satisfied; and
assumes value $0$ otherwise, where $\bs{\ep}_\ell = \max \{ \vep_a,
\; \vep_r \wh{\bs{p}}_\ell \}$ and $w \geq 0$ is a parameter
affecting the shape of the stopping boundary.

Before concluding this subsection, we would like to point out that
it is also possible to modify (\ref{mixgood}) to produce the
following stopping rule: For $\ell = 1, \cd, s$, decision variable
$\bs{D}_\ell$ assumes value $1$ if \be \la{mixsim} \li (
\wh{\bs{p}}_\ell - \f{1}{2} \ri )^2 \geq \f{1}{4} + \f{w}{n_\ell} +
\f{ \bs{\ep}_\ell^2 n_\ell } { 2 \ln (\ze \de) } \ee is satisfied;
and assumes value $0$ otherwise, where $\bs{\ep}_\ell = \max \{
\vep_a, \; \vep_r \wh{\bs{p}}_\ell \}$ and $w \geq 0$ is a parameter
affecting the shape of the stopping boundary.

The sample sizes for the stopping rules associated with
(\ref{mixgood}) and (\ref{mixsim}) can be chosen as suggested in
Section 2.1.  Specifically,  the maximum sample size $n_s$, i.e.,
the sample size of the last stage, should be defined as the smallest
integer such that $\{ \bs{D}_s = 1 \}$ is a sure event. The minimum
sample size $n_1$, i.e., the sample size of the first stage, should
be defined as the smallest integer such that $\{ \bs{D}_1 = 1 \}$ is
an event of a positive probability.  For both stopping rules, we can
show that {\small $\Pr \li \{ \li | \wh{\bs{p}} - p
 \ri | < \vep_a  \; \mrm{or} \; \li | \f{\wh{\bs{p}} - p } {p }
 \ri | < \vep_r \mid p \ri \} > 1 - \de$} for any $p \in (0, 1)$ if $\ze$ is sufficiently
small.

\subsubsection{Asymptotic Analysis of Sampling Schemes}

In this subsection, we shall focus on the asymptotic analysis of
multistage inverse sampling schemes.  Throughout this subsection, we
assume that the multistage sampling schemes follow stopping rules
derived from Chernoff bounds as described in Section
\ref{multiinv_mix}. Moreover, we assume that the sample sizes $n_1,
\cd, n_s$ are chosen as the ascending arrangement of all distinct
elements of  the set defined by (\ref{Bino_mix_ss}).

With regard to the tightness of the double-decision-variable method,
we have

\beT  \la{Bino_mix_DDV_Asp}  Let $\mscr{R}$ be a subset of real
numbers. Define {\small \[ \ovl{P} = \sum_{\ell = 1}^s \Pr \{
\wh{\bs{p}}_\ell \in \mscr{R}, \; \bs{D}_{\ell - 1} = 0, \;
\bs{D}_\ell = 1 \}, \qqu \udl{P} = 1 - \sum_{\ell = 1}^s \Pr \{
\wh{\bs{p}}_\ell \notin \mscr{R}, \; \bs{D}_{\ell - 1} = 0, \;
\bs{D}_\ell = 1 \}.
\]}
Then, $\udl{P} \leq \Pr \{ \wh{\bs{p}} \in \mscr{R} \} \leq \ovl{P}$
and $\lim_{\vep_a \to 0} \li | \Pr \{ \wh{\bs{p}} \in \mscr{R} \} -
\ovl{P} \ri | = \lim_{\vep_a \to 0} \li | \Pr \{ \wh{\bs{p}} \in
\mscr{R} \} - \udl{P} \ri | = 0$ for any $p \in (0, 1)$, where the
limits are taken under the constraint that $\f{\vep_a}{\vep_r}$ is
fixed.

\eeT

See Appendix \ref{App_Bino_mix_DDV_Asp} for a proof.

\bsk

With regard to the asymptotic performance of the sampling scheme as
$\vep_a$ and $\vep_r$ tend to $0$, we have

\beT \la{Bino_mix_Analysis} Let $\mcal{N}_{\mrm{f}} (p, \vep_a,
\vep_r)$ be the minimum sample number $n$ such that {\small
\[
\Pr  \li \{  \li | \f{\sum_{i = 1}^n X_i}{n} - p
 \ri | < \vep_a \; \tx{or} \;  \li | \f{\sum_{i = 1}^n X_i}{n} - p
 \ri | < \vep_r p  \mid p \ri \} > 1 - \ze \de
 \] }
 for a fixed-size sampling procedure. Let {\small $\mcal{N}_{\mrm{m}} (p, \vep_a, \vep_r ) = \f{ \ln (
\ze \de  ) } { \max \{ \mscr{M}_{\mrm{B}} (p, \udl{p} ), \;
\mscr{M}_{\mrm{B}} (p, \ovl{p} ) \}  } $}, where $\udl{p} = \min \{
p - \vep_a, \; \f{p}{1 + \vep_r} \}$ and $\ovl{p}= \max \{ p +
\vep_a, \; \f{p}{1 - \vep_r} \}$.  Define $p^\star =
\f{\vep_a}{\vep_r}, \; d = \sq{ 2 \ln \f{1}{\ze \de} }$,
\[
r(p) = \bec \f{p (1 - p)} {  p^\star ( 1 - p^\star)  } & for
\; p \in (0, p^\star],\\

\f{ p^\star (1 - p)} { p ( 1 - p^\star) }  & for \; p \in (p^\star, 1) \eec \qqu \nu = \bec \f{2}{3} - \f{1}{3} \f{ p (1 - p) (1 - 2 p^\star)} {
p^\star (1 - p^\star) (1 - 2 p)} & \tx{for} \; p \in (0, p^\star],\\
1 - \f{1}{3} \f{p - p^\star}{1 - p^\star} & \tx{for} \; p \in (p^\star, 1). \eec
\]
Let {\small $\ka_p = \f{C_{j_p}}{r (p)} $}, where $j_p$ is the
maximum integer $j$ such that $C_j \geq r(p)$.  Let {\small $\ro_p =
\f{C_{j_p - 1}} { r (p) } - 1$} if $\ka_p = 1, \; j_p > 0$ and
$\ro_p = \ka_p - 1$ otherwise. The following statements hold true
under the condition that $\f{\vep_a}{\vep_r}$ is fixed.

(I): {\small $\Pr \li \{  1 \leq \limsup_{\vep_a \to 0} \f{ \mbf{n}
} { \mcal{N}_{\mrm{m}} (p, \vep_a, \vep_r) } \leq 1 + \ro_p \ri \} =
1$}.  Specially, {\small $\Pr \li \{ \lim_{\vep_a \to 0} \f{ \mbf{n}
} { \mcal{N}_{\mrm{m}} (p, \vep_a, \vep_r) } = \ka_p \ri \} = 1$} if
$\ka_p
> 1$.

(II): $\lim_{\vep_a \to 0} \f{ \bb{E} [ \mbf{n} ] } {
\mcal{N}_{\mrm{f}} (p, \vep_a, \vep_r)} = \li ( \f{ d } {
\mcal{Z}_{\ze \de} } \ri )^2 \times \lim_{\vep_a \to 0} \f{ \bb{E} [
\mbf{n} ] } { \mcal{N}_{\mrm{m}} (p, \vep_a, \vep_r) }$,  where
\[
\lim_{\vep_a \to 0} \f{ \bb{E} [ \mbf{n} ] } { \mcal{N}_{\mrm{m}}
(p, \vep_a, \vep_r) } = \bec
\ka_p & \tx{if} \; \ka_p > 1,\\
1 + \ro_p \Phi (\nu d ) &  \tx{otherwise} \eec
\]
and $1 \leq \lim_{\vep_a \to 0} \f{ \bb{E} [ \mbf{n} ] } {
\mcal{N}_{\mrm{m}} (p, \vep_a, \vep_r) } \leq 1 + \ro_p$.

(III): If $\ka_p > 1$, then $\lim_{\vep_a \to 0} \Pr \{ | \wh{\bs{p}} - p | < \vep_a \; \tx{or} \; | \wh{\bs{p}} - p | < \vep_r p \} = 2 \Phi
\li ( d \sq{\ka_p} \ri ) - 1 > 2 \Phi \li ( d \ri ) - 1 > 1 - 2 \ze \de$. Otherwise, $\Phi \li ( d  \ri )  + \Phi \li ( d \sq{1 + \ro_p}  \ri )
- 1 >  \lim_{\vep_a \to 0} \Pr \{ | \wh{\bs{p}} - p | < \vep_a \; \tx{or} \; | \wh{\bs{p}} - p | < \vep_r p   \} = 1 + \Phi(d) - \Phi(\nu d) -
\Psi (\ro_p, \nu, d) > \Phi \li ( d \ri ) + 2 \Phi \li ( d \sq{1 + \ro_p}  \ri )  - 2  > 1 - 3 \ze \de$.

\eeT

See Appendix \ref{App_Bino_mix_Analysis} for a proof.

\sect{Estimation of Functions of Two Binomial Proportions}

Estimation of functions of two binomial proportions is particularly important in prospective comparative studies such as randomized controlled
clinical trial. More formally, let $X$ and $Y$ be independent Bernoulli random variables such that $\Pr \{ X = 1 \} = 1 - \Pr \{ X = 0 \} = p_x
\in (0, 1)$ and $\Pr \{ Y = 1 \} = 1 - \Pr \{ Y = 0 \} = p_y \in (0, 1)$.  Let $g(.,.)$ be a bivariate function of $p_x$ and $p_y$. It is a
frequent problem to estimate $g(p_x, p_y)$ based on samples of $X$ and $Y$. Typical examples of $g(p_x, p_y)$  are $g(p_x, p_y) = p_x - p_y$ and
$g(p_x, p_y) = \f{p_x}{p_y}$, which are respectively referred to as the {\it difference of population proportions} and {\it ratio of population
proportions}.  In this section, our goal is to develop exact methods for estimating $g(p_x, p_y)$. We will first consider the problem of
estimating $g(p_x, p_y)$ based on fixed-sample-size methods.  Afterward, we will discuss multistage estimation of $g(p_x, p_y)$.   For
simplicity of notations, we denote the tuple $(p_x, p_y)$ by $\bs{p}$.  The set of all possible values of $\bs{p}$ is denoted  by $\Se$. Namely,
$\Se = \{ (p_x, p_y): 0 < p_x < 1, \; 0 < p_y < 1 \}$. The bivariate function $g(p_x, p_y)$ is abbreviated as $g (\bs{p})$.

\subsection{Single-stage Estimation}

Let $X_i, \; i = 1, \cd, N_x$ be i.i.d. samples of $X$.  Let $Y_i, \; i = 1, \cd, N_y$ be i.i.d. samples of $Y$.  Assume that the samples of $X$
and $Y$ are independent. Let $K_x = \sum_{i = 1}^{N_x} X_i$ and $K_y = \sum_{i = 1}^{N_y} Y_i$. Let $\wh{p}_x = \f{K_x}{N_x}, \; \wh{p}_y =
\f{K_y}{N_y}$ and $\wh{\bs{p}} = (\wh{p}_x, \wh{p}_y)$.  A general problem for estimating $g(\bs{p})$ is to construct a confidence interval
$(\mcal{L} ( \wh{\bs{p}} ), \; \mcal{U} ( \wh{\bs{p}} ) )$ for $g(\bs{p})$ such that
\[
\Pr \{ \mcal{L} ( \wh{\bs{p}} ) <  g (\bs{p}) < \mcal{U} ( \wh{\bs{p}} )  \mid \bs{p} \} \geq 1 - \de
\]
for all $\bs{p} \in \Se$, where $\de \in (0, 1)$ is a pre-specified confidence parameter. We propose to solve this problem by virtue of the
coverage tuning technique. Our main idea is as follows:

(i) Seek a class of confidence intervals $ ( \mcal{L} ( \wh{\bs{p}} ), \; \mcal{U} ( \wh{\bs{p}} ) )$ such that the coverage probability can be
controlled by the coverage tuning parameter $\ze > 0$. In other words, the coverage probability, denoted by $P(\ze, \bs{p}) \DEF \Pr \{ \mcal{L}
( \wh{\bs{p}} ) <  g (\bs{p}) < \mcal{U} ( \wh{\bs{p}} )  \mid \bs{p} \}$,  tends to $1$ for all $\bs{p} \in \Se$ as $\ze \to 0$.

(ii)  For a given value of coverage tuning parameter $\ze$, apply Adapted Branch and Bound method in Appendix \ref{BBA} to determine whether the
complementary coverage probability $\Psi(\ze, \bs{p}) \DEF 1 - P(\ze, \bs{p})$ of the confidence interval associated with $\ze$ is no greater
than $\de$ for all $\bs{p} \in \Se$.

(iii) Apply bisection coverage tuning method to determine $\ze > 0$ as large as possible such that $P (\ze, \bs{p}) \geq 1 - \de$ for all
$\bs{p} \in \Se$.

Actually, there are many methods to construct confidence intervals satisfying requirement (i) (see, e.g., \cite{Brown, Newcombe} and the
references therein). As an illustration, consider the confidence interval for $g(\bs{p}) = p_x - p_y$ investigated in \cite{Newcombe}, which is
of the form: \bel &  & \mcal{L} ( \wh{\bs{p}} ) = \wh{p}_x - \wh{p}_y - \mcal{Z}_{\ze \de \sh 2} \sq{
\f{ l_x (1 - l_x) }{N_x} +  \f{ u_y (1 - u_y) }{N_y} }, \la{good88}\\
&  & \mcal{U} ( \wh{\bs{p}} ) = \wh{p}_x - \wh{p}_y + \mcal{Z}_{\ze \de \sh 2} \sq{  \f{ u_x (1 - u_x) }{N_x} +  \f{ l_y (1 - l_y) }{N_y} },
\la{good89} \eel where $l_x, \; u_x$ are the roots for $p$ in the quadratic equation $| \wh{p}_x - p | = \mcal{Z}_{\ze \de \sh 2} \sq{ p (1 - p)
\sh N_x }$; and  $l_y, \; u_y$ are the roots for $p$ in the quadratic equation $| \wh{p}_y - p | = \mcal{Z}_{\ze \de \sh 2} \sq{ p (1 - p) \sh
N_y}$.  Solving these equations yields \bee &  & l_x = \f{ c + 2 N_x \wh{p}_x - \sq{ c^2 + 4 c N_x \wh{p}_x (1 - \wh{p}_x)    }  } { 2 (c + N_x)
}, \qu u_x = \f{ c +
2 N_x \wh{p}_x + \sq{ c^2 + 4 c N_x \wh{p}_x (1 - \wh{p}_x) }  } { 2 (c + N_x) }, \\
&  & l_y = \f{ c + 2 N_y \wh{p}_y - \sq{ c^2 + 4 c N_y \wh{p}_y (1 - \wh{p}_y)    }  } { 2 (c + N_y) }, \qu u_y = \f{ c + 2 N_y \wh{p}_y + \sq{
c^2 + 4 c N_y \wh{p}_y (1 - \wh{p}_y) }  } { 2 (c + N_y) }, \eee where $c = \mcal{Z}_{\ze \de \sh 2}^2$.  Substituting these roots into
(\ref{good88}) and (\ref{good89}) leads to {\small \bee \mcal{L} ( \wh{\bs{p}} )  = \wh{p}_x - \wh{p}_y - \sq{ \li [ \f{ c (1 - 2 \wh{p}_x) -
\sq{ c^2 + 4 c N_x \wh{p}_x (1 - \wh{p}_x) } } { 2 (c + N_x) } \ri ]^2 + \li [ \f{ c (1 - 2 \wh{p}_y) + \sq{ c^2 + 4 c N_y \wh{p}_y (1 -
\wh{p}_y) }  } { 2 (c + N_y) } \ri ]^2
}, &  & \\
\mcal{U} ( \wh{\bs{p}} )  = \wh{p}_x - \wh{p}_y + \sq{ \li [ \f{ c (1 - 2 \wh{p}_x) + \sq{ c^2 + 4 c N_x \wh{p}_x (1 - \wh{p}_x) }  } { 2 (c +
N_x) } \ri ]^2 + \li [ \f{ c (1 - 2 \wh{p}_y) - \sq{ c^2 + 4 c N_y \wh{p}_y (1 - \wh{p}_y) }  } { 2 (c + N_y) } \ri ]^2 }. &  & \eee} Clearly,
the coverage probability of such confidence interval can be controlled by $\ze$, i.e., $P (\ze, \bs{p}) \to 1$ for all $\bs{p} \in \Se$ as $ \ze
\to 0$.

In general, a confidence interval for $g(\bs{p}) = p_x - p_y$ can be constructed as follows.  Let $(L_x, U_x)$ and $(L_y, U_y)$ be confidence
intervals for $p_x$ and $p_y$ respectively such that $\Pr \{ L_x < p_x < U_x \mid p_x \} \geq 1 - \ze \de$ and  $\Pr \{ L_y < p_y < U_y \mid p_y
\} \geq 1 - \ze \de$.  Taking $\mcal{L} ( \wh{\bs{p}} ) = L_x - U_y$ and $\mcal{U} ( \wh{\bs{p}} ) = U_x - L_y$ as the lower and upper
confidence limits for $g(\bs{p}) = p_x - p_y$ leads to $P(\ze, \bs{p}) \geq 1 - 2 \ze \de$.  Similarly, the lower and upper confidence limits
for $g(\bs{p}) = \f{p_x}{p_y}$ can be respectively taken as $\mcal{L} ( \wh{\bs{p}} ) = \f{L_x}{U_y}$ and $\mcal{U} ( \wh{\bs{p}} ) =
\f{U_x}{L_y}$, which ensures that the coverage probability $P(\ze, \bs{p}) \geq 1 - 2 \ze \de$.  Of course, one can use confidence intervals
$(L_x, U_x)$ and $(L_y, U_y)$ with confidence levels approximately equal to $1 - \ze \de$ as long as the coverage probability of the resultant
confidence interval for $g(\bs{p})$ can be controlled by $\ze$.

Given that the structure of the confidence interval is determined so that the coverage probability $P(\ze, \bs{p})$ can be controlled by $\ze$,
we can apply the bisection coverage tuning method to obtain $\ze > 0$ as large as possible such that $P (\ze, \bs{p}) \geq 1 - \de$ for all
$\bs{p} \in \Se$. A critical step for coverage tuning is to determine whether a given $\ze > 0$ is small enough to ensure that the complementary
coverage probability $\Psi (\ze, \bs{p})$ of the confidence interval associated with $\ze$ is no greater than $\de$ for all $\bs{p} \in \Se$. We
propose to apply Adapted Branch and Bound algorithms in Appendix \ref{BBA} to accomplish this task. This needs readily computable bounds of
$\Psi(\ze, \bs{p})$ for $\bs{p}$ in a rectangular domain $\mcal{Q} \DEF \{ (p_x, p_y) : 0 \leq \udl{p}_x \leq p_x \leq \ovl{p}_x \leq 1, \; 0
\leq \udl{p}_y \leq p_y \leq \ovl{p}_y \leq 1 \}$, which will be established in the sequel.

Let $\udl{g}$ and $\ovl{g}$ be lower and upper bounds of $g(\bs{p})$ such that $\udl{g} \leq g(\bs{p}) \leq \ovl{g}$ for all $\bs{p} \in
\mcal{Q}$ and that $\ovl{g} - \udl{g} \to 0$ as $\ovl{p}_x - \udl{p}_x \to 0, \; \ovl{p}_y - \udl{p}_y \to 0$.  Specially, we can take $\udl{g}
= \udl{p}_x - \ovl{p}_y, \; \ovl{g} = \ovl{p}_x - \udl{p}_y$ as the lower and upper bounds for $g(\bs{p}) = p_x - p_y$, and $\udl{g} = \udl{p}_x
\sh \ovl{p}_y, \; \ovl{g} = \ovl{p}_x \sh \udl{p}_y$ as the lower and upper bounds for $g(\bs{p}) = \f{p_x}{p_y}$. By virtue of the bounds of
$g(\bs{p})$, we have \be \la{bb1} \Pr \{ \mcal{L} ( \wh{\bs{p}} )  \geq \ovl{g}  \; \tx{or} \;
 \mcal{U} ( \wh{\bs{p}} ) \leq \udl{g} \mid \bs{p} \} \leq \Psi (\ze, \bs{p}) \leq \Pr \{ \mcal{L} ( \wh{\bs{p}} )  \geq \udl{g} \; \tx{or} \;
  \mcal{U} ( \wh{\bs{p}} ) \leq \ovl{g} \mid \bs{p} \} \ee  for all $\bs{p} \in \mcal{Q}$.  Moreover, the lower and upper bounds
of $\Psi (\ze, \bs{p})$ in (\ref{bb1}) converge as $\ovl{p}_x - \udl{p}_x \to 0, \; \ovl{p}_y - \udl{p}_y \to 0$.   To reduce the computational
complexity for evaluating the bounds of $\Psi(\ze, \bs{p})$, we can apply the truncation technique established in \cite{Chen1}. Specifically,
let $\eta \in (0,1)$ and define {\small \bee &  & a_x  = \mrm{T}_{\mrm{lb}} (\udl{p}_x, N_x, \eta),  \qqu b_x = \mrm{T}_{\mrm{ub}} (\ovl{p}_x,
N_x, \eta),
\qqu a_y  = \mrm{T}_{\mrm{lb}} (\udl{p}_y, N_y, \eta),  \qqu  b_y = \mrm{T}_{\mrm{ub}} (\ovl{p}_y, N_y, \eta), \\
&  & c_x  = \mrm{T}_{\mrm{lb}} (\ovl{p}_x, N_x, \eta), \qqu  d_x =  \mrm{T}_{\mrm{ub}} (\udl{p}_x, N_x, \eta), \qqu c_y  = \mrm{T}_{\mrm{lb}}
(\ovl{p}_y, N_y, \eta), \qqu  d_y =  \mrm{T}_{\mrm{ub}} (\udl{p}_y, N_y, \eta),  \eee} where $\mrm{T}_{\mrm{lb}} (., ., .)$ and
$\mrm{T}_{\mrm{ub}} (., ., .)$ are multivariate functions such that {\small \bel &  & \mrm{T}_{\mrm{lb}} (\se, n, \eta) = \max \li \{0, \;
\f{1}{n} \li \lc n \se + \f{ 1 - 2 \se - \sq{ 1 + \f{18 n \se
(1-\se)}{ \ln \f{2}{\eta}} } } {\f{2}{3n} + \f{3}{ \ln \f{2}{\eta}}}  \ri \rc  \ri \}, \la{truna}\\
&  & \mrm{T}_{\mrm{ub}} (\se, n, \eta) = \min \li \{1, \;  \f{1}{n} \li \lf n \se + \f{ 1 - 2 \se + \sq{ 1 + \f{18 n \se (1-\se)}{ \ln
\f{2}{\eta}} } } {\f{2}{3n} + \f{3}{ \ln \f{2}{\eta}}}  \ri \rf  \ri \} \la{trunb}\eel} for $\se \in [0, 1], \; \eta \in (0, 1)$ and $n \in
\bb{N}$. By virtue of (\ref{bb1}), Theorem 3 of \cite{Chen1} and Bonferroni's inequality, we have \bel
 \Psi (\ze, \bs{p}) \leq 2 \eta + \Pr \{ a_x \leq \wh{p}_x \leq b_x, \; a_y \leq \wh{p}_y \leq b_y, \;
 \mcal{L} ( \wh{\bs{p}} ) \geq \udl{g} \; \tx{or} \; \mcal{U} ( \wh{\bs{p}} ) \leq \ovl{g} \mid
\bs{p} \},  \la{bb18} &  & \\
\Psi (\ze, \bs{p}) \geq \Pr \{ c_x \leq \wh{p}_x \leq d_x, \; c_y \leq \wh{p}_y \leq d_y, \; \mcal{L} ( \wh{\bs{p}} ) \geq \ovl{g} \; \tx{or} \;
\mcal{U} ( \wh{\bs{p}} ) \leq \udl{g} \mid \bs{p} \} \qqu \; \; \la{bb28} &  & \eel for all $\bs{p} \in \mcal{Q}$.  Note that the probabilistic
terms in the lower and upper bounds of $\Psi (\ze, \bs{p})$ given by (\ref{bb18}) and (\ref{bb28}) can be expressed as  summations of finite
number of terms of the form \be \la{termvip} \mscr{T} (a, b, \bs{p}) \DEF \sum_{k_x = a}^b \Pr \{ K_x = k_x \mid p_x \} \Pr \{ u(k_x) \leq K_y
\leq v(k_x) \mid p_y \}, \ee where $a$ and $b$ are integers such that $0 \leq a \leq b \leq N_x$, $u(k_x)$ and $v(k_x)$ are integer-valued
functions of integer $k_x$ such that $0 \leq u(k_x) \leq v(k_x) \leq N_y$, and \bee &  & \Pr \{ K_x = k_x \mid p_x \} = \bi{N_x}{ k_x }
p_x^{k_x} (1 - p_x)^{N_x - k_x}, \\
&  &  \Pr \{ u(k_x) \leq K_y \leq v(k_x) \mid p_y \} = \sum_{k_y = u(k_x)}^{v(k_x)}  \bi{N_y}{ k_y } p_y^{k_y} (1 - p_y)^{N_y - k_y}. \eee  We
propose to bound $\mscr{T} (a, b, \bs{p})$ for $\bs{p} \in \mcal{Q}$ as follows.  Define \bee & & \udl{\Up}_x (k_x) = \min \li (  \Pr \{ K_x =
k_x \mid \udl{p}_x \}, \; \Pr \{ K_x = k_x \mid
\ovl{p}_x \} \ri ),\\
&  & \ovl{\Up}_x (k_x) =  \bec \Pr \{ K_x = k_x \mid p_x^* \} & \tx{for} \; p_x^* \in [\udl{p}_x, \ovl{p}_x],\\
\Pr \{ K_x = k_x \mid \udl{p}_x \} & \tx{for} \; p_x^* < \udl{p}_x,\\
\Pr \{ K_x = k_x \mid \ovl{p}_x \} & \tx{for} \; p_x^* > \ovl{p}_x \eec \eee where $p_x^* = \f{k_x}{N_x}$.  Then, $\udl{\Up}_x (k_x) \leq \Pr \{
K_x = k_x \mid p_x \} \leq \ovl{\Up}_x (k_x)$ for all $p_x \in [\udl{p}_x, \ovl{p}_x]$.  Let $u$ and $v$ be abbreviations of $u(k_x)$ and
$v(k_x)$ respectively.  Define \bee  &  & \udl{\Up}_y (k_x) =  \min \li (  \Pr \{ u \leq K_y \leq v \mid \udl{p}_y \}, \; \Pr \{ u \leq K_y \leq
v \mid
\ovl{p}_y \} \ri ),\\
&  & \ovl{\Up}_y (k_x) = \bec \Pr \{ u \leq K_y \leq v \mid p_y^* \} & \tx{for} \; p_y^* \in [\udl{p}_y, \ovl{p}_y],\\
\Pr \{ u \leq K_y \leq v \mid \udl{p}_y \} & \tx{for} \; p_y^* < \udl{p}_y,\\
\Pr \{ u \leq K_y \leq v \mid \ovl{p}_y \} & \tx{for} \; p_y^* > \ovl{p}_y \eec
 \eee where
{\small $p_y^* = 1 - \li \{ 1 + \li [ \f{ v ! (N_y -v-1)! } { (u-1) ! (N_y -u)! } \ri ]^{\f{1}{v-u+1}} \ri \}^{-1}$}.  By differentiation, it
can be shown that  the derivative of $\Pr \{ u \leq K_y \leq v \mid p_y \}$ with respective to $p_y$ is positive for $p_y < p_y^*$ and negative
for $p_y > p_y^*$.  Hence, $\udl{\Up}_y (k_x) \leq \Pr \{ u(k_x) \leq K_y \leq v(k_x) \mid p_y \} \leq \ovl{\Up}_y (k_x)$ for all $p_y \in
[\udl{p}_y, \ovl{p}_y]$.   It follows that
\[
\sum_{k_x = a}^b \udl{\Up}_x (k_x) \; \udl{\Up}_y (k_x)  \leq \mscr{T} (a, b, \bs{p}) \leq \sum_{k_x = a}^b \ovl{\Up}_x (k_x) \; \ovl{\Up}_y
(k_x)
\]
for all $\bs{p} \in \mcal{Q}$.  Clearly, the lower and upper bounds of $\Psi(\ze, \bs{p})$ in (\ref{bb18}) and (\ref{bb28}) can be  respectively
obtained by summing the bounds of terms like $\mscr{T} (a, b, \bs{p})$.

Based on the lower and upper bounds of $\Psi(\ze, \bs{p})$ obtained by the above method,  we can employ Adapted Branch and Bound technique in
Appendix \ref{BBA} to test whether $\Psi (\ze, \bs{p})$ is no greater than $\de$ for all $\bs{p} \in \Se$.  Consequently, we can apply a
bisection search method to determine the coverage tuning parameter $\ze$ as large as possible such that the coverage probability $P(\ze,
\bs{p})$ of the confidence interval associated with $\ze$ is no less than $1 - \de$ for all $\bs{p} \in \Se$.

In the above discussion, we have been focusing on the interval estimation for $g(\bs{p})$  when the sample sizes $N_x$ and $N_y$ are given. In
many applications, it is important to determine appropriate sample sizes $N_x$ and $N_y$ such that the estimator $g ( \wh{\bs{p}} )$ for
$g(\bs{p})$ satisfies some pre-specified requirements of reliability. In general, the problem can be formulated as follows. Let $\de \in (0, 1)$
be a pre-specified confidence parameter.  Let the margin of error be $\ep \DEF \max \{  \vep_a, \; \vep_r | g(\bs{p})| \}$, where $\vep_a \in
[0, 1)$ and $\vep_r \in [0, 1)$. Such a margin of error can be reduced to the margin of absolute error and the margin of relative error by
taking $\vep_r = 0$ and $\vep_a = 0$ respectively.  For pre-specified confidence parameter $\de$ and margins of error $\ep = \max \{ \vep_a, \;
\vep_r |g(\bs{p})| \}$, a problem of practical importance is to determine sample sizes $N_x, \; N_y$ as small as possible such that $\Pr \{ | g
( \wh{\bs{p}} )  - g (\bs{p}) | < \ep  \mid \bs{p} \}  \geq 1 - \de$ for all $\bs{p}$ in a pre-specified subset of $\Se$. By virtue of the
identity (\ref{RI_Indentity}), we can express event $\{ | g ( \wh{\bs{p}} )  - g (\bs{p}) | < \ep \}$ as $\li \{ \mcal{L} (\wh{\bs{p}}) < g
(\bs{p}) < \mcal{U} (\wh{\bs{p}}) \ri \}$, where $\mcal{L} (\wh{\bs{p}})$ and $\mcal{U} (\wh{\bs{p}}) $ are some functions of $\wh{\bs{p}}$.
Hence, $\Pr \{ | g ( \wh{\bs{p}} ) - g (\bs{p}) | < \ep \mid \bs{p} \}  = \Pr \{ \mcal{L} (\wh{\bs{p}}) < g (\bs{p}) < \mcal{U}
(\wh{\bs{p}})\mid \bs{p} \}$. In practices, one can choose $N_x$ and $N_y$ to be decreasing functions of $\ze$ such that both $N_x$ and $N_y$
tends to infinity as $\ze \to 0$.  This implies that the coverage probability $\mcal{P}(\ze, \bs{p}) \DEF \Pr \{ \mcal{L} (\wh{\bs{p}}) < g
(\bs{p}) < \mcal{U} (\wh{\bs{p}})\mid \bs{p} \}$ tends to $1$ for all $\bs{p} \in \Se$ as $\ze$ tends to $0$. Hence, the sample size problem is
equivalent to finding the largest $\ze$ such that the coverage probability $\mcal{P}(\ze, \bs{p})$ of the confidence interval is no less than $1
- \de$ for all $\bs{p} \in \Se$.  We can apply the preceding interval estimation technique to solve this problem.

\subsection{Multistage Estimation}

In the last section, we have developed exact methods for estimating the function $g(\bs{p})$ based on fixed-size sampling.  In this section, we
shall discuss the multistage estimation of $g(\bs{p})$.  We consider three general estimation problems.  The first problem is to construct a
confidence interval for $g(\bs{p})$ based on a given multistage sampling scheme.  The second problem is to construct a multistage sampling
scheme and an estimator for $g(\bs{p})$ guaranteeing prescribed levels of precision and confidence.  The third problem is to construct a
multistage sampling scheme and a bounded-width confidence interval for $g(\bs{p})$.

More formally, the first problem can be described as follows.  Consider a multistage sampling scheme having $s$ stages with sample sizes $N_{1,
x} < N_{2, x} < \cd < N_{s, x}$ and $N_{1, y} < N_{2, y} < \cd < N_{s, y}$.  Let $X_1, X_2, \cd$ be i.i.d. samples of $X$ and $Y_1, Y_2, \cd$ be
i.i.d. samples of $Y$. As before, assume that the samples of $X$ and $Y$ are independent.  For $\ell = 1, \cd, s$, define \bee &  &  K_{\ell, x}
= \sum_{i = 1}^{N_{\ell, x}} X_i, \qqu K_{\ell, y} = \sum_{i = 1}^{N_{\ell, y}} Y_i, \qqu
\bs{K}_\ell = (K_{\ell, x}, K_{\ell, y}), \\
&  & \wh{\bs{p}}_{\ell, x} = \f{K_{\ell, x}}{N_{\ell, x}}, \qqu \qu \wh{\bs{p}}_{\ell, y} = \f{K_{\ell, y}}{N_{\ell, y}}, \qqu \qu
\wh{\bs{p}}_{\ell} = (\wh{\bs{p}}_{\ell, x}, \wh{\bs{p}}_{\ell, y}) \eee  and let $\mscr{B}_\ell$ be a subset of the support of
$\wh{\bs{p}}_\ell$. The sampling process is continued until $\wh{\bs{p}}_\ell \in \mscr{B}_\ell$ for some $\ell \in \{ 1, \cd, s \}$. Assume
that $\{ \wh{\bs{p}}_s \in \mscr{B}_s \}$ is a sure event. For such a given sampling scheme, the problem is to construct a confidence interval
$( \mscr{L}(\wh{\bs{p}}, \bs{l}), \mscr{U}(\wh{\bs{p}}, \bs{l}) )$ such that \be \la{pr89}
 \Pr \{ \mscr{L}(\wh{\bs{p}}, \bs{l}) <  g(\bs{p}) < \mscr{U}(\wh{\bs{p}}, \bs{l}) \mid \bs{p} \} \geq 1 - \de \ee
 for all $\bs{p} \in \Se$.  The lower and upper confidence limits are functions of $\wh{\bs{p}} =
\wh{\bs{p}}_{\bs{l}}$ and $\bs{l}$, where $\bs{l}$ is the index of stage at which the sampling process is terminated.  We propose a coverage
tuning approach to solve such interval estimation problem.  Our main idea is as follows.

(i) For $\ell = 1, \cd, s$, construct a confidence interval $( \mscr{L}(\wh{\bs{p}}_\ell, \ell), \mscr{U}(\wh{\bs{p}}_\ell, \ell) )$ for
$g(\bs{p})$  in terms of estimator $\wh{\bs{p}}_\ell$, sample sizes $N_{\ell, x}, \; N_{\ell, y}$ and the coverage tuning parameter $\ze
> 0$ such that the coverage probability $\Pr \{ \mscr{L}(\wh{\bs{p}}_\ell, \ell) < g(\bs{p}) <
\mscr{U}(\wh{\bs{p}}_\ell, \ell) \mid \bs{p} \}$ can be made arbitrarily close to $1$ if
$\ze > 0$ is sufficiently small.

(ii)  At the termination of the sampling process, $( \mscr{L}(\wh{\bs{p}}, \bs{l}), \mscr{U}(\wh{\bs{p}}, \bs{l}) )$ is taken as the desired
confidence interval for $g(\bs{p})$.

(iii) The bisection coverage tuning technique is applied to determine  $\ze > 0$ as large as possible such that (\ref{pr89}) is satisfied.

As an illustration of  the above interval estimation method, consider the estimation of $g(\bs{p}) = p_x - p_y$, which is the difference of two
binomial proportions.  For $\ell = 1, \cd, s$, the confidence interval $( \mscr{L}(\wh{\bs{p}}_\ell, \ell), \mscr{U}(\wh{\bs{p}}_\ell, \ell) )$
at the $\ell$-th stage for $p_x - p_y$ can be taken as the confidence interval $( \mcal{L} ( \wh{\bs{p}}), \mcal{U} (\wh{\bs{p}}) )$ defined by
(\ref{good88}) and (\ref{good89}) with the estimator $\wh{\bs{p}}$ identified as $\wh{\bs{p}}_\ell$, and the sample sizes $N_x, \; N_y$
identified as $N_{\ell, x}, \; N_{\ell, y}$ respectively. Then, $\Pr \{ \mscr{L}(\wh{\bs{p}}, \bs{l}) < p_x - p_y < \mscr{U}(\wh{\bs{p}},
\bs{l}) \mid \bs{p} \} \ap 1 - \ze \de \to 1$ for all $\bs{p} \in \Se$ as $\ze \to 0$.

For coverage tuning purpose, we need to bound the complementary coverage probability $1 - \Pr \{ \mscr{L}(\wh{\bs{p}}, \bs{l}) < g (\bs{p})  <
\mscr{U}(\wh{\bs{p}}, \bs{l}) \mid \bs{p} \}$ for $\bs{p}$ in a rectangular domain $\mcal{Q} \DEF \{ (p_x, p_y) : 0 \leq \udl{p}_x \leq p_x \leq
\ovl{p}_x \leq 1, \; 0 \leq \udl{p}_y \leq p_y \leq \ovl{p}_y \leq 1 \}$.  This issue will be addressed after we present the method for
constructing sampling schemes for the second and third problems.

Our second problem is to construct a multistage sampling scheme and an estimator for $g(\bs{p})$ to satisfy some prescribed requirements of
reliability. In Section 2.3, we have shown that such a general problem can be cast in the framework of constructing a multistage sampling scheme
and a random interval $( \mscr{L}(\wh{\bs{p}}, \bs{l}), \; \mscr{U}(\wh{\bs{p}}, \bs{l}) )$ such that $\Pr \{ \mscr{L}(\wh{\bs{p}}, \bs{l}) <
g(\bs{p}) < \mscr{U}(\wh{\bs{p}}, \bs{l}) \mid \bs{p} \} \geq 1 - \de$ for all $\bs{p} \in \Se$.  To design the stopping rule, we can apply the
inclusion principle proposed in Section 2.5. Specifically, we seek a class of sampling schemes associated with the coverage tuning parameter
$\ze$ satisfying the following requirements:

(i) The sample sizes $N_{1, x} (\ze)  < N_{2, x} (\ze) < \cd < N_{s, x} (\ze) $ and $N_{1, y} (\ze) < N_{2, y} (\ze) < \cd < N_{s, y} (\ze) $
are decreasing functions of $\ze$.  Namely, the sample sizes are increasing as $\ze$ decreases.  For simplicity of notations, we abbreviate
$N_{\ell, x} (\ze), \; N_{\ell, y} (\ze)$ respectively as $N_{\ell, x}, \; N_{\ell, y}$ for $\ell = 1, \cd, s$.

(ii) For $\ell = 1, \cd, s$, a confidence interval $( \mcal{L}(\wh{\bs{p}}_\ell, \ell), \mcal{U}(\wh{\bs{p}}_\ell, \ell) )$ for $g(\bs{p})$ can
be constructed in terms of estimators $\wh{\bs{p}}_{\ell, x}, \; \wh{\bs{p}}_{\ell, y}$, sample sizes $N_{\ell, x}, \; N_{\ell, y}$ and the
coverage tuning parameter $\ze > 0$ such that $\Pr \{ g(\bs{p}) \in ( \mcal{L}(\wh{\bs{p}}_\ell, \ell), \; \mcal{U}(\wh{\bs{p}}_\ell, \ell) ) \;
\tx{for} \; \ell = 1, \cd, s \mid \bs{p} \}$ can be made arbitrarily close to $1$ for any $\bs{p} \in \Se$ if $\ze
> 0$ is sufficiently small.

(iii) The sampling process is continued until confidence interval $( \mcal{L}(\wh{\bs{p}}_\ell, \ell), \mcal{U}(\wh{\bs{p}}_\ell, \ell) )$ is
included in $( \mscr{L}(\wh{\bs{p}}_\ell, \ell), \mscr{U}(\wh{\bs{p}}_\ell, \ell) )$ for some $\ell \in \{1, \cd, s \}$.  The desired sequential
random interval is taken as $( \mscr{L}(\wh{\bs{p}}, \bs{l}), \; \mscr{U}(\wh{\bs{p}}, \bs{l}) )$ with $\wh{\bs{p}} = \wh{\bs{p}}_{\bs{l}}$,
where $\bs{l}$ is the index of stage at the termination of the sampling process.

(iv) The sampling process is guaranteed to be terminated at or before the $s$-th stage.

\bsk

Recall that our third problem is to construct a multistage sampling scheme and a bounded-width sequential confidence interval for $g(\bs{p})$.
To design the stopping rule, we can apply the principle proposed in Sections 2 and 14. We seek a class of sampling schemes associated with the
coverage tuning parameter $\ze$ satisfying the following requirements:

(i) The sample sizes $N_{1, x} (\ze)  < N_{2, x} (\ze) < \cd < N_{s, x} (\ze) $ and $N_{1, y} (\ze) < N_{2, y} (\ze) < \cd < N_{s, y} (\ze) $
are decreasing functions of $\ze$.

(ii) For $\ell = 1, \cd, s$, a confidence interval $( \mscr{L}(\wh{\bs{p}}_\ell, \ell), \mscr{U}(\wh{\bs{p}}_\ell, \ell) )$ for $g(\bs{p})$ can
be constructed in terms of estimators $\wh{\bs{p}}_{\ell, x}, \; \wh{\bs{p}}_{\ell, y}$, sample sizes $N_{\ell, x}, \; N_{\ell, y}$ and the
coverage tuning parameter $\ze > 0$ such that  $\Pr \{ g(\bs{p})  \in ( \mscr{L}(\wh{\bs{p}}_\ell, \ell), \; \mscr{U}(\wh{\bs{p}}_\ell, \ell) )
\; \tx{for} \; \ell = 1, \cd, s \mid \bs{p} \}$ can be made arbitrarily close to $1$ for any $\bs{p} \in \Se$ if $\ze > 0$ is sufficiently
small.

(iii) The sampling process is continued until the width of interval {\small $( \mscr{L}(\wh{\bs{p}}_\ell, \ell), \mscr{U}(\wh{\bs{p}}_\ell,
\ell) )$} is not exceeding the prescribed length for some $\ell \in \{1, \cd, s \}$.  The desired sequential bounded-width confidence interval
is taken as $( \mscr{L}(\wh{\bs{p}}, \bs{l}), \; \mscr{U}(\wh{\bs{p}}, \bs{l}) )$ with $\wh{\bs{p}} = \wh{\bs{p}}_{\bs{l}}$, where $\bs{l}$ is
the index of stage at the termination of the sampling process.

(iv) The sampling process is guaranteed to be terminated at or before the $s$-th stage.

\bsk

Within both classes of sampling schemes for the second and third problems, we can apply the bisection coverage tuning technique to determine
$\ze > 0$ as large as possible such that $\Pr \{ \mscr{L}(\wh{\bs{p}}, \bs{l}) <  g(\bs{p}) < \mscr{U}(\wh{\bs{p}}, \bs{l}) \mid \bs{p} \} \geq
1 - \de$ for all $\bs{p} \in \Se$. In applications, one can take the sample sizes for $X$ as functions of the sample sizes for $Y$.  Then, for a
given $\ze$, the maximum sample sizes $N_{s,x}, \; N_{s,y}$ can be determined as the minimum sample sizes for $X$ and $Y$ in a full sequential
scheme such that the sampling process is guaranteed to be terminated.

For the three general problems we discussed above, a common approach is to seek a class of sampling schemes associated with $\ze$ such that the
coverage probability of the resultant sequential random interval can be made arbitrarily close to $1$ if $\ze$ is sufficiently small.
Afterward, the bisection coverage tuning technique is applied to determine $\ze$ as large as possible such that the coverage probability $P(\ze,
\bs{p}) \DEF \Pr \{ \mscr{L}(\wh{\bs{p}}, \bs{l}) <  g(\bs{p}) < \mscr{U}(\wh{\bs{p}}, \bs{l}) \mid \bs{p} \}$ is no less than the pre-specified
level $1 - \de$.  A critical subroutine of the coverage tuning method is to apply Adapted Branch and Bound algorithm proposed in Appendix
\ref{BBA} to determine whether the complementary coverage probability
\[
\Psi (\ze, \bs{p}) \DEF 1 - P(\ze, \bs{p}) = \Pr \{ \mscr{L}(\wh{\bs{p}}, \bs{l}) \geq g(\bs{p}) \; \tx{or} \;  \mscr{U}(\wh{\bs{p}}, \bs{l})
\leq g(\bs{p}) \mid \bs{p} \} \] of the sequential random interval  associated with $\ze$ is no greater than $\de$ for all $\bs{p} \in \Se$.
This requires bounding the complementary coverage probability $\Psi (\ze, \bs{p})$ for $\bs{p}$ in a rectangular domain $\mcal{Q} = \{ (p_x,
p_y) : 0 \leq \udl{p}_x \leq p_x \leq \ovl{p}_x \leq 1, \; 0 \leq \udl{p}_y \leq p_y \leq \ovl{p}_y \leq 1 \}$. This problem is addressed in the
sequel.

As described in Section 5.1, let $\udl{g}$ and $\ovl{g}$ be lower and upper bounds of $g(\bs{p})$ such that $\udl{g} \leq g(\bs{p}) \leq
\ovl{g}$ for all $\bs{p} \in \mcal{Q}$ and that $\ovl{g} - \udl{g} \to 0$ as $\ovl{p}_x - \udl{p}_x \to 0, \; \ovl{p}_y - \udl{p}_y \to 0$.
Then, \be \la{vip898}
 \Pr \{ \mscr{L}(\wh{\bs{p}}, \bs{l}) \geq \ovl{g} \; \tx{or} \;  \mscr{U}(\wh{\bs{p}}, \bs{l}) \leq \udl{g} \mid \bs{p} \}
\leq \Psi (\ze, \bs{p}) \leq \Pr \{ \mscr{L}(\wh{\bs{p}}, \bs{l}) \geq \udl{g} \; \tx{or} \;  \mscr{U}(\wh{\bs{p}}, \bs{l}) \leq \ovl{g} \mid
\bs{p} \} \ee for all $\bs{p} \in \mcal{Q}$.  Note that for a given value of the coverage tuning parameter $\ze$, the stopping rule can be
described as ``the sampling process is continued until $\wh{\bs{p}}_\ell \in \mscr{B}_\ell$ for some $\ell \in \{1, \cd, s\}$,'' where
$\mscr{B}_\ell$ is a subset of the support of $\wh{\bs{p}}_\ell$, i.e.,  $\mscr{B}_\ell \subseteq I_{\wh{\bs{p}}_\ell}$ for $\ell = 1, \cd, s$.
Recall that $I_Z$ denotes the support of random variable $Z$.  Let $\mscr{B}_\ell^c$ denote the complementary set of $\mscr{B}_\ell$ such that
$\mscr{B}_\ell^c \cap \mscr{B}_\ell = \emptyset$ and $\mscr{B}_\ell^c \cup \mscr{B}_\ell = I_{\wh{\bs{p}}_\ell}$ for $\ell = 1, \cd, s$.  By the
definition of the stopping rule, \bel &  & \{ \mscr{L}(\wh{\bs{p}}, \bs{l}) \geq g(\bs{p}) \; \tx{or} \;  \mscr{U}(\wh{\bs{p}}, \bs{l}) \leq g
(\bs{p}) \} \nonumber\\
&  & = \bigcup_{i = 1}^s \{ \wh{\bs{p}}_\ell \in \mscr{B}_\ell^c, \; \ell = 1, \cd, i - 1; \; \wh{\bs{p}}_i \in \mscr{B}_i; \;
\mscr{L}(\wh{\bs{p}}_i, i) \geq g(\bs{p}) \; \tx{or} \; \mscr{U}(\wh{\bs{p}}_i, i) \leq g(\bs{p}) \}. \la{hello} \eel To reduce the
computational complexity for bounding $\Psi (\ze, \bs{p})$, we propose to apply the truncation method established in \cite{Chen1}.  To this end,
let $\eta \in (0, 1)$ be the tolerance of truncation. Define {\small \bee & & a_{\ell, x} = \mrm{T}_{\mrm{lb}} (\udl{p}_x, N_{\ell, x}, \eta),
\qu b_{\ell, x} = \mrm{T}_{\mrm{ub}} (\ovl{p}_x, N_{\ell, x}, \eta),
\qu a_{\ell, y}  = \mrm{T}_{\mrm{lb}} (\udl{p}_y, N_{\ell, y}, \eta),  \qu  b_{\ell, y} = \mrm{T}_{\mrm{ub}} (\ovl{p}_y, N_{\ell, y}, \eta), \\
&  & c_{\ell, x}  = \mrm{T}_{\mrm{lb}} (\ovl{p}_x, N_{\ell, x}, \eta), \qu  d_{\ell, x} =  \mrm{T}_{\mrm{ub}} (\udl{p}_x, N_{\ell, x}, \eta),
\qu c_{\ell, y} = \mrm{T}_{\mrm{lb}} (\ovl{p}_y, N_{\ell, y}, \eta), \qu  d_{\ell, y} =  \mrm{T}_{\mrm{ub}} (\udl{p}_y, N_{\ell, y}, \eta),
\eee} for $\ell = 1, \cd, s$, where $\mrm{T}_{\mrm{lb}} (., ., .)$ and $\mrm{T}_{\mrm{ub}} (., ., .)$ are multivariate functions defined by
(\ref{truna}) and (\ref{trunb}). By virtue of (\ref{vip898}),  (\ref{hello}), Theorem 3 of \cite{Chen1}, and Bonferroni's inequality,  we have
{\small \bel \Psi (\ze, \bs{p}) & \leq & 2 s \eta + \Pr \{ a_{\ell, x} \leq \wh{\bs{p}}_{\ell, x} \leq b_{\ell, x}, \; a_{\ell, y} \leq
\wh{\bs{p}}_{\ell, y} \leq b_{\ell, y}, \; \ell = 1, \cd, s; \; \mscr{L}(\wh{\bs{p}}, \bs{l}) \geq \udl{g} \; \tx{or} \; \mscr{U}(\wh{\bs{p}},
\bs{l}) \leq \ovl{g} \mid \bs{p}
\} \nonumber\\
& = & 2 s \eta + \sum_{i = 1}^s \Pr \{ a_{\ell, x} \leq \wh{\bs{p}}_{\ell, x} \leq b_{\ell, x}, \; a_{\ell, y} \leq \wh{\bs{p}}_{\ell, y} \leq
b_{\ell, y}, \; \ell = 1, \cd, s; \; \nonumber \\
&   & \qqu  \wh{\bs{p}}_\ell \in \mscr{B}_\ell^c, \; \ell = 1, \cd, i - 1; \; \wh{\bs{p}}_i \in \mscr{B}_i; \;
\mscr{L}(\wh{\bs{p}}_i, i) \geq \udl{g} \; \tx{or} \; \mscr{U}(\wh{\bs{p}}_i, i) \leq \ovl{g} \mid \bs{p} \} \nonumber\\
& \leq & 2 s \eta + \sum_{i = 1}^s \Pr \{ a_{\ell, x} \leq \wh{\bs{p}}_{\ell, x} \leq b_{\ell, x}, \; a_{\ell, y} \leq \wh{\bs{p}}_{\ell, y}
\leq b_{\ell, y}, \; \ell = 1, \cd, i; \;  \nonumber\\
&   & \qqu  \wh{\bs{p}}_\ell \in \mscr{B}_\ell^c, \; \ell = 1, \cd, i - 1; \; \wh{\bs{p}}_i \in \mscr{B}_i; \; \mscr{L}(\wh{\bs{p}}_i, i) \geq
\udl{g} \; \tx{or} \; \mscr{U}(\wh{\bs{p}}_i, i) \leq \ovl{g} \mid \bs{p} \} \la{hi8} \eel} and {\small \bel \Psi (\ze, \bs{p}) & \geq & \sum_{i
= 1}^s \Pr \{ c_{\ell, x} \leq \wh{\bs{p}}_{\ell, x} \leq d_{\ell, x}, \; c_{\ell, y} \leq \wh{\bs{p}}_{\ell, y}
\leq d_{\ell, y}, \; \ell = 1, \cd, i; \;  \nonumber \\
&   & \qqu  \wh{\bs{p}}_\ell \in \mscr{B}_\ell^c, \; \ell = 1, \cd, i - 1; \; \wh{\bs{p}}_i \in \mscr{B}_i; \; \mscr{L}(\wh{\bs{p}}_i, i) \geq
\ovl{g} \; \tx{or} \; \mscr{U}(\wh{\bs{p}}_i, i) \leq \udl{g} \mid \bs{p} \} \la{hi9} \eel} for all $\bs{p} \in \mcal{Q}$.  It can be seen that
the bounds of $\Psi(\ze, \bs{p})$ in (\ref{hi8}) and (\ref{hi9}) can be expressed as summations of probabilistic terms like $\Pr \{  \bs{K}_i
\in \mscr{K}_i, \; 1 \leq i \leq \ell \mid \bs{p} \}, \; \ell = 1, \cd, s$, where $\mscr{K}_\ell$ does not depend on $\bs{p}$ and is a subset of
the support of $\bs{K}_\ell$ for $\ell = 1, \cd, s$.  Since $\mscr{K}_\ell$ is independent of $\bs{p}$, we can establish recursive technique for
bounding $\Pr \{  \bs{K}_i \in \mscr{K}_i, \; 1 \leq i \leq \ell \mid \bs{p} \}, \; \ell = 1, \cd, s$ over $\mcal{Q}$ as follows.

For $\ell = 1, \cd, s$, let the realization,  $(k_{\ell, x}, k_{\ell, y})$, of $\bs{K}_\ell$ be denoted by $k_\ell$.  As a consequence of the
independence of the samples of $X$ and $Y$, we have that \be \la{rec383}
 \Pr \{ \bs{K}_1 = k_1 \mid \bs{p} \} = \Pr \{ K_{1, x} = k_{1, x} \mid p_x \}  \Pr \{ K_{1, y} = k_{1, y} \mid p_y
\} \ee for any $k_1 \in \mscr{K}_1$ and that  \bel &  & \Pr \{  \bs{K}_i \in \mscr{K}_i, \; 1 \leq i < \ell
+ 1; \; \bs{K}_{\ell +1} = k_{\ell+1} \mid \bs{p} \} \nonumber \\
&  & = \sum_{k_\ell \in \mscr{K}_\ell} \Pr \{  \bs{K}_i \in \mscr{K}_i, \; 1 \leq i < \ell; \; \bs{K}_\ell = k_\ell \mid \bs{p} \} \nonumber \\
&  & \; \times \Pr \{ K_{\ell+1, x} - K_{\ell, x} = k_{\ell+1, x} - k_{\ell, x} \mid p_x \} \times \Pr \{ K_{\ell+1, y} - K_{\ell, y} =
k_{\ell+1, y} - k_{\ell, y} \mid p_y \}  \qqu \qqu \la{rec3388} \eel for $k_{\ell + 1} \in \mscr{K}_{\ell + 1}$ and $\ell = 1, \cd, s - 1$,
where {\small \bee & & \Pr \{ K_{1, x} = k_{1, x} \mid p_x \}  = \bi{N_{1, x}}{k_{1, x}} p_x^{k_{1, x}} (1 -
p_x)^{N_{1, x} - k_{1, x}},\\
&  & \Pr \{ K_{1, y} = k_{1, y} \mid p_y \} = \bi{N_{1, y}}{k_{1, y}} p_y^{k_{1, y}} (1 - p_y)^{N_{1, y} - k_{1, y}},\\
&  & \Pr \{ K_{\ell+1, x} - K_{\ell, x} = k_{\ell+1, x} - k_{\ell, x} \mid p_x \} = \bi{N_{\ell +1, x} - N_{\ell,x}}{ k_{\ell+1, x} - k_{\ell,
x} } p_x^{  k_{\ell+1, x} - k_{\ell, x} } (1 - p_x)^{ N_{\ell +1, x} - N_{\ell,x} - ( k_{\ell+1, x} - k_{\ell, x}  )  },\\
&  & \Pr \{ K_{\ell+1, y} - K_{\ell, y} = k_{\ell+1, y} - k_{\ell, y} \mid p_y \} = \bi{N_{\ell +1, y} - N_{\ell,y}}{ k_{\ell+1, y} - k_{\ell,
y} } p_y^{  k_{\ell+1, y} - k_{\ell, y} } (1 - p_y)^{ N_{\ell +1, y} - N_{\ell,y} - ( k_{\ell+1, y} - k_{\ell, y}  )  }. \eee } In the sequel,
we shall apply (\ref{rec383}) and (\ref{rec3388}) to develop recursively computable lower and upper bounds $\udl{P}_\ell (k_\ell)$ and
$\ovl{P}_\ell (k_\ell)$ for $\Pr \{  \bs{K}_i \in \mscr{K}_i, \; 1 \leq i < \ell; \; \bs{K}_\ell = k_\ell \mid \bs{p} \}$ such that, for all
$k_\ell \in \mscr{K}_\ell$, \be \la{req88}
 \udl{P}_\ell (k_\ell) \leq \Pr \{  \bs{K}_i \in \mscr{K}_i, \;
1 \leq i < \ell; \; \bs{K}_\ell = k_\ell \mid \bs{p} \} \leq \ovl{P}_\ell (k_\ell), \qqu \ell = 1, \cd, s \ee for all $\bs{p} \in \mcal{Q}$. For
this purpose, define {\small \bee &  & \udl{\Up}_{0,x} (k_{1, x}) = \min \li ( \Pr \{ K_{1, x} = k_{1, x} \mid
\udl{p}_x \}, \; \Pr \{ K_{1, x} = k_{1, x}  \mid \ovl{p}_x \} \ri ),\\
&  & \ovl{\Up}_{0,x} (k_{1, x}) =  \bec \Pr \{ K_{1, x}  = k_{1, x} \mid p_{0, x}^* \} & \tx{for} \; p_{0,x}^* \in [\udl{p}_x, \ovl{p}_x],\\
\Pr \{ K_{1, x} = k_{1, x} \mid \udl{p}_x \} & \tx{for} \; p_{0,x}^* < \udl{p}_x,\\
\Pr \{ K_{1, x} = k_{1, x} \mid \ovl{p}_x \} & \tx{for} \; p_{0,x}^* > \ovl{p}_x \eec
 \eee}
where $p_{0, x}^* = \f{k_{1,x} }{N_{1,x}}$.  Define {\small \bee
 \udl{\Up}_{\ell,x} (k_{\ell, x}, k_{\ell+1, x}) = \min \li (  \Pr \{ K_{\ell+1, x} - K_{\ell, x} = k_{\ell+1, x} - k_{\ell, x} \mid \udl{p}_x
\}, \; \Pr \{ K_{\ell+1, x} - K_{\ell, x} = k_{\ell+1, x} - k_{\ell, x} \mid \ovl{p}_x \} \ri ), &  &\\
\ovl{\Up}_{\ell,x} (k_{\ell, x}, k_{\ell+1, x}) =  \bec \Pr \{ K_{\ell+1, x} - K_{\ell, x} = k_{\ell+1, x} - k_{\ell, x} \mid p_{\ell, x}^* \} & \tx{for} \; p_{\ell,x}^* \in [\udl{p}_x, \ovl{p}_x],\\
\Pr \{ K_{\ell+1, x} - K_{\ell, x} = k_{\ell+1, x} - k_{\ell, x} \mid \udl{p}_x \} & \tx{for} \; p_{\ell,x}^* < \udl{p}_x,\\
\Pr \{ K_{\ell+1, x} - K_{\ell, x} = k_{\ell+1, x} - k_{\ell, x} \mid \ovl{p}_x \} & \tx{for} \; p_{\ell,x}^* > \ovl{p}_x \eec \qqu \qqu \qqu \qqu \qqu \; \;  &  & \\
 \eee} with $p_{\ell, x}^* = \f{k_{\ell + 1,x} - k_{\ell, x}}{N_{\ell+1,x} - N_{\ell, x} }$ for $\ell = 1, \cd, s - 1$.
 Define {\small \bee &  & \udl{\Up}_{0,y} (k_{1, y}) = \min \li (  \Pr \{ K_{1, y} = k_{1, y} \mid
\udl{p}_y \}, \; \Pr \{ K_{1, y} = k_{1, y}  \mid \ovl{p}_y \} \ri ),\\
&  & \ovl{\Up}_{0,y} (k_{1, y}) =  \bec \Pr \{ K_{1, y}  = k_{1, y} \mid p_{0, y}^* \} & \tx{for} \; p_{0,y}^* \in [\udl{p}_y, \ovl{p}_y],\\
\Pr \{ K_{1, y} = k_{1, y} \mid \udl{p}_y \} & \tx{for} \; p_{0,y}^* < \udl{p}_y,\\
\Pr \{ K_{1, y} = k_{1, y} \mid \ovl{p}_y \} & \tx{for} \; p_{0,y}^* > \ovl{p}_y \eec \eee} where $p_{0, y}^* = \f{k_{1,y} }{N_{1,y}}$.  Define
{\small \bee
 \udl{\Up}_{\ell,y} (k_{\ell, y}, k_{\ell+1, y}) = \min \li (  \Pr \{ K_{\ell+1, y} - K_{\ell, y} = k_{\ell+1, y} - k_{\ell, y} \mid \udl{p}_y
\}, \; \Pr \{ K_{\ell+1, y} - K_{\ell, y} = k_{\ell+1, y} - k_{\ell, y} \mid \ovl{p}_y \} \ri ),  &  &\\
\ovl{\Up}_{\ell,y} (k_{\ell, y}, k_{\ell+1, y}) =  \bec \Pr \{ K_{\ell+1, y} - K_{\ell, y} = k_{\ell+1, y} - k_{\ell, y} \mid p_{\ell, y}^* \} & \tx{for} \; p_{\ell,y}^* \in [\udl{p}_y, \ovl{p}_y],\\
\Pr \{ K_{\ell+1, y} - K_{\ell, y} = k_{\ell+1, y} - k_{\ell, y} \mid \udl{p}_y \} & \tx{for} \; p_{\ell,y}^* < \udl{p}_y,\\
\Pr \{ K_{\ell+1, y} - K_{\ell, y} = k_{\ell+1, y} - k_{\ell, y} \mid \ovl{p}_y \} & \tx{for} \; p_{\ell,y}^* > \ovl{p}_y \eec \qqu \qqu \qqu \qqu \qqu \; \;  &  & \\
 \eee} with $p_{\ell, y}^* = \f{k_{\ell + 1,y} - k_{\ell, y}}{N_{\ell+1,y} - N_{\ell, y} }$ for $\ell = 1, \cd, s - 1$.  Then, the lower and upper bounds $\udl{P}_\ell (k_\ell)$ and
$\ovl{P}_\ell (k_\ell)$ satisfying (\ref{req88}) can be computed recursively by  \bee &  & \udl{P}_1 (k_1) = \udl{\Up}_{0,x} (k_{1, x}) \; \udl{\Up}_{0,y} (k_{1, y}), \\
&  & \udl{P}_{\ell + 1} (k_{\ell+1} ) = \sum_{k_\ell \in \mscr{K}_\ell} \udl{P}_\ell (k_\ell) \; \udl{\Up}_{\ell,x} (k_{\ell, x}, k_{\ell+1, x})
\; \udl{\Up}_{\ell,y} (k_{\ell, y}, k_{\ell+1, y}), \qu  \ell = 1, \cd, s-1 \eee and
\bee &  & \ovl{P}_1 (k_1) = \ovl{\Up}_{0,x} (k_{1, x}) \; \ovl{\Up}_{0,y} (k_{1, y}), \\
&  & \ovl{P}_{\ell + 1} (k_{\ell+1} ) = \sum_{k_\ell \in \mscr{K}_\ell} \ovl{P}_\ell (k_\ell) \; \ovl{\Up}_{\ell,x} (k_{\ell, x}, k_{\ell+1, x})
\; \ovl{\Up}_{\ell,y} (k_{\ell, y}, k_{\ell+1, y}), \qu  \ell = 1, \cd, s-1, \eee where $k_1 \in \mscr{K}_1$ and $k_{\ell + 1} \in
\mscr{K}_{\ell + 1}$.  Making use of the bounds $\udl{P}_\ell (k_\ell)$ and $\ovl{P}_\ell (k_\ell)$ for $\Pr \{  \bs{K}_i \in \mscr{K}_i, \; 1
\leq i < \ell; \; \bs{K}_\ell = k_\ell \mid \bs{p} \}$, we have \be \la{good9883}
 \sum_{k_\ell \in \mscr{K}_\ell } \udl{P}_\ell (k_\ell) \leq \Pr \{  \bs{K}_i \in
\mscr{K}_i, \; 1 \leq i \leq \ell \mid \bs{p} \} \leq \sum_{k_\ell \in \mscr{K}_\ell } \ovl{P}_\ell (k_\ell), \qqu \ell = 1, \cd, s \ee
 for all $\bs{p} \in \mcal{Q}$.  Since the bounds of $\Psi(\ze, \bs{p})$ in (\ref{hi8}) and (\ref{hi9}) can
 be expressed as summations of probabilistic terms like $\Pr
\{  \bs{K}_i \in \mscr{K}_i, \; 1 \leq i \leq \ell \mid \bs{p} \}, \; \ell = 1, \cd, s$, we can establish lower and upper bounds of $\Psi(\ze,
\bs{p})$ with respect to $\bs{p} \in \mcal{Q}$ by using (\ref{hi8}),  (\ref{hi9}) and (\ref{good9883}).

By virtue of the lower and upper bounds of $\Psi(\ze, \bs{p})$ computed by the above recursive method,  we can employ Adapted Branch and Bound
technique in Appendix \ref{BBA} to determine whether $\Psi (\ze, \bs{p})$ is no greater than $\de$ for any $\bs{p} \in \Se$. The initial
hypercube $Q_{\mrm{init}}$ can be taken as $\{ (p_x, p_y): 0 \leq p_x \leq 1, \; 0 \leq p_y \leq 1 \}$. Given that it is possible to check the
truth of $\Psi (\ze, \bs{p}) \leq \de, \; \fa \bs{p} \in \Se$, we can apply a bisection search method to obtain the coverage tuning parameter
$\ze$ as large as possible such that the coverage probability $P(\ze, \bs{p})$ of the sequential random interval  $( \mscr{L}(\wh{\bs{p}},
\bs{l}), \mscr{U}(\wh{\bs{p}}, \bs{l}) )$ associated with $\ze$ is no less than $1 - \de$ for all $\bs{p} \in \Se$.

It should be noted that our proposed approach can be readily adapted for estimating functions of means of two Poisson populations, and functions
of proportions of two populations of finite sizes.

\sect{Estimation of Multinomial Proportions}

In probability theory, the multinomial distribution is a generalization of the binomial distribution. The binomial distribution is the
probability distribution of the number of ``successes'' in $n$ independent Bernoulli trials, with the same probability of ``success'' on each
trial. In a multinomial distribution, the analog of the Bernoulli distribution is the categorical distribution, where each trial results in
exactly one of some fixed finite number $\ka$ of possible outcomes, with probabilities $p_1, \cd, p_\ka$ (so that $p_\ell \geq 0$ for $\ell = 1,
\cd, \ka$ and $\sum_{\ell = 1}^\ka p_\ell = 1$), and there are $n$ independent trials.  For $\ell = 1, \cd, \ka$, let random variable $X_\ell$
denote the number of times that outcome number $\ell$ was observed over the $n$ trials. The vector $\bs{X} = (X_1, \cd, X_\ka)$ follows a
multinomial distribution with parameters $n$ and $\bs{p}$, where $\bs{p} = (p_1, \cd, p_{\ka - 1})$.  Define $\Se \DEF \{ (p_1, \cd, p_{\ka -
1}): \;  \sum_{\ell = 1}^{\ka - 1} p_\ell \leq 1 \; \tx{and} \;  p_\ell \geq 0, \; \ell = 1, \cd, \ka - 1  \}$.  For $\bs{p} \in \Se$, the
probability mass function of the multinomial distribution is given by
\[
\Pr \{ X_\ell = x_\ell, \; \ell = 1, \cd, \ka \mid \bs{p} \} = n! \prod_{\ell = 1}^{\ka} \f{ p_\ell^{x_\ell} }{x_\ell !},
\]
where $p_\ka \DEF 1 - \sum_{\ell = 1}^{\ka - 1} p_\ell$ and $x_1, \cd, x_k$ are non-negative integers such that $\sum_{\ell = 1}^\ka x_\ell =
n$. A classical problem in statistical inference is to construct a confidence region for $(p_1, \cd, p_{\ka})$.  Recently, Chafa\"{i} et. al.
proposed in \cite{Chafai} a confidence region with guaranteed confidence level and a small volume. However, such confidence region is difficult
to visualize for category number $\ka$ greater than $2$. This is especially true when the category number $\ka$ gets larger. On the other hand,
as a special type of confidence region, simultaneous confidence intervals offer a straightforward, intuitive, and direct assessment of the
reliability of the estimation. In many applications, it is desirable to construct simultaneous confidence intervals $[ \mcal{L}_\ell
(\wh{p}_\ell), \; \mcal{U}_\ell (\wh{p}_\ell) ], \; \ell = 1, \cd, \ka$ such that
\[
\Pr \{ \mcal{L}_\ell (\wh{p}_\ell) \leq  p_\ell \leq  \mcal{U}_\ell (\wh{p}_\ell), \; \ell = 1, \cd, \ka \mid \bs{p} \}  >  1 - \de
\]
for any $\bs{p} \in \Se$, where $\de \in (0, 1)$ is a pre-specified confidence parameter and $\wh{p}_\ell = \f{X_\ell}{n}$ for $\ell = 1, \cd,
\ka$. Here, $\Pr \{ E \mid \bs{p} \}$ denotes the probability of event $E$ which is determined by the parameter $\bs{p}$.  Wang (2008) made an
unsuccessful attempt in \cite{WangH} to solve this problem. Wang's method depends on statement (i) of her Lemma 2 (see, page 899 of
\cite{WangH}), which is actually an unproven claim. Moreover, even if the claim can be eventually proved, Wang's method still suffers from the
curse of dimensionality because the number of parameter points to be checked grows exponentially  with respect to the category number $\ka$. Her
justification for statement (i) of Lemma 2 is based on the following erroneous argument:

{\it Let $f(x)$ and $g(x)$ be two strictly increasing convex functions with respect to $x > 0$ such that $\lim_{x \downarrow 0} f(x) = 0, \;
\lim_{x \downarrow 0} g(x) > 0$ and that both $f(x)$ and $g(x)$ tend to be infinity as $x$ tends to some positive number.  Then, $f(x)$ and
$g(x)$ have at most two intersections}.

This argument is used by Wang in page 908 of her paper \cite{WangH}. Specifically, she applied the argument to functions $B(p_{k-1})$ and
$C(p_{k-1})$ defined in page 908 of \cite{WangH}. In line 15 from the bottom of page 908, she stated that ``$B(p_{k-1})$ and $C(p_{k-1})$ are
two strictly increasing functions''. Afterward, in lines 12-13 from the bottom of page 908, she concluded that ``there are at most two
intersections of $B(p_{k-1})$ and $C(p_{k-1})$''.  Unfortunately, Wang's argument is incorrect.  As a counterexample,
consider functions \bee &  &  f(x) = \bec x & \tx{for} \; 0 < x \leq 2,\\
x^2 - 3x + 4 & \tx{for} \; 2 < x < 4 \eec\\
&  & g(x) = \bec \f{15}{32} x^2 + \f{1}{8} & \tx{for} \; 0 < x \leq 2,\\
\f{15}{8} (x - 2) + 2  & \tx{for} \; 2 < x < 4 \eec \eee It can be readily shown that, $f(x)$ and $g(x)$ have three intersections.  Hence,
Wang's argument is disproved.  Since Wang's conclusion that ``there are at most two intersections of $B(p_{k-1})$ and $C(p_{k-1})$'' is based on
such an incorrect argument, the statement (i) of her Lemma 2 is certainly incorrect. This affects the validity of Wang's method for determining
the exact coefficients for simultaneous confidence intervals of multinomial proportions.

In view of the situation that there exists no exact method for constructing simultaneous confidence intervals for multinomial proportions
$\bs{p}$, we propose to solve this problem by virtue of the coverage tuning technique.  More specifically, our main idea is as follows:

(i) Seek a class of simultaneous confidence intervals $[ \mcal{L}_\ell (\wh{p}_\ell), \; \mcal{U}_\ell (\wh{p}_\ell) ], \; \ell = 1, \cd, \ka$
such that the coverage probability can be controlled by the coverage tuning parameter $\ze > 0$. In other words, the coverage probability,
denoted by $P(\ze, \bs{p}) \DEF \Pr \{ \mcal{L}_\ell (\wh{p}_\ell) \leq  p_\ell \leq \mcal{U}_\ell (\wh{p}_\ell), \; \ell = 1, \cd, \ka \mid
\bs{p} \}$, tends to $1$ for all $\bs{p} \in \Se$ as $\ze \to 0$.

(ii) For simplicity of establishing lower and upper bounds of coverage probability, choose the lower and upper confidence limits $\mcal{L}_\ell
(\wh{p}_\ell)$ and $\mcal{U}_\ell (\wh{p}_\ell)$ to be nondecreasing functions of $\wh{p}_\ell$ for $\ell = 1, \cd, \ka$.

(iii)  For a given value of coverage tuning parameter $\ze$, apply Adapted Branch and Bound method in Appendix \ref{BBA} to determine whether
the complementary coverage probability $\Psi (\ze, \bs{p}) \DEF  1 - P(\ze, \bs{p})$ of the simultaneous confidence intervals associated with
$\ze$ is no greater than $\de$ for any $\bs{p} \in \Se$.

(iv) Apply bisection coverage tuning method to determine $\ze > 0$ as large as possible such that $P(\ze, \bs{p}) \geq 1 - \de$ for any $\bs{p}
\in \Se$.

Actually, there are many methods to construct confidence intervals fulfilling the above requirements (i) and (ii) (see, e.g., \cite{Fitzpatrick,
Quesenberry} and the references therein).  Specially, four classes of simultaneous confidence intervals satisfying (i) and (ii) are described as
follows.

\bed

\item [Class A]: The lower and upper confidence limits are given by  {\small \[
 \mcal{L}_\ell (\wh{p}_\ell) = \f{ c + 2 n \wh{p}_\ell - \sq{ c^2 + 4 c n \wh{p}_\ell (1 - \wh{p}_\ell)    }  } { 2 (c + n) },  \qqu
\mcal{U}_\ell (\wh{p}_\ell) = \f{ c + 2 n \wh{p}_\ell + \sq{ c^2 + 4 c n \wh{p}_\ell (1 - \wh{p}_\ell) }  } { 2 (c + n) } \]} for $\ell = 1,
\cd, \ka$,  where $c \DEF \chi_{\ka - 1, \ze \de}^2$ is the $100(1-\ze \de)\%$ quantile of the chi-square distribution of $\ka - 1$ degrees of
freedom. This class of the simultaneous confidence intervals has been proposed in \cite{Quesenberry}.

\item [Class B]: The lower and upper confidence limits are given by
\[
\mcal{L}_\ell (\wh{p}_\ell) = \wh{p}_\ell - \f{ \mcal{Z}_{\ze \de \sh 2} }{ 2 \sq{n} },  \qqu \mcal{U}_\ell (\wh{p}_\ell) = \wh{p}_\ell + \f{
\mcal{Z}_{\ze \de \sh 2} }{ 2 \sq{n} }
\]
for $\ell = 1, \cd, \ka$. This class of  simultaneous confidence intervals was proposed in \cite{Fitzpatrick}.

\item [Class C]:  Making use of the binomial confidence interval established by Chen et. al. in \cite{Chen_CI}, we propose to  define simultaneous
confidence intervals with lower and upper confidence limits {\small \bel &  & \mcal{L}_\ell (\wh{p}_\ell) = \wh{p}_\ell + \frac{3}{4} \; \frac{
1 - 2 \wh{p}_\ell -
\sqrt{ 1 + \f{9 n} { 2 \ln \f{2}{\ze \de} } \wh{p}_\ell ( 1- \wh{p}_\ell) } } {1 + \f{9 n } { 8 \ln \f{2}{\ze \de} } }, \la{ChenCIA}\\
&  & \mcal{U}_\ell (\wh{p}_\ell) = \wh{p}_\ell + \frac{3}{4} \; \frac{ 1 - 2 \wh{p}_\ell + \sqrt{ 1 + \f{9 n} { 2 \ln \f{2}{\ze \de} }
\wh{p}_\ell ( 1 - \wh{p}_\ell) } } {1 + \f{9 n} { 8 \ln \f{2}{\ze \de} } } \la{ChenCIB} \eel} for $\ell = 1, \cd, \ka$.  By the coverage theory
in \cite{Chen_CI} and Bonferroni's inequality, the simultaneous confidence intervals defined by (\ref{ChenCIA}) and (\ref{ChenCIB}) guarantee
that the coverage probability  $P(\ze, \bs{p})$ is no less than $1 - \ka \ze \de$ for all $\bs{p} \in \Se$.

\item [Class D]:  Let $\mcal{L}_\ell (\wh{p}_\ell)$ and $\mcal{U}_\ell (\wh{p}_\ell)$ be the lower and upper confidence
limits of the binomial confidence interval proposed by Clopper and Pearson (1934) such that $\Pr \{ \mcal{L}_\ell (\wh{p}_\ell) \leq p_\ell \leq
\mcal{U}_\ell (\wh{p}_\ell) \mid p_\ell \} \geq 1 - \ze \de$ for $\ell = 1, \cd, \ka$.  As a consequence of Bonferroni's inequality,  the
simultaneous confidence intervals $[ \mcal{L}_\ell (\wh{p}_\ell), \; \mcal{U}_\ell (\wh{p}_\ell) ], \; \ell = 1, \cd, \ka$ have a coverage
probability no less than $1 - \ka \ze \de$ for all $\bs{p} \in \Se$.

\eed

Given that the structure of the simultaneous confidence intervals is determined so that the coverage probability $P(\ze, \bs{p})$ can be
controlled by $\ze$, we can apply the bisection coverage tuning method to obtain $\ze > 0$ as large as possible such that $P(\ze, \bs{p}) \geq 1
- \de$ for any $\bs{p} \in \Se$.  A critical step for bisection coverage tuning is to determine whether a given $\ze > 0$ is small enough to
ensure that the coverage probability $P(\ze, \bs{p})$ of the simultaneous confidence intervals associated with $\ze$ is no less than $1 - \de$
for any $\bs{p} \in \Se$. We propose to apply Adapted Branch and Bound technique in Appendix \ref{BBA} to accomplish this task. This needs
readily computable bounds of $\Psi(\ze, \bs{p})$ for $\bs{p} \in \Se$ in a hypercube $\mcal{Q} = \{ (p_1, \cd, p_{\ka - 1}): \udl{p}_\ell \leq
p_\ell \leq \ovl{p}_\ell, \; \ell = 1, \cd, \ka - 1 \}$, where $0 \leq \udl{p}_\ell \leq \ovl{p}_\ell \leq 1, \; \ell = 1, \cd, \ka - 1$ and
$\sum_{\ell = 1}^{\ka - 1} \udl{p}_\ell \leq 1$. We are going to establish the desired bounds in the sequel.

Let $\bs{a}= (a_1, \cd, a_\nu)$ and $\bs{b} = (b_1, \cd, b_\nu)$ be integer-valued vectors. Let $\bs{\se} = (\se_1, \cd, \se_\nu)$ be a vector
of nonnegative elements $\se_i, \; i = \cd, \nu$.  Define multivariate function
\[
S (\bs{a}, \bs{b}, \bs{\se}, n) = \sum_{x_1 = \max \{0, a_1\} }^{\min \{n, b_1\}} \cd \sum_{x_\nu = \max \{0, a_\nu \} }^{\min \{n, b_\nu\} } n!
\f{ \se_1^{x_1} \cd \se_\nu^{x_\nu} }{x_1 ! \cd x_\nu !}  \bb{I} \li ( \sum_{\ell = 1}^\nu x_\ell = n \ri ),
\]
where $\bb{I} (\mcal{E})$ denotes the indicator function for the event $\mcal{E}$.

Assume that $\mcal{U}_\ell(1) \geq 1$ and $\mcal{L}_\ell(0) \leq 0$ for $\ell = 1, \cd, \ka$.  For $\ell = 1, \cd, \ka$, since the lower
confidence limit $\mcal{L}_\ell(\wh{p}_\ell)$ is a nondecreasing function of $\wh{p}_\ell$, we can define inverse function $\mcal{L}_\ell^{-1}
(.)$ of $ \mcal{L}_\ell (.)$ such that $\mcal{L}_\ell^{-1} (\se) = \max \{ z \in I_{\wh{p}_\ell}: \mcal{L}_\ell (z) \leq \se \}$ for $\se \in
[0, 1]$, where $I_{\wh{p}_\ell}$ denotes the support of $\wh{p}_\ell$. Similarly, for $\ell = 1, \cd, \ka$, since the upper confidence limit
$\mcal{U}_\ell(\wh{p}_\ell)$ is a nondecreasing function of $\wh{p}_\ell$, we can define inverse function $\mcal{U}_\ell^{-1} (.)$ of $
\mcal{U}_\ell (.)$ such that $\mcal{U}_\ell^{-1} (\se) = \min \{ z \in I_{\wh{p}_\ell}: \mcal{U}_\ell (z) \geq \se \}$ for $\se \in [0, 1]$. Let
$\udl{p}_\ka = \max \{ 0, 1 - \sum_{\ell = 1}^{\ka - 1} \ovl{p}_\ell \}$ and $\ovl{p}_\ka = 1 - \sum_{\ell = 1}^{\ka - 1} \udl{p}_\ell$.  A
natural method is to directly bound $P(\ze, \bs{p})$ by {\small \bel &   & P(\ze, \bs{p}) \leq \Pr \{ \mcal{U}_\ell^{-1} (\udl{p}_\ell) \leq
\wh{p}_\ell \leq \mcal{L}_\ell
^{-1} ( \ovl{p}_\ell ), \; \ell = 1, \cd, \ka \mid \bs{p} \} \leq S(\bs{A}, \bs{B}, \ovl{\bs{p}}, n), \la{upB}\\
&  & P(\ze, \bs{p}) \geq \Pr \{ \mcal{U}_\ell^{-1} (\ovl{p}_\ell) \leq \wh{p}_\ell  \leq \mcal{L}_\ell ^{-1} ( \udl{p}_\ell ), \; \ell = 1, \cd,
\ka \mid \bs{p} \} \geq S(\bs{C}, \bs{D}, \udl{\bs{p}}, n)  \la{lwB} \eel} for any $\bs{p} \in \Se \cap \mcal{Q}$, where $\udl{\bs{p}} =
(\udl{p}_1, \cd, \udl{p}_\ka), \; \ovl{\bs{p}} = (\ovl{p}_1, \cd, \ovl{p}_\ka)$ and $\bs{A} = (A_1, \cd, A_\ka), \; \bs{B} = (B_1, \cd, B_\ka)$
with
\[
A_\ell = n \mcal{U}_\ell^{-1} (\udl{p}_\ell), \qu B_\ell = n \mcal{L}_\ell ^{-1} ( \ovl{p}_\ell ), \qu  C_\ell = n \mcal{U}_\ell^{-1}
(\ovl{p}_\ell), \qu D_\ell = n \mcal{L}_\ell ^{-1} ( \udl{p}_\ell )
\]
for $\ell = 1, \cd, \ka$.  Since $\Psi(\ze, \bs{p}) = 1 - P(\ze, \bs{p})$, the bounds of $\Psi(\ze, \bs{p})$ can be obtained from the bounds of
$P(\ze, \bs{p})$ for $\bs{p} \in \Se \cap \mcal{Q}$.  At the first glance, the above bounds (\ref{upB}) and (\ref{lwB}) are simple and useful.
Unfortunately, the lower bound (\ref{lwB}) of $P(\ze, \bs{p})$ is useless in practice.  The reason is that, as $\sum_{\ell = 1}^{\ka - 1}
\ovl{p}_\ell \to 1$, we have $\udl{p}_\ka \to 0$ and thus the lower bound (\ref{lwB}) of $P(\ze, \bs{p})$ tends to $0$.  This implies that it is
impossible to verify that $P(\ze, \bs{p})$ is no less than $1 - \de$ for those $\bs{p} \in \Se$ with $\sum_{\ell = 1}^{\ka - 1} p_\ell \ap 1$ by
showing that the lower bound of $P(\ze, \bs{p})$ is no less than $1 - \de$.  In view of this situation, we propose to directly bound the
complementary coverage probability $\Psi (\ze, \bs{p}) = 1 - P (\ze, \bs{p})$.  In this direction, we have established the following result.
\beT \la{multibound} Let $\eta \in (0,1)$.  Let $\mrm{T}_{\mrm{lb}} (., ., .)$ and $\mrm{T}_{\mrm{ub}} (., ., .)$ be multivariate functions
defined by (\ref{truna}) and (\ref{trunb}). For $i = 1, \cd, \ka$, define $\nu_i = \min \{\ka, i + 1 \}$, \bee & & \udl{\se}_{i, \ell} =
\udl{p}_\ell, \qqu
\ovl{\se}_{i, \ell} = \ovl{p}_\ell, \qqu \ell = 1, \cd, \nu_i - 1;\\
&  & \udl{\se}_{i, \nu_i} = \max \li \{ 0, \; 1 - \sum_{\ell = 1}^{\nu_i - 1} \ovl{p}_\ell \ri \}, \qqu \ovl{\se}_{i, \nu_i} = 1 - \sum_{\ell =
1}^{\nu_i - 1} \udl{p}_\ell \eee and $\udl{\bs{\se}}_i = (\udl{\se}_{i,1}, \cd, \udl{\se}_{i, \nu_i}), \; \ovl{\bs{\se}}_i = (\ovl{\se}_{i,1},
\cd, \ovl{\se}_{i,\nu_i})$. For $i = 1, \cd, \ka$, define integer-valued vectors \bee \bs{\mcal{A}}_i = (A_{i,1}, \cd, A_{i,\nu_i}), \qu
\bs{\mcal{B}}_i
= (B_{i,1}, \cd, B_{i,\nu_i}), \qu \bs{\mcal{C}}_i = (C_{i,1}, \cd, C_{i,\nu_i}), \qu \bs{\mcal{D}}_i = (D_{i,1}, \cd, D_{i,\nu_i}), & &\\
\bs{\fra{A}}_i = (\fra{A}_{i,1}, \cd, \fra{A}_{i,\nu_i}), \qu \bs{\fra{B}}_i = (\fra{B}_{i,1}, \cd, \fra{B}_{i,\nu_i}), \qu \bs{\fra{C}}_i =
(\fra{C}_{i,1}, \cd, \fra{C}_{i,\nu_i}), \qu \bs{\fra{D}}_i = (\fra{D}_{i,1}, \cd, \fra{D}_{i,\nu_i}) \; &  & \eee with {\small \bee A_{i, \ell}
= \fra{A}_{i, \ell} = n \times \max \li \{ \mcal{U}_\ell^{-1} (\udl{\se}_{i, \ell}), \;  \mrm{T}_{\mrm{lb}} (\udl{\se}_{i, \ell}, n, \eta) \ri
\}, \qu B_{i, \ell} = \fra{B}_{i, \ell} =  n \times \min \li \{ \mcal{L}_\ell ^{-1} ( \ovl{\se}_{i, \ell} ), \;
\mrm{T}_{\mrm{ub}} (\ovl{\se}_{i, \ell}, n, \eta) \ri \}, &  &\\
C_{i, \ell} = \fra{C}_{i, \ell} = n \times \max \li \{ \mcal{U}_\ell^{-1} (\ovl{\se}_{i, \ell}), \; \mrm{T}_{\mrm{lb}} (\ovl{\se}_{i, \ell}, n,
\eta) \ri \}, \qu
 D_{i, \ell} = \fra{D}_{i, \ell} = n \times \min
\li \{ \mcal{L}_\ell ^{-1} ( \udl{\se}_{i, \ell} ), \; \mrm{T}_{\mrm{ub}} (\udl{\se}_{i, \ell}, n, \eta) \ri \}, &  & \eee} for $i = 2, \cd,
\ka$ and $\ell = 1, \cd, i - 1$; {\small \bee &  & \fra{A}_{i, i} =  n  \mrm{T}_{\mrm{lb}} (\udl{\se}_{i, i}, n, \eta), \qu \fra{C}_{i, i} =  n
\mrm{T}_{\mrm{lb}} (\ovl{\se}_{i, i}, n, \eta), \qu B_{i, i} =  n  \mrm{T}_{\mrm{ub}} (\ovl{\se}_{i, i}, n, \eta), \qu D_{i, i} = n
\mrm{T}_{\mrm{ub}} (\udl{\se}_{i, i}, n, \eta),\\
& & A_{i, i} = \max \li \{ n \mcal{L}_i^{-1} (\udl{\se}_{i, i}) + 1 , \; \fra{A}_{i, i} \ri \}, \qqu C_{i, i}  = \max \li \{ n \mcal{L}_i^{-1}
(\ovl{\se}_{i, i}) + 1, \; \fra{C}_{i, i} \ri \}, \\
 & & \fra{B}_{i, i} = \min \li \{ n \mcal{U}_i^{-1} (\ovl{\se}_{i, i}) - 1 , \; B_{i, i} \ri \}, \qqu \fra{D}_{i, i} = \min \li \{ n \mcal{U}_i^{-1} (\udl{\se}_{i, i})
- 1 , \; D_{i, i} \ri \} \eee} for $i = 1, \cd, \ka$; and {\small \bee &   & A_{i, i+1}  = \fra{A}_{i, i+1} = n \mrm{T}_{\mrm{lb}}
(\udl{\se}_{i, i+1}, n, \eta), \qqu
 B_{i, i+1} = \fra{B}_{i, i+1} =  n  \mrm{T}_{\mrm{ub}} (\ovl{\se}_{i, i+1}, n, \eta), \\
 &  &  C_{i, i+1}  = \fra{C}_{i, i+1} = n \mrm{T}_{\mrm{lb}} (\ovl{\se}_{i, i+1}, n, \eta),  \qqu
 D_{i, i+1} = \fra{D}_{i, i+1} = n  \mrm{T}_{\mrm{ub}} (\udl{\se}_{i, i+1}, n, \eta) \eee} for $i = 1, \cd, \ka - 1$. Then, {\small \bee \sum_{i = 1}^\ka [ S (\bs{\mcal{C}}_i, \bs{\mcal{D}}_i,
\udl{\bs{\se}}_i, n) + S (\bs{\fra{C}}_i, \bs{\fra{D}}_i, \udl{\bs{\se}}_i, n)] \leq \Psi (\ze, \bs{p}) \leq (2 \ka - 1) \eta  + \sum_{i =
1}^\ka [ S (\bs{\mcal{A}}_i, \bs{\mcal{B}}_i, \ovl{\bs{\se}}_i, n) + S (\bs{\fra{A}}_i, \bs{\fra{B}}_i, \ovl{\bs{\se}}_i, n)] \eee} for all
$\bs{p} \in \Se \cap \mcal{Q}$.  \eeT

See Appendix \ref{multibound_app} for a proof.   Actually, we have used the truncation method proposed in \cite{Chen1} to reduce the
computational complexity for bounding $\Psi (\ze, \bs{p})$.  The term $(2 \ka - 1) \eta$ in the upper bound of $\Psi (\ze, \bs{p})$ is due to
truncation. In applications, the truncation tolerance $\eta$ can be chosen to be an extremely small number (e.g., $\eta = \f{\de}{1000}$), while
the reduction of computational complexity can be substantial.  It should be noted that the application of the truncation method introduces no
approximation, but only negligible conservatism.

To further reduce the computational complexity of bounding $\Psi(\ze, \bs{p})$ by virtue of Theorem \ref{multibound}, we can use the recursive
method developed by Frey (2009) to evaluate terms like $S (\bs{a}, \bs{b}, \bs{\se}, n)$ with $\bs{\se} = (\se_1, \cd, \se_\nu)$ and
integer-valued vectors $\bs{a} = (a_1, \cd, a_\nu), \; \bs{b} = (b_1, \cd, b_\nu)$, where $\se_i
> 0, \; n \geq b_i \geq a_i \geq 0, \; i = \cd, \nu$. Under the assumption that $\sum_{i = 1}^\nu a_i < n < \sum_{i = 1}^\nu b_i$ and that $b_i
> a_i$ for at least one $i$ among $i = 1, \cd, \nu$, Frey shows in \cite{Frey} that
\[
S (\bs{a}, \bs{b}, \bs{\se}, n) = n! \li ( \prod_{\ell = 1}^\nu \f{\se_\ell^{a_\ell}}{ a_\ell! } \ri ) \times \sum_{i = 1}^\nu \sum_{j = 1}^{b_i
- a_i} P_r (i, j),
\]
where $r = n - \sum_{i = 1}^\nu a_i$ and $P_r (i, j), \; i = 1, \cd, \nu; \; j = 1, \cd, b_i - a_i$ can be recursively computed by
\[
P_1 (i, j) = \bec \f{\se_i}{a_i + 1}  & \tx{if} \; j = 1 \; \tx{and} \; b_i > a_i,\\
0 & \tx{otherwise} \eec
\]
for $1 \leq i \leq k, \; 1 \leq j \leq b_i - a_i$;
\[
P_{t + 1} (i, 1) = \f{\se_i}{a_i + 1}  \sum_{\ell = 1}^{i - 1} \sum_{j = 1}^{ \min \{b_\ell - a_\ell, t\} }  P_t (\ell, j) \qu \tx{for $1 \leq i
\leq k, \; 1 \leq t \leq r - 1$};
\]
 and $P_{t + 1} (i, j) = \f{\se_i}{a_i + j} P_t (i, j - 1)$ for $1 \leq t \leq r - 1, \; 1 \leq i
\leq k, \; 1 < j \leq b_i - a_i$.

By virtue of the lower and upper bounds of $\Psi(\ze, \bs{p})$ given by Theorem \ref{multibound},  we can employ Adapted Branch and Bound
technique in Appendix \ref{BBA} to test if $\Psi (\ze, \bs{p})$ is no greater than $\de$ for any $\bs{p} \in \Se$. The initial hypercube
$Q_{\mrm{init}}$ can be taken as $\{ (p_1, \cd, p_{\ka - 1}): 0 \leq p_\ell \leq 1, \; \ell = 1, \cd, \ka - 1 \}$. In the branching process,
some hypercubes will be generated which have no intersection with the parameter space $\Se$, i.e., $\sum_{\ell = 1}^{\ka - 1} \udl{p}_\ell > 1$,
where $\udl{p}_\ell$ is the lower bound for $p_\ell$ of $\bs{p}$ in the hypercube. Such hypercubes should be eliminated from further
consideration. Given that it is possible to check the truth of $\Psi (\ze, \bs{p}) \leq \de, \; \fa \bs{p} \in \Se$, we can apply a bisection
search method to determine the coverage tuning parameter $\ze$ as large as possible such that the coverage probability $P(\ze, \bs{p})$ of the
simultaneous confidence intervals associated with $\ze$ is no less than $1 - \de$ for all $\bs{p} \in \Se$.

In the above discussion, we have been focusing on the interval estimation for the multinomial proportions $\bs{p}$ when the sample size $n$ is
given. In many applications, it is important to determine appropriate sample size $n$ such that the estimators $\wh{p}_\ell$ for $p_\ell, \;
\ell = 1, \cd, \ka$ satisfy some pre-specified requirements of reliability. In general, the problem can be formulated as follows.  Let $\de \in
(0, 1)$ be a pre-specified confidence parameter.  For $\ell = 1, \cd, \ka$, let the margin of error for $p_\ell$ be $\ep_\ell \DEF \max \{
\vep_{a, \ell}, \; \vep_{r, \ell} \; p_\ell \}$,  where $\vep_{a, \ell} \in [0, 1)$ and $\vep_{r, \ell} \in [0, 1)$. Such a margin of error can
be reduced to the margin of absolute error and the margin of relative error by taking $\vep_{r, \ell} = 0$ and $\vep_{a, \ell} = 0$
respectively. For pre-specified confidence parameter $\de$ and margins of error $\ep_\ell = \max \{  \vep_{a, \ell}, \; \vep_{r, \ell} p_\ell
\}, \; \ell = 1, \cd, \ka$, a problem of practical importance is to determine a sample size $n$ as small as possible such that
\[
\Pr \{ | \wh{p}_\ell - p_\ell | \leq \ep_\ell, \; \ell = 1, \cd, \ka \mid \bs{p} \}  \geq 1 - \de
\]
for all $\bs{p} \in \Se$. We can apply the above technique to solve this problem.  Invoking identity (\ref{RI_Indentity}), we have
\[
\{ | \wh{p}_\ell - p_\ell | \leq \ep_\ell \} = \li \{ \mcal{L}_\ell (\wh{p}_\ell) \leq p_\ell \leq \mcal{U}_\ell  (\wh{p}_\ell) \ri \}, \qqu
\ell = 1, \cd, \ka
\]
where \be \la{defsim}
 \mcal{L}_\ell (\wh{p}_\ell) = \min \li \{  \wh{p}_\ell - \vep_{a, \ell}, \; \f{ \wh{p}_\ell } { 1 + \vep_{r, \ell} }  \ri
\}, \qqu \mcal{U}_\ell (\wh{p}_\ell)  = \max \li \{ \wh{p}_\ell + \vep_{a, \ell}, \; \f{ \wh{p}_\ell } { 1 - \vep_{r, \ell} }  \ri \} \ee for
$\ell = 1, \cd, \ka$. Hence, $\Pr \{ | \wh{p}_\ell - p_\ell | \leq \ep_\ell, \; \ell = 1, \cd, \ka \mid \bs{p} \}  = \Pr \{ \mcal{L}_\ell
(\wh{p}_\ell) \leq p_\ell \leq \mcal{U}_\ell (\wh{p}_\ell), \; \ell = 1, \cd, \ka \mid \bs{p} \}$.   This implies that the sample size problem
is equivalent to finding the smallest $n$ such that the coverage probability, denoted by $\mcal{P}(n, \bs{p})$, of the simultaneous confidence
intervals defined by (\ref{defsim}) is no less than $1 - \de$.  Clearly, the lower and upper limits of the simultaneous confidence intervals
defined by (\ref{defsim}) are increasing functions of $\wh{p}_\ell$ for $\ell = 1, \cd, \ka$.  Moreover, $\mcal{P}(n, \bs{p})$ tends to $1$ for
any $\bs{p} \in \Se$ as $n$ tends to infinity.  Therefore, for a given sample size $n$, we can apply the same methods for bounding $\Psi(\ze,
\bs{p})$ to establish lower and upper bounds for $1 - \mcal{P}(n, \bs{p})$ with respect to $\bs{p}$ in a hypercube. Consequently, as determining
the truth of $\Psi(\ze, \bs{p}) \leq \de, \; \fa \bs{p} \in \Se$, we can determine the truth of $1 - \mcal{P}(n, \bs{p}) \leq \de, \; \fa \bs{p}
\in \Se$. Given that such a routine can be established, we can obtain the smallest sample size $n$ such that $\mcal{P}(n, \bs{p}) \geq 1 - \de,
\; \fa \bs{p} \in \Se$ by checking $n$ from small to large enough.

\sect{Estimation of Bounded-Variable Means}

In Section 4, we have been focusing on the estimation of binomial parameters.
 Actually, some of the ideas  can be generalized to the estimation
 of means of random variables bounded in interval $[0, 1]$.
 Formally, let $X \in [0, 1]$ be a random variable with expectation $\mu =
 \bb{E} [X]$. We can estimate $\mu$ based on i.i.d. random samples
 $X_1, X_2, \cd$ of $X$ by virtue of multistage sampling schemes.

\subsection{Control of Absolute Error}

To estimate the mean of the bounded variable $X \in [0, 1]$ with an
absolute error criterion, we have multistage sampling schemes
described by the following theorems.

\beT \la{Bounded_Mean_abs_Hoeffding} Let $0 < \vep < \f{1}{2}$. Let
$n_1 < n_2 < \cd < n_s$ be a sequence of sample sizes such that
{\small $n_s \geq \f{ \ln \f{2 s} { \de} } { 2 \vep^2 }$}. Define
$\wh{\bs{\mu}}_\ell = \f{ \sum_{i=1}^{n_\ell} X_i }{n_\ell}$ for
$\ell = 1, \cd, s$.  Suppose that sampling is continued until
$\mscr{M}_{\mrm{B}} ( \f{1}{2} - |\f{1}{2} - \wh{\bs{\mu}}_\ell | ,
\f{1}{2} - |\f{1}{2} - \wh{\bs{\mu}}_\ell | + \vep) \leq \f{1} {
n_\ell } \ln \li ( \f{\de}{2 s} \ri )$. Define $\bs{\wh{\mu}} =
\f{\sum_{i=1}^{\mathbf{n}} X_i}{\mathbf{n}}$, where $\mathbf{n}$ is
the sample size when the sampling is terminated.
 Then, $\Pr \li \{ \li | \bs{\wh{\mu}} - \mu \ri | < \vep \ri \}
\geq 1 - \de$. \eeT

See Appendix \ref{App_Bounded_Mean_abs_Hoeffding} for a proof.

\beT \la{Bounded_Mean_ABS_Massart} Let $0 < \vep < \f{1}{2}$. Let
$n_1 < n_2 < \cd < n_s$ be a sequence of sample sizes such that
{\small $n_s \geq \f{ \ln \f{2 s} { \de} } { 2 \vep^2 }$}. Define
$\wh{\bs{\mu}}_\ell = \f{ \sum_{i=1}^{n_\ell} X_i }{n_\ell}$ for
$\ell = 1, \cd, s$.  Suppose that sampling is continued until
{\small $\li ( \li | \wh{\bs{\mu}}_\ell - \f{1}{2} \ri | - \f{2 \vep
}{3} \ri )^2 \geq \f{1}{4} - \f{ \vep^2 n_\ell } {2 \ln (2 s \sh
\de) }$} for some $\ell \in \{1, \cd, s \}$. Define $\bs{\wh{\mu}} =
\f{\sum_{i=1}^{\mathbf{n}} X_i}{\mathbf{n}}$, where $\mathbf{n}$ is
the sample size when the sampling is terminated.
 Then, $\Pr \li \{ \li | \bs{\wh{\mu}} - \mu \ri | < \vep \ri \}
\geq 1 - \de$. \eeT

See Appendix \ref{App_Bounded_Mean_ABS_Massart} for a proof.  Actually, the estimation of $\mu$ can be more efficient by a computational
approach illustrated in the sequel.

Let $\ro \in (0, 1]$ and $\ze \in (0, \f{1}{\de})$.   By replacing the constant $\f{2}{3}$ and the quantity $\ln \f{\de}{2 s}$ in the stopping
rule described in Theorem \ref{Bounded_Mean_ABS_Massart} as $\ro$ and $\ln (\ze \de)$ respectively, we obtain a stopping rule for estimating the
mean value $\mu$ as follows:

\be \la{stmean} \tx{Continue sampling until {\small $\li ( \li | \wh{\bs{\mu}}_\ell - \f{1}{2} \ri | - \ro \vep \ri )^2 \geq \f{1}{4} + \f{
\vep^2 n_\ell } {2 \ln ( \ze \de) }$} for some $\ell \in \{1, \cd, s \}$}. \ee

As in Theorem \ref{Bounded_Mean_ABS_Massart}, at the termination of the sampling process, the estimator for $\mu$ is taken as $\bs{\wh{\mu}} =
\f{\sum_{i=1}^{\mathbf{n}} X_i}{\mathbf{n}}$.  According to the general criteria proposed in Section \ref{gen_structure}, the minimum sample
size $n_1$ and the maximum sample size $n_s$ should satisfy the following constraint: \be \la{constraint} \f{2 \ro ( 1 - \ro \vep) \ln \f{1}{\ze
\de} }{ \vep } \leq n_1 \leq \cd \leq n_{s-1} < \f{ \ln \f{1}{\ze \de} }{ 2 \vep^2} \leq n_s. \ee

To ensure that the coverage probability $\Pr \li \{ \li | \bs{\wh{\mu}} - \mu \ri | < \vep \mid \mu \ri \}$ is no less than the prescribed
confidence level $1 - \de$, we can apply the coverage tuning technique proposed in Section \ref{Computing}.  For this purpose, we need to
develop a computable upper bound for the complementary coverage probability for $\mu$ contained in an interval.  For $\ell = 1, \cd, s$, define
a decision variable $\bs{D}_\ell$ for the $\ell$-th stage such that $\bs{D}_\ell$ assumes value $1$ if {\small $\li ( \li | \wh{\bs{\mu}}_\ell -
\f{1}{2} \ri | - \ro \vep \ri )^2 \geq \f{1}{4} + \f{ \vep^2 n_\ell } {2 \ln ( \ze \de) }$} and otherwise assumes value $0$.   Consequently, the
stopping rule (\ref{stmean}) can be expressed  as ``continue the sampling process until $\bs{D}_\ell = 1$ for some $\ell \in \{1, \cd, s\}$''.
Similar to (\ref{relaxedB}), we propose to bound the complementary coverage probability $\Pr \li \{ \li | \bs{\wh{\mu}} - \mu \ri | \geq \vep
\mid \mu \ri \}$ with respect to $\mu \in [a, b] \subseteq [0, 1]$ as follows: {\small \be \Pr \li \{ \li | \bs{\wh{\mu}} - \mu \ri | \geq \vep
\mid \mu \ri \}  \leq  \sum_{ \ell = 1 }^s \min \li [ \Pr \{ \bs{D}_{\ell - 1} = 0 \mid \mu \}, \; \Pr \{ | \wh{\bs{\mu}}_\ell - \mu | \geq \vep
\mid \mu \}, \; \Pr \{ \bs{D}_{\ell} = 1 \mid \mu \} \ri ].  \qqu \la{relaxedBB} \ee}  To bound the complementary coverage probability based on
$(\ref{relaxedBB})$, we need to have upper bounds of $\Pr \{ \bs{D}_{\ell} = 0 \mid \mu \}, \; \Pr \{ \bs{D}_{\ell} = 1 \mid \mu \}$ and $\Pr \{
| \wh{\bs{\mu}}_\ell - \mu | \geq \vep \mid \mu \}$ with respect to $\mu \in [a, b]$ for $\ell = 1, \cd, s$.   This can be accomplished as
follows.

By the symmetry of the stopping rule, the coverage probability $\Pr \li \{ \li | \bs{\wh{\mu}} - \mu \ri | < \vep \mid \mu \ri \}$ is
symmetrical about $\mu = \f{1}{2}$. Hence, it suffices to consider $\mu \in [0, \f{1}{2}]$ regarding the coverage probability. Define
\[ c_\ell = \f{1}{2} - \ro \vep - \sq{\f{1}{4} + \f{ \vep^2 n_\ell } { 2 \ln (\ze \de) }}, \qqu \ell = 1, \cd, s.
\]
Let $[a, b] \subseteq [0, \f{1}{2}]$.  Clearly, $\Pr \{ \bs{D}_0 = 0 \mid \mu \}  = \Pr \{ \bs{D}_s = 1 \mid \mu \} = 1$ for $\mu \in [a, b]$.
By virtue of Chernoff-Hoeffding bounds \cite{Chernoff, Hoeffding} and the definition of the stopping rule, for $\mu \in [a, b]$ and $\ell = 1,
\cd, s - 1$, we have \bee \Pr \{ \bs{D}_{\ell} = 0 \mid \mu \} & \leq & \bec \Pr \{ \wh{\bs{\mu}}_\ell
> c_\ell \mid \mu \}  & \;
\tx{if} \; \tx{$b <  c_\ell  < \f{1}{2} - \ro \vep$},\\
1 & \; \tx{otherwise} \eec\\
& \leq &  \bec \exp (n_\ell \mscr{M}_{\mrm{B}} ( c_\ell, b  )  )  & \;
\tx{if} \; \tx{$b <  c_\ell  < \f{1}{2} - \ro \vep$},\\
1 & \; \tx{otherwise} \eec\\
\eee and {\small \bee \Pr \{ \bs{D}_\ell = 1 \mid \mu \} & \leq & \bec \Pr \{  \wh{\bs{\mu}}_\ell \geq 1 - c_\ell \mid \mu \} + \Pr \{
\wh{\bs{\mu}}_\ell \leq c_\ell \mid \mu \} & \;
\tx{if} \; \tx{$c_\ell < a$ and $n_\ell < \f{ 2 \ln (\ze \de) [ (\ro \vep)^2 - \f{1}{4} ] }{ \vep^2 }$},\\
1 & \; \tx{otherwise} \eec \\
& \leq & \bec \exp ( n_\ell \mscr{M}_{\mrm{B}} ( 1 - c_\ell, b  ) ) + \exp ( n_\ell \mscr{M}_{\mrm{B}}  (  c_\ell, a ) ) & \;
\tx{if} \; \tx{$c_\ell < a$ and $n_\ell < \f{ 2 \ln (\ze \de) [ (\ro \vep)^2 - \f{1}{4} ] }{ \vep^2 }$},\\
1 & \; \tx{otherwise} \eec \eee} For  $\mu \in [a, b]$ and $\ell = 1, \cd, s$, under the assumption that $b - a < \vep$, it follows from
Chernoff-Hoeffding bounds that  \bee \Pr \{ | \wh{\bs{\mu}}_\ell - \mu | \geq \vep \mid \mu \} & \leq & \Pr \{ \wh{\bs{\mu}}_\ell \geq \mu +
\vep \mid \mu \} + \Pr \{
\wh{\bs{\mu}}_\ell \leq \mu - \vep \mid \mu \}\\
& \leq & \Pr \{ \wh{\bs{\mu}}_\ell \geq a + \vep \mid \mu \} + \Pr \{ \wh{\bs{\mu}}_\ell \leq b - \vep \mid \mu \}\\
& \leq & \exp (n_\ell \mscr{M}_{\mrm{B}} ( a + \vep, b ) ) + \exp ( n_\ell \mscr{M}_{\mrm{B}} ( b - \vep,  a ) ).  \eee  Applying the above
upper bounds for $\Pr \{ \bs{D}_{\ell} = 0 \mid \mu \}, \; \Pr \{ \bs{D}_{\ell} = 1 \mid \mu \}$ and $\Pr \{ | \wh{\bs{\mu}}_\ell - \mu | \geq
\vep \mid \mu \}$ to (\ref{relaxedBB}) gives an upper bound for the complementary coverage probability $\Pr \{ | \wh{\bs{\mu}} - \mu | \geq \vep
\mid \mu \}$ with respect to $\mu \in [a, b]$.  With such interval bounding technique for the complementary coverage probability, we can use
AMCA to check whether the complementary coverage probability associated with given $\ze$ and $\ro$ is no greater than $\de$ for any $\mu \in [0,
1]$. Consequently, for any fixed $\ro > 0$, by virtue of the coverage tuning technique proposed in Section \ref{Computing}, it is possible to
determine an appropriate value for $\ze$ such that the coverage probability $\Pr \li \{ \li | \bs{\wh{\mu}} - \mu \ri | < \vep \mid \mu \ri \}$
is no less than the prescribed confidence level $1 - \de$ for any underlying $\mu \in (0, 1)$.   We can optimize the performance of the
resultant stopping rule by choosing different values of $\ro > 0$ and computing the corresponding values of $\ze >0$.

\subsection{Control of Relative Error}

To estimate the mean of the  bounded variable $X \in [0, 1]$ with a
relative precision, we have multistage inverse sampling schemes
described by the following theorems.

 \beT
\la{Bounded_Mean_Rev_Hoeffding} Let $0 < \vep < 1$. Let $\ga_1 <
\ga_2 < \cd < \ga_s$ be a sequence of real numbers such that $\ga_1
> \f{1}{\vep}$ and {\small $\ga_s \geq \f{ (1 + \vep) \ln \f{2 s} {
\de} } { (1 + \vep) \ln (1 + \vep) - \vep}$}. For $\ell = 1, \cd,
s$, define $\wh{\bs{\mu}}_\ell = \f{ \ga_\ell }{\mbf{n}_\ell}$,
where $\mbf{n}_\ell$ is the minimum sample number such that
$\sum_{i=1}^{\mbf{n}_\ell} X_i \geq \ga_\ell$. Suppose that sampling
is continued until $ \mscr{M}_{\mrm{B}} (
\f{\ga_\ell}{\mbf{n}_\ell}, \f{\ga_\ell}{\mbf{n}_\ell (1 + \vep)}  )
\leq \f{1} {\mbf{n}_\ell } \ln \li ( \f{\de}{2 s} \ri )$ and
$\mscr{M}_{\mrm{B}} ( \f{\ga_\ell}{\mbf{n}_\ell - 1},
\f{\ga_\ell}{\mbf{n}_\ell (1 - \vep)} ) \leq \f{1} {\mbf{n}_\ell - 1
} \ln \li ( \f{\de}{2 s} \ri )$ for some $\ell \in \{1, \cd, s \}$.
Define $\bs{\wh{\mu}} = \f{ \ga_{\bs{l}} }{\mathbf{n}_{\bs{l}}}$,
where $\bs{l}$ is the index of stage when the sampling is
terminated.
 Then, $\Pr \li \{ \li | \bs{\wh{\mu}} - \mu \ri | < \vep \mu \ri \}
\geq 1 - \de$. \eeT

\beT \la{Bounded_Mean_Rev_Massart} Let $0 < \vep < 1$.  Let $\ga_1 <
\ga_2 < \cd < \ga_s$ be a sequence of real numbers such that $\ga_1
> \f{1}{\vep}$ and {\small $\ga_s \geq \f{ 2 (1 + \vep) (3 + \vep)
\ln \f{2 s} { \de} } { 3 \vep^2 }$}.  For $\ell = 1, \cd, s$, define
$\wh{\bs{\mu}}_\ell = \f{ \ga_\ell }{\mbf{n}_\ell}$, where
$\mbf{n}_\ell$ is the minimum sample number such that
$\sum_{i=1}^{\mbf{n}_\ell} X_i \geq \ga_\ell$.  Suppose that
sampling is continued until $ \mscr{M} ( \f{\ga_\ell}{\mbf{n}_\ell},
\f{\ga_\ell}{\mbf{n}_\ell (1 + \vep)}  ) \leq \f{1} {\mbf{n}_\ell }
\ln \li ( \f{\de}{2 s} \ri )$ and $\mscr{M} (
\f{\ga_\ell}{\mbf{n}_\ell - 1}, \f{\ga_\ell}{\mbf{n}_\ell (1 -
\vep)} ) \leq \f{1} {\mbf{n}_\ell - 1 } \ln \li ( \f{\de}{2 s} \ri
)$ for some $\ell \in \{1, \cd, s \}$. Define $\bs{\wh{\mu}} = \f{
\ga_{\bs{l}} }{\mathbf{n}_{\bs{l}}}$, where $\bs{l}$ is the index of
stage when the sampling is terminated.
 Then, $\Pr \li \{ \li | \bs{\wh{\mu}} - \mu \ri | < \vep \mu \ri \}
\geq 1 - \de$. \eeT

The proofs of Theorems \ref{Bounded_Mean_Rev_Hoeffding} and \ref{Bounded_Mean_Rev_Massart} can be completed by using techniques similar to that
of Theorems \ref{Bounded_Mean_abs_Hoeffding} and \ref{Bounded_Mean_ABS_Massart}.

 We would like to point out that the
construction of stopping rules proposed in the above two theorems requires essentially no computation.   In the sequel, to reduce the sampling
cost, we propose a computational approach for constructing multistage inverse sampling procedures such that the resultant estimator
$\bs{\wh{\mu}}$ satisfy  $\Pr \li \{ \li | \bs{\wh{\mu}} - \mu \ri | < \vep \mu \ri \} \geq 1 - \de$ for any underlying $\mu \in (0, 1)$, where
$\vep \in (0, 1)$ is the margin of relative error and $\de \in (0, 1)$ is the confidence parameter as before.

 Let $X_i \in [0, 1], \; i = 1,
2, \cd$ be i.i.d. random samples of $X$  with common mean $\mu \in (0, 1)$. Let $\ro > 0$ and $\ze \in (0, \f{1}{\de})$. Let $s$ be the number
of stages. Let $\ga_1, \cd, \ga_s$ be positive numbers such that
\[
2 \ro \li ( \f{1}{\vep} + \ro \ri ) \ln (\ze \de) = \ga_1 < \ga_2 < \cd < \ga_s = 2 \li ( \f{1}{\vep} + \ro \ri )^2 \ln (\ze \de).
\]
For $\ell = 1, \cd, s$, let $\mbf{n}_\ell$ be the random number such that $\sum_{i = 1}^{\mbf{n}_\ell - 1} X_i < \ga_\ell \leq \sum_{i =
1}^{\mbf{n}_\ell} X_i$.  We propose the following stopping rule:

Continue sampling until $\f{\ga_\ell}{ \mbf{n}_\ell } \geq 1 + \ro \vep - \f{\ga_\ell}{ 1 + \ro \vep  } \f{ \vep^2 }{ 2 \ln (\ze \de)  }$.

For $\ell = 1, \cd, s$, define $\wh{\bs{\mu}}_\ell = \f{\ga_\ell}{ \mbf{n}_\ell }$.  The estimator $\wh{\bs{\mu}}$ is taken as
$\wh{\bs{\mu}}_{\bs{l}}$, where $\bs{l}$ is the index of stage at the termination of the sampling process.

To ensure that the coverage probability $\Pr \li \{ \li | \bs{\wh{\mu}} - \mu \ri | < \vep \mu \mid \mu \ri \}$ is no less than the prescribed
confidence level $1 - \de$, we can apply the coverage tuning technique proposed in Section \ref{Computing}.  For this purpose, we need to
develop a computable upper bound for the complementary coverage probability for $\mu$ contained in an interval.  For $\ell = 1, \cd, s$, define
a decision variable $\bs{D}_\ell$ for the $\ell$-th stage such that $\bs{D}_\ell$ assumes value $1$ if {\small $\f{\ga_\ell}{ \mbf{n}_\ell }
\geq 1 + \ro \vep - \f{\ga_\ell}{ 1 + \ro \vep  } \f{ \vep^2 }{ 2 \ln (\ze \de)  }$} and otherwise assumes value $0$.   Consequently, the
stopping rule can be expressed as ``continue the sampling process until $\bs{D}_\ell = 1$ for some $\ell \in \{1, \cd, s\}$''.  Similar to
(\ref{relaxedB}),  we propose to bound the complementary coverage probability $\Pr \li \{ \li | \bs{\wh{\mu}} - \mu \ri | \geq \vep \mu \mid \mu
\ri \}$ with respect to $\mu \in [a, b] \subseteq [0, 1]$ as follows: {\small \be \Pr \li \{ \li | \bs{\wh{\mu}} - \mu \ri | \geq \vep \mu \mid
\mu \ri \}  \leq  \sum_{ \ell = 1 }^s \min \li [ \Pr \{ \bs{D}_{\ell - 1} = 0 \mid \mu \}, \; \Pr \{ | \wh{\bs{\mu}}_\ell - \mu | \geq \vep \mu
\mid \mu \}, \; \Pr \{ \bs{D}_{\ell} = 1 \mid \mu \} \ri ]. \qqu \la{relaxedBBrev} \ee} To bound the complementary coverage probability based on
$(\ref{relaxedBBrev})$, we need to have upper bounds of $\Pr \{ \bs{D}_{\ell} = 0 \mid \mu \}, \; \Pr \{ \bs{D}_{\ell} = 1 \mid \mu \}$ and $\Pr
\{ | \wh{\bs{\mu}}_\ell - \mu | \geq \vep \mu \mid \mu \}$ with respect to $\mu \in [a, b]$ for $\ell = 1, \cd, s$.   This can be accomplished
as follows.

Define
\[ c_\ell = 1 + \ro \vep - \f{\ga_\ell}{ 1 + \ro \vep  } \f{ \vep^2 }{ 2 \ln (\ze \de)  }, \qqu \ell = 1, \cd, s.
\]
Let $[a, b] \subseteq [0, 1]$.  Clearly, $\Pr \{ \bs{D}_0 = 0 \mid \mu \}  = \Pr \{ \bs{D}_s = 1 \mid \mu \} = 1$ for $\mu \in [a, b]$. By
virtue of Chernoff-Hoeffding bounds \cite{Chernoff, Hoeffding} and the definition of the stopping rule, for $\mu \in [a, b]$ and $\ell = 1, \cd,
s - 1$, we have \bee \Pr \{ \bs{D}_{\ell} = 0 \mid \mu \} & = & \Pr \{ \wh{\bs{\mu}}_\ell < c_\ell \mid \mu \}  = \Pr \li \{ \mbf{n}_\ell >
\f{\ga_\ell}{c_\ell} \mid \mu \ri \}\\
&  = & \Pr \li \{ \mbf{n}_\ell \geq \li \lf \f{\ga_\ell}{c_\ell}  \ri \rf + 1 \mid \mu \ri \}
 = \Pr \li \{ \sum_{i = 1}^{m_\ell} X_i \leq \ga_\ell \mid \mu \ri \} \\
& \leq &  \bec \exp ( m_\ell \mscr{M}_{\mrm{B}} ( \f{\ga_\ell}{m_\ell}, a  ) )  & \;
\tx{if} \; \tx{$\f{\ga_\ell}{m_\ell} \leq a$},\\
1 & \; \tx{otherwise} \eec\\
\eee where $m_\ell = \li \lf \f{\ga_\ell}{c_\ell}  \ri \rf + 1$.  Similarly, by virtue of Chernoff-Hoeffding bounds \cite{Chernoff, Hoeffding}
and the definition of the stopping rule, for $\mu \in [a, b]$ and $\ell = 1, \cd, s - 1$, we have \bee \Pr \{ \bs{D}_{\ell} = 1 \mid \mu \} & =
& \Pr \{ \wh{\bs{\mu}}_\ell \geq c_\ell \mid \mu \}  = \Pr \li \{
\mbf{n}_\ell \leq \f{\ga_\ell}{c_\ell} \mid \mu \ri \}\\
&  = & \Pr \li \{ \mbf{n}_\ell \leq \li \lf \f{\ga_\ell}{c_\ell}  \ri \rf  \mid \mu \ri \}
 = \Pr \li \{ \sum_{i = 1}^{m_\ell} X_i \geq \ga_\ell \mid \mu \ri \} \\
& \leq &  \bec \exp ( m_\ell \mscr{M}_{\mrm{B}} ( \f{\ga_\ell}{m_\ell}, b  ) )  & \;
\tx{if} \; \tx{$\f{\ga_\ell}{m_\ell} \geq b$},\\
1 & \; \tx{otherwise} \eec\\
\eee where $m_\ell = \li \lf \f{\ga_\ell}{c_\ell}  \ri \rf$.

From the arguments in \cite[Appendix A.1]{Chen3}, it follows immediately that for $\ell = 1, \cd, s$, \be \la{chen8a} \Pr \{ \wh{\bs{\mu}}_\ell
\geq ( 1 + \vep ) \mu \mid \mu \} \leq \exp ( \ga_\ell \mscr{M}_{\mrm{I}} ( (1 + \vep) \mu, \mu )), \ee where the upper bound on the right side
of (\ref{chen8a}) is decreasing with respect to $\mu \in (0, \f{1}{1 + \vep})$.   From the arguments in \cite[Appendix A.1]{Chen3}, it follows
also that for $\ell = 1, \cd, s$, \be \la{chen8b} \Pr \{ \wh{\bs{\mu}}_\ell  \leq ( 1 - \vep ) \mu \mid \mu \} \leq \exp ( \ga_\ell
\mscr{M}_{\mrm{I}} ( (1 - \ep_\ell) \mu, \mu )), \ee where $\ep_\ell$ is the unique number such that $\f{1}{1 - \vep} = \f{1}{1 - \ep_\ell} +
\f{1}{\ga_\ell}$ and the upper bound on the right side of (\ref{chen8b}) is decreasing with respect to $\mu \in (0, 1)$.

For  $\mu \in [a, b] \subseteq (0, \f{1}{1 + \vep}]$ and $\ell = 1, \cd, s$,  it follows from (\ref{chen8a}) and (\ref{chen8b}) that \bee \Pr \{
| \wh{\bs{\mu}}_\ell - \mu | \geq \vep \mu \mid \mu \} & \leq & \Pr \{ \wh{\bs{\mu}}_\ell \geq \mu +
\vep \mu \mid \mu \} + \Pr \{ \wh{\bs{\mu}}_\ell \leq \mu - \vep \mu \mid \mu \}\\
& \leq & \exp ( \ga_\ell \mscr{M}_{\mrm{I}} ( (1 + \vep) \mu, \mu )) + \exp ( \ga_\ell \mscr{M}_{\mrm{I}} ( (1 - \ep_\ell) \mu, \mu ))\\
& \leq & \exp ( \ga_\ell \mscr{M}_{\mrm{I}} ( (1 + \vep) a, a )) + \exp ( \ga_\ell \mscr{M}_{\mrm{I}} ( (1 - \ep_\ell) a, a )).  \eee

For  $\mu \in [a, b] \subseteq (\f{1}{1 + \vep}, 1)$ and $\ell = 1, \cd, s$,  it follows from (\ref{chen8b}) that \bee \Pr \{ |
\wh{\bs{\mu}}_\ell - \mu | \geq \vep \mu \mid \mu \} & \leq & \Pr \{ \wh{\bs{\mu}}_\ell \leq \mu - \vep \mu \mid \mu \}\\
& \leq &  \exp ( \ga_\ell \mscr{M}_{\mrm{I}} ( (1 - \ep_\ell) \mu, \mu )) \leq \exp ( \ga_\ell \mscr{M}_{\mrm{I}} ( (1 - \ep_\ell) a, a )).
\eee

For  $\mu \in [\f{1}{1 + \vep}, 1)$ and $\ell = 1, \cd, s$,  it follows from (\ref{chen8b}) that \bee \Pr \{ | \wh{\bs{\mu}}_\ell - \mu | \geq
\vep \mu \mid \mu \} & \leq & \Pr \{ \wh{\bs{\mu}}_\ell \geq \mu +
\vep \mu \mid \mu \} + \Pr \{ \wh{\bs{\mu}}_\ell \leq \mu - \vep \mu \mid \mu \}\\
& \leq &  \li (  \f{1}{1 + \vep}  \ri )^{\ga_\ell} + \exp ( \ga_\ell \mscr{M}_{\mrm{I}} ( (1 - \ep_\ell) \mu, \mu ))\\
& \leq &  \li (  \f{1}{1 + \vep}  \ri )^{\ga_\ell} + \exp \li ( \ga_\ell \mscr{M}_{\mrm{I}} \li ( \f{1 - \ep_\ell}{1 + \vep}, \f{1}{1 + \vep}
\ri ) \ri ). \eee

For  $\mu \in (0, 1)$ and $\ell = 1, \cd, s$,  it follows from (\ref{chen8b}) that \bee \Pr \{ | \wh{\bs{\mu}}_\ell - \mu | \geq \vep \mu \mid
\mu \} & \leq & \Pr \{ \wh{\bs{\mu}}_\ell \geq \mu +
\vep \mu \mid \mu \} + \Pr \{ \wh{\bs{\mu}}_\ell \leq \mu - \vep \mu \mid \mu \}\\
& \leq & \exp ( \ga_\ell \mscr{M}_{\mrm{I}} ( (1 + \vep) \mu, \mu )) + \exp ( \ga_\ell \mscr{M}_{\mrm{I}} ( (1 - \ep_\ell) \mu, \mu ))\\
& \leq &  \exp \li ( \ga_\ell \li [ \ln \f{1}{1 + \vep}  +  \f{\vep}{1 + \vep} \ri ] \ri ) + \exp \li ( \ga_\ell \li [ \ln \f{1}{1 - \ep_\ell} -
\f{\ep_\ell}{1 - \ep_\ell} \ri ] \ri ). \eee

Applying the above upper bounds for $\Pr \{ \bs{D}_{\ell} = 0 \mid \mu \}, \; \Pr \{ \bs{D}_{\ell} = 1 \mid \mu \}$ and $\Pr \{ |
\wh{\bs{\mu}}_\ell - \mu | \geq \vep \mu \mid \mu \}$ to (\ref{relaxedBBrev}) gives an upper bound for the complementary coverage probability
$\Pr \{ | \wh{\bs{\mu}} - \mu | \geq \vep \mu  \mid \mu \}$ with respect to $\mu \in [a, b]$.  With such interval bounding technique for the
complementary coverage probability, we can use AMCA to check whether the complementary coverage probability associated with given $\ze$ and
$\ro$ is no greater than $\de$ for any $\mu \in (0, 1)$. Consequently, for any fixed $\ro > 0$, by virtue of the coverage tuning technique
proposed in Section \ref{Computing}, it is possible to determine an appropriate value for $\ze$  such that the coverage probability $\Pr \li \{
\li | \bs{\wh{\mu}} - \mu \ri | < \vep \mu \mid \mu \ri \}$ is no less than the prescribed confidence level $1 - \de$ for any underlying $\mu
\in (0, 1)$.  By choosing different values of $\ro > 0$ and computing the corresponding values of $\ze >0$, the performance of the resultant
stopping rule can be optimized.

In some situations, the cost of sampling operation may be high since
samples are obtained one by one when inverse sampling is involved.
In view of this fact, it is desirable to develop multistage
estimation methods without using inverse sampling.  In contrast to
the multistage inverse sampling schemes described above, our
noninverse multistage sampling schemes have infinitely many stages
and deterministic sample sizes $n_1 < n_2 < n_3 < \cd$. Moreover,
the confidence parameter for the $\ell$-th stage, $\de_\ell$,  is
dependent on $\ell$ such that $\de_\ell = \de$ for $1 \leq \ell \leq
\tau$ and $\de_\ell = \de 2^{\tau - \ell}$ for $\ell > \tau$,  where
$\tau$ is a positive integer.  As before, define $\wh{\bs{\mu}}_\ell
= \f{ \sum_{i = 1}^{n_\ell} X_i } { n_\ell }$ for $\ell = 1, 2,
\cd$. The stopping rule is that sampling is continued until
$\bs{D}_\ell = 1$ for some stage with index $\ell$. Define estimator
$\wh{\bs{\mu}} = \wh{\bs{\mu}}_{\bs{l}}$, where $\bs{l}$ is the
index of stage at which the sampling is terminated. We propose two
types of multistage sampling schemes with different stopping rules
as follows.

\bed

\item [Stopping Rule (i):] For $\ell = 1, 2, \cd$, decision variable
$\bs{D}_\ell$  assumes value $1$ if {\small $\mscr{M}_{\mrm{B}} (
\wh{\bs{\mu}}_\ell, \f{\wh{\bs{\mu}}_\ell}{1 + \vep} ) \leq \f{ \ln
( \ze \de_\ell) } { n_\ell }$}; and assumes value $0$ otherwise.

\item [Stopping Rule (ii):] For $\ell = 1, 2, \cd$, decision variable
$\bs{D}_\ell$ assumes value $1$  if {\small
\[
\wh{\bs{\mu}}_\ell \geq \f{ 6(1 + \vep) (3 + \vep) \ln (\ze
\de_\ell) } { 2 (3 + \vep)^2 \ln (\ze \de_\ell) - 9 n_\ell \vep^2};
\]} and assumes value $0$ otherwise.

\eed

Stopping rule (i) is derived by virtue of Chernoff-Hoeffding bounds of
the CDF $\&$ CCDF of $\wh{\bs{\mu}}_\ell$.  Stopping rule (ii) is derived
by virtue of Massart's inequality of the CDF $\&$ CCDF of $\wh{\bs{\mu}}_\ell$.

For both types of multistage sampling schemes described above, we have established the following theorem.

\beT \la{Bounded_noinverse}  The follows statements hold true:

(I): $\Pr \{ \mbf{n} < \iy \} = 1$ for any $\mu \in (0, 1)$ provided
that $\inf_{\ell
> 0} \f{n_{\ell + 1}}{n_\ell} > 1$.

(II): $\bb{E} [ \mbf{n} ] < \iy$  for any $\mu \in (0, 1)$ provided
that  $1 < \inf_{\ell
> 0} \f{n_{\ell + 1}}{n_\ell} \leq \sup_{\ell
> 0} \f{n_{\ell + 1}}{n_\ell} < \iy$.

(III): {\small $\Pr \li \{ \li | \f{ \wh{\bs{\mu}} - \mu } { \mu }
\ri | < \vep \mid \mu \ri \} \geq 1 - \de$} for any $\mu \in (0, 1)$
provided that $ \ze \leq \f{1}{2 (\tau + 1 )}$.

\eeT

The proof of Theorem \ref{Bounded_noinverse} can be accomplished by using techniques similar to that of Theorem \ref
{Bino_Rev_noninverse_Chernoff}, which is given at Appendix \ref{App_Bino_Rev_noninverse_Chernoff}.

\subsection{Control of Absolute  and Relative Errors}

In this subsection, we consider the multistage estimation of the
mean of the bounded variable with a mixed error criterion.
Specifically, we wish to construct a multistage sampling scheme and
its associated estimator $\wh{\bs{\mu}}$ for $\mu = \bb{E} [ X]$
such that $\Pr \{ | \wh{\bs{\mu}} - \mu | < \vep_a, \; |
\wh{\bs{\mu}} - \mu | < \vep_r \mu \} > 1 - \de$.  In the special
case that the variable $X$ is bounded in interval $[0, 1]$, our
multistage sampling schemes and their properties are described by
the following theorems.

\beT \la{Bounded_Mean_mix_Hoeffding} Let $0 < \vep_a < \f{35}{94}$
and {\small $\f{70 \vep_a}{35 - 24 \vep_a } < \vep_r < 1$}.  Let
$n_1 < n_2 < \cd < n_s$ be a sequence of sample sizes such that
{\small $n_s \geq \f{ \ln (2s \sh \de) } { \mscr{M}_{\mrm{B}} (
\f{\vep_a}{\vep_r} + \vep_a, \f{\vep_a}{\vep_r} )  }$}. Define
$\wh{\bs{\mu}}_\ell = \f{ \sum_{i=1}^{n_\ell} X_i }{n_\ell}, \;
\mscr{L} (\wh{\bs{\mu}}_\ell) = \min \{ \wh{\bs{\mu}}_\ell - \vep_a,
\; \f{\wh{\bs{\mu}}_\ell}{ 1 + \vep_r } \}$ and $\mscr{U} (
\wh{\bs{\mu}}_\ell ) = \max \{ \wh{\bs{\mu}}_\ell + \vep_a, \;
\f{\wh{\bs{\mu}}_\ell}{ 1 - \vep_r } \}$ for $\ell = 1, \cd, s$.
Suppose that sampling is continued until {\small  $\max \{
\mscr{M}_{\mrm{B}} (\wh{\bs{\mu}}_\ell, \mscr{L}
(\wh{\bs{\mu}}_\ell) ), \; \mscr{M}_{\mrm{B}} (\wh{\bs{\mu}}_\ell,
\mscr{U} (\wh{\bs{\mu}}_\ell) ) \} \leq \f{1}{n_\ell} \ln \li (
\f{\de}{2 s} \ri )$}. Define $\bs{\wh{\mu}} =
\f{\sum_{i=1}^{\mathbf{n}} X_i}{\mathbf{n}}$, where $\mathbf{n}$ is
the sample size when the sampling is terminated. Then, $\Pr  \{ |
\bs{\wh{\mu}} - \mu | < \vep_a \; \tx{or} \;  | \bs{\wh{\mu}} - \mu
| < \vep_r \mu \} \geq 1 - \de$. \eeT

\bsk

See Appendix \ref{App_Bounded_Mean_mix_Hoeffding} for a proof.

\beT \la{Bounded_Mean_mix_Massart} Let $0 < \vep_a < \f{3}{8}$ and
{\small $\f{6 \vep_a}{3 - 2 \vep_a } < \vep_r < 1$}. Let $n_1 < n_2
< \cd < n_s$ be a sequence of sample sizes such that {\small $n_s
\geq 2  \li ( \f{1}{ \vep_r } + \frac{1} {3} \ri ) \li ( \f{1}{
\vep_a } - \f{1}{ \vep_r } - \frac{1} {3} \ri ) \ln \li (
\f{2s}{\de} \ri )$}.  Define $\wh{\bs{\mu}}_\ell = \f{
\sum_{i=1}^{n_\ell} X_i }{n_\ell}$ and {\small \[ \bs{D}_\ell = \bec
0 & \mrm{for} \; \f{1}{2} - \f{2}{3} \vep_a - \sq{ \f{1}{4} + \f{
n_\ell \vep_a^2 } {2 \ln (\ze \de) } } < \wh{\bs{\mu}}_\ell < \f{
6(1 - \vep_r) (3 - \vep_r) \ln (\ze
\de) } { 2 (3 - \vep_r)^2 \ln (\ze \de) - 9 n_\ell \vep_r^2} \; \mrm{or}\\
  &  \qu \;\; \f{1}{2} + \f{2}{3} \vep_a - \sq{ \f{1}{4} + \f{
n_\ell \vep_a^2 } {2 \ln (\ze \de) } } < \wh{\bs{\mu}}_\ell < \f{
6(1 + \vep_r) (3 + \vep_r) \ln (\ze \de) } { 2 (3 +
\vep_r)^2 \ln (\ze \de) - 9 n_\ell  \vep_r^2},\\
1 & \mrm{else}
 \eec
\]}
for $\ell = 1, \cd, s$, where $\ze = \f{1}{2 s}$.   Suppose that sampling is continued until $\bs{D}_\ell = 1$ for some $\ell \in \{1, \cd, s
\}$. Define $\bs{\wh{\mu}} = \f{\sum_{i=1}^{\mathbf{n}} X_i}{\mathbf{n}}$, where $\mathbf{n}$ is the sample size when the sampling is
terminated.
 Then, $\Pr \{ | \bs{\wh{\mu}} - \mu | < \vep_a \;
 \tx{or} \; | \bs{\wh{\mu}} - \mu | < \vep_r \mu \}
\geq 1 - \de$. \eeT

\bsk

See Appendix \ref{App_Bounded_Mean_mix_Massart} for a proof.

In the general case that $X$ is a random variable bounded in $[a,
b]$, it is useful to estimate the mean $\mu = \bb{E} [X]$  with a
mixed criterion based on i.i.d. samples of $X$ and prior information
that $a \leq \udl{\mu} \leq \mu \leq \ovl{\mu} \leq b$. Let $\vep_a
\in (0, \iy)$ and $\vep_r \in (0, 1)$ be margins of absolute and
relative errors. To describe our multistage estimation methods,
define functions $v(z) = \f{z - a}{b - a}$, {\small \bee & &  g(z) =
\bec \f{1}{b - a}  \min \li \{ z - \vep_a, \; \f{z}{ 1 + \vep_r }
\ri \} - \f{a}{b - a} & \tx{if
$\udl{\mu} > 0$},\\
\f{1}{b - a}  \min \li \{ z - \vep_a, \; \f{z}{ 1 - \vep_r } \ri \}
- \f{a}{b - a} & \tx{if
$\ovl{\mu} < 0$},\\
\f{1}{b - a}  \min \li \{ z - \vep_a, \; \f{z}{ 1 + \mrm{sgn} (z)
\vep_r } \ri \} - \f{a}{b - a} & \tx{if $0 \in [\udl{\mu},
\ovl{\mu}]$} \eec\\
&  & h(z) =  \bec \f{1}{b - a}  \max \li \{ z + \vep_a, \; \f{z}{ 1
- \vep_r } \ri \} - \f{a}{b - a} & \tx{if
$\udl{\mu} > 0$},\\
\f{1}{b - a}  \max \li \{ z + \vep_a, \; \f{z}{ 1 + \vep_r } \ri \}
- \f{a}{b - a} & \tx{if
$\ovl{\mu} < 0$},\\
\f{1}{b - a}  \max \li \{ z + \vep_a, \; \f{z}{ 1 - \mrm{sgn} (z)
\vep_r } \ri \} - \f{a}{b - a} & \tx{if $0 \in [\udl{\mu},
\ovl{\mu}]$} \eec \eee} and \bee &  & \mcal{W}_{\mrm{B}} (z) = \max
\li \{ \mscr{M}_{\mrm{B}} \li ( v(z),
g(z) \ri ), \; \mscr{M}_{\mrm{B}} \li ( v(z),  h(z) \ri ) \ri \},\\
&  & \mcal{W}(z) = \max \li \{ \mscr{M} \li ( v(z), g(z) \ri ), \;
\mscr{M} \li ( v(z), h(z) \ri ) \ri \} \eee for $z \in [a, b]$. By
virtue of such functions and Theorem \ref{Monotone_second_CI_ST}, we
have established multistage sampling schemes as described by
Theorems \ref{Bounded_mix_general_Hoeffding} and
\ref{Bounded_Mean_mix_general_Massart}.

\beT \la{Bounded_mix_general_Hoeffding} Let $n_1 < n_2 < \cd < n_s$
be a sequence of sample sizes such that {\small $n_s \geq  \f{\ln
\f{\de}{2 s} } { \max_{z \in [a, b]} \mcal{W}_{\mrm{B}} (z) }$}.
Define $\wh{\bs{\mu}}_\ell = \f{ \sum_{i=1}^{n_\ell} X_i }{n_\ell}$
for $\ell = 1, \cd, s$.  Suppose that sampling is continued until
$\mcal{W}_{\mrm{B}} (\wh{\bs{\mu}}_\ell) \leq \f{1}{n_\ell} \ln
\f{\de}{2s}$ for some $\ell \in \{1, \cd, s \}$. Define
$\bs{\wh{\mu}} = \f{\sum_{i=1}^{\mathbf{n}} X_i}{\mathbf{n}}$, where
$\mathbf{n}$ is the sample size when the sampling is terminated.
 Then, $\Pr \{ | \bs{\wh{\mu}} - \mu | < \vep_a \; \tx{or} \;
 | \bs{\wh{\mu}} - \mu | < \vep_r |\mu| \}
\geq 1 - \de$ for any $\mu \in [\udl{\mu}, \ovl{\mu}]$. \eeT

For efficiency, the sample sizes in the sampling scheme described in
Theorem \ref{Bounded_mix_general_Hoeffding} can be chosen as a
geometric sequence $n_1, \cd, n_s$.  The minimum sample size $n_1$
can be chosen as the minimum integer no less than $ \f{\ln \f{\de}{2
s} } { \min_{z \in [a, b]} \mcal{W}_{\mrm{B}} (z) }$. The maximum
sample size $n_s$ can be chosen as the minimum integer no less than
$ \f{\ln \f{\de}{2 s} } { \max_{z \in [a, b]} \mcal{W}_{\mrm{B}} (z)
}$.

\beT \la{Bounded_Mean_mix_general_Massart} Let $n_1 < n_2 < \cd <
n_s$ be a sequence of sample sizes such that {\small $n_s \geq
\f{\ln \f{\de}{2 s} } { \max_{z \in [a, b]} \mcal{W} (z) }$}. Define
$\wh{\bs{\mu}}_\ell = \f{ \sum_{i=1}^{n_\ell} X_i }{n_\ell}$ for
$\ell = 1, \cd, s$.  Suppose that sampling is continued until
$\mcal{W} (\wh{\bs{\mu}}_\ell) \leq \f{1}{n_\ell} \ln \f{\de}{2s}$
for some $\ell \in \{1, \cd, s \}$. Define $\bs{\wh{\mu}} =
\f{\sum_{i=1}^{\mathbf{n}} X_i}{\mathbf{n}}$, where $\mathbf{n}$ is
the sample size when the sampling is terminated.
 Then, $\Pr \li \{ \li | \bs{\wh{\mu}} - \mu \ri | < \vep_a \;
 \tx{or} \; \li | \bs{\wh{\mu}} - \mu \ri | < \vep_r |\mu| \ri \}
\geq 1 - \de$ for any $\mu \in [\udl{\mu}, \ovl{\mu}]$. \eeT

In the sampling scheme described in Theorem
\ref{Bounded_Mean_mix_general_Massart}, the minimum  sample size
$n_1$ can be chosen as the minimum integer no less than $ \f{\ln
\f{\de}{2 s} } { \min_{z \in [a, b]} \mcal{W} (z) }$. The maximum
sample size $n_s$ can be chosen as the minimum integer no less than
$ \f{\ln \f{\de}{2 s} } { \max_{z \in [a, b]} \mcal{W} (z) }$.

It should be noted that the minimum and maximum of
$\mcal{W}_{\mrm{B}} (z)$ and $\mcal{W} (z)$ over $[a, b]$ can be
exactly computed by using {\it Branch and Bound} method. For this
purpose, we need to have the upper and lower bounds of
$\mcal{W}_{\mrm{B}} (z)$ and $\mcal{W} (z)$ for $z$ in a subset
$[\udl{z}, \ovl{z}]$ of $[a, b]$. Note that one can partition $[a,
b]$ as subintervals such that the numbers in each subinterval have
the same sign. Without loss of generality, assume that $z$ has the
same sign over $[\udl{z}, \ovl{z}]$. Then,
\[
g(\udl{z}) \leq g(z) \leq g(\ovl{z}),  \qqu h(\udl{z}) \leq h(z)
\leq h(\ovl{z}), \qqu v(\udl{z})  \leq v(z) \leq v(\ovl{z})
\]
for $z \in [\udl{z}, \ovl{z}]$. Since $g(z) + \f{\vep_a}{b - a} \leq
v(z) \leq h(z) - \f{\vep_a}{b - a}$, we have that $g(z) \leq
g(\ovl{z}) \leq v(\udl{z}) \leq v(z) \leq v(\ovl{z}) \leq h(\udl{z})
\leq h(z)$ provided that the subinterval $[\udl{z}, \ovl{z}]$ is
sufficiently narrow.  It follows that \bee &  & \mcal{W}_{\mrm{B}}
(z) \leq \max \li \{ \mscr{M}_{\mrm{B}} \li ( v(\udl{z}), g(\ovl{z})
\ri ), \; \mscr{M}_{\mrm{B}} \li
( v(\ovl{z}), h(\udl{z}) \ri ) \ri \},\\
&  & \mcal{W}_{\mrm{B}} (z) \geq \max \li \{ \mscr{M}_{\mrm{B}} \li
( v(\ovl{z}), g(\udl{z}) \ri ), \; \mscr{M}_{\mrm{B}} \li
( v(\udl{z}), h(\ovl{z}) \ri ) \ri \},\\
&  & \mcal{W} (z) \leq \max \li \{ \mscr{M} \li ( v(\udl{z}) ,
g(\ovl{z}) \ri ),
\; \mscr{M} \li ( v(\ovl{z}), h(\udl{z}) \ri ) \ri \},\\
&  & \mcal{W} (z) \geq \max \li \{ \mscr{M} \li ( v(\ovl{z}),
g(\udl{z}) \ri ), \; \mscr{M} \li ( v(\udl{z}), h(\ovl{z}) \ri ) \ri
\} \eee for $z \in [\udl{z}, \ovl{z}]$, where $[\udl{z}, \ovl{z}]$
is a subset of $[a, b]$ such that $g(\ovl{z}) \leq v(\udl{z}) \leq
v(\ovl{z}) \leq h(\udl{z})$ and that the numbers in $[\udl{z},
\ovl{z}]$ have the same sign.

We would like to point out that the results and methods presented in Sections 7.1, 7.2 and 7.3 are still valid if the assumption that $X_1, X_2,
\cd$ are i.i.d. samples of $X$ are relaxed as follows:

\bee &  & X_k  \; \tx{is bounded in the same interval as $X$ almost surely for $k \in
\bb{N}$},  \\
&  & \bb{E} [ X_k \mid \mscr{F}_{k-1} ] = \bb{E} [ X] = \mu  \qu \tx{almost surely for $k \in \bb{N}$},  \eee where $\{ \mscr{F}_k, \; k = 0, 1,
\cd, \iy \}$ is a sequence of $\si$-subalgebra such that $\{ \emptyset, \Om \} = \mscr{F}_0 \subset \mscr{F}_1 \subset \mscr{F}_2 \subset \cd
\subset \mscr{F}$, with $\mscr{F}_k$ being generated by $X_1, \cd, X_k$.

\subsection{Using the Link between Binomial and Bounded Variables}

Recently, Chen \cite{link_chen} has discovered the following
inherent connection between a binomial parameter and the mean of a
bounded variable.

\beT \la{Link} Let $X$ be a random variable bounded in $[0, 1]$. Let
$U$ a random variable uniformly distributed over $[0, 1]$. Suppose
$X$ and $U$ are independent. Then, $\bb{E} [X] = \Pr \{ X \geq U
\}$. \eeT

To see why Theorem \ref{Link} reveals a relationship between the
mean of a bounded variable and a binomial parameter, we define
\[
Y = \bec 1 & \tx{for} \; X \geq U,\\
0 & \tx{otherwise}. \eec
\]
Then, by Theorem \ref{Link}, we have $\Pr \{ Y = 1 \} = 1 - \Pr \{ Y
= 0 \} = \bb{E} [X]$. This implies that $Y$ is a Bernoulli random
variable and $\bb{E} [X]$ is actually a binomial parameter.  For a
sequence of i.i.d. random samples $X_1, X_2, \cd$ of bounded
variable $X$ and a sequence of i.i.d. random samples $U_1, U_2, \cd$
of uniform variable $U$ such that that $X_i$ is independent with
$U_i$ for all $i$, we can define a sequence of i.i.d. random samples
$Y_1, Y_2, \cd$  of Bernoulli random variable $Y$ by
\[
Y_i = \bec 1 & \tx{for} \; Y_i \geq U_i,\\
0 & \tx{otherwise}. \eec
\]
As a consequence, the techniques of estimating a binomial parameter
can be useful for estimating the mean of a bounded variable.

\section{Estimation of Poisson Parameters}

In this section, we shall consider the multistage estimation of the
mean, $\lm$, of a Poisson random variable $X$ based on its i.i.d.
random samples $X_1, X_2, \cd $.

For $\ell = 1, 2, \cd$, define $K_\ell = \sum_{i = 1}^{n_\ell} X_i,
\; \wh{\bs{\lm}}_\ell = \f{K_\ell} {n_\ell }$, where $n_\ell$ is
deterministic and stands for the sample size at the $\ell$-th stage.
As described in the general structure of our multistage estimation
framework, the stopping rule is that sampling is continued until
$\bs{D}_\ell = 1$ for some $\ell \in \{1, \cd, s\}$. Define
estimator $\wh{\bs{\lm}} = \wh{\bs{\lm}}_{\bs{l}}$, where $\bs{l}$
is the index of stage at which the sampling is terminated. Clearly,
the sample number at the completion of sampling is $\mbf{n} =
n_{\bs{l}}$.

\subsection{Control of Absolute Error}

In this subsection, we shall focus on the design of multistage
sampling schemes for estimating the Poisson parameter $\lm$ with an
absolute error criterion.  Specifically, for $\vep > 0$, we wish to
construct a multistage sampling scheme and its associated estimator
$\wh{\bs{\lm}}$ for $\lm$ such that $\Pr \{ | \wh{\bs{\lm}} - \lm |
< \vep \mid \lm \} > 1 - \de$ for any $\lm \in (0, \iy)$.  As will
be seen below, our multistage sampling procedures have infinitely
many stages and deterministic sample sizes $n_1 < n_2 < n_3 < \cd$.
Moreover, the confidence parameter for the $\ell$-th stage,
$\de_\ell$,  is dependent on $\ell$ such that $\de_\ell = \de$ for
$1 \leq \ell \leq \tau$ and $\de_\ell = \de 2^{\tau - \ell}$ for
$\ell > \tau$,  where $\tau$ is a positive integer.

\subsubsection{Stopping Rule from CDF $\&$ CCDF}

By virtue of the CDF $\&$ CCDF of $\wh{\bs{\lm}}_\ell$, we propose a
class of multistage sampling schemes as follows.

\beT \la{Pos_abs_CDF} Suppose that, for $\ell = 1, 2, \cd$, decision
variable $\bs{D}_\ell$ assumes values $1$ if {\small
$F_{\wh{\bs{\lm}}_\ell} ( \wh{\bs{\lm}}_\ell, \wh{\bs{\lm}}_\ell +
\vep  ) \leq \ze \de_\ell, \; G_{\wh{\bs{\lm}}_\ell} (
\wh{\bs{\lm}}_\ell, \wh{\bs{\lm}}_\ell - \vep  ) \leq \ze
\de_\ell$}; and assumes $0$ otherwise. The following statements hold
true.

(I): $\Pr \{ \mbf{n} < \iy \} = 1$  provided that $\inf_{\ell > 0}
\f{n_{\ell + 1}}{n_\ell} > 1$.

(II): $\bb{E} [ \mbf{n} ] < \iy$ provided that $1 < \inf_{\ell > 0}
\f{n_{\ell + 1}}{n_\ell} \leq \sup_{\ell > 0} \f{n_{\ell +
1}}{n_\ell} < \iy$.

(III): {\small $\Pr  \{  | \wh{\bs{\lm}} - \lm  | < \vep \mid \lm
 \} \geq 1 - \de$} for any $\lm  > 0$ provided that $
\ze \leq \f{1}{ 2 (\tau + 1 ) }$.

(IV):  Let $0 < \eta < \ze \de$ and {\small ${\ell^\star} = \tau + 1
+ \li \lc \f{ \ln ( \ze \de \sh \eta )  } { \ln 2 } \ri \rc$}. Then,
$\Pr  \{ | \wh{\bs{\lm}} - \lm | \geq \vep  \mid \lm  \} < \de$ for
any $\lm \in ( \ovl{\lm}, \iy)$, where $\ovl{\lm}$ is a number such
that {\small $\ovl{\lm} > z_\ell, \; \; \ell = 1, \cd,
{\ell^\star}$} and that $\sum_{\ell = 1}^{\ell^\star} \exp ( n_\ell
\mscr{M}_{\mrm{P}} ( z_\ell, \ovl{\lm}) ) < \de - \eta$ with $z_\ell
= \min \{  z \in I_{\wh{\bs{\lm}}_\ell}: F_{ \wh{\bs{\lm}}_\ell }
\li (z, z +  \vep \ri ) > \ze \de_\ell \; \tx{or} \; G_{
\wh{\bs{\lm}}_\ell } \li ( z, z - \vep \ri ) > \ze \de_\ell \}$,
where $I_{\wh{\bs{\lm}}_\ell}$ represents the support of
$\wh{\bs{\lm}}_\ell$, for $\ell = 1, 2, \cd$.  Moreover, {\small
\bee & & \Pr \li \{ b \leq \wh{\bs{\lm}} -  \vep, \; \bs{l} \leq
\ell^\star \mid a \ri \} \leq \Pr \li \{ \lm \leq \wh{\bs{\lm}} -
\vep \mid \lm \ri \} \leq \f{\eta}{2} + \Pr
\li \{ a \leq \wh{\bs{\lm}} - \vep, \; \bs{l} \leq \ell^\star \mid b \ri \}, \\
&   & \Pr \li \{ a \geq \wh{\bs{\lm}} + \vep, \; \bs{l} \leq
\ell^\star \mid b \ri \}  \leq \Pr \li \{ \lm \geq \wh{\bs{\lm}} +
\vep \mid \lm   \ri \} \leq \f{\eta}{2} + \Pr \li \{ b \geq
\wh{\bs{\lm}} + \vep, \; \bs{l} \leq \ell^\star \mid a \ri \} \eee}
for any $\lm \in [a, b]$, where $a$ and $b$ are numbers such that $0
< b <  a + \vep$.

(V): Let the sample sizes of the multistage sampling scheme be a
sequence $n_\ell = \li \lc m \ga^{\ell - 1} \ri \rc, \; \ell = 1, 2,
\cd$, where $\ga \geq 1 + \f{1}{m} > 1$. Let $\ep > 0, \; 0 < \eta <
1$ and $c = - \mscr{M}_{\mrm{P}} ( \f{\lm}{\eta}, \lm )$. Let $\ka$
be an integer such that {\small $\ka
> \max \li \{ \tau, \;  \f{1}{\ln \ga} \ln \li ( \f{1}{ c m } \ri
) + 1, \; \f{1}{\ln \ga} \ln \li ( \f{1}{ c m } \ln \f{\ga}{c \ep}
\ri ) + 1, \; \tau + \f{1}{\ga - 1} + \f{ \ln (\ze \de) } { \ln 2 }
\ri \} $} and {\small $\mscr{M}_{\mrm{P}}  ( \f{\lm}{\eta},
\f{\lm}{\eta}  + \vep  ) < \f{ \ln (\ze \de_\ka )}{n_\ka}$}. Then,
$\bb{E} [\mbf{n}] < \ep + n_1 + \sum_{\ell = 1}^{\ka} (n_{\ell + 1}
- n_\ell) \Pr \{ \bs{l}
> \ell \}$.

\eeT

The proof of Theorem \ref{Pos_abs_CDF} is similar to that of Theorem \ref{Pos_abs_Chernoff}, which is given at Appendix \ref
{App_Pos_abs_Chernoff}.

\subsubsection{Stopping Rule from Chernoff Bounds}

By virtue of Chernoff bounds of the CDF $\&$ CCDF of
$\wh{\bs{\lm}}_\ell$, we propose a class of multistage sampling
schemes as follows.

\beT \la{Pos_abs_Chernoff}

Suppose that, for $\ell = 1, 2, \cd$, decision variable
$\bs{D}_\ell$ assumes values $1$ if {\small $\mscr{M}_{\mrm{P}}
(\wh{\bs{\lm}}_\ell, \wh{\bs{\lm}}_\ell + \vep )  \leq \f{ \ln ( \ze
\de_\ell ) } { n_\ell }$}; and assumes $0$ otherwise. The following
statements hold true.

(I): $\Pr \{ \mbf{n} < \iy \} = 1$  provided that $\inf_{\ell > 0}
\f{n_{\ell + 1}}{n_\ell} > 1$.

(II): $\bb{E} [ \mbf{n} ] < \iy$ provided that $1 < \inf_{\ell > 0}
\f{n_{\ell + 1}}{n_\ell} \leq \sup_{\ell > 0} \f{n_{\ell +
1}}{n_\ell} < \iy$.

(III): {\small $\Pr \{  | \wh{\bs{\lm}} - \lm  | < \vep \mid \lm \}
\geq 1 - \de$} for any $\lm  > 0$ provided that $\ze \leq \f{1}{ 2
(\tau + 1 ) }$.

(IV):  Let $0 < \eta < \ze \de$ and {\small ${\ell^\star} = \tau + 1
+ \li \lc \f{ \ln ( \ze \de \sh \eta )  } { \ln 2 } \ri \rc$}. Then,
$\Pr \{  | \wh{\bs{\lm}} - \lm  | \geq \vep  \mid \lm  \} < \de$ for
any $\lm \in ( \ovl{\lm}, \iy)$, where $\ovl{\lm}$ is a number such
that {\small $\ovl{\lm} > z_\ell, \; \; \ell = \tau, \cd,
{\ell^\star}$} and that $\sum_{\ell = 1}^{\ell^\star} \exp ( n_\ell
\mscr{M}_{\mrm{P}} ( z_\ell, \ovl{\lm}) ) < \de - \eta$ with
$z_\ell$ satisfying {\small $\mscr{M}_{\mrm{P}} \li (z_\ell, z_\ell
+ \vep \ri  )  = \f{ \ln ( \ze \de_\ell ) } { n_\ell }$} for $\ell =
1, 2, \cd$.  Moreover, {\small \bee &  & \Pr \li \{ b \leq
\wh{\bs{\lm}} -  \vep, \; \bs{l} \leq \ell^\star \mid a \ri \} \leq
\Pr \li \{ \lm \leq \wh{\bs{\lm}} - \vep \mid \lm \ri \} \leq
\f{\eta}{2} + \Pr
\li \{ a \leq \wh{\bs{\lm}} - \vep, \; \bs{l} \leq \ell^\star \mid b \ri \}, \\
&   & \Pr \li \{ a \geq \wh{\bs{\lm}} + \vep, \; \bs{l} \leq
\ell^\star \mid b \ri \}  \leq \Pr \li \{ \lm \geq \wh{\bs{\lm}} +
\vep \mid \lm   \ri \} \leq \f{\eta}{2} + \Pr \li \{ b \geq
\wh{\bs{\lm}} + \vep, \; \bs{l} \leq \ell^\star \mid a \ri \} \eee}
for any $\lm \in [a, b]$, where $a$ and $b$ are numbers such that $0
< b <  a + \vep$.

(V): Let the sample sizes of the multistage sampling scheme be a
sequence $n_\ell = \li \lc m \ga^{\ell - 1} \ri \rc, \; \ell = 1, 2,
\cd$, where $\ga \geq 1 + \f{1}{m} > 1$. Let $\ep > 0, \; 0 < \eta <
1$ and $c = - \mscr{M}_{\mrm{P}} ( \f{\lm}{\eta}, \lm )$. Let $\ka$
be an integer such that {\small $\ka
> \max \li \{ \tau, \;  \f{1}{\ln \ga} \ln \li ( \f{1}{ c m } \ri
) + 1, \; \f{1}{\ln \ga} \ln \li ( \f{1}{ c m } \ln \f{\ga}{c \ep}
\ri ) + 1, \; \tau + \f{1}{\ga - 1} + \f{ \ln (\ze \de) } { \ln 2 }
\ri \} $} and {\small $\mscr{M}_{\mrm{P}}  ( \f{\lm}{\eta},
\f{\lm}{\eta}  + \vep ) < \f{ \ln (\ze \de_\ka )}{n_\ka}$}. Then,
$\bb{E} [\mbf{n}] < \ep + n_1 + \sum_{\ell = 1}^{\ka} (n_{\ell + 1}
- n_\ell) \Pr \{ \bs{l}
> \ell \}$.

\eeT

See Appendix \ref{App_Pos_abs_Chernoff} for a proof.

\subsubsection{Asymptotic Analysis of Multistage Sampling Schemes}

In this subsection, we shall focus on the asymptotic analysis of the
multistage sampling schemes which follow stopping rules derived from
Chernoff bounds of CDF $\&$ CCDF of $\wh{\bs{\lm}}_\ell$ as
described in Theorem \ref{Pos_abs_Chernoff}.

Let $\lm^* > 0$.  We assume that the sample sizes $n_1, n_2, \cd$ are chosen as the ascending arrangement of all distinct elements of the set
\be \la{defposabs} \li \{ \li \lc  \f{ C_{\tau - \ell} \; \ln (\ze \de) }{ \mscr{M}_{\mrm{P}} (\lm^*, \lm^* + \vep ) } \ri \rc : \ell = 1, 2,
\cd \ri \}, \ee where $\tau$ is the maximum integer such that $\f{ C_{\tau - 1} \; \ln (\ze \de) }{ \mscr{M}_{\mrm{P}} (\lm^*, \lm^* + \vep ) }
\geq \f{ \ln \f{1}{\ze \de} }{ \vep }$, i.e., $C_{\tau - 1} \geq - \f{ \mscr{M}_{\mrm{P}} (\lm^*, \lm^* + \vep ) } { \vep }$. With regard to the
asymptotic performance of the sampling scheme, we have

\beT \la{Pos_Asp_Analysis_Abs} Let {\small $\mcal{N}_{\mrm{a}} (\lm, \vep) = \f{ \ln ( \ze \de  ) } { \mscr{M}_{\mrm{P}} ( \lm , \lm + \vep) }
$}. Let $\mcal{N}_{\mrm{f}} (\lm, \vep)$ be the minimum sample number $n$ such that {\small $\Pr \{ | \f{\sum_{i = 1}^n X_i}{n} - \lm  | < \vep
\mid \lm  \} > 1 - \ze \de$} for a fixed-size sampling procedure. Let $j_\lm$ be the largest integer $j$ such that $C_{j} \geq \f{\lm}{\lm^*}$.
Let {\small $\nu = \f{2}{3} (1 - \f{\lm} {\lm^*} ), \; d = \sq{ 2 \ln \f{1}{\ze \de} }$} and {\small $\ka_\lm = \f{\lm^* }{\lm} C_{j_\lm} $}.
Let {\small $\ro_\lm = \f{\lm^* } { \lm } C_{j_\lm - 1} - 1$} if $\ka_\lm = 1$ and $\ro_\lm = \ka_\lm - 1$ otherwise. For $\lm \in (0, \lm^*)$,
the following statements hold true:

(I): {\small $\Pr \li \{  1 \leq \limsup_{\vep \to 0} \f{ \mbf{n} } { \mcal{N}_{\mrm{a}} (\lm, \vep) } \leq 1 + \ro_\lm  \ri \} = 1$}.
Specially, {\small $\Pr \li \{  \lim_{\vep \to 0} \f{ \mbf{n} } { \mcal{N}_{\mrm{a}} (\lm, \vep) } = \ka_\lm  \ri \} = 1$} if $\ka_\lm
> 1$.

(II): {\small $\lim_{\vep \to 0} \f{ \bb{E} [ \mbf{n} ] } { \mcal{N}_{\mrm{f}} (\lm, \vep)} = \li ( \f{ d } { \mcal{Z}_{\ze \de} } \ri )^2
\times \lim_{\vep \to 0} \f{ \bb{E} [ \mbf{n} ] } { \mcal{N}_{\mrm{a}} (\lm, \vep) }$}, where
\[ \lim_{\vep \to 0} \f{ \bb{E} [ \mbf{n} ] }  {
\mcal{N}_{\mrm{a}} (\lm, \vep) } = \bec \ka_\lm
& \tx{if} \; \ka_\lm > 1,\\
1  + \ro_\lm \Phi (\nu d ) &  \tx{otherwise} \eec
\]
and $1 \leq \lim_{\vep \to 0} \f{ \bb{E} [ \mbf{n} ] }  { \mcal{N}_{\mrm{a}} (\lm, \vep) } \leq 1 + \ro_\lm$.

(III): If $\ka_\lm > 1$, then $\lim_{\vep \to 0} \Pr \{ | \wh{\bs{\lm}} - \lm | < \vep \} = 2 \Phi \li ( d \sq{\ka_\lm} \ri ) - 1 > 2 \Phi ( d )
- 1 > 1 - 2 \ze \de$. Otherwise, $\Phi \li ( d  \ri )  + \Phi \li ( d \sq{1 + \ro_\lm}  \ri )  - 1  > \lim_{\vep \to 0} \Pr \{ | \wh{\bs{\lm}} -
\lm | < \vep \} = 1 + \Phi(d) - \Phi(\nu d) - \Psi (\ro_\lm, \nu, d) > \Phi \li ( d \ri ) + 2 \Phi \li ( d \sq{1 + \ro_\lm}  \ri )  - 2  > 1 - 3
\ze \de$.

 \eeT

\bsk

See Appendix \ref{App_Pos_Asp_Analysis_Abs} for a proof.

\subsection{Control of Relative Error}

In this subsection, we shall focus on the design of multistage
sampling schemes for estimating the Poisson parameter $\lm$ with a
relative error criterion.  Specifically, for $\vep \in (0, 1)$, we
wish to construct a multistage sampling scheme and its associated
estimator $\wh{\bs{\lm}}$ for $\lm$ such that $\Pr \{ |
\wh{\bs{\lm}} - \lm | < \vep \lm \mid \lm \} > 1 - \de$ for any $\lm
\in (0, \iy)$.  As will be seen below, our multistage sampling
procedures have infinitely many stages and deterministic sample
sizes $n_1 < n_2 < n_3 < \cd$. Moreover, the confidence parameter
for the $\ell$-th stage, $\de_\ell$,  is dependent on $\ell$ such
that $\de_\ell = \de$ for $1 \leq \ell \leq \tau$ and $\de_\ell =
\de 2^{\tau - \ell}$ for $\ell > \tau$,  where $\tau$ is a positive
integer.

\subsubsection{Stopping Rule from CDF $\&$ CCDF}

By virtue of the CDF $\&$ CCDF of $\wh{\bs{\lm}}_\ell$, we propose a
class of multistage sampling schemes as follows.

\beT \la{Pos_rev_CDF}  Suppose that, for $\ell = 1, 2, \cd$,
decision variable $\bs{D}_\ell$ assumes values $1$ if {\small
$F_{\wh{\bs{\lm}}_\ell} ( \wh{\bs{\lm}}_\ell,
\f{\wh{\bs{\lm}}_\ell}{1 - \vep}  ) \leq \ze \de_\ell, \;
G_{\wh{\bs{\lm}}_\ell}  ( \wh{\bs{\lm}}_\ell,
\f{\wh{\bs{\lm}}_\ell}{1 + \vep} ) \leq \ze \de_\ell$}; and assumes
$0$ otherwise. The following statements hold true.

(I): $\Pr \{ \mbf{n} < \iy \} = 1$  provided that $\inf_{\ell > 0}
\f{n_{\ell + 1}}{n_\ell} > 1$.

(II): $\bb{E} [ \mbf{n} ] < \iy$ provided that $1 < \inf_{\ell > 0}
\f{n_{\ell + 1}}{n_\ell} \leq \sup_{\ell > 0} \f{n_{\ell +
1}}{n_\ell} < \iy$.

(III): {\small $\Pr \li \{ \li | \f{ \wh{\bs{\lm}} - \lm } { \lm }
\ri | < \vep \mid \lm \ri \} \geq 1 - \de$} for any $\lm  > 0$
provided that $ \ze \leq \f{1} { 2 (\tau + 1 ) }$.

(IV):  Let $0 < \eta < \ze \de$ and {\small ${\ell^\star} = \tau + 1
+ \li \lc \f{ \ln ( \ze \de \sh \eta )  } { \ln 2 } \ri \rc$}. Then,
$\Pr \{ | \wh{\bs{\lm}} - \lm  | \geq \vep \lm \mid \lm  \} < \de$
for any $\lm \in (0, \udl{\lm})$, where $\udl{\lm}$ is a number such
that {\small $0 < \udl{\lm} < z_\ell, \; \; \ell = 1, \cd,
{\ell^\star}$} and that $\sum_{\ell = 1}^{\ell^\star} \exp ( n_\ell
\mscr{M}_{\mrm{P}} ( z_\ell, \udl{\lm}) ) < \de - \eta$ with $z_\ell
= \min \{  z \in I_{\wh{\bs{\lm}}_\ell}: F_{ \wh{\bs{\lm}}_\ell }
(z, \f{z}{1 - \vep} ) > \ze \de_\ell \; \tx{or} \; G_{
\wh{\bs{\lm}}_\ell } ( z, \f{z}{1 + \vep}  ) > \ze \de_\ell \}$,
where $I_{\wh{\bs{\lm}}_\ell}$ represents the support of
$\wh{\bs{\lm}}_\ell$, for $\ell = 1, 2, \cd$.  Moreover, {\small
\bee &  & \Pr \li \{ b \leq \f{ \wh{\bs{\lm}} }{ 1 + \vep}, \;
\bs{l} \leq \ell^\star \mid a \ri \} \leq \Pr \li \{ \lm \leq \f{
\wh{\bs{\lm}} }{ 1 + \vep} \mid \lm \ri \} \leq \f{\eta}{2} + \Pr
\li \{ a \leq \f{ \wh{\bs{\lm}} }{ 1 + \vep}, \; \bs{l} \leq \ell^\star \mid b \ri \}, \\
&   & \Pr \li \{ a \geq \f{ \wh{\bs{\lm}} }{ 1 - \vep}, \; \bs{l}
\leq \ell^\star \mid b \ri \}  \leq \Pr \li \{ \lm \geq \f{
\wh{\bs{\lm}} }{ 1 - \vep} \mid \lm   \ri \} \leq \f{\eta}{2} + \Pr
\li \{ b \geq \f{ \wh{\bs{\lm}} }{ 1 - \vep}, \; \bs{l} \leq
\ell^\star \mid a \ri \} \eee} for any $\lm \in [a, b]$, where $a$
and $b$ are numbers such that $0 < b < ( 1 + \vep ) a$.

(V): $\Pr \{ | \wh{\bs{\lm}} - \lm  | \geq \vep \lm \mid \lm  \} <
\de$ for any $\lm \in (\ovl{\lm}, \iy)$, where $\ovl{\lm}$ is a
number such that {\small $\ovl{\lm} > z_1$} and that $2 \exp ( n_1
\mscr{M}_{\mrm{P}} ( (1 + \vep) \ovl{\lm}, \ovl{\lm}) ) + \exp ( n_1
\mscr{M}_{\mrm{P}} ( z_1, \ovl{\lm}) ) < \de$.

(VI): Let the sample sizes of the multistage sampling scheme be a
sequence $n_\ell = \li \lc m \ga^{\ell - 1} \ri \rc, \; \ell = 1, 2,
\cd$, where $\ga \geq 1 + \f{1}{m} > 1$. Let $\ep > 0, \; 0 < \eta <
1$ and $c = - \mscr{M}_{\mrm{P}} ( \eta \lm, \lm )$. Let $\ka$ be an
integer such that {\small $\ka
> \max \li \{ \tau, \;  \f{1}{\ln \ga} \ln \li ( \f{1}{ c m } \ri
) + 1, \; \f{1}{\ln \ga} \ln \li ( \f{1}{ c m } \ln \f{\ga}{c \ep}
\ri ) + 1, \; \tau + \f{1}{\ga - 1} + \f{ \ln (\ze \de) } { \ln 2 }
\ri \} $} and {\small $\mscr{M}_{\mrm{P}} ( \eta \lm, \f{\eta \lm}{1
+ \vep} ) < \f{ \ln (\ze \de_\ka )}{n_\ka}$}. Then, $\bb{E}
[\mbf{n}] < \ep + n_1 + \sum_{\ell = 1}^{\ka} (n_{\ell + 1} -
n_\ell) \Pr \{ \bs{l}
> \ell \}$.

\eeT

\subsubsection{Stopping Rule from Chernoff Bounds}

By virtue of Chernoff bounds of the CDF $\&$ CCDF of
$\wh{\bs{\lm}}_\ell$, we propose a class of multistage sampling
schemes as follows.

\beT \la{Pos_rev_Chernoff}  Suppose that, for $\ell = 1, 2, \cd$,
decision variable $\bs{D}_\ell$ assumes values $1$ if {\small
$\wh{\bs{\lm}}_\ell \geq \f{ \ln (\ze \de_\ell) } { n_\ell } \f{1 +
\vep} { \vep - (1 + \vep) \ln (1 + \vep)}$}; and assumes $0$
otherwise. The following statements hold true.

(I): $\Pr \{ \mbf{n} < \iy \} = 1$  provided that $\inf_{\ell > 0}
\f{n_{\ell + 1}}{n_\ell} > 1$.

(II): $\bb{E} [ \mbf{n} ] < \iy$ provided that $1 < \inf_{\ell > 0}
\f{n_{\ell + 1}}{n_\ell} \leq \sup_{\ell > 0} \f{n_{\ell +
1}}{n_\ell} < \iy$.

(III): {\small $\Pr \{  | \wh{\bs{\lm}} - \lm  | < \vep \lm \mid \lm
\} \geq 1 - \de$} for any $\lm  > 0$ provided that $ \ze \leq \f{1}
{ 2 (\tau + 1 ) }$.

(IV):  Let $0 < \eta < \ze \de$ and {\small ${\ell^\star} = \tau + 1
+ \li \lc \f{ \ln ( \ze \de \sh \eta )  } { \ln 2 } \ri \rc$}. Then,
$\Pr \{ | \wh{\bs{\lm}} - \lm | \geq \vep \lm \mid \lm \} < \de$ for
any $\lm \in (0, \udl{\lm})$, where $\udl{\lm}$ is a number such
that {\small $0 < \udl{\lm} < z_\ell, \; \; \ell = \tau, \cd,
{\ell^\star}$} and that $\sum_{\ell = 1}^{\ell^\star} \exp ( n_\ell
\mscr{M}_{\mrm{P}} ( z_\ell, \udl{\lm}) ) < \de - \eta$ with {\small
$z_\ell = \f{ \ln (\ze \de_\ell) } { n_\ell } \f{1 + \vep} { \vep -
(1 + \vep) \ln (1 + \vep)}$} for $\ell = 1, 2, \cd$. Moreover,
{\small \bee &  & \Pr \li \{ b \leq \f{ \wh{\bs{\lm}} }{ 1 + \vep},
\; \bs{l} \leq \ell^\star \mid a \ri \} \leq \Pr \li \{ \lm \leq \f{
\wh{\bs{\lm}} }{ 1 + \vep} \mid \lm \ri \} \leq \f{\eta}{2} + \Pr
\li \{ a \leq \f{ \wh{\bs{\lm}} }{ 1 + \vep}, \; \bs{l} \leq \ell^\star \mid b \ri \}, \\
&   & \Pr \li \{ a \geq \f{ \wh{\bs{\lm}} }{ 1 - \vep}, \; \bs{l}
\leq \ell^\star \mid b \ri \}  \leq \Pr \li \{ \lm \geq \f{
\wh{\bs{\lm}} }{ 1 - \vep} \mid \lm   \ri \} \leq \f{\eta}{2} + \Pr
\li \{ b \geq \f{ \wh{\bs{\lm}} }{ 1 - \vep}, \; \bs{l} \leq
\ell^\star \mid a \ri \} \eee} for any $\lm \in [a, b]$, where $a$
and $b$ are numbers such that $0 < b < ( 1 + \vep ) a$.

(V): $\Pr \{  | \wh{\bs{\lm}} - \lm  | \geq \vep \lm \mid \lm  \} <
\de$ for any $\lm \in (\ovl{\lm}, \iy)$, where $\ovl{\lm}$ is a
number such that {\small $\ovl{\lm} > z_1$} and that $2 \exp ( n_1
\mscr{M}_{\mrm{P}} ( (1 + \vep) \ovl{\lm}, \ovl{\lm}) ) + \exp ( n_1
\mscr{M}_{\mrm{P}} ( z_1,  \ovl{\lm}) ) < \de$.

(VI): Let the sample sizes of the multistage sampling scheme be a
sequence $n_\ell = \li \lc m \ga^{\ell - 1} \ri \rc, \; \ell = 1, 2,
\cd$, where $\ga \geq 1 + \f{1}{m} > 1$. Let $\ep > 0, \; 0 < \eta <
1$ and $c = - \mscr{M}_{\mrm{P}} ( \eta \lm, \lm )$. Let $\ka$ be an
integer such that {\small $\ka
> \max \li \{ \tau, \;  \f{1}{\ln \ga} \ln \li ( \f{1}{ c m } \ri
) + 1, \; \f{1}{\ln \ga} \ln \li ( \f{1}{ c m } \ln \f{\ga}{c \ep}
\ri ) + 1, \; \tau + \f{1}{\ga - 1} + \f{ \ln (\ze \de) } { \ln 2 }
\ri \} $} and {\small $\mscr{M}_{\mrm{P}}  ( \eta \lm, \f{\eta
\lm}{1 + \vep}  ) < \f{ \ln (\ze \de_\ka )}{n_\ka}$}. Then, $\bb{E}
[\mbf{n}] < \ep + n_1 + \sum_{\ell = 1}^{\ka} (n_{\ell + 1} -
n_\ell) \Pr \{ \bs{l}
> \ell \}$.

\eeT

See Appendix \ref{App_Pos_rev_Chernoff} for a proof.

\subsubsection{Asymptotic Analysis of Multistage Sampling Schemes}

In this subsection, we shall focus on the asymptotic analysis of the
multistage sampling schemes which follow stopping rules derived from
Chernoff bounds of CDF $\&$ CCDF of $\wh{\bs{\lm}}_\ell$ as
described in Theorem \ref{Pos_rev_Chernoff}.  We assume that the
sample sizes $n_1, n_2, \cd$ are chosen as the ascending arrangement
of all distinct elements of the set {\small \be \la{defposrev} \li
\{ \li \lc  \f{C_{\tau - \ell} \; \ln (\ze \de) }{
\mscr{M}_{\mrm{P}} \li ( \lm^\prime, \f{\lm^\prime }{1 + \vep} \ri )
} \ri \rc : \ell = 1, 2, \cd \ri \} \ee} with $0 < \lm^\prime <
\lm^{\prime \prime}$, where $\tau$ is the maximum integer such that
$\f{C_{\tau - 1} \; \ln (\ze \de) }{ \mscr{M}_{\mrm{P}} \li (
\lm^\prime, \f{\lm^\prime }{1 + \vep} \ri ) }  \geq \f{ \ln (\ze
\de) }{ \mscr{M}_{\mrm{P}} \li ( \lm^{\prime \prime}, \f{\lm^{\prime
\prime} }{1 + \vep} \ri ) }$, i.e., $C_{\tau - 1} \geq  \f{
\mscr{M}_{\mrm{P}} \li ( \lm^\prime, \f{\lm^\prime }{1 + \vep} \ri )
} { \mscr{M}_{\mrm{P}} \li ( \lm^{\prime \prime}, \f{\lm^{\prime
\prime} }{1 + \vep} \ri ) }$. With regard to the asymptotic
performance of the sampling scheme, we have

\beT \la{Pos_Asp_Analysis_Rev} Let {\small $\mcal{N}_{\mrm{r}} (\lm,
\vep) = \f{ \ln ( \ze \de ) } { \mscr{M}_{\mrm{P}} ( \lm , \f{\lm}{1
+ \vep}) } $}. Let $\mcal{N}_{\mrm{f}} (\lm, \vep)$ be the minimum
sample number $n$ such that {\small $\Pr \{ | \f{\sum_{i = 1}^n
X_i}{n} - \lm  | < \vep \lm \mid \lm \} > 1 - \ze \de$} for a
fixed-size sampling procedure.  Let $j_\lm$ be the largest integer
$j$ such that $C_{j} \geq \f{\lm^\prime}{\lm}$.  Let  {\small $d =
\sq{2 \ln \f{1}{\ze \de} }$} and  {\small $\ka_\lm = \f{\lm}{
\lm^\prime } C_{j_\lm}$}.  Let {\small $\ro_\lm = \f{\lm}{
\lm^\prime } C_{j_\lm - 1} - 1$} if $\ka_\lm = 1$ and $\ro_\lm =
\ka_\lm - 1$ otherwise. For $\lm \in (\lm^\prime, \lm^{\prime
\prime} )$, the following statements hold true:

(I): {\small $\Pr \li \{  1 \leq \limsup_{\vep \to 0} \f{ \mbf{n} }
{ \mcal{N}_{\mrm{r}} (\lm, \vep) } \leq 1 + \ro_\lm  \ri \} = 1$}.
Specially, {\small $\Pr \li \{  \lim_{\vep \to 0} \f{ \mbf{n} } {
\mcal{N}_{\mrm{r}} (\lm, \vep) } = \ka_\lm  \ri \} = 1$} if $\ka_\lm
> 1$.

(II): {\small $\lim_{\vep \to 0} \f{ \bb{E} [ \mbf{n} ] } {
\mcal{N}_{\mrm{f}} (\lm, \vep)} = \li ( \f{ d } { \mcal{Z}_{\ze \de}
} \ri )^2 \times \lim_{\vep \to 0} \f{ \bb{E} [ \mbf{n} ] } {
\mcal{N}_{\mrm{r}} (\lm, \vep) }$}, where
\[ \lim_{\vep \to 0} \f{ \bb{E} [ \mbf{n} ] }  {
\mcal{N}_{\mrm{r}} (\lm, \vep) } = \bec \ka_\lm
& \tx{if} \; \ka_\lm > 1,\\
1  + \f{\ro_\lm}{2} &  \tx{otherwise} \eec
\]
and $1 \leq \lim_{\vep \to 0} \f{ \bb{E} [ \mbf{n} ] }  {
\mcal{N}_{\mrm{r}} (\lm, \vep) } \leq 1 + \ro_\lm$.

(III): $\lim_{\vep \to 0} \Pr \{ | \wh{\bs{\lm}} - \lm | < \vep \lm
\} = 2 \Phi \li ( d \sq{\ka_\lm} \ri ) - 1 \geq 2 \Phi \li ( d \ri )
- 1 > 1 - 2 \ze \de$.

 \eeT

\bsk

See Appendix \ref{App_Pos_Asp_Analysis_Rev} for a proof.

\subsection{Control of Absolute and Relative Errors}

In this section, we shall focus on the design of multistage sampling
schemes for estimating Poisson parameter $\lm$ with a mixed error
criterion. Specifically, for $\vep_a > 0$ and $0 < \vep_r < 1$, we
wish to construct a multistage sampling scheme and its associated
estimator $\wh{\bs{\lm}}$ for $\lm$ such that $\Pr \{ |
\wh{\bs{\lm}} - \lm | < \vep_a, \;  | \wh{\bs{\lm}} - \lm | < \vep_r
\lm \mid \lm \} > 1 - \de$ for any $\lm \in (0, \iy)$. This is
equivalent to the construction of a random interval with lower limit
$\mscr{L} ( \wh{\bs{\lm}} )$ and upper limit $\mscr{U} (
\wh{\bs{\lm}} )$ such that $\Pr \{  \mscr{L} ( \wh{\bs{\lm}} ) < \lm
< \mscr{U} ( \wh{\bs{\lm}} ) \mid \lm \} > 1 - \de$ for any $\lm \in
(0, \iy)$, where $\mscr{L} ( . )$ and $\mscr{U} ( . )$ are functions
such that $\mscr{L} ( z ) = \min \{ z - \vep_a, \; \f{z}{1 + \vep_r}
\}$ and $\mscr{U} ( z ) = \max \{ z + \vep_a, \; \f{z}{1 - \vep_r}
\}$ for $z \in [0, \iy)$.  In the sequel, we shall propose
multistage sampling schemes such that the number of stages, $s$, is
finite and that the sample sizes are deterministic numbers $n_1 <
n_2 < \cd < n_s$.

\subsubsection {Stopping Rules from CDF $\&$ CCDF and Chernoff Bounds} \la{multimixP}

To estimate $\lm$ with a mixed precision criterion, we propose two
types of multistage sampling schemes with different stopping rules
as follows.

\bed

\item [Stopping Rule (i):] For $\ell = 1, \cd, s$, decision variable
$\bs{D}_\ell$  assumes value $1$ if {\small $F_{\wh{\bs{\lm}}_\ell}
(\wh{\bs{\lm}}_\ell, \mscr{U} ( \wh{\bs{\lm}}_\ell ) ) \leq \ze \de,
\; G_{\wh{\bs{\lm}}_\ell} (\wh{\bs{\lm}}_\ell, \mscr{L} (
\wh{\bs{\lm}}_\ell ) ) \leq \ze \de$}; and assumes value $0$
otherwise.

\item [Stopping Rule (ii):] For $\ell = 1, \cd, s$, decision variable
$\bs{D}_\ell$ assumes value $1$  if {\small
\[
\max \{ \mscr{M}_{\mrm{P}} (\wh{\bs{\lm}}_\ell, \mscr{L} (
\wh{\bs{\lm}}_\ell ) ), \; \mscr{M}_{\mrm{P}} (\wh{\bs{\lm}}_\ell,
\mscr{U} ( \wh{\bs{\lm}}_\ell ) ) \} \leq \f{ \ln ( \ze \de ) } {
n_\ell }; \]} and assumes value $0$ otherwise.

\eed

Stopping rule (i) is derived by virtue of the CDF $\&$ CCDF of
$\wh{\bs{\lm}}_\ell$.  Stopping rule (ii) is derived by virtue of
Chernoff bounds of the CDF $\&$ CCDF of $\wh{\bs{\lm}}_\ell$. For
both types of multistage sampling schemes described above, we have
the following results.

\beT \la{Pos_mix_CDF}   Suppose that the sample size for the $s$-th
stage is no less than {\small $\li \lc  \f{ \ln (\ze \de) } {
\mscr{M}_{\mrm{P}} (\f{\vep_a}{\vep_r} + \vep_a, \f{\vep_a}{\vep_r}
) } \ri \rc$}.  Then, \bee &  & \Pr \{ \lm \leq \mscr{L} (
\wh{\bs{\lm}} ) \mid \lm \} \leq \sum_{\ell = 1}^s \Pr \{ \lm \leq
\mscr{L} ( \wh{\bs{\lm}}_\ell ), \;
\bs{D}_\ell = 1 \mid \lm \} \leq s \ze \de,\\
&    & \Pr \{  \lm \geq \mscr{U} ( \wh{\bs{\lm}} ) \mid \lm \} \leq
\sum_{\ell = 1}^s  \Pr \{  \lm \geq \mscr{U} ( \wh{\bs{\lm}}_\ell ),
\; \bs{D}_\ell = 1 \mid \lm \} \leq s \ze \de \eee for any $\lm >
0$. Moreover, {\small $\Pr \{ | \wh{\bs{\lm}} - \lm | < \vep_a  \;
\mrm{or} \;  | \f{\wh{\bs{\lm}} - \lm } {\lm} | < \vep_r \mid \lm \}
> 1 - \de$} for any $\lm > 0$ provided that $\Pr \{  \lm \leq
\mscr{L} ( \wh{\bs{\lm}} ) \mid \lm \} + \Pr \{  \lm \geq \mscr{U} (
\wh{\bs{\lm}} ) \mid \lm \} < \de$ for any $\lm \in (0, \ovl{\lm}
]$, where $\ovl{\lm} > 0$ is the unique number satisfying {\small
$\sum_{\ell = 1}^s \exp ( n_\ell \mscr{M}_{\mrm{P}} ( \ovl{\lm} (1 +
\vep_r), \ovl{\lm}) ) = \f{\de}{2}$. }

\eeT

\bsk

See Appendix \ref{App_Pos_mix_CDF} for a proof.  Based on the
criteria proposed in Section \ref{gen_structure}, the sample sizes
$n_1 < n_2 < \cd < n_s$ can be chosen as the ascending arrangement
of all distinct elements of {\small \be \la{mixpos} \li \{ \li \lc
\f{ C_{\tau - \ell} \; \ln (\ze \de) }{\mscr{M}_{\mrm{P}}
(\f{\vep_a}{\vep_r} + \vep_a, \f{\vep_a}{\vep_r} )} \ri \rc : \ell =
1, \cd, \tau \ri \},  \ee} where $\tau$ is the maximum integer such
that $\f{ C_{\tau - 1} \; \ln (\ze \de) }{\mscr{M}_{\mrm{P}}
(\f{\vep_a}{\vep_r} + \vep_a, \f{\vep_a}{\vep_r} )} \geq \f{ \ln
\f{1}{\ze \de} }{\vep_a}$, i.e., $C_{\tau - 1} \geq - \f{
\mscr{M}_{\mrm{P}} (\f{\vep_a}{\vep_r} + \vep_a, \f{\vep_a}{\vep_r}
)  } { \vep_a }$. For such a choice of sample sizes, as a result of
Theorem \ref{Pos_mix_CDF}, we have that {\small $\Pr \{ |
\wh{\bs{\lm}} - \lm | < \vep_a  \; \mrm{or} \;  | \f{\wh{\bs{\lm}} -
\lm } {\lm} | < \vep_r \mid \lm \}
> 1 - \de$} for any $\lm > 0$ provided that $\ze < \f{1}{2
\tau}$.

To evaluate the coverage probability associated with a multistage
sampling scheme following a stopping rule derived from Chernoff
bounds, we need to express $\{ \bs{D}_\ell = i \}$ in terms of
$K_\ell$.  For this purpose, the following result is useful.

\beT

\la{Range_Pos_mix}

 Let $\lm^\star = \f{ \vep_a} {\vep_r}$.  Then, {\small $\{ \bs{D}_\ell = 0 \}
 = \{ \mscr{M}_{\mrm{P}} (\wh{\bs{\lm}}_\ell, \mscr{L} ( \wh{\bs{\lm}}_\ell ) )  >
\f{ \ln ( \ze \de  ) } { n_\ell } \} \cup \{ \mscr{M}_{\mrm{P}}
(\wh{\bs{\lm}}_\ell, \mscr{U} ( \wh{\bs{\lm}}_\ell ) )
> \f{ \ln ( \ze \de  ) } { n_\ell } \}$} for $\ell = 1, \cd, s - 1$
and the following statements hold true:

(I) {\small $\{  \mscr{M}_{\mrm{P}} (\wh{\bs{\lm}}_\ell, \mscr{L} (
\wh{\bs{\lm}}_\ell ) ) > \f{ \ln ( \ze \de ) } { n_\ell } \} = \{
n_\ell \; z_a^- < K_\ell < n_\ell \; z_r^+ \}$} where $z_r^+$ is the
unique solution of equation {\small $\mscr{M}_{\mrm{P}} (z, \f{z}{1
+ \vep_r} ) = \f{ \ln ( \ze \de  ) } { n_\ell }$} with respect to $z
\in (\lm^\star + \vep_a, \iy )$, and $z_a^-$ is the unique solution
of equation {\small $\mscr{M}_{\mrm{P}} (z, z - \vep_a ) = \f{ \ln (
\ze \de ) } { n_\ell }$} with respect to $z \in (\vep_a, \lm^\star +
\vep_a)$.

(II) {\small \[ \li \{ \mscr{M}_{\mrm{P}} (\wh{\bs{\lm}}_\ell,
\mscr{U} ( \wh{\bs{\lm}}_\ell ) )
> \f{ \ln ( \ze \de  ) } { n_\ell } \ri \} = \bec  \{ 0
\leq K_\ell < n_\ell \; z_r^- \} &
\tx{for} \;  n_\ell < \f{ \ln \f{1}{\ze \de} } { \vep_a },\\
\{ n_\ell \; z_a^+ < K_\ell < n_\ell \; z_r^- \} & \tx{for} \; \f{
\ln \f{1}{\ze \de} } { \vep_a } \leq n_\ell < \f{ \ln (\ze \de) }
{ \mscr{M}_{\mrm{P}} (\lm^\star - \vep_a, \lm^\star) },\\
\emptyset & \tx{for} \; n_\ell \geq \f{ \ln (\ze \de) } {
\mscr{M}_{\mrm{P}} (\lm^\star - \vep_a, \lm^\star) }
 \eec
\]}
where $z_r^-$ is the unique solution of equation {\small
$\mscr{M}_{\mrm{P}} (z, \f{z}{1 - \vep_r} ) = \f{ \ln ( \ze \de  ) }
{ n_\ell }$} with respect to $z \in (\lm^\star - \vep_a, \iy )$, and
$z_a^+$ is the unique solution of equation {\small
$\mscr{M}_{\mrm{P}} (z, z + \vep_a ) = \f{ \ln ( \ze \de  ) } {
n_\ell }$} with respect to $z \in [0, \lm^\star - \vep_a)$.

\eeT

\bsk

Theorem \ref{Range_Pos_mix} can be shown by a variation of the
argument for Theorem \ref{Bino_range_mix}.

\subsubsection{Asymptotic Stopping Rule}

It should be noted that, for small $\vep_a$ and $\vep_r$, we can
simplify, by using Taylor's series expansion formula $\ln (1 + x) =
x - \f{x^2}{2} + o (x^2)$, the sampling schemes as described in
Section \ref{multimixP} as follows:

(i) The sequence of sample sizes $n_1, \cd, n_s$ is defined as the
ascending arrangement of all distinct elements of {\small $ \li \{
\li \lc C_{\tau - \ell} \li ( \f{2}{\vep_r} \ri ) \ln \f{1}{\ze \de}
\ri \rc : \ell = 1, \cd, \tau \ri \}$}, where $\tau$ is the maximum
integer such that $C_{\tau - 1} \geq \f{ \vep_r}{2}$.

(ii) The decision variables are defined such that $\bs{D}_\ell = 1$
if $n_\ell \geq \f{ \wh{\bs{\lm}}_\ell  \; 2 \ln \f{1}{\ze \de} }{
\max \{ \vep_a^2, \; (\vep_r \wh{\bs{\lm}}_\ell )^2 \} }$; and
$\bs{D}_\ell = 0$ otherwise.

\bsk

For such a simplified sampling scheme, we have \bel \sum_{\ell =
1}^s \Pr \li \{ | \wh{\bs{\lm}}_\ell - \lm | \geq \max \{\vep_a,
\vep_r \lm \}, \; \bs{D}_\ell = 1 \ri \} & \leq & \sum_{\ell = 1}^s
\Pr \li \{ |
\wh{\bs{\lm}}_\ell - \lm | \geq \max \{\vep_a, \vep_r \lm \} \ri \} \nonumber\\
& \leq & \sum_{\ell = 1}^\tau \Pr \li \{ | \wh{\bs{\lm}}_\ell - \lm
| \geq \max \{\vep_a, \vep_r \lm \} \ri \} \nonumber\\
& \leq &  \sum_{\ell = 1}^\tau 2 \exp \li ( n_\ell
\mscr{M}_{\mrm{P}}
\li ( \f{\vep_a}{\vep_r} + \vep_a, \f{\vep_a}{\vep_r} \ri )  \ri ) \la{PineCH}\\
& < & 2 \tau \exp \li ( n_1 \mscr{M}_{\mrm{P}} \li (
\f{\vep_a}{\vep_r} + \vep_a, \f{\vep_a}{\vep_r} \ri )  \ri ),
\la{PLastB} \eel where (\ref{PineCH}) is due to Theorem 1 of
\cite{Chen4}. As can be seen from (\ref{PLastB}), the last bound is
independent of $\lm$ and can be made smaller than $\de$ if $\ze$ is
sufficiently small. This establishes the claim and it follows that
{\small $\Pr \li \{ \li | \wh{\bs{\lm}} - \lm \ri | < \vep_a \;
\mrm{or} \; \li | \f{\wh{\bs{\lm}} - \lm } {\lm }
 \ri | < \vep_r \mid \lm \ri \} > 1 - \de$} for any $\lm \in (0, \iy)$ if $\ze$ is sufficiently
small.

\subsubsection{Asymptotic Analysis of Multistage Sampling Schemes}

In this subsection, we shall focus on the asymptotic analysis of
multistage inverse sampling schemes.  Throughout this subsection, we
assume that the multistage sampling schemes follow stopping rules
derived from Chernoff bounds as described in Section
\ref{multimixP}. Moreover, we assume that the sample sizes $n_1,
\cd, n_s$ are chosen as the ascending arrangement of all distinct
elements of  the set defined by (\ref{mixpos}).

With regard to the tightness of the double-decision-variable method,
we have

\beT

\la{Pos_mix_DDV_Asp}

Let $\mscr{R}$ be a subset of real numbers.  Define {\small \[
\ovl{P}  = \sum_{\ell = 1}^s \Pr \{ \wh{\bs{\lm}}_\ell \in \mscr{R},
\; \bs{D}_{\ell - 1} = 0, \; \bs{D}_\ell = 1 \}, \qqu \udl{P} = 1 -
\sum_{\ell = 1}^s \Pr \{ \wh{\bs{\lm}}_\ell \notin \mscr{R}, \;
\bs{D}_{\ell - 1} = 0, \; \bs{D}_\ell = 1 \}.
\]}
Then, $\udl{P} \leq \Pr \{ \wh{\bs{\lm}} \in \mscr{R} \} \leq
\ovl{P}$ and {\small $\lim_{\vep_a \to 0} | \Pr \{ \wh{\bs{\lm}} \in
\mscr{R} \} - \ovl{P} | = \lim_{\vep_a \to 0} | \Pr \{ \wh{\bs{\lm}}
\in \mscr{R} \} - \udl{P} | = 0$} for any $\lm \in (0, \iy)$,  where
the limits are taken under the constraint that $\f{\vep_a}{\vep_r}$
is fixed.

\eeT

\bsk

See Appendix \ref{App_Pos_mix_DDV_Asp} for a proof.

With regard to the asymptotic performance of the sampling scheme as
$\vep_a$ and $\vep_r$ tend to $0$, we have

\beT \la{Pos_mix_Asp_Analysis}  Let $\mcal{N}_{\mrm{f}} (\lm,
\vep_a, \vep_r)$ be the minimum sample number $n$ such that {\small
\[
\Pr \li \{ \li | \f{\sum_{i = 1}^n X_i}{n} - \lm
\ri | < \vep_a \; \tx{or} \; \li | \f{\sum_{i = 1}^n X_i}{n} - \lm
\ri | < \vep_r \lm \mid \lm \ri \} > 1 - \ze \de
\]}
for a fixed-size sampling procedure.  Let {\small
$\mcal{N}_{\mrm{m}} (\lm, \vep_a, \vep_r) = \f{ \ln ( \ze \de  ) } {
\max \{ \mscr{M}_{\mrm{P}} (\lm, \udl{\lm} ), \; \mscr{M}_{\mrm{P}}
(\lm, \ovl{\lm} ) \}  } $}, where $\udl{\lm} = \min \{ \lm - \vep_a,
\; \f{\lm}{1 + \vep_r} \}$ and $\ovl{\lm}= \max \{ \lm + \vep_a, \;
\f{\lm}{1 - \vep_r} \}$.  Define $\lm^\star = \f{\vep_a}{\vep_r}, \;
d = \sq{ 2 \ln \f{1}{\ze \de} }$,
\[
r(\lm) = \bec \f{\lm} { \lm^\star }  & for
\; \lm \in (0, \lm^\star],\\

\f{\lm^\star}{ \lm  }  & for \; \lm \in (\lm^\star, \iy) \eec \qqu \nu = \bec \f{1}{3} \li ( 2 -
\f{\lm}{\lm^\star} \ri )  & \tx{for} \; \lm \in (0, \lm^\star],\\
1  & \tx{for} \; \lm \in (\lm^\star, \iy). \eec
\]
Let {\small $\ka_\lm = \f{C_{j_\lm}}{ r (\lm) } $}, where $j_\lm$ is
the maximum integer $j$ such that $C_j \geq r (\lm)$. Let {\small
$\ro_\lm = \f{C_{j_\lm - 1}} { r (\lm) } - 1$} if $\ka_\lm = 1, \;
j_\lm > 0$ and $\ro_\lm = \ka_\lm - 1$ otherwise.  The following
statements hold true under the condition that $\f{\vep_a}{\vep_r}$
is fixed.

(I): {\small $\Pr \li \{  1 \leq \limsup_{\vep_a \to 0} \f{ \mbf{n}
} { \mcal{N}_{\mrm{m}} (\lm, \vep_a, \vep_r) } \leq 1 + \ro_\lm \ri
\} = 1$}.  Specially, {\small $\Pr \li \{ \lim_{\vep_a \to 0} \f{
\mbf{n} } { \mcal{N}_{\mrm{m}} (\lm, \vep_a, \vep_r) } = \ka_\lm \ri
\} = 1$} if $\ka_\lm > 1$.

(II): $\lim_{\vep_a \to 0} \f{ \bb{E} [ \mbf{n} ] } {
\mcal{N}_{\mrm{f}}  (\lm, \vep_a, \vep_r)} = \li ( \f{ d } {
\mcal{Z}_{\ze \de} } \ri )^2 \times \lim_{\vep_a \to 0} \f{ \bb{E} [
\mbf{n} ] } { \mcal{N}_{\mrm{m}} (\lm, \vep_a, \vep_r) }$, where
\[
\lim_{\vep_a \to 0} \f{ \bb{E} [ \mbf{n} ] } { \mcal{N}_{\mrm{m}}
(\lm, \vep_a, \vep_r) } = \bec
\ka_\lm & \tx{if} \; \ka_\lm > 1,\\
1 + \ro_\lm \Phi (\nu d ) &  \tx{otherwise} \eec
\]
and $1 \leq \lim_{\vep_a \to 0} \f{ \bb{E} [ \mbf{n} ] } {
\mcal{N}_{\mrm{m}} (\lm, \vep_a, \vep_r) } \leq 1 + \ro_\lm$.

(III): If $\ka_\lm > 1$, then $\lim_{\vep_a \to 0} \Pr \{ | \wh{\bs{\lm}} - \lm | < \vep_a \; \tx{or} \; | \wh{\bs{\lm}} - \lm | < \vep_r \lm \}
= 2 \Phi \li ( d \sq{\ka_\lm} \ri ) - 1 > 2 \Phi \li ( d \ri ) - 1 > 1 - 2 \ze \de$. Otherwise,  $\Phi \li ( d  \ri )  + \Phi \li ( d \sq{1 +
\ro_\lm}  \ri )  - 1  > \lim_{\vep_a \to 0} \Pr \{ | \wh{\bs{\lm}} - \lm | < \vep_a \; \tx{or} \; | \wh{\bs{\lm}} - \lm | < \vep_r \lm   \} = 1
+ \Phi(d) - \Phi(\nu d) - \Psi (\ro_\lm, \nu, d)
> \Phi \li ( d \ri ) + 2 \Phi \li ( d \sq{1 + \ro_\lm}  \ri )  - 2  > 1 - 3 \ze \de$.

\eeT

\bsk

See Appendix \ref{App_Pos_mix_Asp_Analysis} for a proof.

\section{Estimation of Finite Population Proportion}

In this section, we consider the problem of estimating the
proportion of a finite population, which has been discussed in
Section \ref{finite_proportion}.  We shall focus on multistage
sampling schemes with deterministic sample sizes $n_1 < n_2 < \cd <
n_s$.  Our methods are described in the sequel.

Define {\small $K_\ell = \sum_{i=1}^{n_\ell} X_i, \;
\wh{\bs{p}}_\ell = \f{K_\ell}{n_\ell}$} for $\ell = 1, \cd, s$.
Suppose the stopping rule is that sampling without replacement is
continued until $\bs{D}_\ell = 1$ for some $\ell \in \{1, \cd, s\}$.
Define {\small $\wh{\bs{p}} = \wh{\bs{p}}_{\bs{l}}$},  where
$\bs{l}$ is the index of stage at which the sampling is terminated.

By using various functions to define random intervals, we can unify
the estimation problems associated with absolute, relative  and
mixed precision.  Specifically, for estimating $p$ with margin of
absolute error $\vep \in (0,1)$, we have $\Pr \{ | \wh{\bs{p}} - p |
\leq \vep \} = \Pr \{  \mscr{L} ( \wh{\bs{p}} ) < p < \mscr{U} (
\wh{\bs{p}} ) \}$,  where $\mscr{L} (.)$ and $\mscr{U} (.)$ are
functions such that $\mscr{L} (z) = \f{1}{N} \li \lc N ( z - \vep)
\ri \rc - \f{1}{N}$ and $\mscr{U} (z) = \f{1}{N} \li \lf N ( z +
\vep) \ri \rf + \f{1}{N}$ for $z \in [0, 1]$.  For estimating $p$
with margin of relative error $\vep \in (0,1)$, we have $\Pr \{ |
\wh{\bs{p}} - p | \leq \vep p \} = \Pr \{  \mscr{L} ( \wh{\bs{p}} )
< p < \mscr{U} ( \wh{\bs{p}} ) \}$,  where $\mscr{L} (.)$ and
$\mscr{U} (.)$ are functions such that $\mscr{L} (z) = \f{1}{N} \li
\lc N z \sh (1 + \vep) \ri \rc - \f{1}{N}$ and $\mscr{U} (z) =
\f{1}{N} \li \lf N z \sh ( 1 - \vep) \ri \rf + \f{1}{N}$ for $z \in
[0, 1]$. For estimating $p$ with margin of absolute error $\vep_a
\in (0,1)$ and margin of relative error $\vep_r \in (0,1)$, we have
$\Pr \{ | \wh{\bs{p}} - p |  \leq \vep_a \; \tx{or} \; | \wh{\bs{p}}
- p | \leq \vep_r p \} = \Pr \{ \mscr{L} ( \wh{\bs{p}} ) < p <
\mscr{U} ( \wh{\bs{p}} ) \}$, where $\mscr{L} (.)$ and $\mscr{U}
(.)$ are functions such that
\[
\mscr{L} ( z ) = \f{1}{N} \li \lc N \min \li ( z - \vep_a, \; \f{ z
} { 1 + \vep_r  }  \ri ) \ri \rc - \f{1}{N}, \qqu \mscr{U} ( z ) =
\f{1}{N} \li \lf N \max \li ( z + \vep_a, \; \f{ z } { 1 - \vep_r }
\ri ) \ri \rf + \f{1}{N}
\]
for $z \in [0, 1]$.  Therefore, multistage estimation problems
associated with absolute, relative  and mixed precision can be cast
as the general problem of constructing a random interval with lower
limit $\mscr{L} ( \wh{\bs{p}} )$ and upper limit $\mscr{L} (
\wh{\bs{p}} )$ such that $\Pr \{ \mscr{L} ( \wh{\bs{p}} ) < p <
\mscr{U} ( \wh{\bs{p}} ) \} \geq 1 - \de$.  For this purpose, making
use of Theorems \ref{Monotone_second} and \ref{FiniteULE}, we
immediately obtain the following result.

\beC \la{Fthem_abs_ave}

Suppose the sample size of the $s$-th stage is no less than the
minimum number $n$ such that $1 - S_N ( k - 1, n, \mscr{L}
(\f{k}{n}) ) \leq \ze \de$ and $S_N (k , n, \mscr{U} (\f{k}{n}) )
\leq \ze \de$ for $0 \leq k \leq n$. For $\ell = 1, \cd, s$, define
$\bs{D}_\ell$ such that $\bs{D}_\ell$ assumes value $1$ if $1 - S_N
( K_\ell - 1, n_\ell, \mscr{L} (\wh{\bs{p}}_\ell) ) \leq \ze \de, \;
S_N (K_\ell , n_\ell, \mscr{U} (\wh{\bs{p}}_\ell) ) \leq \ze \de$;
and assumes value $0$ otherwise.  Then, \bee &  & \Pr \{ p \leq
\mscr{L} ( \wh{\bs{p}} ) \mid p \} \leq \sum_{\ell = 1}^s \Pr \{ p
\leq \mscr{L} ( \wh{\bs{p}}_\ell ), \;
\bs{D}_\ell = 1 \mid p \} \leq s \ze \de,\\
&    & \Pr \{  p \geq \mscr{U} ( \wh{\bs{p}} ) \mid p \} \leq
\sum_{\ell = 1}^s  \Pr \{  p \geq \mscr{U} ( \wh{\bs{p}}_\ell ), \;
\bs{D}_\ell = 1 \mid p \} \leq s \ze \de \eee and $\Pr \{ \mscr{L} (
\wh{\bs{p}} ) < p < \mscr{U} ( \wh{\bs{p}} ) \mid p \}  \geq 1 - 2 s
\ze \de$ for any $p \in \Se$.

\eeC

\bsk

Let {\small \bee & & n_{\mrm{min}} = 1 +
 \max \li \{n:  1 - S_N \li ( k - 1, n, \mscr{L} \li (\f{k}{n} \ri )
 \ri ) > \ze \de \; \tx{or} \; S_N \li (k , n, \mscr{U} \li (\f{k}{n} \ri ) \ri ) > \ze
\de  \;  \tx{for $0 \leq k \leq n$} \ri \},\\
&  &  n_{\mrm{max}} = \min \li \{n : 1 - S_N \li ( k - 1, n,
\mscr{L} \li (\f{k}{n} \ri ) \ri ) \leq \ze \de \; \tx{and} \; S_N
\li (k , n, \mscr{U} \li ( \f{k}{n} \ri ) \ri ) \leq \ze \de \;
\tx{for $0 \leq k \leq n$} \ri \}. \eee} Based on the criteria
proposed in Section \ref{gen_structure}, the sample sizes $n_1 < n_2
< \cd < n_s$ can be chosen as the ascending arrangement of all
distinct elements of the set {\small $\li \{ \li \lc C_{\tau - \ell}
\; n_{\mrm{max}} \ri \rc : 1 \leq \ell \leq \tau \ri \}$}, where
$\tau$ is the maximum integer such that $C_{\tau - 1} \geq \f{
n_{\mrm{min}} } { n_{\mrm{max}} }$.

\bsk

Now, define \be \la{defc} \mcal{C} (z, p, n, N) = \bec \f{
\bi{Np}{n} } { \bi{N}{n} } &
\tx{for} \; z = 1,\\
\f{ \bi{Np}{n z} \bi{N - Np}{n - n z} } { \bi{ \lf (N+1) z \rf }{n
z} \bi{N - \lf (N+1) z \rf }{n - n z} } & \tx{for} \; z \in \{
\f{k}{n} : k \in \bb{Z}, \; 0 \leq k < n \}
 \eec
\ee where $p \in \Se$. In order to develop multistage sampling
schemes with simple stopping boundaries, we have the following
results.

\beC \la{Fthem_abs_ave_LR}

Suppose the sample size of the $s$-th stage is no less than the
minimum number $n$ such that $\mcal{C} ( \f{k}{n}, \mscr{L}
(\f{k}{n}), n, N ) \leq \ze \de$ and $\mcal{C} ( \f{k}{n}, \mscr{U}
(\f{k}{n}), n, N ) \leq \ze \de$ for $0 \leq k \leq n$. For $\ell =
1, \cd, s$, define $\bs{D}_\ell$ such that $\bs{D}_\ell$ assumes
value $1$ if $\mcal{C}( \wh{\bs{p}}_\ell, \mscr{L}
(\wh{\bs{p}}_\ell), n_\ell, N ) \leq \ze \de, \; \mcal{C}(
\wh{\bs{p}}_\ell, \mscr{U} (\wh{\bs{p}}_\ell), n_\ell, N ) \leq \ze
\de$; and assumes value $0$ otherwise. Then, \bee &  & \Pr \{ p \leq
\mscr{L} ( \wh{\bs{p}} ) \mid p \} \leq \sum_{\ell = 1}^s \Pr \{ p
\leq \mscr{L} ( \wh{\bs{p}}_\ell ), \;
\bs{D}_\ell = 1 \mid p \} \leq s \ze \de,\\
&    & \Pr \{  p \geq \mscr{U} ( \wh{\bs{p}} ) \mid p \} \leq
\sum_{\ell = 1}^s  \Pr \{  p \geq \mscr{U} ( \wh{\bs{p}}_\ell ), \;
\bs{D}_\ell = 1 \mid p \} \leq s \ze \de \eee and $\Pr \{ \mscr{L} (
\wh{\bs{p}} ) < p < \mscr{U} ( \wh{\bs{p}} ) \mid p \}  \geq 1 - 2 s
\ze \de$ for any $p \in \Se$.

\eeC

Corollary \ref{Fthem_abs_ave_LR} can be shown by using Theorems
\ref{Monotone_second}, \ref{FiniteULE} and the inequalities obtained
by Chen \cite{ChenR} as follows: \bel &  & \Pr \li \{ \f{ \sum_{i =
1}^n X_i }{n} \geq z \mid p \ri \} \leq \mcal{C} (z, p, n, N) \qu
\tx{for} \; z \in
\li \{ \f{k}{n}: k \in \bb{Z}, \; n p \leq k \leq n \ri \}, \la{hyp1}\\
&  & \Pr \li \{ \f{ \sum_{i = 1}^n X_i }{n}  \leq z \mid p \ri \}
\leq \mcal{C} (z, p, n, N) \qu \tx{for} \; z \in \li \{ \f{k}{n}: k
\in \bb{Z}, \; 0 \leq k \leq n p \ri \} \la{hyp2} \eel where $p \in
\Se$.  Since $\sum_{i = 1}^n X_i$ has a hypergeometric distribution,
the above inequalities (\ref{hyp1}) and (\ref{hyp2}) provide simple
bounds for the tail probabilities of hypergeometric distribution,
which are substaintially less conservative than Hoeffding's
inequalities \cite{Hoeffding}.

\bsk

It is well known that, for a sampling without replacement with size
$n$, to guarantee that the estimator {\small $\wh{p} = \f{ \sum_{i =
1}^n X_i }{n}$} of the proportion $p = \f{M}{N}$ satisfy $\Pr
\left\{ \left | \wh{p} - p \right | \leq \vep \right\} \geq 1 -
\de$, it suffices to have {\small $n  \geq \f{N p ( 1 - p)}{p ( 1 -
p) + (N - 1) \vep^2 \sh \mcal{Z}_{\de \sh 2}^2 }$},  or
equivalently, {\small $\mcal{Z}_{\de \sh 2}^2 \li ( \f{ N } { n } -
1 \ri ) p ( 1 - p ) \leq (N - 1) \vep^2$} (see formula (1) in page
41 of \cite{Thom}). Therefore, for a very small margin of absolute
error $\vep$, we can develop simple multistage sampling schemes
based normal approximation  as follows.

To estimate the population proportion $p \in \Se$ with margin of
absolute error $\vep \in (0, 1)$, we can choose the sample sizes
$n_1 < n_2 < \cd < n_s$ as the ascending arrangement of all distinct
elements of the set {\small $\li \{ \li \lc \f{N C_{\tau - \ell} }{
1 + 4 (N - 1) \vep^2 \sh \mcal{Z}_{\ze \de}^2 } \ri \rc : \ell = 1,
\cd, \tau \ri \}$}, where $\tau$ is a positive integer.  With such a
choice of sample sizes, we define a stopping rule such that sampling
is continued until
\[
\mcal{Z}_{\ze \de}^2 \li ( \f{  N } { n_\ell } - 1 \ri )
\wh{\bs{p}}_\ell ( 1 - \wh{\bs{p}}_\ell ) \leq (N - 1) \vep^2
\]
is satisfied at some stage with index $\ell$.   Then, $\Pr \left\{
\left | \wh{\bs{p}} - p \right | \leq \vep \mid p \right\} \geq 1 -
\de$ for any $p \in \Se$ provided that the coverage tuning parameter
$\ze$ is sufficiently small.  In order to improve performance,
following the similar idea as  the stopping rule associated with
(\ref{simplegeneral}), we propose a more general stopping rule such
that   sampling is continued until \be \la{simplegeneral_finite} \li
( \li | \wh{\bs{p}}_\ell - \f{1}{2} \ri | - w \vep \ri )^2 \geq
\f{1}{4} + \f{ \vep^2 n_\ell } { 2 \ln (\ze \de) } \times \f{N-1}{N
- n_\ell} \ee is satisfied at some stage with index $\ell \in \{1,
\cd, s \}$. Here, $w \geq 0$ is a parameter affecting the shape of
the stopping boundary.  The factor $\f{N-1}{N - n_\ell}$ is
introduced in consideration of finite population size.  Under the
restriction that $0 < w \vep \leq \f{1}{4}$, the minimum sample size
$n_1$ should be chosen as the smallest integer such that $\li (
\f{1}{2} - w \vep \ri )^2 \geq \f{1}{4} + \f{ \vep^2 n_1 } { 2 \ln
(\ze \de) } \times \f{N-1}{N - n_1}$.  The maximum sample size $n_s$
should be chosen as the smallest integer such that $\f{1}{4} + \f{
\vep^2 n_s } { 2 \ln (\ze \de) } \times \f{N-1}{N - n_s} \leq 0$.
Clearly, for $0 < w \leq \f{1}{4 \vep}$, $\Pr \left\{ \left |
\wh{\bs{p}} - p \right | \leq \vep \mid p \right\} \geq 1 - \de$ for
any $p \in \Se$ provided that the coverage tuning parameter $\ze$ is
sufficiently small. By virtue of the bisection coverage tuning
technique and optimization over $w$, stopping rules of excellence
performance can be obtained.

To estimate the population proportion $p \in \Se$ with margin of
relative error $\vep \in (0, 1)$, we can choose the sample sizes
$n_1 < n_2 < \cd < n_s$ as the ascending arrangement of all distinct
elements of the set $\li \{ \li \lc N C_{\tau - \ell}  \ri \rc :
\ell = 1, \cd, \tau \ri \}$. The stopping rule is that sampling is
continued until
\[
\mcal{Z}_{\ze \de}^2 \li ( \f{  N } { n_\ell } - 1 \ri ) ( 1 -
\wh{\bs{p}}_\ell ) \leq (N - 1) \vep^2 \wh{\bs{p}}_\ell
\]
is satisfied at some stage with index $\ell$.  Then, $\Pr \left\{
\left | \wh{\bs{p}} - p \right | \leq \vep p \mid p \right\} \geq 1
- \de$ for any $p \in \Se$ provided that the coverage tuning
parameter $\ze$ is sufficiently small.  In order to improve the
performance of coverage probability, we propose to revise
(\ref{simplegeneral_finite}) to produce the following stopping rule:
Continue sampling until \be \la{simplegeneral_finiterev} \li ( \li |
\wh{\bs{p}}_\ell - \f{1}{2} \ri | - w \bs{\ep}_\ell \ri )^2 \geq
\f{1}{4} + \f{ \bs{\ep}_\ell^2 n_\ell } { 2 \ln (\ze \de) } \times
\f{N-1}{N - n_\ell} \ee is satisfied at some stage with index $\ell
\in \{1, \cd, s \}$, where $\bs{\ep}_\ell = \vep \wh{\bs{p}}_\ell$
and $w \geq 0$ is a parameter affecting the shape of the stopping
boundary. As suggested in Section 2.1, the maximum sample size $n_s$
should be defined as the smallest integer such that  the sampling
process is sure to be terminated at or before the $s$-th stage. The
minimum sample size $n_1$ should be defined as the smallest integer
such that the sampling process  has a positive probability to be
terminated at the first stage.  For this revised stopping rule, it
can be shown that the coverage probability $\Pr \left\{ \left |
\wh{\bs{p}} - p \right | \leq \vep p \mid p \right\}$ is no less
than $1 - \de$ for any $p \in \Se$ provided that the coverage tuning
parameter $\ze$ is sufficiently small.

To estimate the population proportion $p \in \Se$ with margin of
absolute error $\vep_a \in (0, 1)$ and margin of relative error
$\vep_r \in (0, 1)$, we can choose the sample sizes $n_1 < n_2 < \cd
< n_s$ as the ascending arrangement of all distinct elements of the
set $\li \{ \li \lc n^\star C_{\tau - \ell} \ri \rc : \ell = 1, \cd,
\tau \ri \}$, where {\small $n^\star = \f{N p^\star ( 1 - p^\star)}{
p^\star ( 1 - p^\star) + (N - 1) \vep_a^2 \sh \mcal{Z}_{\ze \de}^2
}$} with $p^\star = \f{\vep_a}{\vep_r} < \f{1}{2}$. The stopping
rule is that sampling is continued until
\[
\mcal{Z}_{\ze \de}^2 \li ( \f{  N } { n_\ell } - 1 \ri )
\wh{\bs{p}}_\ell ( 1 - \wh{\bs{p}}_\ell ) \leq (N - 1) \max \{
\vep_a^2, \; (\vep_r \wh{\bs{p}}_\ell)^2 \}
\]
is satisfied at some stage with index $\ell$.  Then, $\Pr \{  |
\wh{\bs{p}} - p  | \leq \vep_a \; \tx{or} \; | \wh{\bs{p}} - p  |
\leq \vep_r p \mid p \} \geq 1 - \de$ for any $p \in \Se$ provided
that the coverage tuning parameter $\ze$ is sufficiently small.  In
order to further reduce sampling cost, we propose to revise
(\ref{simplegeneral_finite}) to produce the following stopping rule:
Continue sampling until \be \la{simplegeneral_finitemix} \li ( \li |
\wh{\bs{p}}_\ell - \f{1}{2} \ri | - w \bs{\ep}_\ell \ri )^2 \geq
\f{1}{4} + \f{ \bs{\ep}_\ell^2 n_\ell } { 2 \ln (\ze \de) } \times
\f{N-1}{N - n_\ell} \ee is satisfied at some stage with index $\ell
\in \{1, \cd, s \}$, where $\bs{\ep}_\ell = \max \{ \vep_a, \;
\vep_r \wh{\bs{p}}_\ell \}$ and $w \geq 0$ is a parameter affecting
the shape of the stopping boundary.  The minimum sample size $n_1$
and the maximum sample size $n_s$ can be chosen as suggested in
Section 2.1.  For this revised stopping rule, it can be shown that
$\Pr \{  | \wh{\bs{p}} - p  | \leq \vep_a \; \tx{or} \; |
\wh{\bs{p}} - p  | \leq \vep_r p \mid p \} \geq 1 - \de$ for any $p
\in \Se$ provided that the coverage tuning parameter $\ze$ is
sufficiently small.

Before concluding this section, we want to propose another  class of
stopping rules as follows:

\bei

\item [(i):] Continue sampling until $ \li ( \wh{\bs{p}}_\ell - \f{1}{2} \ri )^2
\geq \f{1}{4} + \f{w}{n_\ell} \f{N - n_\ell}{N - 1} + \f{ \vep^2
n_\ell } { 2 \ln (\ze \de) } \times \f{N-1}{N - n_\ell} $ is
satisfied at some stage with index $\ell \in \{1, \cd, s\}$.  This
stopping rule ensures $\Pr \{ | \wh{\bs{p}} - p  | \leq \vep  \mid p
\} \geq 1 - \de$ for any $p \in \Se$ provided that the coverage
tuning parameter $\ze > 0$ is sufficiently small.

\item  [(ii):] Continue sampling until $\li ( \wh{\bs{p}}_\ell - \f{1}{2} \ri )^2
\geq \f{1}{4} + \f{w}{n_\ell} \f{N - n_\ell}{N - 1} + \f{
\bs{\ep}_\ell^2 n_\ell } { 2 \ln (\ze \de) } \times \f{N-1}{N -
n_\ell} $ is satisfied at some stage with index $\ell \in \{1, \cd,
s\}$, where $\bs{\ep}_\ell = \vep \wh{\bs{p}}_\ell$ and $w \geq 0$
is a parameter affecting the shape of the stopping boundary. This
stopping rule ensures $\Pr \{ | \wh{\bs{p}} - p  | \leq \vep p \mid
p \} \geq 1 - \de$ for any $p \in \Se$ provided that the coverage
tuning parameter $\ze > 0$ is sufficiently small.

\item  [(iii):] Continue sampling until $\li ( \wh{\bs{p}}_\ell - \f{1}{2} \ri )^2
\geq \f{1}{4} + \f{w}{n_\ell} \f{N - n_\ell}{N - 1} + \f{
\bs{\ep}_\ell^2 n_\ell } { 2 \ln (\ze \de) } \times \f{N-1}{N -
n_\ell} $ is satisfied at some stage with index $\ell \in \{1, \cd,
s\}$, where $\bs{\ep}_\ell = \max \{ \vep_a, \; \vep_r
\wh{\bs{p}}_\ell \}$ and $w \geq 0$ is a parameter affecting the
shape of the stopping boundary.  This stopping rule ensures $\Pr \{
| \wh{\bs{p}} - p  | \leq \vep_a \; \tx{or} \; | \wh{\bs{p}} - p  |
\leq \vep_r p \mid p \} \geq 1 - \de$ for any $p \in \Se$ provided
that the coverage tuning parameter $\ze > 0$ is sufficiently small.

\eei

\section{Estimation of Multivariate-Hypergeometric Proportions}

In this section, we propose to construct simultaneous confidence intervals for the parameters of a multivariate hypergeometric distribution by
using techniques similar to the methods for estimating the multinomial proportions as described in Section 6.

 Consider a population of $N$ units, among which there are
 $M_1$ units of category $C_1$,  $M_2$ units of category $C_2$, $\cd$, and
 $M_\ka$ units of category $C_\ka$, where $\sum_{\ell = 1}^\ka M_\ell = N$.
 For $\ell = 1, \cd, \ka$, the proportion of the $\ell$-th category is $p_\ell = \f{M_\ell}{N}$.  It is a frequent
 problem to construct simultaneous confidence intervals for $p_1, \cd, p_\ka$ based on sampling without replacement.
This is equivalent to construct simultaneous confidence intervals for $M_1, \cd, M_\ka$.  For $\ell = 1, \cd, \ka$, let $X_\ell$ denote the
number of units of the $\ell$-th category among $n$ units drawn from the population by sampling without replacement.  The vector $\bs{X} = (X_1,
\cd, X_\ka)$ follows a multivariate hypergeometric distribution with parameters $n, \; N$ and $\bs{p}$, where $\bs{p} \DEF ( \f{M_1}{N}, \cd,
\f{M_{\ka - 1}}{N} )$. Define
\[
\Se \DEF \li \{ \li ( \f{M_1}{N}, \cd, \f{M_{\ka - 1}}{N} \ri ): \;  M_1, \cd, M_{\ka - 1} \; \tx{are nonnegative integers such that} \;
\sum_{\ell = 1}^{\ka - 1} M_\ell \leq N  \ri \}. \] For $\bs{p} \in \Se$, the probability mass function of the multivariate hypergeometric
distribution is given by \bee \Pr \{ X_\ell = x_\ell, \; \ell = 1, \cd \ka \} =  \f{ \prod_{\ell=1}^\ka \bi{N p_\ell}{x_\ell}  } { \bi{N}{n} },
\eee where $x_1, \cd, x_\ka$ are integers such that $\sum_{\ell = 1}^\ka x_\ell = n$ and $0 \leq x_\ell \leq N p_\ell, \; \ell = 1, \cd, \ka$.
In many applications, it is desirable to construct simultaneous confidence intervals $[ \mcal{L}_\ell (\wh{p}_\ell), \; \mcal{U}_\ell
(\wh{p}_\ell) ], \; \ell = 1, \cd, \ka$ such that
\[
\Pr \{ \mcal{L}_\ell (\wh{p}_\ell) \leq  p_\ell \leq  \mcal{U}_\ell (\wh{p}_\ell), \; \ell = 1, \cd, \ka \mid \bs{p} \}  >  1 - \de
\]
for any $\bs{p} \in \Se$, where $\de \in (0, 1)$ is a pre-specified confidence parameter and $\wh{p}_\ell \DEF \f{X_\ell}{n}$ for $\ell = 1,
\cd, \ka$.

As in the case of estimating multinomial proportions, we propose to solve this problem by virtue of the coverage tuning technique with main
ideas as follows:

(i) Seek a class of simultaneous confidence intervals $[ \mcal{L}_\ell (\wh{p}_\ell), \; \mcal{U}_\ell (\wh{p}_\ell) ], \; \ell = 1, \cd, \ka$
such that the coverage probability, denoted by $P(\ze, \bs{p}) \DEF \Pr \{ \mcal{L}_\ell (\wh{p}_\ell) \leq  p_\ell \leq \mcal{U}_\ell
(\wh{p}_\ell), \; \ell = 1, \cd, \ka \mid \bs{p} \}$, can be controlled by the coverage tuning parameter $\ze > 0$.

(ii) Choose the lower and upper confidence limits $\mcal{L}_\ell (\wh{p}_\ell)$ and $\mcal{U}_\ell (\wh{p}_\ell)$ to be nondecreasing functions
of $\wh{p}_\ell$ for $\ell = 1, \cd, \ka$.

(iii) Use Adapted Branch and Bound method in Appendix \ref{BBA} to determine whether the complementary coverage probability $\Psi (\ze, \bs{p})
\DEF  1 - P(\ze, \bs{p})$ of the simultaneous confidence intervals associated with $\ze$ is no greater than $\de$ for any $\bs{p} \in \Se$.

(iv) Apply bisection coverage tuning method to determine $\ze > 0$ as large as possible such that $P(\ze, \bs{p}) \geq 1 - \de$ for any $\bs{p}
\in \Se$.

In order to construct confidence intervals fulfilling the above requirements (i) and (ii), one typical method is based on normal approximation.
Another well-known method is similar to constructing binomial confidence intervals proposed by Clopper and Pearson (1934). The later approach
can be described as follows.

For $\ell = 1, \cd, \ka$, define \bee &  & \mcal{L}_\ell (\wh{p}_\ell) = \min \li \{ p \in \bb{P}: \f{ \sum_{i = X_\ell}^n \bi{N p}{i} \bi{ N -
N p}{ n - i } }{ \bi{N}{n} }
> \f{\ze \de}{2} \ri \},\\
&  & \mcal{U}_\ell (\wh{p}_\ell) = \max \li \{ p \in \bb{P}: \f{ \sum_{i = 0}^{X_\ell} \bi{N p}{i} \bi{ N - N p}{ n - i } }{ \bi{N}{n} }
> \f{\ze \de}{2} \ri \},
\eee where $\bb{P} \DEF \{ \f{i}{N}: i = 0, 1, 2, \cd, N \}$.   Clearly, $\Pr \{ \mcal{L}_\ell (\wh{p}_\ell) \leq p_\ell \leq \mcal{U}_\ell
(\wh{p}_\ell) \mid p_\ell \} \geq 1 - \ze \de$ for $\ell = 1, \cd, \ka$.  It follows from Bonferroni's inequality that the simultaneous
confidence intervals $[ \mcal{L}_\ell (\wh{p}_\ell), \; \mcal{U}_\ell (\wh{p}_\ell) ], \; \ell = 1, \cd, \ka$ have a coverage probability no
less than $1 - \ka \ze \de$ for all $\bs{p} \in \Se$.

Since the structure of the simultaneous confidence intervals can be determined so that the coverage probability $P(\ze, \bs{p})$ is controlled
by $\ze$, we can apply the bisection coverage tuning method to obtain $\ze > 0$ as large as possible such that $P(\ze, \bs{p}) \geq 1 - \de$ for
any $\bs{p} \in \Se$.  A critical subroutine for bisection coverage tuning is to determine, via Adapted Branch and Bound technique in Appendix
\ref{BBA}, whether a given $\ze > 0$ is sufficiently small to guarantee that the coverage probability $P(\ze, \bs{p})$ of the simultaneous
confidence intervals associated with $\ze$ is no less than $1 - \de$ for any $\bs{p} \in \Se$.  This subroutine requires computable bounds of
$\Psi(\ze, \bs{p})$ for $\bs{p} \in \Se$ in a hypercube $\mcal{Q} = \{ (p_1, \cd, p_{\ka - 1}): p_\ell \in [\udl{p}_\ell, \ovl{p}_\ell] \cap
\bb{P}, \; \ell = 1, \cd, \ka - 1 \}$, where $\udl{p}_\ell, \; \ovl{p}_\ell \in \bb{P}, \; 0 \leq \udl{p}_\ell \leq \ovl{p}_\ell \leq 1, \; \ell
= 1, \cd, \ka - 1$ and $\sum_{\ell = 1}^{\ka - 1} \udl{p}_\ell \leq 1$.  In the sequel, we will establish such desired bounds.

Let $\bs{a}= (a_1, \cd, a_\nu)$ and $\bs{b} = (b_1, \cd, b_\nu)$ be integer-valued vectors. Let $\bs{\se} = (\se_1, \cd, \se_\nu)$ be a vector
of nonnegative elements $\se_i \in \bb{P}, \; i = \cd, \nu$.  Define multivariate function
\[
S (\bs{a}, \bs{b}, \bs{\se}, n, N) = \sum_{x_1 = \max \{0, a_1\} }^{\min \{n, b_1\}} \cd \sum_{x_\nu = \max \{0, a_\nu \} }^{\min \{n, b_\nu\} }
 \f{ \prod_{\ell=1}^\ka \bi{N \se_\ell}{x_\ell}  } { \bi{N}{n} } \; \bb{I} \li ( \sum_{\ell=1}^n x_\ell = n \ri ),
\]
where $\bb{I} (\mcal{E})$ denotes the indicator function for the event $\mcal{E}$.  Assume that $\mcal{U}_\ell(1) \geq 1$ and $\mcal{L}_\ell(0)
\leq 0$ for $\ell = 1, \cd, \ka$.  For $\ell = 1, \cd, \ka$, since the lower confidence limit $\mcal{L}_\ell(\wh{p}_\ell)$ is a nondecreasing
function of $\wh{p}_\ell$, we can define inverse function $\mcal{L}_\ell^{-1} (.)$ of $ \mcal{L}_\ell (.)$ such that $\mcal{L}_\ell^{-1} (\se) =
\max \{ z \in I_{\wh{p}_\ell}: \mcal{L}_\ell (z) \leq \se \}$ for $\se \in [0, 1]$, where $I_{\wh{p}_\ell}$ denotes the support of
$\wh{p}_\ell$. Similarly, for $\ell = 1, \cd, \ka$, since the upper confidence limit $\mcal{U}_\ell(\wh{p}_\ell)$ is a nondecreasing function of
$\wh{p}_\ell$, we can define inverse function $\mcal{U}_\ell^{-1} (.)$ of $ \mcal{U}_\ell (.)$ such that $\mcal{U}_\ell^{-1} (\se) = \min \{ z
\in I_{\wh{p}_\ell}: \mcal{U}_\ell (z) \geq \se \}$ for $\se \in [0, 1]$. Let $\udl{p}_\ka = \max \{ 0, 1 - \sum_{\ell = 1}^{\ka - 1}
\ovl{p}_\ell \}$ and $\ovl{p}_\ka = 1 - \sum_{\ell = 1}^{\ka - 1} \udl{p}_\ell$.   We propose to bound the complementary coverage probability
$\Psi (\ze, \bs{p}) = 1 - P (\ze, \bs{p})$ by the following result. \beT \la{multibound_finite} Let $\eta \in (0,1)$.  Let $\mrm{T}_{\mrm{lb}}
(., ., .)$ and $\mrm{T}_{\mrm{ub}} (., ., .)$ be multivariate functions defined by (\ref{truna}) and (\ref{trunb}). For $i = 1, \cd, \ka$,
define $\nu_i = \min \{\ka, i + 1 \}$, \bee & & \udl{\se}_{i, \ell} = \udl{p}_\ell, \qqu
\ovl{\se}_{i, \ell} = \ovl{p}_\ell, \qqu \ell = 1, \cd, \nu_i - 1;\\
&  & \udl{\se}_{i, \nu_i} = \max \li \{ 0, \; 1 - \sum_{\ell = 1}^{\nu_i - 1} \ovl{p}_\ell \ri \}, \qqu \ovl{\se}_{i, \nu_i} = 1 - \sum_{\ell =
1}^{\nu_i - 1} \udl{p}_\ell \eee and $\udl{\bs{\se}}_i = (\udl{\se}_{i,1}, \cd, \udl{\se}_{i, \nu_i}), \; \ovl{\bs{\se}}_i = (\ovl{\se}_{i,1},
\cd, \ovl{\se}_{i,\nu_i})$. For $i = 1, \cd, \ka$, define integer-valued vectors {\small \bee \bs{\mcal{A}}_i = (A_{i,1}, \cd, A_{i,\nu_i}), \qu
\bs{\mcal{B}}_i
= (B_{i,1}, \cd, B_{i,\nu_i}), \qu \bs{\mcal{C}}_i = (C_{i,1}, \cd, C_{i,\nu_i}), \qu \bs{\mcal{D}}_i = (D_{i,1}, \cd, D_{i,\nu_i}), & &\\
\bs{\fra{A}}_i = (\fra{A}_{i,1}, \cd, \fra{A}_{i,\nu_i}), \qu \bs{\fra{B}}_i = (\fra{B}_{i,1}, \cd, \fra{B}_{i,\nu_i}), \qu \bs{\fra{C}}_i =
(\fra{C}_{i,1}, \cd, \fra{C}_{i,\nu_i}), \qu \bs{\fra{D}}_i = (\fra{D}_{i,1}, \cd, \fra{D}_{i,\nu_i}) \; &  & \eee} with {\small \bee A_{i,
\ell} = \fra{A}_{i, \ell} = n \times \max \li \{ \mcal{U}_\ell^{-1} (\udl{\se}_{i, \ell}), \;  \mrm{T}_{\mrm{lb}} (\udl{\se}_{i, \ell}, n, \eta)
\ri \}, \qu B_{i, \ell} = \fra{B}_{i, \ell} =  n \times \min \li \{ \mcal{L}_\ell ^{-1} ( \ovl{\se}_{i, \ell} ), \;
\mrm{T}_{\mrm{ub}} (\ovl{\se}_{i, \ell}, n, \eta) \ri \}, &  &\\
C_{i, \ell} = \fra{C}_{i, \ell} = n \times \max \li \{ \mcal{U}_\ell^{-1} (\ovl{\se}_{i, \ell}), \; \mrm{T}_{\mrm{lb}} (\ovl{\se}_{i, \ell}, n,
\eta) \ri \}, \qu
 D_{i, \ell} = \fra{D}_{i, \ell} = n \times \min
\li \{ \mcal{L}_\ell ^{-1} ( \udl{\se}_{i, \ell} ), \; \mrm{T}_{\mrm{ub}} (\udl{\se}_{i, \ell}, n, \eta) \ri \}, &  & \eee} for $i = 2, \cd,
\ka$ and $\ell = 1, \cd, i - 1$; {\small \bee &  & \fra{A}_{i, i} =  n  \mrm{T}_{\mrm{lb}} (\udl{\se}_{i, i}, n, \eta), \qu \fra{C}_{i, i} =  n
\mrm{T}_{\mrm{lb}} (\ovl{\se}_{i, i}, n, \eta), \qu B_{i, i} =  n  \mrm{T}_{\mrm{ub}} (\ovl{\se}_{i, i}, n, \eta), \qu D_{i, i} = n
\mrm{T}_{\mrm{ub}} (\udl{\se}_{i, i}, n, \eta),\\
& & A_{i, i} = \max \li \{ n \mcal{L}_i^{-1} (\udl{\se}_{i, i}) + 1 , \; \fra{A}_{i, i} \ri \}, \qqu C_{i, i}  = \max \li \{ n \mcal{L}_i^{-1}
(\ovl{\se}_{i, i}) + 1, \; \fra{C}_{i, i} \ri \}, \\
 & & \fra{B}_{i, i} = \min \li \{ n \mcal{U}_i^{-1} (\ovl{\se}_{i, i}) - 1 , \; B_{i, i} \ri \}, \qqu \fra{D}_{i, i} = \min \li \{ n \mcal{U}_i^{-1} (\udl{\se}_{i, i})
- 1 , \; D_{i, i} \ri \} \eee} for $i = 1, \cd, \ka$; and {\small \bee &   & A_{i, i+1}  = \fra{A}_{i, i+1} = n \mrm{T}_{\mrm{lb}}
(\udl{\se}_{i, i+1}, n, \eta), \qqu
 B_{i, i+1} = \fra{B}_{i, i+1} =  n  \mrm{T}_{\mrm{ub}} (\ovl{\se}_{i, i+1}, n, \eta), \\
 &  &  C_{i, i+1}  = \fra{C}_{i, i+1} = n \mrm{T}_{\mrm{lb}} (\ovl{\se}_{i, i+1}, n, \eta),  \qqu
 D_{i, i+1} = \fra{D}_{i, i+1} = n  \mrm{T}_{\mrm{ub}} (\udl{\se}_{i, i+1}, n, \eta) \eee} for $i = 1, \cd, \ka - 1$.
 Then, {\small \bee &  & \Psi (\ze, \bs{p}) \geq \sum_{i = 1}^\ka [ S (\bs{\mcal{C}}_i, \bs{\mcal{D}}_i,
\udl{\bs{\se}}_i, n, N) + S (\bs{\fra{C}}_i, \bs{\fra{D}}_i, \udl{\bs{\se}}_i, n, N)],\\
&  & \Psi (\ze, \bs{p}) \leq (2 \ka - 1) \eta  + \sum_{i = 1}^\ka [ S (\bs{\mcal{A}}_i, \bs{\mcal{B}}_i, \ovl{\bs{\se}}_i, n, N) + S
(\bs{\fra{A}}_i, \bs{\fra{B}}_i, \ovl{\bs{\se}}_i, n, N)] \eee} for all $\bs{p} \in \Se \cap \mcal{Q}$.  \eeT

Theorem \ref{multibound_finite} can be proved by mimicking the argument of Theorem \ref{multibound} in Appendix \ref{multibound_app}. As in the
case of estimating multinomial proportions, we have used the truncation method proposed in \cite{Chen1} to reduce the computational complexity
for bounding $\Psi (\ze, \bs{p})$.  To further reduce the computational complexity of the bounds of $\Psi(\ze, \bs{p})$ given by Theorem
\ref{multibound_finite}, we can develop recursive method for computing terms like $S (\bs{a}, \bs{b}, \bs{\se}, n, N)$ with $\bs{\se} = (\se_1,
\cd, \se_\nu)$ and integer-valued vectors $\bs{a} = (a_1, \cd, a_\nu), \; \bs{b} = (b_1, \cd, b_\nu)$, where $\se_i \in \bb{P}, \; n \geq b_i
\geq a_i \geq 0, \; i = \cd, \nu$.  To avoid trivial cases, assume that $\sum_{i = 1}^\nu a_i < n < \sum_{i = 1}^\nu b_i$ and that $b_i
> a_i$ for at least one $i$ among $i = 1, \cd, \nu$.  Notice that \bee  S(\bs{a}, \bs{b}, \bs{\se},
n, N )  &  =  &  \sum_{x_1 = a_1}^{b_1} \cd \sum_{x_\nu = a_\nu}^{b_\nu}  \f{ \prod_{i=1}^\nu \bi{N \se_i}{x_i}  } { \bi{N}{n} } \; \bb{I} \li (
\sum_{i=1}^n x_i = n \ri ) \\
&  = & \f{1}{\bi{N}{n}} \prod_{i=1}^\nu \f{  (N \se_i) !  } {  a_i! \; (N \se_i - a_i)! } \sum_{x_1 = a_1}^{b_1} \cd \sum_{x_\nu =
a_\nu}^{b_\nu}
\prod_{i=1}^\nu  \f{ a_i! \; (N \se_i - a_i)!  }{  x_i! \; (N \se_i - x_i)! } \; \bb{I} \li ( \sum_{i=1}^n x_i = n \ri )\\
&  = &  \f{ \prod_{i=1}^\nu \bi{N \se_i}{a_i}  } { \bi{N}{n} } \sum_{x_1 = a_1}^{b_1} \cd \sum_{x_\nu = a_\nu}^{b_\nu} \prod_{i=1}^\nu \prod_{j
= 1}^{x_i - a_i} \li (  \f{N \se_i + 1}{a_i + j} - 1 \ri )  \; \bb{I} \li ( \sum_{i=1}^n x_i = n \ri ). \eee  By a similar method as used in
\cite{Frey} to derive a recursive algorithm for computing multinomial distributions, it can be readily shown that
\[
S (\bs{a}, \bs{b}, \bs{\se}, n, N) = \f{ \prod_{i=1}^\nu \bi{\se_i N}{a_i}  } { \bi{N}{n} } \sum_{i = 1}^\nu \sum_{j = 1}^{b_i - a_i} P_r (i,
j),
\]
where $r = n - \sum_{i = 1}^\nu a_i$ and $P_r (i, j), \; i = 1, \cd, \nu; \; j = 1, \cd, b_i - a_i$ can be recursively computed by
\[
P_1 (i, j) = \bec \f{\se_i N + 1}{a_i + 1} - 1  & \tx{if} \; j = 1 \; \tx{and} \; b_i > a_i,\\
0 & \tx{otherwise} \eec
\]
for $1 \leq i \leq k, \; 1 \leq j \leq b_i - a_i$;
\[
P_{t + 1} (i, 1) = \li ( \f{\se_i N + 1}{a_i + 1} - 1 \ri )  \sum_{\ell = 1}^{i - 1} \sum_{j = 1}^{ \min \{b_\ell - a_\ell, t\} }  P_t (\ell, j)
\qu \tx{for $1 \leq i \leq k, \; 1 \leq t \leq r - 1$};
\]
 and $P_{t + 1} (i, j) = \li ( \f{\se_i N + 1}{a_i + j} - 1 \ri )  P_t (i, j - 1)$ for $1 \leq t \leq r - 1, \; 1 \leq i
\leq k, \; 1 < j \leq b_i - a_i$.

By virtue of the lower and upper bounds of $\Psi(\ze, \bs{p})$ given by Theorem \ref{multibound_finite},  we can employ Adapted Branch and Bound
technique in Appendix \ref{BBA} to test if $\Psi (\ze, \bs{p})$ is no greater than $\de$ for any $\bs{p} \in \Se$.   Given that it is possible
to check the truth of $\Psi (\ze, \bs{p}) \leq \de, \; \fa \bs{p} \in \Se$, we can apply a bisection search method to determine the coverage
tuning parameter $\ze$ as large as possible such that the coverage probability $P(\ze, \bs{p})$ of the simultaneous confidence intervals
associated with $\ze$ is no less than $1 - \de$ for all $\bs{p} \in \Se$.

As in the case of estimating multinomial proportions, we consider the determination of sample size for estimating the parameters of the
multivariate distribution.  Let $\de \in (0, 1)$ be a pre-specified confidence parameter. For $\ell = 1, \cd, \ka$, let the margin of error for
$p_\ell$ be $\ep_\ell \DEF \max \{ \vep_{a, \ell}, \; \vep_{r, \ell} \; p_\ell \}$,  where $\vep_{a, \ell} \in [0, 1)$ and $\vep_{r, \ell} \in
[0, 1)$.  For pre-specified confidence parameter $\de$ and margins of error $\ep_\ell = \max \{ \vep_{a, \ell}, \; \vep_{r, \ell} \; p_\ell \},
\; \ell = 1, \cd, \ka$, it is desirable to determine a sample size $n$ as small as possible such that $\Pr \{ | \wh{p}_\ell - p_\ell | \leq
\ep_\ell, \; \ell = 1, \cd, \ka \mid \bs{p} \} \geq 1 - \de$ for all $\bs{p} \in \Se$.  As a consequence of identity (\ref{RI_Indentity}), it is
true that $\Pr \{ | \wh{p}_\ell - p_\ell | \leq \ep_\ell, \; \ell = 1, \cd, \ka \mid \bs{p} \}  = \Pr \{ \mcal{L}_\ell (\wh{p}_\ell) \leq p_\ell
\leq \mcal{U}_\ell (\wh{p}_\ell), \; \ell = 1, \cd, \ka \mid \bs{p} \}$, where the simultaneous confidence intervals $[ \mcal{L}_\ell
(\wh{p}_\ell), \; \mcal{L}_\ell (\wh{p}_\ell) ]$ are defined by (\ref{defsim}).  So, the sample size problem is equivalent to finding the
smallest $n$ such that the coverage probability, denoted by $\mcal{P}(n, \bs{p})$, of the simultaneous confidence intervals is no less than $1 -
\de$. For a given sample size $n$, we can apply the above method for bounding $\Psi(\ze, \bs{p})$ to establish lower and upper bounds for $1 -
\mcal{P}(n, \bs{p})$ with respect to $\bs{p}$ in a hypercube. Hence, as checking the truth of $\Psi(\ze, \bs{p}) \leq \de, \; \fa \bs{p} \in
\Se$, we can determine the truth of $1 - \mcal{P}(n, \bs{p}) \leq \de, \; \fa \bs{p} \in \Se$.  Since such a routine can be established,  it is
possible to obtain the smallest sample size $n$ such that $\mcal{P}(n, \bs{p}) \geq 1 - \de, \; \fa \bs{p} \in \Se$ by checking $n$ from small
to sufficiently large.

\section{Taking into Account Prior Information of Parameters}

In many situations, the parameter to be estimated is known to be included in some interval.  This motivates us to propose a general method for
constructing sampling schemes for estimating $\se$ based on prior information that $\se \in \varTheta = [\udl{\se}, \ovl{\se}] \cap \Se$. We
consider the problems of estimating $\se$ for both the absolute error criterion and the relative error criterion.  In many situations, the
sequential method directly working with a single error criterion may lead to extremely large maximum sample size, or unbounded sample size and
stage number.  As will be seen in the sequel, our main idea is to convert the problems associated with these error criteria into an estimation
problem with a mixed error criterion. Such conversion may lead to sampling schemes of {\it finite and substantially smaller maximum sample size
and stage number}. This is the primary reason we propose the method of conversion.

\subsection{Control of Absolute Error}

Given prior information that $\se \in \varTheta = [\udl{\se},
\ovl{\se}] \cap \Se$, where the interval does not contain $0$, it is
a frequent problem to construct a multistage sampling scheme such
that the corresponding estimator $\wh{\bs{\se}}$ for $\se$ ensures
 the absolute error
criterion (\ref{abs888}).  Let $\vep_a = \vep$ and $\vep_r = \f{\vep}{\max(|\ovl{\se}|, |\udl{\se}|)}$, where $\vep > 0$ and $\de \in (0, 1)$
are specified for the absolute error criterion (\ref{abs888}). Then, for arbitrary sampling schemes used to construct the estimator
$\wh{\bs{\se}}$ for $\se$, we have that $\{ | \wh{\bs{\se}} - \se | < \vep \} = \{ | \wh{\bs{\se}} - \se | < \vep_a \; \tx{or} \; |
\wh{\bs{\se}} - \se | < \vep_r |\se| \}$ for any $\se \in \varTheta$. Therefore, the problem of designing a sampling scheme to ensure the
absolute error criterion (\ref{abs888}) can be converted into the problem of constructing a sampling scheme to ensure the  mixed error criterion
(\ref{mix999}), where the later problem has been solved in preceding sections via the construction of sampling schemes with {\it absolutely
bounded and smaller maximum sample sizes and stage numbers} for the parameters of binomial, Poisson, exponential, normal, Gamma, hypergeometric
distributions and the means of bounded random variables.

\subsection{Control of Relative Error}

Given prior information that $\se \in \varTheta = [\udl{\se}, \ovl{\se}] \cap \Se$, where the interval does not contain $0$, it is desirable to
construct a multistage sampling scheme such that  the corresponding estimator $\wh{\bs{\se}}$ for $\se$ guarantees the relative error criterion
(\ref{rev888}).  Let $\vep_a = \vep \min (|\udl{\se}|, |\ovl{\se}|)$ and $\vep_r = \vep$, where $\vep > 0$ and $\de \in (0, 1)$ are specified
for the relative error criterion (\ref{rev888}). Then, for arbitrary sampling schemes used to construct the estimator $\wh{\bs{\se}}$ for $\se$,
we have that $\{ | \wh{\bs{\se}} - \se | < \vep  |\se| \} = \{ | \wh{\bs{\se}} - \se | < \vep_a \; \tx{or} \; | \wh{\bs{\se}} - \se | < \vep_r
|\se| \}$ for any $\se \in \varTheta$. This implies that the problem of designing a sampling scheme to ensure the relative error criterion
(\ref{rev888}) is equivalent to the problem of constructing a sampling scheme to ensure the mixed error criterion (\ref{mix999}), which has been
addressed in preceding sections and thus are not repeated here.

Before concluding this section, we would like to emphasis that as compared to the direct method, the approach of using prior information to
convert the estimation problem with a single error criterion into an equivalent estimation problem with a mixed error criterion always results
in a reduction of maximum sample size.  Frequently, such reduction is of practical importance.

\section{Estimation of Normal Mean}

Let $X$ be a normal random variable of mean $\mu$ and variance
$\si^2$.  In many situations, the variance $\si^2$ is unknown and it
is desirable to estimate $\mu$ with predetermined margin of error
and confidence level based on a sequence of i.i.d. random samples
$X_1, X_2, \cd$ of $X$.

\subsection{Control of Absolute Error}

For {\it a priori} $\vep > 0$, it is useful to construct an
estimator $\wh{\bs{\mu}}$ for $\mu$ such that $\Pr \{ |
\wh{\bs{\mu}} - \mu | < \vep \} > 1 - \de$ for any $\mu \in (- \iy,
\iy)$ and $\si \in (0, \iy)$.

\subsubsection{New Structure of Multistage Sampling}

 Our new multistage sampling method as follows.  Define
\[
\ovl{X}_{n} = \f{ \sum_{i = 1}^{n} X_i } { n}, \qqu S_n =
\sum_{i=1}^{ n } \li ( X_i - \ovl{X}_{n} \ri )^2
\]
for $n = 2, 3, \cd, \iy$.  Let $s$ be a positive number. The
sampling consists of $s + 1$ stages, of which the sample sizes for
the first $s$ stages are chosen as odd numbers $n_\ell = 2 k_\ell +
1, \; \ell = 1, \cd, s$ with $k_1 < k_2 < \cd < k_s$.  Define
$\wh{\bs{\si}}_\ell = \sq{ \f{S_{n_\ell}}{n_\ell - 1} }$ for $\ell =
1, \cd, s$.  Let the coverage tuning parameter $\ze$ be a positive
number less than $\f{1}{2}$. The stopping rule is as follows:

If {\small $n_\ell < ( \wh{\bs{\si}}_\ell \; t _{n_{\ell} - 1, \ze
\de} )^2 \sh \vep^2, \; \ell = 1, \cd, i - 1$} and {\small $n_i \geq
( \wh{\bs{\si}}_i \; t _{n_{i} - 1, \ze \de} )^2 \sh \vep^2$} for
some $i \in \{ 1, \cd, s\}$, then the sampling is stopped at the
$i$-th stage. Otherwise, {\small $\li \lc ( \wh{\bs{\si}}_s  \; t
_{n_{s} - 1, \ze \de} )^2 \sh \vep^2 \ri \rc - n_s$} more samples of
$X$ needs to be taken after the $s$-th stage. The estimator of $\mu$
is defined as $\wh{\bs{\mu}} = \f{ \sum_{i = 1}^{\mbf{n}} X_i } {
\mbf{n} }$, where $\mbf{n}$ is the sample size when the sampling is
terminated.

It should be noted that, in the special case of $s = 1$, the above
sampling scheme reduces to Stein's two-stage procedure \cite{stein}.
Ghosh and Mukhopadhyay have made improvements for the two-stage
procedures (see,  \cite{Gosh},  \cite{Gosh_1980}, \cite{Mukho_1981},
and the references therein).

\beT

\la{Normal_Analytic_Thm}  The following statements hold true.

(I) $\Pr \{ | \wh{\bs{\mu}} - \mu | < \vep \} > 1 - 2 s \ze \de$ for
any $\mu$ and $\si$.

(II) $\lim_{\vep \to 0} \Pr \{ | \wh{\bs{\mu}} - \mu | < \vep \} = 1
- 2 \ze \de$.

(III) $ \bb{E} [ \mbf{n} ] \leq \f{ (\si \; t _{n_{s} - 1, \ze \de}
)^2 }{ \vep^2 } +  n_s$.

(IV) $\limsup_{\vep \to 0}  \bb{E} \li [  \f{\mbf{n}}{ C} \ri ] \leq
\li ( \f{ t _{n_{s} - 1, \ze \de}}{\mcal{Z}_{\ze \de}} \ri )^2$,
where $C = \li ( \f{ \si \; \mcal{Z}_{\ze \de} } { \vep } \ri )^2$.

\eeT

See Appendix \ref{App_Normal_Analytic_Thm} for a proof.

As can be seen from statement (II) of Theorem
\ref{Normal_Analytic_Thm}, to ensure $\Pr \{ | \wh{\bs{\mu}} - \mu |
< \vep \} > 1 - \de$, it suffices to choose the coverage tuning
parameter $\ze$ to be less than $\f{1}{2 s}$.  However, such a
choice is too conservative. To reduce sampling cost, it is possible
to obtain a value of  $\ze$ much greater than $\f{1}{2 s}$ by our
coverage tuning technique.  Such an approach is explored in the
sequel.

\subsubsection{Exact Construction of Sampling Schemes}

To develop an exact computational approach for the determination of
an appropriate value of coverage tuning parameter $\ze$, we need
some preliminary results as follows.

\beT \la{lemchen} Let $1 = k_0 < k_1 < k_2 < \cd$ be a sequence of
positive integers.  Let $0 = z_0 < z_1 < z_2 < \cd$ be a sequence of
positive numbers. Define $h(0, 1) = 1$ and \[
 h(\ell, 1) =
1, \qu h(\ell, m) = \sum_{i = 1}^{k_r} \f{ h(r, i) \; (z_{\ell} -
z_r)^{m - i} } { (m - i)!  }, \qu k_r < m \leq k_{r + 1}, \qu r = 0,
1, \cd, \ell - 1 \] for $\ell = 1, 2, \cd$.  Let $Z_1, Z_2, \cd$ be
i.i.d. exponential random variables with common mean unity.  Then,
\be \la{GRB} \Pr \li \{\sum_{m = 1}^{k_j} Z_m > z_j \; \tx{for} \; j
= 1, \cd, \ell \ri \} =  e^{- z_\ell} \sum_{m = 1}^{k_\ell} h(\ell,
m) \ee for $\ell = 1, 2, \cd$.  Moreover, the following statements
hold true.

(I)  {\small \bee &   & \Pr \li \{ a_j < \sum_{m = 1}^{k_j} Z_m <
b_j \; \tx{for} \; j = 1, \cd, \ell \ri \} \\
& = & \li [ \sum_{i = 1}^{2^{\ell - 1}} \Pr \li \{ \sum_{m =
1}^{k_j} Z_m
> [A_\ell]_{i, j} \; \tx{for} \; j = 1, \cd, \ell \ri \} \ri ] -
\li [ \sum_{i = 1}^{2^{\ell - 1}} \Pr \li \{ \sum_{m = 1}^{k_j} Z_m
> [B_\ell]_{i, j} \; \tx{for} \; j = 1, \cd, \ell \ri \} \ri ],
\eee} where $A_1 = [ a_1 ], \;  B_1 = [ b_1 ]$ and
\[ A_{r + 1} = \bem A_r & a_{r + 1} I_{2^{r - 1} \times 1}\\
B_r & b_{r + 1} I_{2^{r - 1} \times 1} \eem, \qqu
B_{r + 1} = \bem B_r & a_{r + 1} I_{2^{r - 1} \times 1}\\
A_r & b_{r + 1} I_{2^{r - 1} \times 1} \eem, \qqu r = 1, 2, \cd
\] where $I_{2^{r - 1} \times 1}$ represents a column matrix with all $2^{r - 1}$
elements assuming value $1$.

(II)

{\small \bee &  & \Pr \li \{ a_j < \sum_{m = 1}^{k_j} Z_m < b_j \;
\tx{for} \; j = 1, \cd, \ell, \; \sum_{m = 1}^{k_{\ell + 1}} Z_m >
b_{\ell + 1} \ri
\}\\
 & = & \li [ \sum_{i = 1}^{2^{\ell-1}} \Pr \li \{ \sum_{m =
1}^{k_j} Z_m > [E]_{i, j} \; \tx{for} \; j = 1, \cd, \ell + 1 \ri \}
\ri ]  - \li [ \sum_{i = 1}^{2^{\ell-1}} \Pr \li \{ \sum_{m =
1}^{k_j} Z_m > [F]_{i, j} \; \tx{for} \; j = 1, \cd, \ell + 1 \ri \}
\ri ], \qu \eee} where $E = \bem A_\ell & b_{\ell + 1} I_{2^{\ell -
1} \times 1} \eem$ and $F = \bem B_\ell  & b_{\ell + 1} I_{2^{\ell -
1} \times 1} \eem$.

(III) {\small \bee &  & \Pr \li \{ a_j < \sum_{m = 1}^{k_j} Z_m <
b_j \; \tx{for} \; j = 1, \cd, \ell, \; \sum_{m = 1}^{k_{\ell + 1}}
Z_m <
b_{\ell + 1} \ri \}\\
 & = & \Pr \li \{ a_j < \sum_{m = 1}^{k_j} Z_m < b_j
\; \tx{for} \; j = 1, \cd, \ell \ri \} - \Pr \li \{ a_j < \sum_{m =
1}^{k_j} Z_m < b_j \; \tx{for} \; j = 1, \cd, \ell, \; \sum_{m =
1}^{k_{\ell + 1}} Z_m > b_{\ell + 1} \ri \}. \eee}

 \eeT

It should be noted that (\ref{GRB}) is a generalization of the
recursive formulae (47) and (48) of \cite{Robbin} for the case of
multistage estimation.

 For the purpose of computing
appropriate coverage tuning parameter $\ze$, the following results
are useful.

\beT \la{Normal_Main_Thm} Let the sample sizes of the sampling
scheme be odd numbers $n_\ell = 2 k_\ell + 1, \; \ell = 1, \cd, s$,
where $1 = k_0 < k_1 < k_2 < \cd < k_s$. Let $b_0 = 0$ and {\small
$b_\ell = \f{ k_\ell (2 k_\ell + 1) \vep^2 } { (\si \; t _{2 k_\ell,
\ze \de} )^2 }$} for $\ell = 1, \cd, s$.   Define $h(0, 1) = 1, \;
h(\ell, 1) = 1$,
\[
h(\ell, m) = \sum_{i = 1}^{k_r} \f{ h(r, i) \; (b_{\ell} - b_r)^{m -
i} } { (m - i)!  }, \qu k_r < m \leq k_{r + 1}, \qu r = 0, 1, \cd,
\ell - 1
\]
and $H_\ell (\si) = e^{- b_{\ell}} \sum_{m = 1}^{k_\ell} h(\ell, m)$
for $\ell = 1, \cd, s$. Define {\small $c = \f{ n k_s \;  \vep^2 } {
(\si \; t _{2 k_s, \ze \de} )^2 }, \; h^\star (1) = 1$},
\[
h^\star(m) = \sum_{i = 1}^{k_r} \f{ h(r, i) \; (c - b_r)^{m - i} } {
(m - i)!  }, \qu k_r < m \leq k_{r + 1}, \qu r = 0, 1, \cd, s - 1
\]
and $H^\star (\si, n) =  e^{- c} \sum_{m = 1}^{k_s} h^\star (m)$ for
$n \geq n_s$. Then, the following statements hold true.

(I):  $\Pr \{ | \wh{\bs{\mu}} - \mu | \geq \vep \} = 2 \sum_{n \in
\mscr{S}} \li [  1 - \Phi \li (  \f{ \vep \sq{n} } { \si } \ri ) \ri
] \Pr \{ \mbf{n} = n \}$, where $\mscr{S} = \{ n_\ell: 1 \leq \ell
\leq s \} \cup \{ n \in \bb{N}: n > n_s \}$.

(II): $\Pr \{ \mbf{n} = n \} = \bec
H_{\ell - 1} (\si) - H_\ell (\si)  & \tx{for} \; n = n_\ell, \; 1 \leq \ell \leq s,\\
H^\star (\si, n - 1) - H^\star (\si, n)  & \tx{for} \; n > n_s
 \eec$

 where $H_0(\si) \equiv 1$.

(III): For any $\si \in [a, b]$, {\small \bee  &   & \Pr \{ |
\wh{\bs{\mu}} - \mu | \geq \vep \} > 2 \sum_{n \in \mscr{S} \atop{ n
\leq m } } \li [ 1 - \Phi \li ( \f{ \vep \sq{n} } { a } \ri ) \ri ]
\udl{P}_n,\\
&   &  \Pr \{ | \wh{\bs{\mu}} - \mu | \geq \vep \}  < 2 \sum_{n \in
\mscr{S} \atop{ n \leq m } } \li [ 1 - \Phi \li (  \f{ \vep \sq{n} }
{ b } \ri ) \ri ] \ovl{P}_n +  2 \li [ 1 - \Phi \li ( \f{ \vep
\sq{m} } { b } \ri ) \ri ] S_{\mrm{P}} \li (  k_s - 1, \f{ m k_s
\vep^2 } { ( a \; t_{n_s - 1, \ze \de} )^2 } \ri ),
 \eee}
 where
 \bee
 &   & \ovl{P}_n  = \bec
H_{\ell - 1} (b) - H_\ell (a)  & \tx{for} \; n = n_\ell, \; 1 \leq \ell \leq s,\\
H^\star (b, n - 1) - H^\star (a, n)  & \tx{for} \; n > n_s
 \eec\\
 &  & \udl{P}_n  = \bec
H_{\ell - 1} (a) - H_\ell (b)  & \tx{for} \; n = n_\ell, \; 1 \leq \ell \leq s,\\
H^\star (a, n - 1) - H^\star (b, n)  & \tx{for} \; n > n_s
 \eec
\eee and $m > n_s$.

(IV): {\small \bee \bb{E} [ \mathbf{n} ]  & = &  n_1 + \sum_{\ell =
1}^{s-1} (n_{\ell + 1} - n_\ell ) H_\ell(\si) + \sum_{n = n_s}^\iy
H^\star (\si, n)\\
& < &  n_1 + \sum_{\ell = 1}^{s-1} (n_{\ell + 1} - n_\ell ) H_\ell
(\si) + \sum_{n = n_s}^m H^\star (\si, n) + \f{3 (m \ga e)^\up}{\ga
\sq{\up} \; e^{m \ga \up} },  \eee} where $\ga = \f{ \vep^2 }{ ( \si
\;  t_{n_s - 1, \ze \de} )^2}, \; \up = \f{n_s - 1}{2}$ and $m >
\max  \{ \f{1}{\ga}, n_s \}$.

\eeT

See Appendix \ref{App_Normal_Main_Thm} for a proof.

The coverage tuning process requires evaluation of the coverage
probability $\Pr \{ | \wh{\bs{\mu}} - \mu | < \vep \}$ for various
values of $\si$. To reduce the evaluation of coverage probability
with respect to $\si$ to a finite range of $\si$, we have the
following results. \beT \la{Normal_Cut} Let the sample sizes of the
sampling scheme be odd numbers $n_\ell = 2 k_\ell + 1, \; \ell = 1,
\cd, s$, where $1 < k_1 < k_2 < \cd < k_s$. Suppose the coverage
tuning parameter $\ze$ is a positive number less than $\f{1}{2}$.
Then, there exists a unique number $\ovl{\si}$ such that {\small
\[
\sum_{\ell = 1}^{s - 1} \li [ 1 - S_{\mrm{P}} \li (  k_\ell - 1, \f{
n_\ell \; k_\ell \; \vep^2 } { ( \ovl{\si} \; t_{n_s - 1, \ze \de}
)^2 } \ri ) \ri ]  = (1 - 2 \ze) \de \] } and that $\Pr \{ |
\wh{\bs{\mu}} - \mu | \geq \vep \} < \de$ for $\si > \ovl{\si}$.
Similarly, there exists a unique number $\udl{\si}$ such that
\[
1 - \Phi \li ( \f{ \vep \sq{n_1} } { \udl{\si} }  \ri )  +
\sum_{\ell = 1}^{s - 2} \li [ 1 - \Phi \li ( \f{ \vep \sq{n_{\ell +
1} } } { \udl{\si} }  \ri ) \ri ]  S_{\mrm{P}} \li (  k_\ell - 1,
\f{ n_\ell \; k_\ell \; \vep^2 } { ( \udl{\si} \; t_{n_s - 1, \ze
\de} )^2 } \ri ) = \li ( \f{1}{2} - \ze \ri ) \de
\]
and that $\Pr \{ | \wh{\bs{\mu}} - \mu | \geq \vep \} < \de$ for
$\si < \udl{\si}$.  \eeT

\bsk

See Appendix \ref{App_Normal_Cut} for a proof.

\subsection{Control of Relative Error}  For {\it a priori} $\vep > 0$, it is a frequent problem to
construct an estimator $\wh{\bs{\mu}}$ for $\mu$ such that $\Pr \{ |
\wh{\bs{\mu}} - \mu | \leq \vep | \mu | \} \geq 1 - \de$ for any
$\mu \in (- \iy, 0) \cup (0, \iy)$ and $\si \in (0, \iy)$.  For this
purpose, we would like to propose a new sampling method as follows.

\beT  \la{normal_mean_rev}  Define $\de_\ell = \de$ for $1 \leq \ell
\leq \tau$ and $\de_\ell = \de 2^{\tau - \ell}$ for $\ell > \tau$,
where $\tau$ is a positive integer. For $\ell = 1, 2, \cd$, let
$\wh{\bs{\mu}}_\ell = \f{ \sum_{i = 1}^{n_\ell} X_i } { n_\ell }$
and {\small $\wh{\bs{\si}}_\ell = \sq{ \f{1}{n_\ell - 1}
\sum_{i=1}^{ n_\ell } \li ( X_i - \wh{\bs{\mu}}_\ell  \ri )^2}$},
where $n_\ell$ is deterministic and stands for the sample size at
the $\ell$-th stage. Suppose that sampling is continued until $\li |
\wh{\bs{\mu}}_\ell \ri | \geq \f{ t_{n_\ell - 1, \; \ze \de_\ell } }
{ \sq{n_\ell} } \li ( 1 + \f{1}{\vep} \ri ) \wh{\bs{\si}}_\ell$ for
some stage with index $\ell$. Define estimator $\wh{\bs{\mu}} =
\wh{\bs{\mu}}_{\bs{l}}$, where $\bs{l}$ is the index of stage at
which the sampling is terminated. Then, $\Pr \{ \bs{l} < \iy \} = 1$
and {\small $\Pr \li \{ \li | \wh{\bs{\mu}} - \mu \ri | \leq \vep
|\mu|  \ri \} \geq 1 - \de$} for any $\mu \in (- \iy, 0) \cup (0,
\iy)$ and $\si \in (0, \iy)$ provided that $2 (\tau + 1 ) \ze \leq
1$ and $\inf_{\ell > 0} \f{n_{\ell + 1}}{n_\ell} > 1$. \eeT

See Appendix \ref{App_normal_mean_rev} for a proof.

\subsection{Control of Relative and Absolute Errors}

In some situations, it may be appropriate to estimate $\mu$ with a
mixed error criterion specified by $\vep_a
> 0$ and $\vep_r > 0$. In this respect, we have

\beT \la{normalmix}  Define $\de_\ell = \de$ for $1 \leq \ell \leq
\tau$ and $\de_\ell = \de 2^{\tau - \ell}$ for $\ell > \tau$, where
$\tau$ is a positive integer. For $\ell = 1, 2, \cd$, let
$\wh{\bs{\mu}}_\ell = \f{ \sum_{i = 1}^{n_\ell} X_i } { n_\ell }$
and {\small $\wh{\bs{\si}}_\ell = \sq{ \f{1}{n_\ell - 1}
\sum_{i=1}^{ n_\ell } \li ( X_i - \wh{\bs{\mu}}_\ell  \ri )^2}$},
where $n_\ell$ is deterministic and stands for the sample size at
the $\ell$-th stage. Suppose that sampling is continued until $\max
\li ( \vep_a, \f{ \vep_r \li | \wh{\bs{\mu}}_\ell \ri | } {  1 +
\vep_r} \ri ) \geq \f{ t_{n_\ell - 1, \; \ze \de_\ell } } {
\sq{n_\ell} } \; \wh{\bs{\si}}_\ell$ for some stage with index
$\ell$. Define estimator $\wh{\bs{\mu}} = \wh{\bs{\mu}}_{\bs{l}}$,
where $\bs{l}$ is the index of stage at which the sampling is
terminated.  Then, {\small $\Pr \li \{ \li | \bs{\wh{\mu}} - \mu \ri
| < \vep_a \; \tx{or} \; \li | \bs{\wh{\mu}} - \mu \ri | < \vep_r
|\mu| \ri \} \geq 1 - \de$} for any $\mu \in (- \iy, \iy)$ and $\si
\in (0, \iy)$ provided that $2 (\tau + 1 ) \ze \leq 1$ and
$\inf_{\ell > 0} \f{n_{\ell + 1}}{n_\ell} > 1$. \eeT

See Appendix \ref{App_normalmix} for a proof.

\sect{Estimation of Scale Parameters of Gamma Distributions}

In this section, we shall discuss the estimation of the scale
parameter of a Gamma distribution.   In probability theory and
statistics, a random variable $X$ is said to have a gamma
distribution if its density function is of the form
\[
f_X (x) = \frac{x^{k - 1}} { \Gamma(k) \se ^{ k }    }  \exp \li ( -
\frac{x}{\se} \ri ) \;\;\; {\rm for} \;\;\; 0 < x < \infty
\]
where $\se > 0, \; k > 0$ are referred to as the scale parameter and
shape parameter respectively.  Let $X_1, X_2, \cd$ be i.i.d. samples
of $X$. The MLE of the scale parameter $\se$ can be defined as
\[
\wh{\bs{\se}} = \f{ \sum_{i = 1}^n X_i }{n k }.
\]
Let $0 < \vep < 1$ and $0 < \de < 1$.  The goal is determine the
minimum sample size $n$ such that \be \la{expest} \Pr \li \{ \li |
\f{ \wh{\bs{\se}} - \se } { \se  } \ri | < \vep \mid \se \ri \}
> 1 - \de
\ee for any $\se > 0$.  For simplicity of notations, define $Y = n k
\wh{\bs{\se}} = \sum_{i = 1}^n X_i$. Note that $Y$ has a Gamma
distribution of shape parameter $nk$ and scale parameter $\se$.  It
follows that \bee &  &  \Pr \li \{ \li | \f{ \wh{\bs{\se}} - \se } {
\se } \ri | < \vep \mid \se \ri \} = \Pr \{  Y \geq (1 + \vep) n k
\se \mid \se \} + \Pr \{  Y \leq (1 - \vep) n k \se \mid \se \}\\
&  & = \int_{ (1 + \vep) n k \se }^\iy  \frac{x^{n k - 1}} {
\Gamma(nk) \se ^{ n k }    }  \exp \li ( - \frac{x}{\se} \ri ) d x +
\int_0^{ (1 - \vep) n k \se }  \frac{x^{n k - 1}} { \Gamma(nk) \se
^{ n k }    }  \exp \li ( - \frac{x}{\se} \ri ) d x\\
&  & = \int_{ (1 + \vep) n k  }^\iy  \frac{x^{n k - 1}} { \Gamma(nk)
}  \exp \li ( - x \ri ) d x + \int_0^{ (1 - \vep) n k } \frac{x^{n k
- 1}} { \Gamma(nk) }  \exp \li ( - x \ri ) d x \eee for any $\se >
0$.  Therefore, the minimum sample size to ensure (\ref{expest}) is
the minimum integer $n$ such that {\small $\int_{ (1 + \vep) n k
}^\iy \frac{x^{n k - 1}} { \Gamma(nk) }  \exp \li ( - x \ri ) d x +
\int_0^{ (1 - \vep) n k } \frac{x^{n k - 1}} { \Gamma(nk) }  \exp
\li ( - x \ri ) d x > 1 - \de$}, which can be easily computed.

\sect{Exact Bounded-Width Confidence Intervals}

A classical problem in sequential analysis is to construct a
bounded-width confidence interval with a prescribed level of
coverage probability.  Tanaka developed a non-asymptotic method for
constructing bounded-width confidence intervals for the parameters
of the binomial and Poisson distributions \cite{Tanaka}. Although no
approximation is involved, the method is very conservative due to
the bounding techniques employed in the derivation of sequential
confidence intervals. Franz\'{e}n studied the construction of
bounded-width confidence intervals for binomial parameters in his
paper \cite{frazen}. However, no effective method for defining
stopping rules is proposed. In his later paper \cite{frazenB}, he
proposed to construct fixed width confidence intervals based on
sequential probability ratio tests (SPRT). His method can generate
fixed-sample-size confidence intervals based on SPRT. Unfortunately,
he made a fundamental flaw by mistaking that if the width of the
fixed-sample-size confidence interval decreases to be smaller than
the pre-specified length as the number of samples is increasing,
then the fixed-sample-size confidence interval at the termination of
sampling is the overall sequential confidence interval guaranteeing
the desired confidence level.

In this section, we will demonstrate that the general problem of
constructing fixed-width confidence intervals can be solved in our framework of
multistage estimation described in Section \ref{gen_structure}.
Specifically, the problem of constructing a bounded-width confidence interval
can be formulated as the problem of constructing a sequential
random interval with lower limit $\mscr{L} (\wh{\bs{\se}}, \mbf{n})$ and upper
limit $\mscr{U} (\wh{\bs{\se}}, \mbf{n})$ such that
$\mscr{U} (\wh{\bs{\se}}, \mbf{n}) - \mscr{L} (\wh{\bs{\se}}, \mbf{n}) \leq 2 \vep$ and that
$\Pr \{ \mscr{L} (\wh{\bs{\se}}, \mbf{n}) < \se < \mscr{U} (\wh{\bs{\se}}, \mbf{n}) \mid \se \} > 1 - \de$ for any $\se \in \Se$.  For this
purpose, our computational machinery such as bisection coverage tuning and AMCA can be extremely useful.

\subsection{Construction via Coverage Tuning}

As an application of Theorem \ref{Monotone_second}, our general
theory for constructing bounded-width confidence intervals based on
multistage sampling is as follows.

\beC  \la{CorT}

Suppose a multistage sampling scheme satisfies the following
requirements.

(i) For $\ell = 1, \cd, s$,  $\wh{\bs{\se}}_\ell$ is a ULE of $\se$.

(ii) For $\ell = 1, \cd, s$, $\{ \mscr{L}  (\wh{\bs{\se}}_\ell,
\mbf{n}_\ell) \leq \wh{\bs{\se}}_\ell \leq \mscr{U}
(\wh{\bs{\se}}_\ell, \mbf{n}_\ell) \}$ is a sure event.

(iii) For $\ell = 1, \cd, s$, decision variable $\bs{D}_\ell$
assumes value $1$ if $\mscr{U} (\wh{\bs{\se}}_\ell, \mbf{n}_\ell) -
\mscr{L} (\wh{\bs{\se}}_\ell, \mbf{n}_\ell) \leq 2 \vep$ and assumes
value $0$ otherwise.

(iv) {\small $\li \{ \bs{D}_\ell = 1 \ri \} \subseteq \li \{
F_{\wh{\bs{\se}}_\ell} \li ( \wh{\bs{\se}}_\ell, \mscr{U}
(\wh{\bs{\se}}_\ell, \mbf{n}_\ell) \ri ) \leq \ze \de_\ell , \;
G_{\wh{\bs{\se}}_\ell} \li ( \wh{\bs{\se}}_\ell, \mscr{L}
(\wh{\bs{\se}}_\ell, \mbf{n}_\ell) \ri ) \leq \ze \de_\ell \ri \}$}
for $\ell = 1, \cd, s$.

(v) $\{ \mscr{U} (\wh{\bs{\se}}_s, \mbf{n}_s) - \mscr{L}
(\wh{\bs{\se}}_s, \mbf{n}_s) \leq 2 \vep \}$ is a sure event.

Define $\mscr{L} (\wh{\bs{\se}}, \mbf{n}) = \mscr{L}
(\wh{\bs{\se}}_{\bs{l}}, \mbf{n}_{\bs{l}} )$ and $\mscr{U}
(\wh{\bs{\se}}, \mbf{n}) = \mscr{U} (\wh{\bs{\se}}_{\bs{l}},
\mbf{n}_{\bs{l}} )$, where $\bs{l}$ is the index of stage when the
sampling is terminated. Then, $\mscr{U} (\wh{\bs{\se}}, \mbf{n}) -
\mscr{L} (\wh{\bs{\se}}, \mbf{n}) \leq 2 \vep$ and {\small \bee &  &
\Pr \{
 \mscr{L} (\wh{\bs{\se}}, \mbf{n}) \geq \se \mid \se \}
 \leq \sum_{\ell = 1}^s \Pr \{ \mscr{L} (\wh{\bs{\se}}_\ell,
\mbf{n}_\ell) \geq \se, \;  \bs{D}_\ell = 1 \mid \se \} \leq \ze
\sum_{\ell =
1}^s \de_\ell, \\
&  & \Pr \{ \mscr{U} (\wh{\bs{\se}}, \mbf{n}) \leq \se \mid \se \}
\leq \sum_{\ell = 1}^s \Pr \{ \mscr{U} (\wh{\bs{\se}}_\ell,
\mbf{n}_\ell) \leq \se, \; \bs{D}_\ell = 1 \mid \se \} \leq \ze
\sum_{\ell = 1}^s \de_\ell \eee} and $ \Pr \{  \mscr{L}
(\wh{\bs{\se}}, \mbf{n}) < \se  < \mscr{U} (\wh{\bs{\se}}, \mbf{n})
\mid \se \} \geq 1 - 2 \ze \sum_{\ell = 1}^s \de_\ell$ for any $\se
\in \Se$.

\eeC

Corollary \ref{CorT} indicates that if the coverage tuning parameter $\ze > 0$ is sufficiently small, then the coverage probability of the
bounded-width confidence interval described above can be adjusted to be above the desired level.  Actually, in Corollary \ref{CorT}, we have
proposed a general method to construct bounded-width confidence intervals.  Our general stopping rule is: Continue sampling until the difference
between the upper and lower confidence limits is less than the prescribed width.  The confidence limits at the termination of sampling are taken
as the lower and upper bounds of the desired confidence interval.  The coverage probability of the desired confidence interval is guaranteed via
bisection coverage tuning. It should be noted that, in order to simply the stopping boundary, we can use approximate confidence limits of simple
forms. For example, to construct a bounded-width interval of $100 (1 - \de) \%$ confidence level for the binomial parameter $p$ based on i.i.d.
samples $X_1, X_2, \cd$ of Bernoulli variable $X$ of mean $p$ with a multistage sampling scheme of deterministic sample sizes $n_1 < n_2 < \cd <
n_s$, we can use lower and confidence limits \be \la{biap} \mscr{L} (\wh{\bs{p}}_\ell, n_\ell)  =  \wh{\bs{p}}_\ell - \mcal{Z}_{\ze \de} \sq{
\f{ \wh{\bs{p}}_\ell ( 1 - \wh{\bs{p}}_\ell) }{n_\ell} + \f{w}{n_\ell^2} }, \qqu \mscr{U} (\wh{\bs{p}}_\ell, n_\ell)  = \wh{\bs{p}}_\ell +
\mcal{Z}_{\ze \de} \sq{  \f{ \wh{\bs{p}}_\ell ( 1 - \wh{\bs{p}}_\ell) }{n_\ell} + \f{w}{n_\ell^2}  }, \ee where $\ze > 0$ is the coverage tuning
parameter, $w$ is a positive parameter used to adjust the coverage probability, and $\wh{\bs{p}}_\ell = \f{ \sum_{i = 1}^{n_\ell} X_i }{n_\ell}$
for $\ell  = 1, \cd, s$.  Moreover, we can use the confidence interval given by (\ref{best_CI_A}) and (\ref{best_CI_B}).

\subsection{Bounded-Width Confidence Intervals for Binomial Parameters}

In this subsection, we provide concrete multistage sampling schemes
for the construction of bounded-width confidence intervals for
binomial parameters.

\subsubsection{Construction from Clopper-Pearson Intervals}

Making use of Corollary \ref{CorT} and the Clopper-Pearson
confidence interval \cite{Clopper}, we have established the
following sampling scheme. \beC \la{FW1}
   Let $0 < \vep < \f{1}{2}$.  For
$\ell = 1, \cd, s$, let $\mscr{L} ( \wh{\bs{p}}_\ell, n_\ell )$ be
the largest number such that $0 \leq \mscr{L} ( \wh{\bs{p}}_\ell,
n_\ell ) \leq \wh{\bs{p}}_\ell, \; 1 - S_{\mrm{B}} (n_\ell
\wh{\bs{p}}_\ell - 1, n_\ell, \mscr{L} ( \wh{\bs{p}}_\ell, n_\ell ))
\leq \ze \de$ and let $\mscr{U} ( \wh{\bs{p}}_\ell, n_\ell )$ be the
smallest number such that $\wh{\bs{p}}_\ell \leq \mscr{U} (
\wh{\bs{p}}_\ell, n_\ell ) \leq 1, \; S_{\mrm{B}} (n_\ell
\wh{\bs{p}}_\ell, n_\ell, \mscr{U} ( \wh{\bs{p}}_\ell, n_\ell ) )
\leq \ze \de$, where $\wh{\bs{p}}_\ell = \f{ \sum_{i = 1}^{n_\ell}
X_i } {n_\ell}$. For $\ell = 1, \cd, s$, define $\bs{D}_\ell$ such
that $\bs{D}_\ell = 1$ if $\mscr{U} (\wh{\bs{p}}_\ell, n_\ell) -
\mscr{L} (\wh{\bs{p}}_\ell, n_\ell ) \leq 2 \vep$; and $\bs{D}_\ell
= 0$ otherwise. Suppose the stopping rule is that sampling is
continued until $\bs{D}_\ell = 1$ for some $\ell \in \{1, \cd, s\}$.
Suppose that $\{ \mscr{U} (\wh{\bs{p}}_s, n_s) - \mscr{L}
(\wh{\bs{p}}_s, n_s) \leq 2 \vep \}$ is a sure event. Define {\small
$\mscr{L} ( \wh{\bs{p}}, \mbf{n}  ) = \mscr{L} (
\wh{\bs{p}}_{\bs{l}}, n_{\bs{l}} )$} and {\small $\mscr{U} (
\wh{\bs{p}}, \mbf{n} ) = \mscr{U} ( \wh{\bs{p}}_{\bs{l}}, n_{\bs{l}}
)$} with $\wh{\bs{p}} = \wh{\bs{p}}_{\bs{l}}$ and $\mbf{n} =
n_{\bs{l}}$, where $\bs{l}$ is the index of stage when the sampling
is terminated.   Then, $\mscr{U} (\wh{\bs{p}}, \mbf{n}) - \mscr{L}
(\wh{\bs{p}}, \mbf{n}) \leq 2 \vep$,
 {\small \bee &  & \Pr \{
 \mscr{L} (\wh{\bs{p}}, \mbf{n}) \geq p \mid p \}
 \leq \sum_{\ell = 1}^s \Pr \{ \mscr{L} (\wh{\bs{p}}_\ell,
n_\ell) \geq p, \;  \bs{D}_\ell = 1 \mid p \} \leq s \ze \de, \\
&  & \Pr \{ \mscr{U} (\wh{\bs{p}}, \mbf{n}) \leq p \mid p \} \leq
\sum_{\ell = 1}^s \Pr \{ \mscr{U} (\wh{\bs{p}}_\ell, n_\ell) \leq p,
\; \bs{D}_\ell = 1 \mid p \} \leq s \ze \de
 \eee}
and $ \Pr \{  \mscr{L} (\wh{\bs{\se}}, \mbf{n}) < p  < \mscr{U}
(\wh{\bs{\se}}, \mbf{n}) \mid p \} \geq 1 - 2 s \ze \de$ for any $p
\in (0, 1)$.

\eeC

Based on the criteria proposed in Section \ref{gen_structure}, the
maximum sample size $n_s$ can be defined as the smallest integer
such that $\{ \mscr{U} (\wh{\bs{p}}_s, n_s) - \mscr{L}
(\wh{\bs{p}}_s, n_s) \leq 2 \vep \}$ is a sure event.

\subsubsection{Construction from Fishman's Confidence Intervals}

Making use of Corollary \ref{CorT} and Chernoff-Hoeffding
inequalities \cite{Chernoff, Hoeffding}, we have established the
following sampling scheme.

\beC \la{FW2}  Let $0 < \vep < \f{1}{2}$. Suppose the sample size at
the $s$-th stage is no less than {\small $ \li \lc \f{ \ln \f{1}{\ze
\de} } { 2 \vep^2 } \ri \rc$}. For $\ell = 1, \cd, s$, let $\mscr{L}
( \wh{\bs{p}}_\ell, n_\ell )$ be the largest number such that $0
\leq \mscr{L} ( \wh{\bs{p}}_\ell, n_\ell ) \leq \wh{\bs{p}}_\ell, \;
\mscr{M}_{\mrm{B}} \li ( \wh{\bs{p}}_\ell, \mscr{L} (
\wh{\bs{p}}_\ell, n_\ell ) \ri )  \leq \f{ \ln (\ze \de) } {
n_\ell}$ and let $\mscr{U} ( \wh{\bs{p}}_\ell, n_\ell )$ be the
smallest number such that $\wh{\bs{p}}_\ell \leq \mscr{U} (
\wh{\bs{p}}_\ell, n_\ell ) \leq 1, \; \mscr{M}_{\mrm{B}} \li (
\wh{\bs{p}}_\ell, \mscr{U} ( \wh{\bs{p}}_\ell, n_\ell ) \ri ) \leq
\f{ \ln (\ze \de) } { n_\ell}$, where $\wh{\bs{p}}_\ell = \f{
\sum_{i = 1}^{n_\ell} X_i } {n_\ell}$.  For $\ell = 1, \cd, s$,
define $\bs{D}_\ell$ such that $\bs{D}_\ell = 1$ if $\mscr{U}
(\wh{\bs{p}}_\ell, n_\ell) - \mscr{L} (\wh{\bs{p}}_\ell, n_\ell)
\leq 2 \vep$; and $\bs{D}_\ell = 0$ otherwise.  Suppose the stopping
rule is that sampling is continued until $\bs{D}_\ell = 1$ for some
$\ell \in \{1, \cd, s\}$. Define {\small $\mscr{L} ( \wh{\bs{p}},
\mbf{n}  ) = \mscr{L} ( \wh{\bs{p}}_{\bs{l}}, n_{\bs{l}} )$} and
{\small $\mscr{U} ( \wh{\bs{p}}, \mbf{n}  ) = \mscr{U} (
\wh{\bs{p}}_{\bs{l}}, n_{\bs{l}} )$} with $\wh{\bs{p}} =
\wh{\bs{p}}_{\bs{l}}$ and $\mbf{n} = n_{\bs{l}}$, where $\bs{l}$ is
the index of stage when the sampling is terminated.  Then, $\mscr{U}
(\wh{\bs{p}}, \mbf{n}) - \mscr{L} (\wh{\bs{p}}, \mbf{n}) \leq 2
\vep$,
 {\small \bee &  & \Pr \{
 \mscr{L} (\wh{\bs{p}}, \mbf{n}) \geq p \mid p \} \leq \sum_{\ell = 1}^s
 \Pr \{ \mscr{L} (\wh{\bs{p}}_\ell,
n_\ell) \geq p, \;  \bs{D}_\ell = 1 \mid p \} \leq s \ze \de, \\
&  & \Pr \{ \mscr{U} (\wh{\bs{p}}, \mbf{n}) \leq p \mid p \} \leq
\sum_{\ell = 1}^s \Pr \{ \mscr{U} (\wh{\bs{p}}_\ell, n_\ell) \leq p,
\; \bs{D}_\ell = 1 \mid p \} \leq s \ze \de
 \eee}
and $ \Pr \{  \mscr{L} (\wh{\bs{\se}}, \mbf{n}) < p  < \mscr{U}
(\wh{\bs{\se}}, \mbf{n}) \mid p \} \geq 1 - 2 s \ze \de$ for any $p
\in (0, 1)$.

\eeC

Based on the criteria proposed in Section \ref{gen_structure}, the
sample sizes $n_1 < n_2 < \cd < n_s$ can be chosen as the ascending
arrangement of all distinct elements of {\small \be \la{BCID}
 \li \{ \li \lc  \f{ C_{\tau - \ell} \; \ln
\f{1}{\ze \de} } { 2 \vep^2 } \ri \rc : \ell = 1, \cd, \tau \ri \},
\ee} where $\tau$ is the maximum integer such that {\small $\f{
C_{\tau - 1} \; \ln \f{1}{\ze \de} } { 2 \vep^2 } \geq \f{ \ln (\ze
\de) } { \ln (1 - 2 \vep) }$}, i.e., {\small $C_{\tau - 1} \geq \f{
2 \vep^2 } { \ln \f{1}{1 - 2 \vep} }$}.

\subsubsection{Construction from Explicit Confidence Intervals of Chen et al.}

The following sampling scheme is developed based on Corollary
\ref{CorT} and the explicit confidence intervals due to Chen et al
\cite{Chen_CI}.

 \beC
\la{FW3} Let $0 < \vep < \f{3}{4}$. Suppose the sample size at the
$s$-th stage is no less than {\small $\lc \f{8}{9} ( \frac{3}{4
\vep} + 1  ) ( \frac{3}{4 \vep} - 1 ) \ln \frac{1}{\ze \de} \rc$}.
For $\ell = 1, \cd, s$, define $\wh{\bs{p}}_\ell = \f{\sum_{i =
1}^{n_\ell} X_i}{n_\ell}$ and $\bs{D}_\ell$ such that $\bs{D}_\ell =
1$ if {\small $1 - \frac{9 n_\ell}{ 2 \ln (\ze \de) } \;
\wh{\bs{p}}_\ell \li ( 1 - \wh{\bs{p}}_\ell \ri ) \leq \vep^2 \li [
\f{4}{3} - \frac{3 n_\ell}{ 2 \ln (\ze \de) } \ri ]^2$}, and
$\bs{D}_\ell = 0$ otherwise. Suppose the stopping rule is that
sampling is continued until $\bs{D}_\ell = 1$ for some $\ell \in
\{1, \cd, s\}$. Define {\small \bee  &  &  \mscr{L}
(\wh{\bs{p}}_\ell, n_\ell) = \max \li \{ 0, \; \wh{\bs{p}}_\ell +
\frac{3}{4} \; \frac{ 1 - 2 \wh{\bs{p}}_\ell - \sqrt{ 1 - \frac{9
n_\ell}{ 2 \ln (\ze \de) } \; \wh{\bs{p}}_\ell (
1 - \wh{\bs{p}}_\ell) } } {1 - \frac{9 n_\ell}{ 8 \ln (\ze \de) } } \ri \}, \\
&  & \mscr{U} (\wh{\bs{p}}_\ell, n_\ell) = \min \li \{ 1, \;
\wh{\bs{p}}_\ell + \frac{3}{4} \; \frac{ 1 - 2 \wh{\bs{p}}_\ell +
\sqrt{ 1 - \frac{9 n_\ell}{ 2 \ln (\ze \de) } \; \wh{\bs{p}}_\ell (
1 - \wh{\bs{p}}_\ell) } } {1 - \frac{9 n_\ell}{ 8 \ln (\ze \de) } }
\ri \} \eee}  for $\ell = 1, \cd, s$ and $\wh{\bs{p}} =
\wh{\bs{p}}_{\bs{l}}$ and $\mbf{n} = n_{\bs{l}}$, where $\bs{l}$ is
the index of stage when the sampling is terminated. Then, $\mscr{U}
(\wh{\bs{p}}, \mbf{n}) - \mscr{L} (\wh{\bs{p}}, \mbf{n}) \leq 2
\vep$ and {\small \bee &  & \Pr \{
 \mscr{L} (\wh{\bs{p}}, \mbf{n}) \geq p \mid p \}
 \leq \sum_{\ell = 1}^s \Pr \{ \mscr{L} (\wh{\bs{p}}_\ell,
 n_\ell) \geq p, \;  \bs{D}_\ell = 1 \mid p \} \leq s \ze \de, \\
&  & \Pr \{ \mscr{U} (\wh{\bs{p}}, \mbf{n}) \leq p \mid p \} \leq
\sum_{\ell = 1}^s \Pr \{ \mscr{U} (\wh{\bs{p}}_\ell, n_\ell) \leq p,
\; \bs{D}_\ell = 1 \mid p \} \leq s \ze \de
 \eee}
for any $p \in (0, 1)$.  \eeC

Based on the criteria proposed in Section \ref{gen_structure}, the
sample sizes $n_1 < n_2 <\cd < n_s$ can be chosen as the ascending
arrangement of all distinct elements of {\small $\li \{ \li \lc
C_{\tau - \ell} \li ( \frac{1}{2 \vep^2} - \f{8}{9} \ri ) \ln
\frac{1}{\ze \de} \ri \rc : 1 \leq \ell \leq \tau \ri \}$}, where
$\tau$ is the maximum integer such that {\small $C_{\tau - 1} \li (
\frac{1}{2 \vep^2} - \f{8}{9} \ri ) \ln \frac{1}{\ze \de} \geq \li (
\frac{2}{3 \vep} - \f{8}{9} \ri ) \ln \frac{1}{\ze \de}$}, i.e.,
$C_{\tau - 1} \geq \f{4 \vep}{3 + 4 \vep}$.

\subsection{Bounded-Width Confidence Intervals for Finite Population Proportion}

In this subsection, we consider the construction of bounded-width
confidence intervals for finite population proportion, $p$, based on
multistage sampling.  Within the general framework described in
Sections \ref{gen_structure} and \ref{finite_proportion}, we have
established the following method by virtue of Corollary \ref{CorT}
for bounded-width interval estimation.

\beC \la{Bino_Bounded_CI} For $z \in \{ \f{k}{n}: 0 \leq k \leq n
\}$, define $\mscr{L} (z, n) = \min \{ z, L(z, n) \}$ and $\mscr{U}
(z, n) = \max \{ z, U(z, n) \}$, where {\small $L (z, n) = \min \{
\se \in \Se: 1 - S_N (n z - 1, n, \se) > \ze \de \}$} and {\small $U
( z, n ) = \max \{ \se \in \Se: S_N ( n z, n, \se) > \ze \de \}$}.
Suppose the sample size at the $s$-th stage is no less than the
smallest number $n$ such that $\mscr{U} ( z, n ) - \mscr{L} ( z, n )
\leq 2 \vep$ for all $z \in \{ \f{k}{n}: 0 \leq k \leq n \}$.  For
$\ell = 1, \cd, s$, define $\wh{\bs{p}}_\ell = \f{ \sum_{i =
1}^{n_\ell} X_i } {n_\ell}$ and decision variable $\bs{D}_\ell$
which assumes values $1$ if $\mscr{U} ( \wh{\bs{p}}_\ell, n_\ell ) -
\mscr{L} ( \wh{\bs{p}}_\ell, n_\ell ) \leq  2 \vep$ and value $0$
otherwise. Suppose the stopping rule is that sampling is continued
until $\bs{D}_\ell = 1$ for some $\ell \in \{1, \cd, s\}$. Define
$\wh{\bs{p}} = \wh{\bs{p}}_{\bs{l}}$ and $\mbf{n} = n_{\bs{l}}$,
where $\bs{l}$ is the index of stage when the sampling is
terminated. Then, $\mscr{U} (\wh{\bs{p}}, \mbf{n}) - \mscr{L}
(\wh{\bs{p}}, \mbf{n}) \leq 2 \vep$,  {\small \bee &  & \Pr \{
 \mscr{L} (\wh{\bs{p}}, \mbf{n}) > p \mid p \} =
 \Pr \li \{ \mscr{L} (\wh{\bs{p}}, \mbf{n}) - \f{1}{N} \geq  p \mid p \ri \} \leq \sum_{\ell = 1}^s
 \Pr \li \{ \mscr{L} ( \wh{\bs{p}}_\ell, n_\ell ) - \f{1}{N} \geq p, \; \bs{D}_\ell = 1 \mid p \ri \} \leq s \ze \de, \\
&  & \Pr \{ \mscr{U} ( \wh{\bs{p}}, \mbf{n} ) < p \mid p \} = \Pr
\li \{ \mscr{U} ( \wh{\bs{p}}, \mbf{n} ) + \f{1}{N} \leq p \mid p
\ri \} \leq \sum_{\ell = 1}^s \Pr \li \{ \mscr{U} (
\wh{\bs{p}}_\ell, n_\ell )  + \f{1}{N} \leq p, \; \bs{D}_\ell = 1
\mid p \ri \} \leq s \ze \de
 \eee}
and $ \Pr \{  \mscr{L} (\wh{\bs{p}}, \mbf{n}) \leq p \leq  \mscr{U}
( \wh{\bs{p}}, \mbf{n} ) \} \geq 1 - 2 s \ze \de$ for all $p \in
\Se$.  \eeC

Let $n_{\mrm{max}} $ be the smallest number $n$ such that $\mscr{U}
( z, n ) - \mscr{L} ( z, n ) \leq 2 \vep$ for all $z \in \{
\f{k}{n}: 0 \leq k \leq n \}$. Let $n_{\mrm{min}} $ be the largest
number $n$ such that $\mscr{U} ( z, n) - \mscr{L} ( z, n) > 2 \vep$
for all $z \in \{ \f{k}{n}: 0 \leq k \leq n \}$.  Based on the
criteria proposed in Section \ref{gen_structure}, the sample sizes
$n_1 < n_2 < \cd < n_s$ can be chosen as the ascending arrangement
of all distinct elements of {\small $ \li \{ \li \lc C_{\tau - \ell}
\; n_{\mrm{max}} \ri \rc : \ell = 1, \cd, \tau \ri \}$}, where
$\tau$ is the maximum integer such that $C_{\tau - 1} \geq
\f{n_{\mrm{min}}}{ n_{\mrm{max}} }$.

In order to develop multistage sampling schemes with simple stopping
boundaries, we have the following results.

 \beC
\la{Bino_Bounded_CI_LR} For $z \in \{ \f{k}{n}: 0 \leq k \leq n \}$,
define $\mscr{L} (z, n) = \min \{ z, L(z, n) \}$ and $\mscr{U} (z,
n) = \max \{ z, U(z, n) \}$, where {\small $L (z, n) = \min \{ \se
\in \Se: \mcal{C} (z, \se, n, N) > \ze \de \}$} and {\small $U ( z,
n ) = \max \{ \se \in \Se: \mcal{C} ( z, \se, n, N ) > \ze \de \}$},
where $\mcal{C} (z, \se, n, N)$ is defined by (\ref{defc}). Suppose
the sample size at the $s$-th stage is no less than the smallest
number $n$ such that $\mscr{U} ( z, n ) - \mscr{L} ( z, n ) \leq 2
\vep$ for all $z \in \{ \f{k}{n}: 0 \leq k \leq n \}$. For $\ell =
1, \cd, s$, define $\wh{\bs{p}}_\ell = \f{ \sum_{i = 1}^{n_\ell} X_i
} {n_\ell}$ and decision variable $\bs{D}_\ell$ which assumes values
$1$ if $\mscr{U} ( \wh{\bs{p}}_\ell, n_\ell ) - \mscr{L} (
\wh{\bs{p}}_\ell, n_\ell ) \leq  2 \vep$ and value $0$ otherwise.
Suppose the stopping rule is that sampling is continued until
$\bs{D}_\ell = 1$ for some $\ell \in \{1, \cd, s\}$. Define
$\wh{\bs{p}} = \wh{\bs{p}}_{\bs{l}}$ and $\mbf{n} = n_{\bs{l}}$,
where $\bs{l}$ is the index of stage when the sampling is
terminated. Then, $\mscr{U} (\wh{\bs{p}}, \mbf{n}) - \mscr{L}
(\wh{\bs{p}}, \mbf{n}) \leq 2 \vep$,  {\small \bee &  & \Pr \{
 \mscr{L} (\wh{\bs{p}}, \mbf{n}) > p \mid p \} =
 \Pr \li \{ \mscr{L} (\wh{\bs{p}}, \mbf{n}) - \f{1}{N} \geq  p \mid p \ri \} \leq \sum_{\ell = 1}^s
 \Pr \li \{ \mscr{L} ( \wh{\bs{p}}_\ell, n_\ell ) - \f{1}{N} \geq p, \; \bs{D}_\ell = 1 \mid p \ri \} \leq s \ze \de, \\
&  & \Pr \{ \mscr{U} ( \wh{\bs{p}}, \mbf{n} ) < p \mid p \} = \Pr
\li \{ \mscr{U} ( \wh{\bs{p}}, \mbf{n} ) + \f{1}{N} \leq p \mid p
\ri \} \leq \sum_{\ell = 1}^s \Pr \li \{ \mscr{U} (
\wh{\bs{p}}_\ell, n_\ell )  + \f{1}{N} \leq p, \; \bs{D}_\ell = 1
\mid p \ri \} \leq s \ze \de
 \eee}
and $ \Pr \{  \mscr{L} (\wh{\bs{p}}, \mbf{n}) \leq p \leq  \mscr{U}
( \wh{\bs{p}}, \mbf{n} ) \} \geq 1 - 2 s \ze \de$ for all $p \in
\Se$.  \eeC

Corollary \ref{Bino_Bounded_CI_LR} can be shown by virtue of
Corollary \ref{CorT} and inequalities (\ref{hyp1}) and (\ref{hyp2}).

\section{Interval Estimation Based on a Given Sampling Scheme}

In some situations, the sampling scheme is given and the primary
task is to construct a confidence interval for the parameter $\se$
of the underlying distribution. This class of problems fall into the
category of post-experimental analysis. A typical example is the
construction of confidence interval for $\se$ following hypothesis
testing.   In this direction, we have established general interval
estimation methods in the sequel.

\subsection{Confidence Intervals from Inverting Sequential Hypothesis Tests} \la{CPCI}

Define cumulative distribution functions  $F_{\wh{\bs{\se}}} (z,
\se)$ and $G_{\wh{\bs{\se}}} (z, \se)$ as (\ref{defFG}). To
construct a confidence interval of Clopper-Pearson type for a given
sampling scheme, we have the following results.

\beT \la{Classic_CI} Let $\wh{\bs{\se}} = \varphi (X_1, \cd,
X_{\mbf{n}})$ be a ULE of $\se$, where $\mbf{n}$ is the sample
number at the termination of the sampling process. Define confidence
limits $\mscr{L} (\wh{\bs{\se}}, \mbf{n})$ and $\mscr{U}
(\wh{\bs{\se}}, \mbf{n})$ as functions of $(\wh{\bs{\se}}, \mbf{n})$
such that $\{F_{\wh{\bs{\se}}} ( \wh{\bs{\se}},  \mscr{U}
(\wh{\bs{\se}}, \mbf{n}) ) \leq \f{\de}{2}, \; G_{\wh{\bs{\se}}}
(\wh{\bs{\se}}, \mscr{L} (\wh{\bs{\se}}, \mbf{n}) ) \leq \f{\de}{2}
\}$ is a sure event. Then, $\Pr \{ \mscr{L} (\wh{\bs{\se}}, \mbf{n}
) < \se < \mscr{U} (\wh{\bs{\se}}, \mbf{n} ) \mid \se \} \geq 1 -
\de$ for any $\se \in \Se$. \eeT

See Appendix \ref{App_Classic_CI} for a proof.  Armitage had
considered the problem of interval estimation following hypothesis
testing for binomial case \cite{Armitage_1958}. It should be noted
that, by virtue of our computational machinery, exact computation of
confidence intervals is possible for common distributions.

To construct a confidence interval for the proportion of a finite
population after a multistage test in the general framework
described in Sections \ref{gen_structure} and
\ref{finite_proportion}, we propose the following approach.

\beC \la{CI_Finite_proportion} Let $\wh{\bs{p}} = \f{ \sum_{i =
1}^{\mbf{n}} X_i }{ \mbf{n}  }$, where $\mbf{n}$ is the sample
number at the termination of the sampling process.  Define
confidence limits $\mscr{L} (\wh{\bs{p}}, \mbf{n})$ and $\mscr{U}
(\wh{\bs{p}}, \mbf{n})$ as functions of $(\wh{\bs{p}}, \mbf{n})$
such that, for any observation $( \wh{p}, n )$ of $( \wh{\bs{p}},
\mbf{n} )$, $\mscr{L} (\wh{p}, n)$ is the smallest number in $\Se$
satisfying $\Pr \{ \wh{\bs{p}} \geq \wh{p} \mid \mscr{L} (\wh{p}, n)
\} > \f{\de}{2}$ and that $\mscr{U} (\wh{p}, n)$ is the largest
number in $\Se$ satisfying $\Pr \{ \wh{\bs{p}} \leq \wh{p} \mid
\mscr{U} (\wh{p}, n)  \} > \f{\de}{2}$. Then, $\Pr \{ \mscr{L}
(\wh{\bs{p}}, \mbf{n} ) \leq p \leq \mscr{U} (\wh{\bs{p}}, \mbf{n} )
\mid p \} \geq 1 - \de$ for any $p \in \Se$. \eeC

To show Corollary \ref{CI_Finite_proportion}, it suffices to make
use of Theorem \ref{Classic_CI} and the following observations:

(i)  $\wh{\bs{p}}$ is a ULE of $p$;

(ii)  The procedure for constructing the confidence interval ensures
that $\{F_{\wh{\bs{\se}}} ( \wh{\bs{\se}}, \mscr{U} (\wh{\bs{\se}},
\mbf{n}) + \f{1}{N} ) \leq \f{\de}{2}, \; G_{\wh{\bs{\se}}}
(\wh{\bs{\se}}, \mscr{L} (\wh{\bs{\se}}, \mbf{n}) - \f{1}{N} ) \leq
\f{\de}{2} \}$ is a sure event.

(iii) $\{ \mscr{L} (\wh{\bs{p}}, \mbf{n} ) \leq p \leq \mscr{U}
(\wh{\bs{p}}, \mbf{n} ) \} = \{ \mscr{L} (\wh{\bs{p}}, \mbf{n} ) -
\f{1}{N} < p <  \mscr{U} (\wh{\bs{p}}, \mbf{n} ) + \f{1}{N} \}$ for
$p \in \Se$.

\subsection{Confidence Intervals from Coverage Tuning}  \la{CI_fixed}

The method of interval estimation described in Section \ref{CPCI}
suffers from two drawbacks: (i) It is conservative due to the
discrete nature of the underlying variable.  (ii) There is no
closed-form formula for the confidence interval.  In light of this
situation, we shall propose an alternative approach as follows.

Actually, it is possible to define an expression for the confidence
interval such that the lower confidence limit $\mscr{L}$ and upper
confidence limit $\mscr{U}$ are functions of confidence parameter
$\de$, coverage tuning parameter $\ze$ and $\wh{\bs{\se}} =
\wh{\bs{\se}}_{\bs{l}}$, where $\bs{l}$ is the index of stage when
the sampling is terminated and $\wh{\bs{\se}}_\ell, \; \ell = 1,
\cd, s$ are ULEs as defined in Theorem \ref{Classic_CI}. Suppose
$\mscr{L} (\wh{\bs{\se}}, \mbf{n} ) < \wh{\bs{\se}} < \mscr{U}
(\wh{\bs{\se}}, \mbf{n} )$ and
\[
\Pr \{  \se \leq \mscr{L} (\wh{\bs{\se}}_\ell, \mbf{n}_\ell ) \mid
\se \} \leq  \ze \de_\ell, \qqu \Pr \{  \se \geq  \mscr{U}
(\wh{\bs{\se}}_\ell, \mbf{n}_\ell ) \mid \se \} \leq \ze \de_\ell
\]
for $\ell = 1, \cd, s$.  Then,
\[
\Pr \{  \se \leq  \mscr{L} (\wh{\bs{\se}}, \mbf{n} )  \mid \se \}
\leq \sum_{\ell = 1}^s \Pr \{  \se \leq  \mscr{L}
(\wh{\bs{\se}}_\ell, \mbf{n}_\ell ), \; \bs{D}_{\ell } = 1 \mid \se
\} \leq \ze \sum_{\ell = 1}^s  \de_\ell,
\]
\[
\Pr \{  \se \geq  \mscr{U} (\wh{\bs{\se}}, \mbf{n} )  \mid \se \} \leq \sum_{\ell = 1}^s \Pr \{  \se \geq  \mscr{U} (\wh{\bs{\se}}_\ell,
\mbf{n}_\ell ), \; \bs{D}_{\ell } = 1 \mid \se \} \leq \ze \sum_{\ell = 1}^s \de_\ell. \] This implies that it is possible to apply a bisection
search method to obtain a number $\ze$ such that the coverage probability is no less than $1 - \de$. For the purpose of searching $\ze$, we have
established tight bounds for $\Pr \{ \mscr{L} (\wh{\bs{\se}}, \mbf{n} ) < \se < \mscr{U} (\wh{\bs{\se}}, \mbf{n} ) \mid \se \}$ for $\se \in [a,
b] \subseteq \Se$ as in Section \ref{ITVB}. By virtue of such bounds, adaptive maximum checking algorithm described in Section \ref{AMCA} can be
used to determine an appropriate value of $\ze$.  We would like to point out that, for simplicity, we can use approximate confidence limits of
simple forms like (\ref{best_CI_A}) and (\ref{best_CI_B}).  For example, to construct an interval of $100 (1 - \de) \%$ confidence level for the
binomial parameter $p$ based on i.i.d. samples $X_1, X_2, \cd$ of Bernoulli variable $X$ of mean $p$ following a multistage hypothesis testing
scheme of deterministic sample sizes $n_1 < n_2 < \cd < n_s$, we can use the confidence limits given by (\ref{biap}).  It should be noted that,
although approximate confidence limits are used, the coverage probability of the desired confidence interval is rigorously guaranteed by virtue
of bisection coverage tuning.

\subsubsection{Poisson Mean}

At the first glance, it seems that the approach described at the
beginning of Section \ref{CI_fixed} cannot be adapted to Poisson
variables because the parameter space is not bounded.   To overcome
such difficulty, our strategy is to design a confidence interval
such that, for a large number $\lm^* > 0$,  the coverage probability
is always guaranteed for $\lm \in (\lm^*, \iy)$ without tuning the
confidence parameter and that the coverage probability for $\lm \in
(0, \lm^*]$ can be tuned to be no less than $1 - \de$. Such method
is described in more details as follows.

Suppose the multistage testing plan can be put in the general
framework described in Section \ref{gen_structure}.  Let $\al \in
(0, 1)$ and  $\wh{\bs{\lm}}_\ell = \f{\sum_{i = 1}^{\mbf{n}_\ell}
X_i}{\mbf{n}_\ell}$. For every realization, $(\wh{\lm}_\ell,
n_\ell)$, of $(\wh{\bs{\lm}}_\ell, \mbf{n}_\ell)$, let $L =
L(\wh{\lm}_\ell, n_\ell, \al)$ be the largest number such that
$L(\wh{\lm}_\ell, n_\ell, \al) \leq \wh{\lm}_\ell$ and $\Pr \{
\wh{\bs{\lm}}_\ell \geq \wh{\lm}_\ell \mid L \} \leq \al$. Let $U =
U(\wh{\lm}_\ell, n_\ell, \al)$ be the smallest number such that
$U(\wh{\lm}_\ell, n_\ell, \al) \geq \wh{\lm}_\ell$ and $\Pr \{
\wh{\bs{\lm}}_\ell \leq \wh{\lm}_\ell \mid U \} \leq \al$.  One
possible construction of $L$ and $U$ can be found in \cite{Garwood}.
To eliminate the necessity of evaluating the coverage probability of
confidence interval for an infinitely wide range of parameter $\lm$
in the course of coverage tuning, the following result is crucial.

\beT

\la{Pos_CI_Test}
 Define
\[
\mscr{L} (\wh{\bs{\lm}}_\ell, \mbf{n}_\ell) = \bec
L(\wh{\bs{\lm}}_\ell,
\mbf{n}_\ell, \ze \de ) & \tx{if} \; U(\wh{\bs{\lm}}_\ell, \mbf{n}_\ell, \f{\de}{2s} ) \leq \lm^*,\\
L(\wh{\bs{\lm}}_\ell, \mbf{n}_\ell, \f{\de}{2s} ) & \tx{if} \;
U(\wh{\bs{\lm}}_\ell, \mbf{n}_\ell, \f{\de}{2 s} )
> \lm^* \eec
\]
and
\[
\mscr{U} (\wh{\bs{\lm}}_\ell, \mbf{n}_\ell) = \bec  U
(\wh{\bs{\lm}}_\ell,
\mbf{n}_\ell, \ze \de ) & \tx{if} \; U(\wh{\bs{\lm}}_\ell, \mbf{n}_\ell, \f{\de}{2s} ) \leq \lm^*,\\
U (\wh{\bs{\lm}}_\ell, \mbf{n}_\ell, \f{\de}{2s} ) & \tx{if} \;
U(\wh{\bs{\lm}}_\ell, \mbf{n}_\ell, \f{\de}{2 s} )
> \lm^*. \eec
\]
Let the lower and upper confidence limits be,  respectively, defined
as $\mscr{L} (\wh{\bs{\lm}}, \mbf{n}) = \mscr{L}
(\wh{\bs{\lm}}_{\bs{l}}, \mbf{n}_{\bs{l}})$ and $\mscr{U}
(\wh{\bs{\lm}}, \mbf{n}) = \mscr{U} (\wh{\bs{\lm}}_{\bs{l}},
\mbf{n}_{\bs{l}})$, where $\bs{l}$ is the index of stage when the
sampling is terminated. Then, \be \la{Pos_CI_T}
 \Pr \{  \mscr{L}
(\wh{\bs{\lm}}, \mbf{n}) < \lm  < \mscr{U} (\wh{\bs{\lm}}, \mbf{n})
\mid \lm \} \geq 1 - \de \ee for any $\lm \in (0, \iy)$ provided
that (\ref{Pos_CI_T}) holds for any $\lm \in (0, \lm^*]$. \eeT

See Appendix \ref{App_Pos_CI_Test} for a proof.

\subsubsection{Normal Variance}

A wide class of test plans for the variance of a normal distribution
can be described as follows:

   Choose appropriate sample sizes $n_1 < n_2 <
\cd < n_s$ and numbers $a_\ell < b_\ell, \; \ell = 1, \cd, s$. Let
$\wt{\bs{\si}}_\ell = \sqrt{ \f{1}{n_\ell} \sum_{i = 1}^{n_\ell}
(X_i - \ovl{X}_{n_\ell} )^2}$ for $\ell = 1, \cd, s$.   Continue
sampling until $\wt{\bs{\si}}_\ell \leq a_\ell$ or
$\wt{\bs{\si}}_\ell > b_\ell$. When the sampling is terminated,
accept $\mscr{H}_0$ if $\wt{\bs{\si}}_\ell \leq a_\ell$; reject
$\mscr{H}_0$ if $\wt{\bs{\si}}_\ell > b_\ell$.

To construct a confidence interval for $\si$ after the test, we can
use a ULE of $\si$, which is given by $\wt{\bs{\si}} =
\wt{\bs{\si}}_{\bs{l}}$, where $\bs{l}$ is the index of stage when
the test is completed.  Accordingly,  $\mbf{n} = n_{\bs{l}}$ is the
sample number when the test is completed.   A confidence interval
with lower limit $\mscr{L} ( \wt{\bs{\si}}, \mbf{n} )$ and upper
limit $\mscr{U} ( \wt{\bs{\si}}, \mbf{n} )$ can be constructed as
follows:

If $\wt{\bs{\si}}$ assumes value $\wt{\si}$ at the termination of
test, the realization of the upper confidence limit is equal to a
certain value $\si$ such that $\Pr \{ \wt{\bs{\si}} \leq \wt{\si}
\mid \si \} = \f{ \de }{2}$. Similarly, the realization of the lower
confidence limit is equal to a certain value $\si$ such that $\Pr \{
\wt{\bs{\si}} \geq \wt{\si} \mid \si \} = \f{ \de }{2}$.

To find the value of $\si$ such that $\Pr \{ \wt{\bs{\si}} \leq
\wt{\si} \mid \si \} = \f{ \de }{2}$, it is equivalent to find $\si$
such that \be \la{good189} \Pr \{ \wt{\bs{\si}} \leq \wt{\si} \mid
\si \} = \sum_{\ell = 1}^s \Pr \li \{ \wt{\bs{\si}}_\ell \leq
\wt{\si}, \; a_j < \wt{\bs{\si}}_j \leq b_j, \; 1 \leq j < \ell \mid
\si \ri \}. \ee Similarly, to find the value of $\si$ such that $\Pr
\{ \wt{\bs{\si}} \geq \wt{\si} \mid \si \} = \f{ \de }{2}$, it is
equivalent to find $\si$ such that \be \la{good289} \Pr \{
\wt{\bs{\si}} \geq \wt{\si} \mid \si \} = \sum_{\ell = 1}^s \Pr \li
\{ \wt{\bs{\si}}_\ell \geq \wt{\si}, \; a_j < \wt{\bs{\si}}_j \leq
b_j, \; 1 \leq j < \ell \mid \si \ri \}. \ee If we choose the sample
sizes to be odd numbers $n_\ell = 2 k_\ell + 1, \; \ell = 1, \cd,
s$, we can rewrite (\ref{good189}) and (\ref{good289}) respectively
as {\small \be \la{ssen8} \Pr \{ \wt{\bs{\si}} \leq \wt{\si} \mid
\si \} = \sum_{\ell = 1}^s \Pr \li \{ \sum_{m = 1}^{k_\ell} Z_m \leq
\f{n_\ell}{2} \li ( \f{ \wt{\si} }{\si} \ri )^2, \; \f{n_j}{2} \li (
\f{ a_j }{\si} \ri )^2 < \sum_{m = 1}^{k_j} Z_m \leq \f{n_j}{2} \li
( \f{ b_j }{\si} \ri )^2 \; \tx{for} \; 1 \leq j < \ell \mid \si \ri
\} \ee} and {\small \be \la{ssen28} \Pr \{ \wt{\bs{\si}} \geq
\wt{\si} \mid \si \} = \sum_{\ell = 1}^s \Pr \li \{ \sum_{m =
1}^{k_\ell} Z_m \geq \f{n_\ell}{2} \li ( \f{ \wt{\si} }{\si} \ri
)^2, \; \f{n_j}{2} \li ( \f{ a_j }{\si} \ri )^2 < \sum_{m = 1}^{k_j}
Z_m \leq \f{n_j}{2} \li ( \f{ b_j }{\si} \ri )^2 \; \tx{for} \; 1
\leq j < \ell \mid \si \ri \}, \qu \ee} where $Z_1, Z_2, \cd$ are
i.i.d. exponential random variables with common mean unity. As can
be seen from (\ref{ssen8}) and (\ref{ssen28}), the determination of
confidence interval for $\si$ requires the exact computation of the
probabilities in the right-hand sides of (\ref{ssen8}) and
(\ref{ssen28}).  For such computational purpose, we can use Theorem
\ref{lemchen}.

\subsubsection{Exponential Parameters}

Let $X$ be a random variable of density function $f_X (x) =
\f{1}{\se} \exp \li ( - \f{x}{\se} \ri )$. Let $X_1, X_2, \cd$ be
i.i.d. samples of the exponential random variable $X$.  A wide class
of test plans for the parameter $\se$ of the exponential
distribution can be described as follows:

Choose appropriate sample sizes $n_1 < n_2 < \cd < n_s$ and numbers
$a_\ell < b_\ell, \; \ell = 1, \cd, s$.  Define {\small
$\wh{\bs{\se}}_\ell = \f{ \sum_{i = 1}^{n_\ell} X_i }{n_\ell}$} for
$\ell = 1, \cd, s$. Continue sampling until $\wh{\bs{\se}}_\ell \leq
a_\ell$ or $\wh{\bs{\se}}_\ell > b_\ell$.  When the sampling is
terminated, accept $\mscr{H}_0$ if $\wh{\bs{\se}}_\ell \leq a_\ell$;
reject $\mscr{H}_0$ if $\wh{\bs{\se}}_\ell > b_\ell$.

To construct a confidence interval for $\se$ after the test, we can
use a ULE of $\se$, which is given by $\wh{\bs{\se}} =
\wh{\bs{\se}}_{\bs{l}}$, where $\bs{l}$ is the index of stage when
the test is completed.  Accordingly, $\mbf{n} = n_{\bs{l}}$ is the
sample number when the test is completed.   A confidence interval
with lower limit $\mscr{L} ( \wh{\bs{\se}}, \mbf{n} )$ and upper
limit $\mscr{U} ( \wh{\bs{\se}}, \mbf{n} )$ can be constructed as
follows:

If $\wh{\bs{\se}}$ assumes value $\wh{\se}$ when the test is
completed, the realization of the upper confidence limit is equal to
a certain value $\se$ such that $\Pr \{ \wh{\bs{\se}} \leq \wh{\se}
\mid \se \} = \f{ \de }{2}$. Similarly, the realization of the lower
confidence limit is equal to a certain value $\se$ such that $\Pr \{
\wh{\bs{\se}} \geq \wh{\se} \mid \se \} = \f{ \de }{2}$.

To find the value of $\se$ such that $\Pr \{ \wh{\bs{\se}} \leq
\wh{\se} \mid \se \} = \f{ \de }{2}$, it is equivalent to find $\se$
such that \be \la{good189exp} \Pr \{ \wh{\bs{\se}} \leq \wh{\se}
\mid \se \} = \sum_{\ell = 1}^s \Pr \li \{ \wh{\bs{\se}}_\ell \leq
\wh{\se}, \; a_j < \wh{\bs{\se}}_j \leq b_j, \; 1 \leq j < \ell \mid
\se \ri \}. \ee Similarly, to find the value of $\se$ such that $\Pr
\{ \wh{\bs{\se}} \geq \wh{\se} \mid \se \} = \f{ \de }{2}$, it is
equivalent to find $\se$ such that \be \la{good289exp} \Pr \{
\wh{\bs{\se}} \geq \wh{\se} \mid \se \} = \sum_{\ell = 1}^s \Pr \li
\{ \wh{\bs{\se}}_\ell \geq \wh{\se}, \; a_j < \wh{\bs{\se}}_j \leq
b_j, \; 1 \leq j < \ell \mid \se \ri \}. \ee Let $Z_1, Z_2, \cd$ be
i.i.d. exponential random variables with common mean unity. Then, we
can rewrite (\ref{good189exp}) and (\ref{good289exp}) respectively
as {\small \be \la{ssen8exp} \Pr \{ \wh{\bs{\se}} \leq \wh{\se} \mid
\se \} = \sum_{\ell = 1}^s \Pr \li \{ \sum_{m = 1}^{n_\ell} Z_m \leq
n_\ell \li ( \f{ \wh{\se} }{\se} \ri ), \; n_j \li ( \f{ a_j }{\se}
\ri ) < \sum_{m = 1}^{n_j} Z_m \leq n_j \li ( \f{ b_j }{\se} \ri )
\; \tx{for} \; 1 \leq j < \ell \mid \se \ri \} \ee} and {\small \be
\la{ssen28exp} \Pr \{ \wh{\bs{\se}} \geq \wh{\se} \mid \se \} =
\sum_{\ell = 1}^s \Pr \li \{ \sum_{m = 1}^{n_\ell} Z_m \geq n_\ell
\li ( \f{ \wh{\se} }{\se} \ri ), \; n_j \li ( \f{ a_j }{\se} \ri ) <
\sum_{m = 1}^{n_j} Z_m  \leq n_j \li ( \f{ b_j }{\se} \ri ) \;
\tx{for} \; 1 \leq j < \ell \mid \se \ri \}. \ee} As can be seen
from (\ref{ssen8exp}) and (\ref{ssen28exp}), the determination of
confidence interval for $\si$ requires the exact computation of the
probabilities in the right-hand sides of (\ref{ssen8exp}) and
(\ref{ssen28exp}).  For such computational purpose, we can make use
of the results in Theorem \ref{lemchen}.

\sect{Exact Confidence Sequences}

The construction of confidence sequences is a classical problem in
statistics.  The problem has been studied  by Darling and Robbin
\cite{darling, darling2}, Lai \cite{Lai2}, Jennsion and Turnbull
\cite{Jennsion}, and many other researchers.  In this section, we
shall develop a computational approach for the problem in a general
setting as follows.

Let $X_1, X_2, \cd$ be a sequence of samples of random variable $X$
parameterized by $\se \in \Se$.   Consider a multistage sampling
procedure of $s$ stages such that the number of available samples at
the $\ell$-th stage is a random number $\mbf{n}_\ell$ for $\ell = 1,
\cd, s$.  Let $\wh{\bs{\se}}_\ell$ be a function of random tuple
$X_1, \cd, X_{\mbf{n}_\ell}$ for $\ell = 1, \cd, s$.  The objective
is to construct intervals with lower limits $\mscr{L} (
\wh{\bs{\se}}_\ell, \mbf{n}_\ell )$ and upper limits $\mscr{U} (
\wh{\bs{\se}}_\ell, \mbf{n}_\ell )$ such that
\[
\Pr \{ \mscr{L} ( \wh{\bs{\se}}_\ell, \mbf{n}_\ell ) < \se <
\mscr{U} ( \wh{\bs{\se}}_\ell, \mbf{n}_\ell ), \; \ell = 1, \cd, s
\mid \se \}
> 1 - \de
\]
for any $\se \in \Se$.

\subsection{Construction via Coverage Tuning}

Assume that $\wh{\bs{\se}}_\ell$ is a ULE for $\ell = 1, \cd, s$.
For simplicity of notations, let
\[ L_\ell = \mscr{L} ( \wh{\bs{\se}}_\ell, \mbf{n}_\ell ), \qqu
U_\ell = \mscr{U} ( \wh{\bs{\se}}_\ell, \mbf{n}_\ell ), \qqu \ell =
1, \cd, s.
\]
As mentioned earlier, our objective is to construct a sequence of
confidence intervals $(L_\ell, U_\ell), \; 1 \leq \ell \leq s$ such
that $\Pr \{ L_\ell < \se < U_\ell, \;  1 \leq \ell \leq s \mid \se
\} \geq 1 - \de$ for any $\se \in \Se$.  Suppose
\[
\Pr \{  L_\ell < \se < U_\ell \mid \se \} \geq 1 - \ze \de, \qqu 1
\leq \ell \leq s
\]
for any $\se \in \Se$.  By Bonferroni's inequality, we have $\Pr \{
L_\ell < \se < U_\ell, \;  1 \leq \ell \leq s \mid \se \} \geq 1 - s
\ze \de$ for any $\se \in \Se$.  This implies that it is possible to
find an appropriate value of coverage tuning parameter $\ze$ such
that $\Pr \{ L_\ell < \se < U_\ell, \;  1 \leq \ell \leq s \mid \se
\} \geq 1 - \de$ for any $\se \in \Se$.

For this purpose, it suffices to bound the complementary probability
$1 - \Pr \{  L_\ell < \se < U_\ell, \;  1 \leq \ell \leq s \mid \se
\}$ and apply the adaptive maximum checking algorithm described in
Section \ref{AMCA} to find an appropriate value of the coverage
tuning parameter $\ze$ such that $1 - \Pr \{  L_\ell < \se < U_\ell,
\;  1 \leq \ell \leq s \mid \se \} \leq \de$ for any $\se \in [ a, b
] \subseteq \Se$. In this respect, we have

\beT  \la{RCI_Bounds} Let $X_1, X_2, \cd$ be a sequence of identical samples of discrete random variable $X$ which is parameterized by $\se \in
\Se$. For $\ell = 1, \cd, s$, let $\wh{\bs{\se}}_\ell =\varphi (X_1, \cd, X_{\mbf{n}_\ell})$ be a ULE of $\se$.  Let $L_\ell =
\mscr{L}(\wh{\bs{\se}}_\ell, \mbf{n}_\ell)$ and $U_\ell = \mscr{U}(\wh{\bs{\se}}_\ell, \mbf{n}_\ell)$ be bivariate functions of
$\wh{\bs{\se}}_\ell$ and $\mbf{n}_\ell$ such that $\{ L_\ell \leq \wh{\bs{\se}}_\ell \leq U_\ell  \}, \; \ell = 1, \cd, s$ are sure events. Let
$[a, b]$ be a subset of $\Se$. Let $I_\mscr{L}$ denote the intersection of $[a, b]$ and the union of the supports of $L_\ell, \; \ell = 1, \cd,
s$. Let $I_\mscr{U}$ denote the intersection of $[a, b]$ and the union of the supports of $U_\ell, \; \ell = 1, \cd, s$.  Define \bee &   & P_L
(\se) = \sum_{k = 1}^s \Pr \{ L_k \geq \se, \; L_\ell < \se < U_\ell, \; 1 \leq
\ell < k \mid \se \},\\
&  & P_U (\se) = \sum_{k = 1}^s \Pr \{  U_k \leq \se, \; L_\ell <
\se < U_\ell, \; 1 \leq \ell < k \mid \se \}. \eee The following
statements hold true:

(I): $1 - \Pr \{  L_\ell < \se < U_\ell, \;  1 \leq \ell \leq s \mid
\se \} = P_L (\se) + P_U (\se)$.

(II): $P_L (\se)$ is non-decreasing with respect to $p \in \Se$ in
any interval with endpoints being consecutive distinct elements of
$I_\mscr{L} \cup \{a, b \}$.  The maximum of $P_L (\se)$ over $[a,
b]$ is achieved at $I_\mscr{L} \cup \{a, b \}$. Similarly, $P_U
(\se)$ is non-increasing with respect to $p \in \Se$ in any interval
with endpoints being consecutive distinct elements of $I_\mscr{U}
\cup \{a, b \}$. The maximum of $P_U (\se)$ over $[a, b]$ is
achieved at $I_\mscr{U} \cup \{a, b \}$.

(III): Suppose that $\{ L_\ell \geq a \} \subseteq \{
\wh{\bs{\se}}_\ell \geq b \}$ and $\{ U_\ell \leq b \} \subseteq \{
\wh{\bs{\se}}_\ell \leq a \}$ for $\ell = 1, \cd, s$. Then, \bee P_L
(\se) & \leq & \sum_{k = 1}^s \Pr \{ L_k \geq a, \; L_\ell < b, \;
U_\ell
> a, \; 1 \leq \ell < k
\mid b \}, \\
P_U (\se) &  \leq &  \sum_{k = 1}^s \Pr \{  U_k \leq b, \; L_\ell <
b, \; U_\ell
> a, \;  1 \leq \ell < k \mid a \},\\
P_L (\se) & \geq & \sum_{k = 1}^s \Pr \{  L_k \geq b, \;
L_\ell < a, \; U_\ell > b, \; 1 \leq \ell < k \mid a \}, \\
P_U (\se) &  \geq & \sum_{k = 1}^s \Pr \{  U_k \leq a, \; L_\ell <
a, \; U_\ell > b, \;  1 \leq \ell < k \mid b \} \eee for any $\se
\in [ a, b ] \subseteq \Se$.

\eeT

Theorem \ref{RCI_Bounds} can be established by a similar argument as that of Theorem \ref{Main_Bound_Gen}.  It should be noted that no need to
compute $s$ terms in the summation independently. The recursive method described in Section 3.6 can be used.

We would like to point out that, for simplicity, we can use approximate confidence limits of simple forms like (\ref{best_CI_A}) and
(\ref{best_CI_B}).  For example, to construct a confidence sequence of $100 (1 - \de) \%$ confidence level for the binomial parameter $p$ based
on i.i.d. samples $X_1, X_2, \cd$ of Bernoulli variable $X$ of mean $p$ with a multistage sampling scheme of deterministic sample sizes $n_1 <
n_2 < \cd < n_s$, we can use the confidence limits given by (\ref{biap}). It should be noted that, although approximate confidence limits are
used, the coverage probability of the desired confidence sequence is rigorously guaranteed by virtue of bisection coverage tuning.

\subsection{Finite Population Proportion}

To construct a confidence sequence for the proportion, $p$, of a
finite population described in Section \ref{gen_structure}, we have
the following results.

\beT  \la{Finite_RCI_Bounds} Let $L_\ell =
\mscr{L}(\wh{\bs{p}}_\ell, \mbf{n}_\ell)$ and $U_\ell =
\mscr{U}(\wh{\bs{p}}_\ell, \mbf{n}_\ell)$ be bivariate functions of
{\small $\wh{\bs{p}}_\ell = \f{ \sum_{i = 1}^{\mbf{n}_\ell} X_i } {
\mbf{n}_\ell }$} and $\mbf{n}_\ell$ such that $L_\ell \leq
\wh{\bs{p}}_\ell \leq U_\ell$ and that both $N L_\ell$ and $N
U_\ell$ are integer-valued random variables for $\ell = 1, \cd, s$.
Let $a \leq b$ be two elements of $\Se = \{ \f{m}{N} : m = 0, 1,
\cd, N \}$. Let $I_\mscr{L}$ denote the intersection of interval
$(a, b)$ and the union of the supports of $L_\ell - \f{1}{N}, \;
\ell = 1, \cd, s$. Let $I_\mscr{U}$ denote the intersection of
interval $(a, b)$ and the union of the supports of $U_\ell +
\f{1}{N}, \; \ell = 1, \cd, s$. Define \bee & & P_L (p) = \sum_{k =
1}^s \Pr \{ L_k > p, \; L_\ell \leq p \leq U_\ell, \;
1 \leq \ell < k \mid p \},\\
&  & P_U (p) = \sum_{k = 1}^s \Pr \{  U_k < p, \; L_\ell \leq p \leq
U_\ell, \; 1 \leq \ell < k \mid p \}. \eee The following statements
hold true.

(I): $1 - \Pr \{  L_\ell \leq p \leq U_\ell, \;  1 \leq \ell \leq s
\mid p \} = P_L (p) + P_U (p)$.

(II): $P_L (p)$ is non-decreasing with respect to $p \in \Se$ in any
interval with endpoints being consecutive distinct elements of
$I_\mscr{L} \cup \{a, b \}$.  The maximum of $P_L (p)$ over $[a, b]$
is achieved at $I_\mscr{L} \cup \{a, b \}$. Similarly, $P_U (p)$ is
non-increasing with respect to $p \in \Se$ in any interval with
endpoints being consecutive distinct elements of $I_\mscr{U} \cup
\{a, b \}$. The maximum of $P_U (p)$ over $[a, b]$ is achieved at
$I_\mscr{U} \cup \{a, b \}$.

(III): Suppose that $\{ L_\ell \geq a \} \subseteq \{
\wh{\bs{p}}_\ell \geq b \}$ and $\{ U_\ell \leq b \} \subseteq \{
\wh{\bs{p}}_\ell \leq a \}$ for $\ell = 1, \cd, s$. Then, \bee P_L
(p) & \leq & \sum_{k = 1}^s \Pr \{ L_k > a, \; L_\ell \leq b, \;
U_\ell \geq a, \; 1 \leq \ell < k \mid b \}, \\
P_U (p) &  \leq &  \sum_{k = 1}^s \Pr \{  U_k < b, \; L_\ell \leq
b, \; U_\ell \geq a, \;  1 \leq \ell < k \mid a \},\\
P_L (p) & \geq & \sum_{k = 1}^s \Pr \{  L_k > b,
\; L_\ell \leq a, \; U_\ell \geq b, \; 1 \leq \ell < k \mid a \}, \\
P_U (p) &  \geq & \sum_{k = 1}^s \Pr \{  U_k < a, \; L_\ell \leq a,
\; U_\ell \geq b, \;  1 \leq \ell < k \mid b  \} \eee for any $p \in
[ a, b ] \cap \Se$.

\eeT

Theorem \ref{Finite_RCI_Bounds} can be established by a similar
argument as that of Theorem \ref{Main_Bound_Gen}. It should be noted
that our computational machinery such as bisection coverage tuning,
AMCA and recursive algorithm can be used.

\subsection{Poisson Mean}

At the first glance, it seems that the adaptive maximum checking
algorithm described in Section \ref{AMCA}  cannot be adapted to
Poisson variables because the parameter space is not bounded.  To
overcome such difficulty, our strategy is to design a confidence
sequence such that, for a large number $\lm^* > 0$, the coverage
probability is always guaranteed for $\lm \in (\lm^*, \iy)$ without
tuning the confidence parameter and that the coverage probability
for $\lm \in (0, \lm^*]$ can be tuned to be no less than $1 - \de$.
Such method is described in more details as follows.

Let $\al \in (0, 1)$ and $\wh{\bs{\lm}}_\ell = \f{\sum_{i =
1}^{\mbf{n}_\ell} X_i}{\mbf{n}_\ell}$.  For every realization,
$(\wh{\lm}_\ell, n_\ell)$, of $(\wh{\bs{\lm}}_\ell, \mbf{n}_\ell)$,
let $L = L(\wh{\lm}_\ell, n_\ell, \al)$ be the largest number such
that $L(\wh{\lm}_\ell, n_\ell, \al) \leq \wh{\lm}_\ell$ and $\Pr \{
\wh{\bs{\lm}}_\ell \geq \wh{\lm}_\ell \mid L \} \leq \al$.  Let $U =
U(\wh{\lm}_\ell, n_\ell, \al)$ be the smallest number such that
$U(\wh{\lm}_\ell, n_\ell, \al) \geq \wh{\lm}_\ell$ and $\Pr \{
\wh{\bs{\lm}}_\ell \leq \wh{\lm}_\ell \mid U \} \leq \al$.  One
possible construction of $L$ and $U$ can be found in \cite{Garwood}.
To eliminate the necessity of evaluating the coverage probability of
confidence interval for an infinitely wide range of parameter $\lm$
in the course of coverage tuning, the following result is critical.

\beT \la{goal_RCI_Pos}

Define
\[
\mscr{L} (\wh{\bs{\lm}}_\ell, \mbf{n}_\ell) = \bec
L(\wh{\bs{\lm}}_\ell,
\mbf{n}_\ell, \ze \de ) & \tx{if} \; U(\wh{\bs{\lm}}_\ell, \mbf{n}_\ell, \f{\de}{2s} ) \leq \lm^*,\\
L(\wh{\bs{\lm}}_\ell, \mbf{n}_\ell, \f{\de}{2s} ) & \tx{if} \;
U(\wh{\bs{\lm}}_\ell, \mbf{n}_\ell, \f{\de}{2 s} )
> \lm^* \eec
\]
and
\[
\mscr{U} (\wh{\bs{\lm}}_\ell, \mbf{n}_\ell) = \bec  U
(\wh{\bs{\lm}}_\ell,
\mbf{n}_\ell, \ze \de ) & \tx{if} \; U(\wh{\bs{\lm}}_\ell, \mbf{n}_\ell, \f{\de}{2s} ) \leq \lm^*,\\
U (\wh{\bs{\lm}}_\ell, \mbf{n}_\ell, \f{\de}{2s} ) & \tx{if} \;
U(\wh{\bs{\lm}}_\ell, \mbf{n}_\ell, \f{\de}{2 s} )
> \lm^*. \eec
\]
Then, \be \la{goal_RCI_Pos_Ineq} \Pr \{ \mscr{L}
(\wh{\bs{\lm}}_\ell, \mbf{n}_\ell) < \lm < \mscr{U}
(\wh{\bs{\lm}}_\ell, \mbf{n}_\ell), \; \ell = 1, \cd, s \mid \lm \}
\geq 1 - \de \ee for any $\lm \in (0, \iy)$ provided that
(\ref{goal_RCI_Pos_Ineq}) holds for any $\lm \in (0, \lm^*]$.

\eeT

See Appendix \ref{App_goal_RCI_Pos} for a proof.

\subsection{Normal Mean}

Let $X_1, X_2, \cd$ be i.i.d. samples of Gaussian variable $X$ with mean $\mu$ and variance $\si^2$.  By Bonferroni's inequality,
\[
\Pr \{  \ovl{X}_{n_\ell} - \mcal{Z}_{\ze \de} \; \si \sh \sq{n_\ell}
< \mu < \ovl{X}_{n_\ell} + \mcal{Z}_{\ze \de} \; \si \sh \sq{n_\ell}
, \; 1 \leq \ell \leq s \}
> 1 - s \ze \de,
\]
where  $\ovl{X}_{n_\ell} = \f{\sum_{i=1}^{n_\ell} X_i } {n_\ell}$ for $\ell = 1, \cd, s$. It can be seen that if the coverage tuning parameter
$\ze > 0$ is chosen to be small enough, then
\[
\Pr \{  \ovl{X}_{n_\ell} - \mcal{Z}_{\ze \de} \; \si \sh \sq{n_\ell}
 < \mu < \ovl{X}_{n_\ell} + \mcal{Z}_{\ze \de} \; \si \sh \sq{n_\ell} , \; 1 \leq \ell
\leq s \} = 1 - \de.
\]
To compute the coverage probability of the repeated confidence intervals, there is no loss of generality to assume that $\mu = 0$ and $\si^2 =
1$. Hence, it suffices to compute $\Pr \{ | \ovl{X}_{n_\ell} | <   \mcal{Z}_{\ze \de}  \sh \sq{n_\ell}, \; 1 \leq \ell \leq s \}$.   We shall
evaluate the complementary probability {\small \bee  1 - \Pr \{ | \ovl{X}_{n_\ell} | < \mcal{Z}_{\ze \de}  \sh \sq{n_\ell}, \; 1 \leq \ell \leq
s \} & = & \Pr \{  | \ovl{X}_{n_\ell} | \geq
\mcal{Z}_{\ze \de}  \sh \sq{n_\ell} \; \tx{for some} \; \ell \; \tx{among} \; 1, \cd, s \}\\
& = & \sum_{r = 1}^s  \Pr \{  | \ovl{X}_{n_r} | \geq \mcal{Z}_{\ze
\de}  \sh \sq{n_r} \; \tx{and} \; |\ovl{X}_{n_\ell}
| < \mcal{Z}_{\ze \de}  \sh \sq{n_\ell}, \; 1 \leq \ell < r \}\\
& = & 2 \sum_{r = 1}^s  \Pr \{ \ovl{X}_{n_r} \geq \mcal{Z}_{\ze \de} \sh \sq{n_r} \; \tx{and} \; |\ovl{X}_{n_\ell} | < \mcal{Z}_{\ze \de} \sh
\sq{n_\ell}, \; 1 \leq \ell < r \}. \eee}  The bounding method based on consecutive decision variables described in Section \ref{consec} can be
used. Specifically, \bee & & 1 - \Pr \{ | \ovl{X}_{n_\ell} | < \mcal{Z}_{\ze \de} \sh \sq{n_\ell}, \;
1 \leq \ell \leq s \}\\
& \leq & 2 \sum_{r = 1}^s  \Pr \{ \ovl{X}_{n_r} \geq \mcal{Z}_{\ze
\de} \sh \sq{n_r} \; \tx{and} \; |\ovl{X}_{n_\ell} | < \mcal{Z}_{\ze
\de} \sh \sq{n_\ell}, \; \max(1, r - k) \leq \ell < r \} \eee for $1
\leq k < s$.   Such method can be used for the problem of testing
the equality of the mean response of two treatments (see,
\cite{Pocock}, \cite{TRM} and the references therein).  It can also
be applied to the repeated significance tests established by
Armitage, McPherson, and Rowe \cite{AMR}.

\subsection{Normal Variance}

In this section, we shall discuss the construction of confidence
sequence for the variance of a normal distribution.   Let $X_1, X_2,
\cd$ be i.i.d. samples of a normal random variable $X$ of mean $\mu$
and variance $\si^2$.  Our method of constructing a confidence
sequence is follows.

Choose the sample sizes to be odd numbers $n_\ell = 2 k_\ell + 1, \;
\ell = 1, \cd, s$.  Define {\small $\ovl{X}_{n_\ell} = \f{\sum_{i =
1}^s X_i }{n_\ell}$} and $S_{n_\ell} = \sum_{i = 1}^s (X_i -
\ovl{X}_{n_\ell})^2$ for $\ell = 1, \cd, s$.  Note that
\[
\Pr \li \{  \f{ S_{n_\ell}   } {  \chi_{n_\ell - 1, 1 - \ze \de}^2 }
< \si^2 < \f{ S_{n_\ell}   } {  \chi_{n_\ell - 1, \ze \de}^2 }, \; 1
\leq \ell \leq s \ri \}  > 1 - 2 s \ze \de
\]
and {\small \bee \Pr \li \{  \f{ S_{n_\ell} } { \chi_{n_\ell - 1, 1
- \ze \de  }^2 } <  \si^2 < \f{ S_{n_\ell} } { \chi_{n_\ell - 1, \ze
\de }^2 }, \; 1 \leq \ell \leq s \ri \} & = & \Pr \li \{
\chi_{n_\ell - 1, \ze \de  }^2  < \f{ S_{n_\ell} } {\si^2} <
\chi_{n_\ell - 1, 1 - \ze \de }^2, \; 1 \leq \ell \leq s
\ri \}\\
& = & \Pr \li \{  \chi_{n_\ell - 1, \ze \de  }^2  < \sum_{m =
1}^{k_\ell} Z_m <  \chi_{n_\ell - 1, 1 - \ze \de }^2, \; 1 \leq \ell
\leq s \ri \},  \eee} where $Z_1, \; Z_2, \; \cd $ are i.i.d.
exponential random variables with common mean unity.  Therefore, the
coverage probability {\small $\Pr  \li \{  \f{ S_{n_\ell} } {
\chi_{n_\ell - 1, 1 - \ze \de  }^2 } <  \si^2 < \f{ S_{n_\ell} } {
\chi_{n_\ell - 1, \ze \de }^2 }, \; 1 \leq \ell \leq s \ri \}$} can
be exactly computed by virtue of Theorem \ref{lemchen}.
Consequently, we can obtain, via a bisection search method, an
appropriate value of $\ze$ such that
\[
\Pr \li \{  \f{ S_{n_\ell}   } {  \chi_{n_\ell - 1, 1 - \ze \de }^2
} < \si^2 < \f{ S_{n_\ell}   } {  \chi_{n_\ell - 1, \ze \de }^2 },
\; 1 \leq \ell \leq s \ri \}  = 1 - \de.
\]

\subsection{Exponential Parameters}

In this section, we shall consider the construction of confidence
sequences for the parameter $\se$ of a random variable $X$ of
density function $f_X (x) = \f{1}{\se} \exp \li ( - \f{x}{\se} \ri
)$. Let $X_1, X_2, \cd$ be i.i.d. samples of a normal random
variable $X$.   Let $n_1 < n_2 < \cd < n_s$ be a sequence of sample
sizes. Since $\f{2 n \ovl{X}_{n}}{\se}$ has a chi-square
distribution of $2 n$ degrees of freedom,  we have
\[
\Pr \li \{  \chi_{2 n_\ell, \ze \de}^2  < \f{2 n_\ell
\ovl{X}_{n_\ell}}{\se} < \chi_{2 n_\ell, 1 - \ze \de}^2, \; 1 \leq
\ell \leq s \ri \} > 1 - 2 s \ze \de,
\]
or equivalently,
\[
\Pr \li \{ \f{2 \sum_{i = 1}^{n_\ell} X_i}{ \chi_{2 n_\ell, 1 - \ze
\de }^2 } < \se < \f{2 \sum_{i = 1}^{n_\ell} X_i}{\chi_{2 n_\ell,
\ze \de }^2}, \; 1 \leq \ell \leq s \ri \} > 1 - 2 s \ze \de.
\]
Note that
\[
\Pr \li \{ \f{2 \sum_{i = 1}^{n_\ell} X_i}{ \chi_{2 n_\ell, 1 - \ze
\de }^2 } < \se < \f{2 \sum_{i = 1}^{n_\ell} X_i}{\chi_{2 n_\ell,
\ze \de }^2}, \; 1 \leq \ell \leq s \ri \} = \Pr \li \{  \f{ \chi_{2
n_\ell, \ze \de}^2 }{2} < \sum_{i = 1}^{n_\ell} Z_i < \f{ \chi_{2
n_\ell, 1 - \ze \de}^2 }{2}, \; 1 \leq \ell \leq s \ri \},
\]
where $Z_1, \; Z_2, \; \cd $ are i.i.d. exponential random variables
with common mean unity.  Therefore, the coverage probability {\small
$\Pr \li \{ \f{2 \sum_{i = 1}^{n_\ell} X_i}{ \chi_{2 n_\ell, 1 - \ze
\de }^2 } < \se < \f{2 \sum_{i = 1}^{n_\ell} X_i}{\chi_{2 n_\ell,
\ze \de }^2}, \; 1 \leq \ell \leq s \ri \}$} can be exactly computed
by virtue of Theorem \ref{lemchen}. Consequently,  we can obtain,
via a bisection search method, an appropriate value of $\ze$ such
that
\[
\Pr \li \{ \f{2 \sum_{i = 1}^{n_\ell} X_i}{ \chi_{2 n_\ell, 1 - \ze
\de }^2 } < \se < \f{2 \sum_{i = 1}^{n_\ell} X_i}{\chi_{2 n_\ell,
\ze \de }^2}, \; 1 \leq \ell \leq s \ri \} = 1 - \de.
\]

\subsection{Confidence Sequences for Infinite Stages of Sampling}

It should be noted that the theory and techniques of constructing confidence sequences for finite stages of sampling described in preceding
discussion may be applied to the case that the number of stages is infinite, i.e., $s = \iy$.  Note that $\Pr \{ L_\ell < \se < U_\ell \;
\tx{for all} \; \ell = 1, 2, \cd \mid \se \} \leq \Pr \{ L_\ell < \se < U_\ell, \; 1 \leq \ell \leq
\tau  \mid \se \}$ and {\small \bel  &  & \Pr \{ L_\ell < \se < U_\ell \; \tx{for all} \; \ell = 1, 2, \cd \mid \se \} \nonumber \\
&  & \geq \Pr \{ L_\ell < \se < U_\ell, \; 1 \leq \ell \leq \tau  \mid \se \} - \Pr \{ \se \notin (L_k, U_k) \; \tx{for some} \; k > \tau  \mid
\se \}  \la{cool} \eel} for $\tau \geq 1$. This implies that if we can bound $\Pr \{ \se \notin (L_k, U_k) \; \tx{for some} \; k > \tau  \mid
\se \}$ and make it extremely small as compared to $\de$, then the construction of the confidence sequence for infinite stages of sampling can
be treated as the construction of the confidence sequence for a sampling process of finite stages.

As our first illustration, consider the construction of a confidence sequence for the binomial parameter $p$. More formally, let $X_1, X_2, \cd$
be i.i.d. samples of Bernoulli random variable $X$ such that $\Pr \{ X = 1 \} = 1 - \Pr \{ X = 0 \} = p \in (0, 1)$.  Let $\{ n_\ell, \; \ell =
1, 2, \cd \}$ be a sequence of sample sizes such that $n_1 < n_2 < n_3 < \cd$. Define $\wh{p}_\ell = \f{\sum_{i = 1}^{n_\ell} X_i}{n_\ell}$ for
$\ell = 1, 2, \cd$.  Let $0 < \vep < 1$ and $0 < \de < 1$.  We are interested in determining the tightest lower bound, $m$,  for the sample
size, $n_1$, of the first stage such that
\[
\Pr \{ | \wh{p}_\ell - p | < \vep \; \tx{for all} \; \ell = 1, 2, \cd  \mid p \} > 1 - \de
\]
for any $p \in (0, 1)$ provided that $n_1 \geq m$.   Applying (\ref{cool})
with $L_\ell = \wh{p}_\ell - \vep$ and $U_\ell = \wh{p}_\ell + \vep$,
we have \bee &  & 1 - \Pr \{ | \wh{p}_\ell - p | < \vep \; \tx{for all} \; \ell = 1, 2, \cd  \mid p \}\\
&  & \leq \Pr \{ | \wh{p}_k - p | \geq  \vep \; \tx{for some} \; k > \tau \mid p \} + \sum_{k = 1}^\tau \Pr \{ | \wh{p}_k - p | \geq  \vep,  \;
| \wh{p}_\ell - p | < \vep, \; 1 \leq \ell < k  \mid p \}. \eee for $\tau \geq 1$.  From martingale theory, it is well known that $\Pr \{ |
\wh{p}_k - p | \geq \vep \; \tx{for some} \; k > \tau \mid p \} \leq 2 \exp ( - 2 n_{\tau + 1}  \vep^2 )$.  Note that for any $\eta \in (0, 1)$,
it is possible to find stage index $\tau$ large enough such that $2 \exp ( - 2 n_{\tau + 1} \vep^2 )  < \eta$.  Hence, it suffices to ensure
that
\[
\sum_{k = 1}^\tau \Pr \{ | \wh{p}_k - p | \geq  \vep,  \; | \wh{p}_\ell - p | < \vep, \; 1 \leq \ell < k  \mid p \} < \de - \eta
\]
for all $p \in (0, 1)$.  This requirement can be fulfilled by choosing $n_1 \geq m$ with a sufficiently large $m$.  Since the number of stages
involved in the above summation is finite, the computational complexity is manageable by modern computers.  Making use of statement (III) of
Theorem \ref{RCI_Bounds}, we can obtain recursively computable bounds of $\sum_{k = 1}^\tau \Pr \{ | \wh{p}_k - p | \geq  \vep,  \; |
\wh{p}_\ell - p | < \vep, \; 1 \leq \ell < k \mid p \}$ for $p$ in an interval $[a, b] \subseteq (0, 1)$.  Using such bounds, the adapted Branch
and Bound algorithm or AMCA can be used to quickly determine whether $n_1 \geq m$ with a given $m$ is sufficient to guarantee that the coverage
probability is greater than $1 - \de$.  Therefore, we can determine the smallest $m$ such that the desired confidence level is rigorously
guaranteed by virtue of bisection coverage tuning.

As another illustration, consider the construction of a confidence sequence for the mean of a normal random variable.  Let $X_1, X_2, \cd$ be
i.i.d. samples of  a Gaussian random variable $X$ with unknown mean $\mu$ and known variance $\si^2$.  Without loss of generality, assume that
the variance $\si^2$ is equal to $1$.  Let $\{ n_\ell, \; \ell = 1, 2, \cd \}$ be a sequence of sample sizes such that $n_1 < n_2 < n_3 < \cd$.
Define $\wh{\mu}_\ell = \f{\sum_{i = 1}^{n_\ell} X_i}{n_\ell}$ for $\ell = 1, 2, \cd$. Let $\vep > 0$ and $0 < \de < 1$.  We are interested in
determining the tightest lower bound, $m$, for the sample size, $n_1$, of the first stage such that
\[
\Pr \{ | \wh{\mu}_\ell - \mu | < \vep \; \tx{for all} \; \ell = 1, 2, \cd  \mid \mu \} > 1 - \de
\]
for any $\mu \in (- \iy, \iy)$ provided that $n_1 \geq m$.  Making use of (\ref{cool})
with $L_\ell = \wh{\mu}_\ell - \vep$ and $U_\ell = \wh{\mu}_\ell + \vep$,
we have \bee &  & 1 - \Pr \{ | \wh{\mu}_\ell - \mu | < \vep \; \tx{for all} \; \ell = 1, 2, \cd  \mid \mu \}\\
&  & \leq \Pr \{ | \wh{\mu}_k - \mu | \geq  \vep \; \tx{for some} \; k > \tau \mid \mu \} + \sum_{k = 1}^\tau \Pr \{ | \wh{\mu}_k - \mu | \geq
\vep,  \; | \wh{\mu}_\ell - \mu | < \vep, \; 1 \leq \ell < k  \mid \mu \} \eee for $\tau \geq 1$.  From martingale theory, it is well known that
$\Pr \{ | \wh{\mu}_k - \mu | \geq  \vep \; \tx{for some} \; k > \tau \mid \mu \} \leq 2 \exp ( - \f{n_{\tau + 1}  \vep^2}{2} )$.  Note that for
any $\eta \in (0, 1)$, it is possible to find stage index $\tau$ large enough such that $2 \exp ( - \f{n_{\tau + 1}  \vep^2}{2} )   < \eta$.
Hence, it suffices to ensure that
\[
\sum_{k = 1}^\tau \Pr \{ | \wh{\mu}_k - \mu | \geq  \vep,  \; | \wh{\mu}_\ell - \mu | < \vep, \; 1 \leq \ell < k  \mid \mu \} < \de - \eta
\]
for all $\mu \in (- \iy, \iy)$.  Therefore, as in the construction of confidence sequence for the binomial parameter $p$, we can determine the
smallest $m$ such that the desired confidence level is rigorously guaranteed by virtue of bisection coverage tuning.

\section{Multistage Linear Regression}

Regression analysis is a statistical technique for investigating and
modeling the relationship between variables. Applications of
regression are numerous and occur in almost every field, including
engineering, physical sciences, social sciences, economics,
management, life and biological sciences, to name but a few.
Consider a linear model
\[
y = \ba_1 x_1 + \ba_2 x_2 + \cd + \ba_m x_m + w \qu \tx{with} \qu
x_1 \equiv 1,
\]
where $\ba_1, \cd, \ba_m$ are deterministic parameters and $w$ is a
Gaussian random variable of zero mean and variance $\sigma^2$.    A
major task of linear regression is to estimate parameters $\sigma$
and $\ba_i$ based on observations of $y$ for various values of
$x_i$. In order to strictly control estimation error and uncertainty
of inference with as few observations as possible, we shall develop
multistage procedures. To this end, we shall first define some
variables.  Let $\boldsymbol{\beta} = [\beta_1, \cdots,
\beta_m]^\intercal$, where the notation ``$\intercal$'' stands for
the transpose operation. Let $w_1, w_2, \cd$ be a sequence of i.i.d.
samples of $w$. Define
\[
y_i = \ba_1 x_{i1} + \ba_2 x_{i2} + \cd + \ba_m x_{im} + w_i \qu
\tx{with} \qu x_{i1} \equiv 1
\]
for $i = 1, 2, \cd$.  Let $n_\ell, \; \ell = 1, 2, \cd$ be a
sequence of positive integers which is ascending with respect to
$\ell$.  Define
\[
\boldsymbol{Y}_\ell = \bem y_1\\
y_2\\
 \vdots\\
  y_{n_\ell} \eem,
\qqu  \boldsymbol{X}_\ell  = \bem
x_{11} & x_{12} & \cdots & x_{1m}\\
x_{21} & x_{22} & \cdots & x_{2m}\\
\vdots & \vdots & \ddots & \vdots\\
x_{n_\ell 1} & x_{n_\ell 2} & \cdots & x_{n_\ell m} \eem \qqu
\tx{for} \; \ell = 1, 2, \cd.
\]
Assume that $\boldsymbol{X}_\ell^\intercal \boldsymbol{X}_\ell$ is
of rank $m$ for all $\ell$. Define
\[
\boldsymbol{B}_\ell = (\boldsymbol{X}_\ell^\intercal
\boldsymbol{X}_\ell)^{-1} \boldsymbol{X}_\ell^\intercal
\boldsymbol{Y}_\ell, \qqu \widehat{\bs{\si}}_\ell = \sq{
\frac{1}{n_\ell - m} \; \left [ \boldsymbol{Y}_\ell^\intercal
\boldsymbol{Y}_\ell - \boldsymbol{B}_\ell^\intercal
(\boldsymbol{X}^\intercal \boldsymbol{Y}_\ell )\right ] }
\]
for $\ell = 1, 2, \cd$.   For $i = 1, \cd, m$, let $\bs{B}_{i, \ell}
$ denote the $i$-th entry of $\bs{B}_{\ell}$ and let $\left[
(\boldsymbol{X}_\ell^\intercal \boldsymbol{X}_\ell)^{-1}
\right]_{ii}$ denote the $(i, i)$-th entry of
$(\boldsymbol{X}_\ell^\intercal \boldsymbol{X}_\ell)^{-1}$.

\subsection{Control of Absolute Error}

For the purpose of estimating the variance $\si$ and the parameters
$\ba_i$ with an absolute error criterion, we have

\beT \la{Regression_Abs}  Let $\vep > 0$ and $\vep_i > 0$ for $i =
1, \cd, m$.  Let $\tau$ be a positive integer. Suppose the process
of observing $y$ with respect to $x_i$ and $w$ is continued until
{\small $t_{n_\ell - m, \; \ze \de_\ell} \; \widehat{\bs{\si}}_\ell
\; \sqrt{\left[ (\boldsymbol{X}_\ell^\intercal
\boldsymbol{X}_\ell)^{-1} \right]_{ii}} \leq \vep_i$} for $i = 1,
\cd, m$, and
\[
\sqrt{ \frac{ n_\ell - m} { \chi_{n_\ell - m, \; \ze \de_\ell }^2 }
} \; \widehat{\bs{\si}}_\ell - \vep \leq \widehat{\bs{\si}}_\ell
\leq \sqrt{ \frac{ n_\ell - m } { \chi_{n_\ell - m, \; 1 - \ze
\de_\ell }^2 } } \; \widehat{\bs{\si}}_\ell + \vep
\]
at some stage with index $\ell$,  where $\de_\ell = \de$ for $1 \leq
\ell \leq \tau$ and $\de_\ell = \de 2^{\tau - \ell}$ for $\ell >
\tau$.  Define $\wh{\bs{\si}} = \widehat{\bs{\si}}_{\bs{l}}$ and
$\wh{\bs{\ba}} = \bs{B}_{\bs{l}}$, where $\bs{l}$ is the index of
stage at which the observation of $y$ is stopped.  For $i = 1, \cd,
m$, let $\wh{\bs{\ba}}_i$ be the $i$-th entry of $\wh{\bs{\ba}}$.
Then, $\Pr \{ \bs{l} < \iy \} = 1$ and $\Pr \{ | \wh{\bs{\si}} -
\sigma | \leq \vep, \;\; | \wh{\bs{\ba}}_i - \beta_i | \leq \vep_i
\; \tx{for} \; i = 1, \cdots, m \} \geq 1 - \delta$ provided that $2
(m+1) (\tau + 1) \ze \leq 1$ and that $\inf_{\ell > 0} \f{n_{\ell +
1}}{n_\ell} > 1$.

\eeT

See Appendix \ref{App_Regression_Abs}  for a proof.

\subsection{Control of Relative Error}

For the purpose of estimating the variance $\si$ and the parameters
$\ba_i$ with a relative error criterion, we have

\beT \la{Regression_Rev} Let $0 < \vep < 1$ and $0 < \vep_i < 1$ for
$i = 1, \cd, m$.  Let $\tau$ be a positive integer. Suppose the
process of observing $y$ with respect to $x_i$ and $w$ is continued
until {\small $t_{n_\ell - m, \; \ze \de_\ell} \;
\widehat{\bs{\si}}_\ell \; \sqrt{[ (\boldsymbol{X}_\ell^\intercal
\boldsymbol{X}_\ell)^{-1}]_{ii}} \leq \f{\vep_i}{1 + \vep_i} |
\bs{B}_{i, \ell} |$} for $i = 1, \cd, m$, and {\small $\f{
\chi_{n_\ell - m, \; 1 - \ze \de_\ell }^2  } { (1 + \vep)^2 } \leq
n_\ell - m \leq \f{ \chi_{n_\ell - m, \; \ze \de_\ell }^2  } { (1 -
\vep)^2 }$} at some stage with index $\ell$, where $\de_\ell = \de$
for $1 \leq \ell \leq \tau$ and $\de_\ell = \de 2^{\tau - \ell}$ for
$\ell > \tau$.  Define $\wh{\bs{\si}} = \widehat{\bs{\si}}_{\bs{l}}$
and $\wh{\bs{\ba}} = \wh{\bs{\ba}}_{\bs{l}}$, where $\bs{l}$ is the
index of stage at which the observation of $y$ is stopped.  For $i =
1, \cd, m$, let $\wh{\bs{\ba}}_i$ be the $i$-th entry of
$\wh{\bs{\ba}}$. Then, $\Pr \{ \bs{l} < \iy \} = 1$ and $\Pr \{ |
\wh{\bs{\si}} - \sigma | \leq \vep \si, \;\; | \wh{\bs{\ba}}_i -
\beta_i | \leq \vep_i |\beta_i|\; \tx{for}
 \; i = 1, \cdots, m \} \geq 1 - \delta$
provided that $2 (m+1) (\tau + 1) \ze \leq 1$ and that $\inf_{\ell >
0} \f{n_{\ell + 1}}{n_\ell} > 1$.

\eeT

See Appendix \ref{App_Regression_Rev}  for a proof.

\section{Multistage Estimation of Quantile}

The estimation of a quantile of a random variable is a fundamental
problem of practical importance. Specially, in control engineering,
the performance of an uncertain dynamic system can be modeled as a
random variable. Hence, it is desirable to estimate the minimum
level of performance such that the probability of achieving it is
greater than a certain percentage.  In general, the problem of
estimating a quantile can be formulated as follows.

Let $X$ be a random variable with cumulative distribution function
$F_X (.)$. Define quantile $\xi_p = \inf \{ x: F_X (x) > p  \}$ for
$p \in (0, 1)$.  The objective is to estimate $\xi_p$ with
prescribed precision and confidence level based on i.i.d. samples
$X_1, X_2, \cd$ of $X$. To make it possible for the rigorous control
of estimation error and uncertainty of inference, we shall propose
multistage procedures. For this purpose, we need to define some
variables.  For an integer $n$, let $X_{i:n}$ denote the $i$-th
order statistics of i.i.d samples $X_1, \cd, X_n$ of $X$ such that
$- \iy = X_{0:n} < X_{1:n} \leq X_{2:n} \leq \cd \leq X_{n:n} < X_{n
+1:n} = \iy$. Let the sample sizes be a sequence of positive
integers $n_\ell, \; \ell = 1, 2, \cd$ such that $n_1 < n_2 < n_3 <
\cd$. At the $\ell$-th stage, the decision of termination or
continuation of sampling is made based on samples $X_1, \cd,
X_{n_\ell}$.

\subsection{Control of Absolute Error}

For estimating $\xi_p$ with a margin of absolute error $\vep > 0$,
our sampling procedure can be described as follows.

\beT \la{quantile_abs}   For $\ell = 1, 2, \cd$, define $\de_\ell =
\de$ for $1 \leq \ell \leq \tau$ and $\de_\ell = \de 2^{\tau -
\ell}$ for $\ell > \tau$, where $\tau$ is a positive integer.  Let
$i_\ell \leq n_\ell$ be the largest integer such that
$\sum_{k=0}^{i_\ell -1} { n_\ell \choose k} p^k (1 - p)^{n_\ell- k}
\leq \ze \delta_\ell$.  Let $j_\ell \geq 0$ be the smallest integer
such that $\sum_{k= j_\ell}^{n_\ell} { n_\ell \choose k} p^k (1 -
p)^{n_\ell- k} \leq \ze \delta_\ell$. Define $\wh{\bs{\xi}}_{p,
\ell}$ such that $\wh{\bs{\xi}}_{p, \ell} = X_{p n_\ell:n_\ell}$ if
$p n_\ell$ is an integer and $\wh{\bs{\xi}}_{p, \ell} = (\lc p
n_\ell \rc - p n_\ell ) X_{\lf p n_\ell \rf: n_\ell} + (p n_\ell -
\lf p n_\ell \rf ) X_{\lc p n_\ell \rc: n_\ell}$ otherwise. Suppose
that sampling is continued until $X_{j_\ell: n_\ell} - \vep \leq
\wh{\bs{\xi}}_{p, \ell} \leq X_{i_\ell: n_\ell} + \vep$ for some
stage with index $\ell$. Define estimator $\wh{\bs{\xi}}_{p} =
\wh{\bs{\xi}}_{p, \bs{l}}$ where $\bs{l}$ is the index of stage at
which the sampling is terminated. Then, $\Pr \{ \bs{l} < \iy \} = 1$
and $\Pr \{ | \wh{\bs{\xi}}_{p} - \xi_p | \leq \vep \} \geq 1 - \de$
provided that $2 (\tau + 1) \ze \leq 1$ and that $\inf_{\ell > 0}
\f{n_{\ell + 1}}{n_\ell} > 1$.

\eeT

See Appendix \ref{App_quantile_abs} for a proof.

\subsection{Control of Relative Error}

For estimating $\xi_p \neq 0$ with a margin of relative error $\vep
\in (0, 1)$, our sampling procedure can be described as follows.

\beT \la{quantile_rev}  For $\ell = 1, 2, \cd$, define $\de_\ell =
\de$ for $1 \leq \ell \leq \tau$ and $\de_\ell = \de 2^{\tau -
\ell}$ for $\ell > \tau$, where $\tau$ is a positive integer.  Let
$i_\ell \leq n_\ell$ be the largest integer such that
$\sum_{k=0}^{i_\ell -1} { n_\ell \choose k} p^k (1 - p)^{n_\ell- k}
\leq \ze \delta_\ell$.  Let $j_\ell \geq 0$ be the smallest integer
such that $\sum_{k= j_\ell}^{n_\ell} { n_\ell \choose k} p^k (1 -
p)^{n_\ell- k} \leq \ze \delta_\ell$. Define $\wh{\bs{\xi}}_{p,
\ell}$ such that $\wh{\bs{\xi}}_{p, \ell} = X_{p n_\ell:n_\ell}$ if
$p n_\ell$ is an integer and $\wh{\bs{\xi}}_{p, \ell} = (\lc p
n_\ell \rc - p n_\ell ) X_{\lf p n_\ell \rf: n_\ell} + (p n_\ell -
\lf p n_\ell \rf ) X_{\lc p n_\ell \rc: n_\ell}$ otherwise. Suppose
that sampling is continued until $[ 1 - \mrm{sgn} (
\wh{\bs{\xi}}_{p, \ell} ) \vep] X_{j_\ell: n_\ell}  \leq
\wh{\bs{\xi}}_{p, \ell} \leq  [ 1 + \mrm{sgn} ( \wh{\bs{\xi}}_{p,
\ell} ) \vep] X_{i_\ell: n_\ell}$ for some stage with index $\ell$.
Define estimator $\wh{\bs{\xi}}_{p} = \wh{\bs{\xi}}_{p, \bs{l}}$
where $\bs{l}$ is the index of stage at which the sampling is
terminated. Then, $\Pr \{ \bs{l} < \iy \} = 1$ and $\Pr \{ |
\wh{\bs{\xi}}_{p} - \xi_p | \leq \vep | \xi_p | \} \geq 1 - \de$
provided that $2 (\tau + 1) \ze \leq 1$ and that $\inf_{\ell > 0}
\f{n_{\ell + 1}}{n_\ell} > 1$.

\eeT

See Appendix \ref{App_quantile_rev} for a proof.

\subsection{Control of Absolute and Relative Errors}

For estimating $\xi_p$ with margin of absolute error $\vep_a > 0$
and margin of relative error $\vep_r \in (0, 1)$, our sampling
procedure can be described as follows.

\beT \la{quantile_mix}  For $\ell = 1, 2, \cd$, define $\de_\ell =
\de$ for $1 \leq \ell \leq \tau$ and $\de_\ell = \de 2^{\tau -
\ell}$ for $\ell > \tau$, where $\tau$ is a positive integer.  Let
$i_\ell \leq n_\ell$ be the largest integer such that
$\sum_{k=0}^{i_\ell -1} { n_\ell \choose k} p^k (1 - p)^{n_\ell- k}
\leq \ze \delta_\ell$. Let $j_\ell \geq 0$ be the smallest integer
such that $\sum_{k= j_\ell}^{n_\ell} { n_\ell \choose k} p^k (1 -
p)^{n_\ell- k} \leq \ze \delta_\ell$. Define $\wh{\bs{\xi}}_{p,
\ell}$ such that $\wh{\bs{\xi}}_{p, \ell} = X_{p n_\ell:n_\ell}$ if
$p n_\ell$ is an integer and $\wh{\bs{\xi}}_{p, \ell} = (\lc p
n_\ell \rc - p n_\ell ) X_{\lf p n_\ell \rf: n_\ell} + (p n_\ell -
\lf p n_\ell \rf ) X_{\lc p n_\ell \rc: n_\ell}$ otherwise.  Suppose
that sampling is continued until $X_{j_\ell: n_\ell} - \max (
\vep_a, \; \mrm{sgn} ( \wh{\bs{\xi}}_{p, \ell} ) \vep_r X_{j_\ell:
n_\ell} ) \leq \wh{\bs{\xi}}_{p, \ell} \leq X_{i_\ell: n_\ell} +
\max ( \vep_a, \; \mrm{sgn} ( \wh{\bs{\xi}}_{p, \ell} ) \vep_r
X_{i_\ell: n_\ell})$ for some stage with index $\ell$. Define
estimator $\wh{\bs{\xi}}_{p} = \wh{\bs{\xi}}_{p, \bs{l}}$ where
$\bs{l}$ is the index of stage at which the sampling is terminated.
Then, $\Pr \{ \bs{l} < \iy \} = 1$ and $\Pr \{ | \wh{\bs{\xi}}_{p} -
\xi_p | \leq \vep_a \; \tx{or} \; | \wh{\bs{\xi}}_{p} - \xi_p | \leq
\vep_r | \xi_p | \} \geq 1 - \de$ provided that $2 (\tau + 1) \ze
\leq 1$ and that $\inf_{\ell > 0} \f{n_{\ell + 1}}{n_\ell} > 1$.

\eeT

See Appendix \ref{App_quantile_mix} for a proof.

\sect{Conclusion}

In this paper, we have proposed a new framework of multistage estimation. Our general approach provides exact solutions for a wide spectrum of
estimation problems. Specific sampling schemes have been developed for parameters of common distributions. It is demonstrated that our new
methods are unprecedentedly efficient in terms of sampling cost, while rigorously guaranteeing prescribed level of confidence.

\appendix

\sect{Preliminary Results}

We need some preliminary results.

\subsection{Proof of Identity (\ref{RI_Indentity}) }
\la{App_RI_Indentity}

Since the identity is clearly true for the case that $[\udl{\se},
\ovl{\se}]$ does not contain $0$, we only consider the case that $0
\in [\udl{\se}, \ovl{\se}]$.  We claim that \be \la{ineq_RI_special}
\li \{ | \wh{ \bs{\se} } - \se | < \vep_r | \se |  \ri \} \subseteq
\li \{ \f{ \wh{ \bs{\se} } }{ 1 + \mrm{sgn} ( \wh{ \bs{\se} } )
\vep_r } < \se < \f{ \wh{ \bs{\se} } }{ 1 - \mrm{sgn} ( \wh{
\bs{\se} } ) \vep_r } \ri \}. \ee Let $\om \in \{ | \wh{ \bs{\se} }
- \se | < \vep_r | \se | \}$ and $\wh{ \se } = \wh{ \bs{\se} }
(\om)$.  Then, $| \wh{ \se } - \se | < \vep_r | \se |$. To show
(\ref{ineq_RI_special}), it suffices to show {\small $\f{ \wh{ \se }
}{ 1 + \mrm{sgn} ( \wh{ \se } ) \vep_r } < \se < \f{ \wh{ \se }  }{
1 - \mrm{sgn} ( \wh{ \se } ) \vep_r }$}.

In the case of $\se \geq 0$, we have $\wh{ \se } > (\se - \vep_r
|\se| ) \geq 0$ as a result of $| \wh{ \se } - \se | < \vep_r | \se
|$.  Moreover, {\small $\f{ \wh{ \se }  }{ 1 + \mrm{sgn} ( \wh{ \se
} ) \vep_r } = \f{ \wh{ \se }  }{ 1 + \vep_r } < \se < \f{ \wh{ \se
}  }{ 1 - \vep_r } = \f{ \wh{ \se }  }{ 1 - \mrm{sgn} ( \wh{ \se } )
\vep_r }$}.  In the case of $\se < 0$, we have $\wh{ \se } < (\se +
\vep_r |\se| ) < 0$ as a result of $| \wh{ \se } - \se | < \vep_r |
\se |$. Moreover, {\small $\f{ \wh{ \se }  }{ 1 + \mrm{sgn} ( \wh{
\se } ) \vep_r } = \f{ \wh{ \se }  }{ 1 - \vep_r } < \se < \f{ \wh{
\se } }{ 1 + \vep_r } = \f{ \wh{ \se }  }{ 1 - \mrm{sgn} ( \wh{ \se
} ) \vep_r }$}.  Therefore, we have established
(\ref{ineq_RI_special}).

In view of (\ref{ineq_RI_special}), it is obvious that $\{ |
\wh{\bs{\se}} - \se | < \vep_a \; \tx{or} \; | \wh{\bs{\se}} - \se |
< \vep_r |\se| \}  \subseteq \{ \mscr{L} (\wh{\bs{\se}}, \mbf{n})  <
\se < \mscr{U} (\wh{\bs{\se}}, \mbf{n}) \}$.  To complete the proof
of identity (\ref{RI_Indentity}), it remains to show $\{ |
\wh{\bs{\se}} - \se | < \vep_a \; \tx{or} \; | \wh{\bs{\se}} - \se |
< \vep_r |\se| \}  \supseteq \{ \mscr{L} (\wh{\bs{\se}}, \mbf{n})  <
\se < \mscr{U} (\wh{\bs{\se}}, \mbf{n}) \}$.  For this purpose, let
$\om \in \{ \mscr{L} (\wh{\bs{\se}}, \mbf{n})  < \se < \mscr{U}
(\wh{\bs{\se}}, \mbf{n}) \}$ and $\wh{ \se } = \wh{ \bs{\se} }
(\om)$.  Then, \be \la{contra} \min \li \{ \wh{ \se } - \vep_a, \;
\f{ \wh{ \se } }{ 1 + \mrm{sgn} ( \wh{ \se } ) \vep_r } \ri \} < \se
< \max \li \{ \wh{ \se } + \vep_a, \; \f{ \wh{ \se }  }{ 1 -
\mrm{sgn} ( \wh{ \se } ) \vep_r } \ri \} \ee Suppose, to get a
contradiction, that $| \wh{\se} - \se | \geq \vep_a$ and $| \wh{\se}
- \se | \geq \vep_r |\se|$.  There are $8$ cases:

(i) $\se \geq 0, \; \wh{\se} \geq \se + \vep_a, \;  \wh{\se} \geq
\se + \vep_r | \se |$.   In this case, we have $\wh{\se} \geq 0, \;
\se \leq \wh{\se} - \vep_a$ and $\se \leq  \f{ \wh{ \se } }{ 1 +
\vep_r } = \f{ \wh{ \se } }{ 1 + \mrm{sgn} ( \wh{ \se } ) \vep_r }$,
which contradicts the first inequality of (\ref{contra}).

(ii) $\se \geq 0, \; \wh{\se} \leq \se - \vep_a, \;  \wh{\se} \geq
\se + \vep_r | \se |$.  In this case, we have  $ \se + \vep_r | \se
| \leq \wh{\se} \leq \se - \vep_a$, which implies that $\vep_a = 0$
and $\wh{\se} \geq 0$.  Therefore, the first inequality of
(\ref{contra}) can be written as $\f{ \wh{ \se } }{ 1 + \vep_r } <
\se$, which contradicts to $\wh{\se} \geq \se + \vep_r | \se | = (1
+ \vep_r) \se$.

(iii) $\se \geq 0, \; \wh{\se} \geq \se + \vep_a, \;  \wh{\se} \leq
\se - \vep_r | \se |$.  In this case, we have $\se + \vep_a \leq
\wh{\se} \leq \se - \vep_r | \se |$, which implies that $\vep_a = 0$
and $\wh{\se} \geq 0$.  Therefore, the second inequality of
(\ref{contra}) can be written as $\f{ \wh{ \se } }{ 1 - \vep_r } >
\se$, which contradicts to $\wh{\se} \leq \se - \vep_r | \se | = (1
- \vep_r) \se$.

(iv) $\se \geq 0, \; \wh{\se} \leq \se - \vep_a, \;  \wh{\se} \leq
\se - \vep_r | \se |$.  In this case, we have $\se \geq \wh{\se} +
\vep_a$ and $\se \geq \f{ \wh{\se} } { 1 - \vep_r }$.  Hence, by the
second inequality of (\ref{contra}), we have $\f{ \wh{\se} } { 1 -
\vep_r } \leq \se < \f{ \wh{ \se }  }{ 1 - \mrm{sgn} ( \wh{ \se } )
\vep_r }$, which implies $\wh{\se} [ 1 - \mrm{sgn} ( \wh{ \se } )
\vep_r ] < \wh{ \se } (1 - \vep_r)$, i.e., $\vep_r | \wh{ \se } | >
\vep_r \wh{ \se }$.  It follows that $\wh{ \se } < 0$ and thus $\se
< 0$, which contradicts to $\se \geq 0$.

(v) $\se < 0, \; \wh{\se} \geq \se + \vep_a, \;  \wh{\se} \geq \se +
\vep_r | \se |$.  In this case, we have $\se \leq \wh{\se} - \vep_a$
and $\se \leq \f{ \wh{\se} } { 1 - \vep_r }$.  Hence, by the first
inequality of (\ref{contra}), we have $\f{ \wh{\se} } { 1 - \vep_r }
\geq \se > \f{ \wh{ \se }  }{ 1 + \mrm{sgn} ( \wh{ \se } ) \vep_r
}$, which implies $\wh{\se} [ 1 + \mrm{sgn} ( \wh{ \se } ) \vep_r ]
> \wh{ \se } (1 - \vep_r)$, i.e., $\vep_r | \wh{ \se } | > - \vep_r
\wh{ \se }$.  It follows that $\wh{ \se } > 0$ and thus $\se > 0$,
which contradicts to $\se < 0$.

(vi) $\se < 0, \; \wh{\se} \leq \se - \vep_a, \;  \wh{\se} \geq \se
+ \vep_r | \se |$.  In this case, we have $\se - \vep_a \geq
\wh{\se} \geq \se + \vep_r | \se |$, which implies that $\vep_a = 0$
and $\wh{\se} < 0$.  Therefore, the first inequality of
(\ref{contra}) can be written as $\f{ \wh{ \se } }{ 1 - \vep_r } <
\se$, which contradicts to $\wh{\se} \geq \se + \vep_r | \se | = (1
- \vep_r) \se$.

(vii) $\se < 0, \; \wh{\se} \geq \se + \vep_a, \;  \wh{\se} \leq \se
- \vep_r | \se |$.   In this case, we have  $ \se - \vep_r | \se |
\geq \wh{\se} \geq \se + \vep_a$, which implies that $\vep_a = 0$
and $\wh{\se} < 0$.  Therefore, the second inequality of
(\ref{contra}) can be written as $\f{ \wh{ \se } }{ 1 + \vep_r } >
\se$, which contradicts to $\wh{\se} \leq \se - \vep_r | \se | = (1
+ \vep_r) \se$.

(viii) $\se < 0, \; \wh{\se} \leq \se - \vep_a, \;  \wh{\se} \leq
\se - \vep_r | \se |$.  In this case, we have $\wh{\se} < 0, \; \se
\geq \wh{\se} + \vep_a$ and $\se \geq  \f{ \wh{ \se } }{ 1 + \vep_r
} = \f{ \wh{ \se } }{ 1 - \mrm{sgn} ( \wh{ \se } ) \vep_r }$, which
contradicts the second inequality of (\ref{contra}).

From the above $8$ cases, we see that the assumption that $|
\wh{\se} - \se | \geq \vep_a$ and $| \wh{\se} - \se | \geq \vep_r
|\se|$ always leads to a contradiction.  Therefore, it must be true
that either $| \wh{\se} - \se | < \vep_a$ or $| \wh{\se} - \se | <
\vep_r |\se|$. This proves $\{ | \wh{\bs{\se}} - \se | < \vep_a \;
\tx{or} \; | \wh{\bs{\se}} - \se | < \vep_r |\se| \}  \supseteq \{
\mscr{L} (\wh{\bs{\se}})  < \se < \mscr{U} (\wh{\bs{\se}}) \}$ and
consequently completes the proof of identity (\ref{RI_Indentity}).

\subsection{Probability Transform Inequalities}

The well-known probability transform theorem asserts that $\Pr \{ F_Z(Z) \leq \al \} = \Pr \{ G_Z(Z) \leq \al \} = \al$ for any continuous
random variable $Z$ and positive number $\al \in [0, 1]$. In the general case that $Z$ is not necessarily continuous, the probability transform
equalities may not be true. Fortunately, their generalizations, referred to as ``probability transform inequalities'' in this paper,  has been
established in the literature as follows.

 \beL \la{ProbTrans} $\Pr \{ F_Z(Z) \leq \al \}
\leq \al$ and $\Pr \{ G_Z(Z) \leq \al \} \leq \al$ for any random
variable $Z$ and positive number $\al$. \eeL

Lemma \ref{ProbTrans} can be shown as follows.  Let $I_Z$ denote the support of $Z$.  If $\{ z \in I_Z: F_Z(z) \leq \al \}$ is empty, then, $\{
F_Z(Z) \leq \al \}$ is an impossible event and thus $\Pr \{ F_Z(Z) \leq \al \} = 0$.  Otherwise, we can define $z^\star = \max \{ z \in I_Z:
F_Z(z) \leq \al \}$. It follows from the definition of $z^\star$ that $F_Z (z^\star ) \leq \al$. Since $F_Z(z)$ is non-decreasing with respect
to $z$, we have $\{ F_Z(Z) \leq \al \} = \{ Z \leq z^\star \}$. Therefore, $\Pr \{ F_Z(Z) \leq \al  \}  = \Pr \{ Z \leq z^\star \} = F_Z
(z^\star ) \leq \al $ for any $\al > 0$. By a similar method, one can show $\Pr \{ G_Z(Z) \leq \al \} \leq \al$ for any $\al > 0$.

\subsection{Property of ULE} \la{App_Basic_ULE}

\beL

\la{ULE_Basic}

Let $\mbf{m}$ be a stopping time such that for any positive integer $m$, event $\{\mbf{m} = m \}$ depends only on $X_1, \cd, X_m$. Let
$\mscr{E}$ be an event dependent only on random tuple $(X_1, \cd, X_{\mbf{m}})$. Let $\varphi(X_1, \cd, X_{\mbf{m}})$ be a ULE of $\se$. Then,

(i) $\Pr \{ \mscr{E} \mid \se \}$ is non-increasing with respect to
$\se \in \Se$ no less than $z$ provided that $\mscr{E} \subseteq \{
\varphi(X_1, \cd, X_{\mbf{m}}) \leq z \}$.

(ii) $\Pr \{ \mscr{E} \mid \se \}$ is non-decreasing with respect to
$\se \in \Se$ no greater than $z$ provided that $\mscr{E} \subseteq
\{ \varphi(X_1, \cd, X_{\mbf{m}}) \geq z \}$.

\eeL

\bpf

We first consider the case that $X_1, X_2, \cd$ are discrete random
variables. Let $I_{\mbf{m}}$ denote the support of $\mbf{m}$, i.e.,
$I_{\mbf{m}} = \{ \mbf{m} (\om): \om \in \Om \}$. Define $\mscr{X}_m
= \{ (X_1 (\om), \cd, X_m(\om) ): \om \in \mscr{E}, \; \mbf{m} (\om)
= m \}$ for $m \in I_{\mbf{m}}$. Then, \be \la{eeewhole}
 \Pr \{ \mscr{E} \mid \se \}  = \sum_{m \in I_{\mbf{m}}} \;
 \sum_{ (x_1, \cd, x_m) \in \mscr{X}_m } \Pr \{ X_i = x_i, \; i =
1, \cd, m \mid \se \}. \ee

To show statement (i), using the assumption that $\mscr{E} \subseteq
\{ \varphi(X_1, \cd, X_{\mbf{m}}) \leq z \}$, we have $\varphi(x_1,
\cd, x_m) \leq z$ for $(x_1, \cd, x_m) \in \mscr{X}_m$ with $m \in
I_{\mbf{m}}$.  Since $\varphi(X_1, \cd, X_{\mbf{m}})$ is a ULE of
$\se$, we have that $\Pr \{ X_i = x_i, \; i = 1, \cd, m \mid \se \}$
is non-increasing with respect to $\se \in \Se$ no less than $z$. It
follows immediately from (\ref{eeewhole}) that statement (i) is
true.

To show statement (ii), using the assumption that $\mscr{E}
\subseteq \{ \varphi(X_1, \cd, X_{\mbf{m}}) \geq z \}$, we have
$\varphi(x_1, \cd, x_m) \geq z$ for $(x_1, \cd, x_m) \in \mscr{X}_m$
with $m \in I_{\mbf{m}}$.  Since $\varphi(X_1, \cd, X_{\mbf{m}})$ is
a ULE of $\se$, we have that $\Pr \{ X_i = x_i, \; i = 1, \cd, m
\mid \se \}$ is non-decreasing with respect to $\se \in \Se$ no
greater than $z$. It follows immediately from (\ref{eeewhole}) that
statement (ii) is true.

For the case that $X_1, X_2, \cd$ are continuous random variables,
we can also show the lemma by modifying the argument for the
discrete case. Specially, the summation of likelihood function $\Pr
\{ X_i = x_i, \; i = 1, \cd, m \mid \se \}$ over the set of tuple
$(x_1, \cd, x_m )$ is replaced by the integration of the joint
probability density function $f_{X_1, \cd, X_m} (x_1, \cd, x_m,
\se)$ over the set of $(x_1, \cd, x_m)$. This concludes the proof of
Lemma \ref{ULE_Basic}.

\epf

\sect{Proof of Theorem \ref{Monotone_second} }
\la{App_Monotone_second}

Making use of assumptions (ii)-(iii), the definition of the sampling
scheme and the monotonicity of $F_{\wh{\bs{\se}}_\ell}(z, \se)$ as
asserted by Lemma \ref{ULE_Basic}, we have \bee \Pr \{ \se \geq
\mscr{U} (\wh{\bs{\se}}, \mbf{n}) \mid \se \} & = & \sum_{\ell =
1}^s \Pr \{ \se \geq \mscr{U}
(\wh{\bs{\se}}_\ell, \mbf{n}_\ell), \; \bs{l} = \ell \mid \se \}\\
& \leq & \sum_{\ell = 1}^s \Pr \{ \se \geq \mscr{U}
(\wh{\bs{\se}}_\ell, \mbf{n}_\ell), \; \bs{D}_\ell = 1 \mid \se \}\\
& \leq &  \sum_{\ell = 1}^s \Pr \li \{ \se \geq \mscr{U}
(\wh{\bs{\se}}_\ell, \mbf{n}_\ell) \geq \wh{\bs{\se}}_\ell, \;
F_{\wh{\bs{\se}}_\ell} (\wh{\bs{\se}}_\ell,
\mscr{U} (\wh{\bs{\se}}_\ell, \mbf{n}_\ell)) \leq \ze \de_\ell \mid \se \ri \}\\
& \leq &  \sum_{\ell = 1}^s \Pr \li \{ F_{\wh{\bs{\se}}_\ell}
(\wh{\bs{\se}}_\ell, \se )  \leq \ze \de_\ell \mid \se \ri \} \leq
\ze \sum_{\ell = 1}^s \de_\ell \eee for any $\se \in \Se$,  where
the last inequality follows from Lemma \ref{ProbTrans}.

Similarly, we can show that $\Pr \{ \se \leq \mscr{L}
(\wh{\bs{\se}}, \mbf{n}) \mid \se \} \leq \sum_{\ell = 1}^s \Pr \{
\se \leq \mscr{L} (\wh{\bs{\se}}_\ell, \mbf{n}_\ell), \; \bs{D}_\ell
= 1 \mid \se \} \leq \ze \sum_{\ell = 1}^s \de_\ell$.  Hence, $\Pr
\{ \mscr{L} (\wh{\bs{\se}}, \mbf{n})  < \se < \mscr{U}
(\wh{\bs{\se}}, \mbf{n}) \mid \se \} \geq 1 - \Pr \{ \se \leq
\mscr{L} (\wh{\bs{\se}}, \mbf{n}) \mid \se \} - \Pr \{ \se \geq
\mscr{U} (\wh{\bs{\se}}, \mbf{n}) \mid \se \} \geq 1 - 2  \ze
\sum_{\ell = 1}^s \de_\ell$. This concludes the proof of Theorem
\ref{Monotone_second}.

\sect{Proof of Theorem \ref{General Inclusion Principle} } \la{General Inclusion Principle_App}

Making use of the assumption that $\Pr \{ \bs{\tau} < \iy  \} = 1$,  we have \bee
 &   & \Pr \{ \bs{\se} \notin
\bs{\mcal{B}}_{\bs{\tau}}  \} = \Pr \li \{ \bigcup_{\ell \in I_{\bs{\tau}}} \{ \bs{\se} \notin \bs{\mcal{B}}_{\bs{\tau}}, \; \bs{\tau} = \ell \}
\ri \}. \eee By the assumptions on the random regions and stopping time, we have
 {\small \bee
&  & \bigcup_{\ell \in I_{\bs{\tau}} } \{ \bs{\se} \notin \bs{\mcal{B}}_{\bs{\tau}}, \; \bs{\tau} = \ell \} = \bigcup_{\ell \in I_{\bs{\tau}} }
\{ \bs{\se} \notin \bs{\mcal{B}}_\ell, \; \bs{\tau} = \ell \} \subseteq  \bigcup_{\ell \in I_{\bs{\tau}} } \{ \bs{\se} \notin
\bs{\mcal{B}}_\ell, \; \bs{\mcal{A}}_\ell
\subseteq \bs{\mcal{B}}_\ell \} \\
&  &  \subseteq \bigcup_{\ell \in I_{\bs{\tau}} } \{ \bs{\se} \notin \bs{\mcal{A}}_\ell, \; \bs{\mcal{A}}_\ell  \subseteq \bs{\mcal{B}}_\ell \}
\subseteq \bigcup_{\ell \in I_{\bs{\tau}} } \{ \bs{\se} \notin \bs{\mcal{A}}_\ell \} = \{ \bs{\se} \notin \bs{\mcal{A}}_\ell  \; \tx{for some
index $\ell \in I_{\bs{\tau}}$} \}.
 \eee}
Therefore,
 \bee
&  & \Pr \{ \bs{\se} \in \bs{\mcal{B}}_{\bs{\tau}}  \}   = 1 - \Pr \{ \bs{\se} \notin
\bs{\mcal{B}}_{\bs{\tau}}  \} = 1 -  \Pr \li \{ \bigcup_{\ell \in I_{\bs{\tau}} } \{ \bs{\se} \notin \bs{\mcal{B}}_{\bs{\tau}}, \; \bs{\tau} = \ell \} \ri \}\\
&  & \geq 1 - \Pr \{ \bs{\se} \notin \bs{\mcal{A}}_\ell  \; \tx{for some index $\ell \in I_{\bs{\tau}}$} \}\\
&  & = \Pr \{ \bs{\se} \in \bs{\mcal{A}}_\ell  \; \tx{ for $\ell \in I_{\bs{\tau}} $} \mid \se \} \geq 1 - \sum_{\ell \in I_{\bs{\tau}} } \Pr \{
\bs{\se} \notin \bs{\mcal{A}}_\ell \}, \eee  where the last inequality follows from Bonferroni's inequality and the continuity of the
probability measure. This completes the proof of
 the theorem.

\sect{Proof of Theorem \ref{explicit_Chernoff} }
\la{explicit_Chernoff_app}

By the independence of samples, the likelihood function can be
written as {\small $\prod_{i = 1}^n f_X (x_i, \se) = \li [ \prod_{i
= 1}^n c (x_i) \ri ] \times \exp \li ( \eta (\se) \sum_{i = 1}^n x_i
- n \psi (\se) \ri )$}. By the assumption that $\f{d \eta (\se)}{d
\se}
> 0$ and $\f{ \psi^\prime ( \se ) }{ \eta^\prime ( \se ) } = \se$
for $\se \in \Se$, we have that
\[
\f{ d \exp \li ( \eta (\se) z - \psi (\se)  \ri )  } { d \se } =  (z
- \se) \exp \li ( \eta (\se) z - \psi (\se)  \ri ) \f{d \eta
(\se)}{d \se},
\]
which is positive for $\se < z$ and negative for $\se > z$.  This
implies that $\exp \li ( \eta (\se) z - \psi (\se)  \ri )$ is
monotonically increasing with respect to $\se$ less than $z$ and
monotonically decreasing with respect to $\se$ greater than $z$.
This proves that $\ovl{X}_n$ is a ULE of $\se$.  It remains to show
the probabilistic inequalities regarding $\ovl{X}_n$.

Let $\vse (.)$ be the inverse function of $\eta(.)$ such that $\eta
(\vse (\ze) ) = \ze$ for $\ze \in \{ \eta (\se): \se \in \Se \}$.
Define compound function $\phi(.)$ such that $\phi (\zeta) = \psi
(\vse (\zeta))$ for $\ze \in \{ \eta (\se): \se \in \Se \}$. For
simplicity of notations, we abbreviate $\vse (\zeta)$ as $\vse$ when
this can be done without causing confusion.   Using the definition
that $\eta (\vse (\ze) ) = \ze$, the assumption that $\f{ d \psi (
\se ) }{ d \se } = \se \f{ d \eta ( \se ) }{d \se}$, and the chain
rule of differentiation, we have \be \la{use8} \f{ d \phi (\zeta )
}{d \zeta} = \f{ d \psi (\vse) }{d \vse} \f{ d \vse }{ d \zeta} =
\vse \f{ d \eta (\vse) }{d \vse} \f{ d \vse }{ d \zeta} = \vse  \f{
d \eta (\vse) }{d \ze} = \vse \f{ d \ze }{d \ze} = \vse (\ze). \ee
Putting $\zeta = \eta (\se)$, we have {\small $\bb{E} \li [ \exp \li
( t \sum_{i = 1}^n  X_i \ri )  \ri ] =   \exp \li ( n \phi (\zeta +
t) - n \phi (\zeta) \ri )$}.   By virtue of (\ref{use8}), the
derivative of $n \phi (\zeta + t) - n \phi (\zeta)$ with respect to
$t$ is $n \f{d \phi (\zeta + t) }{d t} = n \vse ( \zeta + t)$,
which is equal to $n \vse ( \zeta ) = n \se$ for $t = 0$.  Thus,
$\bb{E} [ \ovl{X}_n ] = \se$, which implies that $\ovl{X}_n$ is also
an unbiased estimator of $\se$.

Again by virtue of (\ref{use8}), the derivative of $- t n z + n \phi
(\zeta + t) - n \phi (\zeta)$ with respect to $t$ is
\[
-  n z +  n \f{d \phi (\zeta + t) }{d t} = - n z + n \vse ( \zeta +
t),
\]
which is equal to $0$ for $t$ such that $\vse ( \zeta + t ) = z$ or
equivalently, $\zeta + t = \eta (z)$, which implies $t = \eta (z) -
\eta (\se )$.  Since $\bb{E} \li [ \exp \li ( n t ( \ovl{X}_n - z)
\ri ) \ri ]$ is a convex function of $t$, its infimum with respect
to $t \in \bb{R}$ is attained at $t = \eta (z) - \eta (\se )$. It
follows that \bee & & \inf_{t \in \bb{R} } \bb{E} \li [ \exp \li ( n
t ( \ovl{X}_n - z) \ri ) \ri ] = \inf_{t \in
\bb{R}} \exp \li ( - t n z +  n \phi (\zeta + t) - n \phi (\zeta) \ri )\\
&  =  & \exp \li ( - [ \eta (z) - \eta (\se ) ]  n z +  n \phi (
\eta
(z) ) - n \phi (\zeta) \ri )\\
&  =  & \exp \li ( - [ \eta (z) - \eta (\se ) ]  n z +  n \psi ( z)
- n
\psi (\se) \ri )\\
& = &  \li [ \f{ \exp \li ( \eta (\se )   z -  \psi (\se) \ri ) } {
\exp \li ( \eta (z ) z - \psi (z) \ri ) } \ri ]^n = [ w (z, \se )
]^n. \eee Finally, the probabilistic inequalities regarding
$\ovl{X}_n$ can be established by virtue of the above results, the
Chernoff bound and the assumption that $\eta (\se)$ is increasing
with respect to $\se \in \Se$.

\sect{Proof of Theorem \ref{fullcontrol} } \la{fullcontrol_app}

In the case of $\sum_{\ell = 1}^\iy \de_\ell < \iy$, by Theorem
\ref{Monotone_second} and Corollary \ref{Monotone_third}, we have
$\Pr \{ \mscr{L} (\wh{\bs{\se}}) < \se < \mscr{U} (\wh{\bs{\se}})
\mid \se \} \geq 1 - 2 \ze \sum_{\ell = 1}^\iy \de_\ell \to 1$ as
$\ze \to 0$.   It remains to show the theorem for the case that $0 <
\al < \de_\ell < \ba$ for all $\ell$.  By Theorem
\ref{explicit_Chernoff}, we have $\mcal{F} (z, \se) = \mcal{G} (z,
\se) = w (z, \se)$ for $\se \in \Se$.  Using $\f{ \psi^\prime ( \se
) }{ \eta^\prime ( \se ) } = \se$ for $\se \in \Se$, we can show
that $\f{ \pa \ln w(z, \se)  } {\pa \se  } = \eta^\prime (\se ) (z -
\se)$ and $\f{ \pa \ln w(z, \se)  } {\pa z } = \eta (\se) - \eta
(z)$,  which implies that $\mcal{F} (z, \se)$ and $\mcal{G} (z,
\se)$ are less than $1$ for $\se \in \Se$ not equal to $z$.  For
$\se \in \Se$, let $\ep$ be a positive number small enough such that
$(\se - \ep, \se + \ep) \subseteq \Se$.  Let $f = \max_{z \in (\se -
\ep, \se + \ep)} \mcal{F} (z, \mscr{U} (z))$ and $g = \max_{z \in
(\se - \ep, \se + \ep)} \mcal{G} (z, \mscr{L} (z))$.  We claim that
$\max \{ f, g \} < 1$.  To show the claim, note that, if $\{ z \in
(\se - \ep, \se + \ep): \mscr{U} (z) < \sup \Se \}$ is an empty set,
then $f$ is equal to $0$; otherwise $f$ is smaller than $1$ as a
consequence of the assumption that $\mscr{U} (\se) > \se$ for all
$\se \in \Se$ and that $\mcal{F} (z, \se) < 1$ when $\se \in \Se$ is
not equal to $z$. Similarly, if $\{ z \in (\se - \ep, \se + \ep):
\mscr{L} (z)
> \inf \Se \}$ is an empty set, then $g$ is equal to $0$;
otherwise $g$ is smaller than $1$.  This proves the claim. By the
definitions of the sampling schemes, \bee \Pr \{ \mbf{n}
> n \} & \leq & \Pr \{ [ \mcal{F} ( \ovl{X}_n, \mscr{U} ( \ovl{X}_n
) ) ]^n
> \ze \de_\ell \; \tx{or} \;  [ \mcal{G} ( \ovl{X}_n, \mscr{L} (
\ovl{X}_n  ) ) ]^n > \ze \de_\ell \}\\
&  \leq  & \Pr \{  [ \mcal{F} ( \ovl{X}_n, \mscr{U} (  \ovl{X}_n  )
) ]^n > \ze \al \; \tx{or} \;  [ \mcal{G} ( \ovl{X}_n, \mscr{L} (
\ovl{X}_n  ) ) ]^n > \ze \al \}. \eee Hence, in the case of $\max \{
f, g \} = 0$, we have $\Pr \{ \mbf{n}
> n \} = 0$.   In the case of $0 < \max
\{ f, g \} < 1$, let {\small $n = \li \lc \f{ \ln (\ze \al)  } {\ln
\max \{ f, g \}  } \ri \rc + 1$}.  Then, $\{  [ \mcal{F} (
\ovl{X}_n, \mscr{U} (  \ovl{X}_n  ) ) ]^n > \ze \al, \; |\ovl{X}_n -
\se| < \ep \}$ and $\{  [ \mcal{G} ( \ovl{X}_n, \mscr{L} ( \ovl{X}_n
) ) ]^n > \ze \al, \; |\ovl{X}_n - \se| < \ep \}$ are impossible
events.  It follows that {\small \bee \Pr \{ \mbf{n} > n \} & \leq &
\Pr \{  [ \mcal{F} ( \ovl{X}_n, \mscr{U} (  \ovl{X}_n  ) )
]^n > \ze \al, \; |\ovl{X}_n - \se| < \ep \}\\
&  & \qu + \Pr \{  [ \mcal{G} ( \ovl{X}_n, \mscr{L} ( \ovl{X}_n ) )
]^n > \ze \al, \; |\ovl{X}_n - \se| < \ep \} + \Pr \{ |\ovl{X}_n - \se| \geq \ep \} \\
&  = & \Pr \{ |\ovl{X}_n - \se| \geq \ep \} \leq [\mcal{F} (\se +
\ep, \se)]^n +  [\mcal{G} (\se - \ep, \se)]^n  \eee} and thus
{\small \bee \Pr \{ \mscr{L} (\wh{\bs{\se}}) \geq \se \; \tx{or} \;
\mscr{U} (\wh{\bs{\se}}) \leq \se \mid \se \}
 & = & \Pr \{ \mscr{L} (\wh{\bs{\se}}) \geq \se \; \tx{or} \;
\mscr{U} (\wh{\bs{\se}}) \leq \se,  \; \mbf{n} \leq n \mid \se \} +
\Pr \{ \mbf{n } > n \}\\
& \leq & 2 n \ze \ba +  [\mcal{F} (\se + \ep, \se)]^n +  [\mcal{G}
(\se - \ep, \se)]^n  \to 0 \eee} as $\ze \to 0$.  This completes the
proof of the theorem.

\sect{Proof of Theorem \ref{Main_Bound_Gen} }
\la{App_Main_Bound_Gen}

Let $\se^\prime < \se^{\prime \prime}$ be two consecutive distinct
elements of $I_\mscr{L} \cup \{ a, b\}$.  Then, $\{ \se \leq
\mscr{L}( \wh{\bs{\se}}, \mbf{n} ) < \se^{\prime \prime}  \}
\subseteq \{ \se^\prime < \mscr{L}( \wh{\bs{\se}}, \mbf{n} ) <
\se^{\prime \prime}  \} = \emptyset$ and it follows that $\{
\mscr{L}( \wh{\bs{\se}}, \mbf{n} ) \geq \se  \} = \{ \mscr{L}(
\wh{\bs{\se}}, \mbf{n} ) \geq \se^{\prime \prime}  \} \cup \{ \se
\leq \mscr{L}( \wh{\bs{\se}}, \mbf{n} ) < \se^{\prime \prime}  \} =
\{ \mscr{L}( \wh{\bs{\se}}, \mbf{n} ) \geq \se^{\prime \prime}  \}$
for any $\se \in (\se^\prime, \se^{\prime \prime}]$.  Recalling that
$ \{ \wh{\bs{\se}} \geq \mscr{L} ( \wh{\bs{\se}}, \mbf{n} ) \}$ is a
sure event, we have $\{ \mscr{L}( \wh{\bs{\se}}, \mbf{n}  ) \geq
\se^{\prime \prime} \} = \{ \wh{\bs{\se}} \geq \se^{\prime \prime},
\; \mscr{L}( \wh{\bs{\se}}, \mbf{n}  ) \geq \se^{\prime \prime} \}$.
Invoking the second statement of Lemma \ref{ULE_Basic}, we have that
$\Pr \{ \se \leq \mscr{L}( \wh{\bs{\se}}, \mbf{n}  ) \; \tx{and} \;
\mscr{E} \; \tx{occurs} \mid \se \} = \Pr \{ \mscr{L}(
\wh{\bs{\se}}, \mbf{n}  ) \geq \se^{\prime \prime} \; \tx{and} \;
\mscr{E} \; \tx{occurs} \mid \se \} = \Pr \{ \wh{\bs{\se}} \geq
\se^{\prime \prime}, \; \mscr{L}( \wh{\bs{\se}}, \mbf{n}  ) \geq
\se^{\prime \prime} \; \tx{and} \; \mscr{E} \; \tx{occurs} \mid \se
\}$ is non-decreasing with respect to $\se \in (\se^\prime,
\se^{\prime \prime} ]$.  This implies that the maximum of $\Pr \{
\se \leq \mscr{L}( \wh{\bs{\se}}, \mbf{n}  ) \; \tx{and} \; \mscr{E}
\; \tx{occurs} \mid \se \}$ with respect to $\se \in ( \se^\prime,
\se^{\prime \prime} ]$ is equal to $\Pr \{ \wh{\bs{\se}} \geq
\se^{\prime \prime}, \; \mscr{L}( \wh{\bs{\se}}, \mbf{n}  ) \geq
\se^{\prime \prime} \; \tx{and} \; \mscr{E} \; \tx{occurs} \mid
\se^{\prime \prime} \}$. Since the argument holds for arbitrary
consecutive distinct elements of $I_\mscr{L} \cup \{ a, b\}$, we
have established statement (I) regarding $\Pr \{ \se \leq \mscr{L}(
\wh{\bs{\se}}, \mbf{n}  ) \; \tx{and} \; \mscr{E} \; \tx{occurs}
\mid \se \}$ for $\se \in [a, b]$.  To prove the statement regarding
$\Pr \{ \se < \mscr{L}( \wh{\bs{\se}}, \mbf{n} ) \; \tx{and} \;
\mscr{E} \; \tx{occurs} \mid \se \}$, note that $\{ \se < \mscr{L}(
\wh{\bs{\se}}, \mbf{n} ) < \se^{\prime \prime} \} \subseteq \{
\se^\prime < \mscr{L}( \wh{\bs{\se}}, \mbf{n} ) < \se^{\prime
\prime} \} = \emptyset$, which implies that $\{ \mscr{L}(
\wh{\bs{\se}}, \mbf{n} ) > \se \} = \{ \mscr{L}( \wh{\bs{\se}},
\mbf{n} ) \geq \se^{\prime \prime}  \} \cup \{ \se < \mscr{L}(
\wh{\bs{\se}}, \mbf{n} ) < \se^{\prime \prime} \} = \{ \mscr{L}(
\wh{\bs{\se}}, \mbf{n} ) \geq \se^{\prime \prime} \}$ for any $\se
\in [ \se^\prime, \se^{\prime \prime} )$. Hence, $\Pr \{ \se <
\mscr{L}( \wh{\bs{\se}}, \mbf{n} ) \; \tx{and} \; \mscr{E} \;
\tx{occurs} \mid \se \} = \Pr \{ \mscr{L}( \wh{\bs{\se}}, \mbf{n}  )
\geq \se^{\prime \prime} \; \tx{and} \; \mscr{E} \; \tx{occurs} \mid
\se \} = \Pr \{ \wh{\bs{\se}} \geq \se^{\prime \prime}, \; \mscr{L}(
\wh{\bs{\se}}, \mbf{n}  ) \geq \se^{\prime \prime} \; \tx{and} \;
\mscr{E} \; \tx{occurs} \mid \se \}$ is non-decreasing with respect
to $\se \in [ \se^\prime, \se^{\prime \prime} )$.  This implies that
the supremum of $\Pr \{ \se < \mscr{L}( \wh{\bs{\se}}, \mbf{n}  ) \;
\tx{and} \; \mscr{E} \; \tx{occurs} \mid \se \}$ with respect to
$\se \in [\se^\prime, \se^{\prime \prime} )$ is equal to $\Pr \{
\wh{\bs{\se}} \geq \se^{\prime \prime}, \; \mscr{L}( \wh{\bs{\se}},
\mbf{n}  ) \geq \se^{\prime \prime} \; \tx{and} \; \mscr{E} \;
\tx{occurs} \mid \se^{\prime \prime} \}$. Since the argument holds
for arbitrary consecutive distinct elements of $I_\mscr{L} \cup \{
a, b\}$, we have established statement (I) regarding $\Pr \{ \se <
\mscr{L}( \wh{\bs{\se}}, \mbf{n}  ) \; \tx{and} \; \mscr{E} \;
\tx{occurs} \mid \se \}$ for $\se \in [a, b]$.

To prove statement (II) regarding $\Pr \{ \se \geq \mscr{U}(
\wh{\bs{\se}}, \mbf{n} ) \; \tx{and} \; \mscr{E} \; \tx{occurs} \mid
\se \}$, let $\se^\prime < \se^{\prime \prime}$ be two consecutive
distinct elements of $I_\mscr{U} \cup \{ a, b\}$. Then, $\{
\se^\prime < \mscr{U}( \wh{\bs{\se}}, \mbf{n} ) \leq \se  \}
\subseteq \{ \se^\prime < \mscr{U}( \wh{\bs{\se}}, \mbf{n} ) <
\se^{\prime \prime}  \} = \emptyset$ and it follows that $\{
\mscr{U}( \wh{\bs{\se}}, \mbf{n} ) \leq \se  \} = \{ \mscr{U}(
\wh{\bs{\se}}, \mbf{n} ) \leq \se^\prime  \} \cup \{ \se^\prime <
\mscr{U}( \wh{\bs{\se}}, \mbf{n} ) \leq \se  \} = \{ \mscr{U}(
\wh{\bs{\se}}, \mbf{n} ) \leq \se^\prime  \}$ for any $\se \in
[\se^\prime, \se^{\prime \prime} )$.  Recalling that $ \{
\wh{\bs{\se}} \leq \mscr{U} ( \wh{\bs{\se}}, \mbf{n} ) \}$ is a sure
event, we have $\{ \mscr{U}( \wh{\bs{\se}}, \mbf{n} ) \leq
\se^\prime  \} = \{ \wh{\bs{\se}} \leq \se^\prime, \; \mscr{U}(
\wh{\bs{\se}}, \mbf{n} ) \leq \se^\prime \}$.  Consequently, $\Pr \{
\mscr{U}( \wh{\bs{\se}}, \mbf{n} ) \leq \se \; \tx{and} \; \mscr{E}
\; \tx{occurs}  \mid \se \} = \Pr \{ \mscr{U}( \wh{\bs{\se}},
\mbf{n} ) \leq \se^\prime \; \tx{and} \; \mscr{E} \; \tx{occurs}
\mid \se \} = \Pr \{ \wh{\bs{\se}} \leq \se^\prime, \; \mscr{U}(
\wh{\bs{\se}}, \mbf{n} ) \leq \se^\prime \; \tx{and} \; \mscr{E} \;
\tx{occurs} \mid \se \}$ is  non-increasing with respect to $\se \in
[\se^\prime, \se^{\prime \prime} )$ as a result of the first
statement of Lemma \ref{ULE_Basic}.   This implies that the maximum
of $\Pr \{ \mscr{U}( \wh{\bs{\se}}, \mbf{n} ) \leq \se \; \tx{and}
\; \mscr{E} \; \tx{occurs} \mid \se \}$  for $\se \in [ \se^\prime,
\se^{\prime \prime} )$ is equal to $\Pr \{ \wh{\bs{\se}} \leq
\se^\prime, \; \mscr{U}( \wh{\bs{\se}}, \mbf{n} ) \leq \se^\prime \;
\tx{and} \; \mscr{E} \; \tx{occurs} \mid \se^\prime \}$.  Since the
argument holds for arbitrary consecutive distinct elements of
$I_\mscr{U} \cup \{ a, b\}$, we have established statement (II)
regarding $\Pr \{ \se \geq \mscr{U}( \wh{\bs{\se}}, \mbf{n}  ) \;
\tx{and} \; \mscr{E} \; \tx{occurs} \mid \se \}$ for $\se \in [a,
b]$.  To prove the statement regarding $\Pr \{ \se > \mscr{U}(
\wh{\bs{\se}}, \mbf{n}  ) \; \tx{and} \; \mscr{E} \; \tx{occurs}
\mid \se \}$, note that $\{ \se^\prime < \mscr{U}( \wh{\bs{\se}},
\mbf{n} ) < \se  \} \subseteq \{ \se^\prime < \mscr{U}(
\wh{\bs{\se}}, \mbf{n} ) < \se^{\prime \prime}  \} = \emptyset$,
which implies that $\{ \mscr{U}( \wh{\bs{\se}}, \mbf{n} ) < \se  \}
= \{ \mscr{U}( \wh{\bs{\se}}, \mbf{n} ) \leq \se^\prime \} \cup \{
\se^\prime < \mscr{U}( \wh{\bs{\se}}, \mbf{n} ) < \se  \} = \{
\mscr{U}( \wh{\bs{\se}}, \mbf{n} ) \leq \se^\prime  \}$ for any $\se
\in ( \se^\prime, \se^{\prime \prime} ]$.  Hence, $\Pr \{ \mscr{U}(
\wh{\bs{\se}}, \mbf{n} ) < \se \; \tx{and} \; \mscr{E} \;
\tx{occurs}  \mid \se \} = \Pr \{ \mscr{U}( \wh{\bs{\se}}, \mbf{n} )
\leq \se^\prime \; \tx{and} \; \mscr{E} \; \tx{occurs} \mid \se \} =
\Pr \{ \wh{\bs{\se}} \leq \se^\prime, \; \mscr{U}( \wh{\bs{\se}},
\mbf{n} ) \leq \se^\prime \; \tx{and} \; \mscr{E} \; \tx{occurs}
\mid \se \}$ is non-increasing with respect to $\se \in (
\se^\prime, \se^{\prime \prime} ]$.  This implies that the supremum
of $\Pr \{ \mscr{U}( \wh{\bs{\se}}, \mbf{n} ) < \se \; \tx{and} \;
\mscr{E} \; \tx{occurs} \mid \se \}$ for $\se \in ( \se^\prime,
\se^{\prime \prime} ]$ is equal to $\Pr \{ \wh{\bs{\se}} \leq
\se^\prime, \; \mscr{U}( \wh{\bs{\se}}, \mbf{n} ) \leq \se^\prime \;
\tx{and} \; \mscr{E} \; \tx{occurs} \mid \se^\prime \}$.  Since the
argument holds for arbitrary consecutive distinct elements of
$I_\mscr{U} \cup \{ a, b\}$, we have established statement (II)
regarding $\Pr \{ \se > \mscr{U}( \wh{\bs{\se}}, \mbf{n}  ) \;
\tx{and} \; \mscr{E} \; \tx{occurs} \mid \se \}$ for $\se \in [a,
b]$.

To show statement (III), note that $\Pr \{ \se \leq \mscr{L}(
\wh{\bs{\se}}, \mbf{n}  ) \; \tx{and} \; \mscr{E} \; \tx{occurs}
\mid \se \}$ is no greater than $\Pr \{ a \leq \mscr{L}(
\wh{\bs{\se}}, \mbf{n}  ) \; \tx{and} \; \mscr{E} \; \tx{occurs}
\mid \se \}$ for any $\se \in [a, b]$. By the assumption that $\{ a
\leq \mscr{L}( \wh{\bs{\se}}, \mbf{n}  ) \} \subseteq \{
\wh{\bs{\se}} \geq b \}$, we have $\Pr \{ a \leq \mscr{L}(
\wh{\bs{\se}}, \mbf{n} ) \; \tx{and} \; \mscr{E} \; \tx{occurs} \mid
\se \} = \Pr \{\wh{\bs{\se}} \geq b, \; a \leq \mscr{L}(
\wh{\bs{\se}}, \mbf{n} ) \; \tx{and} \; \mscr{E} \; \tx{occurs} \mid
\se \}$ for any $\se \in [a, b]$.  As a result of the second
statement of Lemma \ref{ULE_Basic}, we have that $\Pr
\{\wh{\bs{\se}} \geq b, \; a \leq \mscr{L}( \wh{\bs{\se}}, \mbf{n} )
\; \tx{and} \; \mscr{E} \; \tx{occurs} \mid \se \}$ is
non-decreasing with respect to $\se \in [a, b]$.  It follows that
$\Pr \{\wh{\bs{\se}} \geq b, \; a \leq \mscr{L}( \wh{\bs{\se}},
\mbf{n} ) \; \tx{and} \; \mscr{E} \; \tx{occurs} \mid \se \} \leq
\Pr \{\wh{\bs{\se}} \geq b, \; a \leq \mscr{L}( \wh{\bs{\se}},
\mbf{n} ) \; \tx{and} \; \mscr{E} \; \tx{occurs} \mid b \}$ for any
$\se \in [a, b]$, which implies that $\Pr \{ \se \leq \mscr{L}(
\wh{\bs{\se}}, \mbf{n}  ) \; \tx{and} \; \mscr{E} \; \tx{occurs}
\mid \se \} \leq \Pr \{ a \leq \mscr{L}( \wh{\bs{\se}}, \mbf{n} ) \;
\tx{and} \; \mscr{E} \; \tx{occurs} \mid b \}$ for any $\se \in [a,
b]$.  On the other hand, $\Pr \{ \se \leq \mscr{L}( \wh{\bs{\se}},
\mbf{n}  ) \; \tx{and} \; \mscr{E} \; \tx{occurs} \mid \se \} \geq
\Pr \{ b \leq \mscr{L}( \wh{\bs{\se}}, \mbf{n}  ) \; \tx{and} \;
\mscr{E} \; \tx{occurs} \mid \se \}$ for any $\se \in [a, b]$.
Recalling that $\{ \wh{\bs{\se}} \geq \mscr{L} ( \wh{\bs{\se}},
\mbf{n} ) \}$ is a sure event, we have $\Pr \{ b \leq \mscr{L}(
\wh{\bs{\se}}, \mbf{n} ) \; \tx{and} \; \mscr{E} \; \tx{occurs} \mid
\se \} = \Pr \{ b \leq \mscr{L}( \wh{\bs{\se}}, \mbf{n}  ) \leq
\wh{\bs{\se}} \; \tx{and} \; \mscr{E} \; \tx{occurs} \mid \se \}$
for any $\se \in [a, b]$. Hence, applying the second statement of
Lemma \ref{ULE_Basic},  we have that $\Pr \{ b \leq \mscr{L}(
\wh{\bs{\se}}, \mbf{n}  ) \leq \wh{\bs{\se}} \; \tx{and} \; \mscr{E}
\; \tx{occurs} \mid \se \} \geq \Pr \{ b \leq \mscr{L}(
\wh{\bs{\se}}, \mbf{n}  ) \leq \wh{\bs{\se}} \; \tx{and} \; \mscr{E}
\; \tx{occurs} \mid a \} = \Pr \{ b \leq \mscr{L}( \wh{\bs{\se}},
\mbf{n}  ) \; \tx{and} \; \mscr{E} \; \tx{occurs} \mid a \}$ for any
$\se \in [a, b]$, which implies that $\Pr \{ \se \leq \mscr{L}(
\wh{\bs{\se}}, \mbf{n}  ) \; \tx{and} \; \mscr{E} \; \tx{occurs}
\mid \se \} \geq \Pr \{ b \leq \mscr{L}( \wh{\bs{\se}}, \mbf{n}  )
\; \tx{and} \; \mscr{E} \; \tx{occurs} \mid a \}$ for any $\se \in
[a, b]$.  So, we have established $\Pr \{ b \leq \mscr{L}(
\wh{\bs{\se}}, \mbf{n}  ) \; \tx{and} \; \mscr{E} \; \tx{occurs}
\mid a \} \leq \Pr \{ \se \leq \mscr{L}( \wh{\bs{\se}}, \mbf{n} ) \;
\tx{and} \; \mscr{E} \; \tx{occurs} \mid \se \} \leq \Pr \{ a \leq
\mscr{L}( \wh{\bs{\se}}, \mbf{n} ) \; \tx{and} \; \mscr{E} \;
\tx{occurs} \mid b \}$ for any $\se \in [a, b]$.  In a similar
manner, we can show that $\Pr \{ b < \mscr{L}( \wh{\bs{\se}},
\mbf{n}  ) \; \tx{and} \; \mscr{E} \; \tx{occurs} \mid a \} \leq \Pr
\{ \se < \mscr{L}( \wh{\bs{\se}}, \mbf{n} ) \; \tx{and} \; \mscr{E}
\; \tx{occurs} \mid \se \} \leq \Pr \{ a < \mscr{L}( \wh{\bs{\se}},
\mbf{n} ) \; \tx{and} \; \mscr{E} \; \tx{occurs} \mid b \}$ for any
$\se \in [a, b]$.

To show Statement (IV), applying the first statement of Lemma
\ref{ULE_Basic} to the special case that $\mcal{E} = \{
\wh{\bs{\se}} \leq z \}$, we have that $\Pr \{ \wh{\bs{\se}} \leq z
\mid \se \}$ is non-increasing with respect to $\se \in \Se$ no less
than $z$.   Applying the second statement of Lemma \ref{ULE_Basic}
to the special case that $\mcal{E} = \{ \wh{\bs{\se}} > z \}$, we
have that $\Pr \{ \wh{\bs{\se}} > z \mid \se \}$ is non-decreasing
with respect to $\se \in \Se$ no greater than $z$, which implies
that $\Pr \{ \wh{\bs{\se}} \leq z \mid \se \} = 1 - \Pr \{
\wh{\bs{\se}} > z \mid \se \}$ is non-increasing with respect to
$\se \in \Se$ no greater than $z$.  Therefore,  $\Pr \{
\wh{\bs{\se}} \leq z \mid \se \}$ is non-increasing with respect to
$\se \in \Se$. By a similar argument, $\Pr \{ \wh{\bs{\se}} \geq z
\mid \se \}$ is non-decreasing with respect to $\se \in \Se$. Note
that $\Pr \{ \se \leq \mscr{L}( \wh{\bs{\se}} ) \mid \se \}$ is no
greater than $\Pr \{ a \leq \mscr{L}( \wh{\bs{\se}} ) \mid \se \}$
for any $\se \in [a, b]$.  As a result of the monotonicity of
$\mscr{L} (.)$, we have that $\Pr \{ a \leq \mscr{L}( \wh{\bs{\se}}
) \mid \se \}$ is non-decreasing with respect to $\se \in [a, b]$.
It follows that $\Pr \{ a \leq \mscr{L}( \wh{\bs{\se}} ) \mid \se \}
\leq \Pr \{ a \leq \mscr{L}( \wh{\bs{\se}} ) \mid b \}$ for any $\se
\in [a, b]$, which implies that $\Pr \{ \se \leq \mscr{L}(
\wh{\bs{\se}} ) \mid \se \} \leq \Pr \{ a \leq \mscr{L}(
\wh{\bs{\se}} ) \mid b \}$ for any $\se \in [a, b]$.  Other
inequalities in Statement (IV) can be shown by a similar method.

This concludes the proof of Theorem \ref{Main_Bound_Gen}.

\sect{Proof of Theorem \ref{FiniteULE} }
\la{FiniteULE_Ap}  

It is easy to show that, for $x_i \in \{0, 1\}, \; i = 1, \cd, n$,
\[
\Pr \{ X_1 = x_1, \cd, X_n = x_n  \} = h (M, k) \qu \tx{where} \qu h
(M, k) = \bi{M}{k} \bi{N - M}{n - k}
 \li \slash \li [ \bi{n}{k} \bi{N}{n} \ri. \ri ]
\]
with $M = p N$ and $k = \sum_{i = 1}^n x_i$.  Note that $h (M, k) =
0$ if $M$ is smaller than $k$ or greater than $N - n + k$. For $k <
M \leq N - n + k$, we have {\small $\f{ h (M - 1, k) } { h (M, k) }
= \f{M - k} {M} \f{ N - M + 1 } { N - M - n + k + 1 } \leq 1$} if
and only if $M \leq \f{k}{n} (N + 1)$, or equivalently, $M \leq \lf
\f{k}{n} (N + 1) \rf$. It can be checked that $\f{k}{n} (N + 1) - (
N - n + k + 1)$ is equal to $(\f{k}{n} - 1) (N + 1 - n)$,  which is
negative for $k < n$. Hence, for $k < n$, we have that $\lf \f{k}{n}
(N + 1) \rf \leq N - n + k$ and consequently, the maximum of $h (M,
k )$ with respect to $M \in \{0, 1, \cd, N \}$ is achieved at $\li
\lf (N + 1) \f{k}{n} \ri \rf$. For $k = n$, we have $h (M, k ) = h
(M, n) = \bi{M}{n} \slash \bi{N}{n}$, of which the maximum with
respect to $M$ is attained at $M = N$. Therefore, for any $k \in
\{0, 1, \cd, n\}$, the maximum of $h (M, k )$ with respect to $M \in
\{0, 1, \cd, N \}$ is achieved at $\min \li \{ N, \li \lf (N + 1)
\f{k}{n} \ri \rf \ri \}$.  It follows that {\it $\min \{ 1, \f{1}{N}
\li \lf \f{N + 1}{n} \sum_{i=1}^n X_i \ri \rf \}$} is a MLE and also
a ULE for $p \in \Se$.  For simplicity of notations, let $\wh{p} =
\min \{ 1, \f{1}{N} \li \lf (N + 1) \f{k}{n} \ri \rf \}$.  We claim
that $| \wh{p} - \f{k}{n} | < \f{1}{N}$ for $0 \leq k \leq n$. To
prove such claim, we investigate two cases. In the case of $k = n$,
we have $\wh{p} = \f{k}{n} = 1$.  In the case of $k < n$, we have
$\wh{p} = \f{1}{N} \li \lf (N + 1) \f{k}{n} \ri \rf \leq \f{1}{N} (N
+ 1) \f{k}{n} < \f{k}{n} + \f{1}{N}$ and $\wh{p} > \f{1}{N} \li [ (N
+ 1) \f{k}{n} - 1 \ri ] = \f{k}{n} + \f{1}{N} ( \f{k}{n} - 1) \geq
\f{k}{n} - \f{1}{N}$.  The claim is thus proved.  In view of this
established claim and the fact that the difference between any pair
of values of $p \in \Se$ is no less than $\f{1}{N}$, we have that
{\small $\f{ \sum_{i=1}^n X_i }{n}$} is a ULE for $p \in \Se$. This
completes the proof of the theorem.

\sect{Proof of Theorem  \ref{Unbiased_Gen} } \la{App_Unbiased_Gen}

Define {\small $\wh{\bs{\mu}}_\ell = \f{ \sum_{i = 1}^{n_\ell} X_i }
{ n_\ell }$} and $F_\ell(x) = \Pr \{ \wh{\bs{\mu}}_\ell \leq x, \;
\bs{l} = \ell \}$ for $\ell = 1, \cd, s$, where $\bs{l}$ is the
index of stage when the sampling is terminated.  Let $\si^2$ denote
the variance of $X$.

To show statement (I), note that {\small \bee \li | \bb{E}
[\wh{\bs{\mu}} - \mu] \ri | & \leq &
\bb{E} | \wh{\bs{\mu}} - \mu | = \sum_{\ell = 1}^s  \int_{- \iy}^\iy  |x - \mu| \; d F_\ell (x)\\
& = & \sum_{\ell = 1}^s \li [ \int_{| x - \mu | <
\f{1}{\sqrt{n_\ell}} } | x - \mu | \; d F_\ell (x)  + \int_{| x -
\mu | \geq \f{1}{\sqrt{n_\ell}} } | x - \mu | \; d F_\ell (x) \ri ]\\
& =  & \sum_{\ell = 1}^s  \int_{| x - \mu | < \f{1}{\sqrt{n_\ell}} } | x - \mu | \; d F_\ell (x)  + \sum_{\ell = 1}^s \int_{| x - \mu | \geq
\f{1}{\sqrt{n_\ell}} } | x - \mu | \; d F_\ell (x)\\
& \leq  & \sum_{\ell = 1}^s  \int_{| x - \mu | < \f{1}{\sqrt{n_\ell}} } \f{1}{\sqrt{n_\ell}} \; d F_\ell (x)  + \sum_{\ell = 1}^s \int_{|
x - \mu | \geq \f{1}{\sqrt{n_\ell}} } \sqrt{n_\ell} | x - \mu |^2 \; d F_\ell (x) \\
& \leq & \sum_{\ell = 1}^s  \f{1}{\sqrt{n_\ell}} \; \int_{| x - \mu | < \f{1}{\sqrt{n_\ell}} }  d F_\ell (x)  + \sum_{\ell = 1}^s \sqrt{n_\ell}
\int_{- \iy}^\iy  | x - \mu |^2 \; d F_\ell (x)\\
& \leq & \sum_{\ell = 1}^s  \f{1}{\sqrt{n_\ell}} \Pr \li \{ | \wh{\bs{\mu}}_\ell - \mu | <
\f{1}{\sqrt{n_\ell}}, \; \bs{l} = \ell
\ri \} + \sum_{\ell = 1}^s \sqrt{n_\ell} \; \bb{E} [ | \wh{\bs{\mu}}_\ell  - \mu |^2 ]\\
& \leq & \f{1}{\sqrt{n_1}} \sum_{\ell = 1}^s   \Pr \li \{ \bs{l} = \ell \ri \} + \sum_{\ell = 1}^s \sqrt{n_\ell} \; \f{\si^2}{ n_\ell } =
\f{1}{\sqrt{n_1}} + \si^2 \sum_{\ell = 1}^s \f{1}{\sq{n_\ell}}. \eee} By the assumption that $\inf_{\ell > 0} \f{n_{\ell + 1}}{n_\ell} > 1$, we
have that, there exists a positive number $\ro$ such that $n_\ell \geq (1 + \ro)^{2( \ell - 1 )} n_1$ for all $\ell
> 1$.  Hence, \bee \li | \bb{E} [\wh{\bs{\mu}} - \mu] \ri | & \leq &
\bb{E} |\wh{\bs{\mu}} - \mu| \leq \f{1}{\sqrt{n_1}} + \si^2
\sum_{\ell = 1}^s \f{1}{\sqrt{n_\ell}} \leq  \f{1}{\sqrt{n_1}} +
\si^2 \sum_{\ell =
1}^s \f{1}{ \sq{n_1} (1 + \ro)^{\ell - 1} }\\
& \leq & \f{1}{\sqrt{n_1}} +  \f{\si^2}{\sq{n_1}} \sum_{\ell =
1}^\iy \f{1}{ (1 + \ro)^{\ell - 1} } \leq \f{1}{\sqrt{n_1}} +
\f{\si^2}{\sq{n_1}}  \f{1 + \ro}{\ro} \to 0 \eee as $n_1 \to \iy$.
Moreover,  \bee \bb{E} \li [ | \wh{\bs{\mu}} - \mu |^2 \ri ] & = &
\sum_{\ell = 1}^s \int_{- \iy}^\iy  |x - \mu|^2 \; d F_\ell (x) \leq
\sum_{\ell = 1}^s \bb{E} \li [ | \wh{\bs{\mu}}_\ell - \mu |^2 \ri ] \\
& = & \si^2 \sum_{\ell = 1}^s \f{1}{n_\ell} \leq \si^2 \sum_{\ell =
1}^\iy \f{1}{n_1 (1 + \ro)^{2 (\ell - 1)}} = \f{ \si^2 }{ n_1 } \f{
(1 + \ro)^2}{ \ro (2 + \ro)} \to 0 \eee as $n_1 \to \iy$.  This
completes the proof of statement (I).

Now we shall show statement (II).  Since $X$ is a bounded variable,
there exists a positive number $C$ such that $| X - \mu | < C$. By
Chebyshev's inequality, we have $\Pr \{ | \wh{\bs{\mu}}_\ell - \mu |
\geq \f{1}{\sqrt[4]{n_\ell}} \} \leq \f{ \si^2}{\sqrt{ n_\ell} }$
for $\ell = 1, \cd, s$. Therefore, for $k = 1, 2, \cd$, {\small \bee
\bb{E} \li [ | \wh{\bs{\mu}} - \mu |^k \ri ] & = & \sum_{\ell = 1}^s  \int_{- \iy}^\iy  |x - \mu|^k \; d F_\ell (x)\\
& = & \sum_{\ell = 1}^s \li [ \int_{| x - \mu | <
\f{1}{\sqrt[4]{n_\ell}} } | x - \mu |^k \; d F_\ell (x)  + \int_{| x
- \mu | \geq \f{1}{\sqrt[4]{n_\ell}} } | x - \mu |^k \; d F_\ell (x) \ri ]\\
& =  & \sum_{\ell = 1}^s  \int_{| x - \mu | <
\f{1}{\sqrt[4]{n_\ell}} } | x - \mu |^k \; d F_\ell (x)  +
\sum_{\ell = 1}^s \int_{| x - \mu | \geq \f{1}{\sqrt[4]{n_\ell}} } | x - \mu |^k \; d F_\ell (x) \\
& \leq  & \sum_{\ell = 1}^s  \li ( \f{1}{\sqrt[4]{n_\ell}} \ri )^k
 \int_{| x - \mu | < \f{1}{\sqrt[4]{n_\ell}} } d F_\ell (x) + C^k
\sum_{\ell = 1}^s \int_{| x - \mu | \geq \f{1}{\sqrt[4]{n_\ell}} }  d F_\ell (x) \\
& = & \sum_{\ell = 1}^s  \li ( \f{1}{\sqrt[4]{n_\ell}} \ri )^k \Pr
\li \{ | \wh{\bs{\mu}}_\ell - \mu | < \f{1}{\sqrt[4]{n_\ell}}, \;
\bs{l} = \ell \ri \} + C^k \sum_{\ell = 1}^s \Pr \li \{ |
\wh{\bs{\mu}}_\ell - \mu | \geq \f{1}{\sqrt[4]{n_\ell}}, \; \bs{l} = \ell \ri \}\\
& \leq & \li ( \f{1}{\sqrt[4]{n_1}} \ri )^k \sum_{\ell = 1}^s   \Pr
\li \{ \bs{l} = \ell \ri \} + C^k \sum_{\ell = 1}^s \Pr \li \{ |
\wh{\bs{\mu}}_\ell - \mu | \geq \f{1}{\sqrt[4]{n_\ell}} \ri \}\\
& = & \li ( \f{1}{\sqrt[4]{n_1}} \ri )^k + C^k \sum_{\ell = 1}^s \Pr
\li \{ | \wh{\bs{\mu}}_\ell - \mu | \geq \f{1}{\sqrt[4]{n_\ell}} \ri
\} \leq \li ( \f{1}{\sqrt[4]{n_1}} \ri )^k + C^k \sum_{\ell = 1}^s
\f{ \si^2}{\sqrt{ n_\ell }} \to 0 \eee} as $n_1 \to \iy$. Since $\li
| \bb{E} [\wh{\bs{\mu}} - \mu] \ri | \leq \bb{E} | \wh{\bs{\mu}} -
\mu |$, we have that $\bb{E} [\wh{\bs{\mu}} - \mu] \to 0$ as $n_1
\to \iy$.  This completes the proof of statement (II).

\sect{Proof of Theorem \ref{Thm_CDV}} \la{App_Thm_CDV}

We only show the last statement of Theorem \ref{Thm_CDV}.   Note
that {\small \bee n_s - n_1 \; \Pr \{ \bs{l} = 1 \} & = & n_s \; \Pr
\{ \bs{l} \leq s \} - n_1 \; \Pr \{
\bs{l} \leq 1  \} =
 \sum_{\ell = 2}^s \li ( n_\ell \; \Pr \{ \bs{l} \leq \ell \} - n_{\ell - 1} \; \Pr \{ \bs{l} \leq {\ell - 1} \} \ri )\\
& = & \sum_{\ell = 2}^s n_\ell \; ( \Pr \{ \bs{l} \leq \ell \} - \Pr
\{ \bs{l} \leq {\ell - 1} \} ) + \sum_{\ell = 2}^s
(n_\ell - n_{\ell -1} ) \; \Pr \{ \bs{l} \leq {\ell - 1} \}\\
& = & \sum_{\ell = 2}^s n_\ell \; \Pr \{ \bs{l} = \ell \}  +
\sum_{\ell = 2}^s (n_\ell - n_{\ell -1} ) \; \Pr \{ \bs{l} \leq
{\ell - 1} \}, \eee} from which we obtain $n_s - \sum_{\ell = 1}^s
n_\ell \; \Pr \{ \bs{l} = \ell \}  = \sum_{\ell = 2}^s \; (n_\ell -
n_{\ell - 1}) \; \Pr \{ \bs{l} \leq {\ell - 1} \}$. Observing that
$n_s = n_1 + \sum_{\ell = 2}^s \; (n_\ell - n_{\ell - 1})$, we have
\bee \bb{E} [ \mathbf{n} ] & =  & \sum_{\ell = 1}^s n_\ell \; \Pr \{
\bs{l} = \ell  \}  =  n_s - \li ( n_s - \sum_{\ell = 1}^s n_\ell \;
\Pr \{
\bs{l} = \ell \} \ri ) \\
& = & n_1 + \sum_{\ell = 2}^s \; (n_\ell - n_{\ell - 1}) -
\sum_{\ell = 2}^s \; (n_\ell - n_{\ell - 1}) \; \Pr \{ \bs{l} \leq {\ell - 1} \}\\
& = & n_1 + \sum_{\ell = 2}^s \; (n_\ell - n_{\ell - 1}) \; \Pr \{
\bs{l} > \ell - 1 \} =  n_1 + \sum_{\ell = 1}^{s - 1} \; (n_{\ell +
1} - n_{\ell}) \; \Pr \{ \bs{l} > \ell \}. \eee

\section{Proof of Theorem \ref{discrete_sum} } \la{App_discrete_sum}

To prove Theorem \ref{discrete_sum}, we shall only provide the proof
of statement (I), since the proof of statement (II) is similar.  As
a consequence of the assumption that $f(k+1) - f(k) \leq f(k) -
f(k-1)$  for $a < k < b$, we have {\small $\frac{f(b) -f(k) }{b-k}
\leq f(k+1) - f(k) \leq f(k) - f(k-1) \leq \frac{f(k) -f(a)}{k-a}$}
for $a < k < b$. Hence,
\begin{eqnarray*}
\frac{f(b) -f(a)}{b-a} & = & \frac{ \frac{f(b) -f(k) }{b-k} (b-k) +
\frac{f(k) -f(a)}{k-a} (k-a) } {  b-a  }\\
& \leq  & \frac{ \frac{f(k) -f(a) }{k-a} (b-k) + \frac{f(k)
-f(a)}{k-a} (k-a) } {  b-a  } = \frac{f(k) -f(a)}{k-a},
\end{eqnarray*}
which implies $f(k) \geq f(a) + \frac{f(b) -f(a)}{b-a} (k-a)$ for $a
\leq k \leq b$ and it follows that
\[ \sum_{k=a}^b f(k) \geq  (b-a+1) f(a) + \frac{f(b) -f(a)}{b-a}
\sum_{k=a}^b (k-a) = \f{ (b-a+1) [f(b) + f(a)] }{2}.
\]
Again by virtue of the assumption that $f(k+1) - f(k) \leq f(k) -
f(k-1)$ for $a < k < b$, we have
\[
f(k) - f(a) = \sum_{l=a}^{k-1} [f(l+1)-f(l)]  \leq \sum_{l=a}^{k-1}
[f(a+1)-f(a)]  = (k-a) [f(a+1)-f(a)],
\]
\[
f(k) - f(b) = \sum_{l=k}^{b-1} [f(l)-f(l + 1)]  \leq
\sum_{l=k}^{b-1} [f(b - 1)-f(b)] = (k - b) [f(b)-f(b-1)]
\]
for $a < k < b$. Making use of the above established inequalities,
we have
\begin{eqnarray*} \sum_{k=a}^b f(k) & = & (b-a+1) f(a) +
\sum_{k=a}^{i} [f(k) - f(a)] + \sum_{k=i+1}^b [f(b) - f(a)] +
\sum_{k=i+1}^b [f(k) - f(b)]\\
& \leq & (b-a+1) f(a) + \sum_{k=a}^{i} (k - a) [f(a+1) - f(a)] \\
&   & + (b-i) [f(b) - f(a)] + \sum_{k=i+1}^b (k - b) [f(b) - f(b-1)] \\
& = & \al (i) f(a) + \ba (i)  f(b)
\end{eqnarray*}
for $a < i < b$.  Observing that \[ j = a + \frac{f(b) - f(a) + (a -
b) [ f(b) - f(b-1) ]} { f(a+1) + f(b-1) - f(a) - f(b)} = a + \f{ b -
a - (1 - r_{a,b}) (1 - r_b )^{-1} } { 1 + r_{a, b}  ( 1 - r_a ) ( 1
- r_b )^{-1} }
\]
is the solution of equation $f(a) + (i-a) [f(a+1) - f(a)] =  f(b) -
(b-i) [f(b) - f(b-1)]$ with respect to $i$,  we can conclude based
on a geometric argument that the minimum gap between the lower and
upper bounds in (\ref{concavity}) is achieved at $i$ such that $\lf
j \rf \leq i \leq \lc j \rc$.  This completes the proof of Theorem
\ref{discrete_sum}.

\section{Proof of Theorem \ref{continuous_integration} } \la{App_continuous_integration}

To prove Theorem \ref{continuous_integration}, we shall only provide
the proof of statement (I), since the proof of statement (II) is
similar.  Define $g(x) = f(a) + \f{ f(b) - f(a)  } { b - a } (x -
a)$ and
\[
h(x) = \left\{\begin{array}{ll}
   f(a) + f^\prime(a) \; (x-a) \;&  \;{\rm if}\; x \leq t, \\
   f(b) + f^\prime(b) \; (x-b) \;&  \; {\rm if} \; x > t
\end{array} \right.
\]
for $t \in (a, b)$.  By the assumption that $f(x)$ is concave over
$[a, b]$, we have $g(x) \leq f(x) \leq h(x)$ for $x \in [a, b]$ and
it follows that {\small $\int_{a}^b f(x) dx  \geq \int_{a}^ b g(y) d
y = \f{[f(a) + f(b)](b-a) } { 2 }$} and $\int_{a}^b f(x) dx  \leq
\int_{a}^ b g(y) d y + \int_{a}^ b [ h(y) - g(y) ] d y$ with
$\int_{a}^ b [ h(y) - g(y) ] d y =  \int_{a}^t [ h(y) - g(y) ] d y +
\int_t^b [ h(y) - g(y) ] d y = \vDe (t)$.  It can be shown by
differentiation that $\vDe (t)$ attains its minimum at {\small $t =
\f{ f(b) - f(a) + a f^\prime(a) - b f^\prime(b) } { f^\prime(a) -
f^\prime(b) }$}.  This completes the proof of Theorem
\ref{continuous_integration}.

\section{Proofs of Theorems for Estimation of Binomial Parameters}

\subsection{Proof of Theorem  \ref{Bino_ABS_CDF_CH_MA} } \la{App_Bino_ABS_CDF_CH_MA}

We need some preliminary results.  The following lemma can be
readily derived from Hoeffding's inequalities stated in Lemma
\ref{Hoe_Mas}.

\beL \la{decb}

{\small $S_{\mrm{B}} (k, n, p) \leq \exp ( n \mscr{M}_{\mrm{B}} (
\f{k}{n}, p  )  )$} for $0 \leq k \leq n p$. Similarly, {\small $1 -
S_{\mrm{B}} (k - 1, n, p) \leq \exp  ( n \mscr{M}_{\mrm{B}} (
\f{k}{n}, p  )  )$} for $n p \leq k \leq n$. \eeL

\beL \la{lemax}

$\mscr{M}_{\mrm{B}} (z, z - \vep) \leq - 2 \vep^2$ for $0 < \vep < z
< 1$.  Similarly, $\mscr{M}_{\mrm{B}} (z, z + \vep) \leq - 2 \vep^2$
for $0 < z < 1 - \vep < 1$.

\eeL

\bpf

It can be shown that {\small $\f{\pa \mscr{M}_{\mrm{B}} (\mu + \vep,
\mu) } { \pa \vep } = \ln \li ( \f{\mu}{\mu + \vep} \f{ 1 - \mu -
\vep} {1 - \mu} \ri )$} and {\small $\f{\pa^2 \mscr{M}_{\mrm{B}}
(\mu + \vep, \mu) } { \pa \vep^2 } = \f{1}{ (\mu + \vep) (\mu + \vep
- 1) }$} for $0 < \vep < 1 - \mu < 1$. Observing that
$\mscr{M}_{\mrm{B}} (\mu, \mu) = 0$ and {\small $\f{\pa
\mscr{M}_{\mrm{B}} (\mu + \vep, \mu) } { \pa \vep }|_{\vep = 0} =
0$}, by Taylor's expansion formula, we have that there exists a real
number $\vep^* \in (0, \vep)$ such that {\small $\mscr{M}_{\mrm{B}}
(\mu + \vep, \mu) = \f{\vep^2}{2} \f{1}{ (\mu + \vep^*) (\mu +
\vep^* - 1) }$} where the right side is seen to be no greater than $
- 2 \vep^2$. Hence, letting $z = \mu + \vep$, we have
$\mscr{M}_{\mrm{B}} (z, z - \vep) \leq - 2 \vep^2$ for $0 < \vep < z
< 1$.  This completes the proof of the first statement of  the
lemma.

Similarly, it can be verified that {\small $\f{\pa
\mscr{M}_{\mrm{B}} (\mu - \vep, \mu) } { \pa \vep } = - \ln \li (
\f{\mu}{\mu - \vep} \f{ 1 - \mu + \vep} {1 - \mu} \ri )$} and
{\small $\f{\pa^2 \mscr{M}_{\mrm{B}} (\mu - \vep, \mu) } { \pa
\vep^2 } = \f{1}{ (\mu - \vep) (\mu - \vep - 1) }$} for $0 < \vep <
\mu < 1$. Observing that $\mscr{M}_{\mrm{B}} (\mu, \mu) = 0$ and
{\small $\f{\pa \mscr{M}_{\mrm{B}} (\mu - \vep, \mu)  } { \pa \vep
}|_{\vep = 0} = 0$}, by Taylor's expansion formula, we have that
there exists a real number $\vep^\star \in (0, \vep)$ such that
{\small $\mscr{M}_{\mrm{B}} (\mu - \vep, \mu) = \f{\vep^2}{2} \f{1}{
(\mu - \vep^\star) (\mu - \vep^\star - 1) }$} where the right side
is seen to be no greater than $ - 2 \vep^2$.  Therefore, letting $z
= \mu - \vep$, we have $\mscr{M}_{\mrm{B}} (z, z + \vep) \leq - 2
\vep^2$ for $0 < z < 1 - \vep < 1$. This completes the proof of the
second statement of the lemma. \epf

\beL \la{absDS1}

$\{ F_{\wh{\bs{p}}_s}  \li (\wh{\bs{p}}_s, \wh{\bs{p}}_s + \vep \ri
) \leq \ze \de, \; G_{\wh{\bs{p}}_s} \li (\wh{\bs{p}}_s,
\wh{\bs{p}}_s - \vep \ri ) \leq \ze \de \}$ is a sure event.
 \eeL

 \bpf
By the definition of sample sizes, we have $n_s \geq \li \lc \f{\ln
(\ze \de) }{ - 2 \vep^2}  \ri \rc \geq \f{\ln (\ze \de) }{ - 2
\vep^2}$ and consequently  $\f{\ln (\ze \de) }{n_s} \geq - 2
\vep^2$.  By Lemmas \ref{decb} and \ref{lemax}, we have {\small \bee
\Pr \{ F_{\wh{\bs{p}}_s} \li (\wh{\bs{p}}_s, \wh{\bs{p}}_s + \vep
\ri ) \leq \ze \de \} & = & \Pr \li \{ S_{\mrm{B}} \li ( K_s, n_s,
\wh{\bs{p}}_s + \vep \ri ) \leq \ze
\de \ri \}\\
& \geq & \Pr \li \{ \mscr{M}_{\mrm{B}} \li (\wh{\bs{p}}_s,
\wh{\bs{p}}_s + \vep \ri ) \leq \f{\ln (\ze \de)}{n_s} \ri \} \geq
\Pr \li \{ \mscr{M}_{\mrm{B}} \li (\wh{\bs{p}}_s, \wh{\bs{p}}_s +
\vep \ri ) \leq - 2 \vep^2 \ri \} = 1, \\
\Pr \{ G_{\wh{\bs{p}}_s} \li (\wh{\bs{p}}_s, \wh{\bs{p}}_s - \vep
\ri ) \leq \ze \de \} & = & \Pr \li \{ 1 - S_{\mrm{B}} \li
(K_s - 1, n_s, \wh{\bs{p}}_s - \vep \ri ) \leq \ze \de \ri \}\\
& \geq & \Pr \li \{ \mscr{M}_{\mrm{B}} \li (\wh{\bs{p}}_s,
\wh{\bs{p}}_s - \vep \ri ) \leq \f{\ln (\ze \de)}{n_s} \ri \} \geq
\Pr \li \{ \mscr{M}_{\mrm{B}} \li (\wh{\bs{p}}_s, \wh{\bs{p}}_s -
\vep \ri ) \leq - 2 \vep^2 \ri \} = 1
 \eee} which immediately
implies the lemma.

 \epf

\beL \la{lem888}  Let $0 < \vep < \f{1}{2}$. Then,
$\mscr{M}_{\mrm{B}}(z, z + \vep) \geq  \mscr{M}_{\mrm{B}}(z, z -
\vep)$ for $z \in \li [0, \f{1}{2} \ri ]$, and
$\mscr{M}_{\mrm{B}}(z, z + \vep) <  \mscr{M}_{\mrm{B}}(z, z - \vep)$
for $z \in \li ( \f{1}{2}, 1 \ri ]$.  \eeL

\bpf

By the definition of the function $\mscr{M}_{\mrm{B}} (.,.)$, we
have that $\mscr{M}_{\mrm{B}} (z, \mu) = - \iy$ for $z \in [0, 1]$
and $\mu \notin (0, 1)$. Hence, the lemma is trivially true for $0
\leq z \leq \vep$ or $1 - \vep \leq z \leq 1$.  It remains to show
the lemma for $z \in (\vep, 1 - \vep)$.  This can be accomplished by
noting that $\mscr{M}_{\mrm{B}}(z, z + \vep) - \mscr{M}_{\mrm{B}}(z,
z - \vep) = 0$ for $\vep = 0$ and that
\[
\f{ \pa [\mscr{M}_{\mrm{B}}(z, z + \vep) - \mscr{M}_{\mrm{B}}(z, z -
\vep)] } {\pa \vep} = \f{ 2 \vep^2 (1 - 2 z) } {(z^2 - \vep^2) [(1 -
z)^2 - \vep^2] }, \qu \fa z \in (\vep, 1 - \vep)
\]
where the partial derivative is seen to be positive for $z \in \li (
\vep, \f{1}{2} \ri )$ and negative for $z \in \li ( \f{1}{2}, 1 -
\vep \ri )$.  \epf

\beL \la{DS1}

 $\{ \mscr{M}_{\mrm{B}} \li ( \f{1}{2} - \li |\f{1}{2} - \wh{\bs{p}}_s \ri | ,
\f{1}{2} - \li |\f{1}{2} - \wh{\bs{p}}_s \ri | + \vep \ri ) \leq
\f{\ln ( \ze \de )} { n_s } \}$ is a sure event.

\eeL

\bpf To show the lemma, it suffices to show {\small
$\mscr{M}_{\mrm{B}} \li ( \f{1}{2} - \li |\f{1}{2} - z \ri | ,
\f{1}{2} - \li |\f{1}{2} - z \ri | + \vep \ri ) \leq \f{\ln ( \ze
\de )} { n_s }$} for any $z \in [0, 1]$, since $0 \leq \wh{\bs{p}}_s
(\om) \leq 1$ for any $\om \in \Om$. By the definition of sample
sizes, we have {\small $n_s \geq \li \lc \f{  \ln ( \ze \de ) } { -
2 \vep^2 } \ri \rc \geq \f{ \ln ( \ze \de ) } { - 2 \vep^2 }$} and
thus $\f{\ln ( \ze \de )} { n_s } \geq - 2 \vep^2$. Hence, it is
sufficient to show {\small $\mscr{M}_{\mrm{B}}  ( \f{1}{2} - |
\f{1}{2} - z | , \f{1}{2} - | \f{1}{2} - z | + \vep  ) \leq - 2
\vep^2$} for any $z \in [0, 1]$.   This can be accomplished by
considering four cases as follows.

In the case of $z = 0$, we have $\mscr{M}_{\mrm{B}} \li ( \f{1}{2} -
\li |\f{1}{2} - z \ri | , \f{1}{2} - \li |\f{1}{2} - z \ri | + \vep
\ri ) = \mscr{M}_{\mrm{B}} (0, \vep) = \ln (1 - \vep) < - 2 \vep^2$,
where the last inequality follows from the fact that $\ln (1 - x) <
- 2 x^2$ for any $x \in (0, 1)$.

In the case of $0 < z \leq \f{1}{2}$, we have $\mscr{M}_{\mrm{B}}
\li ( \f{1}{2} - \li |\f{1}{2} - z \ri | , \f{1}{2} - \li |\f{1}{2}
- z \ri | + \vep \ri ) = \mscr{M}_{\mrm{B}} (z, z + \vep) \leq - 2
\vep^2$, where the inequality follows from Lemma \ref{lemax} and the
fact that $0 < z \leq \f{1}{2} < 1 - \vep$.

In the case of $\f{1}{2} < z < 1$, we have $\mscr{M}_{\mrm{B}} \li (
\f{1}{2} - \li |\f{1}{2} - z \ri | , \f{1}{2} - \li |\f{1}{2} - z
\ri | + \vep \ri ) = \mscr{M}_{\mrm{B}} (1 - z, 1 - z + \vep) =
\mscr{M}_{\mrm{B}} (z, z - \vep) \leq - 2 \vep^2$, where the
inequality follows from Lemma \ref{lemax} and the fact that $\vep  <
\f{1}{2} <  z < 1$.

In the case of $z = 1$, we have $\mscr{M}_{\mrm{B}} \li ( \f{1}{2} -
\li |\f{1}{2} - z \ri | , \f{1}{2} - \li |\f{1}{2} - z \ri | + \vep
\ri ) = \mscr{M}_{\mrm{B}} (0, \vep) = \ln (1 - \vep) < - 2 \vep^2$.

The proof of the lemma is thus completed.

\epf

\beL \la{abs13}

$\{ \mscr{M}_{\mrm{B}} \li ( \f{1}{2} - \li |\f{1}{2} -
\wh{\bs{p}}_\ell \ri | , \f{1}{2} - \li |\f{1}{2} - \wh{\bs{p}}_\ell
\ri | + \vep \ri ) \leq \f{\ln ( \ze \de )} { n_s } \} \subseteq \{
\mscr{M}_{\mrm{B}} ( \wh{\bs{p}}_\ell,  \wh{\bs{p}}_\ell + \vep )
\leq \f{\ln (\ze \de) }{n_\ell}, \; \mscr{M}_{\mrm{B}} (
\wh{\bs{p}}_\ell, \wh{\bs{p}}_\ell - \vep  ) \leq \f{\ln (\ze \de)
}{n_\ell} \}$ for $\ell = 1, \cd, s$.

\eeL

\bpf

Let {\small $\om \in \{ \mscr{M}_{\mrm{B}} \li ( \f{1}{2} - \li
|\f{1}{2} - \wh{\bs{p}}_\ell \ri | , \f{1}{2} - \li |\f{1}{2} -
\wh{\bs{p}}_\ell \ri | + \vep \ri ) \leq \f{\ln ( \ze \de )} { n_s }
\}$} and {\small $\wh{p}_\ell = \wh{\bs{p}}_\ell (\om)$}.  To show
the lemma, it suffices to show {\small $\max  \{ \mscr{M}_{\mrm{B}}
( \wh{p}_\ell, \wh{p}_\ell + \vep), \; \mscr{M}_{\mrm{B}} \li (
\wh{p}_\ell, \wh{p}_\ell - \vep \ri ) \} \leq \f{\ln (\ze \de)
}{n_\ell}$} by considering two cases: Case (i) $\wh{p}_\ell \leq
\f{1}{2}$; Case (ii) $\wh{p}_\ell
> \f{1}{2}$.

In Case (i), we have {\small $\mscr{M}_{\mrm{B}} ( \wh{p}_\ell,
\wh{p}_\ell + \vep) = \mscr{M}_{\mrm{B}} \li ( \f{1}{2} - \li
|\f{1}{2} - \wh{p}_\ell \ri | , \f{1}{2} - \li |\f{1}{2} -
\wh{p}_\ell \ri | + \vep \ri )  \leq \f{\ln (\ze \de) }{n_\ell}$}.
Since $ \wh{p}_\ell \leq \f{1}{2}$, by Lemma \ref{lem888}, we have
{\small $\mscr{M}_{\mrm{B}} ( \wh{p}_\ell, \wh{p}_\ell - \vep) \leq
\mscr{M}_{\mrm{B}} ( \wh{p}_\ell, \wh{p}_\ell + \vep) \leq \f{\ln
(\ze \de) }{n_\ell}$}.

In Case (ii), we have {\small $\mscr{M}_{\mrm{B}} ( \wh{p}_\ell,
\wh{p}_\ell - \vep) = \mscr{M}_{\mrm{B}}(1 - \wh{p}_\ell, 1 -
\wh{p}_\ell + \vep) = \mscr{M}_{\mrm{B}} \li ( \f{1}{2} - \li
|\f{1}{2} - \wh{p}_\ell \ri | , \f{1}{2} - \li |\f{1}{2} -
\wh{p}_\ell \ri | + \vep \ri )  \leq \f{\ln (\ze \de) }{n_\ell}$}.
Since $\wh{p}_\ell > \f{1}{2}$, by Lemma \ref{lem888}, we have
{\small $\mscr{M}_{\mrm{B}} ( \wh{p}_\ell, \wh{p}_\ell + \vep) <
\mscr{M}_{\mrm{B}} ( \wh{p}_\ell, \wh{p}_\ell - \vep) \leq \f{\ln
(\ze \de) }{n_\ell}$}.  This completes the proof of the lemma.

\epf

\beL \la{absD}

$ \{ ( | \wh{\bs{p}}_s - \f{1}{2} | - \f{2\vep}{3} )^2 \geq \f{1}{4}
+ \f{  n_s \; \vep^2 } {2 \ln (\ze \de) }  \}$ is a sure event.

\eeL \bpf

By the definition of sample sizes, we have {\small $n_s \geq \li \lc
\f{ \ln \f{1}{\ze \de} } { 2 \vep^2  } \ri \rc \geq \f{ \ln
\f{1}{\ze \de} } { 2 \vep^2  }$}, which implies that $\f{1}{4} + \f{
n_s \; \vep^2 } {2 \ln (\ze \de) } \leq 0$. Since $\{ \li ( \li |
\wh{\bs{p}}_s - \f{1}{2} \ri | - \f{2\vep}{3} \ri )^2 \geq 0 \}$ is
a sure event, it follows that {\small $\{ ( | \wh{\bs{p}}_s -
\f{1}{2} | - \f{2\vep}{3} )^2 \geq \f{1}{4} + \f{  n_s \; \vep^2 }
{2 \ln (\ze \de) } \}$} is a sure event. This completes the proof of
the lemma.

\epf

\beL \la{abs381} {\small $\{ ( | \wh{\bs{p}}_\ell - \f{1}{2} | -
\f{2\vep}{3} )^2 \geq \f{1}{4} + \f{  n_\ell \; \vep^2 } {2 \ln (\ze
\de) }  \} \leu  \{ \mscr{M}_{\mrm{B}} \li ( \wh{\bs{p}}_\ell,
\wh{\bs{p}}_\ell + \vep \ri ) \leq \f{\ln (\ze \de)}{n_\ell}, \;
\mscr{M}_{\mrm{B}} \li ( \wh{\bs{p}}_\ell, \wh{\bs{p}}_\ell - \vep
\ri ) \leq \f{\ln (\ze \de)}{n_\ell} \} $} for $\ell = 1, \cd, s$.

\eeL

\bpf

Let {\small $\om \in \{ ( | \wh{\bs{p}}_\ell - \f{1}{2} | -
\f{2\vep}{3} )^2 \geq \f{1}{4} + \f{  n_\ell \; \vep^2 } {2 \ln (\ze
\de) }  \}$ and $\wh{p}_\ell = \wh{\bs{p}}_\ell (\om)$}. Then, \be
\la{eqD}
 \li ( \li |
\wh{p}_\ell - \f{1}{2} \ri | - \f{2\vep}{3} \ri )^2 \geq
 \f{1}{4} + \f{  n_\ell \vep^2 } {2 \ln (\ze \de) }. \ee
To show the lemma, it suffices to show $\mscr{M} \li (\wh{p}_\ell,
\wh{p}_\ell + \vep \ri )  \leq \f{\ln (\ze \de)}{n_\ell}$ and
$\mscr{M} \li (\wh{p}_\ell, \wh{p}_\ell - \vep \ri )  \leq \f{\ln
(\ze \de)}{n_\ell}$.  For the purpose of proving the first
inequality, we need to show \be \la{eqDD} \li ( \wh{p}_\ell -
\f{1}{2} + \f{2\vep}{3} \ri )^2 \geq \f{1}{4} + \f{ n_\ell \vep^2 }
{2 \ln (\ze \de) }.  \ee Clearly, (\ref{eqDD}) holds if {\small
$\f{1}{4} + \f{ n_\ell \vep^2 } {2 \ln (\ze \de) } \leq 0$}.  It
remains to show (\ref{eqDD}) under the condition that {\small
$\f{1}{4} + \f{ n_\ell \vep^2 } {2 \ln (\ze \de) }
> 0$}.  Note that (\ref{eqD}) implies either \be \la{eqD1} \li |
\wh{p}_\ell - \f{1}{2} \ri | - \f{2\vep}{3} \geq \sq{ \f{1}{4} + \f{
n_\ell \vep^2 } {2 \ln (\ze \de) } } \ee or \be \la{eqD2} \li |
\wh{p}_\ell - \f{1}{2} \ri | - \f{2\vep}{3} \leq - \sq{ \f{1}{4} +
\f{  n_\ell \vep^2 } {2 \ln (\ze \de) } }. \ee Since (\ref{eqD1})
implies either {\small $\wh{p}_\ell - \f{1}{2} + \f{2\vep}{3}  \geq
\f{4\vep}{3} + \sq{ \f{1}{4} + \f{  n_\ell \vep^2 } {2 \ln (\ze \de)
} } > \sq{ \f{1}{4} + \f{  n_\ell \vep^2 } {2 \ln (\ze \de) } }$} or
{\small $\wh{p}_\ell - \f{1}{2} + \f{2\vep}{3} \leq - \sq{ \f{1}{4}
+ \f{  n_\ell \vep^2 } {2 \ln (\ze \de) } }$}, it must be true that
(\ref{eqD1}) implies (\ref{eqDD}).  On the other hand, (\ref{eqD2})
also implies (\ref{eqDD}) because (\ref{eqD2}) implies {\small $\sq{
\f{1}{4} + \f{  n_\ell \vep^2 } {2 \ln (\ze \de) } } \leq
\wh{p}_\ell - \f{1}{2}  + \f{2\vep}{3}$}. Hence, we have established
(\ref{eqDD}).

 In the case of $\wh{p}_\ell + \vep \geq 1$, we have $\mscr{M} \li (\wh{p}_\ell, \wh{p}_\ell + \vep \ri
 ) = - \iy < \f{\ln (\ze \de)}{n_\ell}$. In the case of $\wh{p}_\ell + \vep < 1$,
 we have  $- \f{1}{2} <  \wh{p}_\ell - \f{1}{2} + \f{2 \vep}{3} < 1 -
\vep - \f{1}{2} + \f{2 \vep}{3} < \f{1}{2}$ and thus $\f{1}{4} - \li
( \wh{p}_\ell - \f{1}{2} + \f{2 \vep}{3} \ri )^2 > 0$. By virtue of
(\ref{eqDD}),
\[
\mscr{M} \li (\wh{p}_\ell, \wh{p}_\ell + \vep \ri ) = - \f{ \vep^2 }
{ 2 \li [ \f{1}{4} - \li ( \wh{p}_\ell - \f{1}{2} + \f{2 \vep}{3}
\ri )^2 \ri ] }  \leq \f{\ln (\ze \de)}{n_\ell}.
\]
Now, we shall show the second inequality $\mscr{M} \li (\wh{p}_\ell,
\wh{p}_\ell - \vep \ri )  \leq \f{\ln (\ze \de)}{n_\ell}$.  To this
end, we need to establish \be \la{eqDD8} \li ( \wh{p}_\ell -
\f{1}{2}  - \f{2\vep}{3} \ri )^2 \geq \f{1}{4} + \f{ n_\ell \vep^2 }
{2 \ln (\ze \de) } \ee based on (\ref{eqD}).  It is obvious that
(\ref{eqDD8}) holds if {\small $\f{1}{4} + \f{  n_\ell \vep^2 } {2
\ln (\ze \de) } \leq 0$}.  It remains to show (\ref{eqDD8}) under
the condition that {\small $\f{1}{4} + \f{  n_\ell \vep^2 } {2 \ln
(\ze \de) }
> 0$}.  Since (\ref{eqD1})
implies either {\small $\wh{p}_\ell - \f{1}{2} - \f{2\vep}{3}  \leq
- \f{4\vep}{3} - \sq{ \f{1}{4} + \f{  n_\ell \vep^2 } {2 \ln (\ze
\de) } } < - \sq{ \f{1}{4} + \f{  n_\ell \vep^2 } {2 \ln (\ze \de) }
}$} or {\small $\wh{p}_\ell - \f{1}{2} - \f{2\vep}{3}  \geq  \sq{
\f{1}{4} + \f{  n_\ell \vep^2 } {2 \ln (\ze \de) } }$}, it must be
true that (\ref{eqD1}) implies (\ref{eqDD8}). On the other hand,
(\ref{eqD2}) also implies (\ref{eqDD8}) because (\ref{eqD2}) implies
{\small $\wh{p}_\ell - \f{1}{2} - \f{2\vep}{3} \leq - \sq{ \f{1}{4}
+ \f{  n_\ell \vep^2 } {2 \ln (\ze \de) } }$}.  Hence, we have
established (\ref{eqDD8}).

In the case of $\wh{p}_\ell - \vep \leq 0$, we have $\mscr{M} \li
(\wh{p}_\ell, \wh{p}_\ell - \vep \ri ) = - \iy  \leq \f{\ln (\ze
\de)}{n_\ell}$.  In the case of $\wh{p}_\ell - \vep > 0$, we have $-
\f{1}{2} < \vep - \f{1}{2} - \f{2 \vep}{3} < \wh{p}_\ell - \f{1}{2}
- \f{2 \vep}{3} \leq 1 - \f{1}{2} - \f{2 \vep}{3} < \f{1}{2}$ and
thus $\f{1}{4} - \li ( \wh{p}_\ell - \f{1}{2} - \f{2 \vep}{3} \ri
)^2 > 0$. By virtue of (\ref{eqDD8}),
\[
\mscr{M} \li (\wh{p}_\ell, \wh{p}_\ell - \vep \ri ) = - \f{ \vep^2 }
{ 2 \li [ \f{1}{4} - \li ( \wh{p}_\ell - \f{1}{2} - \f{2 \vep}{3}
\ri )^2 \ri ] }  \leq \f{\ln (\ze \de)}{n_\ell}.
\]
Hence, {\small $ \{ \bs{D}_\ell = 1  \}  \leu \{
 \mscr{M} ( \wh{\bs{p}}_\ell, \wh{\bs{p}}_\ell + \vep )
\leq \f{\ln (\ze \de)}{n_\ell}, \; \mscr{M} ( \wh{\bs{p}}_\ell,
\wh{\bs{p}}_\ell - \vep ) \leq \f{\ln (\ze \de)}{n_\ell} \} \leu \{
\mscr{M}_{\mrm{B}} ( \wh{\bs{p}}_\ell, \wh{\bs{p}}_\ell + \vep  )
\leq \f{\ln (\ze \de)}{n_\ell}, \; \mscr{M}_{\mrm{B}} (
\wh{\bs{p}}_\ell, \wh{\bs{p}}_\ell - \vep ) \leq \f{\ln (\ze
\de)}{n_\ell} \} $} for $\ell = 1, \cd, s$. The proof of the lemma
is thus completed.

\epf

\bsk

Now we are in a position to prove Theorem \ref{Bino_ABS_CDF_CH_MA}.

If the stopping rule derived from CDF $\&$ CCDF is used,  then $\{
\bs{D}_s = 1 \}$ is a sure event as a result of Lemma \ref{absDS1}.
Therefore, the sampling scheme satisfies all the requirements
described in Theorem \ref{Monotone_second}, from which Theorem
\ref{Bino_ABS_CDF_CH_MA} immediately follows.

If the stopping rule derived from Chernoff bounds is used, then $\{
\bs{D}_s = 1 \}$ is a sure event as a result of Lemma \ref{DS1}.
Recall that $\exp( \mscr{M}_{\mrm{B}} (z, p) )$  is equal to
$\mcal{F} (z, p)$ and $\mcal{G} (z, p)$ respectively for the cases
of $z \leq p$ and $z \geq p$.  Moreover, $\wh{\bs{p}}_\ell$ is a ULE
of $p$ for $\ell = 1, \cd, s$.  By virtue of these facts and Lemmas
\ref{DS1} and \ref{abs13}, the sampling scheme satisfies all the
requirements described in Corollary \ref{Monotone_third}, from which
Theorem \ref{Bino_ABS_CDF_CH_MA} immediately follows.

If the stopping rule derived from Massart's inequality is used, then
$\{ \bs{D}_s = 1 \}$ is a sure event as a result of Lemma
\ref{absD}.  By virtue of the fact that $\exp( \mscr{M}_{\mrm{B}}
(z, p) )$  is equal to $\mcal{F} (z, p)$ and $\mcal{G} (z, p)$
respectively for the cases of $z \leq p$ and $z \geq p$, the fact
that $\wh{\bs{p}}_\ell$ is a ULE of $p$ for $\ell \in \{ 1, \cd, s
\}$, and Lemmas \ref{absD} and \ref{abs381}, the sampling scheme
satisfies all the requirements described in Corollary
\ref{Monotone_third}, from which Theorem \ref{Bino_ABS_CDF_CH_MA}
immediately follows.

\subsection{Proof of Theorem \ref{Range_Bino_Chernoff} } \la{App_Range_Bino_Chernoff}

Theorem \ref{Range_Bino_Chernoff} can be shown by applying Lemmas
\ref{lem16} and \ref{lem17} to be established in the sequel.

\beL \la{lem16} For $\ell = 1, \cd, s - 1$, {\small \[ \{
\bs{D}_\ell = 0 \} = \li \{ \mscr{M}_{\mrm{B}}(\wh{\bs{p}}_\ell,
\wh{\bs{p}}_\ell + \vep)
> \f{ \ln (\ze \de) } {n_\ell} \ri \} \bigcup \li \{
\mscr{M}_{\mrm{B}}(\wh{\bs{p}}_\ell, \wh{\bs{p}}_\ell - \vep) > \f{
\ln (\ze \de) } {n_\ell} \ri \}.
\]}
\eeL

\bpf  To show the lemma, by the definition of $\bs{D}_\ell$, it
suffices to show

{\small $\li \{ \mscr{M}_{\mrm{B}} \li ( \f{1}{2} - \li |\f{1}{2} -
\wh{\bs{p}}_\ell \ri | , \f{1}{2} - \li |\f{1}{2} - \wh{\bs{p}}_\ell
\ri | + \vep \ri ) \leq \f{ \ln (\ze \de) } {n_\ell} \ri \} = \li \{
\mscr{M}_{\mrm{B}}(\wh{\bs{p}}_\ell, \wh{\bs{p}}_\ell + \vep) \leq
\f{ \ln (\ze \de) } {n_\ell}, \;
\mscr{M}_{\mrm{B}}(\wh{\bs{p}}_\ell, \wh{\bs{p}}_\ell - \vep) \leq
\f{ \ln (\ze \de) } {n_\ell} \ri \}$ } for $\ell = 1, \cd, s - 1$.
For simplicity of notations, we denote $\wh{\bs{p}}_\ell (\om)$ by
$\wh{p}_\ell$ for $\om \in \Om$.  First, we claim that {\small
$\mscr{M}_{\mrm{B}} \li ( \f{1}{2} - \li |\f{1}{2} - \wh{p}_\ell \ri
| , \f{1}{2} - \li |\f{1}{2} - \wh{p}_\ell \ri | + \vep \ri ) \leq
\f{ \ln ( \ze \de) } {n_\ell}$} implies {\small
$\mscr{M}_{\mrm{B}}(\wh{p}_\ell, \wh{p}_\ell + \vep) \leq \f{ \ln (
\ze \de) } {n_\ell}$} and {\small $\mscr{M}_{\mrm{B}}(\wh{p}_\ell,
\wh{p}_\ell - \vep) \leq \f{ \ln ( \ze \de) } {n_\ell}$}.  To prove
this claim, we need to consider two cases: (i) $\wh{p}_\ell \leq
\f{1}{2}$; (ii) $\wh{p}_\ell > \f{1}{2}$.  In the case of
$\wh{p}_\ell \leq \f{1}{2}$, we have {\small
$\mscr{M}_{\mrm{B}}(\wh{p}_\ell, \wh{p}_\ell - \vep) \leq
\mscr{M}_{\mrm{B}}(\wh{p}_\ell, \wh{p}_\ell + \vep) =
\mscr{M}_{\mrm{B}} \li ( \f{1}{2} - \li |\f{1}{2} - \wh{p}_\ell \ri
| , \f{1}{2} - \li |\f{1}{2} - \wh{p}_\ell \ri | + \vep \ri ) \leq
\f{ \ln ( \ze \de) } {n_\ell}$}, where the first inequality follows
from Lemma \ref{lem888}.  Similarly, in the case of $\wh{p}_\ell >
\f{1}{2}$, we have {\small $\mscr{M}_{\mrm{B}}(\wh{p}_\ell,
\wh{p}_\ell + \vep) < \mscr{M}_{\mrm{B}}(\wh{p}_\ell, \wh{p}_\ell -
\vep) = \mscr{M}_{\mrm{B}}(1 - \wh{p}_\ell, 1 - \wh{p}_\ell + \vep)
= \mscr{M}_{\mrm{B}} \li ( \f{1}{2} - \li |\f{1}{2} - \wh{p}_\ell
\ri |, \f{1}{2} - \li |\f{1}{2} - \wh{p}_\ell \ri | + \vep \ri )
\leq \f{ \ln ( \ze \de) } {n_\ell}$}, where the first inequality
follows from Lemma \ref{lem888}.  The claim is thus established.

Second, we claim that {\small $\mscr{M}_{\mrm{B}}(\wh{p}_\ell,
\wh{p}_\ell + \vep) \leq \f{ \ln ( \ze \de) } {n_\ell}$} and {\small
$\mscr{M}_{\mrm{B}}(\wh{p}_\ell, \wh{p}_\ell - \vep) \leq \f{ \ln (
\ze \de) } {n_\ell}$} together imply {\small $\mscr{M}_{\mrm{B}} (
\f{1}{2} -  |\f{1}{2} - \wh{p}_\ell | , \f{1}{2} - |\f{1}{2} -
\wh{p}_\ell | + \vep ) \leq \f{ \ln ( \ze \de) } {n_\ell}$}.  To
prove this claim, we need to consider two cases: (i) $\wh{p}_\ell
\leq \f{1}{2}$; (ii) $\wh{p}_\ell > \f{1}{2}$.  In the case of
$\wh{p}_\ell \leq \f{1}{2}$, we have {\small $\mscr{M}_{\mrm{B}} \li
( \f{1}{2} - \li |\f{1}{2} - \wh{p}_\ell \ri | , \f{1}{2} - \li
|\f{1}{2} - \wh{p}_\ell \ri | + \vep \ri ) =
\mscr{M}_{\mrm{B}}(\wh{p}_\ell, \wh{p}_\ell + \vep) \leq \f{ \ln (
\ze \de) } {n_\ell}$}. Similarly, in the case of $\wh{p}_\ell >
\f{1}{2}$, we have {\small $\mscr{M}_{\mrm{B}} \li ( \f{1}{2} - \li
|\f{1}{2} - \wh{p}_\ell \ri | , \f{1}{2} - \li |\f{1}{2} -
\wh{p}_\ell \ri | + \vep \ri ) = \mscr{M}_{\mrm{B}}(1 - \wh{p}_\ell,
1 - \wh{p}_\ell + \vep) = \mscr{M}_{\mrm{B}}(\wh{p}_\ell,
\wh{p}_\ell - \vep) \leq \f{ \ln ( \ze \de) } {n_\ell}$}. This
establishes our second claim.

Finally, combining our two established claims leads to {\small $ \{
\mscr{M}_{\mrm{B}} ( \f{1}{2} -  |\f{1}{2} - \wh{\bs{p}}_\ell  | ,
\f{1}{2} - |\f{1}{2} - \wh{\bs{p}}_\ell  | + \vep  ) \leq \f{ \ln
(\ze \de) } {n_\ell}  \} =  \{ \mscr{M}_{\mrm{B}}(\wh{\bs{p}}_\ell,
\wh{\bs{p}}_\ell + \vep) \leq \f{ \ln (\ze \de) } {n_\ell}, \;
\mscr{M}_{\mrm{B}}(\wh{\bs{p}}_\ell, \wh{\bs{p}}_\ell - \vep) \leq
\f{ \ln (\ze \de) } {n_\ell}  \}$ }. This completes the proof of the
lemma.

\epf

\beL \la{lem17} For $\ell = 1, \cd, s - 1$, {\small \bee &   & \li
\{ \mscr{M}_{\mrm{B}}(\wh{\bs{p}}_\ell, \wh{\bs{p}}_\ell + \vep) >
\f{ \ln (\ze \de) } {n_\ell} \ri \}
= \{ n_\ell \; \udl{z} < K_\ell < n_\ell \ovl{z} \}, \\
&   & \li \{ \mscr{M}_{\mrm{B}}(\wh{\bs{p}}_\ell, \wh{\bs{p}}_\ell -
\vep)
> \f{ \ln (\ze \de) } {n_\ell} \ri \} = \{ n_\ell (1 -
\ovl{z}) < K_\ell < n_\ell (1 - \udl{z}) \}.  \eee} \eeL

\bpf

Since {\small $\f{\pa \mscr{M}_{\mrm{B}}(z, z + \vep) } {\pa z}
 = \ln \f{(z + \vep) ( 1 - z) }{z (1 - z - \vep) } - \f{ \vep } {
(z + \vep) (1 - z - \vep) }$} for $z \in (0, 1 - \vep)$, it follows
that the partial derivative {\small $\f{ \pa \mscr{M}_{\mrm{B}}(z, z
+ \vep) } {\pa z}$} is equal to $0$ for $z = z^*$. The existence and
uniqueness of $z^*$ can be established by verifying that {\small
$\f{\pa^2 \mscr{M}_{\mrm{B}}(z, z + \vep) } {\pa z^2} = - \vep^2 \li
[ \f{1}{z (z+\vep)^2 }  + \f{1}{ (1 - z) ( 1 - z - \vep)^2 } \ri ] <
0$} for any $z \in (0, 1 - \vep)$
 and that
{\small  \[
 \li.
\f{ \pa \mscr{M}_{\mrm{B}}(z, z + \vep) } {\pa z}  \ri |_{z =
\f{1}{2}} = \ln \f{ 1 + 2 \vep } { 1 - 2 \vep } - \f{ \vep }
{\f{1}{4} - \vep^2 } < 0, \qqu \li. \f{ \pa \mscr{M}_{\mrm{B}}(z, z
+ \vep) } {\pa z}  \ri |_{z = \f{1}{2} - \vep} = \ln \f{ 1 + 2 \vep
} { 1 - 2 \vep } - 4 \vep > 0.
\]}
Since $\mscr{M}_{\mrm{B}}(z^*, z^* + \vep)$ is negative and {\small
$n_\ell < \f{ \ln (\ze \de) } { \mscr{M}_{\mrm{B}}(z^*, z^* + \vep)
}$}, we have that {\small $\mscr{M}_{\mrm{B}}(z^*, z^* + \vep)
> \f{ \ln (\ze \de) } {n_\ell}$}.  On the other hand, by the definition of sample
sizes, we have {\small $n_\ell \geq n_1 = \li \lc \f{\ln (\ze \de)}{
\ln (1 - \vep) } \ri \rc \geq \f{ \ln (\ze \de) } { \lim_{z \to 0}
\mscr{M}_{\mrm{B}}(z, z + \vep) }$}, which implies {\small $\lim_{z
\to 0} \mscr{M}_{\mrm{B}}(z, z + \vep) \leq \f{ \ln (\ze \de) }
{n_\ell}$}. Noting that $\mscr{M}_{\mrm{B}}(z, z + \vep)$ is
monotonically increasing with respect to $z \in (0, z^*)$, we can
conclude from the intermediate value theorem that there exists a
unique number $\udl{z} \in [0, z^*)$ such that {\small
$\mscr{M}_{\mrm{B}}(\udl{z}, \udl{z} + \vep) = \f{ \ln (\ze \de) }
{n_\ell}$}. Similarly, due to the facts that {\small
$\mscr{M}_{\mrm{B}}(z^*, z^* + \vep)
> \f{ \ln (\ze \de) } {n_\ell}, \; \lim_{z \to 1 - \vep}
\mscr{M}_{\mrm{B}}(z, z + \vep) = - \iy <  \f{ \ln (\ze \de) }
{n_\ell}$} and that $\mscr{M}_{\mrm{B}}(z, z + \vep)$ is
monotonically decreasing with respect to $z \in (z^*, 1 - \vep)$, we
can conclude from the intermediate value theorem that there exists a
unique number $\ovl{z} \in (z^*, 1 - \vep)$ such that {\small
$\mscr{M}_{\mrm{B}}(\ovl{z}, \ovl{z} + \vep) = \f{ \ln (\ze \de) }
{n_\ell}$}.  Therefore, we have $\mscr{M}_{\mrm{B}}(z, z + \vep) >
\f{ \ln (\ze \de) } {n_\ell}$ for $z \in (\udl{z}, \ovl{z} )$,  and
$\mscr{M}_{\mrm{B}}(z, z + \vep) \leq  \f{ \ln (\ze \de) } {n_\ell}$
for $z \in [0, \udl{z}] \cup [\ovl{z}, 1]$. This proves that {\small
$ \{ \mscr{M}_{\mrm{B}}(\wh{\bs{p}}_\ell, \wh{\bs{p}}_\ell + \vep) >
\f{ \ln (\ze \de) } {n_\ell} \} = \{ n_\ell \; \udl{z} < K_\ell <
n_\ell \ovl{z} \}$}.  Noting that $\mscr{M}_{\mrm{B}} \li ( \f{1}{2}
+ \upsilon, \f{1}{2} + \upsilon - \vep \ri ) = \mscr{M}_{\mrm{B}}
\li ( \f{1}{2} - \upsilon, \f{1}{2} - \upsilon + \vep \ri )$ for any
$\upsilon \in \li ( 0, \f{1}{2} \ri )$, we have {\small $ \{
\mscr{M}_{\mrm{B}}(\wh{\bs{p}}_\ell, \wh{\bs{p}}_\ell - \vep) > \f{
\ln (\ze \de) } {n_\ell} \} = \{ n_\ell (1 - \ovl{z}) < K_\ell <
n_\ell (1 - \udl{z}) \}$.}  This completes the proof of the lemma.

\epf

\subsection{Proof of Theorem  \ref{them18889} } \la{them18889_app}

The first statement can be shown by virtue of Chernoff-Hoeffding inequality.  To prove (\ref{good9}), let $\ep = \f{ 3 w \vep }{ 2 }$ and $\xi =
\f{1}{\de} \exp \li ( \f{ 9 w^2 }{4} \ln (\ze \de) \ri )$.  Then the sampling scheme can be restated as follows:

Continue sampling until \[
 \li ( \li | \wh{\bs{p}}_\ell - \f{1}{2} \ri | - \f{2}{3} \ep \ri )^2 \geq \f{1}{4} + \f{ \ep^2 n_\ell } { 2 \ln (\xi \de)
}
\]
for some $\ell \in \{1, 2, \cd, s \}$.

By the definition of $\ep$ and $\xi$, we have $\ep \leq \vep$ and
\[
\f{ \ln \f{1}{\xi \de} }{2 \ep^2}  = \f{ -  \f{ 9 w^2 }{4} \ln (\ze \de)  } { 2  \times \f{ 9 w^2 }{4} \vep^2} = \f{ \ln \f{1}{\ze \de} }{2
\vep^2} \leq n_s, \qqu \xi \de = (\ze \de)^{9 w^2 \sh 4}.
\]
Therefore, by Theorem \ref{Bino_ABS_CDF_CH_MA}, we have
\[
\Pr \li \{ \li | \bs{\wh{p}} - p \ri | < \vep \mid p \ri \} \geq \Pr \li \{ \li | \bs{\wh{p}} - p \ri | < \ep \mid p \ri \} \geq 1 - 2 s \xi \de
= 1 - 2 s (\ze \de)^{9 w^2 \sh 4}
\]
for any $p \in (0, 1)$.

\subsection{Proof of Theorem \ref{Bino_DDV_Asp} } \la{App_Bino_DDV_Asp}

We need some preliminary results.

 \beL
 \la{lem31a}
Let $c$ be a positive number.   Let $\ka(\ell, \ep)$ be a bivariate
function of positive number $\ep$ and integer $\ell$.  Let $r$ be a
positive integer dependent on $\ep$.  Suppose that, for  $\ep
> 0$ small enough, $\ka(\ell + 1, \ep)- \ka(\ell, \ep) \geq 1$ for
any $\ell > 0$. Suppose that $\ka(1, \ep)$ tends to infinity as
$\ep$ tends to $0$. Then, $\lim_{\ep \to 0} \sum_{\ell = 1}^r
\ka(\ell, \ep) \; e^{- c \ka(\ell, \ep)} = 0$. \eeL

\bpf

Choose an $\ep > 0$ small enough such that $c \ka(\ell, \ep)
> 1$ for all $\ell \geq 1$.  Since $x e^{-x}$ is monotonically
decreasing with respect to $x > 1$, we have \bee \sum_{\ell = 1}^r
\ka(\ell, \ep) \; e^{- c \ka(\ell, \ep)} & \leq &  \f{1}{c}
\sum_{\ell = 1}^r \lf c \ka(\ell, \ep) \rf \; e^{- \lf c \ka(\ell,
\ep) \rf} \leq \f{1}{c} \sum_{m = \lf c \ka(1, \ep) \rf }^\iy m
e^{-m} \\
& < & \f{1}{c} \int_{ \lf c \ka(1, \ep) \rf - 1 }^\iy x e^{-x} dx =
\f{\lf c \ka(1, \ep) \rf}{c}  e^{ 1 - \lf c \ka(1, \ep) \rf} \to 0
\eee as $\ep \to 0$.  This completes the proof of the lemma.

\epf

\beL

\la{lem32T}

Let $\psi_\ep$ be a function of $\ep \in (0, 1)$ such that $0 < a
\leq \psi_\ep \leq b < 1$. Then, \bee  &  & \mscr{M}_{\mrm{B}}
(\psi_\ep, \psi_\ep + \ep) = - \f{ \ep^2 } { 2 \psi_\ep ( 1 -
\psi_\ep) } + \f{ \ep^3 } { 3 }  \f{ 1 - 2 \psi_\ep } { \psi_\ep^2
(1 - \psi_\ep)^2 } +  o (\ep^3),\\
&  & \mscr{M}_{\mrm{I}} \li ( \psi_\ep, \f{\psi_\ep}{1 + \ep} \ri )
= - \f{ \ep^2 } { 2 ( 1 - \psi_\ep) }  + \f{\ep^3}{3} \f{ 2 -
\psi_\ep } { (1 - \psi_\ep)^2 } + o(\ep^3),\\
&  & \mscr{M}_{\mrm{B}} \li ( \psi_\ep, \f{\psi_\ep}{1 + \ep} \ri )
= - \f{ \ep^2 \psi_\ep } { 2 ( 1 - \psi_\ep) } + \f{ \ep^3 \psi_\ep
(2 - \psi_\ep) } { 3 (1 - \psi_\ep)^2 } + o (\ep^3). \eee  \eeL

\bpf  Using Taylor's series expansion formula $\ln (1 + x)  = x -
\f{x^2}{2} + \f{x^3}{3} + o(x^3)$ for $|x| < 1$,  we have {\small
\bee \mscr{M}_{\mrm{B}} (\psi_\ep, \psi_\ep + \ep) & = & \psi_\ep
\ln \li ( 1 +
\f{\ep}{\psi_\ep} \ri ) + (1 - \psi_\ep) \ln \li ( 1 - \f{\ep}{1 - \psi_\ep} \ri )\\
& = &  - \f{ \ep^2 } { 2 \psi_\ep ( 1 - \psi_\ep) } +
\f{\psi_\ep}{3} \li ( \f{\ep}{\psi_\ep} \ri )^3 +  \f{1 -
\psi_\ep}{3} \li (  - \f{\ep}{1 - \psi_\ep} \ri )^3 \\
&   & + \psi_\ep \times o \li ( \f{\ep^3}{\psi_\ep^3}  \ri )
+ (1 - \psi_\ep) \times o \li ( - \f{\ep^3}{(1 - \psi_\ep)^3}  \ri )\\
& = &  - \f{ \ep^2 } { 2 \psi_\ep ( 1 - \psi_\ep) } + \f{ \ep^3 } {
3 } \f{ 1 - 2 \psi_\ep } { \psi_\ep^2 (1 - \psi_\ep)^2 } +  o
(\ep^3) \eee} for $\ep < \psi_\ep < 1 - \ep$. Since $\lim_{\ep \to
0} \f{\ep}{1 + \ep} \f{\psi_\ep}{1 -
 \psi_\ep} = 0$ and
 \[
\lim_{\ep \to 0} \f{ \f{1 - \psi_\ep}{\psi_\ep} \times o \li ( \li (
\f{\ep}{1 + \ep} \f{\psi_\ep}{1 - \psi_\ep} \ri )^3  \ri)  }{\ep^3}
= \lim_{\ep \to 0} \f{ \f{1 - \psi_\ep}{\psi_\ep} \times o \li ( \li
( \f{\ep}{1 + \ep} \f{\psi_\ep}{1 - \psi_\ep} \ri )^3 \ri) }{\li (
\f{\ep}{1 + \ep} \f{\psi_\ep}{1 - \psi_\ep} \ri )^3 } \f{ \li (
\f{\ep}{1 + \ep} \f{\psi_\ep}{1 - \psi_\ep} \ri )^3 }{\ep^3} = 0,
 \]
we have {\small \bee \mscr{M}_{\mrm{I}} \li ( \psi_\ep,
\f{\psi_\ep}{1 + \ep} \ri ) & = & - \ln \li ( 1 + \ep
\ri ) + \f{1 - \psi_\ep}{\psi_\ep} \ln \li (  1 + \f{\ep}{1 + \ep} \f{\psi_\ep}{1 - \psi_\ep} \ri )\\
& = &  - \ep + \f{\ep^2}{2} - \f{\ep^3}{ 3 } + \f{1 -
\psi_\ep}{\psi_\ep} \li [ \f{\ep}{1 + \ep} \f{\psi_\ep}{1 -
\psi_\ep} - \f{1}{2} \li ( \f{\ep}{1 + \ep} \f{\psi_\ep}{1 -
\psi_\ep} \ri )^2 + \f{1}{3} \li ( \f{\ep}{1
+ \ep} \f{\psi_\ep}{1 - \psi_\ep} \ri )^3  \ri ]\\
&  & + o (\ep^3) + \f{1 - \psi_\ep}{\psi_\ep} \times o \li ( \li (
\f{\ep}{1
+ \ep} \f{\psi_\ep}{1 - \psi_\ep} \ri )^3  \ri) \\
& = &  \f{\ep^2}{2} - \f{\ep^3}{ 3} - \f{\ep^2}{1 + \ep}  - \f{1}{2}
\li ( \f{\ep}{1 + \ep} \ri )^2 \f{\psi_\ep}{1 - \psi_\ep}  +
\f{1}{3}
\f{\ep^3}{(1 + \ep)^3} \f{\psi_\ep^2}{(1 - \psi_\ep)^2} + o (\ep^3)\\
& = & - \f{\ep^2}{2 (1 - \psi_\ep)} + \f{2 \ep^3}{ 3} + \f{\ep^3
\psi_\ep}{1 -
\psi_\ep} + \f{1}{3} \f{\ep^3 \psi_\ep^2}{(1 - \psi_\ep)^2} + o (\ep^3)\\
& = & - \f{ \ep^2 } { 2 ( 1 - \psi_\ep) }  + \f{\ep^3}{3} \f{ 2 -
\psi_\ep } { (1 - \psi_\ep)^2 } + o(\ep^3). \eee}  Since $\psi_\ep$
is bounded in $[a, b]$, we have
\[ \mscr{M}_{\mrm{B}} \li ( \psi_\ep, \f{\psi_\ep}{1 + \ep} \ri ) = \psi_\ep \mscr{M}_{\mrm{I}}
\li ( \psi_\ep, \f{\psi_\ep}{1 + \ep} \ri ) = - \f{ \ep^2 \psi_\ep }
{ 2 ( 1 - \psi_\ep) } + \f{ \ep^3 \psi_\ep (2 - \psi_\ep) } { 3 (1 -
\psi_\ep)^2 } + o (\ep^3).  \]

\epf

\beL \la{lemm21} Let $0 < \vep < \f{1}{2}$.  Then,
there exists a unique number $z^\star \in (\f{1}{2}, \f{1}{2} + \vep)$
such that $\mscr{M}_{\mrm{B}} (z, z - \vep)$ is monotonically increasing with
respect to $z \in (\vep, z^\star)$ and monotonically decreasing
with respect to $z \in (z^\star, 1)$. Similarly, there exists a unique number $z^* \in (\f{1}{2} - \vep, \f{1}{2})$
such that $\mscr{M}_{\mrm{B}} (z, z + \vep)$ is monotonically increasing with
respect to $z \in (0, z^*)$ and monotonically decreasing with respect to $z \in (z^*, 1 - \vep)$.

\eeL

\bpf

Note that {\small $\li. \f{ \pa \mscr{M}_{\mrm{B}}(z, z - \vep) }
{\pa z} \ri |_{z = \f{1}{2}} = \ln \f{ 1 - 2 \vep } { 1 + 2 \vep } +
\f{ \vep } {\f{1}{4} - \vep^2 } > 0$}  because $\ln \f{ 1 - 2 \vep }
{ 1 + 2 \vep } + \f{ \vep } {\f{1}{4} - \vep^2 }$ equals $0$ for
$\vep = 0$ and its derivative with respect to $\vep$ equals to $\f{
2\vep^2 } { (\f{1}{4} - \vep^2 )^2 }$ which is positive for any
positive $\vep$ less than $\f{1}{2}$.   Similarly, {\small $\li. \f{
\pa \mscr{M}_{\mrm{B}}(z, z - \vep) } {\pa z} \ri |_{z = \f{1}{2} +
\vep} = \ln \f{ 1 - 2 \vep } { 1 + 2 \vep } + 4 \vep < 0$} because
$\ln \f{ 1 - 2 \vep } { 1 + 2 \vep } + 4 \vep$ equals $0$ for $\vep
= 0$ and its derivative with respect to $\vep$ equals to $- \f{ 16
\vep^2 } { 1 - 4 \vep^2 } $ which is negative for any positive
$\vep$ less than $\f{1}{2}$.   In view of the signs of {\small $\f{
\pa \mscr{M}_{\mrm{B}}(z, z - \vep) } {\pa z}$}  at $\f{1}{2},
\f{1}{2} + \vep$ and the fact that {\small $\f{\pa^2
\mscr{M}_{\mrm{B}}(z, z - \vep) } { \pa z^2 } = - \vep^2 \li [
\f{1}{z (z - \vep)^2 }  + \f{1}{ (1 - z) ( 1 - z + \vep)^2  } \ri ]
< 0$} for any $z \in (\vep,1)$, we can conclude from the intermediate value theorem that
there exists a unique number $z^\star \in (\f{1}{2}, \f{1}{2} + \vep)$
such that {\small $\li. \f{ \pa \mscr{M}_{\mrm{B}}(z, z - \vep) }
{\pa z} \ri |_{z = z^\star} = 0$}, which implies that
$\mscr{M}_{\mrm{B}} (z, z - \vep)$ is monotonically increasing with
respect to $z \in (\vep, z^\star)$ and monotonically decreasing
with respect to $z \in (z^\star, 1)$.

To show the second statement of the lemma, note that {\small $\li.
\f{ \pa \mscr{M}_{\mrm{B}}(z, z + \vep) } {\pa z}  \ri |_{z =
\f{1}{2}} = \ln \f{ 1 + 2 \vep } { 1 - 2 \vep } - \f{ \vep }
{\f{1}{4} - \vep^2 } < 0$} because $\ln \f{ 1 + 2 \vep } { 1 - 2
\vep } - \f{ \vep } {\f{1}{4} - \vep^2 }$ equals $0$ for $\vep = 0$
and its derivative with respect to $\vep$ equals to $- \f{ 2\vep^2 }
{ (\f{1}{4} - \vep^2 )^2 }$ which is negative for any positive
$\vep$ less than $\f{1}{2}$.  Similarly, {\small $\li. \f{ \pa
\mscr{M}_{\mrm{B}}(z, z + \vep) } {\pa z}  \ri |_{z = \f{1}{2} -
\vep} = \ln \f{ 1 + 2 \vep } { 1 - 2 \vep } - 4 \vep > 0$} because
$\ln \f{ 1 + 2 \vep } { 1 - 2 \vep } - 4 \vep$ equals $0$ for $\vep
= 0$ and its derivative with respect to $\vep$ equals to $\f{ 16
\vep^2 } { 1 - 4 \vep^2 }$ which is positive for any positive $\vep$
less than $\f{1}{2}$.  In view of the signs of {\small $\f{ \pa
\mscr{M}_{\mrm{B}}(z, z + \vep) } {\pa z}$}  at $\f{1}{2} - \vep,
\f{1}{2}$ and the fact that {\small $\f{\pa^2  \mscr{M}_{\mrm{B}}(z,
z + \vep) } {\pa z^2}  =  - \vep^2 \li [ \f{1}{z (z+\vep)^2 }  +
\f{1}{ (1 - z) ( 1 - z - \vep)^2 } \ri ] < 0$}
 for any $z \in (0, 1 - \vep)$, we can conclude from the intermediate value theorem that
there exists a unique number $z^* \in (\f{1}{2} - \vep, \f{1}{2})$
such that {\small $\li. \f{ \pa \mscr{M}_{\mrm{B}}(z, z + \vep) }
{\pa z} \ri |_{z = z^*} = 0$},  which implies that $\mscr{M}_{\mrm{B}} (z, z + \vep)$
is monotonically increasing with respect to $z \in (0, z^*)$ and monotonically decreasing with respect to $z \in
(z^*, 1 - \vep)$.  This completes the proof of the lemma.

\epf

\beL \la{lem33a} If $\vep$ is sufficiently small, then the following
statements hold true.

(I): For $\ell = 1, 2, \cd, s - 1$, there exists a unique number
$z_{\ell} \in [0, \f{1}{2} - \vep)$ such that $n_{\ell} = \f{ \ln (
\ze \de ) } { \mscr{M}_{\mrm{B}} ( z_{\ell}, \; z_{\ell} + \vep )
}$.

(II): $z_{\ell}$ is monotonically increasing with respect to $\ell$
smaller than $s$.

(III): $\lim_{\vep \to 0} z_{\ell} = \f{ 1 - \sq{ 1 - C_{s - \ell}}
}{2}$, where the limit is taken under the restriction that $s -
\ell$ is fixed with respect to $\vep$.

(IV) For $p \in (0, \f{1}{2})$ such that $C_{j_p} = 4 p (1 - p)$ and
$j_p \geq 1$,
\[ \lim_{\vep \to 0} \f{ z_{\ell_\vep}  - p } { \vep } = - \f{ 2 } {
3 },
\]
where $\ell_\vep = s - j_p$.

(V): $\{ \bs{D}_{\ell} = 0 \} = \{ z_{\ell} < \wh{\bs{p}}_{\ell} < 1
- z_{\ell} \}$ for $\ell = 1, 2, \cd, s - 1$.

\eeL

\bsk

{\bf Proof of Statement (I)}:  By the definition of sample sizes, we
have \be \la{from} 0 < \f{ \ln ( \ze \de )  } { \mscr{M}_{\mrm{B}} (
0, \vep ) }  \leq n_{\ell}  <  \f{(1 + C_1) n_s}{2}  < \f{1 +
C_1}{2} \li (  \f{ \ln \f{1}{ \ze \de } } { 2 \vep^2 }  + 1 \ri )
 \ee
for sufficiently small $\vep > 0$. By (\ref{from}), we have {\small
$\f{ \ln ( \ze \de ) } { n_{\ell} } \geq \mscr{M}_{\mrm{B}} ( 0,
\vep )$} and
\[
\f{ \ln ( \ze \de ) } { n_{\ell} }  < - 2 \vep^2 \li ( \f{2}{1 +
C_1} - \f{1}{ n_{\ell} } \ri ) = \f{ - 2 \vep^2} {
\mscr{M}_{\mrm{B}} ( \f{1}{2} - \vep, \f{1}{2} ) } \f{2}{1 + C_1}
\mscr{M}_{\mrm{B}} \li ( \f{1}{2} - \vep, \f{1}{2} \ri  ) + \f{2
\vep^2}{ n_{\ell} } .
\]
Noting that $\lim_{\vep \to 0} \f{2 \vep^2}{ n_{\ell} } = 0$ and
$\lim_{\vep \to 0}  \f{ - 2 \vep^2} { \mscr{M}_{\mrm{B}} ( \f{1}{2}
- \vep, \f{1}{2} ) } = 1$, we have {\small $\f{ \ln ( \ze \de ) } {
n_{\ell} } < \mscr{M}_{\mrm{B}} \li ( \f{1}{2} - \vep, \f{1}{2} \ri
) < 0$} for sufficiently small $\vep > 0$. In view of the established
fact that $ \mscr{M}_{\mrm{B}} ( 0, \vep ) \leq \f{ \ln ( \ze \de )
} { n_{\ell} } < \mscr{M}_{\mrm{B}} \li ( \f{1}{2} - \vep, \f{1}{2}
\ri )$ for small enough $\vep > 0$ and the fact that
$\mscr{M}_{\mrm{B}} ( z, z + \vep )$ is monotonically increasing
with respect to $z \in (0, \f{1}{2} - \vep)$ as asserted by Lemma \ref{lemm21},
invoking the intermediate value theorem, we have that there exists a
unique number $z_{\ell} \in [0, \f{1}{2} - \vep)$ such that
$\mscr{M}_{\mrm{B}} ( z_{\ell}, z_{\ell} + \vep )  = \f{ \ln ( \ze
\de ) } { n_{\ell} }$.  This proves Statement (I).

\bsk

{\bf Proof of Statement (II)}:  Since $n_{\ell}$ is monotonically
increasing with respect to $\ell$ for sufficiently small $\vep > 0$,
we have that $\mscr{M}_{\mrm{B}} ( z_{\ell}, z_{\ell} + \vep )$ is
monotonically increasing with respect to $\ell$ if $\vep > 0$ is
sufficiently small . Recalling that $\mscr{M}_{\mrm{B}} ( z, z +
\vep )$ is monotonically increasing with respect to $z \in (0,
\f{1}{2} - \vep)$, we have that $z_{\ell}$ is monotonically
increasing with respect to $\ell$. This establishes Statement (II).

\bsk

{\bf Proof of Statement (III)}:   For simplicity of notations, let
{\small $b_{\ell} = \f{ 1 - \sq{ 1 - C_{s - \ell} } }{2}$} for $\ell
= 1, 2, \cd, s - 1$. Then, it can be checked that $4b_{\ell} (1 -
b_{\ell}) = C_{s - \ell}$ and, by the definition of sample sizes, we
have \be \la{goode} \f{\mscr{M}_{\mrm{B}} (z_{\ell}, z_{\ell} +
\vep) } { \vep^2 \sh [ 2 b_{\ell} (b_{\ell} - 1)] } = \f{1}{n_{\ell}
} \times  \f{C_{s - \ell}} { 2 \vep^2}  \ln \f{1}{\ze \de} = 1 +
o(1) \ee for $\ell = 1, 2, \cd, s - 1$.

We claim that $\se < z_{\ell} < \f{1}{2}$
for $\se \in (0, b_{\ell})$ if $\vep > 0$ is small enough.  To prove this claim, we use a contradiction method.
Suppose the claim is not true, then there
exists a set, denoted by $S_\vep$,  of infinite many values of $\vep$ such that $z_{\ell}
\leq \se$ for $\vep \in S_\vep$. For small enough $\vep \in S_\vep$,  we have
$z_{\ell} + \vep \leq \se + \vep < b_\ell + \vep < \f{1}{2}$.  Hence, by (\ref{goode}) and the fact that
$\mscr{M}_{\mrm{B}} (z, z + \vep)$ is
monotonically increasing with respect to
$z \in (0, \f{1}{2} - \vep)$ as asserted by Lemma \ref{lemm21}, we have
\[
1 + o(1) = \f{\mscr{M}_{\mrm{B}} (z_{\ell}, z_{\ell} + \vep) } {
\vep^2 \sh [ 2 b_{\ell} (b_{\ell} - 1)] } \geq \f{
\mscr{M}_{\mrm{B}} (\se, \se + \vep) } { \vep^2 \sh [ 2 b_{\ell} (b_{\ell} - 1 )] }
= \f{ \vep^2 \sh [ 2 \se (1 - \se)] + o (\vep^2) } {
\vep^2 \sh [ 2 b_{\ell} (1 - b_{\ell})]  } = \f{ b_{\ell} (1 -
b_{\ell} ) } { \se (1 - \se) } + o(1)
\]
for small enough $\vep \in S_\vep$,  which implies {\small $\f{
b_{\ell} (1 - b_{\ell} ) } { \se (1 - \se) } \leq 1$}, contradicting
to the fact that {\small $\f{ b_{\ell} (1 - b_{\ell} ) } { \se (1 -
\se) } > 1$}.  By (\ref{goode}) and applying Lemma \ref{lem32T} based on the
established condition that $\se < z_{\ell} < \f{1}{2}$ for small
enough $\vep > 0$, we have {\small $\f{\mscr{M}_{\mrm{B}} (z_{\ell}, z_{\ell} + \vep) } { \vep^2 \sh
[ 2 b_{\ell} (b_{\ell} - 1)] } = \f{ \vep^2 \sh [ 2 z_{\ell} (1 -
z_{\ell})]  + o (\vep^2) } { \vep^2 \sh [ 2 b_{\ell} (1 - b_{\ell}
)] } = 1 + o (1)$}, which implies {\small $\f{1}{z_{\ell} (1 - z_{\ell})} -
\f{1} {b_{\ell} (1 - b_{\ell} )} = o(1)$} and consequently $\lim_{\vep
\to 0} z_{\ell} = b_{\ell}$.  This proves Statement (III).

\bsk

{\bf Proof of Statement (IV)}:

Since {\small $n_{\ell_\vep} = \li \lc \f{ C_{s - {\ell_\vep}} \ln
\f{1}{\ze \de} } { 2 \vep^2}  \ri \rc$} and $C_{s - {\ell_\vep}} = 4
p ( 1 - p)$, we can write
\[
n_{\ell_\vep} =  \li \lc  \f{ 2 p ( 1 - p) \ln \f{1}{\ze \de} } {
\vep^2} \ri \rc = \f{ \ln (\ze \de) } { \mscr{M}_{\mrm{B}}
(z_{\ell_\vep}, z_{\ell_\vep} + \vep)  },
\]
from which we have $\f{1}{n_{\ell_\vep}} = o (\vep)$,
\[
1 - o(\vep) = 1 - \f{1}{n_{\ell_\vep}}  < \f{ - \f{ 2 p ( 1 - p) \ln
(\ze \de) } { \vep^2} } { \f{ \ln (\ze \de) } { \mscr{M}_{\mrm{B}}
(z_{\ell_\vep}, z_{\ell_\vep} + \vep)  }  }  \leq 1
\]
and thus \be \la{goodep} \f{ - \f{ 2 p ( 1 - p) \ln (\ze \de) } {
\vep^2} } { \f{ \ln (\ze \de) } { \mscr{M}_{\mrm{B}} (z_{\ell_\vep},
z_{\ell_\vep} + \vep) } } = \f{-  \mscr{M}_{\mrm{B}} (z_{\ell_\vep},
z_{\ell_\vep} + \vep) } { \vep^2 \sh [ 2 p (1 - p)] } =  1 +
o(\vep). \ee

For $\se \in (0,  p)$, we claim that $\se < z_{\ell_\vep} <
\f{1}{2}$ provided that $\vep$ is sufficiently small. Suppose, to
get a contradiction,  that the claim is not true. Then,  there
exists a set of infinite many values of $\vep$ such that
$z_{\ell_\vep} \leq \se$ if $\vep$ in the set is small enough. For
such $\vep < \f{1}{2} - p$, by (\ref{goodep}) and the monotonicity
of $\mscr{M}_{\mrm{B}} (z, z + \vep)$ with respect to $z$, we have
\[
1 + o(\vep) = \f{-  \mscr{M}_{\mrm{B}} (z_{\ell_\vep}, z_{\ell_\vep}
+ \vep) } { \vep^2 \sh [ 2 p (1 - p)] } \geq \f{-
\mscr{M}_{\mrm{B}} (\se, \se + \vep) } { \vep^2 \sh [ 2 p (1 - p)] }
= \f{ \vep^2 \sh [ 2 \se (1 - \se)] + o (\vep^2) } { \vep^2 \sh [ 2
p (1 - p)]  } = \f{ p (1 - p) } { \se (1 - \se) } + o(1)
\]
for small enough $\vep$ in the set,  which contradicts to the fact
that {\small $\f{ p (1 - p) } { \se (1 - \se) } > 1$}.  This proves
our claim.  Since $\se < z_{\ell_\vep}  < \f{1}{2}$ is established,
by (\ref{goodep}) and Lemma \ref{lem32T}, we have
\[
\f{-  \mscr{M}_{\mrm{B}} (z_{\ell_\vep}, z_{\ell_\vep} + \vep) } {
\vep^2 \sh [ 2 p (1 - p)] } = \f{ \vep^2 \sh [ 2 z_{\ell_\vep} (1 -
z_{\ell_\vep})]  - \vep^3 (1 - 2 z_{\ell_\vep}) \sh [3
z_{\ell_\vep}^2 (1 - z_{\ell_\vep})^2 ] + o (\vep^3) } { \vep^2 \sh
[ 2 p (1 - p)] } = 1 + o (\vep)
\]
and consequently, \be \la{good1p} \f{1}{z_{\ell_\vep} (1 -
z_{\ell_\vep})} - \f{1} {p (1 - p)} - \f{2 \vep (1 - 2
z_{\ell_\vep}) }{3 z_{\ell_\vep}^2 (1 - z_{\ell_\vep})^2} + o (\vep)
= 0.  \ee Since $\se < z_{\ell_\vep} < \f{1}{2}$ for small enough
$\vep
> 0$, by (\ref{good1p}), we have {\small $\f{1}{z_{\ell_\vep} (1 - z_{\ell_\vep})}
- \f{1}{p (1 - p)} = o(1)$}, from which it follows that $\lim_{\vep
\to 0} z_{\ell_\vep} = p$. Noting that (\ref{good1p}) can be written
as
\[
\f{ (z_{\ell_\vep} - p) ( z_{\ell_\vep} +  p - 1  )  } { p (1 - p)
z_{\ell_\vep} (1 - z_{\ell_\vep}) } - \f{2 \vep (1 - 2
z_{\ell_\vep})}{3 z_{\ell_\vep}^2 (1 - z_{\ell_\vep})^2} + o (\vep)
= 0
\]
and using the fact that $\lim_{\vep \to 0} z_{\ell_\vep} = p \in (0,
\f{1}{2})$, we have
\[
\f{ z_{\ell_\vep} - p } { \vep } = \f{2 p (1 - p) (1 - 2
z_{\ell_\vep})} { 3 ( z_{\ell_\vep} +  p - 1  ) z_{\ell_\vep} (1 -
z_{\ell_\vep}) } + o(1)
\]
for small enough $\vep > 0$, which implies that $\lim_{\vep \to 0}
\f{ z_{\ell_\vep} - p } { \vep } = - \f{ 2 } { 3 }$.  This proves
Statement (IV).

\bsk

{\bf Proof of Statement (V)}:  Note that {\small \bee  \{
\bs{D}_{\ell} = 0 \} & = &  \li \{ \mscr{M}_{\mrm{B}}  \li (
\f{1}{2} - \li | \f{1}{2} - \wh{\bs{p}}_{\ell} \ri |,  \f{1}{2} -
\li | \f{1}{2} - \wh{\bs{p}}_{\ell} \ri | + \vep \ri ) > \f{\ln (\ze
\de)} {n_{\ell}}
, \; \wh{\bs{p}}_{\ell} \leq \f{1}{2} \ri \}\\
&   & \bigcup \li \{ \mscr{M}_{\mrm{B}}  \li ( \f{1}{2} - \li |
\f{1}{2} - \wh{\bs{p}}_{\ell} \ri |,  \f{1}{2} - \li | \f{1}{2} -
\wh{\bs{p}}_{\ell} \ri | + \vep \ri ) > \f{\ln (\ze \de)}
{n_{\ell}},
\; \wh{\bs{p}}_{\ell} > \f{1}{2}  \ri \}\\
& = &  \li \{ \mscr{M}_{\mrm{B}}  \li ( \wh{\bs{p}}_{\ell}, \;
\wh{\bs{p}}_{\ell}  + \vep \ri ) > \f{\ln (\ze \de)} {n_{\ell}} , \;
\wh{\bs{p}}_{\ell} \leq \f{1}{2} \ri \} \bigcup \li \{
\mscr{M}_{\mrm{B}} \li ( \wh{\bs{p}}_{\ell}, \; \wh{\bs{p}}_{\ell} -
\vep \ri )  > \f{\ln (\ze \de)} {n_{\ell}}, \; \wh{\bs{p}}_{\ell}
> \f{1}{2}  \ri \},  \eee} where we have used the fact that
$\mscr{M}_{\mrm{B}} (z, z + \vep) = \mscr{M}_{\mrm{B}} (1 - z, 1 - z
- \vep)$.   We claim that \bel &  & \li \{ \mscr{M}_{\mrm{B}} \li (
\wh{\bs{p}}_{\ell}, \; \wh{\bs{p}}_{\ell}  + \vep \ri ) > \f{\ln
(\ze \de)} {n_{\ell}} , \; \wh{\bs{p}}_{\ell} \leq \f{1}{2} \ri \} =
 \li \{ z_{\ell} < \wh{\bs{p}}_{\ell} \leq \f{1}{2} \ri \}, \la{claimabs1}\\
&  &  \li \{ \mscr{M}_{\mrm{B}}  \li ( \wh{\bs{p}}_{\ell}, \;
\wh{\bs{p}}_{\ell}  - \vep \ri ) > \f{\ln (\ze \de)} {n_{\ell}} , \;
\wh{\bs{p}}_{\ell} > \f{1}{2} \ri \} = \li \{ \f{1}{2} <
\wh{\bs{p}}_{\ell}  < 1 - z_{\ell} \ri \} \la{claimabs2} \eel for
small enough $\vep > 0$.

To prove (\ref{claimabs1}), let $\om \in \{ \mscr{M}_{\mrm{B}} \li (
\wh{\bs{p}}_{\ell}, \; \wh{\bs{p}}_{\ell}  + \vep \ri )  > \f{\ln
(\ze \de)} {n_{\ell}}, \; \wh{\bs{p}}_{\ell} \leq \f{1}{2} \}$ and
$\wh{p}_{\ell} = \wh{\bs{p}}_{\ell} (\om)$. Then,
$\mscr{M}_{\mrm{B}} ( \wh{p}_{\ell}, \; \wh{p}_{\ell}   + \vep )
> \f{\ln (\ze \de)} {n_{\ell}}$ and $\wh{p}_{\ell} \leq \f{1}{2}$.
Since $z_{\ell} \in [0, \f{1}{2} - \vep)$ and $\mscr{M}_{\mrm{B}}
\li ( z, \; z +  \vep \ri )$ is monotonically increasing with
respect to $z \in (0, \f{1}{2} - \vep)$, it must be true that
$\wh{p}_{\ell} > z_{\ell}$. Otherwise if $\wh{p}_{\ell} \leq
z_{\ell}$, then $\mscr{M}_{\mrm{B}} \li ( \wh{p}_{\ell}, \;
\wh{p}_{\ell}   + \vep \ri ) \leq \mscr{M}_{\mrm{B}} \li ( z_{\ell},
\; z_{\ell}   + \vep \ri ) = \f{\ln (\ze \de)} {n_{\ell}}$, leading
to a contradiction. This proves {\small $ \{ \mscr{M}_{\mrm{B}} \li
( \wh{\bs{p}}_{\ell}, \; \wh{\bs{p}}_{\ell}  + \vep \ri ) > \f{\ln
(\ze \de)} {n_{\ell}} , \; \wh{\bs{p}}_{\ell} \leq \f{1}{2} \}
\subseteq  \{ z_{\ell} < \wh{\bs{p}}_{\ell} \leq \f{1}{2} \}$}  for small enough $\vep > 0$.

Now let $\om \in \li \{ z_{\ell} < \wh{\bs{p}}_{\ell} \leq \f{1}{2}
\ri \}$ and $\wh{p}_{\ell} = \wh{\bs{p}}_{\ell} (\om)$.
Then, $z_{\ell} < \wh{p}_{\ell} \leq \f{1}{2}$. Invoking Lemma \ref{lemm21} that
 there exists a unique number $z^* \in (\f{1}{2}
- \vep, \f{1}{2})$ such that $\mscr{M}_{\mrm{B}}  \li ( z, \; z +  \vep \ri )$ is monotonically
increasing with respect to $z \in (0,  z^*)$ and  monotonically decreasing with respect to $z
\in (z^*, 1 - \vep)$, we have \be \la{recall}
\mscr{M}_{\mrm{B}} \li ( \wh{p}_{\ell}, \; \wh{p}_{\ell}   + \vep
\ri ) > \min \li \{ \mscr{M}_{\mrm{B}} \li ( z_{\ell}, \; z_{\ell} +
\vep \ri ),  \; \mscr{M}_{\mrm{B}} \li ( \f{1}{2}, \; \f{1}{2} +
\vep \ri ) \ri \}. \ee Noting that {\small $\lim_{\vep \to 0} \f{\ln
(\ze \de)} {n_{s}  \mscr{M}_{\mrm{B}} ( \f{1}{2}, \; \f{1}{2} + \vep
)} = 1$},  we have {\small $\mscr{M}_{\mrm{B}} ( \f{1}{2}, \;
\f{1}{2} + \vep )
> \f{\ln (\ze \de)} {n_{\ell}}$} for $\ell < s$ if $\vep > 0$ is small enough.
  By virtue of (\ref{recall}) and $\mscr{M}_{\mrm{B}}
\li ( z_{\ell}, \; z_{\ell} + \vep \ri ) = \f{\ln (\ze \de)}
{n_{\ell}}$, we have $\mscr{M}_{\mrm{B}} \li ( \wh{p}_{\ell}, \;
\wh{p}_{\ell}   + \vep \ri ) > \f{\ln (\ze \de)} {n_{\ell}}$. This
proves {\small $\{ \mscr{M}_{\mrm{B}} \li ( \wh{\bs{p}}_{\ell}, \;
\wh{\bs{p}}_{\ell} + \vep \ri ) > \f{\ln (\ze \de)} {n_{\ell}} , \;
\wh{\bs{p}}_{\ell} \leq \f{1}{2}  \} \supseteq \{ z_{\ell} <
\wh{\bs{p}}_{\ell} \leq \f{1}{2} \}$} and consequently
(\ref{claimabs1}) is established.

To show (\ref{claimabs2}), let $\om \in \{ \mscr{M}_{\mrm{B}} \li (
\wh{\bs{p}}_{\ell}, \; \wh{\bs{p}}_{\ell}  - \vep \ri )  > \f{\ln
(\ze \de)} {n_{\ell}}, \; \wh{\bs{p}}_{\ell} > \f{1}{2} \}$ and
$\wh{p}_{\ell} = \wh{\bs{p}}_{\ell} (\om)$. Then,
$\mscr{M}_{\mrm{B}} ( \wh{p}_{\ell}, \; \wh{p}_{\ell} - \vep  )
> \f{\ln (\ze \de)} {n_{\ell}}$ and $\wh{p}_{\ell} > \f{1}{2}$.
Since $1 - z_{\ell} \in (\f{1}{2} + \vep, 1]$ and
$\mscr{M}_{\mrm{B}} \li ( z, \; z - \vep \ri )$ is monotonically
decreasing with respect to $z \in (\f{1}{2} + \vep, 1)$, it must be
true that $\wh{p}_{\ell} < 1 - z_{\ell}$. Otherwise if
$\wh{p}_{\ell} \geq 1 - z_{\ell}$, then $\mscr{M}_{\mrm{B}} \li (
\wh{p}_{\ell}, \; \wh{p}_{\ell}   - \vep \ri ) \leq
\mscr{M}_{\mrm{B}} \li ( 1 - z_{\ell}, \; 1 - z_{\ell}   - \vep \ri
) = \mscr{M}_{\mrm{B}} \li ( z_{\ell}, \; z_{\ell} + \vep \ri ) =
\f{\ln (\ze \de)} {n_{\ell}}$, leading to a contradiction. This
proves {\small $\{ \mscr{M}_{\mrm{B}}  ( \wh{\bs{p}}_{\ell}, \;
\wh{\bs{p}}_{\ell} - \vep  ) > \f{\ln (\ze \de)} {n_{\ell}} , \;
\wh{\bs{p}}_{\ell} > \f{1}{2} \} \subseteq \{ \f{1}{2} <
\wh{\bs{p}}_{\ell} < 1 - z_{\ell}  \}$}.

Now let $\om \in \li \{ \f{1}{2} < \wh{\bs{p}}_{\ell} < 1 - z_{\ell}
\ri \}$ and $\wh{p}_{\ell} = \wh{\bs{p}}_{\ell} (\om)$.
Then, $\f{1}{2} < \wh{p}_{\ell} < 1 - z_{\ell}$.  Invoking Lemma \ref{lemm21} that
there exists a unique number $z^\star \in
(\f{1}{2}, \f{1}{2} + \vep)$ such that $\mscr{M}_{\mrm{B}}  \li ( z, \; z -  \vep \ri )$ is monotonically
increasing with respect to $z \in (\vep,  z^\star)$ and monotonically decreasing with
respect to $z \in (z^\star, 1)$, we have \be \la{recall2}
\mscr{M}_{\mrm{B}} \li ( \wh{p}_{\ell}, \; \wh{p}_{\ell}   - \vep
\ri ) > \min \li \{  \mscr{M}_{\mrm{B}} \li ( 1 - z_{\ell}, \; 1 -
z_{\ell} - \vep \ri ),  \;  \mscr{M}_{\mrm{B}} \li ( \f{1}{2}, \;
\f{1}{2} - \vep \ri ) \ri \}. \ee Recalling that {\small
$\mscr{M}_{\mrm{B}} \li ( \f{1}{2}, \; \f{1}{2} - \vep \ri ) =
\mscr{M}_{\mrm{B}} \li ( \f{1}{2}, \; \f{1}{2} + \vep \ri )
> \f{\ln (\ze \de)} {n_{\ell}}$} for small enough $\vep
> 0$,  using (\ref{recall2}) and $\mscr{M}_{\mrm{B}}
( 1 - z_{\ell}, \; 1 - z_{\ell} - \vep ) = \mscr{M}_{\mrm{B}} ( z_{\ell}, \;
z_{\ell} + \vep ) = \f{\ln (\ze \de)} {n_{\ell}}$, we have
$\mscr{M}_{\mrm{B}} \li ( \wh{p}_{\ell}, \; \wh{p}_{\ell} - \vep \ri
) > \f{\ln (\ze \de)} {n_{\ell}}$. This proves $ \{
\mscr{M}_{\mrm{B}} \li ( \wh{\bs{p}}_{\ell}, \; \wh{\bs{p}}_{\ell} -
\vep \ri ) > \f{\ln (\ze \de)} {n_{\ell}} , \; \wh{\bs{p}}_{\ell}
> \f{1}{2}  \} \supseteq \{ \f{1}{2} < \wh{\bs{p}}_{\ell} < 1 -
z_{\ell} \}$ and consequently (\ref{claimabs2}) is established. By
virtue of (\ref{claimabs1}) and (\ref{claimabs2}) of the established
claim, we have {\small $\{ \bs{D}_{\ell} = 0 \} = \{ z_{\ell} <
\wh{\bs{p}}_{\ell} \leq \f{1}{2} \} \cup \{ \f{1}{2} <
\wh{\bs{p}}_{\ell} < 1 - z_{\ell} \} = \{ z_{\ell}  <
\wh{\bs{p}}_{\ell} < 1 - z_{\ell} \}$} for small enough $\vep > 0$.
This proves Statement (V).

\beL \la{lem34a} Let $\ell_\vep = s - j_p$.  Then, \be \la{thenlem}
\lim_{\vep \to 0} \sum_{\ell = 1}^{\ell_\vep - 1} n_\ell \Pr \{
\bs{D}_\ell = 1 \} = 0, \qqu \lim_{\vep \to 0} \sum_{\ell =
\ell_\vep + 1}^s n_\ell \Pr \{ \bs{D}_{\ell} = 0 \} = 0 \ee for $p
\in (0, 1)$.  Moreover, $\lim_{\vep \to 0} n_{\ell_\vep} \Pr \{
\bs{D}_{\ell_\vep} = 0 \} = 0$ if $C_{j_p} > 4 p (1 - p)$. \eeL

\bpf

For simplicity of notations, let $b_\ell = \lim_{\vep \to 0}
z_{\ell}$ for $1 \leq \ell < s$.  The proof consists of three main steps as follows.

First, we shall show that (\ref{thenlem}) holds for $p \in (0, \f{1}{2}]$.  By the definition of $\ell_\vep$,
 we have $4 p (1 - p)  > C_{s - \ell_\vep + 1}$.
 Making use of the first three statements of Lemma \ref{lem33a},  we have that {\small $z_\ell < \f{p + b_{\ell_\vep -
1}}{2} < p$} for all $\ell \leq \ell_\vep - 1$ if $\vep$ is
sufficiently small. By the last statement of Lemma \ref{lem33a} and
 using Chernoff bounds, we have {\small \bee \Pr \{ \bs{D}_{\ell} = 1 \} & =
& \Pr \{ \wh{\bs{p}}_{\ell} \leq z_{\ell}
\} + \Pr \{ \wh{\bs{p}}_{\ell} \geq 1 - z_{\ell} \} \leq  \Pr \li \{ \wh{\bs{p}}_{\ell} \leq \f{p + b_{\ell_\vep -
1}}{2} \ri \} + \Pr \li \{ \wh{\bs{p}}_{\ell} \geq 1 - \f{p + b_{\ell_\vep - 1}}{2} \ri \}\\
& \leq & \exp \li ( - 2 n_\ell \li ( \f{p - b_{\ell_\vep - 1}}{2}
\ri )^2 \ri ) + \exp \li ( - 2 n_\ell \li( \f{2 - 3 p - b_{\ell_\vep
- 1}}{2} \ri )^2 \ri ) \eee} for all $\ell \leq \ell_\vep - 1$
provided that $\vep
> 0$ is small enough.  By the definition of $\ell_\vep$, we have
\[
b_{\ell_\vep - 1} = \f{ 1 - \sq{ 1 - C_{s - \ell_\vep + 1}} }{2} <
\f{ 1 - \sq{ 1 - 4 p (1 - p) } }{2}  = p,
\]
which implies that {\small $\li ( \f{p - b_{\ell_\vep - 1}}{2}
\ri )^2$} and {\small $\li( \f{2 - 3 p - b_{\ell_\vep - 1}}{2} \ri )^2$}
are positive constants independent of $\vep > 0$ provided that $\vep > 0$ is small enough.   Hence,
{\small $\lim_{\vep \to 0} \sum_{\ell = 1}^{\ell_\vep - 1} n_\ell \Pr \{
\bs{D}_\ell = 1 \} = 0$} as a result of Lemma \ref{lem31a}.

Similarly, it can be seen from the definition of $\ell_\vep$ that $4
p (1 - p)  < C_{s - \ell_\vep - 1}$. Making use of the first three
statements of Lemma \ref{lem33a}, we have that {\small $z_\ell
> \f{ p + b_{\ell_\vep + 1}}{2} > p$} for $\ell_\vep + 1 \leq \ell <
s$ if $\vep$ is sufficiently small. By the last statement of Lemma
\ref{lem33a} and using Chernoff bound, we have {\small \[ \Pr \{
\bs{D}_{\ell} = 0 \} = \Pr \{ z_{\ell} < \wh{\bs{p}}_{\ell} < 1 -
z_{\ell} \} \leq \Pr \{ \wh{\bs{p}}_{\ell}
> z_{\ell} \} \leq \Pr \li \{ \wh{\bs{p}}_{\ell} > \f{ p +
b_{\ell_\vep + 1}}{2} \ri \} \leq \exp \li ( - 2 n_\ell \li ( \f{ p
- b_{\ell_\vep + 1}}{2} \ri )^2 \ri )
\]}
for $\ell_\vep + 1 \leq \ell < s$ provided that $\vep
> 0$ is small enough.  As a consequence of the definition of $\ell_\vep$, we have that
$b_{\ell_\vep + 1}$ is greater than $p$ and is independent of $\vep > 0$. In view of this and the fact that
$\Pr \{ \bs{D}_s = 0 \} = 0$, we can apply Lemma \ref{lem31a} to conclude that
{\small $\lim_{\vep \to 0} \sum_{\ell = \ell_\vep + 1}^s n_\ell \Pr \{
\bs{D}_\ell = 0 \} = 0$}.

\bsk

Second, we shall show that (\ref{thenlem}) holds for $p \in (\f{1}{2}, 1)$.
As a direct consequence of the definition of $\ell_\vep$,
 we have $4 p (1 - p)  > C_{s - \ell_\vep + 1}$.
 Making use of the first three statements of Lemma \ref{lem33a},  we have that {\small $z_\ell < \f{1 - p +
b_{\ell_\vep - 1}}{2} < 1 - p$} for all $\ell \leq \ell_\vep - 1$ if
$\vep$ is sufficiently small. By the last statement of Lemma
\ref{lem33a} and using Chernoff bounds, we have {\small \bee \Pr \{
\bs{D}_{\ell} = 1 \} & = & \Pr \{ \wh{\bs{p}}_{\ell} \leq z_{\ell}
\} + \Pr \{ \wh{\bs{p}}_{\ell} \geq 1 - z_{\ell} \} \leq  \Pr \li \{
\wh{\bs{p}}_{\ell} \leq \f{1 - p + b_{\ell_\vep - 1}}{2} \ri \} +
\Pr \li \{ \wh{\bs{p}}_{\ell} \geq \f{1 + p -
b_{\ell_\vep - 1}}{2} \ri \}\\
& \leq & \exp \li ( - 2 n_\ell \li ( \f{3 p - 1 - b_{\ell_\vep -
1}}{2} \ri )^2  \ri ) + \exp \li ( - 2 n_\ell \li ( \f{1 - p -
b_{\ell_\vep - 1}}{2} \ri )^2 \ri ) \eee} for all $\ell \leq
\ell_\vep - 1$ provided that $\vep > 0$ is small enough.  As a result of the definition of $\ell_\vep$, we have that
$b_{\ell_\vep - 1}$ is smaller than $1 - p$ and is independent of $\vep
> 0$. Hence, by virtue of Lemma \ref{lem31a}, we have {\small
$\lim_{\vep \to 0} \sum_{\ell = 1}^{\ell_\vep - 1} n_\ell \Pr \{
\bs{D}_\ell = 1 \}  = 0$}.

In a similar manner, by the definition of $\ell_\vep$, we have $4 p
(1 - p)  < C_{s - \ell_\vep - 1}$. Making use of the first three
statements of Lemma \ref{lem33a}, we have that $z_\ell > \f{ 1 - p +
b_{\ell_\vep + 1} }{2} > 1 - p$ for $\ell_\vep + 1 \leq \ell < s$ if
$\vep$ is sufficiently small. By the last statement of Lemma
\ref{lem33a} and using Chernoff bound, \bee \Pr \{ \bs{D}_{\ell} = 0
\} & = & \Pr \{  z_{\ell} < \wh{\bs{p}}_{\ell} < 1 - z_{\ell} \}
\leq \Pr
\{ \wh{\bs{p}}_{\ell} < 1 - z_{\ell} \}\\
& \leq & \Pr \li \{ \wh{\bs{p}}_{\ell} < \f{ 1 + p - b_{\ell_\vep +
1} }{2} \ri \} \leq \exp \li ( - 2 n_\ell \li ( \f{ 1 - p -
b_{\ell_\vep + 1} }{2} \ri )^2 \ri) \eee for $\ell_\vep + 1 \leq
\ell < s$ provided that $\vep
> 0$ is small enough.  Because of the definition of $\ell_\vep$, we have that
$b_{\ell_\vep + 1}$ is greater than $1 - p$ and is independent of $\vep
> 0$.  Noting that $\Pr \{ \bs{D}_s = 0 \} = 0$ and using Lemma \ref{lem31a}, we have
{\small $\lim_{\vep \to 0} \sum_{\ell = \ell_\vep + 1}^s n_\ell \Pr \{
\bs{D}_\ell = 0 \}  = 0$}.

\bsk

Third, we shall show that $\lim_{\vep \to 0} n_{\ell_\vep} \Pr \{
\bs{D}_{\ell_\vep} = 0 \} = 0$  for $p \in (0, 1)$ such that $4 p (1
- p) < C_{j_p}$.

For $p \in (0, \f{1}{2}]$ such that $4 p (1 - p) < C_{j_p}$, making
use of the first three statements of Lemma \ref{lem33a}, we have
{\small $z_{\ell_\vep} > \f{p + b_{\ell_\vep}}{2} > p$} if $\vep$ is
sufficiently small. By the last statement of Lemma \ref{lem33a} and
using Chernoff bound, we have {\small \[ \Pr \{ \bs{D}_{\ell_\vep} =
0 \} = \Pr \{  z_{\ell_\vep} < \wh{\bs{p}}_{\ell_\vep} < 1 -
z_{\ell_\vep} \} \leq \Pr \{
 \wh{\bs{p}}_{\ell_\vep} > z_{\ell_\vep} \} \leq \Pr \li \{
\wh{\bs{p}}_{\ell_\vep} > \f{p + b_{\ell_\vep}}{2} \ri \}  \leq \exp
\li ( - 2 n_{\ell_\vep} \li ( \f{p - b_{\ell_\vep}}{2} \ri )^2 \ri )
\]} for small enough $\vep > 0$. As a consequence of the definition of $\ell_\vep$, we have that
$b_{\ell_\vep}$ is greater than $p$ and is independent of $\vep
> 0$.  It follows
that $\lim_{\vep \to 0} n_{\ell_\vep} \Pr \{ \bs{D}_{\ell_\vep} = 0 \} = 0$.

Similarly, for $p \in (\f{1}{2}, 1)$ such that $4 p (1 - p) <
C_{j_p}$, by virtue of the first three statements of Lemma
\ref{lem33a}, we have {\small $z_{\ell_\vep} > \f{1 - p +
b_{\ell_\vep}}{2} > 1 - p$} if $\vep$ is sufficiently small. By the
last statement of Lemma \ref{lem33a} and using Chernoff bound,
{\small \bee \Pr \{ \bs{D}_{\ell_\vep} = 0 \} & = & \Pr \{
z_{\ell_\vep} < \wh{\bs{p}}_{\ell_\vep} < 1 - z_{\ell_\vep} \} \leq
\Pr \{
\wh{\bs{p}}_{\ell_\vep} < 1 - z_{\ell_\vep}  \}\\
& \leq & \Pr \li \{ \wh{\bs{p}}_{\ell_\vep}  < \f{1 + p -
b_{\ell_\vep}}{2} \ri \} \leq \exp \li ( - 2 n_{\ell_\vep} \li (
\f{1 - p - b_{\ell_\vep}}{2} \ri )^2 \ri ) \eee} for small enough
$\vep > 0$.  Because of the definition of $\ell_\vep$, we have that
$b_{\ell_\vep}$ is greater than $1 - p$ and is independent of $\vep
> 0$.  Hence,  $\lim_{\vep \to 0} n_{\ell_\vep} \Pr \{
\bs{D}_{\ell_\vep} = 0 \} = 0$.

\epf

\bsk

Now we are in a position to prove Theorem \ref{Bino_DDV_Asp}.  To
show $\lim_{\vep \to 0} | \Pr \{ \wh{\bs{p}} \in \mscr{R} \} -
\ovl{P}  | = \lim_{\vep \to 0} | \Pr \{ \wh{\bs{p}} \in \mscr{R} \}
- \udl{P} | = 0$, it suffices to show \be \la{show1}
 \lim_{\vep \to 0}  \sum_{\ell
= 1}^s \Pr \{ \bs{D}_{\ell - 1} = 0, \; \bs{D}_\ell = 1 \} = 1. \ee
This is because $\udl{P} \leq \Pr \{ \wh{\bs{p}} \in \mscr{R} \}
\leq \ovl{P}$ and $\ovl{P} - \udl{P} = \sum_{\ell = 1}^s \Pr \{
\bs{D}_{\ell - 1} = 0, \; \bs{D}_\ell = 1 \} - 1$.  Observing that
{\small \bee
&   & \sum_{\ell = 1}^{\ell_\vep - 1}  \Pr \{ \bs{D}_{\ell - 1} = 0, \;
\bs{D}_\ell = 1 \} \leq  \sum_{\ell = 1}^{\ell_\vep - 1}  \Pr \{
\bs{D}_\ell = 1 \} \leq \sum_{\ell = 1}^{\ell_\vep - 1} n_\ell \Pr
\{ \bs{D}_\ell = 1 \},\\
&   & \sum_{\ell = \ell_\vep + 2}^s  \Pr \{ \bs{D}_{\ell - 1} = 0, \;
\bs{D}_\ell = 1 \} \leq \sum_{\ell = \ell_\vep + 2}^s  \Pr \{
\bs{D}_{\ell - 1} = 0 \} =  \sum_{\ell = \ell_\vep + 1}^s  \Pr \{
\bs{D}_{\ell} = 0 \} \leq  \sum_{\ell = \ell_\vep + 1}^s n_\ell \Pr
\{ \bs{D}_{\ell} = 0 \}
\eee}
and using Lemma \ref{lem34a}, we have {\small $\lim_{\vep \to 0} \sum_{\ell = 1}^{\ell_\vep - 1}  \Pr \{
\bs{D}_{\ell - 1} = 0, \; \bs{D}_\ell = 1 \} = 0$} and {\small $\lim_{\vep
\to 0} \sum_{\ell = \ell_\vep + 2}^s  \Pr \{ \bs{D}_{\ell - 1} = 0,
\; \bs{D}_\ell = 1 \} = 0$}.
Hence, to show (\ref{show1}), it suffices to show $\lim_{\vep \to 0} [ \Pr \{ \bs{D}_{\ell_\vep - 1} = 0, \;
\bs{D}_{\ell_\vep} = 1 \} + \Pr \{ \bs{D}_{\ell_\vep} = 0, \;
\bs{D}_{\ell_\vep + 1} = 1 \} ]  = 1$.  Noting that

{\small $\Pr \{ \bs{D}_{\ell_\vep - 1} = 0,
\; \bs{D}_{\ell_\vep} = 1 \} + \Pr \{ \bs{D}_{\ell_\vep - 1} =
\bs{D}_{\ell_\vep} = 1 \} +  \Pr \{ \bs{D}_{\ell_\vep} = 0, \;
\bs{D}_{\ell_\vep + 1} = 1 \} + \Pr \{ \bs{D}_{\ell_\vep} =
\bs{D}_{\ell_\vep + 1} = 0 \}$}

$ = \Pr \{ \bs{D}_{\ell_\vep } = 1 \} + \Pr \{ \bs{D}_{\ell_\vep} =
0 \} = 1$,

we have {\small \[ \Pr \{ \bs{D}_{\ell_\vep - 1} = 0, \;
\bs{D}_{\ell_\vep} = 1 \} + \Pr \{ \bs{D}_{\ell_\vep} = 0, \;
\bs{D}_{\ell_\vep + 1} = 1 \} = 1 - \Pr \{ \bs{D}_{\ell_\vep - 1} =
\bs{D}_{\ell_\vep} = 1 \} - \Pr \{ \bs{D}_{\ell_\vep} =
\bs{D}_{\ell_\vep + 1} = 0 \}.
\]}
As a result of Lemma \ref{lem34a}, we have {\small $\lim_{\vep \to
0} \Pr \{ \bs{D}_{\ell_\vep - 1} = \bs{D}_{\ell_\vep} = 1 \}  \leq
\lim_{\vep \to 0} \Pr \{ \bs{D}_{\ell_\vep - 1} =  1 \} = 0$} and
{\small $\lim_{\vep \to 0} \Pr \{ \bs{D}_{\ell_\vep} =
\bs{D}_{\ell_\vep + 1} = 0 \} \leq \lim_{\vep \to 0} \Pr \{
\bs{D}_{\ell_\vep + 1} = 0 \} = 0$}. Therefore, $\lim_{\vep \to 0}
\sum_{\ell = 1}^s \Pr \{ \bs{D}_{\ell - 1} = 0, \; \bs{D}_\ell = 1
\} = 1$.  This completes the proof of Theorem \ref{Bino_DDV_Asp}.

\subsection{Proof of Theorem \ref{Bino_Asp_Analysis} } \la{App_Bino_Asp_Analysis}

To prove Theorem \ref{Bino_Asp_Analysis},  we need some preliminary
results.

\beL

\la{lem35a} $\lim_{\vep \to 0} \f{ n_{\ell_\vep} } {
\mcal{N}_{\mrm{a}}  (p, \vep) } = \ka_p, \; \lim_{\vep \to 0} \f{
\vep} {\sq{p ( 1 - p) \sh n_{\ell_\vep} }}   = d \sq{\ka_p}$.

\eeL

\bpf

By the definition of sample sizes, it can be readily shown that
{\small $\lim_{\vep \to 0} \f{ C_{ s - \ell } \ln \f{1}{\ze \de} } {
2 \vep^2 n_\ell} = 1$} for $1 \leq \ell < s$ and it follows that
{\small \bee &   & \lim_{\vep \to 0} \f{ n_{\ell_\vep} } {
\mcal{N}_{\mrm{a}}  (p, \vep) }   =  \lim_{\vep \to 0}  \f{
\mscr{M}_{\mrm{B}} ( \f{1}{2} - | \f{1}{2} - p|, \f{1}{2} - |
\f{1}{2} - p| + \vep)  }{  \ln (\ze \de) } \times \f{ C_{s -
\ell_\vep} } { 2 \vep^2} \ln \f{1}{\ze \de}\\
&   & \qu \qqu \qqu \qu = \lim_{\vep \to 0}  \li [ \f{ \vep^2 }{ 2 p
(1 - p ) } + o (\vep^2) \ri ] \times \f{ C_{s - \ell_\vep} } { 2
\vep^2} = \f{ C_{s
- \ell_\vep} }{4 p (1 - p)}  = \f{C_{j_p} }{4 p (1 - p)} = \ka_p, \\
&   & \lim_{\vep \to 0} \f{ \vep} {\sq{p ( 1 - p) \sh n_{\ell_\vep}
}} = \lim_{\vep \to 0}  \vep \sq{ \f{ C_{s - \ell_\vep} } { 2 \vep^2
p ( 1 - p) } \ln \f{1}{\ze \de} } =  d \sq{ \f{ C_{s - \ell_\vep}
}{4 p (1 - p)} } =  d \sq{ \f{ C_{j_p} }{4 p (1 - p)} } = d
\sq{\ka_p}. \eee}

\epf

\beL \la{limplem} Let $U$ and $V$ be independent Gaussian random
variables with zero means and unit variances.  Let $\ell_\vep = s -
j_p$. Then, for $p \in (0, \f{1}{2}) \cup (\f{1}{2}, 1)$ such that
$C_{j_p} = 4 p (1 - p)$, \bee &   & \lim_{\vep \to 0} \Pr \{ \bs{l}
= \ell_\vep \} = 1 - \lim_{\vep \to 0} \Pr \{ \bs{l} = \ell_\vep + 1
\} =  1 - \Phi \li (
\nu  d \ri ), \\
&   & \lim_{\vep \to 0} \li [ \Pr \{ | \wh{\bs{p}}_{\ell_\vep} - p |
\geq \vep, \; \bs{l}  = \ell_\vep \} + \Pr \{ |
\wh{\bs{p}}_{\ell_\vep + 1} - p | \geq \vep, \;
\bs{l} =  \ell_\vep + 1 \} \ri ]\\
&  &  \qqu \qqu \qu \qu =  \Pr \li \{ U \geq d \ri \} + \Pr \li \{ |U + \sq{\ro_p} V | \geq (1 + \ro_p) d, \; U < \nu d \ri \}. \eee \eeL

\bpf By symmetry, it suffices to show the lemma for $p \in (0,
\f{1}{2})$.   For simplicity of notations, define
\[
a_\ell = \f{ z_\ell - p}{ \sq{ p ( 1 - p) \sh n_\ell} }, \qqu b_\ell
= \f{ \vep } { \sq{ p ( 1 - p) \sh n_\ell} }, \qqu U_\ell = \f{
\wh{\bs{p}}_\ell - p  }{ \sq{ p ( 1 - p) \sh n_\ell} }
\]
for $\ell = 1, \cd, s$.  Since $C_{j_p} = 4 p (1 - p)$, we have
{\small $n_{\ell_\vep} = \li \lc \f{ 2 p ( 1 - p) \ln \f{1}{\ze \de}
} { \vep^2} \ri \rc$} and
\[
\lim_{\vep \to 0} b_{\ell_\vep} = \lim_{\vep \to 0} \f{ \vep } {
\sq{ p ( 1 - p) } }  \sq{ \li \lc \f{ 2 p ( 1 - p) \ln \f{1}{\ze
\de} } { \vep^2} \ri \rc } = \sq{ 2 \ln \f{1}{\ze \de} } = d.
\]
Hence, by Statement (IV) of Lemma \ref{lem33a},
\[
\lim_{\vep \to 0} a_{\ell_\vep} = \lim_{\vep \to 0} b_{\ell_\vep}
\lim_{\vep \to 0} \f{ z_{\ell_\vep} - p } { \vep } = d \lim_{\vep
\to 0} \f{ z_{\ell_\vep} - p } { \vep } =   - \f{2}{3} d = - \nu d.
\]
Let $\eta > 0$. Noting that $\{ \wh{\bs{p}}_{\ell_\vep} \leq
z_{\ell_\vep} \} =   \{ U_{\ell_\vep} \leq a_{\ell_\vep} \}$ and $\{
| \wh{\bs{p}}_{\ell_\vep} - p | \geq \vep \} =  \{ | U_{\ell_\vep} |
\geq b_{\ell_\vep}  \}$,  we have {\small \bee  & & \Pr \li \{
U_{\ell_\vep} \leq - \nu d - \eta \ri \} \leq \Pr \{
\wh{\bs{p}}_{\ell_\vep} \leq z_{\ell_\vep} \} \leq \Pr
\li \{ U_{\ell_\vep}  \leq - \nu d + \eta \ri \},\\
&   & \Pr \li \{ | U_{\ell_\vep} | \geq d + \eta, \; U_{\ell_\vep}
\leq - \nu d - \eta \ri \} \leq \Pr \{ | \wh{\bs{p}}_{\ell_\vep} - p
| \geq \vep, \; \wh{\bs{p}}_{\ell_\vep} \leq z_{\ell_\vep} \} \leq
\Pr \li \{ | U_{\ell_\vep} | \geq d - \eta, \;    U_{\ell_\vep} \leq
- \nu d + \eta \ri \} \eee} for small enough $\vep
> 0$.  Since $U_{\ell_\vep}$ converges in
distribution to a Gaussian random variable $U$ with zero mean and
unit variance as $\vep \to 0$, it must be true that {\small \bee & &
\Pr \li \{ U \leq - \nu d - \eta \ri \} \leq \lim_{\vep \to 0} \Pr
\{ \wh{\bs{p}}_{\ell_\vep} \leq z_{\ell_\vep} \} \leq \Pr \li \{ U
\leq - \nu d + \eta \ri \},\\
&  & \Pr \li \{ | U | \geq d + \eta, \;  U  \leq - \nu d - \eta \ri
\} \leq \lim_{\vep \to 0} \Pr \{ | \wh{\bs{p}}_{\ell_\vep} - p |
\geq \vep, \; \wh{\bs{p}}_{\ell_\vep} \leq z_{\ell_\vep} \} \leq \Pr
\li \{ | U | \geq d - \eta, \; U  \leq - \nu d + \eta \ri \}. \eee}
Since the above inequalities hold true for arbitrarily small $\eta >
0$, we have \bel &
 & \lim_{\vep \to 0} \Pr \{
\wh{\bs{p}}_{\ell_\vep} \leq z_{\ell_\vep} \} = \Pr \li \{  U \leq
- \nu d \ri \} = \Pr \li \{ U \geq \nu d \ri \} = 1 - \Phi (\nu d), \la{suc1}\\
&   & \lim_{\vep \to 0} \Pr \{ | \wh{\bs{p}}_{\ell_\vep} - p | \geq
\vep, \; \wh{\bs{p}}_{\ell_\vep} \leq z_{\ell_\vep} \} = \Pr \li \{
| U | \geq d, \; U \leq - \nu d \ri \} = \Pr \li \{ U \geq d \ri \}.
\la{suc2} \eel Now, we shall consider $\Pr \{ |
\wh{\bs{p}}_{{\ell_\vep} + 1} - p | \geq \vep, \;
\wh{\bs{p}}_{\ell_\vep} > z_{\ell_\vep} \}$. Note that {\small \[
\Pr \{ | \wh{\bs{p}}_{{\ell_\vep} + 1}  - p | \geq \vep, \;
\wh{\bs{p}}_{\ell_\vep}
> z_{\ell_\vep} \} =  \Pr \li \{ | U_{{\ell_\vep} + 1} | \geq
b_{{\ell_\vep} + 1}, \;  U_{\ell_\vep}  > a_{\ell_\vep}  \ri
\} \]}  and
\[
U_{{\ell_\vep} + 1}  = \sq{ \f{ n_{\ell_\vep} } { n_{{\ell_\vep} +
1} }  } U_{\ell_\vep} + \sq{1 - \f{n_{\ell_\vep} } { n_{{\ell_\vep}
+ 1} } } V_{\ell_\vep}, \qu \tx{where} \qu V_{\ell_\vep} = \f{
\sum_{i = n_{\ell_\vep} + 1}^{n_{{\ell_\vep} + 1}} X_i -
(n_{{\ell_\vep} + 1} - n_{\ell_\vep}) p }{ \sq{ p (1 - p) (
n_{{\ell_\vep} + 1} - n_{\ell_\vep} ) } }.
\]
For small enough $\vep > 0$, we have
\[
\Pr \{ | \wh{\bs{p}}_{{\ell_\vep} + 1}  - p | \geq \vep, \; \wh{\bs{p}}_{\ell_\vep} > z_{\ell_\vep} \} \leq \Pr \li \{ \li | U_{{\ell_\vep} + 1}
\ri | \geq \sq{1 + \ro_p} (d - \eta), \;  U_{\ell_\vep}  > - \nu d - \eta \ri \},
\]
\[
\Pr \{ | \wh{\bs{p}}_{{\ell_\vep} + 1}  - p | \geq \vep, \; \wh{\bs{p}}_{\ell_\vep} > z_{\ell_\vep} \} \geq \Pr \li \{ \li | U_{{\ell_\vep} + 1}
\ri | \geq \sq{1 + \ro_p} (d + \eta), \;  U_{\ell_\vep}  > - \nu d + \eta \ri \}.
\]
Note that $U_{\ell_\vep}$ and $V_{\ell_\vep}$ converge in
distribution respectively to independent Gaussian random variables
$U$ and $V$ with zero means and unit variances. Since the
characteristic function of $U_{{\ell_\vep} + 1}$ tends to the
characteristic function of $(U + \sq{\ro_p} V ) \sh \sq{1 + \ro_p}$,
we have {\small \bee &  & \Pr \li \{ \li | U_{{\ell_\vep} + 1} \ri |
\geq d - \eta, \; U_{\ell_\vep} > - \nu d - \eta \ri \} \to  \Pr \li
\{ |U + \sq{\ro_p} V | \geq (1 + \ro_p) (d - \eta) , \; U > - \nu d - \eta \ri \},\\
&  & \Pr \li \{ \li | U_{{\ell_\vep} + 1}  \ri | \geq d + \eta, \; U_{\ell_\vep}  >  - \nu d + \eta \ri \} \to \Pr \li \{ |U + \sq{\ro_p} V |
\geq (1 + \ro_p) (d + \eta), \; U > - \nu d + \eta \ri \} \eee} as $\vep \to 0$.  Since $\eta$ can be arbitrarily small, we have {\small \bel
\lim_{\vep \to 0} \Pr \{ | \wh{\bs{p}}_{{\ell_\vep} + 1} - p | \geq \vep, \; \wh{\bs{p}}_{\ell_\vep} > z_{\ell_\vep} \} & = & \Pr \li \{ |U +
\sq{\ro_p} V | \geq (1 + \ro_p) d, \; U
> - \nu d \ri \} \nonumber\\
& = & \Pr \li \{ |U + \sq{\ro_p} V | \geq (1 + \ro_p) d, \; U < \nu d \ri \} \la{suc3} \eel} for $p \in (0, \f{1}{2})$ such that $C_{j_p} = 4 p
(1 - p)$.    Noting that {\small \bee \Pr \{ \wh{\bs{p}}_{\ell_\vep} \leq z_{\ell_\vep} \; \tx{or} \; \wh{\bs{p}}_{\ell_\vep} \geq 1 -
z_{\ell_\vep} \} \geq \Pr \{ \bs{l} = {\ell_\vep} \}  \geq  \Pr \{ \wh{\bs{p}}_{\ell_\vep} \leq z_{\ell_\vep} \; \tx{or} \;
\wh{\bs{p}}_{\ell_\vep} \geq 1 -
z_{\ell_\vep} \} - \sum_{\ell = 1}^{ \ell_\vep - 1  } \Pr \{ \bs{D}_\ell = 1 \}, \qqu \qqu &  & \\
\Pr \{ 1 - z_{\ell_\vep} > \wh{\bs{p}}_{\ell_\vep} > z_{\ell_\vep} \} \geq \Pr \{ \bs{l} = \ell_\vep + 1 \}  \geq  \Pr \{ 1 - z_{\ell_\vep}
> \wh{\bs{p}}_{\ell_\vep} > z_{\ell_\vep} \} - \Pr \{ \bs{D}_{\ell_\vep + 1} = 0  \}  - \sum_{\ell = 1}^{ \ell_\vep - 1 } \Pr \{ \bs{D}_\ell = 1
\} \; &  & \eee} and using the result that {\small $\lim_{\vep \to 0} \li [ \sum_{\ell = 1}^{{\ell_\vep} - 1} \Pr \{ \bs{D}_\ell = 1 \} + \Pr \{
\bs{D}_{{\ell_\vep} + 1} = 0 \} \ri ] = 0$} as asserted by Lemma \ref{lem34a}, we have $\lim_{\vep \to 0} \Pr \{  \bs{l} = {\ell_\vep} \} =
\lim_{\vep \to 0} \Pr \{ \wh{\bs{p}}_{\ell_\vep} \leq z_{\ell_\vep} \; \tx{or} \; \wh{\bs{p}}_{\ell_\vep} \geq 1 - z_{\ell_\vep} \}$ and
$\lim_{\vep \to 0} \Pr \{ \bs{l} = \ell_\vep + 1 \} = \lim_{\vep \to 0} \Pr \{ 1 - z_{\ell_\vep} > \wh{\bs{p}}_{\ell_\vep} > z_{\ell_\vep} \}$.
We claim that $\lim_{\vep \to 0} \Pr \{\wh{\bs{p}}_{\ell_\vep} \geq 1 - z_{\ell_\vep} \} = 0$ for $p \in (0, \f{1}{2})$.  To show this claim,
note that $\lim_{\vep \to 0} (1 - z_{\ell_\vep} - p) = 1 - 2 p > 0$ as a result of Statement (III) of Lemma \ref{lem33a}. Therefore, $1 -
z_{\ell_\vep} - p > \f{1}{2} - p$ for small enough $\vep > 0$. By virtue of the Chernoff bound, we have $\Pr \{\wh{\bs{p}}_{\ell_\vep} \geq 1 -
z_{\ell_\vep} \} \leq \exp ( - 2 n_{\ell_\vep} (\f{1}{2} - p)^2 )$ for small enough $\vep > 0$, from which the claim immediately follows.  This
implies that \be \la{suc4} \lim_{\vep \to 0} \Pr \{  \bs{l} = {\ell_\vep} \} = \lim_{\vep \to 0} \Pr \{ \wh{\bs{p}}_{\ell_\vep} \leq
z_{\ell_\vep} \}, \qqu \lim_{\vep \to 0} \Pr \{ \bs{l} = \ell_\vep + 1 \} = \lim_{\vep \to 0} \Pr \{ \wh{\bs{p}}_{\ell_\vep}
> z_{\ell_\vep} \}.
\ee Combining (\ref{suc1}) and (\ref{suc4}) yields
\[
\lim_{\vep \to 0} \Pr \{  \bs{l} = {\ell_\vep} \} = 1 - \Phi (\nu
d), \qqu \lim_{\vep \to 0} \Pr \{ \bs{l} = \ell_\vep + 1 \} = \Phi
(\nu d).
\]
Noting that {\small  \bee &   & \Pr \{ | \wh{\bs{p}}_{\ell_\vep} - p
| \geq \vep, \; \bs{l} = {\ell_\vep} \}  \geq  \Pr \{ |
\wh{\bs{p}}_{\ell_\vep} - p | \geq \vep, \; \wh{\bs{p}}_{\ell_\vep}
\notin ( z_{\ell_\vep}, 1 - z_{\ell_\vep} ) \} - \sum_{\ell = 1}^{
\ell_\vep - 1  } \Pr \{ \bs{D}_\ell = 1 \}, \\
&   & \Pr \{ | \wh{\bs{p}}_{\ell_\vep + 1} - p | \geq \vep, \;
\bs{l} = \ell_\vep + 1 \}  \geq   \Pr \{ | \wh{\bs{p}}_{\ell_\vep +
1} - p | \geq \vep, \; \wh{\bs{p}}_{\ell_\vep} \in ( z_{\ell_\vep},
1 - z_{\ell_\vep} ) \}\\
&   & \qqu \qqu  \qqu \qqu \qqu \qqu \qqu - \Pr \{ \bs{D}_{\ell_\vep
+ 1} = 0 \} - \sum_{\ell = 1}^{ \ell_\vep - 1 } \Pr \{ \bs{D}_\ell =
1 \} \eee} and using the result that {\small $\lim_{\vep \to 0} \li
[ \sum_{\ell = 1}^{{\ell_\vep} - 1} \Pr \{ \bs{D}_\ell = 1 \} + \Pr
\{ \bs{D}_{{\ell_\vep} + 1} = 0 \} \ri ] = 0$}, we have \bee  & &
\liminf_{\vep \to 0} \li [ \Pr \{ | \wh{\bs{p}}_{\ell_\vep} - p |
\geq \vep, \; \bs{l} = {\ell_\vep} \} + \Pr \{ |
\wh{\bs{p}}_{\ell_\vep + 1}
- p | \geq \vep, \; \bs{l} = {{\ell_\vep} + 1} \} \ri ]\\
& \geq  & \lim_{\vep \to 0} \li [ \Pr \{ | \wh{\bs{p}}_{\ell_\vep} -
p | \geq \vep, \; \wh{\bs{p}}_{\ell_\vep} \notin ( z_{\ell_\vep}, 1
- z_{\ell_\vep} ) \} + \Pr \{ | \wh{\bs{p}}_{\ell_\vep + 1} - p |
\geq \vep, \; \wh{\bs{p}}_{\ell_\vep} \in ( z_{\ell_\vep}, 1 -
z_{\ell_\vep} ) \} \ri ]. \eee On the other hand, \bee  &   &
\limsup_{\vep \to 0} \li [ \Pr \{ | \wh{\bs{p}}_{\ell_\vep} - p |
\geq \vep, \; \bs{l} = {\ell_\vep} \} + \Pr \{ |
\wh{\bs{p}}_{\ell_\vep + 1}
- p | \geq \vep, \; \bs{l} = {{\ell_\vep} + 1} \} \ri ]\\
& \leq  & \lim_{\vep \to 0} \li [ \Pr \{ | \wh{\bs{p}}_{\ell_\vep} -
p | \geq \vep, \; \wh{\bs{p}}_{\ell_\vep} \notin ( z_{\ell_\vep}, 1
- z_{\ell_\vep} ) \} + \Pr \{ | \wh{\bs{p}}_{\ell_\vep + 1} - p |
\geq \vep, \; \wh{\bs{p}}_{\ell_\vep} \in ( z_{\ell_\vep}, 1 -
z_{\ell_\vep} ) \} \ri ]. \eee Therefore, \bel  &   & \lim_{\vep \to
0} \li [ \Pr \{ | \wh{\bs{p}}_{\ell_\vep} - p | \geq \vep, \; \bs{l}
= {\ell_\vep} \} + \Pr \{ | \wh{\bs{p}}_{\ell_\vep + 1}
- p | \geq \vep, \; \bs{l} = {{\ell_\vep} + 1} \} \ri ] \nonumber\\
& =  & \lim_{\vep \to 0} \li [ \Pr \{ | \wh{\bs{p}}_{\ell_\vep} - p
| \geq \vep, \; \wh{\bs{p}}_{\ell_\vep} \notin ( z_{\ell_\vep}, 1 -
z_{\ell_\vep} ) \} + \Pr \{ | \wh{\bs{p}}_{\ell_\vep + 1} - p | \geq
\vep, \; \wh{\bs{p}}_{\ell_\vep} \in ( z_{\ell_\vep}, 1
- z_{\ell_\vep} ) \} \ri ] \nonumber\\
& =  & \lim_{\vep \to 0} \li [ \Pr \{ | \wh{\bs{p}}_{\ell_\vep} - p
| \geq \vep, \; \wh{\bs{p}}_{\ell_\vep} \leq z_{\ell_\vep} \} + \Pr
\{ | \wh{\bs{p}}_{\ell_\vep + 1} - p | \geq \vep, \;
\wh{\bs{p}}_{\ell_\vep} > z_{\ell_\vep} \} \ri ].  \la{suc5} \eel
Combing (\ref{suc2}), (\ref{suc3}) and (\ref{suc5}) yields \bee & &
\lim_{\vep \to 0} \li [ \Pr \{ | \wh{\bs{p}}_{\ell_\vep} - p | \geq
\vep, \; \bs{l} = {\ell_\vep} \} + \Pr \{ | \wh{\bs{p}}_{\ell_\vep +
1}
- p | \geq \vep, \; \bs{l} = {{\ell_\vep} + 1} \} \ri ] \nonumber\\
& =  & \Pr \{ U \geq d \} + \Pr \li \{ |U + \sq{\ro_p} V | \geq (1 + \ro_p) d, \; U < \nu d \ri \}.  \eee This completes the proof of the lemma.

 \epf

\beL \la{2dGU}

Let $d > 0, \; \ro > 0$ and $0 < \nu \leq 1$.  Let $U$ and $V$ be independent Gaussian variables with zero mean and variance unity. Then,
{\small \bee \Pr \li \{ U \geq d \ri \} + \Pr \{ | U + \sq{\ro} V | \geq (1 + \ro) d, \; U \leq \nu d \} & \geq & 2 - \Phi (d) - \Phi ( d \sq{1 + \ro} ), \\
\Pr \li \{ U \geq d \ri \} + \Pr \{ | U + \sq{\ro} V | \geq (1 + \ro) d, \; U \leq \nu d \} & = & \Psi (\ro, \nu, d) + \Phi (\nu d) - \Phi (d) \\
 &  <  & 3 - 2 \Phi (d \sq{1 + \ro}) - \Phi (d). \eee} \eeL

\bpf

 Clearly, \bee & & \Pr \{ | U + \sq{\ro} V | \geq (1 + \ro ) d,
\; U \leq \nu d \} <  \Pr \{  | U + \sq{\ro} V|  \geq (1 + \ro ) d \} \\
& = &  \Pr \{  |U| \geq d \sq{1 + \ro} \} = 2 [ 1 - \Phi (d \sq{1 + \ro})] \eee

Since $\nu > 0$, we have \bee \Pr \{ | U + \sq{\ro} V | \geq (1 + \ro ) d, \; U \leq \nu d \} & = & \Pr \{ | U + \sq{\ro} V | \geq (1 + \ro ) d, \; U < 0 \}\\
&  &  + \Pr \{ | U + \sq{\ro} V | \geq (1 + \ro ) d, \; 0 \leq U
\leq \nu d \}\\
& > &  \Pr \{ | U + \sq{\ro} V | \geq (1 + \ro ) d, \; U < 0 \}\\
& = & \f{1}{2}  \Pr \{ | U + \sq{\ro} V | \geq (1 + \ro ) d \}\\
& = & \f{1}{2}  \Pr \{  | U | \geq d \sq{1 + \ro} \}  = 1 - \Phi (d \sq{1 + \ro}). \eee  Note that {\small \bee &   & \Pr \li \{ U \geq d \ri \}
+ \Pr \{ \li | U + \sq{\ro} V \ri | \geq (1 + \ro ) d, \; U <
\nu d \}\\
&  =  & \Pr \li \{ U \geq d \ri \} + \Pr \{ U < \nu d \} - \Pr \{ \li | U + \sq{\ro} V \ri | < (1 + \ro ) d,  \;  U
< \nu d \}\\
&  =  & \Pr \li \{ U \geq d \ri \} + \Pr \{ U < \nu d \} - 1 + \Pr \{ \li | U + \sq{\ro} V \ri | \geq (1 + \ro ) d
\; \mrm{or} \; U \geq \nu d \}\\
&  =  & \Pr \li \{ U \geq d \ri \} - \Pr \{ U \geq \nu d \} + \Pr \{ \li | U + \sq{\ro} V \ri | \geq (1 + \ro ) d
\; \mrm{or} \; U \geq \nu d \}\\
&  = & \Pr \{ \li | U + \sq{\ro} V \ri | \geq (1 + \ro ) d \; \mrm{or} \; U \geq \nu d \} - \Pr \{ \nu d \leq U < d \} \eee} and that $\Pr \li
\{ \li | U + \sq{\ro} V \ri | \geq (1 + \ro ) d \; \mrm{or} \;  U \geq \nu d \ri \}$ is the probability that $(U, V)$ is included in a domain
with a boundary which is visible for an observer on the origin and can be represented in polar coordinates $(r, \phi)$ as
\[
\li \{ (r, \phi): r = \f{ \nu d } { | \cos \phi | }, \; - \phi_L \leq \phi \leq \phi_U \ri \} \cup \li \{ (r, \phi): r = \f{ d \sq{1 + \ro} } {
| \cos ( \phi - \phi_{\ro}) | }, \; \phi_U \leq \phi \leq 2 \pi - \phi_L \ri \}.
\]
Hence, by Theorem 6 of \cite{Chen_TH}, we can show that $\Pr \li \{ \li | U + \sq{\ro} V \ri | \geq (1 + \ro ) d \; \mrm{or} \;  U \geq \nu d
\ri \} = \Psi (\ro, \nu, d)$.  The lemma follows immediately.

\epf

\subsubsection{Proof of Statement (I) }

First, we shall show that Statement (I) holds for $p \in (0,
\f{1}{2}]$ such that $C_{j_p} = 4 p (1 - p)$. For this purpose, we
need to show that {\small \be \la{doit} 1  \leq \limsup_{\vep \to 0}
\f{ \mbf{n} (\om)  } {  \mcal{N}_{\mrm{a}}  (p, \vep)  } \leq 1 +
\ro_p \qqu \tx{ for any} \; \om \in \li \{ \lim_{\vep \to 0}
\wh{\bs{p}} = p \ri \}. \ee}
 To show {\small $\limsup_{\vep \to 0}
\f{ \mbf{n} (\om)  } {  \mcal{N}_{\mrm{a}}  (p, \vep)  } \geq 1$},
note that $C_{s - \ell_\vep + 1} < 4 p ( 1 - p) = C_{s - \ell_\vep}
< C_{s - \ell_\vep - 1}$ as a direct consequence of the definition
of $\ell_\vep$ and the assumption that $C_{j_p} = 4 p (1 - p)$. By
the first three statements of  Lemma \ref{lem33a}, we have
$\lim_{\vep \to 0} z_{{\ell}} < p$ for all $\ell \leq {\ell_\vep} -
1$.  Noting that $\lim_{\vep \to 0} \wh{\bs{p}} (\om) = p \leq
\f{1}{2}$, we have $z_\ell < \wh{\bs{p}} (\om) < 1 - z_\ell$ for all
$\ell \leq {\ell_\vep} - 1$ and it follows from the definition of
the sampling scheme that $\mbf{n} (\om) \geq n_{\ell_\vep}$ if $\vep
> 0$ is small enough. By Lemma \ref{lem35a} and noting that $\ka_p =
1$ if $C_{j_p} = 4 p (1 - p)$, we have
 {\small $\limsup_{\vep \to 0} \f{ \mbf{n} (\om)  } {  \mcal{N}_{\mrm{a}}  (p, \vep)  }
\geq \lim_{\vep \to 0}  \f{ n_{{\ell_\vep}}} { \mcal{N}_{\mrm{a}}
(p, \vep) } = \ka_p = 1$}.   To show {\small $\limsup_{\vep \to 0}
\f{ \mbf{n} (\om)  } { \mcal{N}_{\mrm{a}}  (p, \vep)  } \leq 1 +
\ro_p$}, we shall consider three cases: (i) ${\ell_\vep} = s$; (ii)
${\ell_\vep} = s - 1$; (iii) ${\ell_\vep} < s - 1$.  In the case of
${\ell_\vep} = s$, it must be true that $\mbf{n} (\om) \leq n_s =
n_{\ell_\vep}$.  Hence, $\limsup_{\vep \to 0} \f{ \mbf{n} (\om)  } {
\mcal{N}_{\mrm{a}} (p, \vep)  } \leq \lim_{\vep \to 0} \f{
n_{\ell_\vep} } { \mcal{N}_{\mrm{a}}  (p, \vep) } = \ka_p = 1 = 1 +
\ro_p$.  In the case of ${\ell_\vep} = s - 1$, it must be true that
$\mbf{n} (\om) \leq n_s = n_{{\ell_\vep} + 1}$. Therefore, {\small
$\limsup_{\vep \to 0} \f{ \mbf{n} (\om)  } {  \mcal{N}_{\mrm{a}} (p,
\vep)  } \leq \lim_{\vep \to 0}  \f{ n_{{\ell_\vep} + 1}} {
\mcal{N}_{\mrm{a}} (p, \vep) } = \lim_{\vep \to 0}  \f{
n_{{\ell_\vep} + 1}} { n_{{\ell_\vep} }} \times \lim_{\vep \to 0}
\f{ n_{{\ell_\vep} }} { \mcal{N}_{\mrm{a}}  (p, \vep) } = \f{ C_{j_p
- 1} } { C_{j_p} } = 1 + \ro_p$}.  In the case of ${\ell_\vep} < s -
1$, it follows from Lemma \ref{lem33a} that $\lim_{\vep \to 0}
z_{{\ell_\vep} + 1} > p$, which implies that
 $z_{{\ell_\vep} + 1} > p, \; \wh{\bs{p}} (\om) < z_{{\ell_\vep} + 1}$, and thus $\mbf{n} (\om) \leq
n_{{\ell_\vep} + 1}$ for small enough $\vep
> 0$. Therefore, {\small $\limsup_{\vep \to 0} \f{ \mbf{n} (\om)  } {  \mcal{N}_{\mrm{a}}  (p, \vep)  } \leq
\lim_{\vep \to 0}  \f{ n_{{\ell_\vep} + 1}} { \mcal{N}_{\mrm{a}} (p,
\vep) } = \lim_{\vep \to 0}  \f{ n_{{\ell_\vep} + 1}} {
n_{{\ell_\vep} }} \times \lim_{\vep \to 0}  \f{ n_{{\ell_\vep} }} {
\mcal{N}_{\mrm{a}}  (p, \vep) } = \f{ C_{j_p - 1} } { C_{j_p} }  = 1
+ \ro_p$}.  This establishes (\ref{doit}), which implies $\{ 1 \leq
\limsup_{\vep \to 0} \f{ \mbf{n}  } {  \mcal{N}_{\mrm{a}}  (p, \vep)
} \leq 1 + \ro_p \} \supseteq \li \{ \lim_{\vep \to 0} \wh{\bs{p}} =
p \ri \}$. Applying the strong law of large numbers, we have $1 \geq
\Pr \{ 1 \leq \limsup_{\vep \to 0} \f{ \mbf{n}  } {
\mcal{N}_{\mrm{a}}  (p, \vep) } \leq 1 + \ro_p \} \geq \Pr \li \{
\lim_{\vep \to 0} \wh{\bs{p}} = p \ri \} = 1$.  This proves that
Statement (I) holds for $p \in (0, \f{1}{2} ]$ such that $C_{j_p} =
4 p (1 - p)$.

\bsk

Next, we shall show that Statement (I) for $p \in (0, \f{1}{2}]$
such that $C_{j_p} > 4 p (1 - p)$. Note that $C_{s - \ell_\vep + 1}
< 4 p ( 1 - p) < C_{s - \ell_\vep}$ as a direct consequence of the
definitions of $\ell_\vep$ and $j_p$.  By the first three statements
of Lemma \ref{lem33a}, we have $\lim_{\vep \to 0} z_{{\ell_\vep} -
1} < p \leq \f{1}{2}$. It follows that $z_{{\ell}} < p \leq
\f{1}{2}$ for all $\ell \leq {\ell_\vep} - 1$ provided that $\vep
> 0$ is sufficiently small. Therefore, for any $\om \in \li \{
\lim_{\vep \to 0} \wh{\bs{p}} = p \ri \}$, we have $z_\ell <
\wh{\bs{p}} (\om) < 1 - z_\ell$ for all $\ell \leq {\ell_\vep} - 1$
and consequently, $\mbf{n} (\om) \geq n_{{\ell_\vep}}$ provided that
$\vep > 0$ is sufficiently small.  On the other hand, we claim that
$\mbf{n} (\om) \leq n_{{\ell_\vep}}$ provided that $\vep > 0$ is
sufficiently small.  Clearly, this claim is true if $\ell_\vep = s$.
In the case of $\ell_\vep < s$,  by the first three statements of
Lemma \ref{lem33a}, we have $\lim_{\vep \to 0} z_{{\ell_\vep}} > p$
as a consequence of $4 p ( 1 - p) < C_{s - \ell_\vep}$. Hence,
$\wh{\bs{p}} (\om) < z_{{\ell_\vep}}$ provided that $\vep > 0$ is
sufficiently small, which implies that the claim is also true in the
case of $\ell_\vep < s$.  Therefore, $\mbf{n} (\om) =
n_{{\ell_\vep}}$ provided that $\vep > 0$ is sufficiently small.
Applying Lemma \ref{lem35a},  we have {\small $\lim_{\vep \to 0} \f{
\mbf{n} (\om)}
 {  \mcal{N}_{\mrm{a}}  (p, \vep)  } = \lim_{\vep \to 0}
 \f{ n_{{\ell_\vep}}} { \mcal{N}_{\mrm{a}}  (p, \vep) } = \ka_p$}, which implies that
{\small $\{ \lim_{\vep \to 0} \f{ \mbf{n} } {  \mcal{N}_{\mrm{a}}
(p, \vep)  } = \ka_p \} \supseteq \li \{ \lim_{\vep \to 0}
\wh{\bs{p}} = p \ri \}$}.  It follows from the strong law of large
numbers that {\small $1 \geq \Pr \{ \lim_{\vep \to 0} \f{ \mbf{n} }
{ \mcal{N}_{\mrm{a}}  (p, \vep)  } = \ka_p \} \geq \Pr \{ \lim_{\vep
\to 0} \wh{\bs{p}} = p \}$} and thus {\small $\Pr \{ \lim_{\vep \to
0} \f{ \mbf{n} } {  \mcal{N}_{\mrm{a}}  (p, \vep)  } = \ka_p \} =
1$}.  Since $1 \leq \ka_p \leq 1 + \ro_p$, it is obviously true that
$\Pr \{ 1 \leq \limsup_{\vep \to 0} \f{ \mbf{n}  } {
\mcal{N}_{\mrm{a}} (p, \vep) } \leq 1 + \ro_p \} = 1$.  This proves
that Statement (I) holds for $p \in (0, \f{1}{2}]$ such that
$C_{j_p} > 4 p (1 - p)$.

\bsk

In a similar manner, we can show that Statement (I) holds for $p \in
(\f{1}{2}, 1)$. This concludes the proof of Statement (I).

\subsubsection{Proof of Statement (II)}

In the sequel, we will consider the asymptotic value of $\f{ \bb{E}
[ \mbf{n} ] } { \mcal{N}_{\mrm{a}}  (p, \vep)  }$ in three steps.
First, we shall show Statement (II) for $p \in (0, 1)$ such that
$C_{j_p} = 4 p (1 - p)$ and $j_p \geq 1$. By the definition of the
sampling scheme, we have {\small \bee \bb{E} [ \mbf{n} ] & = &
\sum_{\ell = 1}^{{\ell_\vep} - 1} n_\ell \Pr \{ \bs{l} = \ell \} +
\sum_{\ell = {\ell_\vep} + 2}^s n_\ell \Pr \{ \bs{l} = \ell \} +
n_{\ell_\vep} \Pr \{ \bs{l} = {\ell_\vep} \} +
n_{{\ell_\vep} + 1} \Pr \{ \bs{l} = {{\ell_\vep} + 1} \}\\
& \leq &  \sum_{\ell =  1}^{{\ell_\vep} - 1} n_\ell \Pr \{
\bs{D}_\ell = 1 \} + \sum_{\ell =  {\ell_\vep} + 1}^{s - 1} n_{\ell
+ 1} \Pr \{ \bs{D}_\ell = 0 \} + n_{\ell_\vep} \Pr \{ \bs{l} =
{\ell_\vep} \} + n_{{\ell_\vep} + 1} \Pr \{ \bs{l} = {{\ell_\vep} +
1} \} \eee} and {\small $\bb{E} [ \mbf{n} ] \geq n_{\ell_\vep} \Pr
\{ \bs{l} = {\ell_\vep} \} + n_{{\ell_\vep} + 1}
 \Pr \{ \bs{l} = {{\ell_\vep} + 1} \} $}.
Making use of Lemma \ref{lem34a} and the assumption that $\sup_{\ell
> 0} \f{ n_{\ell + 1} }{ n_{\ell} } < \iy $, we have
{\small \bee &  & \lim_{\vep \to 0} \li [ \sum_{\ell =
1}^{{\ell_\vep} - 1} n_\ell \Pr \{ \bs{D}_\ell = 1 \} + \sum_{\ell =
{\ell_\vep} + 1}^{s - 1} n_{\ell + 1} \Pr \{ \bs{D}_\ell = 0  \} \ri ]\\
&   & \leq  \lim_{\vep \to 0} \li [ \sum_{\ell = 1}^{{\ell_\vep} -
1} n_\ell \Pr \{ \bs{D}_\ell = 1 \} + \sup_{\ell
> 0} \f{ n_{\ell + 1} }{ n_{\ell} } \sum_{\ell = {\ell_\vep} +
1}^{s - 1} n_{\ell} \Pr \{ \bs{D}_\ell = 0  \} \ri ] = 0. \eee}
Therefore, {\small \bee &   & \limsup_{\vep \to 0} \f{ \bb{E} [
\mbf{n} ] } { \mcal{N}_{\mrm{a}}  (p, \vep)  }\\
 & \leq & \lim_{\vep
\to 0} \f{ \sum_{\ell = 1}^{{\ell_\vep} - 1} n_\ell \Pr \{
\bs{D}_\ell = 1 \} + \sum_{\ell = {\ell_\vep} + 1}^{s - 1} n_{\ell +
1} \Pr \{ \bs{D}_\ell = 0  \}  +
n_{\ell_\vep} \Pr \{ \bs{l} = \ell_\vep \} + n_{{\ell_\vep} + 1}
\Pr \{ \bs{l} = \ell_\vep + 1 \} } { \mcal{N}_{\mrm{a}}  (p, \vep) }\\
& = & \lim_{\vep \to 0} \f{ n_{\ell_\vep} \Pr \{ \bs{l} = \ell_\vep
\} + n_{{\ell_\vep} + 1} \Pr \{ \bs{l} = \ell_\vep + 1 \} } {
\mcal{N}_{\mrm{a}}  (p, \vep) }. \eee}   On the other hand,
\[
\liminf_{\vep \to 0} \f{ \bb{E} [ \mbf{n} ] } { \mcal{N}_{\mrm{a}}
(p, \vep)  } \geq \lim_{\vep \to 0} \f{ n_{\ell_\vep} \Pr \{ \bs{l}
= \ell_\vep \} + n_{{\ell_\vep} + 1} \Pr \{ \bs{l} = \ell_\vep + 1
\} } { \mcal{N}_{\mrm{a}}  (p, \vep) }.
\]
It follows that
\[
\lim_{\vep \to 0} \f{ \bb{E} [ \mbf{n} ] } { \mcal{N}_{\mrm{a}} (p,
\vep)  } = \lim_{\vep \to 0} \f{ n_{\ell_\vep} \Pr \{ \bs{l} =
\ell_\vep \} + n_{{\ell_\vep} + 1} \Pr \{ \bs{l} = \ell_\vep + 1 \}
} { \mcal{N}_{\mrm{a}}  (p, \vep) }
\]
for $p \in (0, 1)$ such that $C_{j_p} = 4 p (1 - p)$ and $j_p \geq
1$.   Using Lemma \ref{limplem} and the result $\lim_{\vep \to 0}
\f{ n_{\ell_\vep} } { \mcal{N}_{\mrm{a}} (p, \vep) } = \ka_p$ as
asserted by Lemma \ref{lem35a},  we have \bee \lim_{\vep \to 0} \f{
n_{\ell_\vep} \Pr \{ \bs{l} = \ell_\vep \} + n_{{\ell_\vep} + 1} \Pr
\{ \bs{l} = \ell_\vep + 1 \} } { \mcal{N}_{\mrm{a}}  (p, \vep) } & =
& \lim_{\vep \to 0} \f{ n_{\ell_\vep} [ 1 - \Phi (  \nu d ) ] +
n_{\ell_\vep + 1} \Phi ( \nu d ) } { \mcal{N}_{\mrm{a}}
(p, \vep) } \\
& = & 1  + \ro_p \Phi \li (  \nu d \ri ). \eee

Second, we shall show Statement (II) for $p \in (0, 1)$ such that
$C_{j_p} = 4 p (1 - p)$ and $j_p = 0$. In this case, it must be true
that $p = \f{1}{2}$.  By the definition of the sampling scheme, we
have
\[
\bb{E} [ \mbf{n} ]  = \sum_{\ell =  1}^{{\ell_\vep} - 1} n_\ell \Pr
\{ \bs{l} = \ell \} +  n_{\ell_\vep} \Pr \{ \bs{l} = {\ell_\vep} \}
\leq \sum_{\ell =  1}^{{\ell_\vep} - 1} n_\ell \Pr \{ \bs{D}_\ell =
1 \} +  n_{\ell_\vep}
\]
and {\small $\bb{E} [ \mbf{n} ] \geq n_{\ell_\vep} \Pr \{ \bs{l} =
{\ell_\vep} \}  \geq n_{\ell_\vep}  \li ( 1 -  \sum_{\ell =
1}^{{\ell_\vep} - 1} \Pr \{ \bs{D}_\ell = 1 \} \ri )$}.  Therefore,
by Lemma \ref{lem34a},
\[
\limsup_{\vep \to 0} \f{ \bb{E} [ \mbf{n} ] } { \mcal{N}_{\mrm{a}}
(p, \vep)  } \leq \lim_{\vep \to 0}  \f{ \sum_{\ell =
1}^{{\ell_\vep} - 1} n_\ell \Pr \{ \bs{D}_\ell = 1 \} +
n_{\ell_\vep} } { \mcal{N}_{\mrm{a}}  (p, \vep) } = \lim_{\vep \to
0}  \f{ n_{\ell_\vep} } { \mcal{N}_{\mrm{a}}  (p, \vep) } = \ka_p =
1,
\]
\[
\liminf_{\vep \to 0} \f{ \bb{E} [ \mbf{n} ] } { \mcal{N}_{\mrm{a}}
(p, \vep)  } \geq \lim_{\vep \to 0} \f{ n_{\ell_\vep}  \li ( 1 -
\sum_{\ell = 1}^{{\ell_\vep} - 1} \Pr \{ \bs{D}_\ell = 1 \} \ri ) }
{ \mcal{N}_{\mrm{a}}  (p, \vep) } = \lim_{\vep \to 0}  \f{
n_{\ell_\vep} } { \mcal{N}_{\mrm{a}}  (p, \vep) } = \ka_p = 1
\]
and thus $\lim_{\vep \to 0} \f{ \bb{E} [ \mbf{n} ] } {
\mcal{N}_{\mrm{a}} (p, \vep)  } = 1$ for $p \in (0, 1)$ such that
$C_{j_p} = 4 p (1 - p)$ and $j_p = 0$.

Third, we shall show Statements (II)  for $p \in (0, 1)$ such that
$C_{j_p} > 4 p (1 - p)$ .  Note that \bee \bb{E} [ \mbf{n} ] & = &
\sum_{\ell = 1}^{{\ell_\vep} - 1} n_\ell \Pr \{ \bs{l} = \ell \} +
\sum_{\ell = {\ell_\vep} + 1}^s n_\ell \Pr \{
\bs{l} = \ell \} + n_{\ell_\vep} \Pr \{ \bs{l} = {\ell_\vep} \} \\
& \leq & \sum_{\ell =  1}^{{\ell_\vep} - 1} n_\ell \Pr \{
\bs{D}_\ell = 1 \} + \sum_{\ell =  {\ell_\vep}}^{s - 1} n_{\ell + 1}
\Pr \{ \bs{D}_\ell = 0  \} + n_{\ell_\vep} \eee and {\small $\bb{E}
[ \mbf{n} ] \geq n_{\ell_\vep} \Pr \{ \bs{l} = {\ell_\vep} \} \geq
n_{\ell_\vep}  \li ( 1 -  \sum_{\ell =  1}^{{\ell_\vep} - 1} \Pr \{
\bs{D}_\ell = 1 \} - \Pr \{  \bs{D}_{\ell_\vep} = 0 \} \ri )$}.
Therefore, by Lemma \ref{lem34a}, {\small \bee &   & \limsup_{\vep
\to 0} \f{ \bb{E} [ \mbf{n} ] } { \mcal{N}_{\mrm{a}} (p, \vep)  }
\leq \lim_{\vep \to 0}  \f{ \sum_{\ell = 1}^{{\ell_\vep} - 1} n_\ell
\Pr \{ \bs{D}_\ell = 1 \} + \sum_{\ell = {\ell_\vep}}^{s - 1}
n_{\ell + 1} \Pr \{ \bs{D}_\ell = 0 \} + n_{\ell_\vep} } {
\mcal{N}_{\mrm{a}} (p, \vep) } = \lim_{\vep \to 0}  \f{
n_{\ell_\vep} } {
\mcal{N}_{\mrm{a}}  (p, \vep) } = \ka_p,\\
&   & \liminf_{\vep \to 0} \f{ \bb{E} [ \mbf{n} ] } {
\mcal{N}_{\mrm{a}} (p, \vep)  } \geq \lim_{\vep \to 0} \f{
n_{\ell_\vep}  \li ( 1 - \sum_{\ell = 1}^{{\ell_\vep} - 1} \Pr \{
\bs{D}_\ell = 1 \} - \Pr \{ \bs{D}_{\ell_\vep} = 0 \} \ri ) } {
\mcal{N}_{\mrm{a}}  (p, \vep) } = \lim_{\vep \to 0}  \f{
n_{\ell_\vep} } { \mcal{N}_{\mrm{a}}  (p, \vep) } = \ka_p. \eee} So,
$\lim_{\vep \to 0} \f{ \bb{E} [ \mbf{n} ] } { \mcal{N}_{\mrm{a}} (p,
\vep)  } = \ka_p$ for $p \in (0, 1)$ such that $C_{j_p} > 4 p (1 -
p)$.   From the preceding  analysis, we have shown $\lim_{\vep \to
0} \f{ \bb{E} [ \mbf{n} ] } { \mcal{N}_{\mrm{a}}  (p, \vep) }$
exists for all $p \in (0, 1)$. Hence, statement (II) is established
by making use of this result and the fact that
\[
\lim_{\vep \to 0} \f{ \bb{E} [ \mbf{n} ] } { \mcal{N}_{\mrm{f}} (p,
\vep) } = \lim_{\vep \to 0} \f{ \mcal{N}_{\mrm{a}} (p, \vep) } {
\mcal{N}_{\mrm{f}} (p, \vep) } \times \lim_{\vep \to 0} \f{ \bb{E} [
\mbf{n} ] } { \mcal{N}_{\mrm{a}} (p, \vep) } = \f{ 2 \ln \f{1}{\ze
\de} } { \mcal{Z}_{\ze \de}^2  } \times \lim_{\vep \to 0} \f{ \bb{E}
[ \mbf{n} ] } { \mcal{N}_{\mrm{a}} (p, \vep) }.
\]

\subsubsection{Proof of Statement (III)}

As before, we use the notations {\small $b_{\ell} = \f{ \vep } {
\sq{ p ( 1 - p) \sh n_{\ell}} }$} and {\small $U_{\ell} = \f{
\wh{\bs{p}}_{\ell} - p  }{ \sq{ p ( 1 - p) \sh n_{\ell}} }$}.

First, we shall consider $p \in (0, 1)$ such that $C_{j_p} > 4 p (1
- p)$. Applying Lemma \ref{lem34a} based on
 the assumption that $C_{j_p} > 4 p (1 - p)$, we have  {\small \bee &  &  \lim_{\vep \to 0} \Pr \{ \bs{l} <
{{\ell_\vep}} \} \leq \lim_{\vep \to 0} \sum_{\ell = 1}^{{\ell_\vep}
- 1} \Pr \{ \bs{D}_\ell = 1 \} \leq \lim_{\vep \to 0} \sum_{\ell =
1}^{{\ell_\vep} - 1}
n_\ell \Pr \{ \bs{D}_\ell = 1 \} = 0,\\
&   & \lim_{\vep \to 0} \Pr \{ \bs{l} > {{\ell_\vep}} \} \leq
\lim_{\vep \to 0} \Pr \{ \bs{D}_{\ell_\vep} = 0 \} \leq \lim_{\vep
\to 0} n_{{\ell_\vep}} \Pr \{ \bs{D}_{\ell_\vep} = 0 \} = 0 \eee}
and thus $\lim_{\vep \to 0} \Pr \{ \bs{l} \neq {\ell_\vep}  \} = 0$.
Note that $\Pr \{ | \wh{\bs{p}} - p | \geq \vep  \} = \Pr \{ |
\wh{\bs{p}}_{\ell_\vep} - p | \geq \vep, \; \bs{l} = {\ell_\vep} \}
+ \Pr \{ | \wh{\bs{p}} - p | \geq \vep, \; \bs{l} \neq {\ell_\vep}
\}$ and, as a result of the central limit theorem, $U_{\ell_\vep}$
converges in distribution to a standard Gaussian variable $U$.
Hence,  {\small
\[ \lim_{\vep \to 0} \Pr \{ | \wh{\bs{p}} - p | \geq \vep \} =
\lim_{\vep \to 0} \Pr \{ | \wh{\bs{p}}_{\ell_\vep} - p | \geq \vep
\} = \lim_{\vep \to 0} \Pr \li \{ |U_{\ell_\vep}| \geq b_{\ell_\vep}
\ri \} = \Pr \{ |U| \geq d \sq{\ka_p} \}
\]}
and $\lim_{\vep \to 0} \Pr \{ | \wh{\bs{p}} - p | < \vep  \} = \Pr
\{ |U| < d \sq{\ka_p} \} = 2 \Phi (d \sq{\ka_p}) - 1 > 2 \Phi (d) -
1 > 1 - 2 \ze \de$ for $p \in (0, 1)$ such that $C_{j_p} > 4 p (1 -
p)$.

Second, we shall consider $p \in (0, 1)$ such that $C_{j_p} = 4 p (1
- p)$ and $j_p \geq  1$. In this case, it is evident that $\ell_\vep
< s$. By the definition of the sampling scheme, we have that $\Pr \{
\bs{l} > {{\ell_\vep} + 1} \} \leq \Pr \{ \bs{D}_{{\ell_\vep} + 1} =
0 \}$ and that $\Pr \{ \bs{l} = \ell \} \leq \Pr \{ \bs{D}_\ell = 1
\}$ for $\ell < {\ell_\vep}$.  As a result of Lemma \ref{lem34a}, we
have $\lim_{\vep \to 0} \Pr \{ \bs{l} > {{\ell_\vep} + 1} \} \leq
\lim_{\vep \to 0} \Pr \{ \bs{D}_{{\ell_\vep} + 1} = 0 \} = 0$ and
$\lim_{\vep \to 0} \Pr \{ \bs{l} < {\ell_\vep} \} \leq \lim_{\vep
\to 0} \sum_{\ell =  1}^{ {\ell_\vep} - 1} \Pr \{ \bs{D}_\ell = 1 \}
= 0$.  Since \bee  \limsup_{\vep \to 0} \Pr \{ | \wh{\bs{p}} - p |
\geq \vep \} & \leq & \lim_{\vep \to 0} \li [ \Pr \{ |
\wh{\bs{p}}_{\ell_\vep} - p | \geq \vep, \; \bs{l} = {\ell_\vep} \}
+ \Pr \{ | \wh{\bs{p}}_{\ell_\vep + 1} - p | \geq \vep, \; \bs{l} =
{{\ell_\vep} + 1} \} \ri ]\\
&  & + \lim_{\vep \to 0} \Pr \{ \bs{l} < {\ell_\vep} \} + \lim_{\vep
\to 0} \Pr \{ \bs{l} > {{\ell_\vep} + 1} \} \eee and {\small
$\liminf_{\vep \to 0} \Pr \{ | \wh{\bs{p}} - p | \geq \vep \} \geq
\lim_{\vep \to 0} \li [ \Pr \{ | \wh{\bs{p}}_{\ell_\vep} - p | \geq
\vep, \; \bs{l} = {\ell_\vep} \} + \Pr \{ | \wh{\bs{p}}_{\ell_\vep +
1} - p | \geq \vep, \; \bs{l} = {{\ell_\vep} + 1} \} \ri ]$},  we
have {\small $\lim_{\vep \to 0} \Pr \{ | \wh{\bs{p}} - p | \geq \vep
\}   = \lim_{\vep \to 0} \li [ \Pr \{ | \wh{\bs{p}}_{\ell_\vep} - p
| \geq \vep, \; \bs{l} = {\ell_\vep} \} + \Pr \{ |
\wh{\bs{p}}_{\ell_\vep + 1} - p | \geq \vep, \; \bs{l} =
{{\ell_\vep} + 1} \} \ri ]$}.   By Lemma \ref{limplem}, we have
$\lim_{\vep \to 0} \Pr \{ | \wh{\bs{p}} - p | \geq \vep
 \} =  \Pr \li \{ U \geq d \ri \} + \Pr \li \{ |U + \sq{\ro} V | \geq
  (1 + \ro_p) d, \; U < \nu d \ri \}$ for $p \in (0, 1)$ such that $C_{j_p} = 4 p (1 - p)$ and  $j_p \geq
 1$.  As a consequence of Lemma \ref{2dGU}, Statement (III) must be
 true for $p \in (0, 1)$ such that $C_{j_p} = 4 p (1 - p)$ and  $j_p \geq
 1$.

Third, we shall consider $p \in (0, 1)$ such that $C_{j_p} = 4 p (1
- p)$ and $j_p = 0$.  In this case, it must be true that $p =
\f{1}{2}$.  Clearly, $\ell_\vep = s$. It follows from the definition
of the sampling scheme that $\Pr \{ \bs{l} = \ell \} \leq \Pr \{
\bs{D}_\ell = 1 \}$ for $\ell < {\ell_\vep}$. By Lemma \ref{lem34a},
we have $\lim_{\vep \to 0} \Pr \{ \bs{l} < {\ell_\vep} \} \leq
\lim_{\vep \to 0} \sum_{\ell =  1}^{ {\ell_\vep} - 1} \Pr \{
\bs{D}_\ell = 1 \} = 0$. Therefore, $\lim_{\vep \to 0} \Pr \{ \bs{l}
= {\ell_\vep} \} = 1$ and {\small \bee \lim_{\vep \to 0} \Pr \{ |
\wh{\bs{p}} - p | \geq \vep \} & = & \lim_{\vep \to 0} \Pr \{ |
\wh{\bs{p}} -
p | \geq \vep, \; \bs{l} = {\ell_\vep} \} =  \lim_{\vep \to 0}
\Pr \{ | \wh{\bs{p}}_{\ell_\vep} - p | \geq \vep \}\\
& = & \lim_{\vep \to 0} \Pr \li \{ | U_{\ell_\vep} | \geq
b_{\ell_\vep} \ri \} = \Pr \{ | U | \geq d \sq{\ka_p} \} = 2 - 2
\Phi ( d \sq{\ka_p} ) \eee} for $p \in (0, 1)$ such that $C_{j_p} =
4 p (1 - p)$ and $j_p = 0$.

Note that,  for a positive number $z$ and a Gaussian random variable
$X$ with zero mean and unit variance, it holds true that {\small
$\Phi (z) = 1 - \Pr \{ X > z \} > 1 - \inf_{t > 0} \bb{E} [ e^{t (X
- z)} ] = 1 - \inf_{t > 0} e^{ - t z + \f{t^2}{2}} = 1 -  e^{ -
\f{z^2}{2} }$}. So, {\small $ \Phi(d) = \Phi \li ( \sq{ 2 \ln
\f{1}{\ze \de} } \ri ) > 1 - \ze \de$} and consequently,
$\liminf_{\vep \to 0} \Pr \{ | \wh{\bs{p}} - p | < \vep \} > 1 - 2
\ze \de$.   This establishes Statement (III).

\subsection{Proof of Theorem  \ref{Unbiased_Bino} } \la{App_Unbiased_Bino}

Let $I_{\wh{\bs{p}}_\ell}$ denote the support of $\wh{\bs{p}}_\ell$
for $\ell = 1, \cd, s$.  Then, {\small \bee \bb{E} | \wh{\bs{p}} - p
|^k & = & \sum_{\ell = 1}^s \sum_{ \wh{p}_\ell \in
I_{\wh{\bs{p}}_\ell} } | \wh{p}_\ell - p |^k \Pr \{
\wh{\bs{p}}_\ell = \wh{p}_\ell, \; \bs{l} = \ell \} \\
& = & \sum_{\ell = 1}^s \li [ \sum_{ \wh{p}_\ell \in
I_{\wh{\bs{p}}_\ell} \atop{ | \wh{p}_\ell - p | <
\f{p}{\sqrt[4]{\ga_\ell}} } } | \wh{p}_\ell - p |^k \Pr \{
\wh{\bs{p}}_\ell = \wh{p}_\ell, \; \bs{l} = \ell \} + \sum_{
\wh{p}_\ell \in I_{\wh{\bs{p}}_\ell} \atop{ | \wh{p}_\ell - p | \geq
\f{p}{\sqrt[4]{\ga_\ell}} } } | \wh{p}_\ell -
p |^k \Pr \{ \wh{\bs{p}}_\ell = \wh{p}_\ell, \; \bs{l} = \ell \} \ri ]  \\
& = & \sum_{\ell = 1}^s  \sum_{ \wh{p}_\ell \in I_{\wh{\bs{p}}_\ell}
\atop{ | \wh{p}_\ell - p | < \f{p}{\sqrt[4]{\ga_\ell}} } } |
\wh{p}_\ell - p|^k \Pr \{ \wh{\bs{p}}_\ell = \wh{p}_\ell, \; \bs{l}
= \ell \} + \sum_{\ell = 1}^s \sum_{ \wh{p}_\ell \in
I_{\wh{\bs{p}}_\ell} \atop{ | \wh{p}_\ell - p | \geq
\f{p}{\sqrt[4]{\ga_\ell}} } } |\wh{p}_\ell - p|^k
\Pr \{ \wh{\bs{p}}_\ell = \wh{p}_\ell, \; \bs{l} = \ell \}\\
& \leq & \sum_{\ell = 1}^s  \li ( \f{p}{\sqrt[4]{\ga_\ell}} \ri )^k
\sum_{ \wh{p}_\ell \in I_{\wh{\bs{p}}_\ell} \atop{ | \wh{p}_\ell - p
| < \f{p}{\sqrt[4]{\ga_\ell}} } }  \Pr \{ \wh{\bs{p}}_\ell =
\wh{p}_\ell, \; \bs{l} = \ell \} + \sum_{\ell = 1}^s \sum_{
\wh{p}_\ell \in I_{\wh{\bs{p}}_\ell} \atop{ | \wh{p}_\ell - p | \geq
\f{p}{\sqrt[4]{\ga_\ell}} } }
\Pr \{ \wh{\bs{p}}_\ell = \wh{p}_\ell \}\\
& = & \sum_{\ell = 1}^s  \li ( \f{p}{\sqrt[4]{\ga_\ell}} \ri )^k \Pr
\li \{ | \wh{\bs{p}}_\ell - p | < \f{p}{\sqrt[4]{\ga_\ell}}, \;
\bs{l} = \ell \ri \} + \sum_{\ell = 1}^s \Pr
\li \{ | \wh{\bs{p}}_\ell  - p | \geq \f{p}{\sqrt[4]{\ga_\ell}}  \ri \}\\
& \leq & \li ( \f{p}{\sqrt[4]{\ga_1}} \ri )^k \sum_{\ell = 1}^s \Pr
\li \{ \bs{l} = \ell \ri \} + \sum_{\ell = 1}^s \Pr
\li \{ | \wh{\bs{p}}_\ell  - p | \geq \f{p}{\sqrt[4]{\ga_\ell}}  \ri \}\\
& = & \li ( \f{p}{\sqrt[4]{\ga_1}} \ri )^k + \sum_{\ell = 1}^s \Pr
\li \{ | \wh{\bs{p}}_\ell  - p | \geq \f{p}{\sqrt[4]{\ga_\ell}}  \ri
\} \leq  \li ( \f{p}{\sqrt[4]{\ga_1}} \ri )^k +  2 \sum_{\ell = 1}^s
\exp \li ( - \f{\sqrt{\ga_\ell} }{8} \ri ) \eee} for $k = 1, 2,
\cd$, where the last inequality  is derived from Corollary 1 of
\cite{Chen3},  which asserts that
\[
\Pr \li \{ | \wh{\bs{p}}_\ell  - p | \geq \vep p \ri \} \leq 2 \exp
\li ( - \ga_\ell \li [ \ln ( 1 + \vep ) - \f{\vep}{1 + \vep} \ri ]
\ri ) < 2 \exp \li ( - \f{ \ga_\ell \; \vep^2 }{8} \ri ), \qqu \ell
= 1, \cd, s
\]
for $\vep \in (0, 1)$.  By the assumption that $\ga_{\ell + 1} -
\ga_\ell \geq 1$ for any $\ell > 0$, we have that \[ \sum_{\ell =
1}^s \exp \li ( - \f{\sqrt{\ga_\ell} }{8} \ri ) \leq \sum_{m =
\ga_1}^\iy \exp \li ( - \f{\sqrt{m} }{8} \ri ) < \int_{\ga_1 -
1}^\iy \exp \li ( - \f{\sqrt{x} }{8} \ri ) d x \to 0
\]
as $\ga_1 \to \iy$.  Hence, {\small $\bb{E} \li [ | \wh{\bs{p}} - p
|^k \ri ]  <  \li ( \f{p}{\sqrt[4]{\ga_1}} \ri )^k +  2 \int_{\ga_1
- 1}^\iy \exp \li ( - \f{\sqrt{x} }{8} \ri ) d x  \to 0$} as $\ga_1
\to \iy$.  Since $\li | \bb{E} [\wh{\bs{p}} - p] \ri | \leq \bb{E} |
\wh{\bs{p}} - p |$, we have that $\bb{E} [\wh{\bs{p}} - p] \to 0$ as
$n_1 \to \iy$. This completes the proof of the theorem.

\subsection{Proof of Theorem \ref{Bino_Rev_CDF} } \la{App_Bino_Rev_CDF}

We need some preliminary results.

\beL \la{decr}
 Let $0 < \vep < 1$.  Then,  {\small
$\mscr{M}_{\mrm{I}}  ( z, \f{z}{1 + \vep}  )$} is monotonically
decreasing with respect to $z \in (0, 1)$. \eeL

\bpf To show that $\mscr{M}_{\mrm{I}} ( z, \f{z}{1 + \vep} )$ is
monotonically decreasing with respect to $z \in (0, 1)$, we derive
the partial derivative as {\small $\f{ \pa  } { \pa z }
\mscr{M}_{\mrm{I}}  ( z, \f{z}{1 + \vep}  )  =  \f{1}{z^2}  [ \ln (
1 - \f{ \vep z } { 1 + \vep - z } ) +  \f{ \vep z } { 1 + \vep - z }
]$},  where the right side is negative if {\small $\ln ( 1 - \f{
\vep z } { 1 + \vep - z } ) < - \f{ \vep z } { 1 + \vep - z }$}.
This condition is  seen to be true by virtue of the standard
inequality $\ln (1 - x) < - x, \; \fa x \in (0, 1)$ and the fact
that {\small $0 < \f{ \vep z } { 1 + \vep - z } < 1$} as a
consequence of $0 < z < 1$.  This completes the proof of the lemma.

 \epf

\beL \la{revcom}

$\mscr{M}_{\mrm{I}} ( z, \f{z}{1 + \vep} )  > \mscr{M}_{\mrm{I}} (
z, \f{z}{1 - \vep} )$ for $0 < z < 1 - \vep < 1$.

\eeL

\bpf

The lemma follows from the facts that $\mscr{M}_{\mrm{I}} ( z,
\f{z}{1 + \vep} )  = \mscr{M}_{\mrm{I}} ( z, \f{z}{1 - \vep} )$ for
$\vep = 0$ and that
\[
\f{ \pa   } { \pa \vep  }  \mscr{M}_{\mrm{I}} \li ( z, \f{z}{1 +
\vep} \ri ) = - \f{ \vep } { 1 + \vep } \f{1}{1 + \vep - z} >  \f{
\pa } { \pa \vep  }  \mscr{M}_{\mrm{I}} \li ( z, \f{z}{1 - \vep} \ri
)  = - \f{ \vep } { 1 - \vep } \f{1}{1 - \vep - z}.
\]

\epf

\beL \la{revDS1} $\{ F_{\wh{\bs{p}}_s} (\wh{\bs{p}}_s,
\f{\wh{\bs{p}}_s}{1 - \vep} ) \leq \ze \de, \; G_{\wh{\bs{p}}_s}
(\wh{\bs{p}}_s, \f{\wh{\bs{p}}_s}{1 + \vep} ) \leq \ze \de \}$ is a
sure event.

 \eeL

 \bpf

 By Lemma \ref{decb},
{\small  \bel \Pr \li \{ G_{\wh{\bs{p}}_s} \li ( \wh{\bs{p}}_s,
\f{\wh{\bs{p}}_s}{1 + \vep} \ri ) \leq \ze \de \ri \} & = &
 \Pr \li \{ 1 - S_{\mrm{B}} \li (\ga_s - 1, \mathbf{n}_s,
\f{\wh{\bs{p}}_s}{1 + \vep} \ri ) \leq \ze \de \ri \}\\
 & \geq & \Pr
\li \{  \mathbf{n}_s \mscr{M}_{\mrm{B}} \li ( \f{ \ga_s } {
\mathbf{n}_s }, \f{\wh{\bs{p}}_s}{1 + \vep} \ri )  \leq \ln (\ze
\de) \ri \}  =  \Pr \li \{  \f{\ga_s}{\wh{\bs{p}}_s}
\mscr{M}_{\mrm{B}} \li ( \wh{\bs{p}}_s, \f{\wh{\bs{p}}_s}{1 + \vep}
\ri ) \leq \ln (\ze \de)
\ri \} \nonumber\\
& = & \Pr \li \{  \mscr{M}_{\mrm{I}} \li ( \wh{\bs{p}}_s,
\f{\wh{\bs{p}}_s}{1 + \vep} \ri ) \leq \f{ \ln (\ze \de) }{\ga_s}
\ri \}. \la{eqrev1} \eel} Making use of Lemma \ref{decr} and the
fact {\small $\lim_{z \to 0} \mscr{M}_{\mrm{I}} ( z, \f{z}{1 + \vep}
) = \f{\vep}{1 + \vep} - \ln(1 + \vep)$},  we have {\small
$\mscr{M}_{\mrm{I}} ( z, \f{z}{1 + \vep} ) < \f{\vep}{1 + \vep} -
\ln(1 + \vep)$} for any $z \in (0, 1]$. Consequently, {\small $ \{
\mscr{M}_{\mrm{I}} ( \wh{\bs{p}}_s, \f{\wh{\bs{p}}_s}{1 + \vep} )
\leq \f{\vep}{1 + \vep} - \ln(1 + \vep) \}$} is a sure event because
$0 < \wh{\bs{p}}_s (\om) \leq 1$ for any $\om \in \Om$.  By the
definition of $\ga_s$, we have {\small \[ \ga_s = \li \lc \f{ \ln
(\ze \de) } { \f{\vep}{1 + \vep} - \ln(1 + \vep) } \ri \rc \geq \f{
\ln (\ze \de)  } { \f{\vep}{1 + \vep} - \ln(1 + \vep) }.
\]}
Since $\f{\vep}{1 + \vep} - \ln(1 + \vep) < 0$ for any $\vep \in (0,
1)$, we have $\f{ \ln (\ze \de)  } { \ga_s } \geq \f{\vep}{1 + \vep}
- \ln(1 + \vep)$.  Hence, {\small \be \la{eqrev2} \Pr \li \{
 \mscr{M}_{\mrm{I}} \li ( \wh{\bs{p}}_s, \f{\wh{\bs{p}}_s}{1 + \vep}
\ri ) \leq \f{ \ln (\ze \de) }{\ga_s}  \ri \} \geq \Pr \li \{
\mscr{M}_{\mrm{I}} \li ( \wh{\bs{p}}_s, \f{\wh{\bs{p}}_s}{1 + \vep}
\ri ) \leq \f{\vep}{1 + \vep} - \ln(1 + \vep) \ri \} = 1. \ee}
Combining (\ref{eqrev1}) and (\ref{eqrev2}) yields {\small $\Pr \{
G_{\wh{\bs{p}}_s} (\wh{\bs{p}}_s, \f{\wh{\bs{p}}_s}{1 + \vep} ) \leq
\ze \de \} = 1$}.

Similarly, by Lemmas \ref{decb} and \ref{revcom}, {\small  \bel
 \Pr \li \{ F_{\wh{\bs{p}}_s} \li (
\wh{\bs{p}}_s, \f{\wh{\bs{p}}_s}{1 - \vep} \ri ) \leq \ze \de \ri \}
& \geq & \Pr \li \{ S_{\mrm{B}} \li (\ga_s, \mathbf{n}_s,
\f{\wh{\bs{p}}_s}{1 - \vep} \ri ) \leq \ze \de \ri \} \geq  \Pr \li
\{  \mathbf{n}_s \mscr{M}_{\mrm{B}} \li ( \f{ \ga_s
} { \mathbf{n}_s }, \f{\wh{\bs{p}}_s}{1 - \vep} \ri )  \leq \ln (\ze \de) \ri \} \nonumber\\
& = & \Pr \li \{  \f{\ga_s}{\wh{\bs{p}}_s} \mscr{M}_{\mrm{B}} \li (
\wh{\bs{p}}_s, \f{\wh{\bs{p}}_s}{1 - \vep} \ri ) \leq \ln (\ze \de)
\ri \}  =  \Pr \li \{  \mscr{M}_{\mrm{I}} \li ( \wh{\bs{p}}_s,
\f{\wh{\bs{p}}_s}{1 - \vep} \ri ) \leq \f{ \ln (\ze \de)
}{\ga_s} \ri \} \nonumber\\
& \geq & \Pr \li \{  \mscr{M}_{\mrm{I}} \li ( \wh{\bs{p}}_s,
\f{\wh{\bs{p}}_s}{1 + \vep} \ri ) \leq \f{ \ln (\ze \de) }{\ga_s}
\ri \} = 1. \la{eqrev3} \eel} This completes the proof of the lemma.

 \epf

\bsk

Now we are in a position to prove Theorem \ref{Bino_Rev_CDF}.
Clearly, $\wh{\bs{p}}_\ell$ is a ULE of $p$ for $\ell = 1, \cd, s$.
Define {\small $\mscr{L} ( \wh{\bs{p}}_\ell ) =
\f{\wh{\bs{p}}_\ell}{1 + \vep}$} and {\small $\mscr{U} (
\wh{\bs{p}}_\ell ) = \f{\wh{\bs{p}}_\ell}{1 - \vep}$} for $\ell = 1,
\cd, s$.  Then, $\{ \mscr{L} ( \wh{\bs{p}}_\ell ) \leq
\wh{\bs{p}}_\ell \leq \mscr{U} ( \wh{\bs{p}}_\ell )$ is a sure event
for $\ell = 1, \cd, s$.  By the definition of the stopping rule, we
have $\{ \bs{D}_\ell = 1 \}  =  \{ F_{\wh{\bs{p}}_\ell}  (
\wh{\bs{p}}_\ell, \mscr{U} ( \wh{\bs{p}}_\ell ) ) \leq \ze \de, \;
G_{\wh{\bs{p}}_\ell}  ( \wh{\bs{p}}_\ell, \mscr{L} (
\wh{\bs{p}}_\ell )  ) \leq \ze \de \}$ for $\ell = 1, \cd, s$. By
Lemma \ref{revDS1}, we have that $\{  \bs{D}_s = 1 \}$ is a sure
event. So, the sampling scheme satisfies all the requirements
described in Theorem \ref{Monotone_second}, from which
(\ref{rev1cdf}) and (\ref{rev2cdf}) of Theorem \ref{Bino_Rev_CDF}
immediately follows. The other results of Theorem \ref{Bino_Rev_CDF}
can be shown by a similar method as that of the proof of Theorem
\ref{Bino_Rev_Chernoff}.

\subsection{Proof of Theorem \ref{Bino_Rev_Chernoff}} \la{App_Bino_Rev_Chernoff}

Let $X_1, X_2, \cd$ be a sequence of i.i.d. Bernoulli random
variables such that $\Pr \{ X_i = 1 \} = 1 - \Pr \{ X_i = 0 \} = p
\in (0, 1)$ for $i = 1, 2, \cd$. Let $\bs{n}$ be the minimum integer
such that $\sum_{i=1}^{\bs{n}} X_i = \ga$ where $\ga$ is a positive
integer.  In the sequel, from Lemmas \ref{lem_chen} to \ref{tail},
we shall be focusing on probabilities associated with
$\f{\ga}{\bs{n}}$.

\beL \la{lem_chen} {\small \be \la{lem_chen_a} \Pr \li \{
\f{\ga}{\bs{n}} \leq z \ri \} \leq \exp \li ( \ga \mscr{M}_{\mrm{I}}
(z, p) \ri ) \qqu \fa z \in (0, p),\ee \be \la{lem_chen_b} \Pr \li
\{ \f{\ga}{\bs{n}} \geq z \ri \} \leq \exp \li ( \ga
\mscr{M}_{\mrm{I}} (z, p) \ri ) \qqu \fa z \in (p, 1). \ee}
 \eeL

\bpf To  show (\ref{lem_chen_a}),  note that {\small $\Pr \li \{
\f{\ga}{\bs{n}} \leq z \ri \} = \Pr \{ \bs{n} \geq m \} = \Pr \{ X_1
+ \cd + X_m \leq \ga \} = \Pr \{ \f{ \sum_{i=1}^m X_i } {m} \leq
\f{\ga}{m} \}$} where $m = \lc \f{\ga}{z} \rc$.  Since $0 < z < p$,
we have $0 < \f{\ga}{m} = \ga \sh  \lc \f{\ga}{z} \rc \leq \ga \sh (
\f{\ga}{z} ) = z < p$, we can apply Lemma \ref{Hoe_Mas} to obtain
{\small $\Pr \li \{ \f{ \sum_{i=1}^m X_i } {m} \leq \f{\ga}{m} \ri
\} \leq \exp \li (m \mscr{M}_{\mrm{B}} \li ( \f{\ga}{m}, p \ri ) \ri
) = \exp \li (\ga \mscr{M}_{\mrm{I}} \li ( \f{\ga}{m}, p \ri ) \ri
)$}.  Noting that $0 < \f{\ga}{m} \leq z < p$ and that
$\mscr{M}_{\mrm{I}} \li ( z, p \ri )$ is monotonically increasing
with respect to $z \in (0, p)$ as can be seen from {\small $\f{ \pa
\mscr{M}_{\mrm{I}} \li ( z, p \ri ) } { \pa z } = \f{1}{z^2} \ln
\f{1 - z}{1 - p}$}, we have $\mscr{M}_{\mrm{I}} \li ( \f{\ga}{m}, p
\ri ) \leq \mscr{M}_{\mrm{I}} \li (z, p \ri )$ and thus {\small $\Pr
\li \{ \f{\ga}{\bs{n}} \leq z \ri \} = \Pr \li \{ \f{ \sum_{i=1}^m
X_i } {m} \leq \f{\ga}{m} \ri \} \leq \exp \li (\ga
\mscr{M}_{\mrm{I}} \li ( z, p \ri ) \ri )$}.

To show (\ref{lem_chen_b}), note that {\small $\Pr \li \{
\f{\ga}{\bs{n}} \geq z \ri \} = \Pr \{ \bs{n} \leq m \} = \Pr \{ X_1
+ \cd + X_m \geq \ga \} = \Pr \li \{ \f{ \sum_{i=1}^m X_i } {m} \geq
\f{\ga}{m} \ri \}$} where $m = \lf \f{\ga}{z} \rf$.  We need to
consider two cases: (i) $m = \ga$; (ii) $m > \ga$. In the case of $m
= \ga$, we have {\small $\Pr \li \{ \f{\ga}{\bs{n}} \geq z \ri \} =
\Pr \{ X_i = 1, \; i = 1, \cd, \ga \} = \prod_{i=1}^\ga \Pr \{ X_i =
1 \} = p^\ga$}.  Since $\mscr{M}_{\mrm{I}} \li ( z, p \ri )$ is
monotonically decreasing with respect to $z \in (p, 1)$ and $\lim_{z
\to 1} \mscr{M}_{\mrm{I}} \li ( z, p \ri ) = \ln p$, we have $\Pr
\li \{ \f{\ga}{\bs{n}} \geq z \ri \} = p^\ga  < \exp \li (\ga
\mscr{M}_{\mrm{I}} \li ( z, p \ri ) \ri )$. In the case of $m
> \ga$, we have $1 > \f{\ga}{m} = \ga \sh \lf \f{\ga}{z} \rf \geq
\ga \sh  ( \f{\ga}{z} ) = z > p$.  Hence, applying Lemma
\ref{Hoe_Mas}, we obtain {\small $\Pr \li \{ \f{ \sum_{i=1}^m X_i }
{m} \geq \f{\ga}{m} \ri \} \leq \exp \li (m \mscr{M}_{\mrm{B}} \li (
\f{\ga}{m}, p \ri ) \ri ) = \exp \li (\ga \mscr{M}_{\mrm{I}} \li (
\f{\ga}{m}, p \ri ) \ri )$}. Noting that $\mscr{M}_{\mrm{I}} \li (
z, p \ri )$ is monotonically decreasing with respect to $z \in (p,
1)$ and that $1 > \f{\ga}{m} \geq z
> p$, we have $\mscr{M}_{\mrm{I}} \li ( \f{\ga}{m}, p \ri ) \leq
\mscr{M}_{\mrm{I}} \li (z, p \ri )$ and thus {\small $\Pr \li \{
\f{\ga}{\bs{n}} \geq z \ri \} = \Pr \li \{  \f{ \sum_{i=1}^m X_i }
{m} \geq \f{\ga}{m} \ri \} \leq \exp \li (\ga \mscr{M}_{\mrm{I}} \li
( z, p \ri ) \ri )$}.

\epf

The following result, stated as Lemma \ref{Mendob},  have recently
been established by Mendo and Hernando \cite{Mendo2}.

\beL \la{Mendob} Let $\ga \geq 3$ and $\mu_1 \geq \f{\ga - 1}{ \ga -
\f{1}{2} - \sq{\ga - \f{1}{2} } }$.  Then,
 {\small $\Pr \{ \f{ \ga - 1} { \bs{n} } > p \mu_1 \} < 1 -
S_{\mrm{P}} (  \ga - 1, \f{\ga - 1}{ \mu_1} )$} for any $p \in (0,
1)$. \eeL

Since $\Pr \{ \f{\ga}{\bs{n}} > (1 + \vep) p \} = \Pr \{ \f{ \ga -
1} { \bs{n} } \geq \f{ \ga - 1} { \ga} (1 + \vep) p \} = \Pr \{ \f{
\ga - 1} { \bs{n} } \geq p \mu_1 \}$ with $\mu_1 = \f{ \ga - 1} {
\ga} (1 + \vep)$, we can rewrite Lemma \ref{Mendob} as follows:

\beL \la{lemM2} Let $0 < \vep < 1$ and $\ga \geq 3$.  Then, {\small
$\Pr \{ \f{\ga}{\bs{n}} > (1 + \vep) p \} < 1 - S_{\mrm{P}} (  \ga -
1, \f{\ga}{ 1 + \vep}  )$} for any $p \in (0, 1)$ provided that
{\small $1 + \vep \geq \f{\ga}{ \ga - \f{1}{2} - \sq{\ga - \f{1}{2}
} }$}.
 \eeL

The following result stated as Lemma \ref{Mendoa} is due to Mendo
and Hernando \cite{Mendo}.

 \beL \la{Mendoa} Let $\ga \geq 3$ and $\mu_2 \geq \f{\ga +
\sq{\ga}}{\ga - 1}$. Then, {\small  $\Pr \{ \f{ \ga - 1} { \bs{n} }
\geq \f{p}{\mu_2} \} > 1 - S_{\mrm{P}} (  \ga - 1, (\ga - 1) \mu_2
)$} for any $p \in (0, 1)$. \eeL

Since $\Pr \{ \f{\ga}{\bs{n}} \geq (1 - \vep) p \} = \Pr \{ \f{ \ga
- 1} { \bs{n} } \geq \f{ \ga - 1} { \ga} (1 - \vep) p \} = \Pr \{
\f{ \ga - 1} { \bs{n} } \geq \f{p}{\mu_2} \}$ with $\mu_2 =
\f{\ga}{(\ga - 1) (1 - \vep)}$, we can rewrite Lemma \ref{Mendoa} as
follows:

\beL \la{lemM}
  Let $0 < \vep < 1$ and $\ga \geq 3$.  Then,
  {\small $\Pr \{ \f{\ga}{\bs{n}} \geq (1 - \vep) p \} > 1 - S_{\mrm{P}} (  \ga - 1,
  \f{\ga}{ 1 - \vep}  )$} for any $p \in (0, 1)$ provided that {\small $\f{1}{1 - \vep} \geq 1
+ \f{1}{ \sq{\ga} }$}. \eeL

\beL \la{tail} Let $0 < \vep < 1$ and $\ga \in \bb{N}$.  Then,
{\small $\Pr \li \{  \li |  \f{\ga}{\bs{n}} - p  \ri |
> \vep p \ri \} < 1 - S_{\mrm{P}} (  \ga - 1,
  \f{\ga}{ 1 + \vep}  ) + S_{\mrm{P}} (  \ga - 1,
  \f{\ga}{ 1 - \vep}  )$} for any $p \in (0, 1)$
provided that {\small $\ga \geq \li [ \li ( 1 + \vep + \sq{ 1 + 4
\vep + \vep^2 } \ri ) \sh (2 \vep) \ri ]^2 + \f{1}{2}$}.

\eeL

\bpf For simplicity of notations, let {\small $h(\vep) = \li [ \li (
1 + \vep + \sq{ 1 + 4 \vep + \vep^2 } \ri ) \sh (2 \vep) \ri ]^2 +
\f{1}{2}$}.

Clearly, $\Pr \li \{  \li | \f{\ga}{\bs{n}} - p \ri |
> \vep p \ri \} = \Pr \{ \f{\ga}{\bs{n}} > (1 + \vep) p \} +
1 - \Pr \{ \f{\ga}{\bs{n}} \geq (1 - \vep) p \}$. By virtue of
Lemmas \ref{lemM2} and \ref{lemM}, to prove that {\small $\Pr \li \{
\li | \f{\ga}{\bs{n}} - p \ri |
> \vep p \ri \} < 1 - S_{\mrm{P}} (  \ga - 1,
  \f{\ga}{ 1 + \vep}  ) + S_{\mrm{P}} (  \ga - 1,
  \f{\ga}{ 1 - \vep}  )$} for any $p \in (0, 1)$
provided that $\ga \geq h(\vep)$, it suffices to prove the following
statements:

(i) {\small $1 + \vep \geq \f{\ga}{ \ga - \f{1}{2} - \sq{\ga -
\f{1}{2} }  } $} implies {\small $\f{1}{1 - \vep} \geq 1 + \f{1}{
\sq{\ga} }$};

(ii) {\small $1 + \vep \geq \f{\ga}{ \ga - \f{1}{2} - \sq{\ga -
\f{1}{2} }  }$} is equivalent to $\ga \geq h(\vep)$;

(iii) $\ga \geq h(\vep)$ implies $\ga \geq 3$.

To prove statement (i), note that {\small \[ \f{1}{1 - \vep} \geq 1
+ \f{1}{ \sq{\ga} } \LRA \vep \geq \f{ 1 } { \sq{\ga} + 1 }, \qqu  1
+ \vep \geq \f{\ga}{ \ga - \f{1}{2} - \sq{\ga - \f{1}{2} }  } \LRA
\vep \geq \f{ \f{1}{2} + \sq{\ga - \f{1}{2} } }{ \ga - \f{1}{2} -
\sq{\ga - \f{1}{2} } }.
\]}
Hence, it suffices to show {\small $\li ( \f{1}{2} + \sq{\ga -
\f{1}{2} } \ri ) \sh \li ( \ga - \f{1}{2} - \sq{\ga - \f{1}{2} } \ri
) > \f{ 1 } { \sq{\ga} + 1 }$}, i.e., {\small $\f{ \ga  } { \f{1}{2}
+ \sq{\ga - \f{1}{2} } } - 2 < \sq{\ga}$}. Let {\small $t = \sq{\ga
- \f{1}{2} }$}. Then, $\ga = t^2 + \f{1}{2}$ and the inequality
becomes {\small
\[ \ga
> \li ( \f{ \ga  } { \f{1}{2} + \sq{\ga - \f{1}{2} } }  - 2 \ri )^2
\LRA t^2 + \f{1}{2} > \li ( \f{ t^2 + \f{1}{2}  } { t + \f{1}{2}
 }  - 2 \ri )^2,
\]}
i.e., {\small $5 t^3 - \f{9}{4} t^2 - \f{3}{2} t - \f{1}{8} > 0$}
 under the condition that {\small $\f{ t^2 + \f{1}{2}  } { t + \f{1}{2}
 }  - 2 > 0 \LRA (t - 1)^2 > \f{3}{2} \LRA t > 1 + \sq{ \f{3}{2}
 }$.}  Clearly, {\small $5 t^3 - \f{9}{4} t^2 - \f{3}{2} t - \f{1}{8} > 5
t^3 - \f{9}{4} t^3 - \f{3}{2} t^3 - \f{1}{8} t^3 = \f{9}{8} t^3 >
0$} for $t
> 1 + \sq{ \f{3}{2} }$.  It follows that, for $t > 1 + \sq{ \f{3}{2}
}$, i.e., $\ga > 5.4$, the inequality holds.  It can be checked by
hand calculation that it also holds for $\ga = 1, \cd, 5$. Hence,
the inequality holds for all $\ga \geq 1$.  This establishes
statement (i).

To show statement (ii), we rewrite {\small $1 + \vep \geq \f{\ga}{
\ga - \f{1}{2} - \sq{\ga - \f{1}{2} }  }$} in terms of {\small $t =
\sq{\ga - \f{1}{2} }$} as {\small $1 + \vep \geq \f{ t^2 + \f{1}{2}
} { t^2 - t  }$}, which is equivalent to $t^2 - (1 + \vep) t -
\f{1}{2} \geq 0$. Solving this inequality yields {\small $t \geq \f{
1 + \vep + \sq{ 1 + 4 \vep +  \vep^2 }  } { 2 \vep  } \LRA \ga \geq
h(\vep)$. } This proves statement (ii).

 To show statement (iii), it is sufficient to show that $h(\vep) \geq 3$ for $\vep \in (0,
 1]$.  Note that {\small $h(\vep) = \f{1}{4} [ 1 + g(\vep)]^2 +
 \f{1}{2}$} with {\small $g(\vep) = (1  + \sq{ 1 + 4 \vep +
\vep^2  } ) \sh \vep$}. Since {\small $g^\prime(\vep) = - ( \sq{ 1 +
4 \vep + \vep^2 }  + 1 + 2 \vep ) \sh ( \vep^2 \sq{ 1 + 4 \vep +
\vep^2 } ) < 0$}, the minimum of $h(\vep)$ is achieved at $\vep =
1$, which is {\small $\li ( 1 + \sq{ \f{3}{2} } \ri )^2 + \f{1}{2} >
3$}. Hence, $\ga \geq  h(\vep)$ implies $\ga \geq 3$.  This proves
statement (iii).

\epf

\beL \la{PCH} Let $\ovl{X}_n = \f{ \sum_{i=1}^n X_i } {n}$ where
$X_1, \cd, X_n$ are i.i.d. Poisson random variables with mean $\lm >
0$. Then, $\Pr \{ \ovl{X}_n \geq z \} \leq \exp (n
\mscr{M}_{\mrm{P}} (z, \lm) )$ for any $z \in (\lm, \iy)$.
Similarly, $\Pr \{ \ovl{X}_n \leq z \} \leq \exp (n
\mscr{M}_{\mrm{P}} (z, \lm) )$ for any $z \in (0, \lm)$.
 \eeL

 \bpf

Let $Y = n \ovl{X}_n$. Then, $Y$ is a Poisson random variable with
mean $\se = n \lm$.  Let $r = n z$.  If $z > \lm$, then $r > \se$
and,  by virtue of Chernoff's bound \cite{Chernoff}, we have \bee
\Pr \{ \ovl{X}_n \geq z \}  = \Pr \{ Y \geq r \} & \leq & \inf_{t
> 0} \bb{E} \li [ e^{ t (Y - r) } \ri ] = \inf_{t
> 0} \sum_{i
= 0}^\iy e^{ t (i - r)} \f{ \se^i } { i! } e^{- \se}\\
& = & \inf_{t > 0} e^{\se e^t } e^{- \se} e^{ - r \; t} \sum_{i =
0}^\iy \f{ (\se e^t)^i } { i! } e^{- \se e^t} =  \inf_{t > 0} e^{-
\se} e^{ \se e^t - r \; t}, \eee where the infimum is achieved at $t
= \ln \li ( \f{r}{\se} \ri ) > 0$.  For this value of $t$, we have
$e^{- \se} e^{ \se e^t - t r } = e^{- \se} \li ( \f{ \se e } { r }
\ri )^r$.  Hence, we have $\Pr \{ \ovl{X}_n \geq z \}  \leq e^{-
\se} \li ( \f{ \se e } { r } \ri )^r = \exp (n \mscr{M}_{\mrm{P}}
(z, \lm) )$.

Similarly, for any number $z \in (0, \lm)$, we have $\Pr \{
\ovl{X}_n \leq z \} \leq \exp (n \mscr{M}_{\mrm{P}} (z, \lm) )$.

 \epf

\beL \la{boundg} $1 - S_{\mrm{P}} (  \ga - 1,
  \f{\ga}{ 1 + \vep}  ) + S_{\mrm{P}} (  \ga - 1,
  \f{\ga}{ 1 - \vep}  ) < 2 \li [ e^{\vep} (1 + \vep)^{-
(1 + \vep)} \ri ]^{\ga \sh (1 + \vep)}$.  \eeL

\bpf

Let $K^+$ be a Poisson random variable with mean value $\f{\ga}{ 1 +
\vep}$. Let $K^-$ be a Poisson random variable with mean value
$\f{\ga}{ 1 - \vep}$.  Then, we have {\small $\Pr \{ K^+ \geq \ga \}
= 1 - S_{\mrm{P}} (  \ga - 1, \f{\ga}{ 1 + \vep}  )$} and {\small
$\Pr \{ K^- < \ga \} = S_{\mrm{P}} (  \ga - 1,
  \f{\ga}{ 1 - \vep}  )$}.  Applying Lemma \ref{PCH}, we have
\[
\Pr \{ K^+ \geq \ga \} \leq \li [ e^{\vep} (1 + \vep)^{- (1 + \vep)}
\ri ]^{\ga \sh (1 + \vep)}, \qqu \Pr \{ K^- < \ga \} \leq \li [
e^{-\vep} (1 - \vep)^{- (1 - \vep)} \ri ]^{\ga \sh (1 - \vep)}.
\]
It follows that {\small \bee 1 - S_{\mrm{P}} \li (  \ga - 1,
  \f{\ga}{ 1 + \vep} \ri ) + S_{\mrm{P}} \li (  \ga - 1,
  \f{\ga}{ 1 - \vep}  \ri ) & = & \Pr \{ K^+ \geq \ga
\} + \Pr \{ K^- < \ga \}\\
& \leq & \li [ e^{\vep} (1 + \vep)^{- (1 + \vep)} \ri ]^{\ga \sh (1
+ \vep)} + \li [ e^{-\vep} (1 - \vep)^{- (1 - \vep)} \ri ]^{\ga \sh
(1 - \vep)}\\
& \leq & 2 \li [ e^{\vep} (1 + \vep)^{- (1 + \vep)} \ri ]^{\ga \sh
(1 + \vep)}. \eee}

\epf

 \beL
 \la{zell}
 For $\ell = 1, \cd, s - 1$, there exists a unique number $z_\ell \in (0, 1]$ such that {\small
$\mscr{M}_{\mrm{I}} ( z_\ell, \f{z_\ell}{1 + \vep}  ) = \f{ \ln (\ze
\de)  } { \ga_\ell}$}.  Moreover,  $z_1 > z_2 > \cd > z_{s - 1}$.
 \eeL

 \bpf

 By the definition of $\ga_\ell$, we have
{\small  \[ \li \lc \f{ \ln (\ze \de)  }  { - \ln(1 + \vep) } \ri
\rc \leq \ga_\ell < \ga_s = \li \lc \f{ \ln (\ze \de)  } {
\f{\vep}{1 + \vep} - \ln(1 + \vep)  } \ri \rc,
 \]}
which implies {\small $\f{ \ln (\ze \de)  }  { - \ln(1 + \vep) }
\leq \ga_\ell <  \f{ \ln (\ze \de)  } { \f{\vep}{1 + \vep} - \ln(1 +
\vep)  }$}.  Making use of this inequality and the fact
 \[
 \lim_{z \to 0} \mscr{M}_{\mrm{I}} \li ( z, \f{z}{1 + \vep} \ri ) =
\f{\vep}{1 + \vep} - \ln(1 + \vep) < 0, \qqu \lim_{z \to 1}
\mscr{M}_{\mrm{I}} \li ( z, \f{z}{1 + \vep} \ri ) = - \ln(1 + \vep)
< 0,
\]
we have
\[
\lim_{z \to 1} \mscr{M}_{\mrm{I}} \li ( z, \f{z}{1 + \vep} \ri )
\leq \f{ \ln (\ze \de)  } { \ga_\ell} < \lim_{z \to 0}
\mscr{M}_{\mrm{I}} \li ( z, \f{z}{1 + \vep} \ri ).
\]
By Lemma \ref{decr}, {\small $\mscr{M}_{\mrm{I}} ( z, \f{z}{1 +
\vep}  )$} is monotonically decreasing with respect to $z \in (0,
1]$. Hence, there exists a unique number $z_\ell \in (0, 1]$ such
that {\small $\mscr{M}_{\mrm{I}}  ( z_\ell, \f{z_\ell}{1 + \vep}  )
= \f{ \ln (\ze \de)  } { \ga_\ell}$}.

To show that $z_\ell$ decreases with respect to $\ell$, we introduce
function  {\small $F(z, \ga) = \ga \mscr{M}_{\mrm{I}} ( z, \f{z}{1 +
\vep} ) - \ln (\ze \de)$}.  Clearly, {\small \[ \f{ d z } { d \ga }
= - \f{ \f{ \pa } { \pa \ga} F(z, \ga) } { \f{ \pa } { \pa z } F(z,
\ga) } = - \f{ \mscr{M}_{\mrm{I}} \li ( z, \f{z}{1 + \vep} \ri ) } {
\ga \f{ \pa  } { \pa z  } \mscr{M}_{\mrm{I}} \li ( z, \f{z}{1 +
\vep} \ri )  }.
\]}
As can be seen from Lemma \ref{decr} and the fact {\small $\lim_{z
\to 0} \mscr{M}_{\mrm{I}} ( z, \f{z}{1 + \vep} ) < 0$}, we have
{\small $\mscr{M}_{\mrm{I}}  ( z, \f{z}{1 + \vep} ) < 0$} and
{\small $\f{ \pa } { \pa z  }  \mscr{M}_{\mrm{I}}  ( z, \f{z}{1 +
\vep} ) < 0$} for any $z \in (0, 1]$. It follows that $\f{ d z } { d
\ga }$ is negative and consequently $z_1 > z_2 > \cd
> z_{s - 1}$. The proof of the lemma is thus completed.

\epf

\beL \la{pass_a} $\li \{ \bs{D}_\ell = 1 \ri \} \subseteq \li \{
\mscr{M}_{\mrm{I}} \li ( \wh{\bs{p}}_\ell, \f{\wh{\bs{p}}_\ell}{1 +
\vep} \ri ) \leq \f{\ln (\ze \de) }{\ga_\ell}, \; \mscr{M}_{\mrm{I}}
\li ( \wh{\bs{p}}_\ell, \f{\wh{\bs{p}}_\ell}{1 - \vep} \ri ) \leq
\f{\ln (\ze \de) }{\ga_\ell} \ri \}$ for $\ell = 1, \cd, s$. \eeL

\bpf  The lemma is a direct consequence of  Lemma \ref{revcom}.

\epf

\beL \la{well_ch}
 $\bs{D}_s = 1$.
\eeL

\bpf

To show $\bs{D}_s = 1$, it suffices to show {\small
$\mscr{M}_{\mrm{I}} (z, \f{z}{1 + \vep}  ) \leq \f{ \ln (\ze \de) }
{ \ga_s }$} for any $z \in (0, 1]$.  This is because {\small $\{
\bs{D}_s = 1 \} =  \{ \mscr{M}_{\mrm{I}} ( \wh{\bs{p}}_s ,
\f{\wh{\bs{p}}_s}{1 + \vep}  ) \leq \f{\ln (\ze \de) }{\ga_s} \}$}
and $0 < \wh{\bs{p}}_s (\om) \leq 1$ for any $\om \in \Om$.  By the
definition of sample sizes, we have {\small $\ga_s = \li \lc \f{ (1
+ \vep) \ln (\ze \de)}{ \vep - (1 + \vep) \ln (1 + \vep) } \ri \rc
\geq \f{ (1 + \vep) \ln (\ze \de)}{ \vep - (1 + \vep) \ln (1 + \vep)
}$}. Since {\small $\lim_{z \to 0} \mscr{M}_{\mrm{I}} (z, \f{z}{1 +
\vep} )   = \f{\vep}{1 + \vep} - \ln (1 + \vep) < 0$}, we have
{\small $\lim_{z \to 0} \mscr{M}_{\mrm{I}} (z, \f{z}{1 + \vep} )
\leq \f{ \ln (\ze \de)  } { \ga_s }$}. By Lemma \ref{decr}, we have
that {\small $\mscr{M}_{\mrm{I}} (z, \f{z}{1 + \vep}  )$} is
monotonically decreasing with respect to $z \in (0, 1)$. Hence,
{\small $\mscr{M}_{\mrm{I}} (z, \f{z}{1 + \vep} ) < \lim_{z \to 0}
\mscr{M}_{\mrm{I}} (z, \f{z}{1 + \vep} ) \leq \f{ \ln (\ze \de) } {
\ga_s }$} for any $z \in (0, 1)$. Since {\small $\mscr{M}_{\mrm{I}}
 (z, \f{z}{1 + \vep} )$} is a continuous function with respect
to $z \in (0, 1)$ and {\small $\mscr{M}_{\mrm{I}} (1, \f{1}{1 +
\vep} ) = \lim_{z \to 1} \mscr{M}_{\mrm{I}} (z, \f{z}{1 + \vep} )$},
it must be true that {\small $\mscr{M}_{\mrm{I}}  (1, \f{1}{1 +
\vep} ) \leq \f{ \ln (\ze \de)  } { \ga_s }$}. This completes the
proof of the lemma.

\epf

\beL

\la{EvalD}

 {\small $\li \{
\bs{D}_\ell = 1  \ri \} =  \Pr \{  \wh{\bs{p}}_\ell \geq z_\ell \}$
} for $\ell = 1, \cd, s - 1$.

\eeL

\bpf

By Lemma \ref{zell}, for $\ell = 1, \cd, s - 1$, there exists a
unique number $z_\ell \in (0, 1]$ such that $\mscr{M}_{\mrm{I}}  (
z_\ell, \f{ z_\ell } { 1 + \vep }  ) = \f{ \ln (\ze \de) } {
\ga_\ell }$.  From Lemma \ref{decr}, we know that
$\mscr{M}_{\mrm{I}} ( z, \f{ z} { 1 + \vep }  )$ is monotonically
decreasing with respect to $z \in (0, 1)$.  It follows that
$\mscr{M}_{\mrm{I}} ( z, \f{ z} { 1 + \vep }  ) \leq \f{ \ln (\ze
\de) } { \ga_\ell }$ if and only if $z \geq z_\ell$. This implies
that {\small $\li \{ \bs{D}_\ell = 1 \ri \} =  \{ \mscr{M}_{\mrm{I}}
( \wh{\bs{p}}_\ell, \f{ \wh{\bs{p}}_\ell } { 1 + \vep } ) \leq \f{
\ln (\ze \de) } { \ga_\ell } \} =  \Pr \{ \wh{\bs{p}}_\ell \geq
z_\ell \}$ } for $\ell = 1, \cd, s - 1$. The lemma is thus proved.

\epf

\beL \la{lemtruncate}

If $\ze$ is sufficiently small, then $ 1 - S_{\mrm{P}} ( \ga_s - 1,
\f{\ga_s}{ 1 + \vep}  ) + S_{\mrm{P}} (  \ga_s - 1,
  \f{\ga_s}{ 1 - \vep}  )  < \de$, inequality
(\ref{cona}) is satisfied and {\small $\Pr \li \{ \li | \f{
\wh{\bs{p}} - p } { p } \ri | \leq \vep  \ri \} \geq 1 - \de$} for
any $p \in (0, p^*]$.

\eeL

\bpf

It is obvious that inequality (\ref{cona}) is satisfied if $\ze$ is
sufficiently small.  By Lemma \ref{boundg}, we have {\small $1 -
S_{\mrm{P}} (  \ga_s - 1, \f{\ga_s}{ 1 + \vep}  ) + S_{\mrm{P}} (
\ga_s - 1, \f{\ga_s}{ 1 - \vep}  ) < 2 \li [ e^{\vep} (1 + \vep)^{-
(1 + \vep)} \ri ]^{\ga_s \sh (1 + \vep)}$}.  By the definition of
$\ga_s$, we have {\small $\ga_s = \li \lc \f{ (1 + \vep) \ln (\ze
\de) } { \vep - (1 + \vep) \ln (1 + \vep) } \ri \rc \geq \f{ (1 +
\vep) \ln (\ze \de) } { \vep - (1 + \vep) \ln (1 + \vep) }$}, which
implies $1 - S_{\mrm{P}} (  \ga_s - 1, \f{\ga_s}{ 1 + \vep}  ) +
S_{\mrm{P}} (  \ga_s - 1,
  \f{\ga_s}{ 1 - \vep}  )  < 2 \li [ e^{\vep} (1 + \vep)^{-
(1 + \vep)} \ri ]^{\ga_s \sh (1 + \vep)} \leq 2 \ze \de$. It follows
that $1 - S_{\mrm{P}} (  \ga_s - 1, \f{\ga_s}{ 1 + \vep}  ) +
S_{\mrm{P}} (  \ga_s - 1, \f{\ga_s}{ 1 - \vep}  ) < \de$ if $\ze$ is
sufficiently small. From now on and throughout the proof of the
lemma, we assume that $\ze$ is small enough to guarantee $1 -
S_{\mrm{P}} (  \ga_s - 1, \f{\ga_s}{ 1 + \vep}  ) + S_{\mrm{P}} (
\ga_s - 1, \f{\ga_s}{ 1 - \vep}  )  < \de$ and inequality
(\ref{cona}). Applying Lemma \ref{EvalD} and (\ref{lem_chen_b}) of
Lemma \ref{lem_chen}, we have {\small \be \Pr \li \{ \li | \f{
\wh{\bs{p}} - p } { p } \ri |
> \vep, \; \bs{l} = \ell \ri \}  \leq  \Pr \li \{ \bs{l} =
\ell \ri \} \leq \Pr \li \{ \bs{D}_\ell = 1 \ri \} = \Pr \{
\wh{\bs{p}}_\ell \geq z_\ell \} \leq   \exp ( \ga_\ell
\mscr{M}_{\mrm{I}} ( z_\ell, p ) ) \la{chenbound} \ee} for $0 < p <
z_{s - 1} $ and $\ell = 1, \cd, s - 1$.   On the other hand, noting
that {\small
\[ \Pr \li \{ \li | \f{ \wh{\bs{p}} - p } { p } \ri |
> \vep, \; \bs{l} = s \ri \} = \Pr \li \{ \li | \f{
\f{\ga_s}{\mathbf{n}_s} - p } { p } \ri | > \vep, \; \bs{l} = s \ri
\} \leq \Pr \li \{ \li | \f{ \f{\ga_s}{\mathbf{n}_s} - p } { p } \ri
| > \vep \ri \}
\]}
and that {\small $\ga_s \geq \li [ \li ( 1 + \vep + \sq{ 1 + 4 \vep
+ \vep^2 } \ri ) \sh (2 \vep) \ri ]^2 + \f{1}{2}$} as a consequence
of (\ref{cona}) and the definition of $\ga_s$, we can apply Lemma
\ref{tail} to obtain {\small \be \la{bounds} \Pr \li \{ \li | \f{
\wh{\bs{p}} - p } { p } \ri |
> \vep, \; \bs{l} = s \ri \} < 1 - S_{\mrm{P}} \li (
\ga_s - 1, \f{\ga_s}{ 1 + \vep}  \ri ) + S_{\mrm{P}} \li (  \ga_s -
1, \f{\ga_s}{ 1 - \vep}  \ri ) < \de. \ee} Noting that {\small
$\f{\pa \mscr{M}_{\mrm{I}} ( z, p ) } { \pa p} = \f{z - p}{z p (1 -
p)} > 0$} for any $p \in (0, z)$ and that $\lim_{p \to 0}
\mscr{M}_{\mrm{I}} ( z, p ) = - \iy$, we have that $\sum_{\ell =
1}^{s - 1} \exp ( \ga_\ell \mscr{M}_{\mrm{I}} ( z_\ell, p ) )$
decreases monotonically to $0$ as $p$ decreases from $z_{s - 1}$ to
$0$.  Since $1 - S_{\mrm{P}} (  \ga_s - 1, \f{\ga_s}{ 1 + \vep}  ) +
S_{\mrm{P}} (  \ga_s - 1,
  \f{\ga_s}{ 1 - \vep}  )  < \de$, there exists a unique number
$p^* \in (0,z_{s-1})$ such that $1 - S_{\mrm{P}} (  \ga_s - 1,
\f{\ga_s}{ 1 + \vep}  ) + S_{\mrm{P}} (  \ga_s - 1,
  \f{\ga_s}{ 1 - \vep}  )  + \sum_{\ell = 1}^{s - 1} \exp (
\ga_\ell \mscr{M}_{\mrm{I}} ( z_\ell, p^* ) ) = \de$. It follows
that $1 - S_{\mrm{P}} (  \ga_s - 1, \f{\ga_s}{ 1 + \vep}  ) +
S_{\mrm{P}} (  \ga_s - 1,
  \f{\ga_s}{ 1 - \vep}  )  + \sum_{\ell = 1}^{s - 1} \exp ( \ga_\ell \mscr{M}_{\mrm{I}}
( z_\ell, p^* ) ) \leq \de$ for any $p \in (0, p^*]$.  Combining
(\ref{chenbound}) and (\ref{bounds}), we have {\small $\Pr \li \{
\li | \wh{\bs{p}} - p  \ri | > \vep p \ri \} < 1 - S_{\mrm{P}} (
\ga_s - 1, \f{\ga_s}{ 1 + \vep}  ) + S_{\mrm{P}} (  \ga_s - 1,
\f{\ga_s}{ 1 - \vep}  )  + \sum_{\ell = 1}^{s - 1} \exp ( \ga_\ell
\mscr{M}_{\mrm{I}} ( z_\ell, p ) ) \leq \de$} for any $p \in (0,
p^*]$.  This completes the proof of the lemma. \epf

\bsk

We are now in a position to prove Theorem \ref{Bino_Rev_Chernoff}.
As pointed out after Theorem  \ref{Unbiased_Bino}, there exists an
inherent connection between the multistage inverse sampling scheme
for Bernoulli random variable X and a multistage sampling scheme of
sample sizes $\ga_1 < \ga_2 < \cd < \ga_s$ for a geometric random
variable $Y$ with mean value $\se = \f{1}{p}$.  It can be shown that
$\mcal{F} (z, \se ) = \inf_{t < 0 } \bb{E} [e^{ t ( Y - z) } ] =
\exp ( \mscr{M}_{\mrm{I}} ( \f{1}{z}, \f{1}{\se} ) )$ for $z \leq
\se$ and that $\mcal{G} (z, \se ) = \inf_{t > 0 } \bb{E} [e^{ t ( Y
- z) } ] = \exp ( \mscr{M}_{\mrm{I}} ( \f{1}{z}, \f{1}{\se} ) )$ for
$z \geq \se$. Moreover, {\small $\wh{\bs{\se}}_\ell =
\f{1}{\wh{\bs{p}}_\ell} = \f{ \mbf{n}_\ell } { \sum_{i =
1}^{\mbf{n}_\ell} X_i } = \f{ \mbf{n}_\ell } {\ga_\ell} = \f{
\sum_{i = 1}^{\ga_\ell} Y_i } {\ga_\ell} $} is a ULE of $\se =
\f{1}{p}$ for $\ell = 1, \cd, s$, where $Y_1, Y_2, \cd$ are i.i.d.
samples of $Y$.  Define a random interval with lower limit $\mscr{L}
( \wh{\bs{\se}}_\ell ) = (1 - \vep) \wh{\bs{\se}}_\ell$ and upper
limit $\mscr{U} ( \wh{\bs{\se}}_\ell ) = (1 + \vep)
\wh{\bs{\se}}_\ell$ for $\ell = 1, \cd, s$.  Then, $\{ \mscr{L} (
\wh{\bs{\se}}_\ell ) \leq \wh{\bs{\se}}_\ell \leq \mscr{U} (
\wh{\bs{\se}}_\ell ) \}$ is a sure event for $\ell = 1, \cd, s$.  By
virtue of these facts and Lemmas \ref{pass_a} and \ref{well_ch}, we
have that the sampling scheme satisfies requirements (i) -- (v)
described in Corollary \ref{Monotone_third}, from which
(\ref{rev1ch}) and (\ref{rev2ch}) follow immediately.  By Lemma
\ref{lemtruncate}, there exists a positive number $\ze_0$ such that
$1 - S_{\mrm{P}} ( \ga_s - 1, \f{\ga_s}{ 1 + \vep} ) + S_{\mrm{P}} (
\ga_s - 1, \f{\ga_s}{ 1 - \vep}  )  < \de$, inequality (\ref{cona})
is satisfied and {\small $\Pr \{ | \wh{\bs{p}} - p | \leq \vep  p
\mid p \} \geq 1 - \de$} for any $p \in (0, p^*]$ if $0 < \ze <
\ze_0$. Hence, by restricting $\ze$ to be less than $\ze_0$, we can
guarantee $\Pr \{ | \wh{\bs{p}} - p | \leq \vep  p  \mid p \} \geq 1
- \de$ for any $p \in (0, 1)$ by ensuring $\Pr \{ | \wh{\bs{p}} - p
| \leq \vep p \mid p \} \geq 1 - \de$ for any $p \in [p^*, 1)$. This
completes the proof of Theorem \ref{Bino_Rev_Chernoff}.

\subsection{Proof of Theorem \ref{Bino_Rev_Inverse_Massart} } \la{App_Bino_Rev_Inverse_Massart}

We need some preliminary results.

\beL \la{lem199}
 $\li \{  \bs{D}_\ell = 1 \ri \} \subseteq \li \{
\mscr{M}_{\mrm{I}} \li ( \wh{\bs{p}}_\ell, \f{\wh{\bs{p}}_\ell}{1 +
\vep} \ri ) \leq \f{\ln (\ze \de) }{\ga_\ell}, \; \mscr{M}_{\mrm{I}}
\li ( \wh{\bs{p}}_\ell, \f{\wh{\bs{p}}_\ell}{1 - \vep} \ri ) \leq
\f{\ln (\ze \de) }{\ga_\ell}  \ri \}$ for $\ell = 1, \cd, s$. \eeL

\bpf  For simplicity of notations, define {\small
$\mcal{M}_{\mrm{I}} (z, \mu) = \f{ \mscr{M}(z, \mu) }{z}$}.  By
tedious computation, we can show that $\{  \bs{D}_\ell = 1  \} = \{
\mcal{M}_{\mrm{I}} ( \wh{\bs{p}}_\ell, \f{\wh{\bs{p}}_\ell}{1 +
\vep} ) \leq \f{\ln (\ze \de) }{\ga_\ell} \} $ for $\ell = 1, \cd,
s$. Noting that {\small \[ \mcal{M}_{\mrm{I}} \li ( z, \f{z}{1 +
\vep} \ri ) - \mcal{M}_{\mrm{I}} \li ( z, \f{z}{1 - \vep} \ri ) =
\f{ 2 \vep^3 (2 - z) } {3 \li ( 1 + \f{\vep}{3} \ri ) \li [ 1 - z +
\vep \li ( 1 - \f{ z}{3} \ri ) \ri ] \li ( 1 - \f{\vep}{3} \ri ) \li
[ 1 - z - \vep \li ( 1 - \f{ z}{3} \ri ) \ri ]} > 0
\]}
for $0 < z < 1 - \vep$,  we have  \bee \li \{  \bs{D}_\ell = 1 \ri
\} & = & \li \{ \mcal{M}_{\mrm{I}} \li ( \wh{\bs{p}}_\ell,
\f{\wh{\bs{p}}_\ell}{1 + \vep} \ri ) \leq \f{\ln (\ze \de)
}{\ga_\ell}, \; \mcal{M}_{\mrm{I}} \li ( \wh{\bs{p}}_\ell,
\f{\wh{\bs{p}}_\ell}{1 - \vep} \ri ) \leq \f{\ln (\ze \de)
}{\ga_\ell}  \ri \}\\
& \subseteq & \li \{ \mscr{M}_{\mrm{I}} \li ( \wh{\bs{p}}_\ell,
\f{\wh{\bs{p}}_\ell}{1 + \vep} \ri ) \leq \f{\ln (\ze \de)
}{\ga_\ell}, \; \mscr{M}_{\mrm{I}} \li ( \wh{\bs{p}}_\ell,
\f{\wh{\bs{p}}_\ell}{1 - \vep} \ri ) \leq \f{\ln (\ze \de)
}{\ga_\ell}  \ri \} \eee for $\ell = 1, \cd, s$.  This completes the
proof of the lemma.

\epf

\beL \la{well_mas}
 $\bs{D}_s = 1$.
\eeL

\bpf To show $\bs{D}_s = 1$, it suffices to show {\small
$\mcal{M}_{\mrm{I}} (z, \f{z}{1 + \vep} ) \leq \f{ \ln (\ze \de) } {
\ga_s }$} for any $z \in (0, 1]$.  This is because $0 <
\wh{\bs{p}}_s (\om) \leq 1$ for any $\om \in \Om$ and {\small $\{
\bs{D}_s = 1 \} =  \{ \mcal{M}_{\mrm{I}} ( \wh{\bs{p}}_s ,
\f{\wh{\bs{p}}_s}{1 + \vep}  ) \leq \f{\ln (\ze \de) }{\ga_s} \}$}
as asserted by Lemma \ref{lem199}.

By the definition of sample sizes, we have {\small $\ga_s = \li \lc
2 \li ( 1 + \f{\vep}{3} \ri ) ( 1 + \vep ) \f{\ln \f{1}{\ze \de} }{
\vep^2 } \ri \rc \geq 2 \li ( 1 + \f{\vep}{3} \ri ) ( 1 + \vep ) \f{
\ln \f{1}{\ze \de} }{ \vep^2 }$}. Since {\small $\lim_{z \to 0}
\mcal{M}_{\mrm{I}} (z, \f{z}{1 + \vep}  )   = -  \vep^2 \li [ 2 \li
( 1 + \f{\vep}{3} \ri ) ( 1 + \vep ) \ri ]^{-1}  < 0$}, we have
{\small $\lim_{z \to 0} \mcal{M}_{\mrm{I}} (z, \f{z}{1 + \vep}  )
\leq \f{ \ln (\ze \de)  } { \ga_s }$}.

Note that {\small $\mcal{M}_{\mrm{I}} ( z, \f{z}{1 + \vep} ) = - \f{
\vep^2}{2 \li ( 1 + \f{\vep}{3} \ri ) \li [ 1 + \vep - (1 -
\f{\vep}{3}) z \ri ] }$}, from which it can be seen that {\small
$\mcal{M}_{\mrm{I}}  (z, \f{z}{1 + \vep}  )$} is monotonically
decreasing with respect to $z \in (0, 1)$. Hence, {\small
$\mcal{M}_{\mrm{I}} (z, \f{z}{1 + \vep} ) < \lim_{z \to 0}
\mcal{M}_{\mrm{I}} (z, \f{z}{1 + \vep} ) \leq \f{ \ln (\ze \de) } {
\ga_s }$} for any $z \in (0, 1)$. Since {\small $\mcal{M}_{\mrm{I}}
(z, \f{z}{1 + \vep} )$} is a continuous function with respect to $z
\in (0, 1)$ and {\small $\mcal{M}_{\mrm{I}} (1, \f{1}{1 + \vep} ) =
\lim_{z \to 1} \mcal{M}_{\mrm{I}} (z, \f{z}{1 + \vep} )$}, it must
be true that {\small $\mcal{M}_{\mrm{I}} (1, \f{1}{1 + \vep} ) \leq
\f{ \ln (\ze \de)  } { \ga_s }$}. This completes the proof of the
lemma.

\epf

\bsk

Finally,  by virtue of the above preliminary results and a similar
method as that of Theorem \ref{Bino_Rev_Chernoff}, we can establish
Theorem \ref{Bino_Rev_Inverse_Massart}.

\subsection{Proof of Theorem \ref{ASN_Bino_Inverse}}
\la{App_ASN_Bino_Inverse}

Since $ \Pr \{ \mathbf{n} \geq i \}$ depends only on $X_1, \cd,
X_{i}$ for all $i \geq 1$, we have, by Wald's equation, $\bb{E} [
X_1 + \cd + X_{\mathbf{n}} ] = \bb{E} [ X_i ] \; \bb{E} [
\mathbf{n}] = p \; \bb{E} [ \mathbf{n}]$. By the definition of the
sampling scheme, $X_1 + \cd + X_{\mathbf{n}} = \bs{\ga}$, and it
follows that $\bb{E} [ X_1 + \cd + X_{\mathbf{n}} ] = \bs{\ga}$.
Hence, $p \; \bb{E} [ \mathbf{n}] = \bb{E} [ \bs{\ga}]$, leading to
the first identity.

The second identity is shown as follows.  Let $\bs{l}$ be the index
of stage when the sampling is stopped. Then,  setting $\ga_0 = 0$,
we have {\small \bee \sum_{i = 1}^s ( \ga_i - \ga_{i-1}) \Pr \{
\bs{l} \geq i \} & = & \sum_{i = 1}^s \ga_i \Pr \{ \bs{l} \geq i \}
- \sum_{i = 1}^s
\ga_{i-1} \Pr \{ \bs{l} \geq i \}\\
& = & \sum_{i = 1}^s \ga_i \Pr \{ \bs{l} \geq i \}  - \sum_{j =
0}^{s-1} \ga_{j} \Pr \{ \bs{l} \geq j \} +
 \sum_{j = 0}^{s - 1}
\ga_j \Pr \{ \bs{l} = j \}\\
& = & \ga_s \Pr \{ \bs{l} \geq s \} + \sum_{j = 0}^{s - 1} \ga_j \Pr
\{ \bs{l} = j \} =  \sum_{i = 1}^{s} \ga_i \Pr \{ \bs{l} = i \} =
\bb{E}[\ga_{\bs{l}}] = \bb{E}[\bs{\ga}]. \eee} This completes the
proof of Theorem \ref{ASN_Bino_Inverse}.

\subsection{Proof of Theorem \ref{Bino_Rev_noninverse_Chernoff} }
\la{App_Bino_Rev_noninverse_Chernoff}

We need the following lemma.

\beL

\la{lemm23} $\mscr{M}_{\mrm{B}}(z, \f{z}{1 + \vep} )$ is
monotonically decreasing from $0$ to $\ln \f{1}{1 + \vep}$ as $z$
increases from $0$ to $1$.

\eeL

\bpf

The lemma can be established by verifying that {\small \[ \lim_{z
\to 0} \mscr{M}_{\mrm{B}} \li (z, \f{z}{1 + \vep} \ri ) = 0, \qu
\lim_{z \to 1} \mscr{M}_{\mrm{B}} \li (z, \f{z}{1 + \vep} \ri ) =
\ln \f{1}{1 + \vep}, \qu  \lim_{z \to 0}  \f{\pa  } {\pa z  }
\mscr{M}_{\mrm{B}} \li (z, \f{z}{1 + \vep} \ri ) = \ln \f{1}{1 +
\vep} + \f{\vep}{1 + \vep} < 0
\]}
and {\small $ \f{\pa^2  } {\pa z^2 } \mscr{M}_{\mrm{B}} \li (z,
\f{z}{1 + \vep} \ri )  = \f{ \vep^2 } { (z - 1) (1 + \vep - z)^2 } <
0$} for any $z \in (0,1)$.

 \epf

\subsubsection{Proof of Statement (I)}

Let $0 < \eta < 1$ and $r = \inf_{\ell > 0} \f{n_{\ell +
1}}{n_\ell}$. By the assumption that $r > 1$, we have that there
exists a number $\ell^\prime > \max \{ \tau, \tau + \f{2}{r - 1} +
\f{ \ln (\ze \de) } { \ln 2 } \}$ such that $\f{n_{\ell +
1}}{n_\ell}
> \f{r + 1}{2}$ for any $\ell > \ell^\prime$.   Noting that $\f{ \ln (
\ze \de_\ell ) } { n_\ell }$ is negative for any $\ell > 0$ and that
{\small
\[ \f{ \f{ \ln (\ze \de_{\ell + 1}) } { n_{\ell+1} } } { \f{ \ln (
\ze \de_\ell ) } { n_\ell } } < \f{2}{r + 1} \times \f{ (\ell + 1 -
\tau) \ln 2  - \ln (\ze \de)  } { (\ell - \tau) \ln 2 - \ln (\ze
\de) } = \f{2}{r + 1} \times \li ( 1 + \f{1}{ \ell - \tau - \f{\ln
(\ze \de ) } { \ln 2 } } \ri ) < 1
\]}
for {\small $\ell > \ell^\prime$}, we have that $\f{ \ln ( \ze
\de_\ell ) } { n_\ell }$ is monotonically increasing with respect to
$\ell$ greater than $\ell^\prime$.  In view of such monotonicity and
the fact that {\small $\f{ \ln ( \ze \de_\ell ) } { n_\ell } = \f{
\ln \li ( \ze \de 2^{\tau - \ell} \ri ) } { n_\ell  } \to 0 >
\mscr{M}_{\mrm{B}} ( \eta p, \f{\eta p}{1 + \vep} )$} as $\ell \to
\iy$, we have that there exists an integer $\ka$ greater than
{\small $\ell^\prime$} such that {\small $\mscr{M}_{\mrm{B}} ( \eta
p, \f{\eta p}{1 + \vep} ) < \f{ \ln (\ze \de_\ell) } { n_\ell }$}
for all $\ell \geq \ka$. For $\ell$ no less than such $\ka$, we
claim that $z < \eta p$ if {\small $\mscr{M}_{\mrm{B}} ( z, \f{z}{1
+ \vep} ) > \f{\ln (\ze \de_\ell )}{n_\ell}$} and $z \in [0, 1]$. To
prove this claim, suppose, to get a contradiction, that $z \geq \eta
p$. Then, since {\small $\mscr{M}_{\mrm{B}} ( z, \f{z}{1 + \vep} )$}
is monotonically decreasing with respect to $z \in (0,1)$ as
asserted by Lemma \ref{lemm23}, we have {\small $\mscr{M}_{\mrm{B}}
(z, \f{z}{1 + \vep}  ) \leq \mscr{M}_{\mrm{B}}  ( \eta p, \f{ \eta
p}{1 + \vep}  ) < \f{\ln (\ze \de_\ell)}{n_\ell}$}, which is a
contradiction. Therefore, we have shown the claim and it follows
that {\small $ \{ \mscr{M}_{\mrm{B}} ( \f{K_\ell}{n_\ell},
\f{K_\ell}{(1 + \vep) n_\ell }  ) > \f{\ln (\ze \de_\ell) }{n_\ell}
\} \subseteq \{ K_\ell < \eta p n_\ell \}$} for $\ell \geq \ka$. So,
{\small
\[ \Pr \{ \bs{l}
> \ell \}  \leq  \Pr \li \{ \mscr{M}_{\mrm{B}}
\li ( \f{K_\ell}{n_\ell}, \f{K_\ell}{(1 + \vep) n_\ell } \ri ) >
\f{\ln (\ze \de_\ell) }{n_\ell}  \ri \} \leq  \Pr \{ K_\ell < \eta p
n_\ell \} < \exp \li ( - \f{(1 - \eta)^2 p n_\ell}{2} \ri )
\]} for large enough $\ell$,  where the last inequality is due to the multiplicative Chernoff bound
\cite{Chernoff_tour}.  Since $\Pr \{ \bs{l}
> \ell \}  < \exp ( - \f{(1 - \eta)^2 p n_\ell}{2} )$  for sufficiently large $\ell$ and
$n_\ell \to \iy$ as $\ell \to \iy$, we have $\Pr \{ \bs{l} < \iy \}
= 1$ or equivalently, $\Pr \{ \mbf{n} < \iy \} = 1$.  This completes
the proof of statement (I).

\subsubsection{Proof of Statement (II)}  In the course of proving
Statement (I), we have shown that there exists an integer $\ka$ such
that {\small $\Pr \{ \bs{l} > \ell \} < \exp  ( - c n_\ell  )$} for
any $\ell \geq \ka$, where $c = \f{(1 - \eta)^2 p }{2}$. Note that
\bee \bb{E} [ \mbf{n}] & = &  n_1 + \sum_{\ell = 1}^{\ka} (n_{\ell +
1} - n_\ell) \Pr \{ \bs{l} > \ell \} + \sum_{\ell = \ka + 1}^\iy
(n_{\ell + 1} - n_\ell) \Pr \{ \bs{l} > \ell \}. \eee  Let $R =
\sup_{\ell > 0} \f{n_{\ell + 1}}{n_\ell}$.  Then, $n_{\ell + 1} -
n_{\ell} \leq R n_{\ell}$.  Hence, if we choose $\ka$ large enough
such that $c n_1 r^{\ka}
> 1$, then {\small \bee \sum_{\ell = \ka + 1}^\iy (n_{\ell + 1} - n_{\ell}) \; \Pr \{
\bs{l}
> \ell \} & < &  \sum_{\ell = \ka + 1}^\iy (n_{\ell + 1} -
n_{\ell}) \; e^{- c n_\ell}
 \leq \f{R}{c} \sum_{\ell = \ka + 1}^\iy  c n_{\ell} \; e^{- c n_\ell}
  \leq   \f{R}{c} \; \sum_{\ell = \ka}^\iy c n_1 r^{\ell}
\exp ( - c n_1 r^{\ell} )\\
& <  & \f{R}{c} \int_{\ka - 1}^\iy c n_1 r^{\ell} \exp ( - c n_1
r^{\ell} ) d \ell = \f{R}{c} \f{\exp ( - c n_1 r^{\ka - 1} ) }{\ln
r}, \eee} which implies that $\bb{E} [ \mbf{n}] < \iy$.

\subsubsection{Proof of Statement (III)}

By differentiation with respect to $\vep \in (0, 1)$, we can show
that $\mscr{M}_{\mrm{B}} ( z, \f{z}{1 - \vep} ) < \mscr{M}_{\mrm{B}}
( z, \f{z}{1 + \vep} )$ for $0 \leq z < 1 - \vep$.  It follows that
$\{ \bs{D}_\ell = 1 \} = \{ \mscr{M}_{\mrm{B}} ( \wh{\bs{p}}_\ell,
\f{\wh{\bs{p}}_\ell}{1 + \vep} ) \leq \f{\ln (\ze \de_\ell)
}{n_\ell} \} = \{ \mscr{M}_{\mrm{B}} ( \wh{\bs{p}}_\ell,
\f{\wh{\bs{p}}_\ell}{1 + \vep} ) \leq \f{\ln (\ze \de_\ell)
}{n_\ell}, \; \mscr{M}_{\mrm{B}} ( \wh{\bs{p}}_\ell,
\f{\wh{\bs{p}}_\ell}{1 - \vep} ) \leq \f{\ln (\ze \de_\ell)
}{n_\ell} \}$ for $\ell = 1, \cd, s$.  Hence, by the definition of
the sampling scheme, we have \bee \Pr \li \{ \li | \wh{\bs{p}}_\ell
- p \ri | \geq \vep p, \; \bs{l} = \ell \mid p \ri \} & \leq &
 \Pr \li \{ p \leq
\f{\wh{\bs{p}}_\ell}{1 + \vep}, \;  \mscr{M}_{\mrm{B}} \li (
\wh{\bs{p}}_\ell, \f{\wh{\bs{p}}_\ell}{1 + \vep} \ri  ) \leq
\f{\ln (\ze \de_\ell) }{n_\ell} \mid p \ri \}\\
&  & +  \Pr \li \{ p \geq \f{\wh{\bs{p}}_\ell}{1 - \vep}, \;
\mscr{M}_{\mrm{B}} \li ( \wh{\bs{p}}_\ell, \f{\wh{\bs{p}}_\ell}{1 -
\vep} \ri  ) \leq \f{\ln
(\ze \de_\ell) }{n_\ell}  \mid p \ri \}\\
& \leq &  \Pr \li \{ p \leq \f{\wh{\bs{p}}_\ell}{1 + \vep}, \;
\mscr{M}_{\mrm{B}} \li ( \wh{\bs{p}}_\ell, p  \ri  ) \leq
\f{\ln (\ze \de_\ell) }{n_\ell}  \mid p \ri \}\\
&  & +  \Pr \li \{ p \geq \f{\wh{\bs{p}}_\ell}{1 - \vep}, \;
\mscr{M}_{\mrm{B}} \li ( \wh{\bs{p}}_\ell, p  \ri  ) \leq \f{\ln
(\ze \de_\ell) }{n_\ell} \mid p \ri \} \\
& \leq &   \Pr \li \{ G_{\wh{\bs{p}}_\ell} \li ( \wh{\bs{p}}_\ell, p
\ri  ) \leq \ze \de_\ell \mid p \ri \} + \Pr \li \{
F_{\wh{\bs{p}}_\ell} \li (
\wh{\bs{p}}_\ell, p  \ri  ) \leq \ze \de_\ell \mid p \ri \} \\
& \leq & 2 \ze  \de_\ell \eee for any $p \in (0, 1)$ and $\ell = 1,
2, \cd$.  So, $\sum_{\ell = \ell^\star + 1}^\iy \Pr \li \{ \li |
\wh{\bs{p}}_\ell - p \ri | \geq \vep p, \; \bs{l} = \ell \mid p \ri
\} \leq 2 \ze \sum_{\ell = \ell^\star + 1}^\iy \de_\ell \leq 2 (\tau
+ 1) \ze \de$, which implies that $\Pr \li \{ \li | \wh{\bs{p}} - p
\ri | < \vep p  \mid p \ri \}  \geq 1 - \de$ provided that $\ze \leq
\f{1}{ 2 (\tau + 1)}$.

\subsubsection{Proof of Statement (IV)}

Recall that in the course of proving statement (III), we have shown
that $\Pr \li \{ \li | \wh{\bs{p}}_\ell - p \ri | \geq \vep p, \;
\bs{l} = \ell \mid p \ri \} \leq 2 \ze \de_\ell$ for any $\ell > 0$.
Making use of such result,  we have $\sum_{\ell = \ell^\star +
1}^\iy \Pr \li \{ \li | \wh{\bs{p}}_\ell - p \ri | \geq \vep p, \;
\bs{l} = \ell \mid p \ri \} \leq  2 \ze \sum_{\ell = \ell^\star +
1}^\iy \de_\ell \leq \eta$ for any $p \in (0, 1)$. It follows that
{\small \bee \Pr \li \{ \li | \wh{\bs{p}} - p \ri | \geq \vep p \mid
p \ri \} & = & \sum_{\ell = 1}^{\ell^\star} \Pr \li \{ \li |
\wh{\bs{p}}_\ell - p \ri | \geq \vep p, \; \bs{l} = \ell \mid p \ri
\} + \sum_{\ell = \ell^\star + 1}^\iy \Pr \li \{ \li |
\wh{\bs{p}}_\ell
- p \ri | \geq \vep p, \; \bs{l} = \ell \mid p \ri \}\\
& \leq & \sum_{\ell = 1}^{\ell^\star} \Pr \li \{ \li |
\wh{\bs{p}}_\ell - p \ri | \geq \vep p, \; \bs{l} = \ell \mid
p \ri \} + \eta\\
& \leq & \sum_{\ell = 1}^{\ell^\star} \Pr \li \{ \bs{l} = \ell \mid
p \ri \} + \eta \leq  \sum_{\ell = 1}^{\ell^\star} \Pr \li \{
\wh{\bs{p}}_\ell \geq z_\ell \mid p \ri \} + \eta\\
& \leq & \sum_{\ell = 1}^{\ell^\star} \exp ( n_\ell
\mscr{M}_{\mrm{B}} ( z_\ell, p)  ) + \eta < \sum_{\ell =
1}^{\ell^\star} \exp ( n_\ell \mscr{M}_{\mrm{B}} ( z_\ell, p^*)  ) +
\eta < \de \eee} for any $p \in (0, p^*)$.

Now we shall bound $\Pr  \{ p \leq \f{\wh{\bs{p}} }{1 + \vep} \}$
and $\Pr \{ p \geq \f{\wh{\bs{p}} }{1 - \vep}  \}$ for $p \in [a, b]
\subseteq (0, 1)$.  Observing that $\{ a \leq \f{\wh{\bs{p}}_\ell}{1
+ \vep} \} \subseteq \{ \wh{\bs{p}} \geq b \}$ as a consequence of
$b < a (1 + \vep)$, by statement (III) of Theorem
\ref{Main_Bound_Gen}, we have
\[
\Pr \li \{ b \leq \f{ \wh{\bs{p}} }{ 1 + \vep},  \; \bs{l} \leq
\ell^\star \mid a \ri \} \leq \Pr \li \{ p \leq \f{ \wh{\bs{p}} }{ 1
+ \vep},  \; \bs{l} \leq \ell^\star \mid p \ri \} \leq \Pr \li \{ a
\leq \f{ \wh{\bs{p}} }{ 1 + \vep},  \; \bs{l} \leq \ell^\star \mid b
\ri \}
\]
for any $p \in [a, b]$.  On the other hand, \bee \Pr \li \{ p \leq
\f{ \wh{\bs{p}} }{ 1 + \vep},  \; \bs{l} > \ell^\star \mid p \ri \}
& \leq  & \sum_{\ell = \ell^\star + 1}^\iy  \Pr \li \{ p \leq
\f{\wh{\bs{p}}_\ell}{1 + \vep}, \;  \mscr{M}_{\mrm{B}} \li (
\wh{\bs{p}}_\ell, \f{\wh{\bs{p}}_\ell}{1 + \vep} \ri  ) \leq
\f{\ln (\ze \de_\ell) }{n_\ell} \mid p \ri \}\\
& \leq & \sum_{\ell = \ell^\star + 1}^\iy  \Pr \li \{ p \leq
\f{\wh{\bs{p}}_\ell}{1 + \vep}, \;  \mscr{M}_{\mrm{B}} \li (
\wh{\bs{p}}_\ell, p  \ri  ) \leq
\f{\ln (\ze \de_\ell) }{n_\ell}  \mid p \ri \}\\
& \leq &  \sum_{\ell = \ell^\star + 1}^\iy  \Pr \li \{
G_{\wh{\bs{p}}_\ell} \li ( \wh{\bs{p}}_\ell, p  \ri  ) \leq \ze
\de_\ell \mid p \ri \} \leq  \ze \sum_{\ell = \ell^\star + 1}^\iy
\de_\ell \leq \f{\eta}{2} \eee for any $p \in [a, b]$.   Therefore,
$\Pr  \{ b \leq \f{ \wh{\bs{p}} }{ 1 + \vep},  \; \bs{l} \leq
\ell^\star \mid a \} \leq \Pr  \{ p \leq \f{\wh{\bs{p}} }{1 + \vep}
\mid p  \} = \Pr  \{ a \leq \f{ \wh{\bs{p}} }{ 1 + \vep},  \; \bs{l}
\leq \ell^\star \mid b \} + \Pr
 \{ p \leq \f{ \wh{\bs{p}} }{ 1 + \vep},  \; \bs{l} > \ell^\star
\mid p \} \leq \Pr  \{ a \leq \f{ \wh{\bs{p}} }{ 1 + \vep},  \;
\bs{l} \leq \ell^\star \mid b \} + \f{\eta}{2}$ for any $p \in [a,
b]$.

Similarly, observing that $\{ b \geq \f{\wh{\bs{p}}_\ell}{1 - \vep}
\} \subseteq \{ \wh{\bs{p}} \leq a \}$ as a consequence of $b < a (1
+ \vep)$, by statement (IV) of Theorem \ref{Main_Bound_Gen}, we have
\[
\Pr \li \{ a \geq \f{ \wh{\bs{p}} }{ 1 - \vep},  \; \bs{l} \leq
\ell^\star \mid b \ri \} \leq \Pr \li \{ p \geq \f{ \wh{\bs{p}} }{ 1
- \vep},  \; \bs{l} \leq \ell^\star \mid p \ri \} \leq \Pr \li \{ b
\geq \f{ \wh{\bs{p}} }{ 1 - \vep},  \; \bs{l} \leq \ell^\star \mid a
\ri \}
\]
for any $p \in [a, b]$.  On the other hand, \bee \Pr \li \{ p \geq
\f{ \wh{\bs{p}} }{ 1 - \vep},  \; \bs{l} > \ell^\star \mid p \ri \}
& \leq  & \sum_{\ell = \ell^\star + 1}^\iy  \Pr \li \{ p \geq
\f{\wh{\bs{p}}_\ell}{1 - \vep}, \;  \mscr{M}_{\mrm{B}} \li (
\wh{\bs{p}}_\ell, \f{\wh{\bs{p}}_\ell}{1 - \vep} \ri  ) \leq
\f{\ln (\ze \de_\ell) }{n_\ell} \mid p \ri \}\\
& \leq & \sum_{\ell = \ell^\star + 1}^\iy  \Pr \li \{ p \geq
\f{\wh{\bs{p}}_\ell}{1 - \vep}, \;  \mscr{M}_{\mrm{B}} \li (
\wh{\bs{p}}_\ell, p  \ri  ) \leq
\f{\ln (\ze \de_\ell) }{n_\ell}  \mid p \ri \}\\
& \leq &  \sum_{\ell = \ell^\star + 1}^\iy  \Pr \li \{
F_{\wh{\bs{p}}_\ell} \li ( \wh{\bs{p}}_\ell, p  \ri  ) \leq \ze
\de_\ell \mid p \ri \} \leq  \ze \sum_{\ell = \ell^\star + 1}^\iy
\de_\ell \leq \f{\eta}{2} \eee for any $p \in [a, b]$.  Therefore,
$\Pr  \{ a \geq \f{ \wh{\bs{p}} }{ 1 - \vep},  \; \bs{l} \leq
\ell^\star \mid b \} \leq \Pr  \{ p \geq \f{\wh{\bs{p}} }{1 - \vep}
\mid p  \} = \Pr  \{ b \geq \f{ \wh{\bs{p}} }{ 1 - \vep},  \; \bs{l}
\leq \ell^\star \mid a \} + \Pr
 \{ p \geq \f{ \wh{\bs{p}} }{ 1 - \vep},  \; \bs{l} > \ell^\star
\mid p \} \leq \Pr  \{ b \geq \f{ \wh{\bs{p}} }{ 1 - \vep},  \;
\bs{l} \leq \ell^\star \mid a \} + \f{\eta}{2}$ for any $p \in [a,
b]$.  This completes the proof of statement (IV).

\subsubsection{Proof of Statement (V)}  We
need a preliminary result.

\beL \la{lema} Let $p \in (0,1)$ and $\eta \in (0,1)$.  Let $\ka$ be
an integer greater than {\small $\max \{ \tau, \; \tau + \f{1}{\ga -
1} + \f{ \ln (\ze \de ) } { \ln 2 } \}$} such that {\small
$\mscr{M}_{\mrm{B}} ( \eta p, \f{\eta p}{1 + \vep} ) < \f{ \ln ( \ze
\de_\ka ) } {n_\ka}$}. Then,  {\small $\Pr \{ \bs{l} > \ell \} <
\exp  ( - \f{(1 - \eta)^2 p n_\ell}{2}  )$} for any $\ell \geq \ka$.
\eeL

\bpf

Let $m_\ell = m \ga^{\ell - 1}$ for $\ell = 1, 2, \cd$. Noting that
{\small
\[ \f{ \f{ \ln (\ze \de_{\ell + 1}) } { m_{\ell+1} } } { \f{ \ln (
\ze \de_\ell ) } { m_\ell } } = \f{1}{\ga} \times \f{ (\ell + 1 -
\tau) \ln 2  - \ln (\ze \de)  } { (\ell - \tau) \ln 2 - \ln (\ze
\de) } = \f{1}{\ga} \times \li ( 1 + \f{1}{ \ell - \tau - \f{\ln
(\ze \de ) } { \ln 2 } } \ri ) < 1
\]}
for {\small $\ell > \max \{ \tau, \tau +  \f{1}{\ga - 1} + \f{ \ln
(\ze \de) } { \ln 2 } \}$} and that {\small $\f{ \ln ( \ze \de_\ell
) } { m_\ell } = \f{ \ln \li ( \ze \de 2^{\tau - \ell} \ri ) } { m
\ga^{\ell - 1} } \to 0 > \mscr{M}_{\mrm{B}} ( \eta p, \f{\eta p}{1 +
\vep} )$} as $\ell \to \iy$, we have that there exists an integer
$\ka$ greater than {\small $\max \{ \tau, \tau + \f{1}{\ga - 1} +
\f{ \ln (\ze \de) } { \ln 2 } \}$} such that {\small
$\mscr{M}_{\mrm{B}} ( \eta p, \f{\eta p}{1 + \vep} ) < \f{ \ln (\ze
\de_\ell) } { m_\ell }$} for all $\ell \geq \ka$. Since $m_\ell \leq
n_\ell$ and $\mscr{M}_{\mrm{B}} ( \eta p, \f{\eta p}{1 + \vep} ) <
0$, we have that there exists an integer $\ka$ greater than {\small
$\max \{ \tau, \tau + \f{1}{\ga - 1} + \f{ \ln (\ze \de) } { \ln 2 }
\}$} such that {\small $\mscr{M}_{\mrm{B}} ( \eta p, \f{\eta p}{1 +
\vep} ) < \f{ \ln (\ze \de_\ell) } { n_\ell }$} for all $\ell \geq
\ka$.  For $\ell$ greater than such $\ka$, we claim that $z < \eta
p$ if {\small $\mscr{M}_{\mrm{B}} ( z, \f{z}{1 + \vep} ) > \f{\ln
(\ze \de_\ell )}{n_\ell}$} and $z \in [0, 1]$. To prove this claim,
suppose, to get a contradiction, that $z \geq \eta p$. Then, since
{\small $\mscr{M}_{\mrm{B}} ( z, \f{z}{1 + \vep} )$} is
monotonically decreasing with respect to $z \in (0,1)$ as asserted
by Lemma \ref{lemm23}, we have {\small $\mscr{M}_{\mrm{B}} (z,
\f{z}{1 + \vep}  ) \leq \mscr{M}_{\mrm{B}}  ( \eta p, \f{ \eta p}{1
+ \vep}  ) < \f{\ln (\ze \de_\ell)}{n_\ell}$}, which is a
contradiction. Therefore, we have shown the claim and it follows
that {\small $ \{ \mscr{M}_{\mrm{B}} ( \f{K_\ell}{n_\ell},
\f{K_\ell}{(1 + \vep) n_\ell }  ) > \f{\ln (\ze \de_\ell) }{n_\ell}
\} \subseteq \{ K_\ell < \eta p n_\ell \}$} for $\ell \geq \ka$. So,
{\small
\[ \Pr \{ \bs{l}
> \ell \}  \leq  \Pr \li \{ \mscr{M}_{\mrm{B}}
\li ( \f{K_\ell}{n_\ell}, \f{K_\ell}{(1 + \vep) n_\ell } \ri ) >
\f{\ln (\ze \de_\ell) }{n_\ell}  \ri \} \leq  \Pr \{ K_\ell < \eta p
n_\ell \} < \exp \li ( - \f{(1 - \eta)^2 p n_\ell}{2} \ri ),
\]} where the last inequality is due to the multiplicative Chernoff bound \cite{Chernoff_tour}.

\epf

\bsk

We are now in position to prove statement (V) of the theorem. Note
that \bee \bb{E} [ \mbf{n}] & = &  n_1 + \sum_{\ell = 1}^{\ka}
(n_{\ell + 1} - n_\ell) \Pr \{ \bs{l} > \ell \} + \sum_{\ell = \ka +
1}^\iy (n_{\ell + 1} - n_\ell) \Pr \{ \bs{l} > \ell \}. \eee By the
definition of $n_\ell$, we have $n_{\ell + 1} - n_{\ell} \leq (\ga -
1) n_{\ell}$.  By the assumption of $\ep, \; \eta$ and $\ka$, we
have $\ln \f{\ga}{c \ep} > 1$ and thus $\ka > \f{1}{\ln \ga} \ln \li
( \f{1}{ c m } \ln \f{\ga}{c \ep} \ri ) + 1 > \f{1}{\ln \ga} \ln \li
( \f{1}{ c m } \ri ) + 1$, which implies that $c m \ga^{\ka - 1}
> 1$ and $\f{\ga}{c} \exp ( - c m \ga^{\ka - 1} ) < \ep$.
Hence, by Lemma \ref{lema}, we have \bee \sum_{\ell = \ka + 1}^\iy
(n_{\ell + 1} - n_{\ell}) \; \Pr \{ \bs{l}
> \ell \} & < &  \sum_{\ell = \ka + 1}^\iy (n_{\ell + 1} -
n_{\ell}) \; e^{- c n_\ell}
 \leq \f{\ga - 1}{c} \sum_{\ell = \ka + 1}^\iy  c n_{\ell} \; e^{- c n_\ell}\\
& \leq &  \f{\ga - 1}{c} \; \sum_{\ell = \ka}^\iy c m \ga^{\ell}
\exp ( - c m \ga^{\ell} ) <  \f{\ga - 1}{c} \int_{\ka - 1}^\iy c m
\ga^{\ell} \exp ( - c m \ga^{\ell} ) d \ell. \eee Making a change of
variable $x  = c m \ga^{\ell}$, we have $d \ell = \f{1}{\ln \ga}
\f{d x}{x}$ and
\[
\int_{\ka - 1}^\iy  c m \ga^{\ell} \exp ( - c m \ga^{\ell} ) d \ell
= \f{1}{\ln \ga} \int_{ c m \ga^{\ka - 1} }^\iy e^{- x} d x =
\f{\exp ( - c m \ga^{\ka - 1} ) }{\ln \ga}.
\]
It follows that $\sum_{\ell = \ka + 1}^\iy (n_{\ell + 1} - n_{\ell})
\; \Pr \{ \bs{l}
> \ell \} < \f{\ga - 1}{c} \f{\exp ( - c m \ga^{\ka - 1} )
}{\ln \ga} < \f{\ga}{c} \exp ( - c m \ga^{\ka - 1} ) < \ep$. This
completes the proof of statement (V) of Theorem
\ref{Bino_Rev_noninverse_Chernoff}.

\subsection{Proof of Theorem \ref{Bino_Inverse_DDV_Asp} } \la{App_Bino_Inverse_DDV_Asp}

We need some preliminary results.

\beL

\la{lem79rev}

If $\vep$ is sufficiently small, then the following statements hold
true.

(I): For $\ell = 1, \cd, s - 1$, there exists a unique number
$z_\ell \in (0, 1]$ such that {\small $\ga_\ell = \f{ \ln ( \ze \de
) } { \mscr{M}_{\mrm{I}} \li ( z_\ell, \f{z_\ell}{1 + \vep} \ri )
}$}.

(II): $z_\ell$ is monotonically decreasing with respect to $\ell$.

(III): $\lim_{\vep \to 0} z_\ell = 1 - C_{s - \ell}$, where the
limit is taken under the restriction that $\ell - s$ is fixed with
respect to $\vep$.

(IV): For $p \in (0, 1)$ such that $C_{j_p} = 1 - p$,
\[
\lim_{\vep \to 0} \f{ p - z_{\ell_\vep}  } { \vep z_{\ell_\vep} } =
- \f{2}{3},
\]
where $\ell_\vep = s - j_p$.

(V): $\{ \bs{D}_\ell = 0 \} = \{ \wh{\bs{p}}_\ell < z_\ell \}$.

\eeL

\bsk

{\bf Proof of Statement (I)}: By the definition of $\ga_\ell$, we
have \be \la{from1}
 0 < \f{ \ln ( \ze \de )  } { \mscr{M}_{\mrm{I}} ( 1, \f{1}{1 + \vep}
) } \leq \ga_1 \leq \ga_{\ell} < \f{(1 + C_1) \ga_s}{2} < \f{(1 +
C_1) }{2} \li [ \f{ (1 + \vep) \ln \f{1}{ \ze \de } } { (1 + \vep)
\ln (1 + \vep) - \vep } + 1 \ri ]
 \ee
for sufficiently small $\vep > 0$. By (\ref{from1}), we have {\small
$\f{ \ln ( \ze \de ) } { \ga_{\ell} } \geq \mscr{M}_{\mrm{I}} ( 1,
\f{1}{1 + \vep} )$} and {\small \[ \f{ \ln ( \ze \de ) } {
\ga_{\ell} } < \li [ \f{\vep}{1 + \vep} - \ln (1 + \vep)  \ri ] \li
( \f{2}{1 + C_1} - \f{1}{ \ga_{\ell} } \ri ) =  \f{ \f{\vep}{1 +
\vep} - \ln (1 + \vep) }{\mscr{M}_{\mrm{I}} ( 0, 0 )}  \f{2
\mscr{M}_{\mrm{I}} ( 0, 0 ) }{1 + C_1}  + \li [ \ln (1 + \vep) -
\f{\vep}{1 + \vep} \ri ] \f{1}{\ga_\ell} .
\]}
Noting that $\lim_{\vep \to 0} \f{{ (1 + \vep) \ln (1 + \vep) - \vep
}}{(1 + \vep) \ga_\ell}  = 0$ and $\lim_{\vep \to 0}  \f{{ \vep - (1
+ \vep) \ln (1 + \vep) }}{(1 + \vep) \mscr{M}_{\mrm{I}} ( 0, 0 )}  =
1$, we have {\small $\f{ \ln ( \ze \de ) } { \ga_{\ell} } <
\mscr{M}_{\mrm{I}} \li ( 0, 0 \ri )$} for sufficiently small $\vep
> 0$. In view of the established fact that $ \mscr{M}_{\mrm{I}} ( 1,
\f{1}{1 + \vep} ) \leq \f{ \ln ( \ze \de ) } { \ga_{\ell} } <
\mscr{M}_{\mrm{I}} ( 0, 0)$ for small enough $\vep > 0$ and the fact
that $\mscr{M}_{\mrm{I}} ( z, \f{z}{1 + \vep} )$ is monotonically
decreasing with respect to $z \in (0, 1)$ as asserted by Lemma
\ref{decr}, invoking the intermediate value theorem, we have that
there exists a unique number $z_{\ell} \in (0, 1]$ such that
$\mscr{M}_{\mrm{I}} ( z_{\ell}, \f{z_{\ell}}{1 + \vep} ) = \f{ \ln (
\ze \de ) } { \ga_{\ell} }$, which implies Statement (I).

\bsk

{\bf Proof of Statement (II)}: Since $\ga_{\ell}$ is monotonically
increasing with respect to $\ell$ for sufficiently small $\vep > 0$,
we have that $\mscr{M}_{\mrm{I}} ( z_{\ell}, \f{z_{\ell}}{1 + \vep}
)$ is monotonically increasing with respect to $\ell$ for
sufficiently small $\vep > 0$. Recalling that $\mscr{M}_{\mrm{I}} (
z, \f{z}{1 + \vep} )$ is monotonically decreasing with respect to $z
\in (0, 1)$, we have that $z_{\ell}$ is monotonically decreasing
with respect to $\ell$. This establishes Statement (II).

\bsk

{\bf Proof of Statement (III)}:  For simplicity of notations, let
$b_{\ell} = 1 - C_{s - \ell}$ for $\ell = 1, 2, \cd, s - 1$. Then,
it can be checked that $1 - b_{\ell} = C_{s - \ell}$ and, by the
definition of $\ga_\ell$, we have  \be \la{goode1} \f{ (1 -
b_{\ell}) (1 + \vep) \mscr{M}_{\mrm{I}} (z_\ell, \f{z_\ell}{1 +
\vep} ) } { \vep - (1 + \vep) \ln (1 + \vep) } = \f{1}{\ga_{\ell} }
\times \f{ C_{s - \ell} \; (1 + \vep) \ln \f{1}{ \ze \de } } { (1 +
\vep) \ln (1 + \vep) - \vep } = 1 + o(1)  \ee for $\ell = 1, 2, \cd,
s - 1$.

We claim that $z_\ell < \se$ for $\se \in (b_\ell, 1)$ if $\vep > 0$
is small enough.  To prove this claim, we use a contradiction
method. Suppose the claim is not true, then there exists a set,
denoted by $S_\vep$,  of infinitely many values of $\vep$ such that
$z_\ell \geq \se$ for $\vep \in S_\vep$.
  By (\ref{goode1}) and the fact that {\small $\mscr{M}_{\mrm{I}}
(z, \f{z}{1 + \vep} )$} is monotonically decreasing with respect to
$z \in (0, 1)$ as asserted by Lemma \ref{decr}, we have {\small \bee
1 + o (1)  =  \f{ (1 - b_\ell) (1 + \vep) \mscr{M}_{\mrm{I}}
(z_\ell, \f{z_\ell}{1 + \vep} ) } { \vep - (1 + \vep) \ln (1 + \vep)
} \geq \f{ (1 - b_\ell) (1 + \vep) \mscr{M}_{\mrm{I}} (\se,
\f{\se}{1 + \vep} ) } { \vep - (1 + \vep) \ln (1 + \vep) } = \f{ 1 -
b_\ell } { 1 - \se } + o(1) \eee} for small enough $\vep \in
S_\vep$,  which implies {\small $\f{ 1 - b_\ell } { 1 - \se } \leq
1$}, contradicting to the fact that {\small $\f{ 1 - b_\ell } { 1 -
\se } > 1$}.  The claim is thus established.  Similarly, we can show
that $z_\ell > \se^\prime$ for $\se^\prime \in (0, b_\ell)$ if
$\vep$ is small enough. Now we restrict $\vep$ to be small enough so
that $\se^\prime < z_\ell < \se$. Applying Lemma \ref{lem32T} based
on such restriction, we have \be \f{ (1 - b_\ell) (1 + \vep)
\mscr{M}_{\mrm{I}} (z_\ell, \f{z_\ell}{1 + \vep} ) } { \vep - (1 +
\vep) \ln (1 + \vep) }  = \f{ (1 - b_\ell)  \li [ - \f{ \vep^2 } { 2
( 1 - z_\ell) } + o(\vep^2) \ri ]  } { - \f{\vep^2}{2} + o(\vep^2) }
=  \f{ \f{ 1 - b_\ell } { 1 - z_\ell } + o(1) } { 1  + o(1) }.
\la{imm} \ee Combining (\ref{goode1}) and (\ref{imm}) yields  $\f{
b_\ell  - z_\ell} { 1 - z_\ell }  =  o(1)$,  which implies
$\lim_{\vep \to 0} z_\ell = b_\ell$.  This proves Statement (III).

\bsk

{\bf Proof of Statement (IV)}:

Since {\small $\ga_{\ell_\vep} =  \li \lc  \f{ C_{s - \ell_\vep } \;
(1 + \vep) \ln ( \ze \de ) } { \vep - (1 + \vep) \ln (1 + \vep)  }
\ri \rc$} and $C_{s - {\ell_\vep}} = 1 - p$, we can write
\[
\ga_{\ell_\vep} =  \li \lc  \f{ (1 - p) \; (1 + \vep) \ln ( \ze \de
) } { \vep - (1 + \vep) \ln (1 + \vep)  } \ri \rc = \f{ \ln (\ze
\de) } { \mscr{M}_{\mrm{I}} (z_{\ell_\vep}, z_{\ell_\vep} \sh ( 1 +
\vep) )},
\]
from which we have $\f{1}{\ga_{\ell_\vep}} = o (\vep)$,
\[
1 - o(\vep) = 1 - \f{1}{\ga_{\ell_\vep}}  < \f{ \f{ (1 - p) \; (1 +
\vep) \ln ( \ze \de ) } { \vep - (1 + \vep) \ln (1 + \vep)  }  } {
\f{ \ln (\ze \de) } { \mscr{M}_{\mrm{I}} (z_{\ell_\vep},
z_{\ell_\vep} \sh ( 1 + \vep) ) } } \leq 1
\]
and thus \[ \f{ \f{ (1 - p) \; (1 + \vep) \ln ( \ze \de ) } { \vep -
(1 + \vep) \ln (1 + \vep)  } } { \f{ \ln (\ze \de) } {
\mscr{M}_{\mrm{I}} (z_{\ell_\vep}, z_{\ell_\vep} \sh ( 1+ \vep) ) }
} = 1 + o(\vep). \] For $\se \in (p,  1)$, we claim that
$z_{\ell_\vep} < \se$ if $\vep$ is sufficiently small.  Suppose, to
get a contradiction that the claim is not true.  Then, there exists
a set of infinitely many values of $\vep$ such that $z_{\ell_\vep}
\geq \se$ if $\vep$ in the set is small enough. For such $\vep$, by
the monotonicity of $\mscr{M}_{\mrm{I}} (., .)$, we have {\small
\bel 1 + o (\vep) = \f{ \f{ (1 - p) (1 + \vep) \ln \f{1}{\ze \de} }
{ (1 + \vep) \ln (1 + \vep) - \vep}  } { \f{ \ln (\ze \de) } {
\mscr{M}_{\mrm{I}} (z_{\ell_\vep}, z_{\ell_\vep} \sh (1 + \vep) )  }
} & = & \f{ (1 - p) (1 + \vep) \mscr{M}_{\mrm{I}} (z_{\ell_\vep},
z_{\ell_\vep} \sh (1 + \vep) ) } { \vep - (1 +
\vep) \ln (1 + \vep) } \la{qute}\\
 & > & \f{ (1 - p) (1 + \vep)
\mscr{M}_{\mrm{I}} (\se, \se \sh (1 + \vep) ) } { \vep - (1 + \vep)
\ln (1 + \vep) } = \f{ 1 - p } { 1 - \se } + o(1) \nonumber \eel}
for small enough $\vep$ in the set,  which contradicts to the fact
that {\small $\f{ 1 - p } { 1 - \se } > 1$}.  This proves the claim.
 Similarly, we can show that $z_{\ell_\vep} \geq \se^\prime$ for any
$\se^\prime \in (0, p)$. Now we restrict $\vep$ to be small enough
so that $\se^\prime < z_{\ell_\vep} < \se$. By virtue of such
restriction, we have \bel \f{ (1 - p) (1 + \vep) \mscr{M}_{\mrm{I}}
(z_{\ell_\vep}, z_{\ell_\vep} \sh (1 + \vep) ) } { \vep - (1 + \vep)
\ln (1 + \vep) }  & = & \f{ (1 - p)  \li [ - \f{ \vep^2 } { 2 ( 1 -
z_{\ell_\vep}) } + \f{ \vep^3 ( 2 - z_{\ell_\vep} ) } { 3 (1 -
z_{\ell_\vep})^2 } + o(\vep^3) \ri ] }
{ \vep \sh (1 + \vep) - \ln (1 + \vep) } \nonumber\\
& = & \f{ (1 - p)  \li [ - \f{ \vep^2 } { 2 ( 1 - z_{\ell_\vep}) } +
\f{ \vep^3 ( 2 -  z_{\ell_\vep} ) } { 3 (1 - z_{\ell_\vep})^2 } +
o(\vep^3) \ri ] } {
\vep [1 - \vep + \vep^2 + o(\vep^2)]  - [\vep - \f{\vep^2}{2} + \f{\vep^3}{3} + o(\vep^3)] } \nonumber\\
& = & \f{ (1 - p)  \li [ - \f{ \vep^2 } { 2 ( 1 - z_{\ell_\vep}) } +
\f{ \vep^3 ( 2 -  z_{\ell_\vep} ) } { 3 (1 - z_{\ell_\vep})^2 } +
o(\vep^3) \ri ] }
{ - \f{\vep^2}{2} + \f{2 \vep^3}{3} + o(\vep^3) } \nonumber\\
& = & \f{ \f{ 1 - p } { 1 - z_{\ell_\vep} } - \f{ 2 \vep (1 - p) ( 2
- z_{\ell_\vep} ) } { 3 (1 - z_{\ell_\vep})^2 } + o(\vep) } { 1 -
\f{4 \vep}{3} + o(\vep) }.  \la{imm} \eel Combining (\ref{qute}) and
(\ref{imm}) yields {\small $\f{ 1 - p } { 1 - z_{\ell_\vep} } - \f{
2 \vep (1 - p) ( 2 - z_{\ell_\vep} ) } { 3 (1 - z_{\ell_\vep})^2 } =
1 - \f{4 \vep}{3} + o(\vep)$},  i.e.,
\[
\f{ p  - z_{\ell_\vep}} { 1 - z_{\ell_\vep} }  =  \f{4 \vep}{3} -
\f{ 2 \vep (1 - p) ( 2 - z_{\ell_\vep} ) } { 3 (1 - z_{\ell_\vep})^2
} + o(\vep),
\]
i.e.,
\[
\f{ p  - z_{\ell_\vep}} { \vep z_{\ell_\vep} }  =  \f{4  (1 -
z_{\ell_\vep})}{3 z_{\ell_\vep} } - \f{ 2 (1 - p) ( 2 -
z_{\ell_\vep} ) } { 3 z_{\ell_\vep} (1 - z_{\ell_\vep}) } + o(1),
\]
which implies that $\lim_{\vep \to 0} \f{ p  - z_{\ell_\vep}} { \vep
z_{\ell_\vep} }  =  \f{4  (1 - p)}{3 p} - \f{ 2 ( 2 - p ) } { 3 p }
= - \f{2}{3}$. \bsk

{\bf Proof of Statement (V)}: Noting that {\small
$\mscr{M}_{\mrm{I}} (z, \f{z}{1 + \vep} )$} is monotonically
decreasing with respect to $z \in (0, 1)$ as asserted by Lemma
\ref{decr}, we have {\small $\{ \bs{D}_{\ell} = 0 \}  =  \{
\mscr{M}_{\mrm{I}} ( \wh{\bs{p}}_{\ell}, \; \f{\wh{\bs{p}}_{\ell}
}{1  + \vep}  ) > \f{\ln (\ze \de)} {\ga_{\ell}} \} = \li \{
\wh{\bs{p}}_{\ell} < z_\ell \ri \}$} as claimed by statement (V).

\beL \la{Defrev} Let $\ell_\vep = s  - j_p$.  Then, \be \la{revequ}
\lim_{\vep \to 0} \sum_{\ell = 1}^{\ell_\vep - 1} \ga_\ell \Pr \{
\bs{D}_\ell = 1 \} = 0, \qqu \lim_{\vep \to 0} \sum_{\ell =
\ell_\vep + 1}^s \ga_\ell \Pr \{ \bs{D}_{\ell} = 0 \} = 0 \ee for $p
\in (0, 1)$.  Moreover, $\lim_{\vep \to 0} \ga_{\ell_\vep} \Pr \{
\bs{D}_{\ell_\vep} = 0 \} = 0$ if $C_{j_p} > 1 - p$. \eeL

\bpf

For simplicity of notations, let $b_\ell = \lim_{\vep \to 0}
z_{\ell}$ for $1 \leq \ell < s$.  The proof consists of two main
steps as follows.

First, we shall show that (\ref{revequ}) holds for any $p \in (0,
1)$. By the definition of $\ell_\vep$, we have $1 - p  > C_{s -
\ell_\vep + 1}$. Making use of the first three statements of Lemma
\ref{lem79rev}, we have
 {\small $z_\ell > \f{p  + b_{\ell_\vep - 1}}{2} > p$}
 for all $\ell \leq \ell_\vep - 1$ if $\vep$ is sufficiently
small. By the last statement of Lemma \ref{lem79rev} and using Lemma
\ref{lem_chen}, we have \bee \Pr \{ \bs{D}_{\ell} = 1 \}  =  \Pr \{
\wh{\bs{p}}_{\ell} \geq z_{\ell} \} \leq  \Pr \li \{
\wh{\bs{p}}_{\ell} \geq \f{p  + b_{\ell_\vep - 1}}{2} \ri  \} \leq
\exp \li (  \ga_\ell \mscr{M}_{\mrm{I}} \li ( \f{p  + b_{\ell_\vep -
1}}{2}, p \ri )  \ri ) \eee for all $\ell \leq \ell_\vep - 1$ if
$\vep > 0$ is sufficiently small. Since $b_{\ell_\vep - 1}$ is
greater than $p$ and is independent of $\vep > 0$ as a consequence
of the definition of $\ell_\vep$, it follows from Lemma \ref{lem31a}
that {\small $\lim_{\vep \to 0} \sum_{\ell = 1}^{\ell_\vep - 1}
\ga_\ell \Pr \{ \bs{D}_\ell = 1 \}  = 0$}.

Similarly, it can be seen from the definition of $\ell_\vep$ that $1
- p  < C_{s - \ell_\vep - 1}$.  Making use of the first three
statements of Lemma \ref{lem79rev}, we have that {\small $z_\ell <
\f{p + b_{\ell_\vep + 1}}{2} < p$} for $\ell_\vep + 1 \leq \ell < s$
if $\vep$ is sufficiently small.  By the last statement of Lemma
\ref{lem79rev} and using Lemma \ref{lem_chen}, we have
\[
\Pr \{ \bs{D}_{\ell} = 0 \} = \Pr \{  \wh{\bs{p}}_{\ell} < z_{\ell}
\} \leq \Pr \li \{ \wh{\bs{p}}_{\ell} < \f{p + b_{\ell_\vep + 1}}{2}
\ri \} \leq  \exp \li (  \ga_\ell \mscr{M}_{\mrm{I}} \li ( \f{p  +
b_{\ell_\vep + 1}}{2}, p \ri )  \ri )
\]
for $\ell_\vep + 1 \leq \ell < s$ if $\vep
> 0$ is small enough.  By virtue of the definition of $\ell_\vep$, we have that $b_{\ell_\vep + 1}$ is smaller
than $p$ and is independent of $\vep > 0$. In view of this and the
fact that $\Pr \{ \bs{D}_s = 0 \} = 0$, we can use Lemma
\ref{lem31a} to conclude that {\small $\lim_{\vep \to 0} \sum_{\ell
= \ell_\vep + 1}^s \ga_\ell \Pr \{\bs{D}_\ell = 0 \}  = 0$}.

\bsk

Next, we shall show that $\lim_{\vep \to 0} \ga_{\ell_\vep} \Pr \{
\bs{D}_{\ell_\vep} = 0 \} = 0$ for $p \in (0, 1)$ such that $C_{j_p}
> 1 - p$.  Note that $1 - p < C_{s - \ell_\vep}$
because of the definition of $\ell_\vep$. Making use of the first
three statements of Lemma \ref{lem79rev}, we have that {\small
$z_{\ell_\vep} < \f{p + b_{\ell_\vep}}{2}  < p$} if $\vep
> 0$ is small enough. By the last statement of Lemma \ref{lem79rev} and using Lemma \ref{lem_chen}, we have
\[
\Pr \{ \bs{D}_{\ell_\vep} = 0 \} = \Pr \{ \wh{\bs{p}}_{\ell_\vep} <
z_{\ell_\vep} \} \leq \Pr \li \{ \wh{\bs{p}}_{\ell_\vep} < \f{p +
b_{\ell_\vep}}{2} \ri \} \leq  \exp \li (  \ga_{\ell_\vep}
\mscr{M}_{\mrm{I}} \li ( \f{p  + b_{\ell_\vep}}{2}, p \ri )  \ri )
\] for small enough $\vep
> 0$. By virtue of the definition of $\ell_\vep$, we have that $b_{\ell_\vep}$ is smaller
than $p$ and is independent of $\vep > 0$. It follows that
$\lim_{\vep \to 0} \ga_{\ell_\vep} \Pr \{ \bs{D}_{\ell_\vep} = 0 \}
= 0$.  This completes the proof of the lemma.

\epf

\bsk

Finally, we would like to note that Theorem
\ref{Bino_Inverse_DDV_Asp} can be shown by employing Lemma
\ref{Defrev} and a similar argument as the proof of Theorem
\ref{Bino_DDV_Asp}.

\subsection{Proof of Theorem \ref{Bino_Inverse_Asp_Analysis}}  \la{App_Bino_Inverse_Asp_Analysis}

We need some preliminary results.

\beL

\la{lem81rev}

$\lim_{\vep \to 0} \f{\ga_{\ell_\vep}}{ \ga (p, \vep)} = \ka_p, \;
\lim_{\vep \to 0} \vep \sq{ \f{\ga_{{\ell_\vep}} }{1 - p} } = d
\sq{\ka_p}$. \eeL

\bpf

By the definition of $\ga_\ell$, we have
\[
\lim_{\vep \to 0} \f{ C_{s - \ell} \; (1 + \vep) \ln (\ze \de) } {
\ga_\ell [ \vep - (1 + \vep) \ln (1 + \vep)  ]} = 1
\]
for $1 \leq \ell < s$.  It follows that \bee \lim_{\vep \to 0}
\f{\ga_{\ell_\vep}}{ \ga (p, \vep)} & = & \lim_{\vep \to 0} \f{
\mscr{M}_{\mrm{I}} (p, \f{p}{1 + \vep} )} { \ln (\ze \de) } \times
\f{ C_{s - \ell_\vep} \; (1 + \vep) \ln (\ze \de) } { \vep - (1 +
\vep) \ln (1 + \vep) } = \lim_{\vep \to 0} \f{ C_{s - \ell_\vep} \;
(1 + \vep) \mscr{M}_{\mrm{I}} (p,
\f{p}{1 + \vep} ) } { \vep - (1 + \vep) \ln (1 + \vep) }\\
& = & \lim_{\vep \to 0} \f{ C_{s - \ell_\vep} \; (1 + \vep) \li (
\vep^2 \sh [2 (p - 1) ] + o (\vep^2) \ri )} { \vep - (1 + \vep) \ln
(1 + \vep) } = \f{ C_{s - \ell_\vep}  } { 1 - p } = \f{ C_{j_p} }{1
- p}  = \ka_p \eee and \bee \lim_{\vep \to 0} \vep \sq{
\f{\ga_{{\ell_\vep}} }{1 - p} } & = & \lim_{\vep \to 0} \vep \sq{
\f{1}{1 - p}  \f{ C_{s - \ell_\vep} \; (1 + \vep) \ln \f{1}{\ze \de}
} { (1 + \vep) \ln (1 + \vep) - \vep} } = d \sq{\f{ C_{s -
\ell_\vep} } { 1 - p }} = d \sq{\ka_p}.  \eee

\epf

\beL \la{limpleminv} Let $U$ and $V$ be independent Gaussian random
variables with zero means and unit variances.  Then, for $p \in (0,
1)$ such that $C_{j_p} = 1 - p$, \bee  &  & \lim_{\vep \to 0} \Pr \{
\bs{l} = \ell_\vep \} = 1 - \lim_{\vep \to 0} \Pr \{ \bs{l} =
\ell_\vep + 1 \} =  1 - \Phi \li (  \nu  d \ri ),\\
&  & \lim_{\vep \to 0} \li [ \Pr \{ | \wh{\bs{p}}_{\ell_\vep} - p |
\geq \vep p, \; \bs{l} = \ell_\vep \} + \Pr \{ |
\wh{\bs{p}}_{\ell_\vep + 1} - p | \geq \vep p, \;
\bs{l} =  \ell_\vep + 1 \} \ri ]\\
&   & \qqu \qqu \qqu \qqu  = \Pr \li \{ U \geq d \ri \} + \Pr \li \{ |U + \sq{\ro_p} V | \geq (1 + \ro_p) d, \; U < \nu d \ri \}. \eee

\eeL

\bpf

By Statement (V) of Lemma \ref{lem79rev}, we have {\small \bee & &
\Pr \{ \wh{\bs{p}}_{\ell_\vep} \geq z_{\ell_\vep} \} \geq \Pr \{
\bs{l} = {\ell_\vep} \} \geq  \Pr \{ \wh{\bs{p}}_{\ell_\vep} \geq
z_{\ell_\vep} \} - \sum_{\ell = 1}^{ \ell_\vep - 1  } \Pr \{
\bs{D}_\ell = 1 \}, \\
&  & \Pr \{ \wh{\bs{p}}_{\ell_\vep} < z_{\ell_\vep} \} \geq \Pr \{
\bs{l} = \ell_\vep + 1 \}  \geq  \Pr \{ \wh{\bs{p}}_{\ell_\vep} <
z_{\ell_\vep} \} - \Pr \{ \bs{D}_{\ell_\vep + 1} = 0  \} -
\sum_{\ell = 1}^{ \ell_\vep - 1  } \Pr \{ \bs{D}_\ell = 1 \}. \eee}
Making use of this result and the fact that $\lim_{\vep \to 0} \li [
\sum_{\ell = 1}^{{\ell_\vep} - 1} \Pr \{ \bs{D}_\ell = 1 \} + \Pr \{
\bs{D}_{{\ell_\vep} + 1} = 0 \} \ri ] = 0$ as asserted by  Lemma
\ref{Defrev}, we have
\[
\lim_{\vep \to 0} \Pr \{  \bs{l} = {\ell_\vep} \} = \lim_{\vep \to
0} \Pr \{ \wh{\bs{p}}_{\ell_\vep} \geq z_{\ell_\vep} \}, \qqu
\lim_{\vep \to 0} \Pr \{ \bs{l} = \ell_\vep + 1 \} = \lim_{\vep \to
0} \Pr \{ \wh{\bs{p}}_{\ell_\vep} < z_{\ell_\vep} \}.
\]
Noting that {\small \bee &  & \Pr \{ | \wh{\bs{p}}_{\ell_\vep} - p |
\geq \vep p, \; \bs{l} = {\ell_\vep} \}  \geq  \Pr \{ |
\wh{\bs{p}}_{\ell_\vep} - p | \geq \vep p, \;
\wh{\bs{p}}_{\ell_\vep} \geq z_{\ell_\vep} \} - \sum_{\ell = 1}^{
\ell_\vep - 1  } \Pr \{ \bs{D}_\ell = 1 \}, \\
&  & \Pr \{ | \wh{\bs{p}}_{\ell_\vep + 1} - p | \geq \vep p, \;
\bs{l} = \ell_\vep + 1 \} \geq \Pr \{ | \wh{\bs{p}}_{\ell_\vep + 1}
- p | \geq \vep p, \; \wh{\bs{p}}_{\ell_\vep} < z_{\ell_\vep} \} -
\Pr \{ \bs{D}_{\ell_\vep + 1} = 0  \} - \sum_{\ell = 1}^{ \ell_\vep
- 1  } \Pr \{ \bs{D}_\ell = 1 \} \qqu \eee} and using the result
that $\lim_{\vep \to 0} \li [ \sum_{\ell = 1}^{{\ell_\vep} - 1} \Pr
\{ \bs{D}_\ell = 1 \} + \Pr \{ \bs{D}_{{\ell_\vep} + 1} = 0 \} \ri ]
= 0$, we have \bee  &   & \liminf_{\vep \to 0} \li [ \Pr \{ |
\wh{\bs{p}}_{\ell_\vep} - p | \geq \vep p, \; \bs{l} = {\ell_\vep}
\} + \Pr \{ | \wh{\bs{p}}_{\ell_\vep + 1}
- p | \geq \vep p, \; \bs{l} = {{\ell_\vep} + 1} \} \ri ]\\
& \geq  & \lim_{\vep \to 0} \li [ \Pr \{ | \wh{\bs{p}}_{\ell_\vep} -
p | \geq \vep p, \; \wh{\bs{p}}_{\ell_\vep} \geq z_{\ell_\vep} \} +
\Pr \{ | \wh{\bs{p}}_{\ell_\vep + 1} - p | \geq \vep p, \;
\wh{\bs{p}}_{\ell_\vep} < z_{\ell_\vep} \} \ri ].  \eee On the other
hand, \bee  &   & \limsup_{\vep \to 0} \li [ \Pr \{ |
\wh{\bs{p}}_{\ell_\vep} - p | \geq \vep p, \; \bs{l} = {\ell_\vep}
\} + \Pr \{ | \wh{\bs{p}}_{\ell_\vep + 1}
- p | \geq \vep p, \; \bs{l} = {{\ell_\vep} + 1} \} \ri ]\\
& \leq  & \lim_{\vep \to 0} \li [ \Pr \{ | \wh{\bs{p}}_{\ell_\vep} -
p | \geq \vep p, \; \wh{\bs{p}}_{\ell_\vep} \geq z_{\ell_\vep} \} +
\Pr \{ | \wh{\bs{p}}_{\ell_\vep + 1} - p | \geq \vep p, \;
\wh{\bs{p}}_{\ell_\vep} < z_{\ell_\vep} \} \ri ].  \eee Therefore,
\bee  &   & \lim_{\vep \to 0} \li [ \Pr \{ | \wh{\bs{p}}_{\ell_\vep}
- p | \geq \vep p, \; \bs{l} = {\ell_\vep} \} + \Pr \{ |
\wh{\bs{p}}_{\ell_\vep + 1}
- p | \geq \vep p, \; \bs{l} = {{\ell_\vep} + 1} \} \ri ]\\
& =  & \lim_{\vep \to 0} \li [ \Pr \{ | \wh{\bs{p}}_{\ell_\vep} - p
| \geq \vep p, \; \wh{\bs{p}}_{\ell_\vep} \geq z_{\ell_\vep} \} +
\Pr \{ | \wh{\bs{p}}_{\ell_\vep + 1} - p | \geq \vep p, \;
\wh{\bs{p}}_{\ell_\vep} < z_{\ell_\vep} \} \ri ].  \eee Since $\ka_p
= 1$, by Lemma \ref{lem81rev} and Statement (IV) of Lemma
\ref{lem79rev}, we have
\[
\lim_{\vep \to 0}  \f{p - z_{\ell_\vep}}{ z_{\ell_\vep}} \sq{
\f{\ga_{{\ell_\vep}} }{1 - p} } = \lim_{\vep \to 0} \vep \sq{
\f{\ga_{{\ell_\vep}} }{1 - p} } \lim_{\vep \to 0}  \f{p -
z_{\ell_\vep}}{ \vep z_{\ell_\vep}} = d \lim_{\vep \to 0} \f{p -
z_{\ell_\vep}}{ \vep z_{\ell_\vep}} = - \f{2}{3} d = - \nu d.
\]
Note that \bee \f{1}{ \wh{\bs{p}}_{{\ell_\vep} + 1} } - \f{1}{p}  =
\f{ \mbf{n}_{{\ell_\vep} + 1} } {\ga_{{\ell_\vep} + 1}} - \f{1}{p} =
\f{ \ga_{{\ell_\vep}} } {\ga_{{\ell_\vep} + 1}} \sq{ \f{1 - p}{p^2
\ga_{\ell_\vep}} } U_{\ell_\vep} + \f{\ga_{{\ell_\vep} + 1} -
\ga_{{\ell_\vep}} } {\ga_{{\ell_\vep} + 1}} \sq{ \f{1 - p}{p^2
(\ga_{{\ell_\vep} + 1} - \ga_{\ell_\vep}) }  } V_{\ell_\vep}  \eee
where
\[
U_{\ell_\vep} = \li ( \f{1} {\wh{\bs{p}}_{{\ell_\vep}}}  - \f{1}{p}
\ri ) \sq{ \f{p^2 \ga_{\ell_\vep}}{1 - p}  }, \qqu V_{\ell_\vep} =
\li ( \f{ \mbf{n}_{{\ell_\vep} + 1} - \mbf{n}_{{\ell_\vep}}  } {
\ga_{{\ell_\vep} + 1} - \ga_{{\ell_\vep}}} - \f{1}{p} \ri )  \sq{
\f{p^2 (\ga_{{\ell_\vep} + 1} - \ga_{\ell_\vep}) }{1 - p}  }.
\]
By the central limit theorem, $U_{\ell_\vep} \to U$ and
$V_{\ell_\vep}  \to V$ as $\vep \to 0$.  Hence, {\small \bee
U_{{\ell_\vep} + 1} = \li ( \f{1}{ \wh{\bs{p}}_{{\ell_\vep} + 1} } -
\f{1}{p} \ri ) \sq{ \f{p^2 \ga_{{\ell_\vep} + 1}}{1 - p}  } & = &
\li [ \f{ \ga_{{\ell_\vep}} } {\ga_{{\ell_\vep} + 1}} \sq{ \f{1 -
p}{p^2 \ga_{\ell_\vep}} } U_{\ell_\vep} + \f{\ga_{{\ell_\vep} + 1} -
\ga_{{\ell_\vep}} } {\ga_{{\ell_\vep} + 1}} \sq{ \f{1 - p}{p^2
(\ga_{{\ell_\vep} + 1} - \ga_{\ell_\vep}) }  }
V_{\ell_\vep}  \ri ] \sq{ \f{p^2 \ga_{{\ell_\vep} + 1}}{1 - p}  } \\
& = &  \sq{ \f{ \ga_{{\ell_\vep}} } {\ga_{{\ell_\vep} + 1}} }
U_{\ell_\vep} + \sq{ \f{\ga_{{\ell_\vep} + 1} - \ga_{{\ell_\vep}} }
{\ga_{{\ell_\vep} + 1}} } V_{\ell_\vep} \to \sq{ \f{1}{1 + \ro_p} }
U + \sq{ \f{\ro_p}{1 + \ro_p} } V \eee} as $\vep \to 0$.  It can be
seen that $\Pr \{ \wh{\bs{p}}_{{\ell_\vep}} \geq z_{\ell_\vep} \} =
\Pr \{ U_{{\ell_\vep}} \leq \f{p - z_{\ell_\vep}}{p z_{\ell_\vep}}
\sq{ \f{\ga_{{\ell_\vep}} }{1 - p} } \} $, {\small \bee & & \Pr \{ |
\wh{\bs{p}}_{{\ell_\vep} + 1} - p
| \geq \vep p, \; \wh{\bs{p}}_{{\ell_\vep}} < z_{\ell_\vep}  \}\\
& = &  \Pr \{ \wh{\bs{p}}_{{\ell_\vep} + 1} \geq (1 + \vep) p, \;
\wh{\bs{p}}_{{\ell_\vep}} < z_{\ell_\vep}  \} + \Pr \{
\wh{\bs{p}}_{{\ell_\vep} + 1}
\leq (1 - \vep) p, \; \wh{\bs{p}}_{{\ell_\vep}} < z_{\ell_\vep}  \}\\
& = &  \Pr \li \{ \f{1} {\wh{\bs{p}}_{{\ell_\vep} + 1} } - \f{1}{p}
\leq - \f{\vep}{(1 + \vep) p}, \; \f{1}{\wh{\bs{p}}_{{\ell_\vep}} }
- \f{1}{p}
> \f{p - z_{\ell_\vep}}{p z_{\ell_\vep}} \ri \} + \Pr \li \{ \f{1} {\wh{\bs{p}}_{{\ell_\vep} + 1} } - \f{1}{p} \geq
\f{\vep}{(1 - \vep) p}, \; \f{1}{\wh{\bs{p}}_{{\ell_\vep}} } -
\f{1}{p}
> \f{p - z_{\ell_\vep}}{p z_{\ell_\vep}} \ri \}\\
& = & \Pr \li \{ U_{{\ell_\vep} + 1} \leq - \f{\vep}{(1 + \vep) }
\sq{ \f{\ga_{{\ell_\vep} + 1} }{1 - p} }, \; U_{{\ell_\vep}} >  \f{p
- z_{\ell_\vep}}{
z_{\ell_\vep}} \sq{ \f{\ga_{{\ell_\vep}} }{1 - p} } \ri \}\\
&  & +  \Pr \li \{ U_{{\ell_\vep} + 1} \geq  \f{\vep}{(1 - \vep) }
\sq{ \f{\ga_{{\ell_\vep} + 1} }{1 - p} }, \; U_{{\ell_\vep}} >  \f{p
- z_{\ell_\vep}}{ z_{\ell_\vep}} \sq{ \f{\ga_{{\ell_\vep}} }{1 - p}
} \ri \} \eee} and {\small \bee \Pr \{ | \wh{\bs{p}}_{{\ell_\vep}} -
p | \geq \vep p, \; \wh{\bs{p}}_{{\ell_\vep}} \geq z_{\ell_\vep} \}
& = & \Pr \li \{ U_{{\ell_\vep}} \leq - \f{\vep}{(1 + \vep) } \sq{
\f{\ga_{{\ell_\vep}} }{1 - p} }, \; U_{{\ell_\vep}} \leq  \f{p -
z_{\ell_\vep}}{
z_{\ell_\vep}} \sq{ \f{\ga_{{\ell_\vep}} }{1 - p} } \ri \}\\
&  & +  \Pr \li \{ U_{{\ell_\vep}} \geq  \f{\vep}{(1 - \vep) } \sq{
\f{\ga_{{\ell_\vep}} }{1 - p} }, \; U_{{\ell_\vep}} \leq  \f{p -
z_{\ell_\vep}}{ z_{\ell_\vep}} \sq{ \f{\ga_{{\ell_\vep}} }{1 - p} }
\ri \}. \eee} Therefore, for $p \in (0, 1)$ such that $C_{j_p} = 1 -
p$, we have $\lim_{\vep \to 0} \Pr \{ \bs{l} = \ell_\vep \} = 1 -
\lim_{\vep \to 0} \Pr \{ \bs{l} = \ell_\vep + 1 \} =  1 - \Phi \li (
\nu  d \ri )$ and \bee &  & \lim_{\vep \to 0} \li [ \Pr \{ |
\wh{\bs{p}}_{\ell_\vep} - p | \geq \vep p, \; \bs{l} = \ell_\vep \}
+ \Pr \{ | \wh{\bs{p}}_{\ell_\vep + 1} - p | \geq \vep
p, \; \bs{l} =  \ell_\vep + 1 \} \ri ] \\
 & \to & \Pr \li \{  | U |
\geq d, \; U \leq - \nu d
\ri \} + \Pr \li \{ \li | U + \sq{\ro_p}V \ri | \geq (1 + \ro_p) d, \; U > - \nu d \ri \}\\
& = & \Pr \li \{ U \geq d \ri \} + \Pr \li \{ \li |  U + \sq{\ro_p} V \ri | \geq ( 1 + \ro_p ) d, \; U < \nu d \ri \}  \eee as $\vep \to 0$.
This completes the proof of the lemma.

 \epf

\subsubsection{Proof of Statement (I)}

First, we shall show that Statement (I) holds for $p \in (0, 1)$
such that $C_{j_p} = 1 - p$. For this purpose,  we need to show that
{\small \be \la{doitrev} 1  \leq \limsup_{\vep \to 0} \f{ \bs{\ga}
(\om)  } {  \ga  (p, \vep)  } \leq 1 + \ro_p \qqu \tx{ for any} \;
\om \in \li \{ \lim_{\vep \to 0} \wh{\bs{p}} = p \ri \}. \ee}
 To show {\small $\limsup_{\vep \to 0} \f{ \bs{\ga} (\om)  } {
\ga (p, \vep)  } \geq 1$}, note that $C_{{s - \ell_\vep} + 1} < 1 -
p = C_{s - \ell_\vep} < C_{s - \ell_\vep - 1}$ as a direct
consequence of the definition of $\ell_\vep$ and the assumption that
$C_{j_p} = 1 - p$. By the first three statements of  Lemma
\ref{lem79rev}, we have $\lim_{\vep \to 0} z_{{\ell}} > p$ for all
$\ell \leq {\ell_\vep} - 1$.  Noting that $\lim_{\vep \to 0}
\wh{\bs{p}} (\om) = p$,  we have $\wh{\bs{p}} (\om) < z_\ell$ for
all $\ell \leq {\ell_\vep} - 1$ and it follows from the definition
of the sampling scheme that $\bs{\ga} (\om) \geq \ga_{\ell_\vep}$ if
$\vep > 0$ is small enough.
 By Lemma \ref{lem81rev} and noting that $\ka_p = 1$ if $C_{j_p} = 1 - p$, we have
 {\small $\limsup_{\vep \to 0} \f{ \bs{\ga} (\om)  } {  \ga  (p, \vep)  }
\geq \lim_{\vep \to 0}  \f{ \ga_{{\ell_\vep}}} { \ga  (p, \vep) } =
\ka_p = 1$}.

To show {\small $\limsup_{\vep \to 0} \f{ \bs{\ga} (\om)  } {  \ga
(p, \vep)  } \leq 1 + \ro_p$}, we shall consider two cases: (i)
${\ell_\vep} = s - 1$; (ii) ${\ell_\vep} < s - 1$. In the case of
${\ell_\vep} = s - 1$, it must be true that $\bs{\ga} (\om) \leq
\ga_s = \ga_{\ell_\vep + 1}$.  Hence,  $\limsup_{\vep \to 0} \f{
\bs{\ga} (\om) } {  \ga  (p, \vep)  } \leq \lim_{\vep \to 0}  \f{
\ga_{\ell_\vep + 1} } { \ga  (p, \vep) } = \lim_{\vep \to 0}  \f{
\ga_{\ell_\vep } } { \ga  (p, \vep) } \times \lim_{\vep \to 0}  \f{
\ga_{\ell_\vep + 1} } { \ga_{\ell_\vep } } = 1 + \ro_p$. In the case
of ${\ell_\vep} < s - 1$, it follows from the first three statements
of Lemma \ref{lem79rev}
 that $\lim_{\vep \to 0} z_{{\ell_\vep} + 1} < p$, which implies that
 $z_{{\ell_\vep} + 1} < p, \; \wh{\bs{p}} (\om) > z_{{\ell_\vep} + 1}$, and thus $\bs{\ga} (\om) \leq
\ga_{{\ell_\vep} + 1}$ for small enough $\vep
> 0$. Therefore, {\small $\limsup_{\vep \to 0} \f{ \bs{\ga} (\om)  } {  \ga  (p, \vep)  } \leq
\lim_{\vep \to 0}  \f{ \ga_{{\ell_\vep} + 1}} { \ga  (p, \vep) } = 1
+ \ro_p$}.  This establishes (\ref{doitrev}) and it follows that $\{
1 \leq  \limsup_{\vep \to 0} \f{ \bs{\ga}  } {  \ga  (p, \vep)  }
\leq 1 + \ro_p \} \supseteq \li \{ \lim_{\vep \to 0} \wh{\bs{p}} = p
\ri \}$.  According to the strong law of large numbers, we have $1
\geq \Pr \{ 1 \leq  \limsup_{\vep \to 0} \f{ \bs{\ga}  } {  \ga  (p,
\vep)  } \leq 1 + \ro_p \} \geq \Pr \li \{ \lim_{\vep \to 0}
\wh{\bs{p}} = p \ri \} = 1$.  This proves that Statement (I) holds
for $p \in (0,1)$ such that $C_{j_p} = 1 - p$.

\bsk

Next, we shall show that Statement (I) holds for $p \in (0, 1)$ such
that $C_{j_p} > 1 - p$. Note that $C_{s - \ell_\vep + 1} < 1 - p <
C_{s - \ell_\vep}$ as a direct consequence of the definition of
$\ell_\vep$ and the assumption that $C_{j_p} > 1 - p$.  By the first
three statements of Lemma \ref{lem79rev}, we have $\lim_{\vep \to 0}
z_{{\ell_\vep} - 1} > p$ and thus $z_{{\ell}} > p$ for all $\ell
\leq {\ell_\vep} - 1$ provided that $\vep > 0$ is sufficiently
small. Therefore, for any $\om \in \li \{ \lim_{\vep \to 0}
\wh{\bs{p}} = p \ri \}$, we have $\wh{\bs{p}} (\om) < z_\ell$ for
all $\ell \leq {\ell_\vep} - 1$ and consequently, $\bs{\ga} (\om)
\geq \ga_{{\ell_\vep}}$ provided that $\vep > 0$ is sufficiently
small. On the other hand, we claim that $\bs{\ga} (\om) \leq
\ga_{{\ell_\vep}}$. Such claim can be justified by investigating two
cases.  In the case of ${\ell_\vep} = s$, it is trivially true that
$\bs{\ga} (\om) \leq \ga_{{\ell_\vep}}$.  In the case of $\ell_\vep
< s$, we have $ p > \lim_{\vep \to 0} z_{{\ell_\vep}}$ and thus $p >
z_{{\ell_\vep}}$  provided that $\vep > 0$ is sufficiently small.
Therefore, for any $\om \in \li \{ \lim_{\vep \to 0} \wh{\bs{p}} = p
\ri \}$, we have $\wh{\bs{p}} (\om)
> z_{{\ell_\vep}}$ and consequently, $\bs{\ga} (\om) \leq \ga_{{\ell_\vep}}$ provided that
$\vep > 0$ is sufficiently small.  This proves the claim and it
follows that $\bs{\ga} (\om) = \ga_{{\ell_\vep}}$ if $\vep
> 0$ is small enough. Applying Lemma \ref{lem81rev}, we have
{\small $\lim_{\vep \to 0} \f{ \bs{\ga} (\om)  }
 {  \ga  (p, \vep)  } = \lim_{\vep \to 0}  \f{ \ga_{{\ell_\vep}}}
 { \ga  (p, \vep) } = \ka_p$}, which implies that
{\small $\{ \lim_{\vep \to 0} \f{ \bs{\ga} } {  \ga  (p, \vep)  } =
\ka_p \} \supseteq \li \{ \lim_{\vep \to 0} \wh{\bs{p}} = p \ri
\}$}. It follows from the strong law of large numbers that {\small
$1 \geq \Pr \{ \lim_{\vep \to 0} \f{ \bs{\ga} } {  \ga  (p, \vep)  }
= \ka_p \} \geq \Pr \{ \lim_{\vep \to 0} \wh{\bs{p}} = p \} = 1$}
and thus {\small $\Pr \{ \lim_{\vep \to 0} \f{ \bs{\ga} } {  \ga (p,
\vep) } = \ka_p \} = 1$}.  Since $1 \leq \ka_p \leq 1 + \ro_p$, we
have that $\Pr \{ 1 \leq  \limsup_{\vep \to 0} \f{ \bs{\ga}  } { \ga
(p, \vep)  } \leq 1 + \ro_p \} = 1$ is of course true.  This proves
that Statement (I) also holds for $p \in (0, 1)$ such that $C_{j_p}
> 1 - p$. The proof of Statement (I) is thus completed.

\subsubsection{Proof of Statement (III)}

First, we shall consider $p \in (0, 1)$ such that $C_{j_p} = 1 - p$.
In this case, it is evident that $\ell_\vep < s$. It follows from
Lemma \ref{Defrev} and the definition of the sampling scheme that
$\lim_{\vep \to 0} \Pr \{ \bs{l} > {{\ell_\vep} + 1} \} \leq
\lim_{\vep \to 0} \Pr \{ \bs{D}_{{\ell_\vep} + 1} = 0 \} = 0$ and
$\lim_{\vep \to 0} \Pr \{ \bs{l} < {\ell_\vep} \} \leq \lim_{\vep
\to 0} \sum_{\ell =  1}^{ {\ell_\vep} - 1} \Pr \{ \bs{D}_\ell = 1 \}
= 0$. Since \bee  \limsup_{\vep \to 0} \Pr \{ | \wh{\bs{p}} - p |
\geq \vep p \} & \leq & \lim_{\vep \to 0} \li [ \Pr \{ |
\wh{\bs{p}}_{\ell_\vep} - p | \geq \vep p, \; \bs{l} = {\ell_\vep}
\} + \Pr \{ | \wh{\bs{p}}_{\ell_\vep + 1} - p | \geq \vep p, \;
\bs{l} =
{{\ell_\vep} + 1} \} \ri ]\\
&  & + \lim_{\vep \to 0} \Pr \{ \bs{l} < {\ell_\vep} \} + \lim_{\vep
\to 0} \Pr \{ \bs{l} > {{\ell_\vep} + 1} \} \eee and
\[
\liminf_{\vep \to 0} \Pr \{ | \wh{\bs{p}} - p | \geq \vep p \}  \geq
\lim_{\vep \to 0} \li [ \Pr \{ | \wh{\bs{p}}_{\ell_\vep} - p | \geq
\vep p, \; \bs{l} = {\ell_\vep} \} + \Pr \{ | \wh{\bs{p}}_{\ell_\vep
+ 1} - p | \geq \vep p, \; \bs{l} = {{\ell_\vep} + 1} \} \ri ],
\]
we have {\small \[ \lim_{\vep \to 0} \Pr \{ | \wh{\bs{p}} - p | \geq
\vep p \}   =  \lim_{\vep \to 0} \li [ \Pr \{ |
\wh{\bs{p}}_{\ell_\vep} - p | \geq \vep p, \; \bs{l} = {\ell_\vep}
\} + \Pr \{ | \wh{\bs{p}}_{\ell_\vep + 1} - p | \geq \vep p, \;
\bs{l} = {{\ell_\vep} + 1} \} \ri ].
 \]}
By Lemma \ref{limpleminv}, we have  \[ \lim_{\vep \to 0} \Pr \{ |
\wh{\bs{p}} - p | \geq \vep p
 \} =  \Pr \li \{ U \geq d \ri \} + \Pr \li \{ |U + \sq{\ro_p} V | \geq (1 + \ro_p) d, \; U < \nu d \ri \}
 \] for $p \in (0, 1)$ such that $C_{j_p} = 1 - p$.
 As a consequence of Lemma \ref{2dGU}, Statement (III) must be
 true for $p \in (0, 1)$ such that $C_{j_p} = 1 - p$.

Next, we shall consider $p \in (0, 1)$ such that $C_{j_p}
> 1 - p$.  Note that $C_{s - {\ell_\vep} + 1} < 1 - p < C_{s -
{\ell_\vep}}$. Since {\small $U_{\ell_\vep} = \li ( \f{p}
{\wh{\bs{p}}_{{\ell_\vep}}}  - 1 \ri ) \sq{ \f{ \ga_{\ell_\vep}}{1 -
p} }$} converges in distribution to a standard Gaussian variable
$U$, {\small $\lim_{\vep \to 0} \vep \sq{\f{\ga_{\ell_\vep}}{1 - p}
} = d \sq{\ka_p}$} and $\lim_{\vep \to 0} \Pr \{ \bs{\ga} =
\ga_{\ell_\vep} \} = 1$ as can be seen from Statement (I),  we have
{\small \bee \lim_{\vep \to 0} \Pr \{ | \wh{\bs{p}} - p | \geq \vep
p \} & = & \lim_{\vep \to 0} \Pr \{ |
\wh{\bs{p}}_{\ell_\vep} - p | \geq \vep p \}\\
& = & \lim_{\vep \to 0} \Pr \li \{ U_{\ell_\vep}  \geq \f{ \vep} {1
- \vep} \sq{\f{\ga_{\ell_\vep}}{1 - p} } \ri \} + \lim_{\vep \to 0}
\Pr \li \{ U_{\ell_\vep} \leq - \f{ \vep} {1 + \vep}
\sq{\f{\ga_{\ell_\vep}}{1 - p} }
\ri \}\\
& = & \lim_{\vep \to 0} \Pr \li \{ | U_{\ell_\vep} | \geq \vep
\sq{\f{\ga_{\ell_\vep}}{1 - p} } \ri \} =  \Pr \{ | U | \geq d
\sq{\ka_p} \} \eee} and consequently, $\lim_{\vep \to 0} \Pr \{ |
\wh{\bs{p}} - p | < \vep p \} \geq 2 \Phi \li ( d \sq{\ka_p} \ri ) -
1 > 1 - 2 \ze \de$ for $p \in (0, 1)$ such that $C_{j_p} > 1 - p$.
This proves Statement (III).

\bsk

Finally, we would like to note that Statement (II) can be shown by
employing Lemma \ref{Defrev} and similar argument as the proof of
Statement (II) of Theorem \ref{Bino_Asp_Analysis}.

\subsection{Proof of Theorem  \ref{Bino_Asp_Analysis_Noninverse} }
\la{App_Bino_Asp_Analysis_Noninverse}

We need some preliminary results.

\beL \la{bino_noninverse_lem} If $\vep$ is sufficiently small, then
the following statements hold true.

(I): For $\ell = 1, \cd, \tau$, there exists a unique number
$z_{\ell} \in [0, 1]$ such that $n_{\ell} = \f{ \ln ( \ze \de_\ell )
} { \mscr{M}_{\mrm{B}} ( z_{\ell}, \; \f{z_{\ell}}{1 + \vep} ) }$.

(II): $z_{\ell}$ is monotonically decreasing with respect to $\ell$
no greater than $\tau$.

(III): $\lim_{\vep \to 0} z_{\ell} = \li [ 1 + ( 1 - \f{1}{p^*} )
C_{\tau - \ell} \ri ]^{-1}$ for $1 \leq \ell \leq \tau$, where the
limit is taken under the restriction that $\ell - \tau$ is fixed
with respect to $\vep$.

(IV): $\{ \bs{D}_{\ell} = 0 \} = \{ \wh{\bs{p}}_{\ell} < z_{\ell}
\}$ for $\ell = 1, \cd, \tau$.

\eeL

\bsk

{\bf Proof of Statement (I)}: By the definition of $n_\ell$, we have
{\small $0 < \f{ \ln ( \ze \de )  } { \mscr{M}_{\mrm{B}} ( 1,
\f{1}{1 + \vep} ) } \leq \li \lc \f{ \ln ( \ze \de )  } {
\mscr{M}_{\mrm{B}} ( 1, \f{1}{1 + \vep} ) } \ri \rc = n_1 \leq
n_{\ell}$} for sufficiently small $\vep > 0$.  Hence, $\ln \f{1}{1 +
\vep} = \mscr{M}_{\mrm{B}} ( 1, \f{1}{1 + \vep} ) \leq \f{ \ln ( \ze
\de ) } { n_{\ell} } < 0$ for small enough $\vep
> 0$. Recall that $\mscr{M}_{\mrm{B}} ( z, \f{z}{1 + \vep} )$ is
monotonically decreasing with respect to $z \in (0, 1)$ as asserted
by Lemma \ref{lemm23}.  Invoking the intermediate value theorem, we
have that there exists a unique number $z_{\ell} \in (0, 1]$ such
that $\mscr{M}_{\mrm{B}} ( z_{\ell}, \f{z_{\ell}}{1 + \vep} ) = \f{
\ln ( \ze \de ) } { n_{\ell} }$, which implies Statement (I).

\bsk

{\bf Proof of Statement (II)}: Since $n_{\ell}$ is monotonically
increasing with respect to $\ell$ for sufficiently small $\vep > 0$,
we have that $\mscr{M}_{\mrm{B}} ( z_{\ell}, \f{z_{\ell}}{1 + \vep}
)$ is monotonically increasing with respect to $\ell \leq \tau$ for
sufficiently small $\vep > 0$. Recalling that $\mscr{M}_{\mrm{B}} (
z, \f{z}{1 + \vep} )$ is monotonically decreasing with respect to $z
\in (0, 1)$ as asserted by Lemma \ref{lemm23}, we have that
$z_{\ell}$ is monotonically decreasing with respect to $\ell \leq
\tau$. This establishes Statement (II).

\bsk

{\bf Proof of Statement (III)}:

For simplicity of notations, let {\small $b_{\ell} = \li [ 1 + ( 1 -
\f{1}{p^*} ) C_{\tau - \ell} \ri ]^{-1}$} for $\ell = 1, 2, \cd,
\tau$. Then, it can be checked that $\f{p^* (1 - b_{\ell}) }{ b_\ell
(1 - p^*)} = C_{\tau - \ell}$ for $1 \leq \ell \leq \tau$. By the
definition of sample sizes, we have \be \la{goodenoninv} \f{
\mscr{M}_{\mrm{B}} (z_{\ell}, \f{z_{\ell}}{1 + \vep}) } { \vep^2
b_{\ell} \sh [ 2 (b_{\ell} - 1)] } = \f{\ln ( \ze \de) }{n_{\ell} }
\times \f{2 (p^* - 1) C_{\tau - \ell}} { p^* \vep^2} = 1 + o(1) \ee
for $\ell = 1, \cd, \tau$, where
\[
n_\ell = \f{ \ln (\ze \de)  } { \mscr{M}_{\mrm{B}} (z_{\ell},
\f{z_{\ell}}{1 + \vep}) } = \f{[1 + o(1) ] C_{\tau - \ell} \ln (\ze
\de)}{ \mscr{M}_{\mrm{B}} (p^*, \f{p^*}{1 + \vep}) }.
\]
We claim that $\se < z_{\ell} < 1$ for $\se \in (0, b_{\ell})$ if
$\vep > 0$ is small enough.  To prove this claim, we use a
contradiction method. Suppose the claim is not true, then there
exists a set, denoted by $S_\vep$,  of infinite many values of
$\vep$ such that $z_{\ell} \leq \se$ for $\vep \in S_\vep$.  Hence,
by (\ref{goodenoninv}) and the fact that $\mscr{M}_{\mrm{B}} (z,
\f{z}{1 + \vep})$ is monotonically decreasing with respect to $z \in
(0, 1)$ as asserted by Lemma \ref{lemm23}, we have
\[
1 + o(1) = \f{\mscr{M}_{\mrm{B}} (z_{\ell}, \f{z_{\ell}}{1 + \vep})
} { \vep^2 b_{\ell} \sh [ 2 (b_{\ell} - 1)] } \leq \f{
\mscr{M}_{\mrm{B}} (\se, \f{\se}{1 + \vep}) } { \vep^2 b_{\ell} \sh
[ 2 (b_{\ell} - 1 )] } = \f{ \vep^2 \se \sh [ 2 (1 - \se)] + o
(\vep^2) } { \vep^2 b_{\ell} \sh [ 2 (1 - b_{\ell})]  } = \f{ \se (1
- b_{\ell} ) } { b_{\ell} (1 - \se) } + o(1)
\]
for small enough $\vep \in S_\vep$,  which implies {\small $\f{ \se
(1 - b_{\ell} ) } { b_{\ell} (1 - \se) } \geq 1$}, contradicting to
the fact that {\small $\f{ \se (1 - b_{\ell} ) } { b_{\ell} (1 -
\se) } < 1$}.  This proves our claim.  In a similar manner, we can
show that $0 < z_{\ell} < \se^\prime$ for $\se^\prime \in (
b_{\ell}, 1)$ if $\vep > 0$ is small enough.  By (\ref{goodenoninv})
and applying Lemma \ref{lem32T} based on the established condition
that $\se < z_{\ell} < \se^\prime$ for small enough $\vep > 0$, we
have {\small $\f{\mscr{M}_{\mrm{B}} (z_{\ell}, \f{z_{\ell}}{1 +
\vep}) } { \vep^2 b_{\ell} \sh [ 2 (b_{\ell} - 1)] } = \f{ \vep^2
z_{\ell} \sh [ 2 (1 - z_{\ell})]  + o (\vep^2) } { \vep^2 b_{\ell}
\sh [ 2 (1 - b_{\ell} )] } = 1 + o (1)$}, which implies {\small
$\f{z_{\ell} }{ 1 - z_{\ell} } - \f{b_{\ell}} { 1 - b_{\ell} } =
o(1)$} and consequently $\lim_{\vep \to 0} z_{\ell} = b_{\ell}$.
This proves Statement (III).

\bsk

{\bf Proof of Statement (IV)}: Noting that {\small
$\mscr{M}_{\mrm{B}} (z, \f{z}{1 + \vep} )$} is monotonically
decreasing with respect to $z \in (0, 1)$ as asserted by Lemma
\ref{lemm23}, we have {\small $\{ \bs{D}_{\ell} = 0 \}  =  \{
\mscr{M}_{\mrm{B}} ( \wh{\bs{p}}_{\ell}, \; \f{\wh{\bs{p}}_{\ell}
}{1 + \vep}  ) > \f{\ln (\ze \de_\ell)} {n_{\ell}} \} = \li \{
\wh{\bs{p}}_{\ell} < z_\ell \ri \}$} as claimed by statement (IV).

\beL \la{lem34a_noninv} Let $\ell_\vep = \tau  - j_p$.  Then, \be
\la{thenlemninv} \lim_{\vep \to 0} \sum_{\ell = 1}^{\ell_\vep - 1}
n_\ell \Pr \{ \bs{D}_\ell = 1 \} = 0, \qqu \lim_{\vep \to 0}
\sum_{\ell = \ell_\vep + 1}^\tau n_\ell \Pr \{ \bs{D}_{\ell} = 0 \}
= 0 \ee for $p \in (p^*, 1)$. Moreover, $\lim_{\vep \to 0}
n_{\ell_\vep} \Pr \{ \bs{D}_{\ell_\vep} = 0 \} = 0$ for $p \in (p^*,
1)$ such that $C_{j_p} > r (p)$. \eeL

\bpf

For simplicity of notations, let $b_\ell = \lim_{\vep \to 0}
z_{\ell}$ for $1 \leq \ell \leq \tau$.

First, we shall show that (\ref{thenlemninv}) holds for $p \in (p^*,
1)$. By the definition of $\ell_\vep$,
 we have $r (p)  > C_{\tau - \ell_\vep + 1}$.
 Making use of the first three statements of Lemma \ref{bino_noninverse_lem},
 we have that {\small $z_\ell > \f{p + b_{\ell_\vep -
1}}{2} > p$} for all $\ell \leq \ell_\vep - 1$ if $\vep$ is
sufficiently small. By the last statement of Lemma
\ref{bino_noninverse_lem} and
 using Chernoff bound, we have {\small \bee \Pr \{ \bs{D}_{\ell} = 1 \}  =
 \Pr \{ \wh{\bs{p}}_{\ell} \geq z_{\ell} \} \leq  \Pr \li \{
\wh{\bs{p}}_{\ell} > \f{p + b_{\ell_\vep - 1}}{2} \ri \} \leq \exp
\li ( - 2 n_\ell \li ( \f{p - b_{\ell_\vep - 1}}{2} \ri )^2 \ri )
\eee} for all $\ell \leq \ell_\vep - 1$ provided that $\vep
> 0$ is small enough.  By the definition of $\ell_\vep$, we have
\[
b_{\ell_\vep - 1} = \li [ 1 + \li ( 1 - \f{1}{p^*} \ri ) C_{\tau -
\ell_\vep + 1} \ri ]^{-1}  > p,
\]
which implies that {\small $\li ( \f{p - b_{\ell_\vep - 1}}{2} \ri
)^2$} is a positive constant independent of $\vep > 0$ provided that
$\vep > 0$ is small enough.   Hence, {\small $\lim_{\vep \to 0}
\sum_{\ell = 1}^{\ell_\vep - 1} n_\ell \Pr \{ \bs{D}_\ell = 1 \} =
0$} as a result of Lemma \ref{lem31a}.

Similarly, it can be seen from the definition of $\ell_\vep$ that $r
(p) < C_{\tau - \ell_\vep - 1}$. Making use of the first three
statements of Lemma \ref{bino_noninverse_lem}, we have that {\small
$z_\ell < \f{ p + b_{\ell_\vep + 1}}{2} < p$} for $\ell_\vep + 1
\leq \ell \leq \tau$ if $\vep$ is sufficiently small. By the last
statement of Lemma \ref{bino_noninverse_lem} and using Chernoff
bound, we have {\small
\[ \Pr \{ \bs{D}_{\ell} = 0 \} = \Pr \{ \wh{\bs{p}}_{\ell}
< z_{\ell} \} \leq \Pr \li \{ \wh{\bs{p}}_{\ell} < \f{ p +
b_{\ell_\vep + 1}}{2} \ri \} \leq \exp \li ( - 2 n_\ell \li ( \f{ p
- b_{\ell_\vep + 1}}{2} \ri )^2 \ri )
\]}
for $\ell_\vep + 1 \leq \ell \leq \tau$ provided that $\vep
> 0$ is small enough.  As a consequence of the definition of $\ell_\vep$, we have that
$b_{\ell_\vep + 1}$ is smaller than $p$ and is independent of $\vep
> 0$. Therefore, we can apply Lemma \ref{lem31a} to conclude that {\small
$\lim_{\vep \to 0} \sum_{\ell = \ell_\vep + 1}^\tau n_\ell \Pr \{
\bs{D}_\ell = 0 \} = 0$}.

\bsk

Second, we shall show that $\lim_{\vep \to 0} n_{\ell_\vep} \Pr \{
\bs{D}_{\ell_\vep} = 0 \} = 0$  for $p \in (p^*, 1)$ such that
$C_{j_p} > r (p)$.  Clearly, $r (p) < C_{\tau - \ell_\vep}$ because
of the definition of $\ell_\vep$. Making use of the first three
statements of Lemma \ref{bino_noninverse_lem}, we have {\small
$z_{\ell_\vep} < \f{p + b_{\ell_\vep}}{2} < p$} if $\vep$ is
sufficiently small. By the last statement of Lemma
\ref{bino_noninverse_lem} and using Chernoff bound, we have {\small
\[ \Pr \{ \bs{D}_{\ell_\vep} = 0 \} =  \Pr \{
 \wh{\bs{p}}_{\ell_\vep} < z_{\ell_\vep} \} \leq \Pr \li \{
\wh{\bs{p}}_{\ell_\vep} < \f{p + b_{\ell_\vep}}{2} \ri \}  \leq \exp
\li ( - 2 n_{\ell_\vep} \li ( \f{p - b_{\ell_\vep}}{2} \ri )^2 \ri )
\]} for small enough $\vep > 0$. As a consequence of the definition of $\ell_\vep$, we have that
$b_{\ell_\vep}$ is smaller than $p$ and is independent of $\vep
> 0$.  It follows
that $\lim_{\vep \to 0} n_{\ell_\vep} \Pr \{ \bs{D}_{\ell_\vep} = 0
\} = 0$.

\epf

\beL

\la{ineqinv}
 $\lim_{\vep \to 0} \sum_{\ell =  \tau +
1}^{\iy} n_\ell \Pr \{ \bs{l} = \ell \} = 0$ for any $p \in (p^*,
1)$.

\eeL

\bpf

Recalling that the sample sizes $n_1, n_2, \cd$ are chosen as the
ascending arrangement of all distinct elements of the set defined by
(\ref{defss}), we have that
\[
n_\ell = \li \lc  \f{ C_{\tau - \ell} \ln (\ze \de) }{
\mscr{M}_{\mrm{B}} (p^*, \f{p^*}{1 + \vep} ) }  \ri \rc, \qqu \ell =
1, 2, \cd
\]
for small enough $\vep \in (0, 1)$.   By the assumption that
$\inf_{i \in \bb{Z} } \f{C_{i - 1}}{C_{i} } = 1 + \udl{\ro} > 1$, we
have that
\[
n_\ell >  (1 + \udl{\ro} )^{\ell - \tau - 1} \f{ \ln(\ze \de) }{
\mscr{M}_{\mrm{B}} (p^*, \f{p^*}{1 + \vep} ) }, \qqu \ell = \tau +
1, \tau + 2, \cd
\]
for small enough $\vep \in (0, 1)$. So, we have shown that there
exists a number $\vep^* \in (0, 1)$ such that
\[
n_\ell \mscr{M}_{\mrm{B}} \li ( p^*, \f{p^*}{1 + \vep} \ri ) < (1 +
\udl{\ro} )^{\ell - \tau - 1} \ln (\ze \de), \qqu \ell = \tau + 1,
\tau + 2, \cd
\]
for any $\vep \in (0, \vep^*)$. Observing that there exist a
positive integer $\ka^*$ such that $(1 + \udl{\ro} )^{\ell - \tau -
1} \ln (\ze \de) < \ln (\ze \de) - (\ell - \tau) \ln 2 = \ln (\ze
\de_\ell)$ for any $\ell \geq \tau + \ka^*$, we have that there
exists a positive integer $\ka^*$ independent of $\vep$ such that
{\small $\mscr{M}_{\mrm{B}} ( p^*, \f{p^*}{1 + \vep} ) < \f{\ln (\ze
\de_\ell)}{n_\ell}$} for $\ell \geq \tau + \ka^*$ and $0 < \vep <
\vep^*$.  Recall that $\mscr{M}_{\mrm{B}} (z, \f{z}{1 + \vep} )$ is
monotonically decreasing with respect to $z \in (0, 1)$ as asserted
by Lemma \ref{lemm23}.  For $\ell \geq \tau + \ka^*$ and $0 < \vep <
\vep^*$, as a result of {\small $\f{\ln (\ze \de_\ell)}{n_\ell} >
\mscr{M}_{\mrm{B}} ( p^*, \f{p^*}{1 + \vep} ) > \mscr{M}_{\mrm{B}} (
1, \f{1}{1 + \vep} ) = \ln \f{1}{1 + \vep} $}, there exists a unique
number $z_\ell \in [0, 1]$ such that {\small $\mscr{M}_{\mrm{B}}
(z_\ell, \f{z_\ell}{1 + \vep} ) = \f{\ln (\ze \de_\ell)}{n_\ell} >
\mscr{M}_{\mrm{B}} ( p^*, \f{p^*}{1 + \vep} )$}.   Moreover, it must
be true that $z_\ell < p^*$ for $\ell \geq \tau + \ka^*$ and $\vep
\in (0, \vep^*)$. Therefore, for small enough $\vep \in (0,
\vep^*)$, we have {\small \bee \sum_{\ell = \tau + 1}^{\iy} n_\ell
\Pr \{ \bs{l} = \ell \} & = & \sum_{\ell = \tau + 1}^{\tau + \ka^*}
n_\ell \Pr \{ \bs{l} = \ell \} + \sum_{\ell = \tau + \ka^* + 1}^\iy
n_\ell \Pr
\{ \bs{l} = \ell \}\\
& \leq & \sum_{\ell =  \tau + 1}^{\tau + \ka^*} n_\ell \Pr \{
\bs{D}_\tau  = 0 \} + \sum_{\ell =  \tau + \ka^* + 1}^\iy  n_\ell
\Pr \{ \bs{D}_{\ell - 1} = 0 \} \\
& = & \sum_{\ell =  \tau + 1}^{\tau + \ka^*} n_\ell \Pr \{
\bs{D}_\tau = 0 \} + \sum_{\ell =  \tau + \ka^*}^\iy  n_{\ell+1} \Pr
\{ \bs{D}_{\ell} = 0 \} \\
& < & k^* (1 + \ovl{\ro})^{k^*} n_\tau \Pr \{ \bs{D}_\tau = 0 \} +
(1 + \ovl{\ro}) \sum_{\ell = \tau + \ka^*}^\iy  n_{\ell} \Pr \{
\bs{D}_{\ell} = 0 \} \\
& \leq  & k^* (1 + \ovl{\ro})^{k^*} n_\tau \Pr \{ \wh{\bs{p}}_\tau <
z_\tau \} + (1 + \ovl{\ro}) \sum_{\ell = \tau + \ka^*}^\iy  n_{\ell}
\Pr \{ \wh{\bs{p}}_{\ell} < z_\ell \}\\
& \leq  & k^* (1 + \ovl{\ro})^{k^*} n_\tau \Pr \li \{
\wh{\bs{p}}_\tau < \f{p^* + p}{2} \ri \} + (1 + \ovl{\ro})
\sum_{\ell = \tau + \ka^*}^\iy
n_{\ell} \Pr \{ \wh{\bs{p}}_{\ell} < p^* \}\\
& \leq  & k^* (1 + \ovl{\ro})^{k^*} n_\tau \exp  \li ( -
\f{n_\tau}{2} ( p - p^*  )^2 \ri ) + (1 + \ovl{\ro}) \sum_{\ell =
\tau + \ka^*}^\iy n_{\ell} \exp( - 2 n_\ell (p - p^*)^2) \to 0 \eee}
as $\vep \to 0$, where we have used Chernoff bound and the
assumption that $\sup_{i \in \bb{Z}} \f{C_{i - 1}}{C_i} = 1 +
\ovl{\ro} < \iy$. This completes the proof of the lemma. \epf

\beL

\la{lem35ainv} $\lim_{\vep \to 0} \f{ n_{\ell_\vep} } {
\mcal{N}_{\mrm{r}}  (p, \vep) } = \ka_p, \; \lim_{\vep \to 0} \f{
\vep p} {\sq{p ( 1 - p) \sh n_{\ell_\vep} }}   = d \sq{\ka_p}$.

\eeL

\bpf

By the definition of sample sizes, it can be readily shown that
{\small $\lim_{\vep \to 0} \f{ 2 ( 1 - p^* ) C_{ \tau - \ell} \ln
\f{1}{\ze \de} } { p^* \vep^2 n_\ell} = 1$} for any $\ell \geq 1$
and it follows that {\small \bee \lim_{\vep \to 0} \f{ n_{\ell_\vep}
} { \mcal{N}_{\mrm{r}}  (p, \vep) }  & = & \lim_{\vep \to 0}  \f{
\mscr{M}_{\mrm{B}} ( p, \f{p}{1 + \vep} )  }{  \ln (\ze \de) }
\times \f{ 2 ( 1 - p^* ) C_{\tau - \ell_\vep} } {
p^* \vep^2} \ln \f{1}{\ze \de}\\
& = & \lim_{\vep \to 0}  \li [ \f{ p \vep^2 }{ 2 (1 - p ) } + o
(\vep^2) \ri ] \times \f{ 2 ( 1 - p^* ) C_{ \tau - {\ell_\vep}  } } { p^* \vep^2}\\
&  = & \f{p (1 - p^*) C_{\tau - \ell_\vep}}{p^* (1 - p)} = \f{ p (1
- p^*) C_{j_p} }{p^* (1 - p)} = \ka_p, \eee }
 {\small
\bee \lim_{\vep \to 0} \f{ \vep p } {\sq{p ( 1 - p) \sh
n_{\ell_\vep} }}  & = & \lim_{\vep \to 0}  \vep p \sq{ \f{ 2 ( 1 -
p^* ) C_{\tau - \ell_\vep} } { p ( 1 - p) p^* \vep^2} \ln \f{1}{\ze
\de} } = d \sq{ \f{p ( 1 - p^*) C_{\tau - \ell_\vep}}{p^* (1 - p)} } \\
& =  & d \sq{ \f{p (1 - p^*) C_{j_p} }{p^* (1 - p)} } = d
\sq{\ka_p}. \eee}

\epf

\beL \la{limplemnoninv} Let $U$ and $V$ be independent Gaussian
random variables with zero means and unit variances.  Then, for $p
\in (p^*, 1)$ such that $C_{j_p} = r(p)$, \bee &  & \lim_{\vep \to
0} \Pr \{ \bs{l} = \ell_\vep \} = 1 - \lim_{\vep \to 0} \Pr \{
\bs{l} = \ell_\vep + 1 \} =  1 - \Phi \li (  \nu  d \ri ), \\ &  &
\lim_{\vep \to 0} \li [ \Pr \{ | \wh{\bs{p}}_{\ell_\vep} - p | \geq
\vep p, \; \bs{l}  = \ell_\vep \} + \Pr \{ | \wh{\bs{p}}_{\ell_\vep
+ 1} - p | \geq \vep p, \;
\bs{l} =  \ell_\vep + 1 \} \ri ]\\
&   & \qqu \qqu \qqu \qu = \Pr \li \{ U \geq d \ri \} + \Pr \li \{ |U + \sq{\ro_p} V | \geq (1 + \ro_p) d, \; U < \nu d \ri \}. \eee \eeL

Lemma \ref{limplemnoninv} can be shown by a similar method as that
of Lemma \ref{limpleminv}.

\subsubsection{Proof of Statement (I)}

First, we shall show that Statement (I) holds for $p \in (p^*, 1)$
such that {\small $C_{ j_p} = r (p)$}. For this purpose,  we need to
show that {\small \be \la{doitinv} 1 \leq \limsup_{\vep \to 0} \f{
\mbf{n} (\om)  } {  \mcal{N}_{\mrm{r}} (p, \vep)  } \leq 1 + \ro_p
\qqu \tx{ for any} \; \om \in \li \{ \lim_{\vep \to 0} \wh{\bs{p}} =
p \ri \}. \ee}
 To show {\small $\limsup_{\vep \to 0}
\f{ \mbf{n} (\om)  } {  \mcal{N}_{\mrm{r}}  (p, \vep)  } \geq 1$},
note that $C_{\tau - \ell_\vep + 1} < r (p) = C_{\tau - \ell_\vep} <
C_{\tau - \ell_\vep - 1}$ as a direct consequence of the definitions
of $\ell_\vep$ and $j_p$. By the first three statements of  Lemma
\ref{bino_noninverse_lem}, we have $\lim_{\vep \to 0} z_{{\ell}} >
p$ for all $\ell \leq {\ell_\vep} - 1$.  Noting that $\lim_{\vep \to
0} \wh{\bs{p}} (\om) = p$,  we have $\wh{\bs{p}} (\om) < z_\ell$ for
all $\ell \leq {\ell_\vep} - 1$ and it follows from the definition
of the sampling scheme that
 $n_{\ell_\vep} \leq \mbf{n} (\om)$ if $\vep > 0$ is small enough.
 By Lemma \ref{lem35ainv} and noting that $\ka_p = 1$ if {\small $C_{ j_p} = r (p)$}, we have
 {\small $\limsup_{\vep \to 0} \f{ \mbf{n} (\om)  } {  \mcal{N}_{\mrm{r}}  (p, \vep)  }
\geq \lim_{\vep \to 0}  \f{ n_{{\ell_\vep}}} { \mcal{N}_{\mrm{r}}
(p, \vep) } = \ka_p = 1$}.

To show {\small $\limsup_{\vep \to 0} \f{ \mbf{n} (\om)  } {
\mcal{N}_{\mrm{r}}  (p, \vep)  } \leq 1 + \ro_p$}, note that
${\ell_\vep} + 1 \leq \tau$ as a result of $p^* < p < 1$ and the
assumption that {\small $C_{ j_p} = r (p)$}. By virtue of Lemma
\ref{bino_noninverse_lem}, we have $\lim_{\vep \to 0} z_{{\ell_\vep}
+ 1} < p$, which implies
 $\wh{\bs{p}} (\om) > z_{{\ell_\vep} + 1}$ and thus $\mbf{n} (\om) \leq
n_{{\ell_\vep} + 1}$ for small enough $\vep \in (0, 1)$. Therefore,
{\small $\limsup_{\vep \to 0} \f{ \mbf{n} (\om)  } {
\mcal{N}_{\mrm{r}} (p, \vep)  } \leq \lim_{\vep \to 0}  \f{
n_{{\ell_\vep} + 1}} { \mcal{N}_{\mrm{r}} (p, \vep) } = \lim_{\vep
\to 0}  \f{ n_{{\ell_\vep} + 1}} { n_{{\ell_\vep} }} \times
\lim_{\vep \to 0} \f{ n_{{\ell_\vep} }} { \mcal{N}_{\mrm{r}}  (p,
\vep) } = 1 + \ro_p$}. This establishes (\ref{doitinv}), which
implies $\{ 1 \leq \limsup_{\vep \to 0} \f{ \mbf{n}  } {
\mcal{N}_{\mrm{r}} (p, \vep) } \leq 1 + \ro_p \} \supseteq \li \{
\lim_{\vep \to 0} \wh{\bs{p}} = p \ri \}$. Applying the strong law
of large numbers, we have $1 \geq \Pr \{ 1 \leq \limsup_{\vep \to 0}
\f{ \mbf{n}  } { \mcal{N}_{\mrm{r}}  (p, \vep) } \leq 1 + \ro_p \}
\geq \Pr \li \{ \lim_{\vep \to 0} \wh{\bs{p}} = p \ri \} = 1$.  This
proves that Statement (I) holds for $p \in (p^*, 1)$ such that
{\small $C_{ j_p} = r (p)$}.

\bsk

Next, we shall show that Statement (I) holds for $p \in (p^*, 1)$
such that $C_{ j_p} > r(p)$. Note that $C_{\tau - \ell_\vep + 1} < r
(p) < C_{\tau - \ell_\vep}$ as a direct consequence of the
definition of $\ell_\vep$ and the assumption that {\small $C_{ j_p}
> r (p)$}. By the first three statements of
Lemma \ref{bino_noninverse_lem}, we have $\lim_{\vep \to 0}
z_{{\ell_\vep} - 1} > p > \lim_{\vep \to 0} z_{{\ell_\vep}}$ and
thus $z_{{\ell}}
> p > z_{{\ell_\vep}}$ for all $\ell \leq {\ell_\vep} - 1$ provided
that $\vep \in (0, 1)$ is sufficiently small. Therefore, for any
$\om \in \li \{ \lim_{\vep \to 0} \wh{\bs{p}} = p \ri \}$, we have
$z_\ell > \wh{\bs{p}} (\om) > z_{{\ell_\vep}}$ for all $\ell \leq
{\ell_\vep} - 1$ and consequently, $\mbf{n} (\om) = n_{{\ell_\vep}}$
provided that $\vep \in (0, 1)$ is sufficiently small. Applying
Lemma \ref{lem35ainv},  we have {\small $\lim_{\vep \to 0} \f{
\mbf{n} (\om)  }
 {  \mcal{N}_{\mrm{r}}  (p, \vep)  } = \lim_{\vep \to 0}  \f{ n_{{\ell_\vep}}}
 { \mcal{N}_{\mrm{r}}  (p, \vep) } = \ka_p$}, which implies that
{\small $\{ \lim_{\vep \to 0} \f{ \mbf{n} } {  \mcal{N}_{\mrm{r}}
(p, \vep)  } = \ka_p \} \supseteq \li \{ \lim_{\vep \to 0}
\wh{\bs{p}} = p \ri \}$}.  It follows from the strong law of large
numbers that {\small $1 \geq \Pr \{ \lim_{\vep \to 0} \f{ \mbf{n} }
{ \mcal{N}_{\mrm{r}}  (p, \vep)  } = \ka_p \} \geq \Pr \{ \lim_{\vep
\to 0} \wh{\bs{p}} = p \}$} and thus {\small $\Pr \{ \lim_{\vep \to
0} \f{ \mbf{n} } {  \mcal{N}_{\mrm{r}}  (p, \vep)  } = \ka_p \} =
1$}. Since $1 \leq \ka_p \leq 1 + \ro_p$, we have that $\Pr \{ 1
\leq \limsup_{\vep \to 0} \f{ \mbf{n} } { \mcal{N}_{\mrm{r}}  (p,
\vep) } \leq 1 + \ro_p \}$ is of course true. This proves that
Statement (I) holds for $p \in (p^*, 1)$ such that {\small $C_{ j_p}
> r (p)$}.   The proof of Statement (I) is thus completed.

\subsubsection{Proof of Statement (II)}

In the sequel, we will consider the asymptotic value of $\f{ \bb{E}
[ \mbf{n} ] } { \mcal{N}_{\mrm{r}}  (p, \vep)  }$ in three steps.
First, we shall show Statement (II) for $p \in (p^*, 1)$ such that
$C_{j_p} = r(p)$. Clearly, $\ell_\vep < \tau$. By the definition of
the sampling scheme, we have {\small \bee \bb{E} [ \mbf{n} ] & = &
\sum_{\ell = 1}^{{\ell_\vep} - 1} n_\ell \Pr \{ \bs{l} = \ell \} +
\sum_{\ell = {\ell_\vep} + 2}^\tau n_\ell \Pr \{ \bs{l} = \ell \} +
\sum_{\ell =  \tau + 1}^{\iy} n_\ell \Pr \{ \bs{l} = \ell \}\\
&  & + n_{\ell_\vep} \Pr \{ \bs{l} = {\ell_\vep} \} + n_{{\ell_\vep}
+ 1} \Pr \{ \bs{l} = {{\ell_\vep} + 1} \}\\
 & \leq &  \sum_{\ell =
1}^{{\ell_\vep} - 1} n_\ell \Pr \{ \bs{D}_\ell = 1 \} + \sum_{\ell =
{\ell_\vep} + 1}^{\tau - 1} n_{\ell + 1} \Pr \{ \bs{D}_\ell = 0 \} +
\sum_{\ell =  \tau + 1}^{\iy}
n_\ell \Pr \{ \bs{l} = \ell \}\\
&   & + n_{\ell_\vep} \Pr \{ \bs{l} = {\ell_\vep} \} +
n_{{\ell_\vep} + 1} \Pr \{ \bs{l} = {{\ell_\vep} + 1} \} \eee} and
{\small $\bb{E} [ \mbf{n} ] \geq n_{\ell_\vep} \Pr \{ \bs{l} =
{\ell_\vep} \} + n_{\ell_\vep + 1} \Pr \{ \bs{l} = {{\ell_\vep} + 1}
\}$}.   Making use of Lemmas \ref{lem34a_noninv}, \ref{ineqinv} and
the assumption that $\sup_{\ell > 0} \f{n_{\ell + 1}}{n_{\ell}} <
\iy$ for small enough $\vep > 0$, we have {\small \bee &   &
\lim_{\vep \to 0} \li [ \sum_{\ell = 1}^{{\ell_\vep} - 1} n_\ell \Pr
\{ \bs{D}_\ell = 1 \} + \sum_{\ell = {\ell_\vep} + 1}^{\tau - 1}
n_{\ell + 1} \Pr \{ \bs{D}_\ell = 0  \} + \sum_{\ell =  \tau +
1}^{\iy} n_\ell \Pr \{
\bs{l} = \ell \} \ri ]\\
&   & \leq  \lim_{\vep \to 0} \li [ \sum_{\ell = 1}^{{\ell_\vep} -
1} n_\ell \Pr \{ \bs{D}_\ell = 1 \} + \sup_{\ell > 0} \f{n_{\ell +
1}}{n_{\ell}}  \sum_{\ell = {\ell_\vep} + 1}^{\tau - 1} n_{\ell} \Pr
\{ \bs{D}_\ell = 0  \} + \sum_{\ell = \tau + 1}^{\iy} n_\ell \Pr \{
\bs{l} = \ell \} \ri ] = 0. \eee} Therefore, {\small \[
\limsup_{\vep \to 0} \f{ \bb{E} [ \mbf{n} ] } { \mcal{N}_{\mrm{r}}
(p, \vep)  } \leq  \lim_{\vep \to 0}  \f{ n_{\ell_\vep} \Pr \{
\bs{l} = {\ell_\vep} \} + n_{\ell_\vep + 1} \Pr \{ \bs{l} =
{{\ell_\vep} + 1} \} } { \mcal{N}_{\mrm{r}}  (p, \vep) } \]} and
{\small \[ \liminf_{\vep \to 0} \f{ \bb{E} [ \mbf{n} ] } {
\mcal{N}_{\mrm{r}} (p, \vep)  } \geq  \lim_{\vep \to 0}  \f{
n_{\ell_\vep} \Pr \{ \bs{l} = {\ell_\vep} \} + n_{\ell_\vep + 1} \Pr
\{ \bs{l} = {{\ell_\vep} + 1} \} } { \mcal{N}_{\mrm{r}}  (p, \vep)
}.
\]}
It follows that {\small \[ \lim_{\vep \to 0} \f{ \bb{E} [ \mbf{n} ]
} { \mcal{N}_{\mrm{r}} (p, \vep)  } =  \lim_{\vep \to 0}  \f{
n_{\ell_\vep} \Pr \{ \bs{l} = {\ell_\vep} \} + n_{\ell_\vep + 1} \Pr
\{ \bs{l} = {{\ell_\vep} + 1} \} } { \mcal{N}_{\mrm{r}}  (p, \vep) }
\]}
Using Lemma \ref{limplemnoninv} and the result $\lim_{\vep \to 0}
\f{ n_{\ell_\vep} } { \mcal{N}_{\mrm{r}} (p, \vep) } = \ka_p$ as
asserted by Lemma \ref{lem35ainv},  we have \bee \lim_{\vep \to 0}
\f{ n_{\ell_\vep} \Pr \{ \bs{l} = \ell_\vep \} + n_{{\ell_\vep} + 1}
\Pr \{ \bs{l} = \ell_\vep + 1 \} } { \mcal{N}_{\mrm{r}}  (p, \vep) }
& = & \lim_{\vep \to 0} \f{ n_{\ell_\vep} [ 1 - \Phi (  \nu d ) ] +
n_{\ell_\vep + 1} \Phi ( \nu d ) } { \mcal{N}_{\mrm{r}} (p, \vep) } \\
& = & 1  + \ro_p \Phi \li (  \nu d \ri ). \eee

\bsk

Second, we shall show Statement (II) for $p \in (p^*, 1)$ such that
$C_{j_p} > r(p)$.  Note that \bee \bb{E} [ \mbf{n} ] & = &
\sum_{\ell = 1}^{{\ell_\vep} - 1} n_\ell \Pr \{ \bs{l} = \ell \} +
\sum_{\ell =  {\ell_\vep} + 1}^\tau n_\ell \Pr \{ \bs{l} = \ell \} +
n_{\ell_\vep} \Pr \{ \bs{l} = {\ell_\vep} \} + \sum_{\ell = \tau +
1}^{\iy} n_\ell \Pr \{
\bs{l} = \ell \}\\
& \leq & \sum_{\ell =  1}^{{\ell_\vep} - 1} n_\ell \Pr \{
\bs{D}_\ell = 1 \} + \sum_{\ell =  {\ell_\vep}}^{\tau - 1} n_{\ell +
1} \Pr \{ \bs{D}_\ell = 0  \} + n_{\ell_\vep} + \sum_{\ell =  \tau +
1}^{\iy} n_\ell \Pr \{ \bs{l} = \ell \} \eee and {\small $\bb{E} [
\mbf{n} ] \geq n_{\ell_\vep} \Pr \{ \bs{l} = {\ell_\vep} \} \geq
n_{\ell_\vep}  \li ( 1 -  \sum_{\ell =  1}^{{\ell_\vep} - 1} \Pr \{
\bs{D}_\ell = 1 \} - \Pr \{  \bs{D}_{\ell_\vep} = 0 \} \ri )$}.
Therefore, by Lemma \ref{lem34a_noninv}, {\small \bee \limsup_{\vep
\to 0} \f{ \bb{E} [ \mbf{n} ] } { \mcal{N}_{\mrm{r}} (p, \vep)  } &
\leq & \lim_{\vep \to 0}  \f{ \sum_{\ell = 1}^{{\ell_\vep} - 1}
n_\ell \Pr \{ \bs{D}_\ell = 1 \} + \sum_{\ell = {\ell_\vep}}^{\tau -
1} n_{\ell + 1} \Pr \{ \bs{D}_\ell = 0 \} + n_{\ell_\vep} +
\sum_{\ell =  \tau + 1}^{\iy} n_\ell \Pr \{ \bs{l} = \ell \} } {
\mcal{N}_{\mrm{r}}  (p,
\vep) }\\
& = & \lim_{\vep \to 0}  \f{ n_{\ell_\vep} } { \mcal{N}_{\mrm{r}}
(p, \vep) } = \ka_p, \eee}
\[
\liminf_{\vep \to 0} \f{ \bb{E} [ \mbf{n} ] } { \mcal{N}_{\mrm{r}}
(p, \vep)  } \geq \lim_{\vep \to 0} \f{ n_{\ell_\vep}  \li ( 1 -
\sum_{\ell = 1}^{{\ell_\vep} - 1} \Pr \{ \bs{D}_\ell = 1 \} - \Pr \{
\bs{D}_{\ell_\vep} = 0 \} \ri ) } { \mcal{N}_{\mrm{r}}  (p, \vep) }
= \lim_{\vep \to 0}  \f{ n_{\ell_\vep} } { \mcal{N}_{\mrm{r}}  (p,
\vep) } = \ka_p.
\]
So, $\lim_{\vep \to 0} \f{ \bb{E} [ \mbf{n} ] } { \mcal{N}_{\mrm{r}}
(p, \vep)  } = \ka_p$ for $p \in (p^*, 1)$ such that $C_{j_p} >
r(p)$.

From the preceding  analysis, we have obtained $\limsup_{\vep \to 0}
\f{ \bb{E} [ \mbf{n} ] } { \mcal{N}_{\mrm{r}}  (p, \vep) }$ for all
$p \in (p^*, 1)$. Hence, statement (II) is established by making use
of this result and the fact that
\[
\lim_{\vep \to 0} \f{ \bb{E} [ \mbf{n} ] } { \mcal{N}_{\mrm{f}} (p,
\vep) } = \lim_{\vep \to 0} \f{ \mcal{N}_{\mrm{r}} (p, \vep) } {
\mcal{N}_{\mrm{f}} (p, \vep) } \times \lim_{\vep \to 0} \f{ \bb{E} [
\mbf{n} ] } { \mcal{N}_{\mrm{r}} (p, \vep) } = \f{ 2 \ln \f{1}{\ze
\de} } { \mcal{Z}_{\ze \de}^2  } \times \lim_{\vep \to 0} \f{ \bb{E}
[ \mbf{n} ] } { \mcal{N}_{\mrm{r}} (p, \vep) }.
\]

\subsubsection{Proof of Statement (III)}

First, we shall consider $p \in (p^*, 1)$ such that $C_{j_p} =
r(p)$. In this case, it is evident that $\ell_\vep < \tau$. By the
definition of the sampling scheme, we have that $\Pr \{ \bs{l} >
{{\ell_\vep} + 1} \} \leq \Pr \{ \bs{D}_{{\ell_\vep} + 1} = 0 \}$
and that $\Pr \{ \bs{l} = \ell \} \leq \Pr \{ \bs{D}_\ell = 1 \}$
for $\ell < {\ell_\vep}$.  As a result of Lemma \ref{lem34a_noninv},
we have $\lim_{\vep \to 0} \Pr \{ \bs{l} > {{\ell_\vep} + 1} \} \leq
\lim_{\vep \to 0} \Pr \{ \bs{D}_{{\ell_\vep} + 1} = 0 \} = 0$ and
$\lim_{\vep \to 0} \Pr \{ \bs{l} < {\ell_\vep} \} \leq \lim_{\vep
\to 0} \sum_{\ell =  1}^{ {\ell_\vep} - 1} \Pr \{ \bs{D}_\ell = 1 \}
= 0$.  Since \bee \limsup_{\vep \to 0} \Pr \{ | \wh{\bs{p}} - p |
\geq \vep p \} & \leq & \lim_{\vep \to 0} \li [ \Pr \{ |
\wh{\bs{p}}_{\ell_\vep} - p | \geq \vep p, \; \bs{l} = {\ell_\vep}
\} + \Pr \{ | \wh{\bs{p}}_{\ell_\vep + 1} - p | \geq \vep p, \;
\bs{l} =
{{\ell_\vep} + 1} \} \ri ]\\
&  & + \lim_{\vep \to 0} \Pr \{ \bs{l} < {\ell_\vep} \} + \lim_{\vep
\to 0} \Pr \{ \bs{l} > {{\ell_\vep} + 1} \} \eee and
\[
\liminf_{\vep \to 0} \Pr \{ | \wh{\bs{p}} - p | \geq \vep p \}  \geq
\lim_{\vep \to 0} \li [ \Pr \{ | \wh{\bs{p}}_{\ell_\vep} - p | \geq
\vep p, \; \bs{l} = {\ell_\vep} \} + \Pr \{ | \wh{\bs{p}}_{\ell_\vep
+ 1} - p | \geq \vep p, \; \bs{l} = {{\ell_\vep} + 1} \} \ri ],
\]
we have {\small \[ \lim_{\vep \to 0} \Pr \{ | \wh{\bs{p}} - p | \geq
\vep p \}   =  \lim_{\vep \to 0} \li [ \Pr \{ |
\wh{\bs{p}}_{\ell_\vep} - p | \geq \vep p, \; \bs{l} = {\ell_\vep}
\} + \Pr \{ | \wh{\bs{p}}_{\ell_\vep + 1} - p | \geq \vep p, \;
\bs{l} = {{\ell_\vep} + 1} \} \ri ].
 \]}
By Lemma \ref{limplemnoninv}, we have {\small $\lim_{\vep \to 0} \Pr \{ | \wh{\bs{p}} - p | \geq \vep p \} =  \Pr \li \{ U \geq d \ri \} + \Pr
\li \{ |U + \sq{\ro_p} V | \geq (1 + \ro_p) d, \; U < \nu d \ri \}$} for $p \in (p^*, 1)$ such that $C_{j_p} = r(p)$.  As a consequence of Lemma
\ref{2dGU}, Statement (III) must be true for $p \in (p^*, 1)$ such that $C_{j_p} = r(p)$.

Next, we shall consider $p \in (p^*, 1)$ such that $C_{j_p} > r(p)$.
Applying Lemma \ref{lem34a_noninv}, we have  {\small \bee & &
\lim_{\vep \to 0} \Pr \{ \bs{l} < {{\ell_\vep}} \} \leq \lim_{\vep
\to 0} \sum_{\ell = 1}^{{\ell_\vep} - 1} \Pr \{ \bs{D}_\ell = 1 \}
\leq \lim_{\vep \to 0} \sum_{\ell = 1}^{{\ell_\vep} - 1}
n_\ell \Pr \{ \bs{D}_\ell = 1 \} = 0,\\
&   & \lim_{\vep \to 0} \Pr \{ \bs{l} > {{\ell_\vep}} \} \leq
\lim_{\vep \to 0} \Pr \{ \bs{D}_{\ell_\vep} = 0 \} \leq \lim_{\vep
\to 0} n_{{\ell_\vep}} \Pr \{ \bs{D}_{\ell_\vep} = 0 \} = 0 \eee}
and thus $\lim_{\vep \to 0} \Pr \{ \bs{l} \neq {\ell_\vep}  \} = 0$.
Note that $\Pr \{ | \wh{\bs{p}} - p | \geq \vep p \} = \Pr \{ |
\wh{\bs{p}}_{\ell_\vep} - p | \geq \vep p, \; \bs{l} = {\ell_\vep}
\} + \Pr \{ | \wh{\bs{p}} - p | \geq \vep p, \; \bs{l} \neq
{\ell_\vep} \}$ and, as a result of the central limit theorem,
{\small $U_{\ell_\vep} = \f{ \wh{\bs{p}}_{\ell_\vep} - p }{ \sq{ p (
1 - p) \sh n_{\ell_\vep}} }$} converges in distribution to a
standard Gaussian variable $U$. Hence,  {\small
\[ \lim_{\vep \to 0} \Pr \{ | \wh{\bs{p}} - p | \geq \vep p \} =
\lim_{\vep \to 0} \Pr \{ | \wh{\bs{p}}_{\ell_\vep} - p | \geq \vep p
\} = \lim_{\vep \to 0} \Pr \li \{ |U_{\ell_\vep}| \geq \f{ \vep p} {
\sq{ p ( 1 - p) \sh n_{\ell_\vep}} } \ri \} = \Pr \{ |U| \geq d
\sq{\ka_p} \}
\]}
and $\lim_{\vep \to 0} \Pr \{ | \wh{\bs{p}} - p | < \vep p \} = \Pr
\{ |U| < d \sq{\ka_p} \} = 2 \Phi (d \sq{\ka_p}) - 1 > 2 \Phi (d) -
1 > 1 - 2 \ze \de$ for $p \in (p^*, 1)$.  Here we have used the fact
that  {\small $\Phi (z)
>  1 -  e^{ - \f{z^2}{2} }$} and {\small $ \Phi(d)
= \Phi ( \sq{ 2 \ln \f{1}{\ze \de} } ) > 1 - \ze \de$}.  This proves
Statement (III).

\subsection{Proof of Theorem \ref{Bino_mix_CDF_CH} } \la{App_Bino_mix_CDF_CH}

We need some preliminary results.

\beL \la{CCon}

{\small $\mscr{M}_{\mrm{B}} (z, z-\vep)$} is monotonically
increasing with respect to $z \in (\vep, p + \vep)$ provided that
{\small $0 < \vep < \f{35}{94}$} and {\small $0 < p < \f{1}{2} -
\f{12}{35} \vep$}. \eeL

 \bpf

Define {\small $g(\vep, p) = \f{ \vep } {p(1-p)  } + \ln \f{ p (1 -
p - \vep) } { (p + \vep) (1 - p) }$} for $0 < p < 1$ and $0 < \vep <
1 - p$.  We shall first show that $g(\vep, p) > 0$ if $0 < \vep <
\f{35}{94}$ and {\small $0 < p < \f{1}{2} - \f{12}{35} \vep$}.

Let $\f{1}{3} < k < 1$ and $0 < \vep \leq \f{1}{2 (1 + k)}$. It can
be shown by tedious computation that {\small $\f{ \pa g \li ( \vep,
\f{1}{2} - k \vep \ri ) } { \pa \vep }  = \f{ 16 \vep^2 \li [  3 k
-1 - 4 (1 - k) k^2 \vep^2 \ri ] } { (1 - 4 k^2 \vep^2)^2 [ 1 - 4
(k-1)^2 \vep^2 ] }$}, which implies that $g \li ( \vep, \f{1}{2} - k
\vep \ri )$ is monotonically increasing with respect to {\small
$\vep \in \li ( 0, \; \f{1}{2 k} \sq{ \f{2} { 1 - k } - 3} \ri )$}
and is monotonically decreasing with respect to {\small $\vep \in
\li ( \f{1}{2 k} \sq{ \f{2} { 1 - k } - 3}, \; \f{1}{2 (1 + k)} \ri
]$}. Since $g \li (0, \f{1}{2} \ri ) = 0$, we have that $g \li (
\vep, \f{1}{2} - k \vep \ri )$ is positive for $0 < \vep \leq
\f{1}{2 (1 + k)}$ if {\small $g \li ( \vep, \f{1}{2} - k \vep \ri
)$} is positive for $\vep = \f{1}{2 (1 + k)}$. For $\vep = \f{1}{2
(1 + k)}$ with $k = \f{12}{35}$, we have {\small $g \li ( \vep,
\f{1}{2} - k \vep \ri ) = 1 + \f{1}{2 k +1} - \ln \li (2 + \f{1}{k}
\ri ) = 1 + \f{35}{59} - \ln \li ( 2 + \f{35}{12} \ri )$}, which is
positive because $e \times e^{ \f{35}{59} } > 2.718 \times
\sum_{i=0}^4 \f{1}{i!} \li ( \f{35}{59} \ri )^i  > 2 + \f{35}{12}$.
It follows that $g \li ( \vep, \f{1}{2} - \f{12}{35} \vep \ri )$ is
positive for any $\vep \in \li (0, \f{35}{94} \ri )$. Since {\small
$\f{ \pa g(\vep, p) } { \pa p } = - \vep^2 \li [  \f{ 1  } { (p +
\vep) p^2  } + \f{ 1 } { (1 - p - \vep) (1 - p)^2 } \ri ]$} is
negative, we have that $g(\vep, p)$ is positive for $0 < \vep <
\f{35}{94}$ if $0 < p < \f{1}{2} - \f{12}{35}
 \vep$.

 Finally, the lemma is established by verifying that {\small $\f{\pa^2
\mscr{M}_{\mrm{B}}(z, z - \vep) } { \pa z^2 } = - \vep^2 \li [
\f{1}{z (z - \vep)^2 }  + \f{1}{ (1 - z) ( 1 - z + \vep)^2  } \ri ]
< 0$} for any $z \in (\vep,1)$ and that {\small $\li . \f{\pa
\mscr{M}_{\mrm{B}} (z, z-\vep)} { \pa z } \ri |_{z = p + \vep} =
g(\vep, p)$}.

 \epf

\beL \la{lemm22}

$\mscr{M}_{\mrm{B}} (p - \vep, p) < \mscr{M}_{\mrm{B}} (p + \vep, p)
< - 2 \vep^2$ for $0 < \vep < p < \f{1}{2} < 1 - \vep$.

\eeL

\bpf The lemma follows from the facts that $\mscr{M}_{\mrm{B}} (p-
\vep, p) - \mscr{M}_{\mrm{B}} (p+ \vep, p) = 0$ for $\vep = 0$ and
that
\[
\f{\pa [ \mscr{M}_{\mrm{B}} (p-\vep, p) - \mscr{M}_{\mrm{B}} (p+
\vep, p)  ]} { \pa \vep } =  \ln \li [1 + \f{\vep^2} { p^2 } \f{2 p
- 1}{(1-p)^2 - \vep^2} \ri ],
\]
where the right side is negative for $0 < \vep < p < \f{1}{2} < 1 -
\vep$.  By Lemma \ref{lemax}, we have $\mscr{M}_{\mrm{B}} (p + \vep,
p) < - 2 \vep^2$ for $0 < \vep < p < \f{1}{2} < 1 - \vep$.  This
completes the proof of the lemma.

\epf

\beL

\la{lemm24} $\mscr{M}_{\mrm{B}}(z, \f{z}{1 - \vep} )$ is
monotonically decreasing from $0$ to $- \iy$ as $z$ increases from
$0$ to $1-\vep$.

\eeL

\bpf

The lemma can be shown by verifying that {\small \[ \lim_{z \to 0}
\mscr{M}_{\mrm{B}} \li (z, \f{z}{1 - \vep} \ri ) = 0, \qu \lim_{z
\to 1-\vep} \mscr{M}_{\mrm{B}} \li (z, \f{z}{1 - \vep} \ri ) = -
\iy, \qu \lim_{z \to 0} \f{\pa }{ \pa z} \mscr{M}_{\mrm{B}} \li (z,
\f{z}{1 - \vep} \ri ) = \ln \f{1}{1 - \vep} - \f{\vep}{1 - \vep} < 0
\]}
and {\small $\f{\pa^2 }{ \pa z^2} \mscr{M}_{\mrm{B}} \li (z, \f{z}{1
- \vep} \ri ) = \f{ \vep^2 } { (z - 1) (1 - \vep - z)^2 } < 0$} for
any $z \in (0,1 - \vep)$.

 \epf

 \beL
 \la{lemm25}

 $\mscr{M}_{\mrm{B}} ( z, \f{z}{1 + \vep} ) >
\mscr{M}_{\mrm{B}} ( z, \f{z}{1 - \vep} )$ for $0 < z < 1 - \vep <
1$. \eeL

\bpf

The lemma follows from the facts that $\mscr{M}_{\mrm{B}} ( z,
\f{z}{1 + \vep} ) - \mscr{M}_{\mrm{B}} ( z, \f{z}{1 - \vep} ) = 0$
for $\vep = 0$ and that
\[
\f{ \pa } {\pa \vep} \li [ \mscr{M}_{\mrm{B}} \li ( z, \f{z}{1 +
\vep} \ri ) - \mscr{M}_{\mrm{B}} \li ( z, \f{z}{1 - \vep} \ri ) \ri
] = \f{ 2 \vep^2 z (2 - z) } {(1 - \vep^2) [(1 - z)^2 - \vep^2] }
> 0
\]
for $z \in (0, 1 - \vep)$.

\epf

\beL \la{mixDS1}

$\li \{ \mscr{M}_{\mrm{B}} \li ( \wh{\bs{p}}_s, \mscr{L} (
\wh{\bs{p}}_s ) \ri ) \leq \f{\ln (\ze \de) } { n_s }, \;
\mscr{M}_{\mrm{B}} \li ( \wh{\bs{p}}_s, \mscr{U} ( \wh{\bs{p}}_s )
\ri ) \leq \f{\ln (\ze \de) } { n_s } \ri \}$ is a sure event.

\eeL

\bpf

For simplicity of notations, we denote $p^\star =
\f{\vep_a}{\vep_r}$. In order to show the lemma, it suffices to show
\bel & & \li \{
 \mscr{M}_{\mrm{B}} \li ( \wh{\bs{p}}_s, \f{ \wh{\bs{p}}_s }{ 1 - \vep_r} \ri )
> \f{\ln (\ze \de)}{n_s}, \;
\wh{\bs{p}}_s > p^\star - \vep_a \ri \} = \emptyset, \la{ep81BC}\\
&   & \li \{  \mscr{M}_{\mrm{B}} ( \wh{\bs{p}}_s, \wh{\bs{p}}_s +
\vep_a)
> \f{\ln (\ze \de)}{n_s}, \; \wh{\bs{p}}_s \leq
p^\star - \vep_a \ri \} = \emptyset, \la{ep82BC}\\
&   & \li \{  \mscr{M}_{\mrm{B}} \li ( \wh{\bs{p}}_s, \f{
\wh{\bs{p}}_s }{ 1 + \vep_r} \ri ) > \f{\ln (\ze \de)}{n_s}, \;
\wh{\bs{p}}_s > p^\star + \vep_a \ri \} = \emptyset, \la{ep83BC}\\
&   & \li \{  \mscr{M}_{\mrm{B}} ( \wh{\bs{p}}_s, \wh{\bs{p}}_s -
\vep_a)
> \f{\ln (\ze \de)}{n_s}, \; \wh{\bs{p}}_s \leq
p^\star + \vep_a \ri \} = \emptyset. \la{ep84BC} \eel By the
definition of $n_s$,  we have {\small $n_s \geq \li \lc  \f{ \ln
(\ze \de) } {  \mscr{M}_{\mrm{B}} \li ( p^\star + \vep_a, p^\star
\ri ) } \ri \rc \geq \f{ \ln (\ze \de) } {  \mscr{M}_{\mrm{B}} \li (
p^\star + \vep_a, p^\star \ri ) } $}.  By the assumption on $\vep_a$
and $\vep_r$, we have $0 < \vep_a < p^\star < \f{1}{2} < 1 -
\vep_a$. Hence, by Lemma \ref{lemm22}, we have $ \mscr{M}_{\mrm{B}}
\li ( p^\star - \vep_a, p^\star \ri ) <  \mscr{M}_{\mrm{B}} \li (
p^\star + \vep_a, p^\star \ri ) < 0$ and it follows that \be
\la{rela8BC} \f{\ln (\ze \de)}{n_s} \geq  \mscr{M}_{\mrm{B}} \li (
p^\star + \vep_a, p^\star \ri ) >  \mscr{M}_{\mrm{B}} \li ( p^\star
- \vep_a, p^\star \ri ). \ee By (\ref{rela8BC}), {\small \be
\la{ep989BC} \li \{  \mscr{M}_{\mrm{B}} \li ( \wh{\bs{p}}_s, \f{
\wh{\bs{p}}_s }{ 1 - \vep_r} \ri )
> \f{\ln (\ze \de)}{n_s}, \; \wh{\bs{p}}_s > p^\star - \vep_a \ri \}
\subseteq \li \{  \mscr{M}_{\mrm{B}} \li ( \wh{\bs{p}}_s, \f{
\wh{\bs{p}}_s }{ 1 - \vep_r} \ri ) >  \mscr{M}_{\mrm{B}} \li (
p^\star - \vep_a, p^\star \ri ), \; \wh{\bs{p}}_s > p^\star - \vep_a
\ri \}. \ee} Noting that $ \mscr{M}_{\mrm{B}} \li ( p^\star -
\vep_a, p^\star \ri ) =  \mscr{M}_{\mrm{B}} \li ( p^\star - \vep_a,
\f{p^\star - \vep_a}{1 - \vep_r} \ri )$ and making use of the fact
that $ \mscr{M}_{\mrm{B}}(z, \f{z}{1 - \vep} )$ is monotonically
decreasing with respect to $z \in (0, 1 - \vep)$ as asserted by
Lemma \ref{lemm24}, we have \be \la{ep99BC} \li \{
\mscr{M}_{\mrm{B}} \li ( \wh{\bs{p}}_s, \f{ \wh{\bs{p}}_s }{ 1 -
\vep_r} \ri )
>  \mscr{M}_{\mrm{B}} \li ( p^\star - \vep_a, p^\star \ri ) \ri \} = \{
\wh{\bs{p}}_s < p^\star - \vep_a \}. \ee Combining (\ref{ep989BC})
and (\ref{ep99BC}) yields (\ref{ep81BC}). By (\ref{rela8BC}),
{\small \be \la{eq9998BC} \li \{  \mscr{M}_{\mrm{B}} (
\wh{\bs{p}}_s, \wh{\bs{p}}_s + \vep_a)
> \f{\ln (\ze \de)}{n_s}, \; \wh{\bs{p}}_s \leq
p^\star - \vep_a \ri \} \subseteq \li \{  \mscr{M}_{\mrm{B}} (
\wh{\bs{p}}_s, \wh{\bs{p}}_s + \vep_a) >  \mscr{M}_{\mrm{B}} \li (
p^\star - \vep_a, p^\star \ri ), \; \wh{\bs{p}}_s \leq p^\star -
\vep_a \ri \}. \ee} By the assumption on $\vep_a$ and $\vep_r$, we
have $p^\star - \vep_a < \f{1}{2} - \vep_a$.   Recalling the fact
that $ \mscr{M}_{\mrm{B}}(z, z + \vep)$ is monotonically increasing
with respect to $z \in (0, \f{1}{2} - \vep)$ as asserted by Lemma
\ref{lemm21}, we have that the event in the right-hand side of
(\ref{eq9998BC}) is an impossible event and consequently,
(\ref{ep82BC}) is established. By (\ref{rela8BC}), {\small \be
\la{ep98998BC} \li \{  \mscr{M}_{\mrm{B}} \li ( \wh{\bs{p}}_s, \f{
\wh{\bs{p}}_s }{ 1 + \vep_r} \ri ) > \f{\ln (\ze \de)}{n_s}, \;
\wh{\bs{p}}_s > p^\star + \vep_a \ri \} = \li \{  \mscr{M}_{\mrm{B}}
\li ( \wh{\bs{p}}_s, \f{ \wh{\bs{p}}_s }{ 1 + \vep_r} \ri )
>  \mscr{M}_{\mrm{B}} \li ( p^\star + \vep_a, p^\star \ri ), \; \wh{\bs{p}}_s >
p^\star + \vep_a \ri \}.  \ee} Noting that $ \mscr{M}_{\mrm{B}} \li
( p^\star + \vep_a, p^\star \ri ) =  \mscr{M}_{\mrm{B}} \li (
p^\star + \vep_a, \f{p^\star + \vep_a}{1 + \vep_r} \ri )$ and making
use of the fact that $ \mscr{M}_{\mrm{B}}(z, \f{z}{1 + \vep} )$ is
monotonically decreasing with respect to $z \in (0, 1)$ as asserted
by Lemma \ref{lemm23}, we have \be \la{ep9998BC} \li \{
\mscr{M}_{\mrm{B}} \li ( \wh{\bs{p}}_s, \f{ \wh{\bs{p}}_s }{ 1 +
\vep_r} \ri )
>  \mscr{M}_{\mrm{B}} \li ( p^\star + \vep_a, p^\star \ri ) \ri \} = \{
\wh{\bs{p}}_s < p^\star + \vep_a \}. \ee Combining (\ref{ep98998BC})
and (\ref{ep9998BC}) yields (\ref{ep83BC}). By (\ref{rela8BC}),
{\small \be \la{eq999833BC} \li \{  \mscr{M}_{\mrm{B}} (
\wh{\bs{p}}_s, \wh{\bs{p}}_s - \vep_a)
> \f{\ln (\ze \de)}{n_s}, \; \wh{\bs{p}}_s \leq
p^\star + \vep_a \ri \} \subseteq \li \{  \mscr{M}_{\mrm{B}} (
\wh{\bs{p}}_s, \wh{\bs{p}}_s - \vep_a) >  \mscr{M}_{\mrm{B}} \li (
p^\star + \vep_a, p^\star \ri ), \; \wh{\bs{p}}_s \leq p^\star +
\vep_a \ri \}. \ee} By the assumption on $\vep_a$ and $\vep_r$, we
have that $ \mscr{M}_{\mrm{B}}(z, z - \vep)$ is monotonically
increasing with respect to $z \in (\vep_a, p^\star + \vep_a)$ as a
result of Lemma \ref{CCon}. Hence,  the event in the right-hand side
of (\ref{eq999833BC}) is an impossible event and consequently,
(\ref{ep84BC}) is established. This completes the proof of the
lemma.

\epf

\bsk

Now we are in a position to prove Theorem \ref{Bino_mix_CDF_CH}. If
the multistage sampling scheme follows a stopping rule derived from
Chernoff bounds,  then $\{ \bs{D}_s = 1 \}$ is a sure event as a
result of Lemma \ref{mixDS1}. Recall that $\exp( \mscr{M}_{\mrm{B}}
(z, p) )$  is equal to $\mcal{F} (z, p)$ and $\mcal{G} (z, p)$
respectively for the cases of $z \leq p$ and $z \geq p$.  Moreover,
$\wh{\bs{p}}_\ell$ is a ULE of $p$ for $\ell = 1, \cd, s$. So, the
sampling scheme satisfies all the requirements described in
Corollary \ref{Monotone_third}, from which Theorem
\ref{Bino_mix_CDF_CH} immediately follows.

If the multistage sampling scheme follows a stopping rule derived
from CDF $\&$ CCDF,  then, by Lemmas \ref{decb}, we have {\small
\bee &  & 1 \geq \Pr \{ G_{\wh{\bs{p}}_s} ( \wh{\bs{p}}_s,  \mscr{L}
( \wh{\bs{p}}_s ) ) \leq \ze \de_s \} = \Pr \li \{ 1 - S_{\mrm{B}}
(K_s - 1, n_s, \mscr{L} ( \wh{\bs{p}}_s ) ) \leq \ze \de \ri \}\\
&    & \qqu \qqu  \qqu \qqu  \qqu  \qqu \qqu \geq \Pr \li \{ n_s
\mscr{M}_{\mrm{B}} \li ( \wh{\bs{p}}_s,
\mscr{L} ( \wh{\bs{p}}_s ) \ri ) \leq \ln (\ze \de) \ri \} = 1,\\
&   & 1 \geq \Pr \{ F_{\wh{\bs{p}}_s} ( \wh{\bs{p}}_s,  \mscr{U} (
\wh{\bs{p}}_s ) ) \leq \ze \de_s \} = \Pr \li \{  S_{\mrm{B}} ( K_s,
n_s, \mscr{U} ( \wh{\bs{p}}_s ) ) \leq \ze \de \ri \}\\
&  & \qqu \qqu  \qqu \qqu  \qqu  \qqu \qqu \geq \Pr \li \{  n_s
\mscr{M}_{\mrm{B}} \li ( \wh{\bs{p}}_s, \mscr{U} ( \wh{\bs{p}}_s )
\ri ) \leq \ln (\ze \de) \ri \} = 1 \eee} and thus $\Pr \{
F_{\wh{\bs{p}}_s} ( \wh{\bs{p}}_s, \mscr{U} ( \wh{\bs{p}}_s ) ) \leq
\ze \de_s, \; G_{\wh{\bs{p}}_s} ( \wh{\bs{p}}_s,  \mscr{L} (
\wh{\bs{p}}_s ) ) \leq \ze \de_s \} = 1$, which implies that $\{
\bs{D}_s = 1 \}$ is a sure event.  So, the sampling scheme satisfies
all the requirements described in Theorem \ref{Monotone_second},
from which Theorem \ref{Bino_mix_CDF_CH} immediately follows.

\subsection{Proof of Theorem \ref{Bino_range_mix} }  \la{App_Bino_range_mix}

We need some preliminary results.

\beL \la{lem27} $\li \{ \mscr{M}_{\mrm{B}} (\wh{\bs{p}}_\ell,
\wh{\bs{p}}_\ell - \vep_a )  > \f{ \ln ( \ze \de  ) } { n_\ell }, \;
\wh{\bs{p}}_\ell \leq p^\star + \vep_a \ri \} = \{ z_a^- <
\wh{\bs{p}}_\ell \leq p^\star + \vep_a \}$. \eeL

\bpf

By the definition of sample sizes, we have {\small $n_s = \li \lc
\f{ \ln (\ze \de) } { \mscr{M}_{\mrm{B}} (p^\star + \vep_a, p^\star)
} \ri \rc$} and thus {\small $n_\ell \leq n_s - 1 < \f{ \ln (\ze
\de) } { \mscr{M}_{\mrm{B}} (p^\star + \vep_a, p^\star) } = \f{ \ln
(\ze \de) } { \mscr{M}_{\mrm{B}}(z^\star, z^\star - \vep_a) }$}
where $z^\star = p^\star + \vep_a$. Since
$\mscr{M}_{\mrm{B}}(z^\star, z^\star - \vep_a)$ is negative, we have
{\small $\mscr{M}_{\mrm{B}}(z^\star, z^\star - \vep_a)
> \f{ \ln (\ze \de) } {n_\ell}$}.  Noting that {\small $\lim_{z \to \vep_a} \mscr{M}_{\mrm{B}}(z, z -
\vep_a) = - \iy < \f{ \ln (\ze \de) } {n_\ell}$} and that
$\mscr{M}_{\mrm{B}}(z, z - \vep_a)$ is monotonically increasing with
respect to $z \in (\vep_a, z^\star)$ as asserted by Lemma
\ref{CCon}, we can conclude from the intermediate value theorem that
there exists a unique number $z_a^- \in (\vep_a, p^\star + \vep_a)$
such that {\small $\mscr{M}_{\mrm{B}}(z_a^-, z_a^- + \vep_a) = \f{
\ln (\ze \de) } {n_\ell}$}.  Finally, by virtue of the monotonicity
of $\mscr{M}_{\mrm{B}}(z, z - \vep_a)$ with respect to $z \in
(\vep_a, z^\star)$, the lemma is established. \epf

\beL \la{lem28} $\li \{ \mscr{M}_{\mrm{B}} \li (\wh{\bs{p}}_\ell,
\f{\wh{\bs{p}}_\ell}{1 + \vep_r}  \ri )  > \f{ \ln ( \ze \de  ) } {
n_\ell }, \; \wh{\bs{p}}_\ell > p^\star + \vep_a \ri \} = \{ p^\star
+ \vep_a < \wh{\bs{p}}_\ell  < z_r^+ \}$. \eeL

\bpf Note that {\small $ \mscr{M}_{\mrm{B}}(z^\star, \f{z^\star}{1 +
\vep_r} ) = \mscr{M}_{\mrm{B}}(z^\star, z^\star - \vep_a)
> \f{ \ln (\ze \de) } {n_\ell}$}.  By the definition of sample
sizes, we have {\small $n_1 = \li \lc \f{ \ln (\ze \de) } { \ln (1
\sh (1 + \vep_r) ) } \ri \rc$} and thus {\small $n_\ell \geq n_1
\geq \f{ \ln (\ze \de) } { \ln (1 \sh (1 + \vep_r) ) } = \f{ \ln
(\ze \de) } { \mscr{M}_{\mrm{B}}(1, 1 \sh (1 + \vep_r) ) } = \f{ \ln
(\ze \de) } { \lim_{z \to 1} \mscr{M}_{\mrm{B}}(z, z \sh (1 +
\vep_r) ) }$}, which implies {\small $\lim_{z \to 1}
\mscr{M}_{\mrm{B}} (z, \f{z}{1 + \vep_r} ) \leq \f{ \ln (\ze \de) }
{n_\ell}$}. Noting that $\mscr{M}_{\mrm{B}}(z, \f{z} {1 + \vep_r} )$
is monotonically decreasing with respect to $z \in (z^\star, 1)$, we
can conclude from the intermediate value theorem that there exists a
unique number $z_r^+ \in (z^\star, 1]$ such that {\small
$\mscr{M}_{\mrm{B}}(z_r^+, \f{z_r^+}{1 + \vep_r} ) = \f{ \ln (\ze
\de) } {n_\ell}$}.  Finally, by virtue of the monotonicity of
$\mscr{M}_{\mrm{B}}(z, \f{z}{1 + \vep_r} )$ with respect to $z \in
(z^\star, 1]$, the lemma is established.

\epf

\beL \la{lem29} For $\ell = 1, \cd, s - 1$, {\small \[ \li \{
\mscr{M}_{\mrm{B}} (\wh{\bs{p}}_\ell, \wh{\bs{p}}_\ell + \vep_a )  >
\f{\ln ( \ze \de )}{n_\ell}, \; \wh{\bs{p}}_\ell \leq p^\star -
\vep_a \ri \} = \bec \{ 0 \leq \wh{\bs{p}}_\ell \leq p^\star -
\vep_a \} & \tx{for} \; n_\ell < \f{ \ln (\ze \de) } { \ln (1 - \vep_a) }, \\
\{ z_a^+ < \wh{\bs{p}}_\ell \leq p^\star - \vep_a \} & \tx{for} \;
\f{ \ln (\ze \de) } { \ln (1 - \vep_a) } \leq n_\ell < \f{ \ln (\ze
\de) }
{ \mscr{M}_{\mrm{B}} (p^\star - \vep_a, p^\star) }, \\
\emptyset & \tx{for} \; n_\ell \geq \f{ \ln (\ze \de) } {
\mscr{M}_{\mrm{B}} (p^\star - \vep_a, p^\star) }. \eec \]} \eeL

\bpf

In the case of {\small $n_\ell < \f{ \ln (\ze \de) } { \ln (1 -
\vep_a) }$}, it is obvious that {\small $\ln (1 - \vep_a)
> \f{ \ln (\ze \de) } {n_\ell}$}. Since {\small $\lim_{z \to
0} \mscr{M}_{\mrm{B}} (z, z + \vep_a) = \ln (1 - \vep_a) < 0$}, we
have $\lim_{z \to 0} \mscr{M}_{\mrm{B}}  (z, z + \vep_a) > \f{ \ln
(\ze \de) } {n_\ell}$.  Observing that $\mscr{M}_{\mrm{B}} (z, z +
\vep_a)$ is monotonically increasing with respect to $z \in (0,
p^\star - \vep_a)$, we have {\small $\mscr{M}_{\mrm{B}}  (z, z +
\vep_a) > \f{ \ln (\ze \de) } {n_\ell}$} for any $z \in [0, p^\star
- \vep_a]$. It follows that {\small $\li \{ \mscr{M}_{\mrm{B}}
(\wh{\bs{p}}_\ell, \wh{\bs{p}}_\ell + \vep_a )  > \f{\ln \li ( \ze
\de \ri )}{n_\ell}, \; \wh{\bs{p}}_\ell \leq p^\star - \vep_a \ri \}
= \{ 0 \leq \wh{\bs{p}}_\ell \leq p^\star - \vep_a \}$ }.

In the case of {\small $\f{ \ln (\ze \de) } { \ln (1 - \vep_a) }
\leq n_\ell < \f{ \ln (\ze \de) } { \mscr{M}_{\mrm{B}} (p^\star -
\vep_a, p^\star) }$}, we have {\small $n_\ell < \f{ \ln (\ze \de) }
{ \mscr{M}_{\mrm{B}} (p^\star - \vep_a, p^\star) } = \f{ \ln (\ze
\de) } { \mscr{M}_{\mrm{B}}(z^*, z^* + \vep_a) }$} where $z^* =
p^\star - \vep_a$. Observing that $\mscr{M}_{\mrm{B}}(z^*, z^* +
\vep_a)$ is negative, we have {\small $\mscr{M}_{\mrm{B}}(z^*, z^* +
\vep_a)
> \f{ \ln (\ze \de) } {n_\ell}$}.  On the other hand,  {\small $\lim_{z \to 0} \mscr{M}_{\mrm{B}}(z, z
+ \vep_a) \leq \f{ \ln (\ze \de) } {n_\ell}$} as a consequence of
{$n_\ell \geq  \f{ \ln (\ze \de) } { \ln (1 - \vep_a) } = \f{ \ln
(\ze \de) } { \lim_{z \to 0} \mscr{M}_{\mrm{B}}(z, z + \vep_a) }$}.
Since $\mscr{M}_{\mrm{B}}(z, z + \vep_a)$ is monotonically
increasing with respect to $z \in (0, z^*) \subset (0, \f{1}{2} -
\vep_a)$, we can conclude from the intermediate value theorem that
there exists a unique number $z_a^+ \in [0, p^\star - \vep_a)$ such
that {\small $\mscr{M}_{\mrm{B}}(z_a^+, z_a^+ + \vep_a) = \f{ \ln
(\ze \de) } {n_\ell}$}.  By virtue of the monotonicity of
$\mscr{M}_{\mrm{B}}(z, z + \vep_a)$ with respect to $z \in (0,
z^*)$, we have {\small $\li \{ \mscr{M}_{\mrm{B}} (\wh{\bs{p}}_\ell,
\wh{\bs{p}}_\ell + \vep_a )
> \f{\ln \li ( \ze \de \ri )}{n_\ell}, \; \wh{\bs{p}}_\ell \leq
p^\star - \vep_a \ri \} = \{ z_a^+ < \wh{\bs{p}}_\ell \leq p^\star -
\vep_a \}$ }.

In the case of {\small $n_\ell \geq \f{ \ln (\ze \de) } {
\mscr{M}_{\mrm{B}} (p^\star - \vep_a, p^\star) }$}, we have {\small
$n_\ell \geq \f{ \ln (\ze \de) } { \mscr{M}_{\mrm{B}} (p^\star -
\vep_a, p^\star) } = \f{ \ln (\ze \de) } { \mscr{M}_{\mrm{B}}(z^*,
z^* + \vep_a) }$}. Due to the fact that $\mscr{M}_{\mrm{B}}(z^*, z^*
+ \vep_a)$ is negative, we have {\small $\mscr{M}_{\mrm{B}}(z^*, z^*
+ \vep_a) \leq \f{ \ln (\ze \de) } {n_\ell}$}.  Since
$\mscr{M}_{\mrm{B}}(z, z + \vep_a)$ is monotonically increasing with
respect to $z \in (0, z^*) \subset (0, \f{1}{2} - \vep_a)$, we have
that $\mscr{M}_{\mrm{B}}(z, z + \vep_a) \leq \f{ \ln (\ze \de) }
{n_\ell}$ for any $z \in [0, z^*]$. This implies that {\small $\li
\{ \mscr{M}_{\mrm{B}} (\wh{\bs{p}}_\ell, \wh{\bs{p}}_\ell + \vep_a )
> \f{\ln \li ( \ze \de \ri )}{n_\ell}, \; \wh{\bs{p}}_\ell
\leq p^\star - \vep_a \ri \} = \emptyset$}. This completes the proof
of the lemma.

 \epf

\beL \la{lem30} For $\ell = 1, \cd, s - 1$, {\small \[ \li \{
\mscr{M}_{\mrm{B}} \li (\wh{\bs{p}}_\ell, \f{\wh{\bs{p}}_\ell}{1 -
\vep_r}  \ri )  > \f{ \ln ( \ze \de  ) } { n_\ell }, \;
\wh{\bs{p}}_\ell
> p^\star - \vep_a \ri \} = \bec \{ p^\star - \vep_a <
\wh{\bs{p}}_\ell  < z_r^- \}  & \tx{for} \; n_\ell < \f{ \ln (\ze
\de) }
{ \mscr{M}_{\mrm{B}} (p^\star - \vep_a, p^\star) },\\
\emptyset & \tx{for} \; n_\ell \geq \f{ \ln (\ze \de) } {
\mscr{M}_{\mrm{B}} (p^\star - \vep_a, p^\star) }. \eec
\]} \eeL

\bpf

In the case of {\small $n_\ell < \f{ \ln (\ze \de) } {
\mscr{M}_{\mrm{B}} (p^\star - \vep_a, p^\star) }$}, we have  {\small
$ \mscr{M}_{\mrm{B}}(z^*, \f{z^*}{1 - \vep_r} ) =
\mscr{M}_{\mrm{B}}(z^*, z^* + \vep_a) = \mscr{M}_{\mrm{B}} (p^\star
- \vep_a, p^\star) > \f{ \ln (\ze \de) } {n_\ell}$}. Noting that
{\small $\lim_{z \to 1 - \vep_r} \mscr{M}_{\mrm{B}} (z, \f{z}{1 -
\vep_r} ) = - \iy < \f{ \ln (\ze \de) } {n_\ell}$} and that
$\mscr{M}_{\mrm{B}}(z, \f{z} {1 - \vep_r} )$ is monotonically
decreasing with respect to $z \in (z^*, 1 - \vep_r)$, we can
conclude from the intermediate value theorem that there exists a
unique number $z_r^- \in (z^*, 1 - \vep_r)$ such that {\small
$\mscr{M}_{\mrm{B}}(z_r^-, \f{z_r^-}{1 - \vep_r} ) = \f{ \ln (\ze
\de) } {n_\ell}$}.  By virtue of the monotonicity of
$\mscr{M}_{\mrm{B}}(z, \f{z} {1 - \vep_r} )$ with respect to $z \in
(z^*, 1 - \vep_r)$, we have {\small $ \{ \mscr{M}_{\mrm{B}}
(\wh{\bs{p}}_\ell, \f{\wh{\bs{p}}_\ell}{1 - \vep_r} )  > \f{\ln (\ze
\de )}{n_\ell}, \; \wh{\bs{p}}_\ell
> p^\star - \vep_a  \} = \{ p^\star - \vep_a <
\wh{\bs{p}}_\ell  < z_r^- \}$}.

In the case of {\small $n_\ell \geq \f{ \ln (\ze \de) } {
\mscr{M}_{\mrm{B}} (p^\star - \vep_a, p^\star) }$}, we have  {\small
$ \mscr{M}_{\mrm{B}}(z^*, \f{z^*}{1 - \vep_r} ) \leq \f{ \ln (\ze
\de) } {n_\ell}$}. Noting that $\mscr{M}_{\mrm{B}}(z, \f{z}{1 -
\vep_r} )$ is monotonically decreasing with respect to $z \in (z^*,
1 - \vep_r)$, we can conclude that {\small $\mscr{M}_{\mrm{B}}(z,
\f{z}{1 - \vep_r} ) \leq \f{ \ln (\ze \de) } {n_\ell}$} for any $z
\in [z^*, 1 - \vep_r)$.  This implies that {\small $ \{
\mscr{M}_{\mrm{B}} (\wh{\bs{p}}_\ell, \f{\wh{\bs{p}}_\ell}{1 -
\vep_r} )  > \f{\ln (\ze \de )}{n_\ell}, \; \wh{\bs{p}}_\ell
> p^\star - \vep_a  \} = \emptyset$}.  The proof of the lemma is
thus completed.

\epf

\bsk

We are now in position to prove Theorem \ref{Bino_range_mix}.
Clearly, it follows directly from the definition of $\bs{D}_\ell$
that {\small $\{ \bs{D}_\ell = 0 \}
 = \{ \mscr{M}_{\mrm{B}} (\wh{\bs{p}}_\ell, \mscr{L} ( \wh{\bs{p}}_\ell ) )  >
\f{ \ln ( \ze \de  ) } { n_\ell } \} \cup \{ \mscr{M}_{\mrm{B}}
(\wh{\bs{p}}_\ell, \mscr{U} ( \wh{\bs{p}}_\ell ) )  > \f{ \ln ( \ze
\de  ) } { n_\ell } \}$}.  It remains to show statements (I) and
(II).

With regard to statement (I), invoking the definition of $\mscr{L} (
\wh{\bs{p}}_\ell )$, we have {\small \bee \li \{ \mscr{M}_{\mrm{B}}
(\wh{\bs{p}}_\ell, \mscr{L} ( \wh{\bs{p}}_\ell ) )
> \f{ \ln ( \ze \de  ) } { n_\ell } \ri \} & = & \li \{
\mscr{M}_{\mrm{B}} (\wh{\bs{p}}_\ell, \wh{\bs{p}}_\ell - \vep_a )  >
\f{ \ln ( \ze \de  ) } { n_\ell }, \;  \wh{\bs{p}}_\ell \leq
p^\star + \vep_a \ri \}\\
&   &  \bigcup \li \{ \mscr{M}_{\mrm{B}} \li (\wh{\bs{p}}_\ell,
\f{\wh{\bs{p}}_\ell}{1 + \vep_r}  \ri )  > \f{ \ln ( \ze \de  ) } {
n_\ell }, \; \wh{\bs{p}}_\ell > p^\star + \vep_a \ri
\}\\
& = &  \{ z_a^- < \wh{\bs{p}}_\ell \leq p^\star + \vep_a \} \cup \{
p^\star + \vep_a < \wh{\bs{p}}_\ell  < z_r^+\}\\
& = &  \{ z_a^- < \wh{\bs{p}}_\ell  < z_r^+\} = \{ n_\ell \; z_a^- <
K_\ell  < n_\ell \; z_r^+\} \eee } where the second equality is due
to Lemma \ref{lem27} and Lemma \ref{lem28}.  This establishes
statement (I).

The proof of statement (II) can be completed by applying Lemma
\ref{lem29}, Lemma \ref{lem30} and observing that {\small \bee \li
\{ \mscr{M}_{\mrm{B}} (\wh{\bs{p}}_\ell, \mscr{U} ( \wh{\bs{p}}_\ell
) )
> \f{ \ln ( \ze \de  ) } { n_\ell } \ri \} & = & \li \{
\mscr{M}_{\mrm{B}} (\wh{\bs{p}}_\ell, \wh{\bs{p}}_\ell + \vep_a )  >
\f{ \ln ( \ze \de  ) } { n_\ell }, \;  \wh{\bs{p}}_\ell
\leq p^\star - \vep_a \ri \}\\
&   &  \bigcup \li \{ \mscr{M}_{\mrm{B}} \li (\wh{\bs{p}}_\ell,
\f{\wh{\bs{p}}_\ell}{1 - \vep_r}  \ri )  > \f{ \ln ( \ze \de  ) } {
n_\ell }, \; \wh{\bs{p}}_\ell > p^\star - \vep_a \ri \}. \eee} This
completes the proof of Theorem \ref{Bino_range_mix}.

\subsection{Proof of Theorem \ref{Bino_mix_Massart} }  \la{App_Bino_mix_Massart}

 We need some preliminary results, especially some properties
 of function $\mscr{M}(z, \mu)$.

\beL \la{decrea} $\mscr{M} (z, z + \vep )$ is monotonically
increasing with respect to  $z \in (0, \f{1}{2} - \f{2 \vep}{3})$,
and is monotonically decreasing with respect to $z \in (\f{1}{2} -
\f{2 \vep}{3}, 1 - \vep)$.  Similarly, $\mscr{M} (z, z - \vep )$ is
monotonically increasing with respect to $z \in (\vep, \f{1}{2} +
\f{2 \vep}{3})$, and is monotonically decreasing with respect to $z
\in ( \f{1}{2} + \f{2 \vep}{3}, 1)$. \eeL

\bpf The lemma can be established by checking the partial
derivatives
\[
\f{\pa \mscr{M} (z, z + \vep )} { \pa z } = \f{ \vep^2 } { \li [ \li
(z + \f{2 \vep}{3} \ri ) \li ( 1- z - \f{2 \vep}{3} \ri ) \ri ]^2 }
\li ( \f{1}{2} - \f{2 \vep}{3} - z \ri ),
\]
\[
\f{\pa \mscr{M} (z, z - \vep )} { \pa z } = \f{ \vep^2 } { \li [ \li
(z - \f{2 \vep}{3} \ri ) \li ( 1- z + \f{2 \vep}{3} \ri ) \ri ]^2 }
\li ( \f{1}{2} + \f{2 \vep}{3} - z \ri ).
\]

\epf

\beL \la{lem888m}  Let $0 < \vep < \f{1}{2}$. Then, $\mscr{M}(z, z -
\vep) \leq \mscr{M}(z, z + \vep) \leq - 2 \vep^2$ for $z \in \li [0,
\f{1}{2} \ri ]$, and $\mscr{M}(z, z + \vep) < \mscr{M}(z, z - \vep)
\leq - 2 \vep^2$ for $z \in \li ( \f{1}{2}, 1 \ri ]$. \eeL

\bpf

By the definition of the function $\mscr{M} (.,.)$, we have that
$\mscr{M} (z, \mu) = - \iy$ for $z \in [0, 1]$ and $\mu \notin (0,
1)$. Hence, the lemma is trivially true for $0 \leq z \leq \vep$ or
$1 - \vep \leq z \leq 1$.  It remains to show the lemma for $z \in
(\vep, 1 - \vep)$.  This can be accomplished by noting that
\[
\mscr{M} (z, z + \vep ) - \mscr{M} (z, z - \vep ) = \f{ 2 \vep^3 ( 1
- 2 z) } {3 \li (z + \f{2 \vep}{3} \ri ) \li ( 1- z - \f{2 \vep}{3}
\ri )  \li (z - \f{2 \vep}{3} \ri ) \li ( 1- z + \f{2 \vep}{3} \ri )
}.
\]
where the right-hand side is seen to be positive for $z \in \li (
\vep, \f{1}{2} \ri )$ and negative for $z \in \li ( \f{1}{2}, 1 -
\vep \ri )$.  By Lemma \ref{decrea}, the maximums of $\mscr{M}(z, z
+ \vep)$ and $\mscr{M}(z, z - \vep)$ are shown to be $- 2 \vep^2$.
This completes the proof of the lemma.

\epf

\beL \la{lem888m2B}

$\mscr{M} ( z, \f{z}{1 - \vep} ) < \mscr{M} ( z, \f{z}{1 + \vep} ) <
0$ for $0 < z < 1 - \vep < 1$.

\eeL

\bpf

It can be verified that
\[
\mscr{M} \li ( z, \f{z}{1 + \vep} \ri ) - \mscr{M} \li ( z, \f{z}{1
- \vep} \ri ) = \f{ 2 \vep^3 z (2 - z) } {3 \li ( 1 + \f{\vep}{3}
\ri ) \li [ 1 - z + \vep \li ( 1 - \f{ z}{3} \ri ) \ri ]  \li ( 1 -
\f{\vep}{3} \ri ) \li [ 1 - z - \vep \li ( 1 - \f{ z}{3} \ri ) \ri
]},
\]
from which it can be seen that {\small $\mscr{M} ( z, \f{z}{1 -
\vep} ) < \mscr{M} ( z, \f{z}{1 + \vep} ) < 0$} for $z \in (0, 1 -
\vep)$.

\epf

\beL \la{lemm22m}

$\mscr{M} (\mu - \vep, \mu) < \mscr{M} (\mu + \vep, \mu) \leq - 2
\vep^2$ for $0 < \vep < \mu < \f{1}{2} < 1 - \vep$.

\eeL

\bpf The lemma follows from Lemma \ref{lem888m} and the fact that
\[
\mscr{M} (\mu - \vep, \mu) - \mscr{M} (\mu + \vep, \mu) = \f{ \vep^3
(2 \mu - 1) } { 3 \li ( \mu - \f{\vep}{3} \ri ) \li ( 1 - \mu +
\f{\vep}{3} \ri )  \li ( \mu + \f{\vep}{3} \ri ) \li ( 1 - \mu -
\f{\vep}{3} \ri )  },
\]
where the right-hand side is negative for $0 < \vep < \mu < \f{1}{2}
< 1 - \vep$.

\epf

\beL \la{decrev} {\small $\mscr{M} ( z, \f{z}{1 + \vep} )$} is
monotonically decreasing with respect to $z \in (0, 1)$. Similarly,
{\small $\mscr{M} ( z, \f{z}{1 - \vep} )$} is monotonically
decreasing with respect to $z \in (0, 1 - \vep)$. \eeL

\bpf

The lemma can be shown by verifying that
\[
\f{\pa } { \pa z} \mscr{M} \li ( z, \f{z}{1 + \vep} \ri ) = - \f{
\vep^2}{2 \li ( 1 + \f{\vep}{3} \ri )} \times \f{ 1 + \vep } { \li [
(1 + \vep) ( 1 - z) + \f{2 \vep z}{3} \ri ]^2 } < 0
\]
for $z \in (0, 1)$ and that
\[
\f{\pa } { \pa z} \mscr{M} \li ( z, \f{z}{1 - \vep} \ri ) = - \f{
\vep^2}{2 \li ( 1 - \f{\vep}{3} \ri ) } \times \f{ 1 - \vep } { \li
[ (1 - \vep) ( 1 - z) - \f{2 \vep z}{3} \ri ]^2 } < 0
\]
for $z \in (0, 1 - \vep)$.

\epf

\beL

\la{Ddec} For any fixed $z \in (0, 1)$, $\mscr{M}(z,\mu)$ is
monotonically increasing with respect to $\mu \in (0, z)$, and is
monotonically decreasing with respect to $\mu \in (z, 1)$.
Similarly, for any fixed $\mu \in (0, 1)$, $\mscr{M}(z,\mu)$ is
monotonically increasing with respect to $z \in (0, \mu)$, and is
monotonically decreasing with respect to $z \in (\mu, 1)$.

\eeL

\bpf

The lemma can be shown by checking the following partial
derivatives: \bee \f{\pa \mscr{M}(z,\mu)}{\pa \mu} & = & \f{(z -
\mu) \li [ \mu(1 - z) + z ( 1 - \mu) + z ( 1 -z) \ri ] } { 3 \li [
\li ( \f{2\mu}{3} + \f{z}{3} \ri ) \li ( 1 - \f{2\mu}{3} - \f{z}{3}
\ri ) \ri ]^2
},\\
\f{\pa \mscr{M}(z,\mu)}{\pa z} & = & \f{(\mu - z) \li [ \mu(1 -
\f{2\mu}{3} - \f{z}{3} ) + \f{z - \mu}{6} \ri ] } {
 \li [ \li ( \f{2\mu}{3} + \f{z}{3} \ri ) \li ( 1 - \f{2\mu}{3} - \f{z}{3} \ri ) \ri ]^2 }
  = \f{(\mu - z) \li [ (1 - \mu) (\f{2\mu}{3} + \f{z}{3} ) + \f{\mu - z}{6} \ri ] } {
 \li [ \li ( \f{2\mu}{3} + \f{z}{3} \ri ) \li ( 1 - \f{2\mu}{3} - \f{z}{3} \ri ) \ri ]^2 }.
\eee

\epf

\beL

\la{good8a}

{\small $\{  \bs{D}_\ell = 1 \} \subseteq \{  \mscr{M}_{\mrm{B}} (
\wh{\bs{p}}_\ell, \mscr{U} ( \wh{\bs{p}}_\ell )  ) \leq \f{\ln (\ze
\de)}{n_\ell}, \; \mscr{M}_{\mrm{B}} ( \wh{\bs{p}}_\ell, \mscr{L} (
\wh{\bs{p}}_\ell )  ) \leq \f{\ln (\ze \de)}{n_\ell} \}$} for $\ell
= 1, \cd, s$.

\eeL

\bpf

By the definition of $n_s$, we can show that $n_s \leq \f{ \ln
\f{1}{\ze \de} }{ 2\vep_a^2 }$, which implies that $\f{1}{4} + \f{
n_\ell \vep_a^2 } {2 \ln (\ze \de) } \geq 0$ for $\ell = 1, \cd, s$.
It can be shown by tedious computation that {\small \bel &   & \li
\{ \mscr{M} \li ( \wh{\bs{p}}_\ell, \wh{\bs{p}}_\ell + \vep_a \ri )
> \f{\ln (\ze \de)}{n_\ell} \ri \} = \li \{  \f{1}{2} - \f{2}{3}
\vep_a - \sq{ \f{1}{4} + \f{ n_\ell \vep_a^2 } {2 \ln (\ze \de) } }
< \wh{\bs{p}}_\ell < \f{1}{2} - \f{2}{3} \vep_a + \sq{ \f{1}{4} +
\f{ n_\ell \vep_a^2 } {2 \ln (\ze \de) } } \ri \}, \qqu \la{eqabsa} \qqu\\
&   & \li \{ \mscr{M} \li ( \wh{\bs{p}}_\ell,  \wh{\bs{p}}_\ell -
\vep_a \ri ) > \f{\ln (\ze \de)}{n_\ell} \ri \} = \li \{  \f{1}{2} +
\f{2}{3} \vep_a - \sq{ \f{1}{4} + \f{ n_\ell \vep_a^2 } {2 \ln (\ze
\de) } } < \wh{\bs{p}}_\ell < \f{1}{2} + \f{2}{3} \vep_a + \sq{
\f{1}{4} + \f{ n_\ell \vep_a^2 } {2 \ln (\ze \de) } } \ri \}, \qqu  \la{eqabsb}\\
&   & \li \{ \mscr{M} \li ( \wh{\bs{p}}_\ell, \f{\wh{\bs{p}}_\ell}{1
+ \vep_r} \ri ) > \f{ \ln (\ze \de) }{n_\ell} \ri \} =  \li \{
\wh{\bs{p}}_\ell < \f{ 6(1 + \vep_r) (3 + \vep_r) \ln (\ze \de) } {
2 (3 + \vep_r)^2 \ln (\ze \de) - 9 n_\ell \vep_r^2 } \ri \}, \la{eqreva}\\
&   & \li \{ \mscr{M} \li ( \wh{\bs{p}}_\ell, \f{\wh{\bs{p}}_\ell}{1
- \vep_r} \ri ) > \f{ \ln (\ze \de) }{n_\ell} \ri \} =  \li \{
\wh{\bs{p}}_\ell < \f{ 6(1 - \vep_r) (3 - \vep_r) \ln (\ze \de) } {
2 (3 - \vep_r)^2 \ln (\ze \de) -  9 n_\ell \vep_r^2} \ri \}
\la{eqrevb} \eel} for $\ell = 1, \cd, s$.  By (\ref{eqrevb}), we
have {\small \be \la{part1a} \li \{ \mscr{M} \li ( \wh{\bs{p}}_\ell,
\f{ \wh{\bs{p}}_\ell }{ 1 - \vep_r} \ri ) > \f{\ln (\ze
\de)}{n_\ell}, \; \wh{\bs{p}}_\ell > \f{\vep_a}{\vep_r} - \vep_a \ri
\} = \li \{ \f{\vep_a}{\vep_r} - \vep_a < \wh{\bs{p}}_\ell < \f{ 6(1
- \vep_r) (3 - \vep_r) \ln (\ze \de) } { 2 (3 - \vep_r)^2 \ln (\ze
\de) -  9 n_\ell \vep_r^2}  \ri \}. \ee} By the assumption that $0 <
\vep_a < \f{3}{8}$ and $\f{ 6 \vep_a } { 3 - 2 \vep_a } < \vep_r <
1$, we have $\f{\vep_a}{\vep_r} - \vep_a < \f{1}{2} - \f{4
\vep_a}{3}$. Hence, by virtue of (\ref{eqabsa}), we have {\small \be
\la{part2b} \li \{ \mscr{M} ( \wh{\bs{p}}_\ell, \wh{\bs{p}}_\ell +
\vep_a ) > \f{\ln (\ze \de)}{n_\ell}, \; \wh{\bs{p}}_\ell \leq
\f{\vep_a}{\vep_r} - \vep_a \ri \} = \li \{ \f{1}{2} - \f{2}{3}
\vep_a - \sq{ \f{1}{4} + \f{ n_\ell \vep_a^2 } {2 \ln (\ze \de) } }
< \wh{\bs{p}}_\ell  \leq \f{\vep_a}{\vep_r} - \vep_a \ri \}.  \ee}
Therefore, making use of (\ref{part1a}) and (\ref{part2b}), we have
{\small \bel \li \{ \mscr{M} ( \wh{\bs{p}}_\ell, \mscr{U} (
\wh{\bs{p}}_\ell ) )
> \f{\ln (\ze \de)}{n_\ell} \ri  \} & = & \li \{ \mscr{M} \li (
\wh{\bs{p}}_\ell, \f{ \wh{\bs{p}}_\ell }{ 1 - \vep_r} \ri ) > \f{\ln
(\ze \de)}{n_\ell},
\; \wh{\bs{p}}_\ell > \f{\vep_a}{\vep_r} - \vep_a \ri \} \la{latera}\\
&  & \cup \li \{ \mscr{M} ( \wh{\bs{p}}_\ell, \wh{\bs{p}}_\ell +
\vep_a ) > \f{\ln (\ze \de)}{n_\ell}, \; \wh{\bs{p}}_\ell \leq
\f{\vep_a}{\vep_r} - \vep_a \ri \} \nonumber\\
& = &  \li \{  \f{1}{2} - \f{2}{3} \vep_a - \sq{ \f{1}{4} + \f{
n_\ell \vep_a^2 } {2 \ln (\ze \de) } } < \wh{\bs{p}}_\ell < \f{ 6(1
- \vep_r) (3 - \vep_r) \ln (\ze \de) } { 2 (3 - \vep_r)^2 \ln (\ze
\de) -  9 n_\ell \vep_r^2}  \ri \}.   \la{imp1} \qqu \eel} By
(\ref{eqreva}), we have {\small \be \la{part1a88} \li \{ \mscr{M}
\li ( \wh{\bs{p}}_\ell, \f{ \wh{\bs{p}}_\ell }{ 1 + \vep_r} \ri ) >
\f{\ln (\ze \de)}{n_\ell}, \; \wh{\bs{p}}_\ell > \f{\vep_a}{\vep_r}
+ \vep_a \ri \} = \li \{ \f{\vep_a}{\vep_r} + \vep_a <
\wh{\bs{p}}_\ell < \f{ 6(1 + \vep_r) (3 + \vep_r) \ln (\ze \de) } {
2 (3 + \vep_r)^2 \ln (\ze \de) -  9 n_\ell \vep_r^2}  \ri \}. \ee}
By the assumption that $0 < \vep_a < \f{3}{8}$ and $\f{ 6 \vep_a } {
3 - 2 \vep_a } < \vep_r < 1$, we have $\f{\vep_a}{\vep_r} + \vep_a <
\f{1}{2} + \f{2 \vep_a}{3}$. Hence, by virtue of (\ref{eqabsb}), we
have {\small \be \la{part2b88} \li \{ \mscr{M} ( \wh{\bs{p}}_\ell,
\wh{\bs{p}}_\ell - \vep_a ) > \f{\ln (\ze \de)}{n_\ell}, \;
\wh{\bs{p}}_\ell \leq \f{\vep_a}{\vep_r} + \vep_a \ri \} = \li \{
\f{1}{2} + \f{2}{3} \vep_a - \sq{ \f{1}{4} + \f{ n_\ell \vep_a^2 }
{2 \ln (\ze \de) } } < \wh{\bs{p}}_\ell  \leq \f{\vep_a}{\vep_r} +
\vep_a \ri \}.  \ee} Therefore, making use of (\ref{part1a88}) and
(\ref{part2b88}), we have {\small \bel \li \{  \mscr{M} (
\wh{\bs{p}}_\ell, \mscr{L} ( \wh{\bs{p}}_\ell ) ) > \f{\ln (\ze
\de)}{n_\ell} \ri  \} & = & \li \{ \mscr{M} \li ( \wh{\bs{p}}_\ell,
\f{ \wh{\bs{p}}_\ell }{ 1 + \vep_r} \ri ) > \f{\ln (\ze
\de)}{n_\ell},
\; \wh{\bs{p}}_\ell > \f{\vep_a}{\vep_r} + \vep_a \ri \} \la{laterb}\\
&  & \cup \li \{ \mscr{M} ( \wh{\bs{p}}_\ell, \wh{\bs{p}}_\ell -
\vep_a ) > \f{\ln (\ze \de)}{n_\ell}, \; \wh{\bs{p}}_\ell \leq
\f{\vep_a}{\vep_r} + \vep_a \ri \} \nonumber\\
& = &  \li \{  \f{1}{2} + \f{2}{3} \vep_a - \sq{ \f{1}{4} + \f{
n_\ell \vep_a^2 } {2 \ln (\ze \de) } } < \wh{\bs{p}}_\ell < \f{ 6(1
+ \vep_r) (3 + \vep_r) \ln (\ze \de) } { 2 (3 + \vep_r)^2 \ln (\ze
\de) -  9 n_\ell \vep_r^2}  \ri \}. \qqu  \la{imp2}  \eel} It
follows from (\ref{imp1}) and (\ref{imp2}) that \be \la{laterc} \{
\bs{D}_\ell = 0 \} = \li \{  \mscr{M} ( \wh{\bs{p}}_\ell, \mscr{U} (
\wh{\bs{p}}_\ell )  ) > \f{\ln (\ze \de)}{n_\ell} \ri  \} \cup \li
\{ \mscr{M} ( \wh{\bs{p}}_\ell, \mscr{L} ( \wh{\bs{p}}_\ell )  ) >
\f{\ln (\ze \de)}{n_\ell} \ri  \}, \ee which implies that $\{
\bs{D}_\ell = 1 \} =  \{  \mscr{M} ( \wh{\bs{p}}_\ell, \mscr{U} (
\wh{\bs{p}}_\ell )  ) \leq \f{\ln (\ze \de)}{n_\ell}, \; \mscr{M} (
\wh{\bs{p}}_\ell, \mscr{L} ( \wh{\bs{p}}_\ell )  ) \leq \f{\ln (\ze
\de)}{n_\ell}  \}$ for $\ell = 1, \cd, s$.  So, {\small \bee \{
\bs{D}_\ell = 1 \} & = & \li \{  \mscr{M} ( \wh{\bs{p}}_\ell,
\mscr{U} ( \wh{\bs{p}}_\ell )  ) \leq \f{\ln (\ze \de)}{n_\ell}, \;
\mscr{M} ( \wh{\bs{p}}_\ell, \mscr{L} ( \wh{\bs{p}}_\ell )  ) \leq
\f{\ln (\ze \de)}{n_\ell} \ri  \}\\
& \subseteq  & \li \{  \mscr{M}_{\mrm{B}} ( \wh{\bs{p}}_\ell,
\mscr{U} ( \wh{\bs{p}}_\ell )  ) \leq \f{\ln (\ze \de)}{n_\ell}, \;
\mscr{M}_{\mrm{B}} ( \wh{\bs{p}}_\ell, \mscr{L} ( \wh{\bs{p}}_\ell )
) \leq \f{\ln (\ze \de)}{n_\ell} \ri  \} \eee} for $\ell = 1, \cd,
s$.  This completes the proof of the lemma.

\epf

\beL

\la{bb3800}

$\bs{D}_s = 1$.

\eeL

\bpf

For simplicity of notations, we denote $p^\star =
\f{\vep_a}{\vep_r}$.  In view of (\ref{latera}), (\ref{laterb}) and
(\ref{laterc}), we have that, in order to show $\bs{D}_s = 1$, it
suffices to show \bel &   & \li \{ \mscr{M} \li ( \wh{\bs{p}}_s, \f{
\wh{\bs{p}}_s }{ 1 - \vep_r} \ri ) > \f{\ln (\ze \de)}{n_s}, \;
\wh{\bs{p}}_s > p^\star - \vep_a \ri \} = \emptyset, \la{ep81}\\
&   & \li \{ \mscr{M} ( \wh{\bs{p}}_s, \wh{\bs{p}}_s + \vep_a)
> \f{\ln (\ze \de)}{n_s}, \; \wh{\bs{p}}_s \leq
p^\star - \vep_a \ri \} = \emptyset, \la{ep82}\\
&   & \li \{ \mscr{M} \li ( \wh{\bs{p}}_s, \f{ \wh{\bs{p}}_s }{ 1 +
\vep_r} \ri ) > \f{\ln (\ze \de)}{n_s}, \;
\wh{\bs{p}}_s > p^\star + \vep_a \ri \} = \emptyset, \la{ep83}\\
&   & \li \{ \mscr{M} ( \wh{\bs{p}}_s, \wh{\bs{p}}_s - \vep_a)
> \f{\ln (\ze \de)}{n_s}, \; \wh{\bs{p}}_s \leq
p^\star + \vep_a \ri \} = \emptyset. \la{ep84} \eel By the
definition of $n_s$,  we have {\small $n_s \geq \li \lc  \f{ \ln
(\ze \de) } { \mscr{M} \li ( p^\star + \vep_a, p^\star \ri ) } \ri
\rc \geq \f{ \ln (\ze \de) } { \mscr{M} \li ( p^\star + \vep_a,
p^\star \ri ) } $}.  By the assumption on $\vep_a$ and $\vep_r$, we
have $0 < \vep_a < p^\star < \f{1}{2} < 1 - \vep_a$.  Hence, by
Lemma \ref{lemm22m}, we have $\mscr{M} \li ( p^\star - \vep_a,
p^\star \ri ) < \mscr{M} \li ( p^\star + \vep_a, p^\star \ri ) < 0$
and it follows that \be \la{rela8} \f{\ln (\ze \de)}{n_s} \geq
\mscr{M} \li ( p^\star + \vep_a, p^\star \ri ) >  \mscr{M} \li (
p^\star - \vep_a, p^\star \ri ). \ee By (\ref{rela8}), {\small \be
\la{ep989} \li \{ \mscr{M} \li ( \wh{\bs{p}}_s, \f{ \wh{\bs{p}}_s }{
1 - \vep_r} \ri )
> \f{\ln (\ze \de)}{n_s}, \; \wh{\bs{p}}_s > p^\star - \vep_a \ri \}
\subseteq \li \{ \mscr{M} \li ( \wh{\bs{p}}_s, \f{ \wh{\bs{p}}_s }{
1 - \vep_r} \ri ) > \mscr{M} \li ( p^\star - \vep_a, p^\star \ri ),
\; \wh{\bs{p}}_s > p^\star - \vep_a \ri \}. \ee} Noting that
$\mscr{M} \li ( p^\star - \vep_a, p^\star \ri ) = \mscr{M} \li (
p^\star - \vep_a, \f{p^\star - \vep_a}{1 - \vep_r} \ri )$ and making
use of the fact that $\mscr{M}(z, \f{z}{1 - \vep} )$ is
monotonically decreasing with respect to $z \in (0, 1 - \vep)$ as
asserted by Lemma \ref{decrev}, we have \be \la{ep99} \li \{
\mscr{M} \li ( \wh{\bs{p}}_s, \f{ \wh{\bs{p}}_s }{ 1 - \vep_r} \ri )
> \mscr{M} \li ( p^\star - \vep_a, p^\star \ri ) \ri \} = \{
\wh{\bs{p}}_s < p^\star - \vep_a \}. \ee Combining (\ref{ep989}) and
(\ref{ep99}) yields (\ref{ep81}). By (\ref{rela8}), {\small \be
\la{eq9998} \li \{ \mscr{M} ( \wh{\bs{p}}_s, \wh{\bs{p}}_s + \vep_a)
> \f{\ln (\ze \de)}{n_s}, \; \wh{\bs{p}}_s \leq
p^\star - \vep_a \ri \} \subseteq \li \{ \mscr{M} ( \wh{\bs{p}}_s,
\wh{\bs{p}}_s + \vep_a) > \mscr{M} \li ( p^\star - \vep_a, p^\star
\ri ), \; \wh{\bs{p}}_s \leq p^\star - \vep_a \ri \}. \ee} By the
assumption on $\vep_a$ and $\vep_r$, we have $p^\star - \vep_a <
\f{1}{2} - \f{2 \vep_a}{3}$.   Recalling the fact that $\mscr{M}(z,
z + \vep)$ is monotonically increasing with respect to $z \in (0,
\f{1}{2} - \f{2 \vep}{3})$ as asserted by Lemma \ref{decrea}, we
have that the event in the right-hand side of (\ref{eq9998}) is an
impossible event and consequently, (\ref{ep82}) is established. By
(\ref{rela8}), {\small \be \la{ep98998} \li \{ \mscr{M} \li (
\wh{\bs{p}}_s, \f{ \wh{\bs{p}}_s }{ 1 + \vep_r} \ri ) > \f{\ln (\ze
\de)}{n_s}, \; \wh{\bs{p}}_s > p^\star + \vep_a \ri \} = \li \{
\mscr{M} \li ( \wh{\bs{p}}_s, \f{ \wh{\bs{p}}_s }{ 1 + \vep_r} \ri )
> \mscr{M} \li ( p^\star + \vep_a, p^\star \ri ), \; \wh{\bs{p}}_s >
p^\star + \vep_a \ri \}.  \ee} Noting that $\mscr{M} \li ( p^\star +
\vep_a, p^\star \ri ) = \mscr{M} \li ( p^\star + \vep_a, \f{p^\star
+ \vep_a}{1 + \vep_r} \ri )$ and making use of the fact that
$\mscr{M}(z, \f{z}{1 + \vep} )$ is monotonically decreasing with
respect to $z \in (0, 1)$ as asserted by Lemma \ref{decrev}, we have
\be \la{ep9998} \li \{ \mscr{M} \li ( \wh{\bs{p}}_s, \f{
\wh{\bs{p}}_s }{ 1 + \vep_r} \ri )
> \mscr{M} \li ( p^\star + \vep_a, p^\star \ri ) \ri \} = \{
\wh{\bs{p}}_s < p^\star + \vep_a \}. \ee Combining (\ref{ep98998})
and (\ref{ep9998}) yields (\ref{ep83}).  By (\ref{rela8}), {\small
\be \la{eq999833} \li \{ \mscr{M} ( \wh{\bs{p}}_s, \wh{\bs{p}}_s -
\vep_a) > \f{\ln (\ze \de)}{n_s}, \; \wh{\bs{p}}_s \leq p^\star +
\vep_a \ri \} \subseteq \li \{ \mscr{M} ( \wh{\bs{p}}_s,
\wh{\bs{p}}_s - \vep_a) > \mscr{M} \li ( p^\star + \vep_a, p^\star
\ri ), \; \wh{\bs{p}}_s \leq p^\star + \vep_a \ri \}. \ee} By the
assumption on $\vep_a$ and $\vep_r$, we have $p^\star + \vep_a <
\f{1}{2} + \f{2 \vep_a}{3}$.   Recalling the fact that $\mscr{M}(z,
z - \vep)$ is monotonically increasing with respect to $z \in (0,
\f{1}{2} + \f{2 \vep}{3})$ as asserted by Lemma \ref{decrea}, we
have that the event in the right-hand side of (\ref{eq999833}) is an
impossible event and consequently, (\ref{ep84}) is established. This
completes the proof of the lemma.

\epf

\bsk

Now we are in a position to prove Theorem \ref{Bino_mix_Massart}.
Recall that $\exp( \mscr{M}_{\mrm{B}} (z, p) )$  is equal to
$\mcal{F} (z, p)$ and $\mcal{G} (z, p)$ respectively for the cases
of $z \leq p$ and $z \geq p$.  Moreover, $\wh{\bs{p}}_\ell$ is a ULE
of $p$ for $\ell = 1, \cd, s$.  Furthermore, $\{ \bs{D}_s = 1 \}$ is
a sure event as a result of Lemma \ref{bb3800}. So, the sampling
scheme satisfies all the requirements described in Corollary
\ref{Monotone_third}, from which Theorem \ref{Bino_mix_Massart}
immediately follows.

\subsection{Proof of Theorem \ref{Bino_mix_DDV_Asp} }  \la{App_Bino_mix_DDV_Asp}

We need some preliminary results.

\beL

\la{lem52m} If $\vep_a$ is sufficiently small, then the following
statements hold true.

(I): For $1 \leq \ell < s$, there exists a unique number $z_\ell \in
[0, p^\star - \vep_a)$ such that $n_\ell = \f{ \ln ( \ze \de ) } {
\mscr{M}_{\mrm{B}} ( z_\ell, \; z_\ell + \vep_a ) }$ for $n_\ell
\geq \f{ \ln (\ze \de) } { \ln ( 1 - \vep_a) } $.

(II): For $1 \leq \ell < s$, there exists a unique number $y_\ell
\in (p^\star + \vep_a, 1]$ such that $n_\ell = \f{ \ln ( \ze \de ) }
{ \mscr{M}_{\mrm{B}} ( y_\ell, \; \f{y_\ell }{1 + \vep_r} ) }$.

(III): $z_\ell$ is monotonically increasing with respect to $\ell$;
$y_\ell$ is monotonically decreasing with respect to $\ell$.

(IV): $\lim_{\vep_a \to 0} z_\ell = \f{ 1 - \sq{ 1 - 4 p^\star (1 -
p^\star) C_{s - \ell}} }{2}$ and $\lim_{\vep_a \to 0} y_\ell =
\f{1}{1 + \li ( \f{1}{p^\star} - 1  \ri ) C_{s - \ell}}$, where the
limits are taken under the constraint that $\f{\vep_a}{\vep_r}$ and
$s - \ell$ are fixed with respect to $\vep_a$.

(V): Let $\ell_\vep = s - j_p$.  For $p \in (p^\star, 1)$ such that
$C_{j_p} = r (p)$,
\[
\lim_{\vep_r \to 0} \f{ z_{\ell_\vep} - p}{\vep_r p} = 1 - \f{1}{3} \f{p - p^\star}{1 - p^\star}.
\]
For $p \in (0, p^\star)$ such that $C_{j_p} = r (p)$,
\[
\lim_{\vep_a \to 0} \f{ z_{\ell_\vep} - p}{\vep_a} = \f{ p ( 1 - p) (1 - 2 p^\star) }{ 3 p^\star (1 - p^\star) (1 - 2 p) } - \f{2}{3}.
\]

(VI):
\[
\{ \bs{D}_\ell = 0 \} = \bec  \{ z_\ell < \wh{\bs{p}}_\ell < y_\ell
\} & \tx{for} \;  n_\ell \geq \f{ \ln (\ze
\de) } { \ln ( 1 - \vep_a) };\\
\{ 0 < \wh{\bs{p}}_\ell < y_\ell \} & \tx{for} \; n_\ell < \f{ \ln
(\ze \de) } { \ln ( 1 - \vep_a) }. \eec
\]
\eeL

\bsk

{\bf Proof of Statement (I)}: By the definition of sample sizes, we
have {\small $\f{ \ln ( \ze \de ) } { n_{\ell} } \geq
\mscr{M}_{\mrm{B}} ( 0, \vep_a )$} and \be \la{from2pp}
 n_{\ell} < \f{(1 + C_1) n_s}{2} < \f{(1 + C_1)
}{2} \li [ \f{ \ln ( \ze \de ) } { \mscr{M}_{\mrm{B}} ( p^\star +
\vep_a, \; p^\star  ) } + 1 \ri ]
 \ee
for sufficiently small $\vep_a > 0$. As a consequence of
(\ref{from2pp}), we have {\small \[ \f{ \ln ( \ze \de ) } { n_{\ell}
} < \mscr{M}_{\mrm{B}} ( p^\star + \vep_a, \; p^\star  )  \li (
\f{2}{1 + C_1} - \f{1}{ n_{\ell} } \ri ) = \f{\mscr{M}_{\mrm{B}} (
p^\star + \vep_a, \; p^\star  ) } { \mscr{M}_{\mrm{B}} ( p^\star -
\vep_a, \; p^\star ) } \li ( \f{2}{1 + C_1}  \ri )
\mscr{M}_{\mrm{B}} ( p^\star - \vep_a, \; p^\star ) -
\f{\mscr{M}_{\mrm{B}} ( p^\star + \vep_a, \; p^\star  ) }{ n_{\ell}
}
\]}
provided that $\vep_a > 0$ is sufficiently small.  Noting that
\[
\lim_{\vep_a \to 0} \f{\mscr{M}_{\mrm{B}} ( p^\star + \vep_a, \; p^\star  ) } { \mscr{M}_{\mrm{B}} ( p^\star - \vep_a, \; p^\star ) } =
1, \qqu \lim_{\vep_a \to 0} \f{\mscr{M}_{\mrm{B}} ( p^\star + \vep_a, \; p^\star  ) }{ n_{\ell} } = 0,
\]
we have that $\f{ \ln ( \ze \de ) } { n_{\ell} } <
\mscr{M}_{\mrm{B}} ( p^\star - \vep_a, \; p^\star )$ for small
enough $\vep_a > 0$.   In view of the established fact that $ \mscr{M}_{\mrm{B}} ( 0,
\vep_a ) \leq \f{ \ln ( \ze \de ) } { n_{\ell} } <
\mscr{M}_{\mrm{B}} \li ( p^\star - \vep_a, \; p^\star \ri )$ and the
fact that $\mscr{M}_{\mrm{B}} ( z, z + \vep_a )$ is monotonically
increasing with respect to $z \in (0, p^\star - \vep_a)$ as asserted
by Lemma \ref{lemm21}, invoking the intermediate value theorem, we have that
there exists a unique number $z_{\ell} \in [0, p^\star - \vep_a)$
such that $\mscr{M}_{\mrm{B}} ( z_{\ell}, z_{\ell} + \vep_a )  = \f{
\ln ( \ze \de ) } { n_{\ell} }$, which implies Statement (I).

\bsk

{\bf Proof of Statement (II)}: By the definition of sample sizes, we
have \be \la{from3} \f{ \ln (\ze \de) } { \mscr{M}_{\mrm{B}} ( 1,
\f{1}{1 + \vep_r} )} \leq n_1 \leq n_{\ell} < \f{(1 + C_1) n_s}{2} <
\f{(1 + C_1) }{2} \li [ \f{ \ln ( \ze \de ) } { \mscr{M}_{\mrm{B}} (
p^\star + \vep_a, \; p^\star  ) } + 1 \ri ]
 \ee
and consequently, {\small $\f{ \ln ( \ze \de ) } { n_{\ell} } \geq
\mscr{M}_{\mrm{B}} ( 1, \f{1}{1 + \vep_r} )$},  {\small \bee \f{ \ln
( \ze \de ) } { n_{\ell} } <  \mscr{M}_{\mrm{B}} ( p^\star + \vep_a,
\; p^\star  )  \li ( \f{2}{1 + C_1}  - \f{1}{ n_{\ell} } \ri )  =
\li ( \f{2}{1 + C_1} \ri ) \mscr{M}_{\mrm{B}} \li ( p^\star +
\vep_a, \; \f{p^\star + \vep_a}{ 1 + \vep_r } \ri ) -
\f{\mscr{M}_{\mrm{B}} ( p^\star + \vep_a, \; p^\star  ) }{ n_{\ell}
}  \eee} for sufficiently small $\vep_a > 0$. Noting that {\small
$\lim_{\vep_a \to 0} \f{\mscr{M}_{\mrm{B}} ( p^\star + \vep_a, \;
p^\star  ) }{ n_{\ell} } = 0$}, we have {\small $\f{ \ln ( \ze \de )
} { n_{\ell} } < \mscr{M}_{\mrm{B}} ( p^\star + \vep_a, \;
\f{p^\star + \vep_a}{ 1 + \vep_r } )$} for small enough $\vep_a >
0$.  In view of the established fact that {\small
$\mscr{M}_{\mrm{B}} ( 1, \f{1}{1 + \vep_r} ) \leq \f{ \ln ( \ze \de
) } { n_{\ell} } < \mscr{M}_{\mrm{B}}  ( p^\star + \vep_a, \;
\f{p^\star + \vep_a}{ 1 + \vep_r }  )$} and the fact that
$\mscr{M}_{\mrm{B}} ( z, \f{z}{1 + \vep_r} )$ is monotonically
decreasing with respect to $z \in (0, 1)$ as asserted by Lemma
\ref{lemm23},  invoking the intermediate value theorem, we have that
there exists a unique number $y_{\ell} \in (p^\star + \vep_a, 1]$
such that $\mscr{M}_{\mrm{B}} ( y_{\ell}, \f{y_{\ell}}{1 + \vep_r} )
= \f{ \ln ( \ze \de ) } { n_{\ell} }$, which implies Statement (II).

\bsk

{\bf Proof of Statement (III)}: Since $n_{\ell}$ is monotonically
increasing with respect to $\ell$ if $\vep_a > 0$ is sufficiently
small, we have that $\mscr{M}_{\mrm{B}} ( z_{\ell}, z_{\ell} +
\vep_a )$ is monotonically increasing with respect to $\ell$ for
small enough $\vep_a > 0$.  Recalling that $\mscr{M}_{\mrm{B}} ( z,
z + \vep_a )$ is monotonically increasing with respect to $z \in (0,
p^\star - \vep_a)$, we have that $z_{\ell}$ is monotonically
increasing with respect to $\ell$. Similarly, $\mscr{M}_{\mrm{B}} (
y_{\ell}, \f{y_{\ell}}{1 + \vep_r} )$ is monotonically increasing
with respect to $\ell$ for sufficiently small $\vep_a > 0$.
Recalling that $\mscr{M}_{\mrm{B}} ( z, \f{z}{1 + \vep_r} )$ is
monotonically decreasing with respect to $z \in (0, 1)$, we have
that $y_{\ell}$ is monotonically decreasing with respect to $\ell$.
This establishes Statement (III).

\bsk

{\bf Proof of Statement (IV)}:  We first consider $\lim_{\vep_a \to
0} z_\ell$.  For simplicity of notations, define $b_\ell = \f{ 1 -
\sq{ 1 - 4 p^\star (1 - p^\star) C_{s - \ell}} }{2}$ for $\ell < s$
such that $n_\ell \geq \f{\ln (\ze \de)}{\ln ( 1 - \vep_a)}$. Then,
it can be checked that $\f{ b_\ell ( 1 - b_\ell) }{ p^\star (1 -
p^\star) } = C_{s - \ell}$ and, by the definition of sample sizes,
 we have \be \la{notiabs} \f{ b_\ell ( 1 -
b_\ell) }{ p^\star (1 - p^\star) }
 \f{ \mscr{M}_{\mrm{B}}
(z_\ell, z_\ell  + \vep_a ) } {\mscr{M}_{\mrm{B}} ( p^\star +
\vep_a,   p^\star  )} = \f{1}{n_{\ell} } \times \f{C_{s - \ell} \ln
(\ze \de) } { \mscr{M}_{\mrm{B}} ( p^\star + \vep_a,   p^\star ) } =
1 + o(1) \ee for $\ell < s$ such that $n_\ell \geq \f{\ln (\ze
\de)}{\ln ( 1 - \vep_a)}$.

We claim that $z_\ell
> \se$ for $\se \in (0, b_\ell)$ provided that $\vep_a > 0$ is sufficiently small.
Such a claim can be shown by a contradiction method as follows.  Suppose this claim is not true, then
there is a set, denoted by $S_{\vep_a}$,  of infinitely many values
of $\vep_a$ such that $z_\ell \leq \se$
for any $\vep_a \in S_{\vep_a}$.  For small enough $\vep_a \in S_{\vep_a}$,
it is true that  $z_\ell \leq \se < b_\ell < \f{1}{2} - \vep_a$.
By (\ref{notiabs}) and the fact that $\mscr{M}_{\mrm{B}}
(z, z  + \vep )$ is monotonically increasing with
respect to $z \in (0, \f{1}{2} - \vep)$ as asserted by Lemma \ref{lemm21}, we have
\[
\f{ b_\ell ( 1 - b_\ell) }{ p^\star (1 - p^\star) }
 \f{ \mscr{M}_{\mrm{B}}
(z_\ell, z_\ell  + \vep_a ) } {\mscr{M}_{\mrm{B}} ( p^\star + \vep_a,   p^\star  )}
= 1 + o(1)  \geq \f{ b_\ell ( 1 - b_\ell) }{ p^\star (1
- p^\star) } \f{ \mscr{M}_{\mrm{B}} (\se, \se  + \vep_a ) }
{\mscr{M}_{\mrm{B}} ( p^\star + \vep_a,   p^\star  )} = \f{b_\ell (1 -
b_\ell)}{\se (1 - \se)} + o(1)
\]
for small enough $\vep_a \in S_{\vep_a}$, which implies
$\f{b_\ell(1 - b_\ell)}{\se (1 - \se)} \leq 1$, contradicting to the
fact that $\f{b_\ell(1 - b_\ell)}{\se (1 - \se)} > 1$.  This proves
the claim. Now we restrict $\vep_a$ to be small enough so that $\se
< z_\ell < p^\star$. Making use of  (\ref{notiabs}) and applying Lemma \ref{lem32T}
based on the condition that $z_\ell \in ( \se, p^\star ) \subset (0, 1)$,
we have
\[
\f{ b_\ell ( 1 - b_\ell) }{ p^\star (1 - p^\star) }  \times  \f{
\vep_a^2 \sh [2 z_\ell  ( z_\ell - 1) ] + o (\vep_a^2)  } { \vep_a^2
\sh [ 2 p^\star (p^\star - 1)] + o (\vep_a^2)  } = 1 + o(1),
\]
which implies  {\small $\f{ b_\ell (1 - b_\ell) } { z_\ell ( 1 - z_\ell) }  = 1 +
o(1)$} and thus $\lim_{\vep_a \to 0} z_\ell = b_\ell$.

We now consider $\lim_{\vep_a \to 0} y_\ell$.  For simplicity of
notations, define $a_\ell = \f{1}{1 + \li ( \f{1}{p^\star} - 1  \ri
) C_{s - \ell}}$ for $1 \leq \ell < s$. Then, it can be checked that
$\f{ p^\star }{ 1 - p^\star } \f{1 - a_\ell}{a_\ell} = C_{s - \ell}$
and, by the definition of sample sizes,  \be \la{noti} \f{ p^\star
}{ 1 - p^\star } \f{1 - a_\ell}{a_\ell} \f{ \mscr{M}_{\mrm{B}}
(y_\ell, \f{y_\ell}{ 1 + \vep_r} ) } {\mscr{M}_{\mrm{B}} ( p^\star +
\vep_a,   p^\star  )} = \f{1}{n_{\ell} } \times \f{ C_{s - \ell} \ln
(\ze \de) } {\mscr{M}_{\mrm{B}} ( p^\star + \vep_a, p^\star )} = 1 +
o(1). \ee

We claim that $y_\ell < \se$ for $\se \in (a_\ell, 1)$ if
$\vep_r > 0$ is small enough.  To prove this claim, we use a contradiction method.
Suppose this claim is not true, then there is a set,
denoted by $S_{\vep_r}$,  of infinitely many values of $\vep_r$ such
that $y_\ell \geq \se$ for any $\vep_r \in S_{\vep_r}$.
By (\ref{noti}) and the fact that {\small $\mscr{M}_{\mrm{B}} (z, \f{z}{1 + \vep} )$}
is monotonically decreasing with respect to $z \in (0, 1)$ as asserted by Lemma \ref{lemm23}, we have
\[
\f{ p^\star }{ 1 - p^\star } \f{1 - a_\ell}{a_\ell}  \f{
\mscr{M}_{\mrm{B}} (y_\ell, \f{y_\ell}{1 + \vep_r} ) }
{\mscr{M}_{\mrm{B}} ( p^\star + \vep_a,   p^\star  )} = 1 + o(1)
\geq \f{ p^\star }{ 1 - p^\star } \f{1 - a_\ell}{a_\ell}  \f{
\mscr{M}_{\mrm{B}} (\se, \f{\se}{1 + \vep_r} ) }
{\mscr{M}_{\mrm{B}} ( p^\star + \vep_a,   p^\star  )} = \f{\se(1
- a_\ell)}{a_\ell (1 - \se)} + o(1)
\]
for small enough $\vep_r \in S_{\vep_r}$, which implies $\f{\se(1 -
a_\ell)}{a_\ell (1 - \se)} \leq 1$, contradicting to the fact that
$\f{\se(1 - a_\ell)}{a_\ell (1 - \se)} > 1$.  This proves the claim.
Now we restrict $\vep_r$ to be small enough so that $p^\star <
y_\ell < \se$.  By (\ref{noti}) and applying Lemma \ref{lem32T}
based on the condition that $y_\ell \in ( p^\star, \se) \subset (0, 1)$, we have
\[
\f{ p^\star }{ 1 - p^\star } \f{1 - a_\ell}{a_\ell} \times  \f{
\vep_r^2 y_\ell \sh [ 2 (y_\ell - 1) ]  + o (\vep_r^2)  } {
\vep_a^2 \sh [ 2 p^\star  (p^\star - 1)] +  o (\vep_a^2)  } = 1 +
o(1),
\]
which implies {\small $\f{ y_\ell  - a_\ell } { a_\ell ( 1 - y_\ell) }  =  o(1)$}
and thus $\lim_{\vep_r \to 0} y_\ell = a_\ell$.

\bsk

{\bf Proof of Statement (V)}:

We shall first consider $p \in (p^\star, 1)$.  For small enough $\vep_r > 0$, there exists $z_{\ell_\vep} \in (p^\star, 1)$ such that\[
n_{\ell_\vep} =  \f{ \ln (\ze \de) } { \mscr{M}_{\mrm{B}} (z_{\ell_\vep}, z_{\ell_\vep} \sh (1 + \vep_r) ) }  = \li \lc \f{ C_{s - {\ell_\vep}}
\; \ln (\ze \de) } {\mscr{M}_{\mrm{B}} (p^\star + \vep_a, p^\star)} \ri \rc = \li \lc \f{ p^\star }{ 1 - p^\star } \f{1 - p}{p}  \f{ \ln (\ze
\de) } {\mscr{M}_{\mrm{B}} (p^\star + \vep_a, p^\star)}  \ri \rc.
\]
For $\se \in (p, 1)$, we claim that $z_{\ell_\vep} < \se$ if $\vep_r$ is sufficiently small.  Suppose, to get a contradiction, that this claim
is not true.  Then,  there exists a set, denoted by $S_{\vep_r}$, of infinitely many values of $\vep_r$ such that $z_{\ell_\vep} \geq \se$ for
any value of $\vep_r$ in $S_{\vep_r}$. Noting that \be \la{noti} \f{ \f{ p^\star }{ 1 - p^\star } \f{1 - p}{p} \f{ \ln (\ze \de) }
{\mscr{M}_{\mrm{B}} (p^\star + \vep_a, p^\star)}  } { \f{ \ln (\ze \de) } { \mscr{M}_{\mrm{B}} (z_{\ell_\vep}, z_{\ell_\vep} \sh (1 + \vep_r) )
}  } =  \f{ p^\star }{ 1 - p^\star } \f{1 - p}{p}  \f{ \mscr{M}_{\mrm{B}} (z_{\ell_\vep}, z_{\ell_\vep} \sh (1 + \vep_r) ) } {\mscr{M}_{\mrm{B}}
(p^\star + \vep_a, p^\star)} = 1 + o(\vep_r), \ee we have
\[
\f{ p^\star }{ 1 - p^\star } \f{1 - p}{p}  \f{ \mscr{M}_{\mrm{B}} (z_{\ell_\vep}, z_{\ell_\vep} \sh (1 + \vep_r) ) } {\mscr{M}_{\mrm{B}}
(p^\star + \vep_a, p^\star)} = 1 + o(\vep_r) \geq \f{ p^\star }{ 1 - p^\star } \f{1 - p}{p}  \f{ \mscr{M}_{\mrm{B}} (\se, \se \sh (1 + \vep_r) )
} {\mscr{M}_{\mrm{B}} (p^\star + \vep_a, p^\star)} = \f{\se(1 - p)}{p (1 - \se)} + o(1)
\]
for any value of $\vep_r$ in $S_{\vep_r}$, which contradicts to the fact that $\f{\se(1 - p)}{p (1 - \se)} > 1$.  This proves the claim. Now we
restrict $\vep_r$ to be small enough so that $p^\star < z_{\ell_\vep} < \se$.  Using $\ln (1 + x) = x - \f{x^2}{2} + \f{x^3}{3} + o (x^3)$ for
$|x| < 1$, we can show that  \be \la{nowadd}  \mscr{M}_{\mrm{B}} (p + \vep, p) = - \f{\vep^2} { 2 p ( 1 - p)}  + \f{ 1 - 2 p } { 6 [p ( 1 -
p)]^2  } \vep^3 + o (\vep^3) \ee for $0 < \vep < 1 - p$. Since $z_{\ell_\vep}$ is bounded with respect to $\vep$, by (\ref{noti}),
(\ref{nowadd})  and Lemma \ref{lem32T}, we have
\[
\f{ p^\star }{ 1 - p^\star } \f{1 - p}{p} \times  \f{ - \vep_r^2 z_{\ell_\vep} \sh [ 2 ( 1 - z_{\ell_\vep}) ] +  \vep_r^3 z_{\ell_\vep} (2 -
z_{\ell_\vep}) \sh [ 3 (1 - z_{\ell_\vep})^2 ] + o (\vep_r^3) } { - \vep_r^2 p^\star \sh [ 2 (1 - p^\star)] + \vep_r^3 p^\star (1 - 2 p^\star)
\sh [ 6 (1 - p^\star)^2 ] + o (\vep_r^3)  } = 1 + o(\vep_r),
\]
i.e.,
\[
\f{ \f{ z_{\ell_\vep} (1 - p) } { p ( 1 - z_{\ell_\vep}) } - \f{ 2 \vep_r z_{\ell_\vep} (1 - p) (2 - z_{\ell_\vep}) } { 3 p (1 -
z_{\ell_\vep})^2 } + o (\vep_r) }{1 -   \vep_r (1 - 2 p^\star) \sh [ 3 (1 - p^\star) ] + o(\vep_r)} = 1 + o(\vep_r),
\]
i.e.,
\[
\f{ z_{\ell_\vep}  - p } { p ( 1 - z_{\ell_\vep}) } - \f{ 2 \vep_r z_{\ell_\vep} (1 - p) (2 - z_{\ell_\vep}) } { 3 p (1 - z_{\ell_\vep})^2 } = -
\f{ \vep_r (1 - 2 p^\star)}{ 3 (1 - p^\star) } + o(\vep_r),
\]
i.e.,
\[
\f{ z_{\ell_\vep}  - p } { p  } - \f{ 2 \vep_r z_{\ell_\vep} (1 - p) (2 - z_{\ell_\vep}) } { 3 p (1 - z_{\ell_\vep}) } = - \f{ \vep_r (1 - 2
p^\star) (1 - z_{\ell_\vep})} { 3 (1 - p^\star) } + o(\vep_r),
\]
which implies that $\lim_{\vep_r \to 0} z_{\ell_\vep} = p$ and consequently,
\[
\lim_{\vep_r \to 0} \f{ z_{\ell_\vep} - p}{\vep_r p} = \f{2 (2 - p)}{3} - \f{ ( 1 - 2 p^\star) (1 - p)} { 3 (1 - p^\star)} = 1 - \f{1}{3} \f{ p
- p^\star } {1 - p^\star } \in \li (\f{2}{3}, 1 \ri ).
\]

Next, we shall consider $p \in (0, p^\star)$.  For small enough $\vep_a > 0$, there exists $z_{\ell_\vep} \in (0, p^\star)$ such that\[
n_{\ell_\vep} =  \f{ \ln (\ze \de) } { \mscr{M}_{\mrm{B}} (z_{\ell_\vep}, z_{\ell_\vep} + \vep_a )  }  = \li \lc \f{ C_{s - {\ell_\vep}} \; \ln
(\ze \de) } {\mscr{M}_{\mrm{B}} ( p^\star + \vep_a, p^\star )} \ri \rc = \li \lc \f{ p (1 - p)}{ p^\star (1 - p^\star) } \f{ \ln (\ze \de) }
{\mscr{M}_{\mrm{B}} (p^\star + \vep_a, p^\star)} \ri \rc.
\]
For $\se \in (0, p)$, we claim that $z_{\ell_\vep} > \se$ if $\vep_a$ is sufficiently small.  Suppose, to get a contradiction, that this claim
is not true. Then,  there exists a set, denoted by $S_{\vep_a}$, of infinitely many values of $\vep_a$ such that $z_{\ell_\vep} \leq \se$ for
any value of $\vep_a$ in $S_{\vep_a}$. Noting that \be \la{notiabs} \f{ \f{ p ( 1 - p) }{ p^\star (1 - p^\star) } \f{ \ln (\ze \de) }
{\mscr{M}_{\mrm{B}} ( p^\star + \vep_a, p^\star )}  } { \f{ \ln (\ze \de) } { \mscr{M}_{\mrm{B}} (z_{\ell_\vep}, z_{\ell_\vep} + \vep_a )  }  }
= \f{ p ( 1 - p) }{ p^\star (1 - p^\star) }
 \f{ \mscr{M}_{\mrm{B}}
(z_{\ell_\vep}, z_{\ell_\vep}  + \vep_a ) } {\mscr{M}_{\mrm{B}} ( p^\star + \vep_a, p^\star )} = 1 + o(\vep_a), \ee we have
\[
\f{ p ( 1 - p) }{ p^\star (1 - p^\star) }
 \f{ \mscr{M}_{\mrm{B}}
(z_{\ell_\vep}, z_{\ell_\vep}  + \vep_a ) } {\mscr{M}_{\mrm{B}} ( p^\star + \vep_a, p^\star )} = 1 + o(\vep_a)  > \f{ p ( 1 - p) }{ p^\star (1 -
p^\star) } \f{ \mscr{M}_{\mrm{B}} (\se, \se  + \vep_a ) } {\mscr{M}_{\mrm{B}} ( p^\star + \vep_a, p^\star )} = \f{p(1 - p)}{\se (1 - \se)} +
o(1)
\]
for any value of $\vep_a$ in $S_{\vep_a}$, which contradicts to the fact that $\f{p(1 - p)}{\se (1 - \se)} > 1$.  This proves the claim. Now we
restrict $\vep_a$ to be small enough so that $\se < z_{\ell_\vep} < p^\star$. Since $z_{\ell_\vep}$ is bounded with respect to $\vep$, by
(\ref{nowadd}),  (\ref{notiabs})  and Lemma \ref{lem32T}, we have
\[
\f{ p ( 1 - p) }{ p^\star (1 - p^\star) }  \times  \f{ - \vep_a^2 \sh [2 z_{\ell_\vep}  ( 1 - z_{\ell_\vep}) ] +  \vep_a^3 (1 - 2 z_{\ell_\vep})
\sh [ 3 z_{\ell_\vep}^2 (1 - z_{\ell_\vep})^2 ] + o (\vep_a^3)  } { - \vep_a^2 \sh [ 2 p^\star (1 - p^\star)] + \vep_a^3 (1 - 2 p^\star) \sh [ 6
(p^\star )^2(1 - p^\star)^2 ] + o (\vep_a^3) } = 1 + o(\vep_a),
\]
i.e.,
\[
\f{ \f{ p (1 - p) } { z_{\ell_\vep} ( 1 - z_{\ell_\vep}) } - \f{ 2 \vep_a p (1 - p) (1 - 2 z_{\ell_\vep}) } { 3 z_{\ell_\vep}^2 (1 -
z_{\ell_\vep})^2 } + o (\vep_a) }{1 -  \vep_a (1 - 2 p^\star) \sh [ 3 p^\star (1 - p^\star) ] + o(\vep_a)} = 1 + o(\vep_a),
\]
i.e.,
\[
\f{ p (1 - p) } { z_{\ell_\vep} ( 1 - z_{\ell_\vep}) } - \f{ 2 \vep_a p (1 - p) (1 - 2 z_{\ell_\vep}) } { 3 z_{\ell_\vep}^2 (1 -
z_{\ell_\vep})^2 } = 1 - \f{ \vep_a (1 - 2 p^\star) }{ 3 p^\star (1 - p^\star) } + o(\vep_a),
\]
i.e.,
\[
\f{(z_{\ell_\vep} - p)(1 - z_{\ell_\vep} - p)}{ z_{\ell_\vep} ( 1 - z_{\ell_\vep}) } = \f{\vep_a (1 - 2 p^\star) }{ 3 p^\star (1 - p^\star) } -
\f{ 2 \vep_a p (1 - p) (1 - 2 z_{\ell_\vep}) } { 3 z_{\ell_\vep}^2 (1 - z_{\ell_\vep})^2 } + o(\vep_a),
\]
i.e.,
\[
\f{z_{\ell_\vep} - p}{ \vep_a } = \f{ z_{\ell_\vep} ( 1 - z_{\ell_\vep}) (1 - 2 p^\star) }{ 3 p^\star (1 - p^\star) (1 - z_{\ell_\vep} - p) } -
\f{ 2 p (1 - p) (1 - 2 z_{\ell_\vep}) } { 3 z_{\ell_\vep} (1 - z_{\ell_\vep}) (1 - z_{\ell_\vep} - p) } + o(1),
\]
which implies that $\lim_{\vep_r \to 0} z_{\ell_\vep} = p$ and consequently,
\[
\lim_{\vep_a \to 0} \f{ z_{\ell_\vep} - p}{\vep_a} = \f{ p ( 1 - p) (1 - 2 p^\star) }{ 3 p^\star (1 - p^\star) (1 - 2 p) } - \f{2}{3} = - \nu
\in \li ( - \f{2}{3}, - \f{1}{3} \ri ).
\]

\bsk

{\bf Proof of Statement (VI)}:  By the definition of the sampling
scheme, {\small \bee  &   & \{ \bs{D}_{\ell} = 0  \}  =  \li \{ \max
\{ \mscr{M}_{\mrm{B}} (\wh{\bs{p}}_\ell, \mscr{L} ( \wh{\bs{p}}_\ell
) ), \; \mscr{M}_{\mrm{B}} (\wh{\bs{p}}_\ell, \mscr{U} (
\wh{\bs{p}}_\ell ) ) \} > \f{ \ln ( \ze \de  ) } { n_\ell }, \;  |
\wh{\bs{p}}_\ell - p^\star | \leq \vep_a \ri \}\\
&  & \qu \bigcup \li \{ \max \{ \mscr{M}_{\mrm{B}}
(\wh{\bs{p}}_\ell, \mscr{L} ( \wh{\bs{p}}_\ell ) ), \;
\mscr{M}_{\mrm{B}} (\wh{\bs{p}}_\ell, \mscr{U} ( \wh{\bs{p}}_\ell )
) \} > \f{ \ln ( \ze \de  ) } { n_\ell }, \;
\wh{\bs{p}}_\ell  < p^\star - \vep_a \ri \}\\
&  & \qu \bigcup \li \{ \max \{ \mscr{M}_{\mrm{B}}
(\wh{\bs{p}}_\ell, \mscr{L} ( \wh{\bs{p}}_\ell ) ), \;
\mscr{M}_{\mrm{B}} (\wh{\bs{p}}_\ell, \mscr{U} ( \wh{\bs{p}}_\ell )
) \} > \f{ \ln ( \ze \de  ) } { n_\ell }, \;
\wh{\bs{p}}_\ell  > p^\star + \vep_a \ri \}\\
&   & \qu = \li \{ \max \li \{ \mscr{M}_{\mrm{B}} (\wh{\bs{p}}_\ell,
\wh{\bs{p}}_\ell - \vep_a), \; \mscr{M}_{\mrm{B}} \li (
\wh{\bs{p}}_\ell, \f{\wh{\bs{p}}_\ell}{1 - \vep_r} \ri ) \ri \} >
\f{\ln ( \ze \de  ) } { n_\ell }, \;  | \wh{\bs{p}}_\ell -  p^\star | \leq \vep_a \ri \}\\
&  & \qu \bigcup \li \{ \mscr{M}_{\mrm{B}} (\wh{\bs{p}}_\ell,
\wh{\bs{p}}_\ell + \vep_a ) > \f{ \ln ( \ze \de  ) } { n_\ell }, \;
\wh{\bs{p}}_\ell  < p^\star - \vep_a \ri \} \bigcup \li \{
\mscr{M}_{\mrm{B}} \li ( \wh{\bs{p}}_\ell, \f{ \wh{\bs{p}}_\ell}{1 +
\vep_r} \ri ) > \f{ \ln ( \ze \de  ) } { n_\ell }, \;
\wh{\bs{p}}_\ell  > p^\star + \vep_a \ri \}. \qqu \eee} We claim
that if $\vep_a > 0$ is sufficiently small, then it is true that
{\small \bel &  &  \li \{ \max \li \{ \mscr{M}_{\mrm{B}}
(\wh{\bs{p}}_\ell, \wh{\bs{p}}_\ell - \vep_a), \; \mscr{M}_{\mrm{B}}
\li ( \wh{\bs{p}}_\ell, \f{\wh{\bs{p}}_\ell}{1 - \vep_r} \ri ) \ri
\} > \f{ \ln ( \ze \de  ) } { n_\ell }, \; | \wh{\bs{p}}_\ell -
p^\star | \leq \vep_a \ri \}  =   \li \{  |
\wh{\bs{p}}_\ell - p^\star | \leq \vep_a \ri \}, \qqu \qu \la{that0}\\
 &   &  \li \{ \mscr{M}_{\mrm{B}} (\wh{\bs{p}}_\ell,
\wh{\bs{p}}_\ell + \vep_a ) > \f{ \ln ( \ze \de ) } { n_\ell }, \;
\wh{\bs{p}}_\ell  < p^\star - \vep_a \ri \} =  \{  z_\ell <
\wh{\bs{p}}_\ell  < p^\star - \vep_a \} \; \; \tx{for} \; \f{ \ln
(\ze \de) } { \ln ( 1 - \vep_a) } \leq n_\ell < n_s, \qqu \qu \la{that}\\
&  &  \li \{ \mscr{M}_{\mrm{B}} (\wh{\bs{p}}_\ell, \wh{\bs{p}}_\ell
+ \vep_a ) > \f{ \ln ( \ze \de ) } { n_\ell }, \; \wh{\bs{p}}_\ell <
p^\star - \vep_a \ri \} = \{  0 < \wh{\bs{p}}_\ell  < p^\star -
\vep_a \} \; \; \tx{for} \;  n_1 \leq
n_\ell < \f{ \ln (\ze \de) } { \ln ( 1 - \vep_a) }, \qqu \qu \la{that2}\\
&    &  \li \{ \mscr{M}_{\mrm{B}} \li ( \wh{\bs{p}}_\ell,
\f{\wh{\bs{p}}_\ell}{ 1  + \vep_r} \ri ) > \f{ \ln ( \ze \de  ) } {
n_\ell }, \; \wh{\bs{p}}_\ell  > p^\star + \vep_a \ri \} =  \li \{
p^\star + \vep_a <  \wh{\bs{p}}_\ell  < y_\ell \ri \}. \la{that3}
 \eel}
To show (\ref{that0}), note that \be \la{since}
 n_{\ell} < \f{(1 + C_1) n_s}{2} < \f{(1 + C_1) }{2} \li [ \f{
\ln ( \ze \de ) } { \mscr{M}_{\mrm{B}} ( p^\star + \vep_a,   p^\star
) } + 1 \ri ],  \ee  which implies that
\[
\f{ \ln ( \ze \de ) } { n_{\ell} } < \f{ \mscr{M}_{\mrm{B}} (
p^\star + \vep_a,   p^\star  ) } {\mscr{M}_{\mrm{B}} ( p^\star -
\vep_a, \; p^\star - \vep_a - \vep_a )}  \li ( \f{2}{1 + C_1} \ri )
\mscr{M}_{\mrm{B}} ( p^\star - \vep_a, \; p^\star - \vep_a - \vep_a
) - \f{\mscr{M}_{\mrm{B}} ( p^\star + \vep_a,   p^\star )}{n_\ell}
\]
if $\vep_a > 0$ is sufficiently small.  Noting that
\[
\lim_{\vep_a \to 0} \f{ \mscr{M}_{\mrm{B}} ( p^\star + \vep_a,   p^\star  )  }
{  \mscr{M}_{\mrm{B}} ( p^\star - \vep_a, \; p^\star -
\vep_a - \vep_a ) } = \lim_{\vep_a \to 0}  \f{ \f{ \vep_a^2 } { 2
p^\star  (  p^\star - 1 )   } + o(\vep_a^2)  } {  \f{ \vep_a^2 } {
2 (p^\star - \vep_a)  ( p^\star -  \vep_a  - 1)   } + o(\vep_a^2) }
= 1
\]
and $\lim_{\vep_a \to 0} \f{\mscr{M}_{\mrm{B}} ( p^\star + \vep_a,   p^\star  )}{n_\ell} = 0$,
we have \be \la{have1} \f{ \ln ( \ze \de )
} { n_{\ell} } < \mscr{M}_{\mrm{B}} ( p^\star - \vep_a, \; p^\star -
\vep_a - \vep_a ) \ee for small enough $\vep_a > 0$.  Again by
(\ref{since}), we have
\[
\f{ \ln ( \ze \de ) } { n_{\ell} } < \f{ \mscr{M}_{\mrm{B}} (
p^\star + \vep_a,   p^\star  ) } { \mscr{M}_{\mrm{B}}  ( p^\star +
\vep_a,  \f{p^\star + \vep_a}{1 -  \vep_r}  ) }  \li ( \f{2}{1 +
C_1} \ri ) \mscr{M}_{\mrm{B}} \li ( p^\star + \vep_a, \; \f{p^\star
+ \vep_a}{1 -  \vep_r} \ri ) - \f{\mscr{M}_{\mrm{B}} ( p^\star +
\vep_a, p^\star )}{n_\ell}
\]
if $\vep_a > 0$ is sufficiently small. Noting that
\[
\lim_{\vep_a \to 0} \f{ \mscr{M}_{\mrm{B}} ( p^\star + \vep_a,   p^\star  ) }
{ \mscr{M}_{\mrm{B}} ( p^\star + \vep_a, \f{p^\star +
\vep_a}{1 -  \vep_r} ) } = \lim_{\vep_a \to 0}  \f{ \f{ \vep_a^2 }
{ 2 p^\star  ( p^\star - 1)   } + o(\vep_a^2)  } { \f{ \vep_a^2
} { 2 (p^\star + \vep_a)  ( p^\star +  \vep_a - 1)   } + o \li (
\f{(p^\star + \vep_a)^2 \vep_r^2}{(1 - \vep_r)^2} \ri ) } = 1
\]
and $\lim_{\vep_a \to 0} \f{\mscr{M}_{\mrm{B}} ( p^\star + \vep_a,   p^\star  )}{n_\ell} = 0$, we have \be \la{have2}
 \f{ \ln ( \ze \de )
} { n_{\ell} } < \mscr{M}_{\mrm{B}} \li ( p^\star + \vep_a, \;
\f{p^\star + \vep_a}{1 -  \vep_r } \ri ) \ee for small enough
$\vep_a
> 0$.  It can be seen from Lemmas \ref{lemm21}
and \ref{lemm24} that, for $z \in [p^\star - \vep_a, \; p^\star +
\vep_a]$, $\mscr{M}_{\mrm{B}} ( z, z - \vep_a)$ is monotonically
increasing with respect to $z$ and $\mscr{M}_{\mrm{B}} ( z, \f{z}{1
- \vep_r} )$ is monotonically decreasing with respect to $z$. By
(\ref{have1}) and (\ref{have2}), we have {\small $\f{ \ln ( \ze \de
) } { n_{\ell} } < \mscr{M}_{\mrm{B}} ( z, z - \vep_a )$} and
{\small $\f{ \ln ( \ze \de ) } { n_{\ell} } < \mscr{M}_{\mrm{B}} (
z, \f{z}{1 - \vep_r} )$} for any $z \in [ p^\star - \vep_a,
p^\star + \vep_a ]$ if  $\vep_a > 0$ is small enough. This proves
(\ref{that0}).

To show (\ref{that}), let $\om \in \{ \mscr{M}_{\mrm{B}} \li (
\wh{\bs{p}}_{\ell}, \wh{\bs{p}}_{\ell}  + \vep_a \ri )  > \f{\ln
(\ze \de)} {n_{\ell}}, \wh{\bs{p}}_{\ell} < p^\star - \vep_a \}$
and $\wh{p}_{\ell} = \wh{\bs{p}}_{\ell} (\om)$. Then,
$\mscr{M}_{\mrm{B}} ( \wh{p}_{\ell}, \; \wh{p}_{\ell}   + \vep_a )
> \f{\ln (\ze \de)} {n_{\ell}}$ and $\wh{p}_{\ell} < p^\star - \vep_a$.
Since $z_{\ell} \in [0, p^\star - \vep_a)$ and $\mscr{M}_{\mrm{B}}
\li ( z, z + \vep_a \ri )$ is monotonically increasing with
respect to $z \in (0, p^\star - \vep_a)$, it must be true that
$\wh{p}_{\ell} > z_{\ell}$. Otherwise if $\wh{p}_{\ell} \leq
z_{\ell}$, then $\mscr{M}_{\mrm{B}} \li ( \wh{p}_{\ell},
\wh{p}_{\ell}   + \vep_a \ri ) \leq \mscr{M}_{\mrm{B}} \li (
z_{\ell}, z_{\ell}   + \vep_a \ri ) = \f{\ln (\ze \de)}
{n_{\ell}}$, leading to a contradiction. This proves {\small $ \{
\mscr{M}_{\mrm{B}} \li ( \wh{\bs{p}}_{\ell}, \wh{\bs{p}}_{\ell} +
\vep_a \ri ) > \f{\ln (\ze \de)} {n_{\ell}} , \; \wh{\bs{p}}_{\ell}
< p^\star - \vep_a \} \subseteq  \{ z_{\ell} < \wh{\bs{p}}_{\ell} <
p^\star - \vep_a \}$}.  Now let $\om \in \li \{ z_{\ell} <
\wh{\bs{p}}_{\ell} < p^\star - \vep_a \ri \}$ and $\wh{p}_{\ell} =
\wh{\bs{p}}_{\ell} (\om)$.  Then, $z_{\ell} <
\wh{p}_{\ell} < p^\star - \vep_a$.  Noting that
$\mscr{M}_{\mrm{B}}  \li ( z, z + \vep_a \ri )$ is monotonically increasing with respect to $z
\in (0,  p^\star - \vep_a)$, we have that {\small
$\mscr{M}_{\mrm{B}} \li ( \wh{p}_{\ell},  \wh{p}_{\ell}   + \vep_a
\ri ) > \mscr{M}_{\mrm{B}} \li ( z_{\ell}, z_{\ell} + \vep_a \ri )
= \f{\ln (\ze \de)} {n_{\ell}}$}, which implies {\small $\{
\mscr{M}_{\mrm{B}} \li ( \wh{\bs{p}}_{\ell}, \wh{\bs{p}}_{\ell} +
\vep_a \ri ) > \f{\ln (\ze \de)} {n_{\ell}} , \wh{\bs{p}}_{\ell} <
p^\star - \vep_a \} \supseteq \{ z_{\ell} < \wh{\bs{p}}_{\ell} <
p^\star - \vep_a \}$}.  This establishes (\ref{that}).

Note that, for any $z \in (0, p^\star - \vep_a)$, we have
$\mscr{M}_{\mrm{B}} ( z, z + \vep_a ) > \mscr{M}_{\mrm{B}} ( 0,
\vep_a ) = \ln ( 1 - \vep_a) \geq \f{\ln (\ze \de)} {n_{\ell}}$, which
implies (\ref{that2}).

To show (\ref{that3}), let {\small $\om \in \{ \mscr{M}_{\mrm{B}}  (
\wh{\bs{p}}_{\ell}, \f{\wh{\bs{p}}_\ell}{ 1  + \vep_r} )
> \f{\ln (\ze \de)} {n_{\ell}}, \; \wh{\bs{p}}_{\ell} > p^\star +
\vep_a \}$} and $\wh{p}_{\ell} = \wh{\bs{p}}_{\ell} (\om)$. Then,
{\small $\mscr{M}_{\mrm{B}} ( \wh{p}_{\ell}, \f{\wh{p}_{\ell}}{ 1 +
\vep_r} ) > \f{\ln (\ze \de)} {n_{\ell}}$ and $\wh{p}_{\ell} >
p^\star + \vep_a$}. Since $y_{\ell} \in (p^\star + \vep_a, 1]$ and
{\small $\mscr{M}_{\mrm{B}} ( z, \f{z}{1 + \vep_r} )$} is
monotonically decreasing with respect to $z \in (p^\star + \vep_a,
1)$, it must be true that $\wh{p}_{\ell} < y_{\ell}$. Otherwise if
$\wh{p}_{\ell} \geq y_{\ell}$, then {\small $\mscr{M}_{\mrm{B}} (
\wh{p}_{\ell}, \f{\wh{p}_{\ell}}{1   + \vep_r} ) \leq
\mscr{M}_{\mrm{B}} ( y_{\ell}, \f{y_{\ell}}{ 1 + \vep_r } ) = \f{\ln
(\ze \de)} {n_{\ell}}$}, leading to a contradiction. This proves
{\small $ \{ \mscr{M}_{\mrm{B}} ( \wh{\bs{p}}_{\ell}, \f{
\wh{\bs{p}}_{\ell}}{1 + \vep_r}  ) > \f{\ln (\ze \de)} {n_{\ell}} ,
\; \wh{\bs{p}}_{\ell} > p^\star + \vep_a \} \subseteq \{ p^\star +
\vep_a < \wh{\bs{p}}_{\ell} < y_\ell \}$}.  Now let $\om \in \li \{
 p^\star + \vep_a < \wh{\bs{p}}_{\ell} < y_{\ell}
\ri \}$ and $\wh{p}_{\ell} = \wh{\bs{p}}_{\ell} (\om)$.
Then, $p^\star + \vep_a < \wh{p}_{\ell} < y_{\ell}$. Noting that
{\small $\mscr{M}_{\mrm{B}} ( z, \f{ z }{1 +  \vep_r} )$} is
monotonically decreasing with respect to $z \in (0,
1)$, we have that {\small $\mscr{M}_{\mrm{B}} ( \wh{p}_{\ell}, \f{
\wh{p}_{\ell} }{ 1 + \vep_r} ) > \mscr{M}_{\mrm{B}} ( y_{\ell}, \f{
y_{\ell} }{1 + \vep_r} ) = \f{\ln (\ze \de)} {n_{\ell}}$},
which implies {\small $\{ \mscr{M}_{\mrm{B}} ( \wh{\bs{p}}_{\ell}, \f{
\wh{\bs{p}}_{\ell} }{ 1 + \vep_r} ) > \f{\ln (\ze \de)} {n_{\ell}} ,
\; \wh{\bs{p}}_{\ell} > p^\star + \vep_a \} \supseteq \{  p^\star +
\vep_a < \wh{\bs{p}}_{\ell} < y_{\ell} \}$}. This establishes
(\ref{that3}).

 \bsk

 \beL \la{Defmix} Let $\ell_\vep = s  - j_p$.  Then, under the constraint that limits are taken
with $\f{\vep_a}{\vep_r}$ fixed, \be \la{goodone} \lim_{\vep_a \to
0} \sum_{\ell = 1}^{\ell_\vep - 1} n_\ell \Pr \{ \bs{D}_\ell = 1 \}
= 0, \qqu \lim_{\vep_a \to 0} \sum_{\ell = \ell_\vep + 1}^s n_\ell
\Pr \{ \bs{D}_{\ell} = 0 \} = 0 \ee for $p \in (0, 1)$.  Moreover,
$\lim_{\vep_a \to 0} n_{\ell_\vep} \Pr \{ \bs{D}_{\ell_\vep} = 0 \}
= 0$ if $C_{j_p} > r (p)$. \eeL

\bpf

For simplicity of notations, let $a_\ell = \lim_{\vep_a \to 0}
y_\ell$ and $b_\ell = \lim_{\vep_a \to 0} z_\ell$.
The proof consists of three main  steps as follows.

First, we shall show that (\ref{goodone}) holds for $p \in (0,
p^\star]$. By the definition of $\ell_\vep$, we have {\small $r (p)
> C_{s - \ell_\vep + 1}$}. Making use of the first four statements
of Lemma \ref{lem52m}, we have that {\small $z_\ell < \f{p +
b_{\ell_\vep - 1}}{2} < p$} for all $\ell \leq \ell_\vep - 1$ with
$n_\ell \geq \f{ \ln (\ze \de) } { \ln (1 - \vep_a) }$ and that
{\small $y_\ell  > \f{ p^\star + a_{s - 1} }{2}
> p^\star$} for $1 \leq \ell < s$ if $\vep_a$ is sufficiently small.
Therefore, by the last statement of Lemma \ref{lem52m} and using
Chernoff bound, we have that {\small \bee \Pr \{ \bs{D}_{\ell} = 1
\} & = & \Pr \{ \wh{\bs{p}}_{\ell} \leq z_{\ell} \} + \Pr \{
\wh{\bs{p}}_{\ell} \geq y_{\ell} \} \leq \Pr \li \{
\wh{\bs{p}}_{\ell} \leq \f{p + b_{\ell_\vep -
1}}{2} \ri \} + \Pr \li \{ \wh{\bs{p}}_{\ell} \geq \f{ p^\star + a_{s - 1}  }{2} \ri \}\\
& \leq & \exp \li ( - 2 n_\ell \li( \f{p - b_{\ell_\vep - 1}}{2}
\ri)^2 \ri ) + \exp \li ( - 2 n_\ell \li ( \f{ p^\star + a_{s - 1}
}{2} - p \ri )^2 \ri ) \eee} for all $\ell \leq \ell_\vep -
1$ with $n_\ell \geq \f{ \ln (\ze \de) } { \ln (1 - \vep_a) }$  and that \bee \Pr \{ \bs{D}_{\ell} = 1 \}  & =  & \Pr \{
\wh{\bs{p}}_{\ell} \geq y_{\ell} \} + \Pr \{
\wh{\bs{p}}_{\ell} = 0 \} \leq  \Pr \li \{
\wh{\bs{p}}_{\ell} \geq \f{ p^\star + a_{s - 1}  }{2} \ri \} + \Pr \{
\wh{\bs{p}}_{\ell} = 0 \}\\
& \leq & \exp \li ( - 2 n_\ell \li ( \f{ p^\star + a_{s - 1}  }{2} -
p \ri )^2 \ri ) + \exp ( - 2 n_\ell p^2) \eee for all $\ell$ with
$n_\ell < \f{ \ln (\ze \de) } { \ln (1 - \vep_a) }$ if $\vep_a > 0$
is small enough.  As a consequence of the definition of $\ell_\vep$,
we have that $b_{\ell_\vep - 1}$ is smaller than $p$ and is
independent of $\vep_a > 0$. Hence, we can apply  Lemma \ref{lem31a}
to conclude that {\small $\lim_{\vep_a \to 0} \sum_{\ell =
1}^{\ell_\vep - 1} n_\ell \Pr \{ \bs{D}_\ell = 1 \}  = 0$}.

 Similarly, it can be seen from the definition of $\ell_\vep$ that
 $r (p)  < C_{s - \ell_\vep - 1}$.
  Making use of the first four statements of
Lemma \ref{lem52m}, we have that {\small $z_\ell > \f{p +
b_{\ell_\vep + 1}}{2} > p$} for $\ell_\vep + 1 \leq \ell < s$ if
$\vep_a$ is sufficiently small. By the last statement of Lemma
\ref{lem52m} and using Chernoff bound, we have {\small \[ \Pr \{
\bs{D}_{\ell} = 0 \} = \Pr \{  z_{\ell} < \wh{\bs{p}}_{\ell} <
y_{\ell} \} \leq \Pr \{ \wh{\bs{p}}_{\ell} > z_{\ell} \} \leq \Pr
\li \{ \wh{\bs{p}}_{\ell} > \f{p + b_{\ell_\vep + 1}}{2} \ri \} \leq
\exp \li ( - 2 n_\ell \li ( \f{p - b_{\ell_\vep + 1}}{2} \ri )^2 \ri
)
\]}
for $\ell_\vep + 1 \leq \ell < s$ if $\vep_a
> 0$ is small enough.  By virtue of the definition of $\ell_\vep$,
we have that $b_{\ell_\vep + 1}$ is greater than $p$ and is
independent of $\vep_a > 0$. In view of this and the fact that $\Pr
\{ \bs{D}_s = 0 \} = 0$, we can use Lemma \ref{lem31a} to arrive at
{\small $\lim_{\vep_a \to 0} \sum_{\ell = \ell_\vep + 1}^s n_\ell
\Pr \{ \bs{D}_\ell = 0 \} = 0$}.  This proves that (\ref{goodone})
holds for $p \in (0, p^\star]$.

\bsk

Second, we shall show that (\ref{goodone}) holds for $p \in
(p^\star, 1)$. As a direct consequence of the definition of
$\ell_\vep$, we have {\small $r (p) > C_{s - \ell_\vep + 1}$}.
Making use of the first four statements of Lemma \ref{lem52m}, we
have that {\small $y_\ell > \f{p + a_{\ell_\vep - 1} }{2} > p$} for
all $\ell \leq \ell_\vep - 1$ and {\small $z_{s - 1} < \f{p^\star +
b_{s - 1}}{2} < p^\star$} if $\vep_a$ is sufficiently small. By the
last statement of Lemma \ref{lem52m} and using Chernoff bound, we
have \bee \Pr \{ \bs{D}_{\ell} = 1 \} & \leq & \Pr \{
\wh{\bs{p}}_{\ell} \geq y_{\ell} \} + \Pr \{ \wh{\bs{p}}_{\ell} \leq
z_{s - 1} \}
 \leq  \Pr \li \{ \wh{\bs{p}}_{\ell} \geq \f{p + a_{\ell_\vep - 1}
}{2}  \ri \} + \Pr \li \{ \wh{\bs{p}}_{\ell} \leq \f{p^\star + b_{s - 1}}{2} \ri \}\\
& \leq & \exp \li ( - 2 n_\ell \li( \f{p - a_{\ell_\vep - 1} }{2}
\ri)^2  \ri ) + \exp \li ( - 2 n_\ell \li (p - \f{p^\star + b_{s -
1}}{2} \ri )^2 \ri ) \eee for all $\ell \leq \ell_\vep - 1$ provided
that $\vep_a > 0$ is small enough.  As a result of the definition of $\ell_\vep$,
we have that $a_{\ell_\vep - 1}$ is greater than $p$ and is independent of $\vep_a > 0$.
 Hence, it follows from Lemma \ref{lem31a} that {\small $\lim_{\vep_a \to 0} \sum_{\ell = 1}^{\ell_\vep - 1} n_\ell \Pr \{
\bs{D}_\ell = 1 \}  = 0$}.

In a similar manner, by the definition of $\ell_\vep$, we have
{\small $r (p) < C_{s - \ell_\vep - 1}$}.  Making use of the first
four statements of Lemma \ref{lem52m}, we have that {\small $y_\ell
< \f{ p + a_{\ell_\vep + 1} }{2} < p$}
 for $\ell_\vep + 1 \leq \ell < s$ if $\vep_a$ is sufficiently small.
 By the last statement of Lemma \ref{lem52m} and using Chernoff bound, we have
{\small \[ \Pr \{ \bs{D}_{\ell} = 0 \} = \Pr \{  z_{\ell} <
\wh{\bs{p}}_{\ell} < y_{\ell} \} \leq \Pr \{ \wh{\bs{p}}_{\ell} <
y_{\ell} \} \leq \Pr \li \{ \wh{\bs{p}}_{\ell} < \f{ p +
a_{\ell_\vep + 1} }{2} \ri \} \leq \exp \li ( - 2 n_\ell \li (  \f{
p - a_{\ell_\vep + 1} }{2} \ri )^2 \ri )
\]}
for $\ell_\vep + 1 \leq \ell < s$ if $\vep_a
> 0$ is small enough.  Clearly, $\Pr \{ \bs{D}_s = 0 \} = 0$.  As a consequence of the definition of $\ell_\vep$,
we have that $a_{\ell_\vep + 1}$ is smaller than $p$ and is independent of $\vep_a > 0$.
 Hence, it follows from Lemma \ref{lem31a} that  {\small $\lim_{\vep_a \to 0} \sum_{\ell = \ell_\vep + 1}^s n_\ell \Pr \{
\bs{D}_\ell = 0 \} = 0$}.  This proves that (\ref{goodone}) holds for $p \in (p^\star, 1)$.

\bsk

Third, we shall show $\lim_{\vep_a \to 0} n_{\ell_\vep} \Pr \{
\bs{D}_{\ell_\vep} = 0 \} = 0$ for $p \in (0, 1)$ such that $C_{j_p}
> r(p)$.

For $p \in (0, p^\star)$ such that $C_{j_p}
> r(p)$, we have {\small $r (p) < C_{s - \ell_\vep}$} because of the
definition of $\ell_\vep$. Making use of the first four statements
of Lemma \ref{lem52m}, we have that {\small $z_{\ell_\vep} > \f{p +
b_{\ell_\vep} }{2} > p$} if $\vep_a > 0$ is small enough. By the
last statement of Lemma \ref{lem52m} and using Chernoff bound, we
have {\small \[ \Pr \{ \bs{D}_{\ell_\vep} = 0 \} = \Pr \{
z_{\ell_\vep} < \wh{\bs{p}}_{\ell_\vep} < y_{\ell_\vep} \} \leq \Pr
\{
 \wh{\bs{p}}_{\ell_\vep} > z_{\ell_\vep} \} \leq \Pr \li \{
\wh{\bs{p}}_{\ell_\vep} > \f{p + b_{\ell_\vep} }{2} \ri \}  \leq
\exp \li ( - 2 n_{\ell_\vep} \li ( \f{p - b_{\ell_\vep} }{2}  \ri
)^2 \ri ).
\]}
Since $b_{\ell_\vep}$ is greater than $p$ and is independent of $\vep_a > 0$ due to the definition of $\ell_\vep$,
it follows that  $\lim_{\vep_a \to 0} n_{\ell_\vep} \Pr \{ \bs{D}_{\ell_\vep} = 0
\} = 0$.

For $p \in (p^\star, 1)$ such that $C_{j_p}
> r(p)$, we have {\small $r (p) < C_{s - \ell_\vep}$} as a result of the
definition of $\ell_\vep$. Making use of the first four statements
of Lemma \ref{lem52m}, we have that {\small $y_{\ell_\vep} < \f{p +
a_{\ell_\vep}}{2} < p$} if $\vep_a > 0$ is sufficiently small. By
the last statement of Lemma \ref{lem52m} and using Chernoff bound,
we have {\small \[ \Pr \{ \bs{D}_{\ell_\vep} = 0 \} = \Pr \{
z_{\ell_\vep} < \wh{\bs{p}}_{\ell_\vep} < y_{\ell_\vep} \} \leq \Pr
\{ \wh{\bs{p}}_{\ell_\vep} < y_{\ell_\vep} \} \leq \Pr \li \{
\wh{\bs{p}}_{\ell_\vep}  < \f{p + a_{\ell_\vep}}{2} \ri \}  \leq
\exp \li ( - 2 n_{\ell_\vep} \li ( \f{p - a_{\ell_\vep}}{2} \ri )^2
\ri ).
\]}
Since $a_{\ell_\vep}$ is smaller than $p$ and is independent of
$\vep_a > 0$ as a consequence of the definition of $\ell_\vep$, it
follows  that $\lim_{\vep_a \to 0} n_{\ell_\vep} \Pr \{
\bs{D}_{\ell_\vep} = 0 \} = 0$. This proves $\lim_{\vep_a \to 0}
n_{\ell_\vep} \Pr \{ \bs{D}_{\ell_\vep} = 0 \} = 0$ for $p \in (0,
1)$ such that $C_{j_p} > r(p)$. The proof of the lemma is thus
completed.

\epf

\bsk

The proof of Theorem \ref{Bino_mix_DDV_Asp} can be accomplished by
employing Lemma \ref{Defmix} and a similar argument as the proof of
Theorem \ref{Bino_DDV_Asp}.

\subsection{Proof of Theorem \ref{Bino_mix_Analysis} }  \la{App_Bino_mix_Analysis}

As a result of the definitions of $\ka_p$ and $r(p)$,
 we have that $\ka_p > 1$ if and only if $r(p)$ is not an integer.
 To prove Theorem \ref{Bino_mix_Analysis}, we need some preliminary results.

\beL
\la{lem54m}

$\lim_{\vep_a \to 0} \f{ n_{\ell_\vep}  } { \mcal{N}_{\mrm{m}} (p,
\vep_a, \vep_r) } = \ka_p, \; \lim_{\vep_a \to 0} \vep_a \sq{\f{
n_{\ell_\vep}}{p(1 - p)} } = d \sq{\ka_p}, \; \lim_{\vep_r \to 0}
\vep_r \sq{\f{ p n_{\ell_\vep}}{1 - p} }  = d \sq{\ka_p}$.

\eeL

\bpf

First, we shall consider $p \in (0, p^\star]$.  By the definition of
sample sizes, we have \be \la{defasn} \lim_{\vep_a \to 0} \f{  C_{ s
- \ell} \; \ln (\ze \de) } {n_\ell \mscr{M}_{\mrm{B}} (p^\star +
\vep_a, p^\star )} = 1 \ee for $1 \leq \ell < s$. It follows that
{\small \bee \lim_{\vep_a \to 0} \f{ n_{\ell_\vep} } {
\mcal{N}_{\mrm{m}} (p, \vep_a, \vep_r) } & = & \lim_{\vep_a \to 0}
\f{ \mscr{M}_{\mrm{B}} (p, p + \vep_a )  } { \ln (\ze \de)  } \times
\f{ C_{s - \ell_\vep} \ln (\ze \de) } {\mscr{M}_{\mrm{B}} (p^\star +
\vep_a, p^\star )} = \lim_{\vep_a \to 0} \f{ C_{s - \ell_\vep}
\mscr{M}_{\mrm{B}} (p, p + \vep_a
) } {\mscr{M}_{\mrm{B}} (p^\star + \vep_a, p^\star )} \\
& = & \lim_{\vep_a \to 0} \f{ C_{s - \ell_\vep} \li ( \vep_a^2 \sh
[2 p (p - 1)] + o (\vep_a^2) \ri ) } { \vep_a^2 \sh
[2 p^\star ( p^\star - 1)] + o (\vep_a^2)} \\
& = &  \f{ p^\star (1 - p^\star) }{ p (1 - p)} C_{s - \ell_\vep}  =
\f{ p^\star (1 - p^\star) }{ p (1 - p)} C_{j_p}  = \ka_p \eee} and
{\small \bee \lim_{\vep_a \to 0} \vep_a \sq{\f{ n_{\ell_\vep}}{p(1 -
p)} } & = & \lim_{\vep_a \to 0} \vep_a \sq{ \f{ 1 } { p(1 - p) } \f{
C_{s - \ell_\vep} \ln
(\ze \de) } {\mscr{M}_{\mrm{B}} ( p^\star + \vep_a, p^\star )} }\\
& = & \lim_{\vep_a \to 0} \vep_a \sq{ \f{ 1 } { p(1 - p) } \times
\f{ C_{s - \ell_\vep} \ln (\ze \de) } { \vep_a^2 \sh [2 p^\star
(p^\star - 1)] + o (\vep_a^2) }  } = d \sq{\f{ p^\star (1 - p^\star)
}{ p (1 - p)}  C_{s - \ell_\vep}  }
 = d \sq{\ka_p}.  \eee}

\bsk

Next,  we shall consider $p \in (p^\star, 1]$.  By virtue of (\ref{defasn}), we have
 \bee \lim_{\vep_r \to 0} \f{n_{\ell_\vep} }{\mscr{N}_{\mrm{m}} (p, \vep_a, \vep_r)} & = &
\lim_{\vep_r \to 0} \f{ \mscr{M}_{\mrm{B}} (p, \f{p}{1 + \vep_r} ) }
{ \ln (\ze \de)  } \times \f{ C_{s - \ell_\vep} \ln (\ze \de) }
{\mscr{M}_{\mrm{B}} ( p^\star + \vep_a,   p^\star  )}  =
\lim_{\vep_r \to 0} \f{ C_{s - \ell_\vep} \mscr{M}_{\mrm{B}}
(p, \f{p}{1 + \vep_r} ) } {\mscr{M}_{\mrm{B}} ( p^\star + \vep_a,   p^\star  )}\\
& = & \lim_{\vep_r \to 0} \f{ C_{s - \ell_\vep} \li ( \vep_r^2 p \sh
[ 2 (p - 1)] + o (\vep_r^2) \ri ) } {
\vep_a^2  \sh [ 2 p^\star (p^\star - 1)] + o (\vep_a^2)  } \\
& = & \f{ p (1 - p^\star) }{ p^\star (1 - p)} C_{s - \ell_\vep}  =
\f{ p (1 - p^\star) }{ p^\star (1 - p)} C_{j_p}  = \ka_p \eee and
\bee \lim_{\vep_r \to 0} \vep_r \sq{\f{ p n_{\ell_\vep}}{1 - p} }  &
= & \lim_{\vep_r \to 0} \vep_r \sq{ \f{ p } { 1 - p } \f{ C_{s -
\ell_\vep} \ln (\ze
\de) } {\mscr{M}_{\mrm{B}} ( p^\star + \vep_a,   p^\star  )}  }\\
& = & \lim_{\vep_r \to 0} \vep_r \sq{ \f{ p } { 1 - p } \times \f{
C_{s - \ell_\vep} \ln (\ze \de) } {  \vep_a^2 \sh [ 2 p^\star
(p^\star - 1)] + o (\vep_a^2) }  } = d \sq{ \f{ p (1 - p^\star) }{
p^\star (1 - p)} C_{s - \ell_\vep} } = d \sq{\ka_p}. \eee

\epf

\beL \la{limplemmix} Let $U$ and $V$ be independent Gaussian random
variables with zero means and unit variances.  Then, for $p \in (0,
1)$ such that $C_{j_p} = r(p)$ and $j_p \geq 1$,
\bee  &  &
\lim_{\vep \to 0} \Pr \{ \bs{l} = \ell_\vep \} = 1 - \lim_{\vep \to
0} \Pr \{ \bs{l} = \ell_\vep + 1 \} =  1 - \Phi \li (
\nu  d \ri ),\\
&  & \lim_{\vep \to 0} \li [ \Pr \{ | \wh{\bs{p}}_{\ell_\vep} - p |
\geq \vep_p, \; \bs{l}  = \ell_\vep \} + \Pr \{ |
\wh{\bs{p}}_{\ell_\vep + 1} - p | \geq \vep_p, \;
\bs{l} =  \ell_\vep + 1 \} \ri ]\\
&   & \qqu \qqu \qqu = \Pr \li \{ U \geq d \ri \} + \Pr \li \{ |U + \sq{\ro_p} V | \geq (1 + \ro_p) d, \; U < \nu d \ri \}, \eee where $\vep_p =
\max \{ \vep_a, \vep_r p \}$.

\eeL

\bpf

First, we shall consider $p \in [p^\star, 1)$.  Since $\ka_p = 1$,
by Statement (V) of Lemma \ref{lem52m},  we have
\[
\lim_{\vep_r \to 0} \f{ z_{\ell_\vep} - p}{ \sq{ p ( 1 - p) \sh
n_{\ell_\vep}} } = \lim_{\vep_r \to 0} \vep_r \sq{\f{ p
n_{\ell_\vep}}{1 - p} } \lim_{\vep_r \to 0} \f{ z_{\ell_\vep} -
p}{\vep_r p} = d  \lim_{\vep_r \to 0} \f{ z_{\ell_\vep} - p}{\vep_r
p} = \nu d.
\]
Note that {\small \[ \Pr \{ | \wh{\bs{p}}_{\ell_\vep} - p | \geq
\vep_r p, \; \wh{\bs{p}}_{\ell_\vep} \geq z_{\ell_\vep} \} = \Pr \li
\{ \f{ | \wh{\bs{p}}_{\ell_\vep} - p | }{ \sq{ p ( 1 - p) \sh
n_{\ell_\vep}} } \geq \vep_r \sq{\f{ p n_{\ell_\vep}}{1 - p} }, \;
\f{ \wh{\bs{p}}_{\ell_\vep} - p }{ \sq{ p ( 1 - p) \sh
n_{\ell_\vep}} } \geq \f{ z_{\ell_\vep} - p}{ \sq{ p ( 1 - p) \sh
n_{\ell_\vep}} } \ri \}.
\]}
Therefore, \bee &  & \Pr \{ | \wh{\bs{p}}_{\ell_\vep} - p | \geq
\vep_p, \; \bs{l}  = \ell_\vep \} + \Pr \{ | \wh{\bs{p}}_{\ell_\vep
+ 1} - p | \geq \vep_p, \; \bs{l} =  \ell_\vep + 1 \} \\  & \to &
\Pr \{  |U| \geq b, \; U \geq \nu d \} + \Pr \li \{  \li |
 U + \sq{\ro_p} V \ri | \geq (1 + \ro_p ) d, \; U < \nu d \ri \}\\
& = & \Pr \{  U \geq d \} + \Pr \li \{  \li |  U + \sq{\ro_p} V \ri | \geq (1 + \ro_p) d, \; U < \nu d \ri \} \eee for $p \in (p^\star, 1)$ such
that $C_{j_p} = r (p)$.

\bsk

Next, we shall consider $p \in (0, p^\star)$.  Since $\ka_p = 1$, by
Statement (V) of Lemma \ref{lem52m},  we have
\[
\lim_{\vep_a \to 0} \f{ z_{\ell_\vep} - p}{ \sq{ p ( 1 - p) \sh
n_{\ell_\vep}} } = \lim_{\vep_a \to 0} \vep_a \sq{\f{
n_{\ell_\vep}}{p(1 - p)} } \lim_{\vep_a \to 0} \f{ z_{\ell_\vep} -
p}{\vep_a} = d \lim_{\vep_a \to 0} \f{ z_{\ell_\vep} - p}{\vep_a} =
- \nu d.
\]
Note that {\small \[ \Pr \{ | \wh{\bs{p}}_{\ell_\vep} - p | \geq
\vep_a, \; \wh{\bs{p}}_{\ell_\vep} \leq z_{\ell_\vep} \} = \Pr \li
\{ \f{ | \wh{\bs{p}}_{\ell_\vep} - p | }{ \sq{ p ( 1 - p) \sh
n_{\ell_\vep}} } \geq \vep_a \sq{\f{ n_{\ell_\vep}}{p(1 - p)} }, \;
\f{ \wh{\bs{p}}_{\ell_\vep} - p }{ \sq{ p ( 1 - p) \sh
n_{\ell_\vep}} } \leq \f{ z_{\ell_\vep} - p}{ \sq{ p ( 1 - p) \sh
n_{\ell_\vep}} } \ri \}.
\]}
Therefore, $\Pr \{ \bs{D}_{\ell_\vep} = 1 \} \to \Pr \{ U \geq \nu d
\}$ and \bee &  & \Pr \{ | \wh{\bs{p}}_{\ell_\vep} - p | \geq
\vep_p, \; \bs{l}  = \ell_\vep \} + \Pr \{ | \wh{\bs{p}}_{\ell_\vep
+ 1} - p | \geq \vep_p, \; \bs{l} =  \ell_\vep + 1 \} \\
 & \to & \Pr \{ |U| \geq d, \; U \leq - \nu d \} + \Pr \li \{ \li |  U +
\sq{\ro_p} V
\ri | \geq (1 + \ro_p) d, \; U > - \nu d \ri \}\\
& = & \Pr \li \{ U \geq d \ri \} + \Pr \li \{ \li |  U + \sq{\ro_p} V \ri | \geq (1 + \ro_p) d  , \; U < \nu d \ri \} \eee for $p \in (0,
p^\star)$ such that $C_{j_p} = r (p)$.

\epf

Now, we shall first show that Statement (I) holds for $p \in (0,
p^\star]$ such that $C_{j_p} = r(p)$. For this purpose, we need to
show that {\small \be \la{doitm} 1 \leq \limsup_{\vep_a \to 0} \f{
\mbf{n} (\om)  } {  \mcal{N}_{\mrm{m}} (p, \vep_a, \vep_r)  } \leq 1
+ \ro_p \qqu \tx{ for any} \; \om \in \li \{ \lim_{\vep_a \to 0}
\wh{\bs{p}} = p \ri \}. \ee} To show {\small $\limsup_{\vep_a \to 0}
\f{ \mbf{n} (\om)  } { \mcal{N}_{\mrm{m}} (p, \vep_a, \vep_r)  }
\geq 1$}, note that $C_{s - \ell_\vep + 1} < r (p) = C_{s -
\ell_\vep} < C_{s - \ell_\vep - 1}$ as a direct consequence of the
definitions of $\ell_\vep$ and $j_p$. By the first four statements
of  Lemma \ref{lem52m}, we have $\lim_{\vep_a \to 0} z_{{\ell}} < p$
for all $\ell \leq {\ell_\vep} - 1$ with {\small $n_\ell \geq \f{
\ln (\ze \de) }{ \ln (1 - \vep_a) }$}. Noting that $\lim_{\vep_a \to
0} \wh{\bs{p}} (\om) = p$,  we have $\wh{\bs{p}} (\om) > z_\ell$ for
all $\ell \leq {\ell_\vep} - 1$ with {\small $n_\ell \geq \f{ \ln
(\ze \de) }{ \ln (1 - \vep_a) }$} and it follows from the definition
of the sampling scheme that
 $\mbf{n} (\om) \geq n_{\ell_\vep}$ if $\vep_a > 0$ is small enough.
 By Lemma \ref{lem54m} and noting that $\ka_p = 1$ if $C_{j_p} = r(p)$, we have
 {\small $\limsup_{\vep_a \to 0} \f{ \mbf{n} (\om)  } {  \mcal{N}_{\mrm{m}}  (p, \vep_a, \vep_r)  }
\geq \lim_{\vep_a \to 0}  \f{ n_{{\ell_\vep}}} { \mcal{N}_{\mrm{m}}
(p, \vep_a, \vep_r) } = \ka_p = 1$}.

To show {\small $\limsup_{\vep_a \to 0} \f{ \mbf{n} (\om)  } {
\mcal{N}_{\mrm{m}}  (p, \vep_a, \vep_r)  } \leq 1 + \ro_p$}, we
shall consider three cases: (i) ${\ell_\vep} = s$; (ii) ${\ell_\vep}
= s - 1$; (iii) ${\ell_\vep} < s - 1$. In the case of ${\ell_\vep} =
s$, it must be true that $\mbf{n} (\om) \leq n_s = n_{\ell_\vep}$.
Hence,  {\small $\limsup_{\vep_a \to 0} \f{ \mbf{n} (\om)  } {
\mcal{N}_{\mrm{m}} (p, \vep_a, \vep_r)  } \leq \lim_{\vep_a \to 0}
\f{ n_{\ell_\vep} } { \mcal{N}_{\mrm{m}}  (p, \vep_a, \vep_r) } =
\ka_p = 1 = 1 + \ro_p$}.  In the case of ${\ell_\vep} = s - 1$, it
must be true that $\mbf{n} (\om) \leq n_s = n_{\ell_\vep + 1}$.
Hence,  {\small $\limsup_{\vep_a \to 0} \f{ \mbf{n} (\om)  } {
\mcal{N}_{\mrm{m}} (p, \vep_a, \vep_r) } \leq \lim_{\vep_a \to 0}
\f{ n_{\ell_\vep + 1} } { \mcal{N}_{\mrm{m}} (p, \vep_a, \vep_r) } =
1 + \ro_p$}.  In the case of ${\ell_\vep} < s - 1$, it follows from
Lemma \ref{lem52m} that $\lim_{\vep_a \to 0} z_{{\ell_\vep} + 1} >
p$, which implies that
 $z_{{\ell_\vep} + 1} > p, \; \wh{\bs{p}} (\om) < z_{{\ell_\vep} + 1}$, and thus $\mbf{n} (\om) \leq
n_{{\ell_\vep} + 1}$ for small enough $\vep_a > 0$. Therefore,
{\small $\limsup_{\vep_a \to 0} \f{ \mbf{n} (\om)  } {
\mcal{N}_{\mrm{m}}  (p, \vep_a, \vep_r)  } \leq \lim_{\vep_a \to 0}
\f{ n_{{\ell_\vep} + 1}} { \mcal{N}_{\mrm{m}} (p, \vep_a, \vep_r) }
= \lim_{\vep_a \to 0}  \f{ n_{{\ell_\vep} + 1}} { n_{{\ell_\vep}}}
\times \lim_{\vep_a \to 0}  \f{ n_{{\ell_\vep} }} {
\mcal{N}_{\mrm{m}} (p, \vep_a, \vep_r) } = 1 + \ro_p$}.  This
establishes (\ref{doitm}) and it follows that
 $\{ 1 \leq  \limsup_{\vep_a \to 0} \f{ \mbf{n}  } {
\mcal{N}_{\mrm{m}}  (p, \vep_a, \vep_r)  } \leq 1 + \ro_p \}
\supseteq \li \{ \lim_{\vep_a \to 0} \wh{\bs{p}} = p \ri \}$.
According to the strong law of large numbers, we have $1 \geq \Pr \{
1 \leq \limsup_{\vep_a \to 0} \f{ \mbf{n}  } {  \mcal{N}_{\mrm{m}}
(p, \vep_a, \vep_r) } \leq 1 + \ro_p \} \geq \Pr \li \{ \lim_{\vep_a
\to 0} \wh{\bs{p}} = p \ri \} = 1$.  This proves that Statement (I)
holds for $p \in (0, p^\star ]$ such that $C_{j_p} = r(p)$.

\bsk

Next, we shall show that Statement (I) holds for $p \in (0,
p^\star]$ such that $C_{j_p} > r(p)$. Note that {\small $C_{s -
\ell_\vep + 1} < r (p) < C_{s - \ell_\vep}$} as a direct consequence
of the definitions of $\ell_\vep$ and $j_p$.  By the first four
statements of Lemma \ref{lem52m}, we have $\lim_{\vep_a \to 0}
z_{{\ell_\vep} - 1} < p$ and thus $z_{{\ell}} < p$ for all $\ell
\leq {\ell_\vep} - 1$ with {\small $n_\ell \geq \f{ \ln (\ze \de) }{
\ln (1 - \vep_a) }$} provided that $\vep_a > 0$ is sufficiently
small. Therefore, for any $\om \in \li \{ \lim_{\vep_a \to 0}
\wh{\bs{p}} = p \ri \}$, we have $ z_\ell < \wh{\bs{p}} (\om) <
y_\ell$ for all $\ell \leq {\ell_\vep} - 1$ with {\small $n_\ell
\geq \f{ \ln (\ze \de) }{ \ln (1 - \vep_a) }$} and consequently,
$\mbf{n} (\om) \geq n_{{\ell_\vep}}$ provided that $\vep_a > 0$ is
sufficiently small. On the other hand, we can show that $\mbf{n}
(\om) \leq n_{{\ell_\vep}}$ if $\vep_a > 0$ is small enough by
investigating two cases.  In the case of $\ell_\vep = s$, it is
trivially true that $\mbf{n} (\om) \leq n_{{\ell_\vep}}$. In the
case of $\ell_\vep < s$,  we have $p < \lim_{\vep_a \to 0}
z_{{\ell_\vep}}$ and thus $p < z_{{\ell_\vep}}$ provided that
$\vep_a > 0$ is sufficiently small. Therefore, for any $\om \in \li
\{ \lim_{\vep_a \to 0} \wh{\bs{p}} = p \ri \}$, we have $\wh{\bs{p}}
(\om) < z_{{\ell_\vep}}$ and consequently, $\mbf{n} (\om) \leq
n_{{\ell_\vep}}$ provided that $\vep_a > 0$ is sufficiently small.
So, we have established that $\mbf{n} (\om) = n_{{\ell_\vep}}$ if
$\vep_a > 0$ is sufficiently small.  Applying Lemma \ref{lem54m}, we
have {\small $\lim_{\vep_a \to 0} \f{ \mbf{n} (\om) } {
\mcal{N}_{\mrm{m}}  (p, \vep_a, \vep_r) } = \lim_{\vep_a \to 0} \f{
n_{{\ell_\vep}}}
 { \mcal{N}_{\mrm{m}}  (p, \vep_a, \vep_r) } = \ka_p$}, which implies that
{\small $\{ \lim_{\vep_a \to 0} \f{ \mbf{n} } {  \mcal{N}_{\mrm{a}}
(p, \vep)  } = \ka_p \} \supseteq \li \{ \lim_{\vep_a \to 0}
\wh{\bs{p}} = p \ri \}$}.  It follows from the strong law of large
numbers that {\small $1 \geq \Pr \{ \lim_{\vep_a \to 0} \f{ \mbf{n}
} { \mcal{N}_{\mrm{m}}  (p, \vep_a, \vep_r)  } = \ka_p \} \geq \Pr
\{ \lim_{\vep_a \to 0} \wh{\bs{p}} = p \} = 1$} and thus {\small
$\Pr \{ \lim_{\vep_a \to 0} \f{ \mbf{n} } {  \mcal{N}_{\mrm{m}}  (p,
\vep_a, \vep_r)  } = \ka_p \} = 1$}.  Since $1 \leq \ka_p \leq 1 +
\ro_p$, it is of course true that $\Pr \{ 1 \leq \limsup_{\vep_a \to
0} \f{ \mbf{n}  } {  \mcal{N}_{\mrm{m}} (p, \vep_a, \vep_r) } \leq 1
+ \ro_p \} = 1$.  This proves that Statement (I) holds true for $p
\in (0, p^\star]$ such that $C_{j_p} > r(p)$. Thus, we have shown
that Statement (I) holds true for $p \in (0, p^\star]$.

\bsk

In a similar manner, we can show that Statement (I) is true for $p
\in (p^\star, 1)$. This concludes the proof for Statement (I) of the
theorem.

\bsk

To show Statements (II) and (III),  we can employ Lemmas
\ref{Defmix}, \ref{lem54m} and mimic the corresponding arguments for
Theorem \ref{Bino_Asp_Analysis} by identifying $\vep_a$ and $\vep_r
p$ as $\vep$ for the cases of $p \leq p^\star$  and $p > p^\star$
respectively in the course of proof. Specially, in order to prove
Statement (III), we need to make use of the following observation:
\[
\Pr \{ | \wh{\bs{p}} - p | \geq \vep_a, \; | \wh{\bs{p}} - p | \geq \vep_r p \} = \bec
\Pr \{ | \wh{\bs{p}} - p | \geq \vep_a \} & \tx{for} \; p \in (0, p^\star],\\
\Pr \{ | \wh{\bs{p}} - p | \geq \vep_r p \} & \tx{for} \; p \in (p^\star, 1)
\eec
\]
\[ \Pr \{ | \wh{\bs{p}}_{\ell} - p | \geq \vep_a \} = \Pr \li \{  |U_\ell|
\geq \vep_a \sq{ \f{n_{\ell_\vep} } {p (1 - p)} } \ri \}, \qqu
\Pr \{ | \wh{\bs{p}}_{\ell} - p | \geq \vep_r p \} = \Pr \li \{  |U_\ell|
\geq \vep_r \sq{ \f{p n_{\ell}}{1 - p}  } \ri \}
\]
where, according to the central limit theorem,  $U_\ell = \f{ |
\wh{\bs{p}}_{\ell} - p | }{ \sq{ p (1 - p) \sh n_{\ell}} }$ converges in distribution
to a Gaussian random variable of zero mean and unit variance as $\vep_a \to 0$.

\subsection{Proof of Theorem \ref{multibound}} \la{multibound_app}

For simplicity of notations, let the complementary event of $\{ \mcal{L}_\ell (\wh{p}_\ell) \leq p_\ell  \leq \mcal{U}_\ell (\wh{p}_\ell), \; 1
\leq \ell \leq \ka \}$ be denoted by $\bs{\mscr{E}}$.  Define \bee &  &  E_i = \{ \mcal{L}_i (\wh{p}_i) > p_i,  \; \mcal{L}_\ell (\wh{p}_\ell)
\leq
p_\ell \leq \mcal{U}_\ell (\wh{p}_\ell), \; 1 \leq \ell < i \},\\
&  & \fra{E}_i = \{ \mcal{U}_i (\wh{p}_i) < p_i,  \; \mcal{L}_\ell (\wh{p}_\ell) \leq p_\ell  \leq \mcal{U}_\ell (\wh{p}_\ell), \; 1 \leq \ell <
i \} \eee and $\mcal{E}_i = E_i \cup \fra{E}_i$ for $i = 1, \cd, \ka$.  By an induction method,  we can show that \be \la{VVIPP}
 \bs{\mscr{E}} =
\cup_{i = 1}^\ka \mcal{E}_i. \ee   By the monotonicity of the confidence limits with respect to $\wh{p}_\ell$, we have \bel & & E_i \subseteq
\{ \wh{p}_i
> \mcal{L}_i^{-1} (\udl{p}_i),  \; \mcal{U}_\ell^{-1} (\udl{p}_\ell) \leq \wh{p}_\ell \leq \mcal{L}_\ell^{-1} (\ovl{p}_\ell), \;
1 \leq \ell < i \}, \la{crt188}\\
&  & E_i \supseteq \{ \wh{p}_i > \mcal{L}_i^{-1} (\ovl{p}_i),  \; \mcal{U}_\ell^{-1} (\ovl{p}_\ell) \leq \wh{p}_\ell \leq \mcal{L}_\ell^{-1}
(\udl{p}_\ell), \; 1 \leq \ell < i  \}, \la{crt288}\\
&  & \fra{E}_i \subseteq  \{ \wh{p}_i < \mcal{U}_i^{-1} (\ovl{p}_i),  \; \mcal{U}_\ell^{-1} (\udl{p}_\ell) \leq \wh{p}_\ell \leq
\mcal{L}_\ell^{-1} (\ovl{p}_\ell), \; 1 \leq \ell < i  \}, \la{crt1889}\\
&  & \fra{E}_i \supseteq \{ \wh{p}_i < \mcal{U}_i^{-1} (\udl{p}_i),  \; \mcal{U}_\ell^{-1} (\ovl{p}_\ell) \leq \wh{p}_\ell \leq
\mcal{L}_\ell^{-1} (\udl{p}_\ell), \; 1 \leq \ell < i \} \la{crt18899} \eel for all $\bs{p} \in \Se \cap \mcal{Q}$ and $i = 1, \cd, \ka$. Define
random variables
\[
Y_{i, \ell} = X_\ell, \qqu \ell = 1, \cd, \nu_i - 1; \qqu Y_{i, \nu_i} = n - \sum_{\ell = 1}^{\nu_i - 1} X_\ell
\]
for $i = 1, \cd, \ka$. Then, $(Y_{i, 1}, \cd, Y_{i, \nu_i})$ follows a multinomial distribution  with parameters $n$ and $(p_1, \cd, p_{\nu_i -
1})$. Define event \bee &  & \mscr{G}  =  \{ n \mrm{T}_{\mrm{lb}} (\udl{p}_\ell, n, \eta) \leq X_\ell \leq n \mrm{T}_{\mrm{ub}} (\ovl{p}_\ell,
n, \eta), \; \ell
= 1, \cd, \ka \}\\
&  & \qqu \cap  \{  n \mrm{T}_{\mrm{lb}} (\udl{\se}_{i, \nu_i}, n, \eta) \leq Y_{i, \nu_i}  \leq n \mrm{T}_{\mrm{ub}} (\ovl{\se}_{i, \nu_i}, n,
\eta), \; i = 1, \cd, \ka - 1 \}. \eee  Note that both $\mrm{T}_{\mrm{lb}} (\se, n, \eta)$ and $\mrm{T}_{\mrm{ub}} (\se, n, \eta)$ are
non-decreasing with respect to $\se \in [0, 1]$.  By (\ref{crt188}) and the definition of $\bs{\mcal{A}}_i, \; \bs{\mcal{B}}_i, \; i = 1, \cd,
\ka$, \be \la{OK1}
 E_i \cap \mscr{G} \subseteq \{ A_{i, \ell} \leq Y_{i, \ell} \leq B_{i, \ell}, \; \ell = 1, \cd, \nu_i \} \ee for all $\bs{p} \in \Se \cap \mcal{Q}$ and
$i = 1, \cd, \ka$.  By (\ref{crt1889}) and the definition of $\bs{\fra{A}}_i, \; \bs{\fra{B}}_i, \; i = 1, \cd, \ka$, \be \la{OK2}
 \fra{E}_i
\cap \mscr{G} \subseteq \{ \fra{A}_{i, \ell} \leq Y_{i, \ell} \leq \fra{B}_{i, \ell}, \; \ell = 1, \cd, \nu_i \} \ee for all $\bs{p} \in \Se
\cap \mcal{Q}$ and $i = 1, \cd, \ka$.  By Bonferroni's inequality and Theorem 3 of \cite{Chen1}, \be \la{OK3}
 \Pr \{ \bs{\mscr{E}} \mid \bs{p} \}
\leq (2 \ka - 1) \eta + \Pr \{ \bs{\mscr{E}} \cap \mscr{G} \mid \bs{p} \} \ee for all $\bs{p} \in \Se \cap \mcal{Q}$.  Making use of
(\ref{VVIPP}), (\ref{OK1}), (\ref{OK2}) and (\ref{OK3}), we have \bee \Pr \{ \bs{\mscr{E}} \mid \bs{p} \} & \leq & (2 \ka - 1) \eta + \sum_{i =
1}^\ka [ \Pr \{ E_i \cap \mscr{G} \mid \bs{p} \} + \Pr \{ \fra{E}_i \cap \mscr{G} \mid \bs{p} \} ]\\
& \leq &  (2 \ka - 1) \eta + \sum_{i = 1}^\ka  \Pr \{  A_{i, \ell} \leq Y_{i, \ell} \leq B_{i, \ell}, \; \ell = 1, \cd, \nu_i \mid \bs{p}  \}\\
&   & + \sum_{i = 1}^\ka  \Pr \{  \fra{A}_{i, \ell} \leq Y_{i, \ell} \leq \fra{B}_{i, \ell}, \; \ell = 1, \cd, \nu_i \mid \bs{p}  \}\\
& \leq & (2 \ka - 1) \eta + \sum_{i = 1}^\ka [S ( \bs{\mcal{A}}_i, \bs{\mcal{B}}_i, \ovl{\bs{\se}}_i, n ) + S ( \bs{\fra{A}}_i, \bs{\fra{B}}_i,
\ovl{\bs{\se}}_i, n ) ] \eee for all $\bs{p} \in \Se \cap \mcal{Q}$.  On the other side, by (\ref{crt288}) and the definition of
$\bs{\mcal{C}}_i, \; \bs{\mcal{D}}_i, \; i = 1, \cd, \ka$, \be \la{OK18}
 E_i \supseteq \{ C_{i, \ell} \leq Y_{i, \ell} \leq D_{i, \ell}, \; \ell = 1, \cd, \nu_i \} \ee for all $\bs{p} \in \Se \cap \mcal{Q}$ and
$i = 1, \cd, \ka$.  By (\ref{crt18899}) and the definition of $\bs{\fra{C}}_i, \; \bs{\fra{D}}_i, \; i = 1, \cd, \ka$, \be \la{OK28}
 \fra{E}_i \supseteq \{ \fra{C}_{i, \ell} \leq Y_{i, \ell} \leq \fra{D}_{i, \ell}, \; \ell = 1, \cd, \nu_i \} \ee for all $\bs{p} \in \Se \cap
\mcal{Q}$ and $i = 1, \cd, \ka$.  Applying (\ref{VVIPP}), (\ref{OK18}) and (\ref{OK28}), we have \bee \Pr \{ \bs{\mscr{E}} \mid \bs{p} \} & =  &
\sum_{i = 1}^\ka [ \Pr \{ E_i  \mid \bs{p} \} + \Pr \{ \fra{E}_i  \mid \bs{p} \} ]\\
& \geq & \sum_{i = 1}^\ka  \Pr \{  C_{i, \ell} \leq Y_{i, \ell} \leq D_{i, \ell}, \; \ell = 1, \cd, \nu_i \mid \bs{p}  \}\\
&  & + \sum_{i = 1}^\ka  \Pr \{  \fra{C}_{i, \ell} \leq Y_{i, \ell} \leq \fra{D}_{i, \ell}, \; \ell = 1, \cd, \nu_i \mid \bs{p}  \}\\
& \geq & \sum_{i = 1}^\ka [S ( \bs{\mcal{C}}_i, \bs{\mcal{D}}_i, \udl{\bs{\se}}_i, n ) + S ( \bs{\fra{C}}_i, \bs{\fra{D}}_i, \udl{\bs{\se}}_i, n
) ] \eee for all $\bs{p} \in \Se \cap \mcal{Q}$.  This completes the proof of the theorem.

\subsection{Proof of Theorem \ref{Bounded_Mean_abs_Hoeffding} } \la{App_Bounded_Mean_abs_Hoeffding}

We need some preliminary results.

\beL \la{absal1} Let $\ovl{X}_n = \f{\sum_{i=1}^n X_i}{n}$, where
$X_1, \; \cd,\; X_n$ are i.i.d.  random variables such that $0 \leq
X_i \leq 1$ and $\bb{E}[ X_i ] = \mu \in (0, 1)$ for $i = 1, \; \cd,
n$. Then,  $\Pr \li \{ \ovl{X}_n \geq \mu, \; \mscr{M}_{\mrm{B}} \li
( \ovl{X}_n ,\mu \ri ) \leq \f{\ln \al}{n} \ri \} \leq \al$ for any
$\al > 0$. \eeL

\bpf  For simplicity of notations, let $F_{\ovl{X}_n} (z) = \Pr \li
\{ \ovl{X}_n \leq z \ri \}$.  By Lemma \ref{Hoe_Mas}, we have that
$\{ \ovl{X}_n \geq \mu \} = \{ \ovl{X}_n \geq \mu, \; F_{\ovl{X}_n}
(\ovl{X}_n) \leq  \exp \li (n \mscr{M}_{\mrm{B}} \li (\ovl{X}_n, \mu
\ri ) \ri ) \}$.  Therefore, {\small \bee \li \{ \ovl{X}_n \geq \mu,
\; \mscr{M}_{\mrm{B}} \li ( \ovl{X}_n ,\mu \ri ) \leq \f{\ln \al}{n}
\ri \}   & = &  \li \{ \ovl{X}_n \geq \mu, \; \mscr{M}_{\mrm{B}} \li
( \ovl{X}_n ,\mu \ri ) \leq \f{\ln \al}{n}, \; F_{\ovl{X}_n}
(\ovl{X}_n) \leq  \exp \li (n \mscr{M}_{\mrm{B}} \li (\ovl{X}_n, \mu
\ri ) \ri ) \ri \}\\
& \subseteq &  \{ F_{\ovl{X}_n} (\ovl{X}_n) \leq \al \} \eee} and
thus Lemma \ref{absal1} follows from Lemma \ref{ProbTrans}.

\epf

\beL  \la{absal2} Let $\ovl{X}_n = \f{\sum_{i=1}^n X_i}{n}$,  where
$X_1, \; \cd,\; X_n$ are i.i.d. random variables such that $0 \leq
X_i \leq 1$ and $\bb{E}[ X_i ] = \mu \in (0, 1)$ for $i = 1, \; \cd,
n$. Then, $\Pr \li \{ \ovl{X}_n \leq \mu, \; \mscr{M}_{\mrm{B}} \li
( \ovl{X}_n,\mu \ri ) \leq \f{\ln \al}{n} \ri \} \leq \al$ for any
$\al > 0$. \eeL

\bpf For simplicity of notations, let $G_{\ovl{X}_n} (z) = \Pr \li
\{ \ovl{X}_n \geq z \ri \}$. By Lemma \ref{Hoe_Mas}, we have that
$\{ \ovl{X}_n \leq \mu \} = \{ \ovl{X}_n \leq \mu, \; G_{\ovl{X}_n}
(\ovl{X}_n) \leq \exp \li (n \mscr{M}_{\mrm{B}} \li (\ovl{X}_n, \mu
\ri ) \ri ) \}$.  Therefore, {\small \bee \li \{ \ovl{X}_n \leq \mu,
\; \mscr{M}_{\mrm{B}} \li ( \ovl{X}_n ,\mu \ri ) \leq \f{\ln \al}{n}
\ri \}   & = &  \li \{ \ovl{X}_n \leq \mu, \; \mscr{M}_{\mrm{B}} \li
( \ovl{X}_n ,\mu \ri ) \leq \f{\ln \al}{n}, \; G_{\ovl{X}_n}
(\ovl{X}_n) \leq  \exp \li (n \mscr{M}_{\mrm{B}} \li (\ovl{X}_n, \mu
\ri ) \ri ) \ri \}\\
& \subseteq &  \{ G_{\ovl{X}_n} (\ovl{X}_n) \leq \al \} \eee} and
thus Lemma \ref{absal2} follows from Lemma \ref{ProbTrans}.

\epf

Now we are in a position to show Theorem
\ref{Bounded_Mean_abs_Hoeffding}.  By a similar method as that of
Lemma \ref{DS1}, we can show that $\{ \mscr{M}_{\mrm{B}} \li (
\f{1}{2} - \li |\f{1}{2} - \wh{\bs{\mu}}_s \ri | , \f{1}{2} - \li
|\f{1}{2} - \wh{\bs{\mu}}_s \ri | + \vep \ri ) \leq \f{\ln (
\f{\de}{2 s} )} { n_s } \}$ is a sure event.  By a similar method as
that of Lemma \ref{abs13}, we can show that $\{ \mscr{M}_{\mrm{B}}
\li ( \f{1}{2} - \li |\f{1}{2} - \wh{\bs{\mu}}_\ell \ri | , \f{1}{2}
- \li |\f{1}{2} - \wh{\bs{\mu}}_\ell \ri | + \vep \ri ) \leq \f{\ln
( \f{\de}{2 s} )} { n_s } \} \subseteq \{ \mscr{M}_{\mrm{B}} (
\wh{\bs{\mu}}_\ell, \wh{\bs{\mu}}_\ell + \vep ) \leq \f{\ln
(\f{\de}{2 s}) }{n_\ell}, \; \mscr{M}_{\mrm{B}} (
\wh{\bs{\mu}}_\ell, \wh{\bs{\mu}}_\ell - \vep  ) \leq \f{\ln
(\f{\de}{2 s}) }{n_\ell} \}$ for $\ell = 1, \cd, s$.  Making use of
these facts and Lemmas \ref{absal1} and \ref{absal2}, we have
{\small \bee \Pr \{ |\wh{\bs{\mu}} - \mu | \geq \vep \} & \leq &
\sum_{\ell = 1}^s  \Pr \li \{ \mu \geq \wh{\bs{\mu}}_\ell + \vep, \;
\mscr{M}_{\mrm{B}} ( \wh{\bs{\mu}}_\ell,
\wh{\bs{\mu}}_\ell + \vep ) \leq \f{\ln (\f{\de}{2 s}) }{n_\ell} \ri \}\\
&  & +  \sum_{\ell = 1}^s  \Pr \li \{ \mu \leq \wh{\bs{\mu}}_\ell -
\vep, \; \mscr{M}_{\mrm{B}} ( \wh{\bs{\mu}}_\ell, \wh{\bs{\mu}}_\ell
- \vep ) \leq \f{\ln (\f{\de}{2 s}) }{n_\ell} \ri \}\\
& \leq & \sum_{\ell = 1}^s  \Pr \li \{ \mu \geq \wh{\bs{\mu}}_\ell,
\; \mscr{M}_{\mrm{B}} ( \wh{\bs{\mu}}_\ell, \mu ) \leq \f{\ln
(\f{\de}{2 s}) }{n_\ell} \ri \} +  \sum_{\ell = 1}^s  \Pr \li \{ \mu
\leq \wh{\bs{\mu}}_\ell, \; \mscr{M}_{\mrm{B}} ( \wh{\bs{\mu}}_\ell,
\mu ) \leq \f{\ln (\f{\de}{2 s}) }{n_\ell} \ri \} \leq \de, \eee}
from which Theorem \ref{Bounded_Mean_abs_Hoeffding} follows.

\subsection{Proof of Theorem \ref{Bounded_Mean_ABS_Massart} } \la{App_Bounded_Mean_ABS_Massart}

We need some preliminary results.

\beL \la{lemmas1} Let $\ovl{X}_n = \f{\sum_{i=1}^n X_i}{n}$, where
$X_1, \cd, X_n$ are i.i.d. random variables such that $0 \leq X_i
\leq 1$ and $\bb{E}[X_i] = \mu \in (0,1)$ for $i = 1, \cd, n$. Then,
{\small $\Pr \li \{ \ovl{X}_n \geq \mu, \; \mscr{M} \li ( \ovl{X}_n
,  \mu \ri ) \leq \f{\ln \al}{n} \ri \} \leq \al$} for any $\al >
0$. \eeL

\bpf

For simplicity of notations, let $F_{\ovl{X}_n} (z) = \Pr \li \{
\ovl{X}_n \leq z \ri \}$.  By Lemma \ref{Hoe_Mas}, we have that $\{
\ovl{X}_n \geq \mu \} = \{ \ovl{X}_n \geq \mu, \; F_{\ovl{X}_n}
(\ovl{X}_n) \leq  \exp \li (n \mscr{M} \li (\ovl{X}_n, \mu \ri ) \ri
) \}$. Therefore, {\small \bee \li \{ \ovl{X}_n \geq \mu, \;
\mscr{M} \li ( \ovl{X}_n ,\mu \ri ) \leq \f{\ln \al}{n} \ri \}   & =
&  \li \{ \ovl{X}_n \geq \mu, \; \mscr{M} \li ( \ovl{X}_n ,\mu \ri )
\leq \f{\ln \al}{n}, \; F_{\ovl{X}_n} (\ovl{X}_n) \leq \exp \li (n
\mscr{M} \li (\ovl{X}_n, \mu
\ri ) \ri ) \ri \}\\
& \subseteq &  \{ F_{\ovl{X}_n} (\ovl{X}_n) \leq \al \} \eee} and
thus  Lemma \ref{lemmas1} follows from Lemma \ref{ProbTrans}.

\epf

\beL \la{lemmas2} Let $\ovl{X}_n = \f{\sum_{i=1}^n X_i}{n}$, where
$X_1, \cd, X_n$ are i.i.d. random variables such that $0 \leq X_i
\leq 1$ and $\bb{E}[X_i] = \mu \in (0,1)$ for $i = 1, \cd, n$. Then,
{\small $\Pr \{ \ovl{X}_n \leq \mu, \;  \mscr{M} \li ( \ovl{X}_n ,
\mu \ri ) \leq \f{\ln \al}{n} \}  \leq \al$} for any $\al > 0$. \eeL

\bpf

For simplicity of notations, let $G_{\ovl{X}_n} (z) = \Pr \li \{
\ovl{X}_n \geq z \ri \}$. By Lemma \ref{Hoe_Mas}, we have that $\{
\ovl{X}_n \leq \mu \} = \{ \ovl{X}_n \leq \mu, \; G_{\ovl{X}_n}
(\ovl{X}_n) \leq \exp \li (n \mscr{M} \li (\ovl{X}_n, \mu \ri ) \ri
) \}$. Therefore, {\small \bee \li \{ \ovl{X}_n \leq \mu, \;
\mscr{M} \li ( \ovl{X}_n ,\mu \ri ) \leq \f{\ln \al}{n} \ri \}   & =
&  \li \{ \ovl{X}_n \leq \mu, \; \mscr{M} \li ( \ovl{X}_n ,\mu \ri )
\leq \f{\ln \al}{n}, \; G_{\ovl{X}_n} (\ovl{X}_n) \leq \exp \li (n
\mscr{M} \li (\ovl{X}_n, \mu
\ri ) \ri ) \ri \}\\
& \subseteq &  \{ G_{\ovl{X}_n} (\ovl{X}_n) \leq \al \} \eee} and
thus Lemma \ref{lemmas2} follows from Lemma \ref{ProbTrans}.

\epf

Now we are in a position to show Theorem
\ref{Bounded_Mean_ABS_Massart}.  By a similar method as that of
Lemma \ref{absD}, we can show that $ \{ ( | \wh{\bs{\mu}}_s -
\f{1}{2} | - \f{2\vep}{3} )^2 \geq \f{1}{4} + \f{  n_s \; \vep^2 }
{2 \ln (\f{\de}{2 s}) }  \}$ is a sure event.  By a similar method
as that of Lemma \ref{abs381}, we can show that {\small $\{ ( |
\wh{\bs{\mu}}_\ell - \f{1}{2} | - \f{2\vep}{3} )^2 \geq \f{1}{4} +
\f{ n_\ell \; \vep^2 } {2 \ln (\f{\de}{2 s}) }  \} \leu  \{
\mscr{M}_{\mrm{B}} \li ( \wh{\bs{\mu}}_\ell, \wh{\bs{\mu}}_\ell +
\vep \ri ) \leq \f{\ln (\f{\de}{2 s})}{n_\ell}, \;
\mscr{M}_{\mrm{B}} \li ( \wh{\bs{\mu}}_\ell, \wh{\bs{\mu}}_\ell -
\vep \ri ) \leq \f{\ln (\f{\de}{2 s})}{n_\ell} \} $} for $\ell = 1,
\cd, s$. Therefore,  by a variation of the argument for Theorem
\ref{Bounded_Mean_abs_Hoeffding} and using Lemmas \ref{lemmas1} and
\ref{lemmas2}, we have $\Pr \{ |\wh{\bs{\mu}} - \mu | \geq \vep \}
\leq  \de$,  from which Theorem \ref{Bounded_Mean_ABS_Massart}
follows.

\subsection{Proof of Theorem \ref{Bounded_Mean_mix_Hoeffding} } \la{App_Bounded_Mean_mix_Hoeffding}

By a similar method as that of Lemma \ref{mixDS1}, we can show that
{\small $\{ \mscr{M}_{\mrm{B}} ( \wh{\bs{p}}_s, \mscr{L} (
\wh{\bs{p}}_s ) ) \leq \f{\ln (\ze \de) } { n_s }, \;
\mscr{M}_{\mrm{B}} ( \wh{\bs{p}}_s, \mscr{U} ( \wh{\bs{p}}_s ) )
\leq \f{\ln (\ze \de) } { n_s } \}$} is a sure event.  By Lemmas
\ref{absal1} and \ref{absal2}, we have {\small \bee \Pr \{
|\wh{\bs{\mu}} - \mu | \geq \vep \} & \leq & \sum_{\ell = 1}^s \Pr
\li \{ \mu \geq \mscr{U} ( \wh{\bs{\mu}}_\ell ), \;
\mscr{M}_{\mrm{B}} ( \wh{\bs{\mu}}_\ell, \mscr{U} (
\wh{\bs{\mu}}_\ell ) ) \leq \f{\ln (\f{\de}{2 s}) }{n_\ell} \ri \}\\
&   & + \sum_{\ell = 1}^s  \Pr \li \{ \mu \leq \mscr{L} (
\wh{\bs{\mu}}_\ell ), \; \mscr{M}_{\mrm{B}} ( \wh{\bs{\mu}}_\ell,
\mscr{L} ( \wh{\bs{\mu}}_\ell ) ) \leq \f{\ln (\f{\de}{2 s}) }{n_\ell} \ri \}\\
& \leq & \sum_{\ell = 1}^s  \Pr \li \{ \mu \geq \wh{\bs{\mu}}_\ell,
\; \mscr{M}_{\mrm{B}} ( \wh{\bs{\mu}}_\ell, \mu ) \leq \f{\ln
(\f{\de}{2 s}) }{n_\ell} \ri \} +  \sum_{\ell = 1}^s  \Pr \li \{ \mu
\leq \wh{\bs{\mu}}_\ell, \; \mscr{M}_{\mrm{B}} ( \wh{\bs{\mu}}_\ell,
\mu ) \leq \f{\ln (\f{\de}{2 s}) }{n_\ell} \ri \} \leq  \de, \eee}
from which Theorem \ref{Bounded_Mean_mix_Hoeffding} follows.

\subsection{Proof of Theorem \ref{Bounded_Mean_mix_Massart} }
\la{App_Bounded_Mean_mix_Massart}

By a similar method as that of Lemma \ref{bb3800}, we can show that
$\{ \bs{D}_s = 1 \}$ is a sure event.  By a similar method as that
of Lemma \ref{good8a}, we can show that {\small $\{  \bs{D}_\ell = 1
\} \subseteq  \{  \mscr{M}_{\mrm{B}} ( \wh{\bs{\mu}}_\ell, \mscr{U}
( \wh{\bs{\mu}}_\ell )  ) \leq \f{\ln (\ze \de)}{n_\ell}, \;
\mscr{M}_{\mrm{B}} ( \wh{\bs{\mu}}_\ell, \mscr{L} (
\wh{\bs{\mu}}_\ell )  ) \leq \f{\ln (\ze \de)}{n_\ell} \}$} for
$\ell = 1, \cd, s$. Therefore,  by a variation of the argument for
Theorem \ref{Bounded_Mean_mix_Hoeffding} and using Lemmas
\ref{lemmas1} and \ref{lemmas2}, we can establish Theorem
\ref{Bounded_Mean_mix_Massart}.

\sect{Proofs of Theorems for Estimation of Poisson Parameters}

\subsection{Proof of Theorem \ref{Pos_abs_Chernoff} }  \la{App_Pos_abs_Chernoff}

First, we shall show statement (I).  Let $0 < \eta < 1$ and $r =
\inf_{\ell > 0} \f{n_{\ell + 1}}{n_\ell}$. By the assumption that $r
> 1$, we have that there exists a number $\ell^\prime > \max \{
\tau, \tau + \f{2}{r - 1} + \f{ \ln (\ze \de) } { \ln 2 } \}$ such
that $\f{n_{\ell + 1}}{n_\ell}
> \f{r + 1}{2}$ for any $\ell > \ell^\prime$.   Noting that {\small
\[ \f{ \f{ \ln (\ze \de_{\ell + 1}) } { n_{\ell+1} } } { \f{ \ln (
\ze \de_\ell ) } { n_\ell } } < \f{2}{r + 1} \times \f{ (\ell + 1 -
\tau) \ln 2  - \ln (\ze \de)  } { (\ell - \tau) \ln 2 - \ln (\ze
\de) } = \f{2}{r + 1} \times \li ( 1 + \f{1}{ \ell - \tau - \f{\ln
(\ze \de ) } { \ln 2 } } \ri ) < 1
\]}
for {\small $\ell > \ell^\prime$} and that {\small $\f{ \ln ( \ze
\de_\ell ) } { n_\ell } = \f{ \ln \li ( \ze \de 2^{\tau - \ell} \ri
) } { n_\ell  } \to 0 > \mscr{M}_{\mrm{P}}  ( \f{\lm}{\eta},
\f{\lm}{\eta} + \vep  )$} as $\ell \to \iy$, we have that there
exists an integer $\ka$ greater than {\small $\ell^\prime$} such
that {\small $\mscr{M}_{\mrm{P}}  ( \f{\lm}{\eta}, \f{\lm}{\eta} +
\vep ) < \f{ \ln (\ze \de_\ell) } { n_\ell }$} for all $\ell \geq
\ka$. For $\ell$ no less than such $\ka$, we claim that $z >
\f{\lm}{\eta}$ if {\small $\mscr{M}_{\mrm{P}} ( z, z + \vep )
> \f{\ln (\ze \de_\ell )}{n_\ell}$} and $z \in [0, \iy)$. To prove
this claim, suppose, to get a contradiction, that $z \leq
\f{\lm}{\eta}$. Then, since {\small $\mscr{M}_{\mrm{P}} ( z,  z +
\vep )$} is monotonically increasing with respect to $z > 0$, we
have {\small $\mscr{M}_{\mrm{P}} (z, z + \vep  ) \leq
\mscr{M}_{\mrm{P}} ( \f{\lm}{\eta}, \f{\lm}{\eta} + \vep ) < \f{\ln
(\ze \de_\ell)}{n_\ell}$}, which is a contradiction. Therefore, we
have shown the claim and it follows that {\small $ \{
\mscr{M}_{\mrm{P}} ( \f{K_\ell}{n_\ell}, \f{K_\ell}{ n_\ell } + \vep
) > \f{\ln (\ze \de_\ell) }{n_\ell}  \} \subseteq \{ \f{ K_\ell }
{n_\ell} > \f{\lm}{\eta} \}$} for $\ell \geq \ka$. So,  {\small
\[ \Pr \{ \bs{l}
> \ell \}  \leq  \Pr \li \{ \mscr{M}_{\mrm{P}}
\li ( \f{K_\ell}{n_\ell}, \f{K_\ell}{n_\ell } + \vep \ri ) > \f{\ln
(\ze \de_\ell) }{n_\ell}  \ri \} \leq  \Pr \li \{ \f{K_\ell}{n_\ell}
> \f{\lm}{\eta} \ri \} < \exp \li ( - c n_\ell  \ri ),
\]} where $c = - \mscr{M}_{\mrm{P}} ( \f{\lm}{\eta}, \lm )$ and
the last inequality is due to Chernoff bounds. Since $\Pr \{ \bs{l}
> \ell \}  < \exp ( - c n_\ell  )$ and
$n_\ell \to \iy$ as $\ell \to \iy$, we have $\Pr \{ \bs{l} < \iy \}
= 1$ or equivalently, $\Pr \{ \mbf{n} < \iy \} = 1$.  This completes
the proof of statement (I).

To show statement (II) of Theorem \ref{Pos_abs_Chernoff}, we can use
an argument similar to the proof of statement (II) of Theorem
\ref{Bino_Rev_noninverse_Chernoff}.

To show statement (III) of Theorem \ref{Pos_abs_Chernoff}, we can
use an argument similar to the proof of statement (III) of Theorem
\ref{Bino_Rev_noninverse_Chernoff}.

To show statement (IV) of Theorem \ref{Pos_abs_Chernoff}, we can use
an argument similar to the proof of statement (IV) of Theorem
\ref{Bino_Rev_noninverse_Chernoff} and make use of the observation
that {\small \bee \Pr \li \{ \li | \wh{\bs{\lm}} - \lm \ri | \geq
\vep \mid \lm \ri \} & = & \sum_{\ell = 1}^{\ell^\star} \Pr \li \{
\li | \wh{\bs{\lm}}_\ell - \lm \ri | \geq \vep, \; \bs{l} = \ell
\mid \lm \ri \} + \sum_{\ell = \ell^\star + 1}^\iy \Pr \li \{ \li |
\wh{\bs{\lm}}_\ell
- \lm \ri | \geq \vep, \; \bs{l} = \ell \mid \lm \ri \}\\
& \leq & \sum_{\ell = 1}^{\ell^\star} \Pr \li \{ \li |
\wh{\bs{\lm}}_\ell - \lm \ri | \geq \vep, \; \bs{l} = \ell \mid
\lm \ri \} + \eta\\
& \leq & \sum_{\ell = 1}^{\ell^\star} \Pr \li \{ \bs{l} = \ell \mid
\lm \ri \} + \eta \leq  \sum_{\ell = 1}^{\ell^\star} \Pr \li \{
\wh{\bs{\lm}}_\ell \leq z_\ell \mid \lm \ri \} + \eta \leq
\sum_{\ell = 1}^{\ell^\star} \exp ( n_\ell \mscr{M}_{\mrm{P}} (
z_\ell,  \lm) ) + \eta. \eee}

To show statement (V) of Theorem \ref{Pos_abs_Chernoff}, we can use
an argument similar to the proof of statement (V) of Theorem
\ref{Bino_Rev_noninverse_Chernoff}.

\subsection{Proof of Theorem  \ref{Pos_Asp_Analysis_Abs} } \la{App_Pos_Asp_Analysis_Abs}

Theorem \ref{Pos_Asp_Analysis_Abs}  can be established by using a
method similar to that of Theorem \ref{Bino_Asp_Analysis_Noninverse}
based on the following preliminary results.

\beL \la{Plemm21} Let $\vep > 0$.  Then, $\mscr{M}_{\mrm{P}} (z, z +
\vep)$ is monotonically increasing with respect to $z > 0$. \eeL

\bpf Note that $\mscr{M}_{\mrm{P}} (z, z + \vep) = - \vep + z \ln
\li ( \f{z + \vep}{z} \ri )$ and \[ \f{ \pa \mscr{M}_{\mrm{P}} (z, z
+ \vep) } { \pa z } = \ln \li ( \f{z + \vep}{z} \ri ) - \f{\vep}{z +
\vep} = - \ln \li ( 1 - \f{\vep}{z + \vep} \ri ) - \f{\vep}{z +
\vep} > 0, \qqu \fa z
> 0
\]
where the inequality follows from $\ln (1 - x) \leq - x, \; \fa x
\in [0, 1)$.

 \epf

\beL \la{Pos_abs_noninverse_lem} If $\vep$ is sufficiently small,
then the following statements hold true.

(I): For $\ell = 1, \cd, \tau$, there exists a unique number
$z_{\ell} \in [0, \iy)$ such that $n_{\ell} = \f{ \ln ( \ze \de_\ell
) } { \mscr{M}_{\mrm{P}} ( z_{\ell}, \; z_{\ell} + \vep ) }$.

(II): $z_{\ell}$ is monotonically increasing with respect to $\ell$
no greater than $\tau$.

(III): $\lim_{\vep \to 0} z_{\ell} =  \lm^* C_{\tau - \ell}$ for $1
\leq \ell \leq \tau$, where the limit is taken under the restriction
that $\ell - \tau$ is fixed with respect to $\vep$.

(IV): $\{ \bs{D}_{\ell} = 0 \} = \{ \wh{\bs{\lm}}_{\ell} > z_{\ell}
\}$ for $\ell = 1, \cd, \tau$.

\eeL

\bpf Lemma \ref{Pos_abs_noninverse_lem} can be shown by a similar
method as that of Lemma \ref{bino_noninverse_lem}. \epf

\beL \la{lPos_abs_noninv} Define $\ell_\vep = \tau  - j_\lm$, where
$j_\lm$ is the largest integer $j$ such that $C_{j} \geq
\f{\lm}{\lm^*}$. Then, \be \la{posabslemninv} \lim_{\vep \to 0}
\sum_{\ell = 1}^{\ell_\vep - 1} n_\ell \Pr \{ \bs{D}_\ell = 1 \} =
0, \qqu \lim_{\vep \to 0} \sum_{\ell = \ell_\vep + 1}^\tau n_\ell
\Pr \{ \bs{D}_{\ell} = 0 \} = 0 \ee for $\lm \in (0, \lm^*)$.
Moreover, $\lim_{\vep \to 0} n_{\ell_\vep} \Pr \{ \bs{D}_{\ell_\vep}
= 0 \} = 0$ for $\lm \in (0, \lm^*)$ such that $C_{j_\lm} >
\f{\lm}{\lm^*}$. \eeL

\bpf

For simplicity of notations, let $b_\ell = \lim_{\vep \to 0}
z_{\ell}$ for $1 \leq \ell \leq \tau$.

First, we shall show that (\ref{posabslemninv}) holds for $\lm \in
(0, \lm^*)$. By the definition of $\ell_\vep$,
 we have $b_{\ell_\vep - 1} = \lm^* C_{\tau - \ell_\vep +
 1} =  \lm^* C_{ j_\lm + 1} < \lm$.
 Making use of the first three statements of Lemma \ref{Pos_abs_noninverse_lem},
 we have that {\small $z_\ell < \f{\lm + b_{\ell_\vep -
1}}{2} < \lm$} for all $\ell \leq \ell_\vep - 1$ if $\vep$ is
sufficiently small. By the last statement of Lemma
\ref{Pos_abs_noninverse_lem}, we have {\small \bee \Pr \{
\bs{D}_{\ell} = 1 \}  =
 \Pr \{ \wh{\bs{\lm}}_{\ell} \leq z_{\ell} \} \leq  \Pr \li \{
\wh{\bs{\lm}}_{\ell} < \f{\lm + b_{\ell_\vep - 1}}{2} \ri \} \leq
\exp \li ( n_\ell \mscr{M}_{\mrm{P}} \li ( \f{\lm + b_{\ell_\vep -
1}}{2}, \lm \ri ) \ri ) \eee} for all $\ell \leq \ell_\vep - 1$
provided that $\vep > 0$ is small enough.   Since {\small $\f{\lm +
b_{\ell_\vep - 1}}{2}$} is independent of $\vep > 0$, we have
{\small $\lim_{\vep \to 0} \sum_{\ell = 1}^{\ell_\vep - 1} n_\ell
\Pr \{ \bs{D}_\ell = 1 \} = 0$} as a result of Lemma \ref{lem31a}.

Similarly, it can be seen from the definition of $\ell_\vep$ that
$b_{\ell_\vep + 1} = \lm^* C_{\tau - \ell_\vep -
 1} =  \lm^* C_{j_\lm - 1} > \lm$.  Making use of
the first three statements of Lemma \ref{Pos_abs_noninverse_lem}, we
have that {\small $z_\ell > \f{ \lm + b_{\ell_\vep + 1}}{2} > \lm$}
for $\ell_\vep + 1 \leq \ell \leq \tau$ if $\vep$ is sufficiently
small. By the last statement of Lemma \ref{Pos_abs_noninverse_lem},
we have {\small
\[ \Pr \{ \bs{D}_{\ell} = 0 \} = \Pr \{ \wh{\bs{\lm}}_{\ell}
> z_{\ell} \} \leq \Pr \li \{ \wh{\bs{\lm}}_{\ell} > \f{ \lm +
b_{\ell_\vep + 1}}{2} \ri \} \leq \exp \li ( n_\ell
\mscr{M}_{\mrm{P}} \li ( \f{ \lm + b_{\ell_\vep + 1}}{2}, \lm \ri )
\ri )
\]}
for $\ell_\vep + 1 \leq \ell \leq \tau$ provided that $\vep
> 0$ is small enough.  Therefore, we can apply Lemma \ref{lem31a} to conclude that {\small
$\lim_{\vep \to 0} \sum_{\ell = \ell_\vep + 1}^\tau n_\ell \Pr \{
\bs{D}_\ell = 0 \} = 0$}.

\bsk

Second, we shall show that $\lim_{\vep \to 0} n_{\ell_\vep} \Pr \{
\bs{D}_{\ell_\vep} = 0 \} = 0$  for $\lm \in (0, \lm^*)$ such that
$C_{j_\lm} > \f{\lm}{\lm^*}$. Clearly, $b_{\ell_\vep} = \lm^*
C_{\tau - \ell_\vep} = \lm^* C_{j_\lm} > \lm$. Making use of the
first three statements of Lemma \ref{Pos_abs_noninverse_lem}, we
have {\small $z_{\ell_\vep}
>\f{\lm + b_{\ell_\vep}}{2} > \lm$} if $\vep$ is sufficiently small.
By the last statement of Lemma \ref{Pos_abs_noninverse_lem}, we have
{\small
\[ \Pr \{ \bs{D}_{\ell_\vep} = 0 \} =  \Pr \{
 \wh{\bs{\lm}}_{\ell_\vep} > z_{\ell_\vep} \} \leq \Pr \li \{
\wh{\bs{\lm}}_{\ell_\vep} > \f{\lm + b_{\ell_\vep}}{2} \ri \}  \leq
\exp \li ( n_{\ell_\vep} \mscr{M}_{\mrm{P}} \li ( \f{\lm +
b_{\ell_\vep}}{2}, \lm \ri ) \ri )
\]} for small enough $\vep > 0$. It follows
that $\lim_{\vep \to 0} n_{\ell_\vep} \Pr \{ \bs{D}_{\ell_\vep} = 0
\} = 0$.

\epf

\beL

\la{Pos_abs_ineqinv}
 $\lim_{\vep \to 0} \sum_{\ell =  \tau +
1}^{\iy} n_\ell \Pr \{ \bs{l} = \ell \} = 0$ for any $\lm \in (0,
\lm^*)$.

\eeL

\bpf

Recalling that the sample sizes $n_1, n_2, \cd$ are chosen as the
ascending arrangement of all distinct elements of the set defined by
(\ref{defposabs}), we have that
\[
n_\ell =  \li \lc  \f{ C_{\tau - \ell} \; \ln (\ze \de) }{
\mscr{M}_{\mrm{P}} \li ( \lm^*, \lm^* + \vep \ri ) } \ri \rc, \qqu
\ell = 1, 2, \cd
\]
for small enough $\vep > 0$.   By the assumption that $\inf_{i \in
\bb{Z} } \f{C_{i - 1}} {C_{i} } = 1 + \udl{\ro} > 1$, we have that
\[
n_\ell  > (1 + \udl{\ro} )^{\ell - \tau - 1} \f{ \ln (\ze \de) } {
\mscr{M}_{\mrm{P}} \li ( \lm^*, \lm^* + \vep \ri ) }, \qqu \ell =
\tau + 1, \tau + 2, \cd
\]
for small enough $\vep > 0$. So, there exists a number $\vep^* > 0$
such that
\[
n_\ell \mscr{M}_{\mrm{P}} \li ( \lm^*, \lm^* + \vep \ri ) < (1 +
\udl{\ro} )^{\ell - \tau - 1} \ln (\ze \de), \qqu \ell = \tau + 1,
\tau + 2, \cd
\]
for any $\vep \in (0, \vep^*)$. Observing that there exist a
positive integer $\ka^*$ such that $(1 + \udl{\ro} )^{\ell - \tau -
1} \ln (\ze \de) < \ln (\ze \de) - (\ell - \tau) \ln 2 = \ln (\ze
\de_\ell)$ for any $\ell \geq \tau + \ka^*$, we have that there
exists a positive integer $\ka^*$ independent of $\vep$ such that
{\small $\mscr{M}_{\mrm{P}} ( \lm^*, \lm^* + \vep ) < \f{\ln (\ze
\de_\ell)}{n_\ell}$} for $\ell \geq \tau + \ka^*$ and $0 < \vep <
\vep^*$.  Note that $\mscr{M}_{\mrm{P}} (z, z + \vep )$ is
monotonically increasing
 with respect to $z \in (0, \iy)$ as asserted by Lemma \ref{Plemm21}.
   For $\ell \geq \tau + \ka^*$ and $0 <
\vep < \vep^*$, as a result of {\small $\f{\ln (\ze
\de_\ell)}{n_\ell} > \mscr{M}_{\mrm{P}} ( \lm^*, \lm^* + \vep ) $},
there exists a unique number $z_\ell \in [0, \iy)$ such that {\small
$\mscr{M}_{\mrm{P}} (z_\ell, z_\ell + \vep ) = \f{\ln (\ze
\de_\ell)}{n_\ell} > \mscr{M}_{\mrm{P}} ( \lm^*, \lm^* + \vep )$}.
Moreover, it must be true that $z_\ell > \lm^*$ for $\ell \geq \tau
+ \ka^*$ and $\vep \in (0, \vep^*)$. Therefore, for small enough
$\vep \in (0, \vep^*)$, we have {\small \bee \sum_{\ell = \tau +
1}^{\iy} n_\ell \Pr \{ \bs{l} = \ell \} & = & \sum_{\ell = \tau +
1}^{\tau + \ka^*} n_\ell \Pr \{ \bs{l} = \ell \} + \sum_{\ell = \tau
+ \ka^* + 1}^\iy n_\ell \Pr
\{ \bs{l} = \ell \}\\
& \leq & \sum_{\ell =  \tau + 1}^{\tau + \ka^*} n_\ell \Pr \{
\bs{D}_\tau  = 0 \} + \sum_{\ell =  \tau + \ka^* + 1}^\iy  n_\ell
\Pr \{ \bs{D}_{\ell - 1} = 0 \} \\
& = & \sum_{\ell =  \tau + 1}^{\tau + \ka^*} n_\ell \Pr \{
\bs{D}_\tau = 0 \} + \sum_{\ell =  \tau + \ka^*}^\iy  n_{\ell+1} \Pr
\{ \bs{D}_{\ell} = 0 \} \\
& \leq & k^* (1 + \ovl{\ro})^{k^*} n_\tau \Pr \{ \bs{D}_\tau = 0 \}
+ (1 + \ovl{\ro}) \sum_{\ell = \tau + \ka^*}^\iy  n_{\ell} \Pr \{
\bs{D}_{\ell} = 0 \} \\
& \leq  & k^* (1 + \ovl{\ro})^{k^*} n_\tau \Pr \{ \wh{\bs{\lm}}_\tau
> z_\tau \} + (1 + \ovl{\ro}) \sum_{\ell = \tau + \ka^*}^\iy
n_{\ell} \Pr \{ \wh{\bs{\lm}}_{\ell} > z_\ell \}\\
& \leq  & k^* (1 + \ovl{\ro})^{k^*} n_\tau \Pr \li \{
\wh{\bs{\lm}}_\tau > \f{\lm^* + \lm}{2} \ri \} + (1 + \ovl{\ro})
\sum_{\ell = \tau + \ka^*}^\iy n_{\ell} \Pr \{ \wh{\bs{\lm}}_{\ell} > \lm^* \}\\
& \leq  & k^* (1 + \ovl{\ro})^{k^*} n_\tau \exp \li ( n_\tau
\mscr{M}_{\mrm{P}} \li ( \f{\lm + \lm^*}{2}, \lm \ri ) \ri)\\
&  & + (1 + \ovl{\ro}) \sum_{\ell = \tau + \ka^*}^\iy n_{\ell} \exp(
n_\ell \mscr{M}_{\mrm{P}} ( \lm^*, \lm) ) \to 0 \eee} as $\vep \to
0$, where we have used  the assumption that $\sup_{i \in \bb{Z}}
\f{C_{i - 1}}{C_i} = 1 + \ovl{\ro} < \iy$. This completes the proof
of the lemma.  \epf

\subsection{Proof of Theorem \ref{Pos_rev_Chernoff} } \la{App_Pos_rev_Chernoff}

To show statement (I) of Theorem \ref{Pos_rev_Chernoff}, we can use
an argument similar to the proof of statement (I) of Theorem
\ref{Bino_Rev_noninverse_Chernoff}.

To show statement (II) of Theorem \ref{Pos_rev_Chernoff}, we can use
an argument similar to the proof of statement (II) of Theorem
\ref{Bino_Rev_noninverse_Chernoff}.

To show statement (III) of Theorem \ref{Pos_rev_Chernoff}, we can
use an argument similar to the proof of statement (III) of Theorem
\ref{Bino_Rev_noninverse_Chernoff}.

To show statement (IV) of Theorem \ref{Pos_rev_Chernoff}, we can use
an argument similar to the proof of statement (IV) of Theorem
\ref{Bino_Rev_noninverse_Chernoff} and make use of the observation
that {\small \bee \Pr \li \{ \li | \wh{\bs{\lm}} - \lm \ri | \geq
\vep \lm \mid \lm \ri \} & = & \sum_{\ell = 1}^{\ell^\star} \Pr \li
\{ \li | \wh{\bs{\lm}}_\ell - \lm \ri | \geq \vep \lm, \; \bs{l} =
\ell \mid \lm \ri \} + \sum_{\ell = \ell^\star + 1}^\iy \Pr \li \{
\li | \wh{\bs{\lm}}_\ell
- \lm \ri | \geq \vep \lm, \; \bs{l} = \ell \mid \lm \ri \}\\
& \leq & \sum_{\ell = 1}^{\ell^\star} \Pr \li \{ \li |
\wh{\bs{\lm}}_\ell - \lm \ri | \geq \vep \lm, \; \bs{l} = \ell \mid
\lm \ri \} + \eta\\
& \leq & \sum_{\ell = 1}^{\ell^\star} \Pr \li \{ \bs{l} = \ell \mid
\lm \ri \} + \eta \leq  \sum_{\ell = 1}^{\ell^\star} \Pr \li \{
\wh{\bs{\lm}}_\ell \geq z_\ell \mid \lm \ri \} + \eta \leq
\sum_{\ell = 1}^{\ell^\star} \exp ( n_\ell \mscr{M}_{\mrm{P}} (
z_\ell,  \lm) ) + \eta. \eee}

To show statement (V) of Theorem \ref{Pos_rev_Chernoff}, we can use
an argument similar to the proof of statement (V) of Theorem
\ref{Bino_Rev_noninverse_Chernoff} and make use of the observation
that \bee \Pr \li \{ \li | \wh{\bs{\lm}} - \lm \ri | \geq \vep \lm
\mid \lm \ri \} & \leq & \Pr \li \{ \li | \wh{\bs{\lm}} - \lm \ri |
\geq \vep \lm, \; \bs{l} = 1 \mid \lm \ri \} + \Pr \li \{ \li |
\wh{\bs{\lm}} - \lm \ri | \geq \vep
\lm, \; \bs{l} > 1 \mid \lm \ri \}\\
& \leq &  \Pr \li \{ \li | \wh{\bs{\lm}}_1 - \lm \ri | \geq \vep \lm
\mid \lm \ri \} + \Pr \li \{ \bs{l} > 1 \mid \lm \ri \}\\
& \leq &  \Pr \li \{ \li | \wh{\bs{\lm}}_1 - \lm \ri | \geq \vep \lm
\mid \lm \ri \} + \Pr \li \{ \wh{\bs{\lm}}_1 < z_1 \mid \lm \ri \}\\
& \leq & 2 \exp ( n_1 \mscr{M}_{\mrm{P}} ( (1 + \vep) \lm,  \lm) ) +
\exp ( n_1 \mscr{M}_{\mrm{P}} ( z_1,  \lm) ). \eee

\subsection{Proof of Theorem \ref{Pos_Asp_Analysis_Rev}}  \la{App_Pos_Asp_Analysis_Rev}

Theorem \ref{Pos_Asp_Analysis_Rev}  can be established by using a
method similar to that of Theorem \ref{Bino_Asp_Analysis_Noninverse}
based on the following preliminary results.

\beL \la{Pos_rev_noninverse_lem} If $\vep$ is sufficiently small,
then the following statements hold true.

(I): For $\ell = 1, \cd, \tau$, there exists a unique number
$z_{\ell} \in [0, \iy)$ such that $n_{\ell} = \f{ \ln ( \ze \de_\ell
) } { \mscr{M}_{\mrm{P}} ( z_{\ell}, \; \f{z_{\ell}}{1 + \vep} ) }$.

(II): $z_{\ell}$ is monotonically decreasing with respect to $\ell$
no greater than $\tau$.

(III): $\lim_{\vep \to 0} z_{\ell} =  \f{\lm^\prime}{C_{\tau -
\ell}}$ for $1 \leq \ell \leq \tau$, where the limit is taken under
the restriction that $\ell - \tau$ is fixed with respect to $\vep$.

(IV): $\{ \bs{D}_{\ell} = 0 \} = \{ \wh{\bs{\lm}}_{\ell} < z_{\ell}
\}$ for $\ell = 1, \cd, \tau$.

\eeL

\bpf Lemma \ref{Pos_rev_noninverse_lem} can be shown by a similar
method as that of Lemma \ref{bino_noninverse_lem}. \epf

\beL \la{lPos_rev_noninv} Define $\ell_\vep = \tau - j_\lm$, where
$j_\lm$ is the largest integer $j$ such that $C_{j} \geq
\f{\lm^\prime}{\lm}$.  Then, \be \la{posrevlemninv} \lim_{\vep \to
0} \sum_{\ell = 1}^{\ell_\vep - 1} n_\ell \Pr \{ \bs{D}_\ell = 1 \}
= 0, \qqu \lim_{\vep \to 0} \sum_{\ell = \ell_\vep + 1}^\tau n_\ell
\Pr \{ \bs{D}_{\ell} = 0 \} = 0 \ee for $\lm \in (\lm^\prime,
\lm^{\prime \prime})$. Moreover, $\lim_{\vep \to 0} n_{\ell_\vep}
\Pr \{ \bs{D}_{\ell_\vep} = 0 \} = 0$ for $\lm \in (\lm^\prime,
\lm^{\prime \prime})$ such that $C_{j_\lm} > \f{\lm^\prime}{\lm}$.
\eeL

\bpf

For simplicity of notations, let $b_\ell = \lim_{\vep \to 0}
z_{\ell}$ for $1 \leq \ell \leq \tau$.

First, we shall show that (\ref{posrevlemninv}) holds for $\lm \in
(\lm^\prime, \lm^{\prime \prime})$. By the definition of
$\ell_\vep$,
 we have $b_{\ell_\vep - 1} = \f{\lm^\prime}{ C_{\tau - \ell_\vep +
 1} } =  \f{\lm^\prime}{C_{j_\lm + 1} } > \lm$.
 Making use of the first three statements of Lemma \ref{Pos_rev_noninverse_lem},
 we have that {\small $z_\ell > \f{\lm + b_{\ell_\vep -
1}}{2} > \lm$} for all $\ell \leq \ell_\vep - 1$ if $\vep$ is
sufficiently small. By the last statement of Lemma
\ref{Pos_rev_noninverse_lem}, we have {\small \bee \Pr \{
\bs{D}_{\ell} = 1 \}  =
 \Pr \{ \wh{\bs{\lm}}_{\ell} \geq z_{\ell} \} \leq  \Pr \li \{
\wh{\bs{\lm}}_{\ell} > \f{\lm + b_{\ell_\vep - 1}}{2} \ri \} \leq
\exp \li ( n_\ell \mscr{M}_{\mrm{P}} \li ( \f{\lm + b_{\ell_\vep -
1}}{2}, \lm \ri ) \ri ) \eee} for all $\ell \leq \ell_\vep - 1$
provided that $\vep > 0$ is small enough.   Since {\small $\f{\lm +
b_{\ell_\vep - 1}}{2}$} is independent of $\vep > 0$, we have
{\small $\lim_{\vep \to 0} \sum_{\ell = 1}^{\ell_\vep - 1} n_\ell
\Pr \{ \bs{D}_\ell = 1 \} = 0$} as a result of Lemma \ref{lem31a}.

Similarly, it can be seen from the definition of $\ell_\vep$ that
$b_{\ell_\vep + 1} = \f{\lm^\prime}{ C_{\tau - \ell_\vep -
 1} }=  \f{ \lm^\prime }{C_{j_\lm - 1} } < \lm$.  Making use of
the first three statements of Lemma \ref{Pos_rev_noninverse_lem}, we
have that {\small $z_\ell < \f{ \lm + b_{\ell_\vep + 1}}{2} < \lm$}
for $\ell_\vep + 1 \leq \ell \leq \tau$ if $\vep$ is sufficiently
small. By the last statement of Lemma \ref{Pos_rev_noninverse_lem},
we have {\small
\[ \Pr \{ \bs{D}_{\ell} = 0 \} = \Pr \{ \wh{\bs{\lm}}_{\ell}
< z_{\ell} \} \leq \Pr \li \{ \wh{\bs{\lm}}_{\ell} < \f{ \lm +
b_{\ell_\vep + 1}}{2} \ri \} \leq \exp \li ( n_\ell
\mscr{M}_{\mrm{P}} \li ( \f{ \lm + b_{\ell_\vep + 1}}{2}, \lm \ri )
\ri )
\]}
for $\ell_\vep + 1 \leq \ell \leq \tau$ provided that $\vep
> 0$ is small enough.  Therefore, we can apply Lemma \ref{lem31a} to conclude that {\small
$\lim_{\vep \to 0} \sum_{\ell = \ell_\vep + 1}^\tau n_\ell \Pr \{
\bs{D}_\ell = 0 \} = 0$}.

\bsk

Second, we shall show that $\lim_{\vep \to 0} n_{\ell_\vep} \Pr \{
\bs{D}_{\ell_\vep} = 0 \} = 0$  for $\lm \in (\lm^\prime,
\lm^{\prime \prime})$ such that $C_{j_\lm}
> \f{\lm^\prime}{\lm}$. Clearly, $b_{\ell_\vep} = \f{\lm^\prime}{C_{\tau - \ell_\vep}} = \f{\lm^\prime}{ C_{j_\lm} } <
\lm$. Making use of the first three statements of Lemma
\ref{Pos_rev_noninverse_lem}, we have {\small $z_{\ell_\vep} <
\f{\lm + b_{\ell_\vep}}{2} < \lm$} if $\vep$ is sufficiently small.
By the last statement of Lemma \ref{Pos_rev_noninverse_lem}, we have
{\small
\[ \Pr \{ \bs{D}_{\ell_\vep} = 0 \} =  \Pr \{
 \wh{\bs{\lm}}_{\ell_\vep} < z_{\ell_\vep} \} \leq \Pr \li \{
\wh{\bs{\lm}}_{\ell_\vep} < \f{\lm + b_{\ell_\vep}}{2} \ri \}  \leq
\exp \li ( n_{\ell_\vep} \mscr{M}_{\mrm{P}} \li ( \f{\lm +
b_{\ell_\vep}}{2}, \lm \ri ) \ri )
\]} for small enough $\vep > 0$. It follows
that $\lim_{\vep \to 0} n_{\ell_\vep} \Pr \{ \bs{D}_{\ell_\vep} = 0
\} = 0$.

\epf

\beL

\la{Pos_rev_ineqinv}
 $\lim_{\vep \to 0} \sum_{\ell =  \tau +
1}^{\iy} n_\ell \Pr \{ \bs{l} = \ell \} = 0$ for any $\lm \in
(\lm^\prime, \lm^{\prime \prime})$.

\eeL

\bpf

Recalling that the sample sizes $n_1, n_2, \cd$ are chosen as the
ascending arrangement of all distinct elements of the set defined by
(\ref{defposrev}), we have that {\small \[ n_\ell = \li \lc \f{
C_{\tau - \ell} \; \ln (\ze \de) }{ \mscr{M}_{\mrm{P}} \li (
\lm^\prime, \f{\lm^\prime}{1 + \vep} \ri ) } \ri \rc, \qqu \ell = 1,
2, \cd
\]}
for small enough $\vep \in (0, 1)$. By the assumption that $\inf_{i
\in \bb{Z} } \f{C_{i - 1}} {C_{i} } = 1 + \udl{\ro} > 1$, we have
that
\[
n_\ell > (1 + \udl{\ro} )^{\ell - \tau - 1} \f{ \ln (\ze \de) } {
\mscr{M}_{\mrm{P}} \li ( \lm^\prime, \f{\lm^\prime}{1 + \vep} \ri )
}, \qqu \ell = \tau + 1, \tau + 2, \cd
\]
for small enough $\vep \in (0, 1)$. So, there exists a number
$\vep^* \in (0, 1)$ such that
\[
n_\ell \mscr{M}_{\mrm{P}} \li ( \lm^\prime, \f{\lm^\prime}{1 + \vep}
\ri ) < (1 + \udl{\ro})^{\ell - \tau - 1} \ln (\ze \de), \qqu \ell =
\tau + 1, \tau + 2, \cd
\]
for any $\vep \in (0, \vep^*)$. Observing that there exist a
positive integer $\ka^*$ such that $(1 + \udl{\ro})^{\ell - \tau -
1} \ln (\ze \de) < \ln (\ze \de) - (\ell - \tau) \ln 2 = \ln (\ze
\de_\ell)$ for any $\ell \geq \tau + \ka^*$, we have that there
exists a positive integer $\ka^*$ independent of $\vep$ such that
{\small $\mscr{M}_{\mrm{P}} ( \lm^\prime, \f{\lm^\prime}{1 + \vep} )
< \f{\ln (\ze \de_\ell)}{n_\ell}$} for $\ell \geq \tau + \ka^*$ and
$0 < \vep < \vep^*$.  Note that $\mscr{M}_{\mrm{P}} (z, \f{z}{1 +
\vep} ) = z
 [ \f{\vep}{1 + \vep} - \ln (1 + \vep)]$ is monotonically decreasing
 with respect to $z \in (0, \iy)$.  For $\ell \geq \tau + \ka^*$ and $0 <
\vep < \vep^*$, as a result of {\small $\f{\ln (\ze
\de_\ell)}{n_\ell} > \mscr{M}_{\mrm{P}} ( \lm^\prime,
\f{\lm^\prime}{1 + \vep} ) $}, there exists a unique number $z_\ell
\in [0, \iy)$ such that {\small $\mscr{M}_{\mrm{P}} (z_\ell,
\f{z_\ell}{1 + \vep} ) = \f{\ln (\ze \de_\ell)}{n_\ell} >
\mscr{M}_{\mrm{P}} ( \lm^\prime, \f{\lm^\prime}{1 + \vep} )$}.
Moreover, it must be true that $z_\ell < \lm^\prime$ for $\ell \geq
\tau + \ka^*$ and $\vep \in (0, \vep^*)$. Therefore, for small
enough $\vep \in (0, \vep^*)$, we have {\small \bee \sum_{\ell =
\tau + 1}^{\iy} n_\ell \Pr \{ \bs{l} = \ell \} & = & \sum_{\ell =
\tau + 1}^{\tau + \ka^*} n_\ell \Pr \{ \bs{l} = \ell \} + \sum_{\ell
= \tau + \ka^* + 1}^\iy n_\ell \Pr
\{ \bs{l} = \ell \}\\
& \leq & \sum_{\ell =  \tau + 1}^{\tau + \ka^*} n_\ell \Pr \{
\bs{D}_\tau  = 0 \} + \sum_{\ell =  \tau + \ka^* + 1}^\iy  n_\ell
\Pr \{ \bs{D}_{\ell - 1} = 0 \} \\
& = & \sum_{\ell =  \tau + 1}^{\tau + \ka^*} n_\ell \Pr \{
\bs{D}_\tau = 0 \} + \sum_{\ell =  \tau + \ka^*}^\iy  n_{\ell+1} \Pr
\{ \bs{D}_{\ell} = 0 \} \\
& \leq & k^* (1 + \ovl{\ro})^{k^*} n_\tau \Pr \{ \bs{D}_\tau = 0 \}
+ (1 + \ovl{\ro}) \sum_{\ell = \tau + \ka^*}^\iy  n_{\ell} \Pr \{
\bs{D}_{\ell} = 0 \} \\
& \leq  & k^* (1 + \ovl{\ro})^{k^*} n_\tau \Pr \{ \wh{\bs{\lm}}_\tau
< z_\tau \} + (1 + \ovl{\ro}) \sum_{\ell = \tau + \ka^*}^\iy
n_{\ell}
\Pr \{ \wh{\bs{\lm}}_{\ell} < z_\ell \}\\
& \leq  & k^* (1 + \ovl{\ro})^{k^*} n_\tau \Pr \li \{
\wh{\bs{\lm}}_\tau < \f{\lm^\prime + \lm}{2} \ri \} + (1 +
\ovl{\ro}) \sum_{\ell = \tau +
\ka^*}^\iy n_{\ell} \Pr \{ \wh{\bs{\lm}}_{\ell} < \lm^\prime \}\\
& \leq  & k^* (1 + \ovl{\ro})^{k^*} n_\tau \exp \li ( n_\tau
\mscr{M}_{\mrm{P}} \li ( \f{\lm + \lm^\prime}{2}, \lm \ri ) \ri)\\
&  & + (1 + \ovl{\ro}) \sum_{\ell = \tau + \ka^*}^\iy n_{\ell} \exp(
n_\ell \mscr{M}_{\mrm{P}} ( \lm^\prime, \lm) ) \to 0 \eee} as $\vep
\to 0$, where we have used the assumption that $\sup_{i \in \bb{Z}}
\f{C_{i - 1}}{C_i} = 1 + \ovl{\ro} < \iy$. This completes the proof
of the lemma.  \epf

\subsection{Proof of Theorem \ref{Pos_mix_CDF} } \la{App_Pos_mix_CDF}

We need some preliminary results.  The following results, stated as
Lemma \ref{lem99}, can be derived from Chernoff bounds.

\beL \la{lem99} {\small $S_{\mrm{P}} (0, k, n \lm) \leq \exp ( n
\mscr{M}_{\mrm{P}} ( \f{k}{n}, \lm ) )$} for $ 0 \leq k \leq n \lm$.
Similarly, {\small $S_{\mrm{P}} (k, \iy, n \lm) \leq \exp  ( n
\mscr{M}_{\mrm{P}} ( \f{k}{n}, \lm ) )$} for $k \geq n \lm$. \eeL

\beL \la{Plemm22} $\mscr{M}_{\mrm{P}} (\lm - \vep, \lm) <
\mscr{M}_{\mrm{P}} (\lm + \vep, \lm) < 0$ for any $\vep \in (0,
\lm]$.

\eeL

\bpf In the case of $\vep = \lm > 0$, we have $\mscr{M}_{\mrm{P}}
(\lm + \vep, \lm) = \vep - 2 \vep \ln 2 > - \vep =
\mscr{M}_{\mrm{P}} (\lm - \vep, \lm)$.  In the case of $0 < \vep <
\lm$, the lemma follows from the facts that $\mscr{M}_{\mrm{P}} (\lm
+ \vep, \lm) = \mscr{M}_{\mrm{P}} (\lm - \vep, \lm)$ for $\vep = 0$
and $\f{ \pa } { \pa \vep }  [ \mscr{M}_{\mrm{P}} (\lm + \vep, \lm)
- \mscr{M}_{\mrm{P}} (\lm - \vep, \lm) ] = \ln \f{\lm^2}{\lm^2 -
\vep^2} > 0$ for any $ \vep \in (0, \lm)$.   To show
$\mscr{M}_{\mrm{P}} (\lm + \vep, \lm) < 0$ for any $\vep \in (0,
\lm]$, note that $\mscr{M}_{\mrm{P}} (\lm + \vep, \lm) = \vep + (\lm
+ \vep ) \ln \f{\lm}{\lm + \vep} < \vep + (\lm + \vep )  \times \f{
- \vep}{\lm + \vep} = 0$.  This completes the proof of the lemma.

\epf

 \beL \la{Plemm21B} Let $\vep > 0$.  Then, $\mscr{M}_{\mrm{P}} (z, z - \vep)$
 is monotonically increasing with respect to $z > \vep$. \eeL

\bpf Note that $\mscr{M}_{\mrm{P}} (z, z - \vep)  =  \vep + z \ln
\li ( \f{z - \vep}{z} \ri )$ and \[ \f{ \pa \mscr{M}_{\mrm{P}} (z, z
- \vep) } { \pa z } = \ln \li (  \f{z - \vep}{z} \ri ) + \f{\vep}{z
- \vep} = - \ln \li ( 1 + \f{\vep}{z - \vep} \ri ) + \f{\vep}{z -
\vep} > 0
\]
where the last inequality follows from $\ln (1 + x) \leq  x, \; \fa
x \in [0, 1)$.

 \epf

\beL \la{Plemm24} Let $0 < \vep < 1$. Then,  {\small
$\mscr{M}_{\mrm{P}} ( z, \f{z}{1 - \vep} ) < \mscr{M}_{\mrm{P}} ( z,
\f{z}{1 + \vep} )$} and {\small $\f{\pa } { \pa z }
\mscr{M}_{\mrm{P}} ( z, \f{z}{1 - \vep}  ) < \f{\pa } { \pa z }
\mscr{M}_{\mrm{P}}  ( z, \f{z}{1 + \vep} ) < 0$} for $z > 0$. \eeL

\bpf Note that {\small $\mscr{M}_{\mrm{P}} ( z, \f{z}{1 + \vep} ) -
\mscr{M}_{\mrm{P}} ( z, \f{z}{1 - \vep}  )  = z \; g (\vep)$} where
$g(\vep) = \f{ \vep } { 1 + \vep } + \f{ \vep } { 1 - \vep } + \ln (
\f{1 - \vep} {1 + \vep} )$. Since $g(0)  = 0$ and $\f{d g (\vep) } {
d \vep  }  = \f{4 \vep^2} { (1 - \vep^2)^2 }
> 0$, we have $g (\vep)  > 0$ for $0 < \vep < 1$.  It follows that
{\small $\mscr{M}_{\mrm{P}} ( z, \f{z}{1 - \vep} ) <
\mscr{M}_{\mrm{P}} ( z, \f{z}{1 + \vep} )$}.

Using the inequality $\ln (1 - x) <  - x, \; \fa x \in (0, 1)$, we
have {\small $\f{\pa } { \pa z } \mscr{M}_{\mrm{P}} ( z, \f{z}{1 +
\vep} ) = \f{ \vep } { 1 + \vep } + \ln ( 1 - \f{\vep}{ 1 + \vep } )
< 0$}.  Noting that {\small $\f{\pa  } { \pa z }  [
\mscr{M}_{\mrm{P}} ( z, \f{z}{1 + \vep} )  - \mscr{M}_{\mrm{P}} ( z,
\f{z}{1 - \vep} ) ]  = g (\vep)
> 0$},  we have {\small $\f{\pa } { \pa z } \mscr{M}_{\mrm{P}} (
z, \f{z}{1 - \vep} ) < \f{\pa } { \pa z } \mscr{M}_{\mrm{P}} ( z,
\f{z}{1 + \vep} ) < 0$}.

\epf

\beL \la{PmixDS1} {\small $\Pr  \{ \mscr{M}_{\mrm{P}}  (
\wh{\bs{\lm}}_s, \mscr{L} ( \wh{\bs{\lm}}_s )  ) \leq \f{\ln (\ze
\de) } { n_s }, \; \mscr{M}_{\mrm{P}} ( \wh{\bs{\lm}}_s, \mscr{U} (
\wh{\bs{\lm}}_s ) ) \leq \f{\ln (\ze \de) } { n_s } \} = 1$}. \eeL

\bpf

For simplicity of notations, we denote $\lm^\star =
\f{\vep_a}{\vep_r}$.  By the definitions of $\mscr{L} (
\wh{\bs{\lm}}_s )$ and $\mscr{U} ( \wh{\bs{\lm}}_s )$, we have that,
in order to show the lemma, it suffices to show \bel & & \li \{
\mscr{M}_{\mrm{P}} \li ( \wh{\bs{\lm}}_s, \f{ \wh{\bs{\lm}}_s }{ 1 -
\vep_r} \ri )
> \f{\ln (\ze \de)}{n_s}, \;
\wh{\bs{\lm}}_s > \lm^\star - \vep_a \ri \} = \emptyset, \la{ep81POS}\\
&   & \li \{ \mscr{M}_{\mrm{P}} ( \wh{\bs{\lm}}_s, \wh{\bs{\lm}}_s +
\vep_a)
> \f{\ln (\ze \de)}{n_s}, \; \wh{\bs{\lm}}_s \leq
\lm^\star - \vep_a \ri \} = \emptyset, \la{ep82POS}\\
&   & \li \{ \mscr{M}_{\mrm{P}} \li ( \wh{\bs{\lm}}_s, \f{
\wh{\bs{\lm}}_s }{ 1 + \vep_r} \ri ) > \f{\ln (\ze \de)}{n_s}, \;
\wh{\bs{\lm}}_s > \lm^\star + \vep_a \ri \} = \emptyset, \la{ep83POS}\\
&   & \li \{ \mscr{M}_{\mrm{P}} ( \wh{\bs{\lm}}_s, \wh{\bs{\lm}}_s -
\vep_a)
> \f{\ln (\ze \de)}{n_s}, \; \wh{\bs{\lm}}_s \leq
\lm^\star + \vep_a \ri \} = \emptyset. \la{ep84POS} \eel By the
definition of $n_s$,  we have {\small $n_s \geq \li \lc  \f{ \ln
(\ze \de) } { \mscr{M}_{\mrm{P}} \li ( \lm^\star + \vep_a, \lm^\star
\ri ) } \ri \rc \geq \f{ \ln (\ze \de) } { \mscr{M}_{\mrm{P}} \li (
\lm^\star + \vep_a, \lm^\star \ri ) } $}.  By the assumption on
$\vep_a$ and $\vep_r$, we have $0 < \vep_a < \lm^\star$.  Hence, by
Lemma \ref{Plemm22}, we have $\mscr{M}_{\mrm{P}} \li ( \lm^\star -
\vep_a, \lm^\star \ri ) < \mscr{M}_{\mrm{P}} \li ( \lm^\star +
\vep_a, \lm^\star \ri ) < 0$ and it follows that \be \la{rela8POS}
\f{\ln (\ze \de)}{n_s} \geq \mscr{M}_{\mrm{P}} \li ( \lm^\star +
\vep_a, \lm^\star \ri ) >  \mscr{M}_{\mrm{P}} \li ( \lm^\star -
\vep_a, \lm^\star \ri ). \ee By (\ref{rela8POS}), {\small \be
\la{ep989POS} \li \{ \mscr{M}_{\mrm{P}} \li ( \wh{\bs{\lm}}_s, \f{
\wh{\bs{\lm}}_s }{ 1 - \vep_r} \ri )
> \f{\ln (\ze \de)}{n_s}, \; \wh{\bs{\lm}}_s > \lm^\star - \vep_a \ri \}
\subseteq \li \{ \mscr{M}_{\mrm{P}} \li ( \wh{\bs{\lm}}_s, \f{
\wh{\bs{\lm}}_s }{ 1 - \vep_r} \ri ) > \mscr{M}_{\mrm{P}} \li (
\lm^\star - \vep_a, \lm^\star \ri ), \; \wh{\bs{\lm}}_s > \lm^\star
- \vep_a \ri \}. \ee} Noting that $\mscr{M}_{\mrm{P}} \li (
\lm^\star - \vep_a, \lm^\star \ri ) = \mscr{M}_{\mrm{P}} \li (
\lm^\star - \vep_a, \f{\lm^\star - \vep_a}{1 - \vep_r} \ri )$ and
making use of the fact that $\mscr{M}_{\mrm{P}}(z, \f{z}{1 - \vep}
)$ is monotonically decreasing with respect to $z \in (0, \iy)$ as
asserted by Lemma \ref{Plemm24}, we have \be \la{ep99POS} \li \{
\mscr{M}_{\mrm{P}} \li ( \wh{\bs{\lm}}_s, \f{ \wh{\bs{\lm}}_s }{ 1 -
\vep_r} \ri ) > \mscr{M}_{\mrm{P}} \li ( \lm^\star - \vep_a,
\lm^\star \ri ) \ri \} = \{ \wh{\bs{\lm}}_s < \lm^\star - \vep_a \}.
\ee Combining (\ref{ep989POS}) and (\ref{ep99POS}) yields
(\ref{ep81POS}).  By (\ref{rela8POS}), {\small \be \la{eq9998POS}
\li \{ \mscr{M}_{\mrm{P}} ( \wh{\bs{\lm}}_s, \wh{\bs{\lm}}_s +
\vep_a)
> \f{\ln (\ze \de)}{n_s}, \; \wh{\bs{\lm}}_s \leq
\lm^\star - \vep_a \ri \} \subseteq \li \{ \mscr{M}_{\mrm{P}} (
\wh{\bs{\lm}}_s, \wh{\bs{\lm}}_s + \vep_a) > \mscr{M}_{\mrm{P}} \li
( \lm^\star - \vep_a, \lm^\star \ri ), \; \wh{\bs{\lm}}_s \leq
\lm^\star - \vep_a \ri \}. \ee} By the assumption on $\vep_a$ and
$\vep_r$, we have $\lm^\star - \vep_a > 0$. Recalling the fact that
$\mscr{M}_{\mrm{P}}(z, z + \vep)$ is monotonically increasing with
respect to $z \in (0, \iy)$ as asserted by Lemma \ref{Plemm21}, we
have that the event in the right-hand side of (\ref{eq9998POS}) is
an impossible event and consequently, (\ref{ep82POS}) is
established. By (\ref{rela8POS}), {\small \be \la{ep98998POS} \li \{
\mscr{M}_{\mrm{P}} \li ( \wh{\bs{\lm}}_s, \f{ \wh{\bs{\lm}}_s }{ 1 +
\vep_r} \ri ) > \f{\ln (\ze \de)}{n_s}, \; \wh{\bs{\lm}}_s >
\lm^\star + \vep_a \ri \} = \li \{ \mscr{M}_{\mrm{P}} \li (
\wh{\bs{\lm}}_s, \f{ \wh{\bs{\lm}}_s }{ 1 + \vep_r} \ri )
> \mscr{M}_{\mrm{P}} \li ( \lm^\star + \vep_a, \lm^\star \ri ), \; \wh{\bs{\lm}}_s >
\lm^\star + \vep_a \ri \}.  \ee} Noting that $\mscr{M}_{\mrm{P}} \li
( \lm^\star + \vep_a, \lm^\star \ri ) = \mscr{M}_{\mrm{P}} \li (
\lm^\star + \vep_a, \f{\lm^\star + \vep_a}{1 + \vep_r} \ri )$ and
making use of the fact that $\mscr{M}_{\mrm{P}}(z, \f{z}{1 + \vep}
)$ is monotonically decreasing with respect to $z \in (0, \iy)$ as
asserted by Lemma \ref{Plemm24},  we have \be \la{ep9998POS} \li \{
\mscr{M}_{\mrm{P}} \li ( \wh{\bs{\lm}}_s, \f{ \wh{\bs{\lm}}_s }{ 1 +
\vep_r} \ri )
> \mscr{M}_{\mrm{P}} \li ( \lm^\star + \vep_a, \lm^\star \ri ) \ri \} = \{
\wh{\bs{\lm}}_s < \lm^\star + \vep_a \}. \ee Combining
(\ref{ep98998POS}) and (\ref{ep9998POS}) yields (\ref{ep83POS}). By
(\ref{rela8POS}), {\small \be \la{eq999833POS} \li \{
\mscr{M}_{\mrm{P}} ( \wh{\bs{\lm}}_s, \wh{\bs{\lm}}_s - \vep_a)
> \f{\ln (\ze \de)}{n_s}, \; \wh{\bs{\lm}}_s \leq
\lm^\star + \vep_a \ri \} \subseteq \li \{ \mscr{M}_{\mrm{P}} (
\wh{\bs{\lm}}_s, \wh{\bs{\lm}}_s - \vep_a) > \mscr{M}_{\mrm{P}} \li
( \lm^\star + \vep_a, \lm^\star \ri ), \; \wh{\bs{\lm}}_s \leq
\lm^\star + \vep_a \ri \}. \ee}  Recalling the fact that
$\mscr{M}_{\mrm{P}}(z, z - \vep)$ is monotonically increasing with
respect to $z \in (\vep, \iy)$ as stated by Lemma \ref{Plemm21B}, we
have that the event in the right-hand side of (\ref{eq999833POS}) is
an impossible event and consequently, (\ref{ep84POS}) is
established. This completes the proof of the lemma.

\epf

\beL \la{cutlm} $\Pr \li \{ \li | \f{\wh{\bs{\lm}} - \lm } {\lm}
 \ri | \geq \vep_r  \mid \lm \ri \} < \de$ for $\lm \in [\ovl{\lm},
 \iy)$.
\eeL

\bpf

 Note that \bel \Pr \li \{ \li | \f{\wh{\bs{\lm}} - \lm } {\lm}
 \ri | \geq \vep_r  \mid \lm \ri \} & = & \sum_{\ell = 1}^s
 \Pr \li \{ \li | \f{\wh{\bs{\lm}}_\ell - \lm } {\lm}
 \ri | \geq \vep_r, \; \bs{l} = \ell  \mid \lm \ri \}
 \leq  \sum_{\ell = 1}^s \Pr \li \{ \li | \f{\wh{\bs{\lm}}_\ell -
\lm } {\lm} \ri | \geq \vep_r \mid \lm \ri \} \nonumber\\
 & \leq & \sum_{\ell = 1}^s \li [ \exp ( n_\ell \mscr{M}_{\mrm{P}} ( \lm
+ \lm \vep_r, \lm) ) + \exp ( n_\ell \mscr{M}_{\mrm{P}} ( \lm -
\lm \vep_r, \lm) )  \ri ] \la{ppp} \\
& < & 2 \sum_{\ell = 1}^s \exp ( n_\ell \mscr{M}_{\mrm{P}} ( \lm (1
+ \vep_r), \lm) ) \nonumber \eel where (\ref{ppp}) follows from
Lemma \ref{PCH}.  Since $\lim_{\lm \to 0} \mscr{M}_{\mrm{P}} ( \lm
(1 + \vep_r), \lm) = 0$ and $\lim_{\lm \to \iy} \mscr{M}_{\mrm{P}} (
\lm (1 + \vep_r), \lm) = - \iy$, there exists a unique number
$\ovl{\lm} > 0$ such that $\sum_{\ell = 1}^s \exp ( n_\ell
\mscr{M}_{\mrm{P}} ( \ovl{\lm} (1 + \vep_r), \ovl{\lm}) ) =
\f{\de}{2}$.  Finally, the lemma is established by noting that
$\mscr{M}_{\mrm{P}} ( \lm (1 + \vep_r), \lm)$ is monotonically
decreasing with respect to $\lm > 0$.

\epf

\bsk

Now we are in a position to prove Theorem \ref{Pos_mix_CDF}. The
second statement of Theorem \ref{Pos_mix_CDF} is a result of Lemma
\ref{cutlm}.

If the multistage sampling scheme follows a stopping rule derived
from Chernoff bounds,  then $\{ \bs{D}_s = 1 \}$ is a sure event as
a result of Lemma \ref{PmixDS1}.  It can be shown that $\exp(
\mscr{M}_{\mrm{P}} (z, \lm) )$  is equal to $\mcal{F} (z, \lm)$ and
$\mcal{G} (z, \lm)$ respectively for the cases of $z \leq \lm$ and
$z \geq \lm$.  Moreover, $\wh{\bs{\lm}}_\ell$ is a ULE of $\lm$ for
$\ell = 1, \cd, s$. So, the sampling scheme satisfies all the
requirements described in Corollary \ref{Monotone_third}, from which
Theorem \ref{Pos_mix_CDF} immediately follows.

If the multistage sampling scheme follows a stopping rule derived
from CDF $\&$ CCDF,  then, by Lemmas \ref{lem99} and \ref{PmixDS1},
we have {\small \bee  &  & \Pr \{ G_{\wh{\bs{\lm}}_s} (
\wh{\bs{\lm}}_s, \mscr{L} ( \wh{\bs{\lm}}_s ) ) \leq \ze \de_s \} =
\Pr \{ 1 - S_{\mrm{P}} (K_s - 1, n_s \mscr{L} ( \wh{\bs{\lm}}_s ) )
\leq \ze
\de \}\\
&  & \qqu \qqu  \qqu \qqu \qqu \qqu \geq \Pr  \{ n_s
\mscr{M}_{\mrm{P}} ( \wh{\bs{\lm}}_s, \mscr{L} ( \wh{\bs{\lm}}_s ) )
\leq \ln (\ze
\de) \} = 1,\\
&  & \Pr \{ F_{\wh{\bs{\lm}}_s} ( \wh{\bs{\lm}}_s, \mscr{U} (
\wh{\bs{\lm}}_s ) ) \leq \ze \de_s \} = \Pr \{ S_{\mrm{P}} (K_s, n_s
\mscr{U} ( \wh{\bs{\lm}}_s ) ) \leq \ze \de \} \geq \Pr \{ n_s
\mscr{M}_{\mrm{P}} ( \wh{\bs{\lm}}_s, \mscr{U} ( \wh{\bs{\lm}}_s ) )
\leq \ln (\ze \de) \} = 1 \eee} and thus $\Pr \{ F_{\wh{\bs{\lm}}_s}
( \wh{\bs{\lm}}_s, \mscr{U} ( \wh{\bs{\lm}}_s ) ) \leq \ze \de_s, \;
G_{\wh{\bs{\lm}}_s} ( \wh{\bs{\lm}}_s, \mscr{L} ( \wh{\bs{\lm}}_s )
) \leq \ze \de_s \} = 1$, which implies that $\{ \bs{D}_s = 1 \}$ is
a sure event.  So, the sampling scheme satisfies all the
requirements described in Theorem \ref{Monotone_second}, from which
Theorem \ref{Pos_mix_CDF} immediately follows.

\subsection{Proof of Theorem \ref{Pos_mix_DDV_Asp} }  \la{App_Pos_mix_DDV_Asp}

We need some preliminary results.

\beL

\la{lem98PS}
If $\vep_a$ is sufficiently small, then the following
statements hold true.

(I): For $1 \leq \ell < s$, there exists a unique number $z_\ell \in
[0, \lm^\star - \vep_a)$ such that $n_\ell = \f{ \ln ( \ze \de ) } {
\mscr{M}_{\mrm{P}} ( z_\ell, z_\ell + \vep_a ) }$.

(II): For $1 \leq \ell < s$, there exists a unique number $y_\ell
\in (\lm^\star + \vep_a, \iy)$ such that $n_\ell = \f{ \ln ( \ze \de
) } { \mscr{M}_{\mrm{P}} \li ( y_\ell, \f{ y_\ell }{1 + \vep_r} \ri
) }$.

(III): $z_\ell$ is monotonically increasing with respect to $\ell$;
$y_\ell$ is monotonically decreasing with respect to $\ell$.

(IV): $\lim_{\vep_a \to 0} z_\ell = \lm^\star C_{s - \ell}$ and
$\lim_{\vep_a \to 0} y_\ell = \f{ \lm^\star}{ C_{ s - \ell} }$,
where the limits are taken under the constraint that
$\f{\vep_a}{\vep_r}$ and $s - \ell$ are fixed with respect to
$\vep_a$.

(V): Let $\ell_\vep = s - j_\lm$.  For $\lm \in (\lm^\star, \iy)$
such that $C_{j_\lm} = r (\lm)$,
\[
\lim_{\vep_r \to 0} \f{ z_{\ell_\vep} - \lm}{\vep_r \lm} = 1.
\]
For $\lm \in (0, \lm^\star)$ such that $C_{j_\lm} = r (\lm)$,
\[
\lim_{\vep_a \to 0} \f{ z_{\ell_\vep} - \lm}{\vep_a} = \f{1}{3} \li ( \f{ \lm}{\lm^\star} - 2 \ri ).
\]

(VI): $\{ \bs{D}_\ell = 0 \} = \{ z_\ell < \wh{\bs{\lm}}_\ell <
y_\ell \}$ for $1 \leq \ell < s$.

\eeL

\bsk

{\bf Proof of Statement (I)}:

By the definition of sample sizes, we have {\small $\f{ \ln ( \ze
\de ) } { n_{\ell} } \geq \mscr{M}_{\mrm{P}} ( 0, \vep_a )$} and \be
\la{from2}
 n_{\ell} < \f{(1 + C_1) n_s}{2} < \f{(1 + C_1)}{2}
\li [ \f{ \ln ( \ze \de ) } { \mscr{M}_{\mrm{P}} ( \lm^\star +
\vep_a,  \lm^\star) } + 1 \ri ]
 \ee
for sufficiently small $\vep_a > 0$. By (\ref{from2}), we have
{\small \[ \f{ \ln ( \ze \de ) } { n_{\ell} }  < \mscr{M}_{\mrm{P}}
( \lm^\star + \vep_a,  \lm^\star)  \li ( \f{2} {1 + C_1} - \f{1}{
n_{\ell} } \ri ) = \f{\mscr{M}_{\mrm{P}} ( \lm^\star + \vep_a,
\lm^\star) } { \mscr{M}_{\mrm{P}} ( \lm^\star - \vep_a, \; \lm^\star
) } \f{2 \mscr{M}_{\mrm{P}} ( \lm^\star - \vep_a, \; \lm^\star ) }
{1 + C_1}   - \f{\mscr{M}_{\mrm{P}} ( \lm^\star + \vep_a, \lm^\star)
}{ n_{\ell} } .
\]}
Noting that
\[
\lim_{\vep_a \to 0} \f{\mscr{M}_{\mrm{P}} ( \lm^\star + \vep_a,  \lm^\star) } { \mscr{M}_{\mrm{P}} ( \lm^\star - \vep_a, \; \lm^\star
) } = 1, \qqu \lim_{\vep_a \to 0} \f{\mscr{M}_{\mrm{P}} ( \lm^\star + \vep_a,  \lm^\star) }{ n_{\ell} } = 0,
\]
we have that $\f{ \ln ( \ze \de ) } { n_{\ell} } <
\mscr{M}_{\mrm{P}} ( \lm^\star - \vep_a, \; \lm^\star )$ for small
enough $\vep_a > 0$.  In view of the established fact that $ \mscr{M}_{\mrm{P}} ( 0,
\vep_a ) \leq \f{ \ln ( \ze \de ) } { n_{\ell} } <
\mscr{M}_{\mrm{P}} \li ( \lm^\star - \vep_a, \; \lm^\star \ri )$ and
the fact that $\mscr{M}_{\mrm{P}} ( z, z + \vep_a )$ is
monotonically increasing with respect to $z  > 0$ as
asserted by Lemma \ref{Plemm21}, invoking the intermediate value
theorem, we have that there exists a unique number $z_{\ell} \in [0,
\lm^\star - \vep_a)$ such that $\mscr{M}_{\mrm{P}} ( z_{\ell},
z_{\ell} + \vep_a )  = \f{ \ln ( \ze \de ) } { n_{\ell} }$, which
implies Statement (I).

\bsk

{\bf Proof of Statement (II)}: By (\ref{from2}), we have {\small
\bee \f{ \ln ( \ze \de ) } { n_{\ell} } <   \mscr{M}_{\mrm{P}} (
\lm^\star + \vep_a,  \lm^\star)  \li ( \f{2} {1 + C_1} - \f{1}{
n_{\ell} } \ri ) = \li ( \f{2} {1 + C_1} \ri ) \mscr{M}_{\mrm{P}}
\li ( \lm^\star + \vep_a, \; \f{\lm^\star + \vep_a}{ 1 + \vep_r }
\ri ) - \f{\mscr{M}_{\mrm{P}} ( \lm^\star + \vep_a,  \lm^\star) }{
n_{\ell} } . \eee} Noting that {\small $\lim_{\vep_a \to 0}
\f{\mscr{M}_{\mrm{P}} ( \lm^\star + \vep_a,  \lm^\star) }{ n_{\ell}
} = 0$}, we have that {\small $\f{ \ln ( \ze \de ) } { n_{\ell} } <
\mscr{M}_{\mrm{P}} ( \lm^\star + \vep_a, \; \f{\lm^\star + \vep_a}{
1 + \vep_r } )$} for small enough $\vep_a > 0$.  In view of the
established fact that {\small $ \f{ \ln ( \ze \de ) } { n_{\ell} } <
\mscr{M}_{\mrm{P}} ( \lm^\star + \vep_a, \; \f{\lm^\star + \vep_a}{
1 + \vep_r } )$} and the fact that $\mscr{M}_{\mrm{P}} ( z, \f{z}{1
+ \vep_r} )$ is monotonically decreasing to $- \iy$ with respect to
$z \in (0, \iy)$ as asserted by Lemma \ref{Plemm24},  invoking the
intermediate value theorem, we have that there exists a unique
number $y_{\ell} \in (\lm^\star + \vep_a, \iy)$ such that
$\mscr{M}_{\mrm{P}} ( y_{\ell}, \f{y_{\ell}}{1 + \vep_r} ) = \f{ \ln
( \ze \de ) } { n_{\ell} }$, which implies Statement (II).

\bsk

{\bf Proof of Statement (III)}:  Since $n_{\ell}$ is monotonically
increasing with respect to $\ell$ if $\vep_a > 0$ is sufficiently
small, we have that $\mscr{M}_{\mrm{P}} ( z_{\ell}, z_{\ell} +
\vep_a )$ is monotonically increasing with respect to $\ell$ for
small enough $\vep_a > 0$.  Recalling that $\mscr{M}_{\mrm{P}} ( z,
z + \vep_a )$ is monotonically increasing with respect to $z  > 0$,
we have that $z_{\ell}$ is monotonically increasing with respect to
$\ell$. Similarly, $\mscr{M}_{\mrm{P}} ( y_{\ell}, \f{y_{\ell}}{1 +
\vep_r} )$ is monotonically increasing with respect to $\ell$ for
sufficiently small $\vep_a > 0$. Recalling that $\mscr{M}_{\mrm{P}}
( z, \f{z}{1 + \vep_r} )$ is monotonically decreasing with respect
to $z > 0$, we have that $y_{\ell}$ is monotonically decreasing with
respect to $\ell$. This establishes Statement (III).

\bsk

{\bf Proof of Statement (IV)}: We first consider $\lim_{\vep_a \to
0} z_\ell$.  For simplicity of notations, define $b_\ell = \lm^\star
C_{s - \ell}$ for $\ell < s$. Then, it can be checked that $\f{
b_\ell }{ \lm^\star } = C_{s - \ell}$ and, by the definition of
sample sizes, we have
 \be \la{notiabsP} \f{ b_\ell }{ \lm^\star }
 \f{ \mscr{M}_{\mrm{P}}
(z_\ell, z_\ell  + \vep_a ) } {\mscr{M}_{\mrm{P}} ( \lm^\star +
\vep_a, \lm^\star )} = \f{1}{n_{\ell} } \times \f{C_{s - \ell} \;
\ln (\ze \de) } { \mscr{M}_{\mrm{P}} ( \lm^\star + \vep_a, \lm^\star
) } = 1 + o(1) \ee for $\ell < s$.

We claim that $z_\ell
> \se$ for $\se \in (0, b_\ell)$ if $\vep_a > 0$ is small enough.
To prove this claim, we use a contradiction method.
Suppose this claim is not true, then there is a set, denoted by $S_{\vep_a}$,  of infinitely many values
of $\vep_a$ such that $z_\ell \leq \se$ for any $\vep_a \in S_{\vep_a}$.
By (\ref{notiabsP}) and the fact that $\mscr{M}_{\mrm{P}} ( z, z + \vep_a )$ is
monotonically increasing with respect to $z  > 0$ as
asserted by Lemma \ref{Plemm21}, we have
\[
\f{ b_\ell  }{ \lm^\star  }
 \f{ \mscr{M}_{\mrm{P}}
(z_\ell, z_\ell  + \vep_a ) } {\mscr{M}_{\mrm{P}} ( \lm^\star + \vep_a, \lm^\star )} = 1 + o(1)  \geq \f{ b_\ell  }{ \lm^\star  }
\f{ \mscr{M}_{\mrm{P}} (\se, \se  + \vep_a ) } {\mscr{M}_{\mrm{P}}
( \lm^\star + \vep_a, \lm^\star )} = \f{b_\ell }{\se } + o(1)
\]
for small enough $\vep_a \in S_{\vep_a}$, which implies
$\f{b_\ell}{\se} \leq 1$, contradicting to the fact that
$\f{b_\ell}{\se} > 1$.  This proves the claim. Now we restrict
$\vep_a$ to be small enough so that $\se < z_\ell < \lm^\star$.
Since $z_\ell$ is bounded in interval $(\se, \lm^\star)$, we have $\mscr{M}_{\mrm{P}}
(z_\ell, z_\ell  + \vep_a ) = - \vep_a^2 \sh (2 z_\ell) + o
(\vep_a^2)$ and by (\ref{notiabsP}), we have
\[
\f{ b_\ell }{ \lm^\star }  \times  \f{ - \vep_a^2 \sh (2 z_\ell) + o
(\vep_a^2)  } { - \vep_a^2 \sh (2 \lm^\star) + o (\vep_a^2)  } = 1 +
o(1),
\]
which implies $\f{ b_\ell } { z_\ell }  = 1 + o(1)$ and thus $\lim_{\vep_a \to 0} z_\ell = b_\ell$.

We now consider $\lim_{\vep_a \to 0} y_\ell$.  For simplicity of
notations, define $a_\ell = \f{\lm^\star}{C_{s - \ell}}$ for $1 \leq
\ell < s$. Then, it can be checked that $\f{ \lm^\star }{a_\ell} =
C_{s - \ell}$ and, by the definition of sample sizes, we have \be
\la{notiPP} \f{ \lm^\star }{a_\ell} \f{ \mscr{M}_{\mrm{P}} (y_\ell,
\f{y_\ell}{1 + \vep_r} ) } {\mscr{M}_{\mrm{P}} ( \lm^\star + \vep_a,
\lm^\star )} = \f{1}{n_{\ell} } \times \f{C_{s - \ell} \;  \ln (\ze
\de) } {\mscr{M}_{\mrm{P}} ( \lm^\star + \vep_a, \lm^\star )} = 1 +
o(1). \ee

We claim that $y_\ell < \se$ for $\se \in (a_\ell, \iy)$ if $\vep_r > 0$ is small
enough.  To prove this claim, we use a contradiction method.
Suppose this claim is not true, then there is a set,
denoted by $S_{\vep_r}$,  of infinitely many values of $\vep_r$ such
that $y_\ell \geq \se$ for any $\vep_r \in S_{\vep_r}$.
By (\ref{notiPP}) and the fact that $\mscr{M}_{\mrm{P}} (
z, \f{z}{1 + \vep_r} )$ is monotonically decreasing  with
respect to $z \in (0, \iy)$ as asserted by Lemma \ref{Plemm24}, we have
\[
\f{ \lm^\star }{a_\ell}  \f{ \mscr{M}_{\mrm{P}} (y_\ell, \f{y_\ell}{1 + \vep_r} ) }
{\mscr{M}_{\mrm{P}} ( \lm^\star + \vep_a, \lm^\star )} = 1 + o(1) \geq \f{ \lm^\star }{a_\ell} \f{
\mscr{M}_{\mrm{P}} (\se, \f{\se}{1 + \vep_r} ) }
{\mscr{M}_{\mrm{P}} ( \lm^\star + \vep_a, \lm^\star )} =
\f{\se}{a_\ell} + o(1)
\]
for small enough $\vep_r \in S_{\vep_r}$, which implies
$\f{\se}{a_\ell} \leq 1$, contradicting to the fact that
$\f{\se}{a_\ell} > 1$.  This proves the claim. Now we restrict
$\vep_r$ to be small enough so that $\lm^\star < y_\ell < \se$.
Since $y_\ell$ is bounded in interval $(\lm^\star, \se)$, we have
$\mscr{M}_{\mrm{P}} (y_\ell, \f{y_\ell}{1 + \vep_r} ) = - \vep_r^2 y_\ell \sh 2   + o
(\vep_r^2)$ and by (\ref{notiPP}), we have
\[
\f{ \lm^\star }{a_\ell} \times  \f{ - \vep_r^2 y_\ell \sh 2   + o
(\vep_r^2)  } { - \vep_a^2  \sh (2 \lm^\star) +  o (\vep_a^2)  } = 1 +
o(1),
\]
which implies $\f{ y_\ell  - a_\ell } { a_\ell  }  =  o(1)$ and thus $\lim_{\vep_r \to 0} y_\ell = a_\ell$.

\bsk

{\bf Proof of Statement (V)}:

We shall first consider $\lm \in (\lm^\star, \iy)$ such that $C_{j_\lm} = \f{ \lm^\star }{ \lm  }$.  Note that \be \la{moreadd}
\mscr{M}_{\mrm{P}} (\lm + \vep, \lm) = \vep - (\lm + \vep) \ln \li ( 1 + \f{\vep}{\lm} \ri ) = - \f{ \vep^2}{ 2 \lm } + \f{ \vep^3}{ 6 \lm^2} +
o(\vep) \ee for $\lm > 0$ and $\vep > 0$.  Let $\psi_\ep$ be a function of $\ep \in (0, 1)$ such that $|\psi_\ep|$ is bounded from above by a
constant independent of $\ep$. Then, by Taylor's series expansion formula, we have \bel \mscr{M}_{\mrm{P}} \li (\psi_\ep, \f{\psi_\ep}{1 + \ep}
\ri ) & = & \f{\ep \; \psi_\ep}{1 + \ep} - \psi_\ep \ln (1 + \ep)  = \ep \; \psi_\ep \li [ 1 - \ep + \ep^2 + o(\ep^2) \ri ] -
\psi_\ep \li [ \ep - \f{\ep^2}{2} + \f{\ep^3}{3} + o (\ep^3) \ri ] \nonumber\\
& = & - \f{\ep^2 \psi_\ep}{2} + \f{2 \ep^3 \psi_\ep}{3} + o (\ep^3) \la{use381} \eel  for $\ep \in (0, 1)$.   By the definition of sample sizes,
for small enough $\vep_r$, there exists $z_{\ell_\vep} \in (\lm^\star, \iy)$ such that \be \la{use382} n_{\ell_\vep} =  \f{ \ln (\ze \de) } {
\mscr{M}_{\mrm{P}} (z_{\ell_\vep}, z_{\ell_\vep} \sh (1 + \vep_r) ) }  = \li \lc \f{ C_{s - {\ell_\vep}} \; \ln (\ze \de) } {\mscr{M}_{\mrm{P}}
( \lm^\star + \vep_a, \lm^\star )} \ri \rc = \li \lc \f{ \lm^\star }{ \lm } \f{ \ln (\ze \de) } {\mscr{M}_{\mrm{P}} ( \lm^\star + \vep_a,
\lm^\star )} \ri \rc, \ee from which we can use an argument similar to the proof of Statement (III) to deduce that $z_{\ell_\vep}$ is smaller
than $\se$ for $\se \in (\lm, \iy)$ if $\vep_r > 0$ is small enough.   Hence, by (\ref{moreadd}), (\ref{use381}) and (\ref{use382}), we have
\[
1 + o (\vep_r) = \f{ \f{ \lm^\star }{ \lm  } \f{ \ln (\ze \de) } {\mscr{M}_{\mrm{P}} ( \lm^\star + \vep_a, \lm^\star )} } { \f{ \ln (\ze \de) }
{ \mscr{M}_{\mrm{P}} (z_{\ell_\vep}, z_{\ell_\vep} \sh (1 + \vep_r) )  }} = \f{ \lm^\star }{ \lm  } \f{ \mscr{M}_{\mrm{P}} (z_{\ell_\vep},
z_{\ell_\vep} \sh (1 + \vep_r) ) } {\mscr{M}_{\mrm{P}} ( \lm^\star + \vep_a, \lm^\star )} = \f{ \lm^\star }{ \lm  }  \f{ - \f{\vep_r^2
z_{\ell_\vep} }{2} + \f{2 \vep_r^3 z_{\ell_\vep}}{3} + o (\vep_r^3)  } { - \f{\vep_r^2 \lm^\star }{2} + \f{\vep_r^3 \lm^\star}{6} + o (\vep_r^3)
},
\]
and consequently,
\[
1 + o (\vep_r) =  \f{ \lm^\star }{ \lm  }  \f{ z_{\ell_\vep} - \f{4 \vep_r z_{\ell_\vep}}{3} + o (\vep_r)  } {  \lm^\star  - \f{ \vep_r
\lm^\star}{3} + o (\vep_r) },
\]
which implies that
\[
\lm^\star \li ( z_{\ell_\vep} - \f{4 \vep_r z_{\ell_\vep}}{3} \ri ) = \lm \li ( \lm^\star - \f{\vep_r \lm^\star}{3} \ri ) + o (\vep_r),
\]
i.e., $z_{\ell_\vep} \li ( 1 - \f{4 \vep_r }{3} \ri )  = \lm \li ( 1 - \f{\vep_r }{3} \ri ) + o (\vep_r)$.  It follows that $\lim_{\vep_r \to 0}
\f{ z_{\ell_\vep} - \lm}{\vep_r \lm} = 1$ and thus
\[
\lim_{\vep_r \to 0} \f{ z_{\ell_\vep} - \lm}{ \sq{ \lm \sh n_{\ell_\vep}} } = b \lim_{\vep_r \to 0} \f{ z_{\ell_\vep} - \lm}{\vep_r \lm} = b.
\]

Next, we shall now consider $\lm \in (0, \lm^\star)$ such that $C_{j_\lm} = \f{\lm}{ \lm^\star }$.  Let $\psi_\ep$ be a function of $\ep \in (0,
\iy)$ such that $\f{1}{|\psi_\ep|}$ is bounded from above by a constant independent of $\ep$.  Then, by Taylor's series expansion formula, we
have {\small \be \la{use881} \mscr{M}_{\mrm{P}} (\psi_\ep, \psi_\ep  + \ep ) = - \ep  + \psi_\ep \ln \li ( 1 + \f{\ep}{ \psi_\ep } \ri ) =  -
\ep + \psi_\ep \li [ \f{\ep}{\psi_\ep} - \f{\ep^2}{2 \psi_\ep^2} + \f{\ep^3}{3 \psi_\ep^3} + o (\ep^3) \ri ] = - \f{\ep^2 }{2 \psi_\ep} +
\f{\ep^3}{3 \psi_\ep^2} + o (\ep^3). \qqu \ee} By the definition of sample sizes, for small enough $\vep_a$, there exists $z_{\ell_\vep} \in (0,
\lm^\star)$ such that \be \la{use882}
 n_{\ell_\vep} =  \f{ \ln (\ze
\de) } { \mscr{M}_{\mrm{P}} (z_{\ell_\vep}, z_{\ell_\vep} + \vep_a ) }  = \li \lc \f{ C_{s - {\ell_\vep}} \; \ln (\ze \de) } {\mscr{M}_{\mrm{P}}
(\lm^\star + \vep_a, \lm^\star )}  \ri \rc = \li \lc \f{\lm}{ \lm^\star } \f{ \ln (\ze \de) } {\mscr{M}_{\mrm{P}} (\lm^\star + \vep_a, \lm^\star
)}  \ri \rc. \ee from which we can use an argument similar to the proof of Statement (III) to deduce that $z_{\ell_\vep}$ is greater than $\se$
for $\se \in (0, \lm)$ if $\vep_a > 0$ is small enough. Hence, by (\ref{moreadd}), (\ref{use881}) and (\ref{use882}), we have
\[
1 + o (\vep_a) = \f{ \f{\lm}{ \lm^\star } \f{ \ln (\ze \de) } {\mscr{M}_{\mrm{P}} (\lm^\star + \vep_a, \lm^\star )} } { \f{ \ln (\ze \de) } {
\mscr{M}_{\mrm{P}} (z_{\ell_\vep}, z_{\ell_\vep} + \vep_a ) } } = \f{\lm}{ \lm^\star } \f{ \mscr{M}_{\mrm{P}} (z_{\ell_\vep}, z_{\ell_\vep} +
\vep_a ) } {\mscr{M}_{\mrm{P}} (\lm^\star + \vep_a, \lm^\star )} = \f{\lm}{ \lm^\star } \f{ - \f{\vep_a^2 }{2 z_{\ell_\vep}} + \f{\vep_a^3}{3
z_{\ell_\vep}^2} + o (\vep_a^3) } {  - \f{\vep_a^2 }{2 \lm^\star} + \f{\vep_a^3}{6 (\lm^\star)^2} + o (\vep_a^3) }
\]
and consequently,
\[
1 + o (\vep_a) =  \f{ \f{\lm }{z_{\ell_\vep}} - \f{2 \vep_a \lm}{3 z_{\ell_\vep}^2} + o (\vep_a) } {  1  - \f{ \vep_a}{3 \lm^\star} + o (\vep_a)
},
\]
which implies that $\f{\lm }{z_{\ell_\vep}} - \f{2 \vep_a \lm}{3 z_{\ell_\vep}^2} = 1  - \f{ \vep_a}{3 \lm^\star} + o (\vep_a)$, i.e.,
$\f{z_{\ell_\vep}  - \lm } {\vep_a} = \f{z_{\ell_\vep}}{3 \lm^\star}  - \f{ 2 \lm} { 3 z_{\ell_\vep} } + z_{\ell_\vep} \f{o(\vep_a)}{\vep_a}$.
So, we have
\[
\lim_{\vep_a \to 0} \f{ z_{\ell_\vep} - \lm}{\vep_a} = \f{1}{3} \li ( \f{ \lm}{\lm^\star} - 2 \ri )  \in \li ( - \f{2}{3}, - \f{1}{3} \ri ).
\]

\bsk

{\bf Proof of Statement (VI)}:  By the definition of the sampling
scheme, we have {\small \bee  \{ \bs{D}_{\ell} = 0  \} & = & \li \{
\max \{ \mscr{M}_{\mrm{P}} (\wh{\bs{\lm}}_\ell, \udl{\bs{\lm}}_\ell
), \; \mscr{M}_{\mrm{P}} (\wh{\bs{\lm}}_\ell, \ovl{\bs{\lm}}_\ell )
\} > \f{ \ln ( \ze \de  )
} { n_\ell }, \;  | \wh{\bs{\lm}}_\ell - \lm^\star | \leq \vep_a \ri \}\\
&  & \bigcup \li \{ \max \{ \mscr{M}_{\mrm{P}} (\wh{\bs{\lm}}_\ell,
\udl{\bs{\lm}}_\ell ), \; \mscr{M}_{\mrm{P}} (\wh{\bs{\lm}}_\ell,
\ovl{\bs{\lm}}_\ell ) \} > \f{ \ln ( \ze \de  ) } { n_\ell }, \;
\wh{\bs{\lm}}_\ell  < \lm^\star - \vep_a \ri \}\\
&  & \bigcup \li \{ \max \{ \mscr{M}_{\mrm{P}} (\wh{\bs{\lm}}_\ell,
\udl{\bs{\lm}}_\ell ), \; \mscr{M}_{\mrm{P}} (\wh{\bs{\lm}}_\ell,
\ovl{\bs{\lm}}_\ell ) \} > \f{ \ln ( \ze \de  ) } { n_\ell }, \;
\wh{\bs{\lm}}_\ell  > \lm^\star + \vep_a \ri \}\\
& = & \li \{ \max \li \{ \mscr{M}_{\mrm{P}}  ( \wh{\bs{\lm}}_\ell,
\wh{\bs{\lm}}_\ell - \vep_a), \; \mscr{M}_{\mrm{P}} \li (
\wh{\bs{\lm}}_\ell, \f{\wh{\bs{\lm}}_\ell}{1 - \vep_r} \ri ) \ri \}
> \f{ \ln ( \ze \de  ) } { n_\ell }, \; | \wh{\bs{\lm}}_\ell - \lm^\star | \leq \vep_a \ri \}\\
&  & \bigcup \li \{ \mscr{M}_{\mrm{P}} (\wh{\bs{\lm}}_\ell,
\wh{\bs{\lm}}_\ell + \vep_a ) > \f{ \ln ( \ze \de  ) } { n_\ell },
\; \wh{\bs{\lm}}_\ell  < \lm^\star - \vep_a \ri \} \\
&   &  \bigcup \li \{ \mscr{M}_{\mrm{P}} \li ( \wh{\bs{\lm}}_\ell,
\f{\wh{\bs{\lm}}_\ell}{1 + \vep_r} \ri ) > \f{ \ln ( \ze \de  ) } {
n_\ell }, \; \wh{\bs{\lm}}_\ell
> \lm^\star + \vep_a \ri \}. \eee}

We claim that, {\small \bel &  &  \li \{ \max \li \{
\mscr{M}_{\mrm{P}} (\wh{\bs{\lm}}_\ell, \wh{\bs{\lm}}_\ell -
\vep_a), \; \mscr{M}_{\mrm{P}} \li ( \wh{\bs{\lm}}_\ell,
\f{\wh{\bs{\lm}}_\ell}{1 - \vep_r} \ri ) \ri \}
> \f{ \ln ( \ze \de ) } { n_\ell }, \;  | \wh{\bs{\lm}}_\ell - \lm^\star | \leq \vep_a \ri \}
 =  \li \{ | \wh{\bs{\lm}}_\ell - \lm^\star | \leq \vep_a \ri
\}, \nonumber\\
 \la{thatP0}\\
&   & \li \{ \mscr{M}_{\mrm{P}} (\wh{\bs{\lm}}_\ell,
\wh{\bs{\lm}}_\ell + \vep_a ) > \f{ \ln ( \ze \de ) } { n_\ell }, \;
\wh{\bs{\lm}}_\ell  < \lm^\star - \vep_a \ri \} = \{  z_\ell <
\wh{\bs{\lm}}_\ell  < \lm^\star - \vep_a \},
\la{thatP}\\
&   & \li \{ \mscr{M}_{\mrm{P}} \li ( \wh{\bs{\lm}}_\ell,
\f{\wh{\bs{\lm}}_\ell}{ 1  + \vep_r} \ri ) > \f{ \ln ( \ze \de  ) }
{ n_\ell }, \; \wh{\bs{\lm}}_\ell  > \lm^\star + \vep_a \ri \} = \li
\{ \lm^\star + \vep_a <  \wh{\bs{\lm}}_\ell < y_\ell \ri \}
\la{that3P} \eel} for $1 \leq \ell < s$ provided that $\vep_a$ is
sufficiently small.

To show (\ref{thatP0}), note that \be \la{sinceP}
 n_{\ell} < \f{(1 + C_1) n_s}{ 2 } < \f{1 + C_1}{2} \li [  \f{
\ln ( \ze \de ) } { \mscr{M}_{\mrm{P}} ( \lm^\star + \vep_a,
\lm^\star ) } + 1 \ri ],  \ee  from which we have
\[
\f{ \ln ( \ze \de ) } { n_{\ell} } < \f{ \mscr{M}_{\mrm{P}} (
\lm^\star + \vep_a, \lm^\star ) } {\mscr{M}_{\mrm{P}} ( \lm^\star -
\vep_a,  \lm^\star - \vep_a - \vep_a )}  \li ( \f{2} {1 + C_1} \ri )
\mscr{M}_{\mrm{P}} ( \lm^\star - \vep_a,  \lm^\star - \vep_a -
\vep_a ) - \f{\mscr{M}_{\mrm{P}} (  \lm^\star + \vep_a, \lm^\star
)}{n_\ell}.
\]
Noting that
\[
\lim_{\vep_a \to 0} \f{ \mscr{M}_{\mrm{P}} (  \lm^\star + \vep_a, \lm^\star )  }
{  \mscr{M}_{\mrm{P}} ( \lm^\star - \vep_a, \lm^\star -
\vep_a - \vep_a ) } = \lim_{\vep_a \to 0}  \f{ - \f{ \vep_a^2 } { 2
\lm^\star  } + o(\vep_a^2)  } { - \f{ \vep_a^2 } { 2 (\lm^\star -
\vep_a)   } + o(\vep_a^2) } = 1
\]
and $\lim_{\vep_a \to 0} \f{\mscr{M}_{\mrm{P}} (
\lm^\star + \vep_a, \lm^\star )}{n_\ell} = 0$, we have \be \la{have1PP} \f{ \ln
( \ze \de ) } { n_{\ell} } < \mscr{M}_{\mrm{P}} ( \lm^\star -
\vep_a, \lm^\star - \vep_a - \vep_a ) \ee for small enough $\vep_a
> 0$. Again by (\ref{sinceP}), we have
\[
\f{ \ln ( \ze \de ) } { n_{\ell} } < \f{ \mscr{M}_{\mrm{P}} (
 \lm^\star + \vep_a, \lm^\star ) } { \mscr{M}_{\mrm{P}} ( \lm^\star
+ \vep_a, \f{\lm^\star + \vep_a}{1 -  \vep_r} ) } \li ( \f{2} {1 +
C_1} \ri ) \mscr{M}_{\mrm{P}} \li ( \lm^\star + \vep_a, \f{\lm^\star
+ \vep_a}{1 -  \vep_r} \ri ) - \f{\mscr{M}_{\mrm{P}} (
 \lm^\star + \vep_a, \lm^\star )}{n_\ell}.
\]
Noting that
\[
\lim_{\vep_a \to 0} \f{ \mscr{M}_{\mrm{P}} (  \lm^\star + \vep_a, \lm^\star ) }
{ \mscr{M}_{\mrm{P}} ( \lm^\star + \vep_a, \f{\lm^\star +
\vep_a}{1 -  \vep_r} ) } = \lim_{\vep_a \to 0} \f{ - \f{ \vep_a^2 }
{ 2 \lm^\star    } + o(\vep_a^2) } {  - \f{ \vep_a^2 } { 2
(\lm^\star + \vep_a) } + o \li ( \f{(\lm^\star + \vep_a)^2
\vep_r^2}{(1 - \vep_r)^2} \ri ) } = 1
\]
and $\lim_{\vep_a \to 0} \f{\mscr{M}_{\mrm{P}} ( \lm^\star + \vep_a, \lm^\star )}{n_\ell} = 0$, we have \be \la{have2PP}
 \f{ \ln ( \ze \de )
} { n_{\ell} } < \mscr{M}_{\mrm{P}} \li ( \lm^\star + \vep_a,
\f{\lm^\star + \vep_a}{1 -  \vep_r } \ri ) \ee for small enough
$\vep_a
> 0$. Note that, for $z \in [\lm^\star - \vep_a, \; \lm^\star +
\vep_a]$, $\mscr{M}_{\mrm{P}} ( z, z - \vep_a)$ is monotonically
increasing with respect to $z$ and {\small $\mscr{M}_{\mrm{P}} ( z,
\f{z}{1 - \vep_r} )$} is monotonically decreasing with respect to
$z$. By (\ref{have1PP}) and (\ref{have2PP}), we have {\small $\f{
\ln ( \ze \de ) } { n_{\ell} } < \mscr{M}_{\mrm{P}} ( z, z - \vep_a
)$ and $\f{ \ln ( \ze \de ) } { n_{\ell} } < \mscr{M}_{\mrm{P}} ( z,
\f{z}{1 - \vep_r} )$} for any $z \in [ \lm^\star - \vep_a, \lm^\star
+ \vep_a ]$ if $\vep_a > 0$ is small enough. This proves
(\ref{thatP0}).

To show (\ref{thatP}), let {\small $\om \in \{ \mscr{M}_{\mrm{P}} (
\wh{\bs{\lm}}_{\ell}, \wh{\bs{\lm}}_{\ell}  + \vep_a  )  > \f{\ln (\ze
\de)} {n_{\ell}}, \; \wh{\bs{\lm}}_{\ell} < \lm^\star - \vep_a \}$}
and $\wh{\lm}_{\ell} = \wh{\bs{\lm}}_{\ell} (\om)$. Then,
$\mscr{M}_{\mrm{P}} ( \wh{\lm}_{\ell},  \wh{\lm}_{\ell}   + \vep_a )
> \f{\ln (\ze \de)} {n_{\ell}}$ and $\wh{\lm}_{\ell} < \lm^\star -
\vep_a$. Since $z_{\ell} \in [0, \lm^\star - \vep_a)$ and
$\mscr{M}_{\mrm{P}} ( z, z + \vep_a  )$ is monotonically increasing
with respect to $z \in (0, \lm^\star - \vep_a)$, it must be true
that $\wh{\lm}_{\ell} > z_{\ell}$. Otherwise if $\wh{\lm}_{\ell}
\leq z_{\ell}$, then $\mscr{M}_{\mrm{P}} ( \wh{\lm}_{\ell},
\wh{\lm}_{\ell}   + \vep_a  ) \leq \mscr{M}_{\mrm{P}}  ( z_{\ell},
z_{\ell}   + \vep_a  ) = \f{\ln (\ze \de)} {n_{\ell}}$, leading to a
contradiction. This proves {\small $ \{ \mscr{M}_{\mrm{P}} (
\wh{\bs{\lm}}_{\ell},  \wh{\bs{\lm}}_{\ell} + \vep_a  ) > \f{\ln
(\ze \de)} {n_{\ell}} , \; \wh{\bs{\lm}}_{\ell} < \lm^\star - \vep_a
\} \subseteq  \{ z_{\ell} < \wh{\bs{\lm}}_{\ell} < \lm^\star -
\vep_a \}$}.  Now let $\om \in \{ z_{\ell} < \wh{\bs{\lm}}_{\ell} <
\lm^\star - \vep_a \}$ and $\wh{\lm}_{\ell} = \wh{\bs{\lm}}_{\ell}
(\om)$.  Then, $z_{\ell} < \wh{\lm}_{\ell} <
\lm^\star - \vep_a$.  Noting that $\mscr{M}_{\mrm{P}}  ( z, z +  \vep_a  )$ is
monotonically increasing with respect to $z  > 0$, we have that {\small $\mscr{M}_{\mrm{P}} (
\wh{\lm}_{\ell},  \wh{\lm}_{\ell}   + \vep_a ) > \mscr{M}_{\mrm{P}} (
z_{\ell},  z_{\ell} + \vep_a ) = \f{\ln (\ze \de)} {n_{\ell}}$},
which implies {\small $\{ \mscr{M}_{\mrm{P}} ( \wh{\bs{\lm}}_{\ell},
\wh{\bs{\lm}}_{\ell} + \vep_a ) > \f{\ln (\ze \de)} {n_{\ell}} , \;
\wh{\bs{\lm}}_{\ell} < \lm^\star - \vep_a \} \supseteq \{ z_{\ell} <
\wh{\bs{\lm}}_{\ell} < \lm^\star - \vep_a \}$}.  This establishes
(\ref{thatP}).

To show (\ref{that3P}), let {\small $\om \in \{ \mscr{M}_{\mrm{P}} (
\wh{\bs{\lm}}_{\ell}, \f{\wh{\bs{\lm}}_\ell}{ 1  + \vep_r}  )
> \f{\ln (\ze \de)} {n_{\ell}}, \; \wh{\bs{\lm}}_{\ell} > \lm^\star +
\vep_a \}$} and $\wh{\lm}_{\ell} = \wh{\bs{\lm}}_{\ell} (\om)$.
Then, {\small $\mscr{M}_{\mrm{P}} ( \wh{\lm}_{\ell}, \f{
\wh{\lm}_{\ell}}{ 1 + \vep_r} ) > \f{\ln (\ze \de)} {n_{\ell}}$} and
$\wh{\lm}_{\ell} > \lm^\star + \vep_a$. Since $y_{\ell} \in
(\lm^\star + \vep_a, \iy)$ and {\small $\mscr{M}_{\mrm{P}} ( z, \f{
z }{1 + \vep_r} )$} is monotonically decreasing with respect to $z > 0$, it must be true that
$\wh{\lm}_{\ell} < y_{\ell}$. Otherwise if $\wh{\lm}_{\ell} \geq
y_{\ell}$, then {\small $\mscr{M}_{\mrm{P}} ( \wh{\lm}_{\ell}, \f{
\wh{\lm}_{\ell} }{1 + \vep_r}  ) \leq \mscr{M}_{\mrm{P}} ( y_{\ell},
\f{ y_{\ell} }{1 + \vep_r} ) = \f{\ln (\ze \de)} {n_{\ell}}$},
leading to a contradiction. This proves {\small $ \{
\mscr{M}_{\mrm{P}}  ( \wh{\bs{\lm}}_{\ell}, \f{ \wh{\bs{\lm}}_{\ell}
}{1 + \vep_r} ) > \f{\ln (\ze \de)} {n_{\ell}} , \;
\wh{\bs{\lm}}_{\ell} > \lm^\star + \vep_a \} \subseteq \{ \lm^\star
+ \vep_a < \wh{\bs{\lm}}_{\ell} < y_\ell \}$}. Now let $\om \in \{
\lm^\star + \vep_a < \wh{\bs{\lm}}_{\ell} < y_{\ell} \}$ and
$\wh{\lm}_{\ell} = \wh{\bs{\lm}}_{\ell} (\om)$.
Then, $\lm^\star + \vep_a < \wh{\lm}_{\ell} < y_{\ell}$. Noting that {\small
$\mscr{M}_{\mrm{P}} ( z, \f{z}{1 +  \vep_r} )$} is monotonically
decreasing with respect to $z > 0$, we
have that {\small $\mscr{M}_{\mrm{P}} ( \wh{\lm}_{\ell}, \f{
\wh{\lm}_{\ell} }{ 1 + \vep_r} ) > \mscr{M}_{\mrm{P}} ( y_{\ell},
\f{ y_{\ell} }{ 1 + \vep_r} ) = \f{\ln (\ze \de)} {n_{\ell}}$},
which implies {\small $\{ \mscr{M}_{\mrm{P}} ( \wh{\bs{\lm}}_{\ell}, \f{
\wh{\bs{\lm}}_{\ell} }{ 1 + \vep_r} )
> \f{\ln (\ze \de)} {n_{\ell}} , \; \wh{\bs{\lm}}_{\ell} > \lm^\star +
\vep_a \} \supseteq  \{  \lm^\star + \vep_a < \wh{\bs{\lm}}_{\ell} <
y_{\ell} \}$}. This establishes (\ref{that3P}).

\beL \la{DefmixP} Let $\ell_\vep = s - j_\lm$.  Then, under the
constraint that limits are taken with $\f{\vep_a}{\vep_r}$ fixed,
\be \la{eqmix} \lim_{\vep_a \to 0} \sum_{\ell = 1}^{\ell_\vep - 1}
n_\ell \Pr \{ \bs{D}_\ell = 1 \} = 0, \qqu \lim_{\vep_a \to 0}
\sum_{\ell = \ell_\vep + 1}^s n_\ell \Pr \{ \bs{D}_{\ell} = 0 \} = 0
\ee for $\lm \in (0, \iy)$.  Moreover, $\lim_{\vep_a \to 0}
n_{\ell_\vep} \Pr \{ \bs{D}_{\ell_\vep} = 0 \} = 0$ if $C_{j_\lm} >
r (\lm)$. \eeL

\bpf

Throughout the proof of the lemma, we restrict $\vep_a$ to be
small enough such that {\small $\f{ \ln \f{1}{\ze \de} } { \vep_a }
< \f{ \ln (\ze \de) } { \mscr{M}_{\mrm{P}} \li ( \lm, \f{ \lm }{1 +
\vep_r} \ri ) }$}.  For simplicity of notations, let $a_\ell = \lim_{\vep_a \to 0}
y_\ell$ and $b_\ell = \lim_{\vep_a \to 0} z_\ell$.  The proof
consists of three main steps as follows.

\bsk

First, we shall show that (\ref{eqmix}) holds for $\lm \in (0, \lm^\star]$.
By the definition of $\ell_\vep$, we have
 {\small $\f{\lm} { \lm^\star } > C_{s - \ell_\vep + 1}$}.
 Making use of the first four statements of Lemma \ref{lem98PS}, we have that
 {\small $z_\ell < \f{ \lm + b_{\ell_\vep -
1} }{2} < \lm$} for all $\ell \leq \ell_\vep - 1$ and {\small $y_{s
- 1} > \f{ \lm^\star + a_{s - 1} }{2} > \lm^\star$} if $\vep_a$ is
sufficiently small. By the
last statement of Lemma \ref{lem98PS} and using Lemma \ref{PCH},
we have \bee \Pr \{ \bs{D}_{\ell} = 1 \} & = & \Pr
\{ \wh{\bs{\lm}}_{\ell} \leq z_{\ell}
\} + \Pr \{ \wh{\bs{\lm}}_{\ell} \geq y_{\ell} \} \leq  \Pr \{ \wh{\bs{\lm}}_{\ell} \leq z_\ell  \} + \Pr \li
\{ \wh{\bs{\lm}}_{\ell} \geq y_{s - 1} \ri \}\\
& \leq & \Pr \li \{ \wh{\bs{\lm}}_{\ell} \leq \f{ \lm + b_{\ell_\vep
- 1} }{2} \ri \} + \Pr \li \{ \wh{\bs{\lm}}_{\ell} \geq \f{
\lm^\star + a_{s - 1} }{2} \ri \}\\
& \leq & \exp \li ( n_{\ell} \mscr{M}_{\mrm{P}} \li ( \f{ \lm +
b_{\ell_\vep - 1} }{2}, \lm \ri ) \ri ) +  \exp \li ( n_{\ell}
\mscr{M}_{\mrm{P}} \li (\f{ \lm^\star + a_{s - 1} }{2}, \lm \ri )
\ri ) \eee for all $\ell \leq \ell_\vep - 1$ if $\vep_a > 0$ is
small enough.  Noting that {\small $b_{\ell_\vep - 1} = \lm^\star
C_{j_\lm + 1}$}, $a_{s - 1} = \f{\lm^\star}{C_1}$,
\[
\f{ \lm + b_{\ell_\vep - 1} }{2} = \f{ \lm + \lm^\star C_{j_\lm + 1}
}{2} < \lm, \qqu \f{ \lm^\star + a_{s - 1} }{2} = \f{ \lm^\star +
\f{\lm^\star}{C_1} }{2} > \lm
\]
which are constants independent of $\vep_a > 0$. Therefore, both
{\small $\mscr{M}_{\mrm{P}} ( \f{ \lm + b_{\ell_\vep - 1} }{2}, \lm
)$} and {\small $\mscr{M}_{\mrm{P}} (\f{ \lm^\star + a_{s - 1} }{2},
\lm )$} are negative constants independent of $\vep_a > 0$.  It
follows from  Lemma \ref{lem31a} that {\small $\lim_{\vep_a \to 0}
\sum_{\ell = 1}^{\ell_\vep - 1} n_\ell \Pr \{ \bs{D}_\ell = 1 \}  =
0$}.

Similarly, it can be seen from the definition of $\ell_\vep$ that
 $\f{\lm} { \lm^\star }  < C_{s - \ell_\vep - 1}$.  Making use of the first four statements of Lemma \ref{lem98PS},
we have that
 {\small $z_\ell > \f{\lm + b_{\ell_\vep + 1}}{2} > \lm$} for
$\ell_\vep + 1 \leq \ell < s$ if $\vep_a$ is sufficiently small. By
the last statement of Lemma \ref{lem98PS} and using Lemma \ref{PCH},
we have {\small \[ \Pr \{ \bs{D}_{\ell} = 0 \}  =  \Pr \{  z_{\ell}
< \wh{\bs{\lm}}_{\ell} < y_{\ell} \} \leq \Pr \{
\wh{\bs{\lm}}_{\ell}
> z_{\ell} \} \leq  \Pr \li \{ \wh{\bs{\lm}}_{\ell} > \f{\lm + b_{\ell_\vep + 1}}{2} \ri \} \leq \exp \li (
n_{\ell} \mscr{M}_{\mrm{P}} \li ( \f{\lm + b_{\ell_\vep + 1}}{2},
\lm \ri ) \ri ) \]} for $\ell_\vep + 1 \leq \ell < s$ if $\vep_a
> 0$ is small enough.  By virtue of the definition of $\ell_\vep$, we have that $b_{\ell_\vep + 1}$ is greater
than $\lm$ and is independent of $\vep_a
> 0$.  In view of this and the fact that $\Pr \{ \bs{D}_s = 0 \} = 0$,
 we can use Lemma \ref{lem31a} to arrive at $\lim_{\vep_a \to 0} \sum_{\ell = \ell_\vep + 1}^s n_\ell \Pr \{
\bs{D}_\ell = 0 \} = 0$.

\bsk

Second, we shall show that (\ref{eqmix}) holds for $\lm \in
(\lm^\star, \iy)$. As a direct consequence of the definition of
$\ell_\vep$, we have $\f{\lm^\star} { \lm } > C_{s - \ell_\vep +
1}$. Making use of the first four statements of Lemma \ref{lem98PS},
we have that $y_\ell > \f{\lm + a_{\ell_\vep - 1}}{2}
> \lm$ for all $\ell \leq \ell_\vep - 1$ and {\small $z_{s-1} < \f{
\lm^\star + b_{s - 1}} {2}$} if $\vep_a$ is sufficiently small. By
the last statement of Lemma \ref{lem98PS} and using Lemma \ref{PCH},
we have \bee \Pr \{ \bs{D}_{\ell} = 1 \} & = & \Pr \{
\wh{\bs{\lm}}_{\ell} \geq y_{\ell} \} + \Pr \{ \wh{\bs{\lm}}_{\ell}
\leq z_\ell \}
\leq  \Pr \{ \wh{\bs{\lm}}_{\ell} \geq y_{\ell} \} + \Pr \{ \wh{\bs{\lm}}_{\ell} \leq z_{s - 1} \}\\
& \leq & \Pr \li \{ \wh{\bs{\lm}}_{\ell} \geq \f{\lm + a_{\ell_\vep - 1}}{2} \ri \} + \Pr
\li \{ \wh{\bs{\lm}}_{\ell} \leq \f{ \lm^\star + b_{s - 1}} {2} \ri \}\\
& \leq & \exp \li ( n_{\ell} \mscr{M}_{\mrm{P}} \li ( \f{\lm + a_{\ell_\vep - 1}}{2},
\lm \ri ) \ri ) + \exp \li ( n_{\ell} \mscr{M}_{\mrm{P}} \li ( \f{ \lm^\star + b_{s - 1}} {2}, \lm \ri ) \ri )
 \eee for all $\ell \leq \ell_\vep - 1$ if $\vep_a > 0$ is small enough.
 By virtue of the definition of $\ell_\vep$, we have that $a_{\ell_\vep - 1}$ is greater
than $\lm$ and is independent of $\vep_a > 0$.  Hence, it follows
from Lemma \ref{lem31a} that $\lim_{\vep_a \to 0} \sum_{\ell =
1}^{\ell_\vep - 1} n_\ell \Pr \{ \bs{D}_\ell = 1 \}  = 0$.

In a similar manner, by the definition of $\ell_\vep$, we have
$\f{\lm^\star} { \lm } < C_{\ell_\vep - 1}$.  Making use of the
first four statements of Lemma \ref{lem98PS}, we have that $y_\ell <
\f{ \lm + a_{ \ell_\vep + 1 } }{2} < \lm$ for $\ell_\vep + 1 \leq
\ell < s$ if $\vep_a$ is sufficiently small. By the last statement
of Lemma \ref{lem98PS} and using Lemma \ref{PCH}, we have {\small \[
 \Pr \{ \bs{D}_{\ell} = 0
\}  =  \Pr \{  z_{\ell} < \wh{\bs{\lm}}_{\ell} < y_{\ell} \} \leq
\Pr \{ \wh{\bs{\lm}}_{\ell} < y_{\ell} \} \leq \Pr \li \{
\wh{\bs{\lm}}_{\ell} < \f{ \lm + a_{ \ell_\vep + 1 } }{2} \ri \}
 \leq  \exp \li ( n_{\ell} \mscr{M}_{\mrm{P}} \li (\f{ \lm + a_{ \ell_\vep + 1 } }{2}, \lm \ri) \ri )
\]} for $\ell_\vep + 1 \leq \ell < s$ if $\vep
> 0$ is small enough.  As a result of the definition of $\ell_\vep$, we have that $a_{\ell_\vep + 1}$ is smaller
than $\lm$ and is independent of $\vep_a
> 0$.  In view of this and the fact that $\Pr \{ \bs{D}_s = 0 \} = 0$,
 we can use Lemma \ref{lem31a} to conclude that $\lim_{\vep_a \to 0} \sum_{\ell = \ell_\vep + 1}^s n_\ell \Pr \{
\bs{D}_\ell = 0 \} = 0$.  This proves that (\ref{eqmix}) holds for $\lm \in (\lm^\star, \iy)$.

\bsk

Third,  we shall show that $\lim_{\vep \to 0} n_{\ell_\vep} \Pr \{
\bs{D}_{\ell_\vep} = 0 \} = 0$ if $C_{j_\lm} > r(\lm)$.

For $\lm \in (0, \lm^\star)$ such that $C_{j_\lm} > r(\lm)$, we have
$\f{ \lm }{\lm^\star} < C_{s - \ell_\vep}$ because of the definition
of $\ell_\vep$.  Making use of the first four statements of Lemma
\ref{lem98PS}, we have that $z_{\ell_\vep}
> \f{\lm  + b_{\ell_\vep}}{2}  > \lm$ if $\vep_a
> 0$ is small enough.  By the last statement of Lemma \ref{lem98PS} and using Lemma \ref{PCH}, we have

{\small $\Pr
\{ \bs{D}_{\ell_\vep} = 0 \}  =  \Pr \{  z_{\ell_\vep} <
\wh{\bs{\lm}}_{\ell_\vep} < y_{\ell_\vep} \} \leq \Pr \{
 \wh{\bs{\lm}}_{\ell_\vep} > z_{\ell_\vep} \} \leq
  \Pr \li \{ \wh{\bs{\lm}}_{\ell_\vep} > \f{\lm  + b_{\ell_\vep}}{2} \ri \}  \leq \exp \li (
n_{\ell_\vep} \mscr{M}_{\mrm{P}} \li (\f{\lm  + b_{\ell_\vep}}{2}, \lm \ri ) \ri )$}.
Since $b_{\ell_\vep}$ is greater than $\lm$ and is independent of $\vep_a > 0$
due to the definition of $\ell_\vep$, it follows that  $\lim_{\vep_a \to 0} n_{\ell_\vep} \Pr \{ \bs{D}_{\ell_\vep} =
0 \} = 0$.

For $\lm \in (\lm^\star, \iy)$ such that $C_{j_\lm} > r(\lm)$, we
have $\f{ \lm^\star }{\lm} < C_{s - \ell_\vep}$ as a result of the
definition of $\ell_\vep$. Making use of the first four statements
of Lemma \ref{lem98PS}, we have that $y_{\ell_\vep} < \f{\lm +
a_{\ell_\vep}}{2} < \lm$ if $\vep_a
> 0$ is small enough.  By the last statement of Lemma \ref{lem98PS} and using Lemma \ref{PCH}, we have

{\small $\Pr \{ \bs{D}_{\ell_\vep} = 0 \}  =
 \Pr \{  z_{\ell_\vep} < \wh{\bs{\lm}}_{\ell_\vep} < y_{\ell_\vep}
\} \leq \Pr \{ \wh{\bs{\lm}}_{\ell_\vep} < y_{\ell_\vep}  \}
\leq  \Pr \li \{ \wh{\bs{\lm}}_{\ell_\vep}  < \f{\lm + a_{\ell_\vep}}{2} \ri  \} \leq
\exp \li ( n_{\ell_\vep} \mscr{M}_{\mrm{P}} \li ( \f{\lm + a_{\ell_\vep}}{2}, \lm \ri ) \ri )$}.
Since $a_{\ell_\vep}$ is smaller than $\lm$ and is independent of $\vep_a > 0$
as a consequence of the definition of $\ell_\vep$,
it follows that $\lim_{\vep_a \to 0} n_{\ell_\vep} \Pr \{ \bs{D}_{\ell_\vep} = 0 \} = 0$.
This concludes the proof of the lemma.

\epf

\bsk

Finally, we would like to note that the proof of Theorem
\ref{Pos_mix_DDV_Asp} can be completed by employing Lemma
\ref{DefmixP} and a similar argument as that of Theorem
\ref{Bino_DDV_Asp}.

\subsection{Proof of Theorem \ref{Pos_mix_Asp_Analysis} }  \la{App_Pos_mix_Asp_Analysis}

As a result of the definitions of $\ka_\lm$ and $r(\lm)$,
 we have that $\ka_\lm > 1$ if and only if $C_{j_\lm} > r(\lm)$.
 To prove Theorem \ref{Pos_mix_Asp_Analysis}, we need some preliminary results.

\beL

\la{lem100PP} $\lim_{\vep_a \to 0} \f{ n_{\ell_\vep}
}{\mcal{N}_{\mrm{m}} (\lm, \vep_a, \vep_r)} = \ka_\lm, \;
\lim_{\vep_a \to 0} \vep_a \sq{ \f{n_{\ell_\vep}}{\lm} } = d
\sq{\ka_\lm}, \; \lim_{\vep_r \to 0} \vep_r \sq{\lm n_{\ell_\vep}} =
d \sq{\ka_\lm}$. \eeL

\bpf

First, we shall consider $\lm \in (0, \lm^\star)$.  Note that \bee
\mscr{M}_{\mrm{P}} (z, z  + \vep )  =  - \vep  + z \ln \li ( 1 +
\f{\vep}{ z } \ri ) =  - \vep + z \li [ \f{\vep}{z} - \f{\vep^2}{2
z^2} +  o (\vep^2) \ri ] =  - \f{\vep^2 }{2 z} + o (\vep^2).  \eee
By the definition of sample sizes, we have \be \la{gook}
\lim_{\vep_a \to 0} \f{ C_{s - \ell} \ln (\ze \de) } { n_\ell
\mscr{M}_{\mrm{P}} (\lm^\star  + \vep_a, \lm^\star)} = 1 \ee for  $1
\leq \ell < s$. It follows that  \bee \lim_{\vep_a \to 0} \f{
n_{\ell_\vep} }{\mcal{N}_{\mrm{m}} (\lm, \vep_a, \vep_r)} & = &
\lim_{\vep_a \to 0} \f{ \mscr{M}_{\mrm{P}} (\lm, \lm + \vep_a )  } {
\ln (\ze \de) } \times \f{ C_{s - \ell_\vep} \ln (\ze \de) }
{\mscr{M}_{\mrm{P}} (\lm^\star  + \vep_a, \lm^\star)} = \lim_{\vep_a
\to 0} \f{ C_{s - \ell_\vep} \mscr{M}_{\mrm{P}} (\lm, \lm + \vep_a )
} {\mscr{M}_{\mrm{P}} (\lm^\star  +
\vep_a, \lm^\star )}\\
& = & \lim_{\vep_a \to 0} \f{ C_{s - \ell_\vep} [ - \f{\vep_a^2 }{2
\lm} +  o (\vep_a^2) ]  } {- \f{\vep_a^2 }{2 \lm^\star} + o
(\vep_a^2)} = \f{\lm^\star}{\lm} C_{s - \ell_\vep} =
\f{\lm^\star}{\lm} C_{j_\lm} = \ka_\lm \eee and {\small \bee
\lim_{\vep_a \to 0} \vep_a \sq{ \f{n_{\ell_\vep}}{\lm} }  & = &
\lim_{\vep_a \to 0} \vep_a \sq{\f{1}{\lm}  \f{ C_{s - \ell_\vep} \ln
(\ze \de) } {\mscr{M}_{\mrm{P}} ( \lm^\star +
\vep_a, \lm^\star )} }\\
& = & \lim_{\vep_a \to 0} \vep_a \sq{\f{1}{\lm} \f{ C_{s -
\ell_\vep} \ln (\ze \de) } {- \f{\vep_a^2 }{2 \lm^\star} + o
(\vep_a^2) } } = d \sq{ \f{\lm^\star}{\lm} C_{s - \ell_\vep} } = d
\sq{\ka_\lm}. \eee}

We shall next consider $\lm \in (\lm^\star, \iy)$.
Note that {\small \bee \mscr{M}_{\mrm{P}} \li (z, \f{z}{1 + \vep} \ri )  =
 \f{\vep z}{1 + \vep} - z \ln (1 + \vep)  =  \vep z \li [ 1 - \vep
+ o(\vep) \ri ] - z \li [ \vep - \f{\vep^2}{2} + o (\vep^2) \ri ]  =
- \f{\vep^2 z}{2} + o (\vep^2).  \eee} By (\ref{gook}), we have
{\bee  \lim_{\vep_r \to 0} \f{ n_{\ell_\vep} } { \mcal{N}_{\mrm{m}}
(\lm, \vep_a, \vep_r) } & = & \lim_{\vep_r \to 0} \f{
\mscr{M}_{\mrm{P}} (\lm, \f{\lm}{1 + \vep_r} ) }{ \ln (\ze \de) }
\f{ C_{s - \ell_\vep} \ln (\ze \de) }
{\mscr{M}_{\mrm{P}} ( \lm^\star + \vep_a, \lm^\star)}\\
& = & \lim_{\vep_r \to 0} \f{ C_{s - \ell_\vep} \mscr{M}_{\mrm{P}}
(\lm, \f{\lm}{1 + \vep_r} ) } {\mscr{M}_{\mrm{P}} ( \lm^\star +
\vep_a, \lm^\star)}
 =  \lim_{\vep_r \to 0} \f{ C_{s - \ell_\vep} [ - \f{\vep_r^2
\lm}{2} + o (\vep_r^2)  ]  } {  - \f{\vep_a^2 }{2 \lm^\star} + o
(\vep_a^2)  }\\
& = & \f{\lm}{ \lm^\star } C_{s - \ell_\vep}
 = \f{\lm}{ \lm^\star } C_{j_\lm} = \ka_\lm \eee} and \bee \lim_{\vep_r \to 0}
\vep_r \sq{\lm n_{\ell_\vep}} & = & \lim_{\vep_r \to 0} \vep_r \sq{
\f{ \lm C_{s - \ell_\vep} \ln (\ze \de) }
{\mscr{M}_{\mrm{P}} ( \lm^\star + \vep_a, \lm^\star)} } \\
& = & \lim_{\vep_r \to 0} \vep_r \sq{ \f{ \lm C_{s - \ell_\vep} \ln
(\ze \de) } {- \f{\vep_a^2}{2 \lm^\star} + o (\vep_a^2) } } = d \sq{
\f{\lm}{ \lm^\star } C_{s - \ell_\vep}  } = d \sq{\ka_\lm}. \eee

\epf

\beL \la{limplempos} Let $U$ and $V$ be independent Gaussian random
variables with zero means and unit variances.  Then, for $\lm \in
(0, \iy)$ such that $C_{j_\lm} = r (\lm)$ and $j_\lm \geq 1$, \bee &
& \lim_{\vep \to 0} \Pr \{ \bs{l} = \ell_\vep \} = 1 - \lim_{\vep
\to 0} \Pr \{ \bs{l} = \ell_\vep + 1 \} =  1 - \Phi \li (  \nu  d
\ri ), \\
&  & \lim_{\vep \to 0} \li [ \Pr \{ |
\wh{\bs{\lm}}_{\ell_\vep} - \lm | \geq \vep_\lm, \; \bs{l}  =
\ell_\vep \} + \Pr \{ | \wh{\bs{\lm}}_{\ell_\vep + 1} - \lm | \geq
\vep_\lm, \;
\bs{l} =  \ell_\vep + 1 \} \ri ]\\
&   & \qqu \qqu \qqu \qqu = \Pr \li \{ U \geq d \ri \} + \Pr \li \{ |U + \sq{\ro_\lm} V | \geq (1 + \ro_\lm) d, \; U < \nu d \ri \}, \eee where
$\vep_\lm = \max \{ \vep_a, \vep_r \lm \}$.

\eeL

\bpf

We shall first consider $\lm \in (\lm^\star, \iy)$  such that
$C_{j_\lm} = r (\lm)$.  Since $\ka_\lm = 1$, by Statement (V) of
Lemma \ref{lem98PS}, we have
\[
\lim_{\vep_r \to 0} \f{ z_{\ell_\vep} - \lm}{ \sq{ \lm \sh n_{\ell_\vep}} } = \lim_{\vep_r \to 0} \vep_r \sq{\lm n_{\ell_\vep}} \lim_{\vep_r \to
0} \f{ z_{\ell_\vep} - \lm }{\vep_r \lm} = d \lim_{\vep_r \to 0} \f{ z_{\ell_\vep} - \lm}{\vep_r \lm} = d.
\]
By a similar argument as in the proof of Lemma \ref{limplem}, we can
show that
\[
\lim_{\vep \to 0} \Pr \{ \bs{l} = {\ell_\vep} \} = 1 - \lim_{\vep
\to 0} \Pr \{ \bs{l} = {\ell_\vep} + 1 \} = \lim_{\vep \to 0} \Pr \{
\wh{\bs{\lm}}_{\ell_\vep} \geq z_{\ell_\vep} \}
\]
\bee &  & \lim_{\vep \to 0} \li [ \Pr \{ |
\wh{\bs{\lm}}_{{\ell_\vep}} - \lm | \geq \vep_\lm, \; \bs{l}  =
{\ell_\vep} \} + \Pr \{ | \wh{\bs{\lm}}_{{\ell_\vep} + 1} - \lm |
\geq \vep_\lm, \;
\bs{l} =  {\ell_\vep} + 1 \} \ri ]\\
& = & \lim_{\vep \to 0} \li [ \Pr \{ | \wh{\bs{\lm}}_{\ell_\vep} - \lm | \geq \vep_r \lm, \; \wh{\bs{\lm}}_{\ell_\vep} \geq z_{\ell_\vep} \} +
\Pr \{ | \wh{\bs{\lm}}_{{\ell_\vep} + 1} - \lm | \geq \vep_r \lm, \; \wh{\bs{\lm}}_{\ell_\vep} < z_{\ell_\vep} \} \ri ]. \eee Note that
\[ \Pr \{ | \wh{\bs{\lm}}_{\ell_\vep} - \lm | \geq \vep_r \lm, \;
\wh{\bs{\lm}}_{\ell_\vep} \geq z_{\ell_\vep} \} = \Pr \li \{ \f{ |
\wh{\bs{\lm}}_{\ell_\vep} - \lm | }{ \sq{ \lm \sh n_{\ell_\vep}} }
\geq \vep_r \sq{\lm n_{\ell_\vep}}, \;    \f{
\wh{\bs{\lm}}_{\ell_\vep} - \lm }{ \sq{ \lm \sh n_{\ell_\vep}} }
\geq \f{ z_{\ell_\vep} - \lm }{ \sq{ \lm \sh n_{\ell_\vep}} } \ri
\}.
\]
Therefore, \bee &  & \Pr \{ | \wh{\bs{\lm}}_{\ell_\vep} - \lm | \geq
\vep_\lm, \; \bs{l} = \ell_\vep \} + \Pr \{ |
\wh{\bs{\lm}}_{\ell_\vep + 1} - \lm | \geq \vep_\lm, \; \bs{l} =
\ell_\vep + 1 \}\\
& \to & \Pr \{  |U| \geq d, \; U \geq d \} + \Pr \li \{  \li | U + \sq{\ro_\lm} V \ri | \geq (1 + \ro_\lm) d, \; U < d \ri \} \\
&  = & \Pr \{  U \geq d \} + \Pr \li \{  \li | U + \sq{\ro_\lm} V \ri | \geq (1 + \ro_\lm) d, \; U < d \ri \} \eee for $\lm \in (\lm^\star,
\iy)$ such that $C_{j_\lm} = r (\lm)$.

Next, we shall now consider $\lm \in (0, \lm^\star)$ such that
$C_{j_\lm} = r (\lm)$.  Since $\ka_\lm = 1$, by Statement (V) of
Lemma \ref{lem98PS}, we have
\[
\lim_{\vep_a \to 0} \f{ z_{\ell_\vep} - \lm}{ \sq{ \lm \sh
n_{\ell_\vep}} } = \lim_{\vep_a \to 0} \vep_a \sq{ \f{ n_{\ell_\vep}
}{ \lm } } \lim_{\vep_a \to 0} \f{ z_{\ell_\vep} - \lm }{\vep_a} = d
\lim_{\vep_a \to 0} \f{ z_{\ell_\vep} - \lm}{\vep_a} = - \nu d.
\]
Clearly,
\[ \Pr \{ | \wh{\bs{\lm}}_{\ell_\vep} - \lm | \geq \vep_a, \;
\wh{\bs{\lm}}_{\ell_\vep} \leq z_{\ell_\vep} \} = \Pr \li \{ \f{ |
\wh{\bs{\lm}}_{\ell_\vep} - \lm | }{ \sq{ \lm \sh n_{\ell_\vep}} }
\geq \vep_a \sq{ \f{n_{\ell_\vep} }{ \lm } }, \;    \f{
\wh{\bs{\lm}}_{\ell_\vep} - \lm }{ \sq{ \lm \sh n_{\ell_\vep}} }
\leq \f{ z_{\ell_\vep} - \lm }{ \sq{ \lm \sh n_{\ell_\vep}} } \ri
\}.
\]
Therefore, \bee &  & \Pr \{ | \wh{\bs{\lm}}_{\ell_\vep} - \lm | \geq
\vep_\lm, \; \bs{l} = \ell_\vep \} + \Pr \{ |
\wh{\bs{\lm}}_{\ell_\vep + 1} - \lm | \geq \vep_\lm, \; \bs{l} =
\ell_\vep + 1 \}\\
 & \to & \Pr \{ |U| \geq d, \; U \leq - \nu d \} +
\Pr \li \{ \li | U + \sq{\ro_\lm} V \ri | \geq (1 + \ro_\lm) d, \; U >
- \nu d \ri \}\\
& = & \Pr \li \{ U \geq d \ri \} + \Pr \li \{ \li |  U + \sq{\ro_\lm} V \ri | \geq (1 + \ro_\lm) d, \; U < \nu d \ri \}. \eee

\epf

\bsk

Finally, we would like to note that the proof of Theorem
\ref{Pos_mix_Asp_Analysis} can be completed by employing Lemma
\ref{lem100PP} and similar arguments as that of Theorem
\ref{Bino_Asp_Analysis}. Specially, we need to restrict $\vep_a$ to
be small enough such that {\small $\f{ \ln \f{1}{\ze \de} } { \vep_a
} < \f{ \ln (\ze \de) } { \mscr{M}_{\mrm{P}} \li ( \lm, \f{ \lm }{1
+ \vep_r} \ri ) }$}.  For the purpose of proving Statement (III), we
need to make use of the following observation:
\[
\Pr \{ | \wh{\bs{\lm}} - \lm | \geq \vep_a, \; | \wh{\bs{\lm}} - \lm | \geq \vep_r \lm \} = \bec
\Pr \{ | \wh{\bs{\lm}} - \lm | \geq \vep_a \} & \tx{for} \; \lm \in (0, \lm^\star],\\
\Pr \{ | \wh{\bs{\lm}} - \lm | \geq \vep_r \lm \} & \tx{for} \; \lm \in (\lm^\star, \iy)
\eec
\]
\[ \Pr \{ | \wh{\bs{\lm}}_{\ell} - \lm | \geq \vep_a \} = \Pr \li \{  |U_\ell|
\geq \vep_a \sq{ \f{n_{\ell_\vep} } {\lm } } \ri \}, \qqu
\Pr \{ | \wh{\bs{\lm}}_{\ell} - \lm | \geq \vep_r \lm \} = \Pr \li \{  |U_\ell|
\geq \vep_r \sq{ \lm n_{\ell}  } \ri \}
\]
where, according to the central limit theorem,  $U_\ell = \f{ |
\wh{\bs{\lm}}_{\ell} - \lm | }{ \sq{ \lm \sh n_{\ell}} }$ converges in distribution
to a Gaussian random variable $U$ of zero mean and unit variance as $\vep_a \to 0$.

\sect{Proofs of Theorems for Estimation of Normal Mean}

\subsection{Proof of Theorem \ref{Normal_Analytic_Thm} } \la{App_Normal_Analytic_Thm}

First, we shall show statement (I) which asserts that $\Pr \{ |
\wh{\bs{\mu}} - \mu | < \vep \} > 1 - 2 s \ze \de$.    Define
$\mbf{m} = \max \{ n_s, \; \lc ( \wh{\bs{\si}}_s  \; t _{n_{s} - 1,
\ze \de} )^2 \sh \vep^2 \rc \}$.  Then, $\{ \sq{ \mbf{m} } \geq (
\wh{\bs{\si}}_s  \; t _{n_{s} - 1, \ze \de} ) \sh \vep \}$ is a sure
event and by the definition of the sampling scheme,  \bel \Pr \{ |
\ovl{X}_{\mbf{n}} - \mu | \geq \vep, \; \mathbf{n} \geq n_s \} & = &
\Pr \{ | \ovl{X}_{\mbf{m}} - \mu | \geq \vep, \; \mathbf{n} \geq n_s
\} \leq \Pr \{ | \ovl{X}_{\mbf{m}}
- \mu | \geq \vep \} \nonumber\\
& = & \Pr \{ | \ovl{X}_{\mbf{m}} - \mu | \geq \vep, \; \sq{ \mbf{m}
} \geq ( \wh{\bs{\si}}_s  \; t _{n_{s} - 1, \ze \de} ) \sh
\vep \} \nonumber\\
&  \leq & \Pr \li \{ \sq{ \mbf{m} } | \ovl{X}_{\mbf{m}} - \mu | \geq
\vep \times \f{ \wh{\bs{\si}}_s \; t _{n_{s} - 1, \ze \de} }{
\vep } \ri \} \nonumber\\
& = & \Pr \li \{ \f{ \sq{\mbf{m}} | \ovl{X}_{\mbf{m}} - \mu |
}{\wh{\bs{\si}}_s } \geq t _{n_{s} - 1, \ze \de}  \ri \}.
\la{try6888} \eel Note that $\sq{\mbf{m}} ( \ovl{X}_{\mbf{m}} - \mu
) \sh \si$ is a standard Gaussian variable and that $\sq{\mbf{m}} (
\ovl{X}_{\mbf{m}} - \mu ) \sh \si$ is independent of
$\wh{\bs{\si}}_s$  because {\small \bee \Pr \li \{ \f{ \sq{\mbf{m}}
( \ovl{X}_{\mbf{m}} - \mu ) }{ \si } \leq u \ri \} & = & \sum_{m =
n_s}^\iy
 \Pr \li \{ \f{ \sq{m} ( \ovl{X}_{m} - \mu ) } { \si } \leq u, \;  \mbf{m} = m \ri \} \\
 & = & \sum_{m = n_s}^\iy \Pr \li \{ \f{ \sq{m} ( \ovl{X}_{m} - \mu )
}{\si} \leq u \ri \} \Pr \{ \mbf{m} = m  \} =  \sum_{m = n_s}^\iy
\Phi (u)  \Pr \{ \mbf{m} = m \} = \Phi (u) \qu \eee} and \bee \Pr
\li \{ \f{ \sq{\mbf{m}} ( \ovl{X}_{\mbf{m}} - \mu ) }{ \si } \leq u,
\; \wh{\bs{\si}}_s \leq v \ri \} & = & \sum_{m = n_s}^\iy \Pr \li \{
\f{ \sq{m} (
\ovl{X}_{m} - \mu ) }{ \si } \leq u, \;  \mbf{m} = m, \; \wh{\bs{\si}}_s \leq v \ri \}\\
& = & \sum_{m = n_s}^\iy \Pr \li \{ \f{ \sq{m} ( \ovl{X}_{m} - \mu )
}{ \si } \leq u \ri \} \Pr \{ \mbf{m} = m, \; \wh{\bs{\si}}_s \leq v \}\\
& = & \sum_{m = n_s}^\iy \Phi (u)  \Pr \{ \mbf{m} = m, \;
\wh{\bs{\si}}_s \leq v \}  =  \Phi (u) \Pr \{ \wh{\bs{\si}}_s \leq v \}\\
& = &  \Pr \{ \sq{\mbf{m}} ( \ovl{X}_{\mbf{m}} - \mu ) \sh \si \leq
u \} \Pr \{ \wh{\bs{\si}}_s \leq v \} \eee for any $u$ and $v$.
Therefore, $\sq{\mbf{m}} ( \ovl{X}_{\mbf{m}} - \mu ) \sh
\wh{\bs{\si}}_s$ has a Student $t$-distribution of $n_s - 1$ degrees
of freedom.  It follows from (\ref{try6888}) that \be \la{try68} \Pr
\{ | \ovl{X}_{\mbf{n}} - \mu | \geq \vep, \; \mathbf{n} \geq n_s \}
\leq  2 \ze \de. \ee By the definition of the sampling scheme, we
have $ \{ \mathbf{n} = n_\ell \} \subset \li \{ \vep \geq \f{
\wh{\bs{\si}}_\ell \; t _{n_{\ell} - 1, \ze \de} }{ \sq{n_\ell} }
\ri \}$  and thus {\small \be \Pr \{ | \ovl{X}_{\mbf{n}} - \mu |
\geq \vep, \; \mathbf{n} = n_\ell \} \leq \Pr \li \{ |
\ovl{X}_{n_\ell} - \mu | \geq \vep \geq \f{ \wh{\bs{\si}}_\ell \; t
_{n_{\ell} - 1, \ze \de} }{ \sq{n_\ell} } \ri \} \leq \Pr \li \{ \f{
\sq{n_\ell} | \ovl{X}_{n_\ell} - \mu | }{ \wh{\bs{\si}}_\ell } \geq
t _{n_{\ell} - 1, \ze \de}  \ri \} = 2 \ze \de \la{comb3398} \ee}
for $\ell = 1, \cd, s - 1$. Combining (\ref{try68}) and
(\ref{comb3398}) yields \be \la{nowuse} \Pr \{ | \wh{\bs{\mu}} - \mu
| \geq \vep \}  =  \Pr \{ | \ovl{X}_{\mbf{n}} - \mu | \geq \vep, \;
\mathbf{n} \geq n_s \} + \sum_{\ell = 1}^{s - 1} \Pr \{ |
\ovl{X}_{\mbf{n}} - \mu | \geq \vep, \; \mathbf{n} = n_\ell \} \leq
2 s \ze \de, \ee  which implies that $\Pr \{ | \wh{\bs{\mu}} - \mu |
< \vep \} > 1 - 2 s \ze \de$ for any $\mu$ and $\si$.  This proves
statement (I).

Second, we shall show statement (II) which asserts that $\lim_{\vep
\to 0} \Pr \{ | \wh{\bs{\mu}} - \mu | < \vep \} = 1 - 2 \ze \de$.
Obviously,  $\lim_{\vep \to 0} \Pr \{ \mathbf{n} < n_s \} = 0$.
Hence, $\lim_{\vep \to 0} \sum_{\ell = 1}^{s - 1} \Pr \{ |
\ovl{X}_{\mbf{n}} - \mu | \geq \vep, \; \mathbf{n} = n_\ell \} = 0$
and \bel \Pr \{ | \wh{\bs{\mu}} - \mu | \geq \vep \} & = & \Pr \{ |
\ovl{X}_{\mbf{n}} - \mu | \geq \vep, \; \mathbf{n} \geq n_s \} +
\sum_{\ell = 1}^{s - 1} \Pr \{ | \ovl{X}_{\mbf{n}} - \mu | \geq
\vep,  \; \mathbf{n} = n_\ell \} \nonumber\\
&  \to  &  \Pr \{ | \ovl{X}_{\mbf{n}} - \mu | \geq \vep, \;
\mathbf{n} \geq n_s \} \la{try98} \eel as $\vep \to 0$.  By virtue
of (\ref{try68}) and (\ref{try98}), we have $\limsup_{\vep \to 0}
\Pr \{ | \wh{\bs{\mu}} - \mu | < \vep \} \leq 2 \ze \de$, which
implies that \be \la{infinq} \liminf_{\vep \to 0} \Pr \{ |
\wh{\bs{\mu}} - \mu | \geq \vep \} \geq 1 - 2 \ze \de. \ee On the
other hand, by (\ref{try98}) and the fact that $\lim_{\vep \to 0}
\Pr \{ \mathbf{n} \geq n_s \} = 1$, we have {\small \bel \Pr \{ |
\wh{\bs{\mu}} - \mu | < \vep \} & \to & \Pr \{ | \ovl{X}_{\mbf{n}} -
\mu | < \vep, \; \mathbf{n} \geq n_s \} = \Pr \{ | \ovl{X}_{\mbf{m}}
- \mu | < \vep, \; \mathbf{n} \geq n_s \} \nonumber\\
&  \to  &  \Pr \{ | \ovl{X}_{\mbf{m}} - \mu | < \vep \}
\nonumber\\
& < & \Pr \li \{ | \ovl{X}_{\mbf{m} } - \mu | < \vep \leq \f{ (1 +
\eta) \wh{\bs{\si}}_s \; t _{n_{s} - 1, \ze \de} }{ \sq{\mbf{m} } }
\ri \} + \Pr \li \{ \f{ (1 + \eta) \wh{\bs{\si}}_s \; t _{n_{s} - 1,
\ze \de} }{ \sq{\mbf{m} } } < \vep
\ri \} \nonumber\\
& \leq & \Pr \li \{ \f{ \sq{\mbf{m} } | \ovl{X}_{\mbf{m} } - \mu | }
{ \wh{\bs{\si}}_s } < (1 + \eta)
 t _{n_{s} - 1, \ze \de} \ri \} + \Pr \li \{ \f{ (1 + \eta)
\wh{\bs{\si}}_s \; t _{n_{s} - 1, \ze \de} }{ \sq{\mbf{m} } } < \vep
\ri \} \nonumber \eel} as $\vep \to 0$,  where $\eta$ is a positive
number.  Noting that {\small \bee \Pr \li \{  \f{ (1 + \eta)
\wh{\bs{\si}}_s \; t _{n_{s} - 1, \ze \de} }{ \sq{\mbf{m} } } < \vep
\ri \} \leq  \Pr \li \{ \f{ (1 + \eta) \wh{\bs{\si}}_s \; t _{n_{s}
- 1, \ze \de} }{ \sq{ ( \wh{\bs{\si}}_s \; t _{n_{s} - 1, \ze \de}
)^2 \sh \vep^2 + n_s} } < \vep \ri \} =  \Pr \li \{
\wh{\bs{\si}}_s^2 < \f{ n_s \; \vep^2 } { \eta(2 + \eta)(t _{n_{s} -
1, \ze \de} )^2 }  \ri \} \eee}  which tends to $0$ as $\vep \to 0$,
we have
\[
\limsup_{\vep \to 0} \Pr \{ | \wh{\bs{\mu}} - \mu | < \vep \} \leq
\Pr \li \{ \f{ \sq{\mbf{m} } | \ovl{X}_{\mbf{m} } - \mu | } {
\wh{\bs{\si}}_s } < (1 + \eta)
 t _{n_{s} - 1, \ze \de} \ri \}.
\]
Since the above argument holds for arbitrarily small $\eta > 0$, we
have \be \la{supinq} \limsup_{\vep \to 0} \Pr \{ | \wh{\bs{\mu}} -
\mu | < \vep \} \leq \Pr \li \{ \f{ \sq{\mbf{m} } | \ovl{X}_{\mbf{m}
} - \mu | } { \wh{\bs{\si}}_s } \leq  t _{n_{s} - 1, \ze \de} \ri \}
= 1 - 2 \ze \de. \ee Combing (\ref{infinq}) and (\ref{supinq})
yields  $\lim_{\vep \to 0} \Pr \{ | \wh{\bs{\mu}} - \mu | < \vep \}
= 1 - 2 \ze \de$.  This proves statement (II).

\bsk

Finally, statements (III) and (IV) can be shown by making use of the
observation that $ \mbf{n} \leq ( \wh{\bs{\si}}_s  \; t _{n_{s} - 1,
\ze \de} )^2 \sh \vep^2  +  n_s$.   This completes the proof of
Theorem \ref{Normal_Analytic_Thm}.

\subsection{Proof of Theorem \ref{Normal_Main_Thm} } \la{App_Normal_Main_Thm}

\subsubsection{Proof of Statement (I)}

Define Helmert transform \be \la{Hel}
 U_i = \f{ X_i - \mu }{\si}, \qqu  V_i = \f{ U_1 + \cd +
U_i - i U_{i + 1}  } { \sq{i (i + 1) } }, \qqu  W_i = \f{ U_1 + \cd
+ U_{i}  } { \sq{i} } \ee for $i = 1, 2, \cd, \iy$. Clearly, the
$U_i$ are independent Gaussian variables with zero mean and variance
unity. Since the transformation from $(U_1, \cd, U_i)$ to $(V_1,
\cd, V_{i - 1}, W_i)$ is orthogonal for any $i \geq 2$, the $V_i$
are independent Gaussian variables with zero mean and variance
unity. It is easily seen that $\sq{n} ( \ovl{X}_{n} - \mu ) \sh \si
 = W_{n}$ and $S_{n} = \si^2 ( \sum_{i = 1}^{n} U_i^2 - W_{n}^2 ) = \si^2(
V_1^2 + \cd + V_{n - 1}^2 )$ for $n = 2, 3, \cd, \iy$. Hence, by the
definition of the sampling scheme, we have that $\{ | \ovl{X}_{n} -
\mu | \geq \vep \}$ is independent of $\{ \mbf{n} = n \}$ for any $n
\in \mscr{S}$.  It follows from  such independency and the
definition of the sampling scheme that {\small \bee   \Pr \{ |
\wh{\bs{\mu}} - \mu | \geq \vep \} & = & \sum_{ n \in \mscr{S} } \Pr
\{ | \wh{\bs{\mu}} - \mu | \geq \vep, \; \mbf{n} = n \} = \sum_{ n
\in \mscr{S} } \Pr \{ | \ovl{X}_{n} - \mu
| \geq \vep, \; \mbf{n} = n \},\\
& = & \sum_{ n \in \mscr{S} } \Pr \{ | \ovl{X}_{n} - \mu | \geq \vep
\} \Pr \{ \mbf{n} = n \} = 2 \sum_{ n \in \mscr{S} } \li [ 1 - \Phi
\li (  \f{ \vep \sq{n} } { \si }
 \ri ) \ri ] \Pr \{ \mbf{n} = n \}.
\eee} This proves statement (I).

\subsubsection{Proof of Statement (II)}

Define $Z_j = \f{ V_{2j - 1} + V_{2 j} } {2}$ for $j = 1, 2, \cd,
\iy$,  where $V_i$ are defined in (\ref{Hel}).  It is easy to see
that $Z_j$ are identical and independent exponential random
variables with density $e^{-z}$. By the definition of
$\wh{\bs{\si}}_\ell$, we have {\small $\wh{\bs{\si}}_\ell = \sq {
\f{ S_{2 k_\ell + 1} } { 2 k_\ell } } = \si \sq { \f{ \sum_{j =
1}^{k_\ell} Z_j } { k_\ell } }, \; \ell = 1, \cd, s$} and thus
{\small \bel &  & \li \{ \f{ ( \wh{\bs{\si}}_\ell \; t _{n_{\ell} -
1, \ze \de} )^2 } { \vep^2 } > n_\ell \ri \} = \li \{ \sum_{j =
1}^{k_\ell} Z_j
> b_\ell \ri \}, \;  \qqu \ell = 1, \cd, s, \la{qu1}\\
&   & \li \{  \f{ ( \wh{\bs{\si}}_s  \; t _{n_{s} - 1, \ze \de} )^2
}{ \vep^2}  > n \ri \} = \li \{  \sum_{j = 1}^{k_s} Z_j
> c \ri \}, \qqu n \geq n_s. \la{qu2}
\eel} It follows from (\ref{qu1}) and the definition of the stopping
rule that \be \la{key1}
 \{\mbf{n} > n_\ell \}  =  \li \{\sum_{j =
1}^{k_i} Z_j
> b_{i} \; \tx{for} \; 1 \leq i \leq \ell \ri \} \ee for $\ell = 1,
\cd, s$.  Making use of (\ref{key1}) and Theorem \ref{lemchen}, we
have \be \la{usa1} \Pr \{\mbf{n} > n_\ell \} =  H_\ell (\si) \ee for
$\ell = 1, \cd, s$.  Similarly, it follows from (\ref{qu2}) and the
definition of the stopping rule that \be \la{key2} \{\mbf{n} > n \}
= \li \{ \sum_{j = 1}^{k_s} Z_j
> c, \; \sum_{j = 1}^{k_\ell} Z_j
> b_\ell \; \tx{for} \; 1 \leq \ell < s \ri \}
\ee for $n \geq n_s$.  Making use of (\ref{key2}) and Theorem
\ref{lemchen}, we have \be \la{usa2}
 \Pr \{\mbf{n} > n \} =  H^\star (\si, n) \ee for $n \geq n_s$.  By
 virtue of (\ref{usa1}), we have $\Pr \{ \mbf{n} = n_1 \} = 1 - \Pr \{ \mbf{n} > n_1 \} = H_0(\si) -
H_1(\si)$ and $\Pr \{ \mbf{n} = n_\ell \} = \Pr \{ \mbf{n}
> n_{\ell - 1} \} - \Pr \{ \mbf{n} > n_\ell \} = H_{\ell - 1} (\si) - H_\ell (\si)$ for $1 <
\ell \leq s$.  In a similar manner, using (\ref{usa2}), we have $\Pr
\{ \mbf{n} = n \} = \Pr \{ \mbf{n} > n - 1 \} - \Pr \{ \mbf{n}
> n \} = H^\star (\si, n - 1) - H^\star (\si, n)$ for $n > n_s$.
This completes the proof of statement (II).

\subsubsection{Proof of Statement (III)}

By the established statement (I), we have \be \la{ineq169} \Pr \{ |
\wh{\bs{\mu}} - \mu | \geq \vep \}  =  2 \sum_{n \in \mscr{S} \atop{
n \leq m } } \li [ 1 - \Phi \li ( \f{\vep \sq{n} } { \si } \ri ) \ri
] \Pr \{ \mbf{n} = n \} + 2 \sum_{ n
> m  } \li [ 1 - \Phi \li ( \f{\vep
\sq{n} } { \si } \ri ) \ri ] \Pr \{ \mbf{n} = n \}. \ee Note that
{\small \bel \sum_{ n > m  } \li [ 1 - \Phi \li ( \f{\vep \sq{n} } {
\si } \ri ) \ri ] \Pr \{ \mbf{n} = n \} & < & \sum_{ n
> m  } \li [ 1 - \Phi \li ( \f{\vep
\sq{m} } { \si } \ri ) \ri ] \Pr \{ \mbf{n} = n \}  =  \li [ 1 -
\Phi \li ( \f{\vep \sq{m} } { \si } \ri ) \ri ] \Pr
\{ \mbf{n} > m \} \nonumber \\
& < & \li [ 1 -  \Phi \li (  \f{ \vep \sq{m} } { \si }
 \ri ) \ri ] \Pr \li \{  \chi_{n_s
- 1}^2 >  \f{ m (n_s - 1) \vep^2 } { ( \si t_{n_s - 1, \ze \de} )^2
} \ri \} \nonumber \\
& \leq & \li [ 1 -  \Phi \li (  \f{ \vep \sq{m} } { b }
 \ri ) \ri ] \Pr \li \{  \chi_{n_s
- 1}^2 >  \f{ m (n_s - 1) \vep^2 } { ( a \; t_{n_s - 1, \ze \de} )^2
} \ri \} \nonumber\\
& = &  \li [ 1 -  \Phi \li (  \f{ \vep \sq{m} } { b }
 \ri ) \ri ]  S_{\mrm{P}} \li (  k_s - 1, \f{ m k_s \vep^2 } { ( a \;
t_{n_s - 1, \ze \de} )^2 } \ri ) \la{ineq170} \eel} for any $\si \in
[ a, b]$, where $\chi_{n_s - 1}^2$ represents a chi-square random
variable of $n_s - 1$ degrees of freedom. Observing that $H_\ell
(\si)$ is monotonically increasing with respect to $\si \in [a, b]$
for $0 \leq \ell \leq s$ and that $H^\star (\si, n)$ is
monotonically increasing with respect to $\si \in [a, b]$ for $n
\geq n_s$, we have $\udl{P}_n \leq \Pr \{ \mbf{n} = n \} \leq
\ovl{P}_n$ for $\si \in [a, b]$. Therefore, \be \sum_{n \in \mscr{S}
\atop{ n \leq m } } \li [ 1 - \Phi \li ( \f{\vep \sq{n} } { a } \ri
) \ri ] \udl{P}_n \leq \sum_{n \in \mscr{S} \atop{ n \leq m } } \li
[ 1 - \Phi \li ( \f{\vep \sq{n} } { \si } \ri ) \ri ] \Pr \{ \mbf{n}
= n \} \leq \sum_{n \in \mscr{S} \atop{ n \leq m } } \li [ 1 - \Phi
\li ( \f{\vep \sq{n} } { b } \ri ) \ri ] \ovl{P}_n \la{ineq171} \ee
for $\si \in [a, b]$.  So, statement (III) follows from
(\ref{ineq169}), (\ref{ineq170}) and (\ref{ineq171}).

\subsubsection{Proof of Statement (IV)}

Applying (\ref{usa1}) and (\ref{usa2}), we have \bel \bb{E} [
\mathbf{n} ] & = & n_1 + \sum_{\ell = 1}^{s-1} (n_{\ell + 1} -
n_\ell ) \Pr \{ \mathbf{n} > n_\ell \} +
\sum_{n = n_s}^\iy \Pr \{ \mathbf{n} > n \} \nonumber\\
& = & n_1 + \sum_{\ell = 1}^{s-1} (n_{\ell + 1} - n_\ell ) H_\ell
(\si) + \sum_{n = n_s}^\iy H^\star (\si, n) \la{come1} \eel and \be
\la{come2} \bb{E} [ \mathbf{n} ] =  n_1 + \sum_{\ell = 1}^{s-1}
(n_{\ell + 1} - n_\ell ) H_\ell (\si) + \sum_{n = n_s}^m H^\star
(\si, n) + \sum_{n = m + 1}^\iy \Pr \{ \mathbf{n} > n \}. \ee Note
that
\[
\Pr \{ \mbf{n} > n \} < \Pr \li \{ \chi_{n_s - 1}^2
> \f{ n (n_s - 1) \vep^2 }{  ( \si t_{n_s - 1, \ze \de} )^2} \ri \} = \Pr \li \{ \chi_{n_s - 1}^2
> (n_s - 1) n \ga \ri \} < \li [ n \ga  e^{- (n \ga - 1)} \ri ]^\up
\]
for $n \geq m$, where the last inequality can be deduced from
Chernoff bounds.  Therefore,
\[
\sum_{n = m + 1}^\iy \Pr \{ \mbf{n} > n \} < \f{e^\up} {\ga} \sum_{n
= m + 1}^\iy g( n \ga ) \; \ga,
\]
where we have introduced function $g(x) = (x e^{-x})^\up$ for
simplicity of notations.   Note that $g(x)$ is monotonically
decreasing with respect to $x$ greater than $1$ because $g^\prime
(x) = \up g(x) \li ( \f{1}{ x } - 1 \ri ) < 0$ for $x > 1$. Making
use of the assumption that $n \ga \geq m \ga > 1$ and the monotone
decreasing property of $g(x)$, we have
 \[ \sum_{n = m + 1}^\iy g( n \ga ) \; \ga  <
\int_{m \ga}^\iy g ( x ) dx = \f{ \up!} { \up^{\up + 1} } \int_{m
\up \ga}^\iy \f{ \lm^\up e^{-\lm} } {\up!} d \lm,
\]  where
\bee \int_{m \up \ga}^\iy \f{ \lm^\up e^{-\lm} } {\up!} d \lm  & = &
e^{-m \up \ga} \sum_{i = 0}^\up \f{ (m \up \ga)^i } { i! } =  \Pr \{ K \leq \up \}\\
& < & \inf_{h > 0} e^{ h \up} \bb{E} [ e^{- h K} ]  =  e^{- m \up
\ga} \li ( \f{m \up \ga e}{\up} \ri )^\up = e^{- m \up \ga} (m \ga
e)^\up \eee with $K$ representing a Poisson random variable with
mean $m \up \ga$. It follows that
\[
\sum_{n = m + 1}^\iy  \Pr \{ \mbf{n} > n \}  < \f{ e^\up \up!} { \ga
\up^{\up + 1}  } e^{- m \up \ga} (m \ga e)^\up =  \f{\up!}{\ga \up}
\li (  \f{m \ga}{\up} \ri )^\up e^{- (m \ga - 2) \up}.
\]
Using inequality $\up! < \sq{2 \pi \up} \; \up^\up e^{- \up +
\f{1}{12 \up}}$, we have {\small \be \la{come3} \sum_{n = m + 1}^\iy
\Pr \{ \mbf{n} > n \} < \f{1}{\ga \up} \sq{2 \pi \up} \; \up^\up
e^{- \up + \f{1}{12 \up}} \li ( \f{m \ga}{\up} \ri )^\up e^{- (m \ga
- 2) \up} = \f{1}{\ga} \sq{\f{2 \pi}{\up} } \; (m \ga)^\up e^{- (m
\ga - 1) \up + \f{1}{12 \up}} < \f{3 (m \ga e)^\up}{\ga \sq{\up} \;
e^{m \ga \up} } . \ee}  So, the proof of statement (IV) can be
completed by combining (\ref{come1}), (\ref{come2}) and
(\ref{come3}).

\subsection{Proof of Theorem \ref{Normal_Cut} } \la{App_Normal_Cut}

By (\ref{try68}) and (\ref{nowuse}), we have \be \la{ineq175} \Pr \{
| \wh{\bs{\mu}} - \mu | \geq \vep \} \leq 2 \ze \de + \sum_{\ell =
1}^{s - 1} \Pr \{ | \ovl{X}_{\mbf{n}} - \mu | \geq \vep, \;
\mathbf{n} = n_\ell \}.  \ee By the definition of the sampling
scheme, we have {\small \bel  \sum_{\ell = 1}^{s - 1} \Pr \{ |
\ovl{X}_{\mbf{n}} - \mu | \geq \vep, \; \mathbf{n} = n_\ell \} &
\leq & \sum_{\ell = 1}^{s - 1} \Pr \li \{ S_{n_\ell} \leq \f{ n_\ell
(n_\ell - 1) \vep^2 }{ t_{n_\ell - 1, \ze \de}^2} \ri \} =
\sum_{\ell = 1}^{s - 1} \Pr \li \{ \chi_{n_\ell - 1}^2 \leq \f{
n_\ell (n_\ell - 1) \vep^2 }{ ( \si \; t_{n_\ell - 1, \ze \de})^2}
\ri \} \nonumber\\
 & = &  \sum_{\ell = 1}^{s - 1} \li [ 1 - S_{\mrm{P}} \li (  k_\ell - 1, \f{ n_\ell \; k_\ell \;
\vep^2 } { ( \si \; t_{n_s - 1, \ze \de} )^2 } \ri ) \ri ]
\la{ineq176} \eel} and {\small \bel &   & \sum_{\ell = 1}^{s - 1}
\Pr \{ | \ovl{X}_{\mbf{n}} - \mu | \geq \vep, \; \mathbf{n} = n_\ell \} \nonumber \\
 & \leq & \Pr \{ | \ovl{X}_{n_1} - \mu | \geq \vep \}  + \sum_{\ell = 1}^{s
- 2} \Pr \li \{ | \ovl{X}_{n_{\ell+1}} - \mu | \geq \vep, \;
S_{n_{\ell}} > \f{
n_\ell (n_\ell - 1) \vep^2 }{ t_{n_\ell - 1, \ze \de}^2}  \ri \} \nonumber\\
& =  & \Pr \{ | \ovl{X}_{n_1} - \mu | \geq \vep \}  +  \sum_{\ell =
1}^{s - 2} \Pr \li \{ | \ovl{X}_{n_{\ell+1}} - \mu | \geq \vep \ri
\} \Pr \li \{ S_{n_{\ell}} > \f{ n_\ell
(n_\ell - 1) \vep^2 }{ t_{n_\ell - 1, \ze \de}^2} \ri \} \nonumber \\
& = & 2 \li [ 1 - \Phi \li ( \f{ \vep \sq{n_1} } {  \si }  \ri ) \ri
] + 2  \sum_{\ell = 1}^{s - 2} \li [ 1 - \Phi \li ( \f{ \vep
\sq{n_{\ell + 1} } } {  \si }  \ri ) \ri ] \Pr \li \{ \chi_{n_\ell -
1}^2 > \f{ n_\ell (n_\ell - 1) \vep^2 }{ (\si \;  t_{n_\ell - 1, \ze
\de})^2}  \ri \} \nonumber\\
& = & 2 \li [ 1 - \Phi \li ( \f{ \vep \sq{n_1} } {  \si }  \ri ) \ri
] + 2  \sum_{\ell = 1}^{s - 2} \li [ 1 - \Phi \li ( \f{ \vep
\sq{n_{\ell + 1} } } {  \si }  \ri ) \ri ] S_{\mrm{P}} \li ( k_\ell
- 1, \f{ n_\ell \; k_\ell \; \vep^2 } { ( \si \; t_{n_s - 1, \ze
\de} )^2 } \ri ). \qqu \la{ineq177} \eel} Combining (\ref{ineq175})
and (\ref{ineq176}) yields \be \la{ineq178} \Pr \{ | \wh{\bs{\mu}} -
\mu | \geq \vep \} \leq 2 \ze \de + \sum_{\ell = 1}^{s - 1} \li [ 1
- S_{\mrm{P}} \li (  k_\ell - 1, \f{ n_\ell \; k_\ell \; \vep^2 } {
( \si \; t_{n_s - 1, \ze \de} )^2 } \ri ) \ri ], \ee where the upper
bound in the right side of (\ref{ineq178}) monotonically decreases
from $s - 1 + 2 \ze \de$ to $2 \ze \de$ as $\si$ increases from $0$
to $\iy$. Since $0 < \ze < \f{1}{2}$, there exists a unique number
$\ovl{\si}$ such that
\[
\sum_{\ell = 1}^{s - 1} \li [ 1 - S_{\mrm{P}} \li (  k_\ell - 1, \f{
n_\ell \; k_\ell \; \vep^2 } { ( \ovl{\si} \; t_{n_s - 1, \ze \de}
)^2 } \ri ) \ri ] = (1 - 2 \ze) \de
\]
and that $\Pr \{ | \wh{\bs{\mu}} - \mu | \geq \vep \} < \de$ for
$\si > \ovl{\si}$.  On the other hand, combining (\ref{ineq175}) and
(\ref{ineq177}) yields {\small \be \la{ineq179} \Pr \{ |
\wh{\bs{\mu}} - \mu | \geq \vep \} \leq 2 \ze \de + 2 \li [ 1 - \Phi
\li ( \f{ \vep \sq{n_1} } {  \si }  \ri ) \ri ] + 2  \sum_{\ell =
1}^{s - 2} \li [ 1 - \Phi \li ( \f{ \vep \sq{n_{\ell + 1} } } {  \si
}  \ri ) \ri ] S_{\mrm{P}} \li (  k_\ell - 1, \f{ n_\ell \; k_\ell
\; \vep^2 } { ( \si \; t_{n_s - 1, \ze \de} )^2 } \ri ), \ee} where
the upper bound in the right side of (\ref{ineq179}) monotonically
increases from $2 \ze \de$ to $s - 1 + 2 \ze \de$ as $\si$ increases
from $0$ to $\iy$. Since $0 < \ze < \f{1}{2}$, there exists a unique
number $\udl{\si}$ such that
\[
1 - \Phi \li ( \f{ \vep \sq{n_1} } { \udl{\si} }  \ri )  +
\sum_{\ell = 1}^{s - 2} \li [ 1 - \Phi \li ( \f{ \vep \sq{n_{\ell +
1} } } { \udl{\si} }  \ri ) \ri ] S_{\mrm{P}} \li (  k_\ell - 1, \f{
n_\ell \; k_\ell \; \vep^2 } { ( \udl{\si} \; t_{n_s - 1, \ze \de}
)^2 } \ri )  = \li ( \f{1}{2} - \ze \ri ) \de
\]
and that $\Pr \{ | \wh{\bs{\mu}} - \mu | \geq \vep \} < \de$ for
$\si < \udl{\si}$.  This completes the proof of Theorem
\ref{Normal_Cut}.

\subsection {Proof of Theorem \ref{normal_mean_rev} } \la{App_normal_mean_rev}

By the definition of the stopping rule, we have \bee  \Pr \li \{ \li
| \wh{\bs{\mu}} - \mu \ri |
> \vep |\mu| \ri \} \leq  \sum_{\ell = 1}^\iy \Pr \li \{ \li |
\wh{\bs{\mu}}_\ell - \mu \ri | > \vep |\mu|, \; \li |
\wh{\bs{\mu}}_\ell \ri | \geq \f{ t_{n_\ell - 1, \; \ze \de_\ell } }
{ \sq{n_\ell} } \li ( 1 + \f{1}{\vep} \ri ) \wh{\bs{\si}}_\ell \ri
\}. \eee By virtue of identity (\ref{RI_Indentity}), we have {\small
\bee &    & \Pr \li \{ \li | \wh{\bs{\mu}}_\ell - \mu \ri |
> \vep |\mu|, \; \li | \wh{\bs{\mu}}_\ell \ri | \geq \f{ t_{n_\ell -
1, \; \ze \de_\ell } } { \sq{n_\ell} } \li ( 1 + \f{1}{\vep} \ri )
\wh{\bs{\si}}_\ell \ri
\}\\
 & = & \Pr \li \{ \mu < \f{ \wh{\bs{\mu}}_\ell } { 1 + \mrm{sgn}
(\wh{\bs{\mu}}_\ell) \vep }, \; \li | \wh{\bs{\mu}}_\ell \ri | \geq
\f{ t_{n_\ell - 1, \; \ze \de_\ell } } { \sq{n_\ell} } \li ( 1 +
\f{1}{\vep} \ri ) \wh{\bs{\si}}_\ell  \ri \}\\
&    & + \Pr \li \{ \mu > \f{ \wh{\bs{\mu}}_\ell } { 1 - \mrm{sgn}
(\wh{\bs{\mu}}_\ell) \vep }, \; \li | \wh{\bs{\mu}}_\ell \ri | \geq
\f{ t_{n_\ell - 1, \; \ze \de_\ell } } { \sq{n_\ell} } \li ( 1 +
\f{1}{\vep} \ri ) \wh{\bs{\si}}_\ell  \ri \}\\
& = & \Pr \li \{ \wh{\bs{\mu}}_\ell - \mu > \f{ \vep |
\wh{\bs{\mu}}_\ell | } { 1 + \mrm{sgn} (\wh{\bs{\mu}}_\ell) \vep },
\; \li | \wh{\bs{\mu}}_\ell \ri | \geq \f{ t_{n_\ell - 1, \; \ze
\de_\ell } } { \sq{n_\ell} } \li ( 1 +
\f{1}{\vep} \ri ) \wh{\bs{\si}}_\ell  \ri \}\\
&    & + \Pr \li \{ \mu - \wh{\bs{\mu}}_\ell  > \f{ \vep
 | \wh{\bs{\mu}}_\ell | } { 1 - \mrm{sgn} (\wh{\bs{\mu}}_\ell) \vep }, \;
\li | \wh{\bs{\mu}}_\ell \ri | \geq \f{ t_{n_\ell - 1, \; \ze
\de_\ell } } { \sq{n_\ell} } \li ( 1 + \f{1}{\vep} \ri )
\wh{\bs{\si}}_\ell  \ri \}\\
& \leq & \Pr \li \{ | \wh{\bs{\mu}}_\ell - \mu | > \f{ \vep |
\wh{\bs{\mu}}_\ell | } { 1 +  \vep } \geq \f{ t_{n_\ell - 1, \; \ze
\de_\ell } } { \sq{n_\ell} } \wh{\bs{\si}}_\ell  \ri \} \leq  \Pr
\li \{ \f{   \sq{n_\ell} | \wh{\bs{\mu}}_\ell - \mu | }{
\wh{\bs{\si}}_\ell } > t_{n_\ell - 1, \; \ze \de_\ell }  \ri \} = 2
\ze \de_\ell \eee} for all $\ell > 0$.   Therefore, $\Pr \li \{ \li
| \wh{\bs{\mu}} - \mu \ri | > \vep |\mu| \ri \} \leq 2 \sum_{\ell =
1}^\iy \ze \de_\ell = 2 (\tau + 1) \ze \de$.

The finite stopping property of the sampling scheme can be shown by
an argument similar to the proof of statement (I) of Theorem
\ref{Bino_Rev_noninverse_Chernoff}.

\subsection {Proof of Theorem  \ref{normalmix} } \la{App_normalmix}

By the definition of the stopping rule, we have {\small \bee  \Pr
\li \{ \li | \wh{\bs{\mu}} - \mu \ri | > \max ( \vep_a, \; \vep_r
|\mu| ) \ri \} \leq  \sum_{\ell = 1}^\iy \Pr \li \{ \li |
\wh{\bs{\mu}}_\ell - \mu \ri | > \max ( \vep_a, \; \vep_r |\mu| ),
\; \max \li ( \vep_a, \f{\vep_r \li | \wh{\bs{\mu}}_\ell \ri |}{1 +
\vep_r}  \ri ) \geq \f{ t_{n_\ell - 1, \; \ze \de_\ell } } {
\sq{n_\ell} } \wh{\bs{\si}}_\ell \ri \}.  \eee} By virtue of
identity (\ref{RI_Indentity}), we have {\small \bee & & \Pr \li \{
\li | \wh{\bs{\mu}}_\ell  - \mu \ri | > \max ( \vep_a, \; \vep_r
|\mu| ), \; \max \li ( \vep_a, \f{\vep_r \li | \wh{\bs{\mu}}_\ell
\ri |}{1 + \vep_r}  \ri ) \geq \f{ t_{n_\ell - 1,
\; \ze \de_\ell } } { \sq{n_\ell} } \wh{\bs{\si}}_\ell \ri \}\\
& = &  \Pr \li \{ \mu < \min \li ( \wh{\bs{\mu}}_\ell - \vep_a, \;
\f{ \wh{\bs{\mu}}_\ell } { 1 + \mrm{sgn} (\wh{\bs{\mu}}_\ell) \vep_r
} \ri ),  \; \max \li ( \vep_a, \f{\vep_r \li | \wh{\bs{\mu}}_\ell
\ri |}{1 + \vep_r}  \ri ) \geq \f{ t_{n_\ell - 1, \; \ze \de_\ell }
} { \sq{n_\ell} } \wh{\bs{\si}}_\ell \ri \}\\
&  & + \Pr \li \{ \mu > \max \li ( \wh{\bs{\mu}}_\ell + \vep_a, \;
\f{ \wh{\bs{\mu}}_\ell } { 1 - \mrm{sgn} (\wh{\bs{\mu}}_\ell) \vep_r
} \ri ),  \; \max \li ( \vep_a, \f{\vep_r \li | \wh{\bs{\mu}}_\ell
\ri |}{1 + \vep_r}  \ri ) \geq \f{ t_{n_\ell - 1, \; \ze \de_\ell }
} { \sq{n_\ell} } \wh{\bs{\si}}_\ell \ri \}\\
& = &  \Pr \li \{ \wh{\bs{\mu}}_\ell - \mu > \max \li ( \vep_a, \;
\f{ \vep_r |\wh{\bs{\mu}}_\ell| } { 1 + \mrm{sgn}
(\wh{\bs{\mu}}_\ell) \vep_r } \ri ),  \; \max \li ( \vep_a,
\f{\vep_r \li | \wh{\bs{\mu}}_\ell \ri |}{1 + \vep_r}  \ri ) \geq
\f{ t_{n_\ell - 1, \; \ze \de_\ell }
} { \sq{n_\ell} } \wh{\bs{\si}}_\ell \ri \}\\
&  & + \Pr \li \{ \mu - \wh{\bs{\mu}}_\ell > \max \li (  \vep_a, \;
\f{ \vep_r | \wh{\bs{\mu}}_\ell | } { 1 - \mrm{sgn}
(\wh{\bs{\mu}}_\ell) \vep_r } \ri ),  \; \max \li ( \vep_a,
\f{\vep_r \li | \wh{\bs{\mu}}_\ell \ri |}{1 + \vep_r}  \ri ) \geq
\f{ t_{n_\ell - 1, \; \ze \de_\ell } } { \sq{n_\ell} }
\wh{\bs{\si}}_\ell \ri \}\\
& \leq &  \Pr \li \{ | \wh{\bs{\mu}}_\ell - \mu | >  \max \li (
\vep_a, \f{\vep_r \li | \wh{\bs{\mu}}_\ell \ri |}{1 + \vep_r}  \ri )
\geq \f{ t_{n_\ell - 1, \; \ze \de_\ell } } { \sq{n_\ell} }
\wh{\bs{\si}}_\ell \ri \} \\
& \leq & \Pr \li \{ \f{   \sq{n_\ell} | \wh{\bs{\mu}}_\ell - \mu |
}{ \wh{\bs{\si}}_\ell } > t_{n_\ell - 1, \; \ze \de_\ell }  \ri \} =
2 \ze \de_\ell \eee} for all $\ell > 0$.  Therefore, $\Pr \li \{ \li
| \wh{\bs{\mu}} - \mu \ri | > \max ( \vep_a, \; \vep_r |\mu| ) \ri
\} \leq 2 \sum_{\ell = 1}^\iy \ze \de_\ell = 2 (\tau + 1) \ze \de$.

The finite stopping property of the sampling scheme can be shown by
an argument similar to the proof of statement (I) of Theorem
\ref{Bino_Rev_noninverse_Chernoff}.

\section{Proof of Theorem \ref{Classic_CI} }  \la{App_Classic_CI}

Since $\wh{\bs{\se}}$ is a ULE of $\se$, by virtue of Lemma
\ref{ULE_Basic}, we have that $\Pr \{ \wh{\bs{\se}} \leq z \mid \se
\}$ is non-increasing with respect to $\se$ no less than $z$ and
that $\Pr \{ \wh{\bs{\se}} > z \mid \se \}$ is non-decreasing with
respect to $\se$ no greater than $z$.  This implies that $\Pr \{
\wh{\bs{\se}} \leq z \mid \se \}$ is non-increasing with respect to
$\se \in \Se$ and that $\Pr \{ \wh{\bs{\se}} \geq z \mid \se \}$ is
non-decreasing with respect to $\se \in \Se$. By the definitions of
$F_{\wh{\bs{\se}}} (z, \se)$ and $G_{\wh{\bs{\se}}} (z, \se)$ given
in Section \ref{defFG}, we have that $F_{\wh{\bs{\se}}} (z, \se)$ is
non-increasing with respect to $\se \in \Se$ and that
$G_{\wh{\bs{\se}}} (z, \se)$ is non-decreasing with respect to $\se
\in \Se$.  Recalling the definition of $\mscr{U} ( \wh{\bs{\se}},
\mbf{n} )$, we have that $\{ F_{\wh{\bs{\se}}} ( \wh{\bs{\se}},
\mscr{U} ( \wh{\bs{\se}}, \mbf{n} ) ) \leq \f{\de}{2} \}$ is a sure
event and consequently {\small \[ \{ \se \geq \mscr{U} (
\wh{\bs{\se}}, \mbf{n} ) \} = \li \{ \se \geq \mscr{U} (
\wh{\bs{\se}}, \mbf{n} ), \; F_{\wh{\bs{\se}}} (\wh{\bs{\se}},
\mscr{U} ( \wh{\bs{\se}}, \mbf{n} ) ) \leq \f{\de}{2} \ri \}
\subseteq \li \{ \se \geq \mscr{U} ( \wh{\bs{\se}}, \mbf{n} ), \;
F_{\wh{\bs{\se}}} (\wh{\bs{\se}}, \se ) \leq \f{\de}{2} \ri \}
\subseteq \li \{ F_{\wh{\bs{\se}}} (\wh{\bs{\se}}, \se ) \leq
\f{\de}{2} \ri \},
\]}
which implies that \be \la{fact1} \Pr \{ \se \geq \mscr{U} (
\wh{\bs{\se}}, \mbf{n} ) \} \leq \Pr \li \{ F_{\wh{\bs{\se}}}
(\wh{\bs{\se}}, \se ) \leq \f{\de}{2} \ri \} \leq \f{\de}{2}, \ee
where the last inequality follows from  Lemma \ref{ProbTrans}.  On
the other hand, recalling the definition of $\mscr{L} (
\wh{\bs{\se}}, \mbf{n} )$, we have that $\{ G_{\wh{\bs{\se}}} (
\wh{\bs{\se}}, \mscr{L} ( \wh{\bs{\se}}, \mbf{n} ) ) \leq \f{\de}{2}
\}$ is a sure event and consequently {\small \[  \{ \se \leq
\mscr{L} ( \wh{\bs{\se}}, \mbf{n} ) \} = \li \{ \se \leq \mscr{L} (
\wh{\bs{\se}}, \mbf{n} ), \; G_{\wh{\bs{\se}}} ( \wh{\bs{\se}},
\mscr{L} ( \wh{\bs{\se}}, \mbf{n} ) ) \leq \f{\de}{2} \ri \}
\subseteq \li \{ \se \leq \mscr{L} ( \wh{\bs{\se}}, \mbf{n} ), \;
G_{\wh{\bs{\se}}} ( \wh{\bs{\se}}, \se ) \leq \f{\de}{2} \ri \}
\subseteq \li \{ G_{\wh{\bs{\se}}} ( \wh{\bs{\se}}, \se ) \leq
\f{\de}{2}  \ri \},
\]}
which implies that \be \la{fact2} \Pr \{ \se \leq \mscr{L} (
\wh{\bs{\se}}, \mbf{n} ) \} \leq \Pr \li \{ G_{\wh{\bs{\se}}} (
\wh{\bs{\se}}, \se ) \leq \f{\de}{2} \ri \} \leq \f{\de}{2}, \ee
where the last inequality follows from  Lemma \ref{ProbTrans}.
Finally, by virtue of (\ref{fact1}) and (\ref{fact2}), we have
\[
\Pr \{  \mscr{L} ( \wh{\bs{\se}}, \mbf{n} ) < \se <  \mscr{U} (
\wh{\bs{\se}}, \mbf{n} )  \mid \se \} \geq 1 - \Pr \{ \se \geq
\mscr{U} ( \wh{\bs{\se}}, \mbf{n} ) \} - \Pr \{ \se \leq \mscr{L} (
\wh{\bs{\se}}, \mbf{n} ) \} \geq 1 - \f{\de}{2} - \f{\de}{2} = 1 -
\de.
\]
This completes the proof of the theorem.

\section{Proof of Theorem \ref{Pos_CI_Test} } \la{App_Pos_CI_Test}

Note that \bee &  & \Pr \{ \mscr{L} (\wh{\bs{\lm}}_\ell,
\mbf{n}_\ell) < \lm < \mscr{U} (\wh{\bs{\lm}}_\ell, \mbf{n}_\ell) \mid \lm  \}\\
& \geq & \Pr \li \{ \mscr{L} (\wh{\bs{\lm}}_\ell, \mbf{n}_\ell) <
\lm < \mscr{U} (\wh{\bs{\lm}}_\ell, \mbf{n}_\ell), \; U \li (
\wh{\bs{\lm}}_\ell, \mbf{n}_\ell,
\f{\de}{2s} \ri ) > \lm^* \mid \lm  \ri \}\\
& = & \Pr \li \{ L \li ( \wh{\bs{\lm}}_\ell, \mbf{n}_\ell,
\f{\de}{2s} \ri ) < \lm < U \li ( \wh{\bs{\lm}}_\ell, \mbf{n}_\ell,
\f{\de}{2s} \ri ), \; U \li ( \wh{\bs{\lm}}_\ell, \mbf{n}_\ell,
\f{\de}{2s} \ri ) > \lm^* \mid \lm  \ri \}\\
& = & \Pr \li \{ L \li ( \wh{\bs{\lm}}_\ell, \mbf{n}_\ell,
\f{\de}{2s} \ri ) < \lm < U \li ( \wh{\bs{\lm}}_\ell, \mbf{n}_\ell,
\f{\de}{2s} \ri ) \mid \lm  \ri \} \geq 1 - \f{\de}{2s}
 \eee
for any $\lm \in [\lm^*, \iy)$.  Therefore, \bee  \Pr \li \{ \lm
\notin \li ( \mscr{L} (\wh{\bs{\lm}}_\ell, \mbf{n}_\ell) , \;
\mscr{U} (\wh{\bs{\lm}}_\ell, \mbf{n}_\ell) \ri ), \; \bs{l} = \ell
\mid \lm \ri \} \leq  \Pr \li \{ \lm \notin \li ( \mscr{L}
(\wh{\bs{\lm}}_\ell, \mbf{n}_\ell) , \; \mscr{U}
(\wh{\bs{\lm}}_\ell, \mbf{n}_\ell ) \ri ) \mid \lm \ri \} \leq
\f{\de}{2s} \eee for $\ell = 1, \cd, s $ and any $\lm \in [\lm^*,
\iy)$. It follows that \bee \Pr \li \{ \lm \notin \li ( \mscr{L}
(\wh{\bs{\lm}}, \mbf{n} ) , \; \mscr{U} (\wh{\bs{\lm}}, \mbf{n} )
\ri )  \mid \lm \ri \}
  =  \sum_{\ell = 1}^s \Pr \li \{ \lm \notin \li ( \mscr{L}
(\wh{\bs{\lm}}_\ell, \mbf{n}_\ell ) , \; \mscr{U}
(\wh{\bs{\lm}}_\ell, \mbf{n}_\ell ) \ri ), \; \bs{l} = \ell \mid \lm
\ri \} \leq \de \eee for any $\lm \in [\lm^*, \iy)$.  The theorem
immediately follows.

\section{Proof of Theorem \ref{goal_RCI_Pos} }
\la{App_goal_RCI_Pos}

Note that \bee &   & \Pr \{ \mscr{L} (\wh{\bs{\lm}}_\ell,
\mbf{n}_\ell) < \lm < \mscr{U} (\wh{\bs{\lm}}_\ell, \mbf{n}_\ell) \mid \lm  \}\\
& \geq & \Pr \li \{ \mscr{L} (\wh{\bs{\lm}}_\ell, \mbf{n}_\ell) <
\lm < \mscr{U} (\wh{\bs{\lm}}_\ell, \mbf{n}_\ell), \; U \li (
\wh{\bs{\lm}}_\ell, \mbf{n}_\ell,
\f{\de}{2s} \ri ) > \lm^* \mid \lm  \ri \}\\
& = & \Pr \li \{ L \li ( \wh{\bs{\lm}}_\ell, \mbf{n}_\ell,
\f{\de}{2s} \ri ) < \lm < U \li ( \wh{\bs{\lm}}_\ell, \mbf{n}_\ell,
\f{\de}{2s} \ri ), \; U \li ( \wh{\bs{\lm}}_\ell, \mbf{n}_\ell,
\f{\de}{2s} \ri ) > \lm^* \mid \lm  \ri \}\\
& = & \Pr \li \{ L \li ( \wh{\bs{\lm}}_\ell, \mbf{n}_\ell,
\f{\de}{2s} \ri ) < \lm < U \li ( \wh{\bs{\lm}}_\ell, \mbf{n}_\ell,
\f{\de}{2s} \ri ) \mid \lm  \ri \} \geq 1 - \f{\de}{2s}
 \eee
for any $\lm \in [\lm^*, \iy)$.  The theorem immediately follows.

\sect{Proofs of Theorems for Multistage Linear Regression}

\subsection{Proof of Theorem \ref{Regression_Abs} } \la{App_Regression_Abs}

By the definition of the stopping rule, \bel \Pr \{ |
\wh{\bs{\ba}}_i - \beta_i |
> \vep_i \} & \leq  & \sum_{\ell = 1}^\iy \Pr \li \{ |
\bs{B}_{i, \ell} - \beta_i | > \vep_i \geq t_{n_\ell - m, \; \ze
\de_\ell} \; \widehat{\bs{\si}}_\ell \; \sqrt{\left[
(\boldsymbol{X}_\ell^\intercal \boldsymbol{X}_\ell)^{-1}
\right]_{ii}} \ri \} \nonumber\\
& \leq  & \sum_{\ell = 1}^\iy \Pr \li \{  \f{ | \bs{B}_{i, \ell} -
\beta_i | } { \widehat{\bs{\si}}_\ell \; \sqrt{\left[
(\boldsymbol{X}_\ell^\intercal \boldsymbol{X}_\ell)^{-1}
\right]_{ii}}  } > t_{n_\ell - m, \; \ze \de_\ell} \ri \}
\la{tbound} \eel for $i = 1, \cd, m$.  From the classical theory of
linear regression, we know that $\bs{B}_{i, \ell} - \beta_i$ is a
Gaussian random variable of zero mean,  variance $\si^2 \left[
(\boldsymbol{X}_\ell^\intercal \boldsymbol{X}_\ell)^{-1}
\right]_{ii}$ and that $(n_\ell - m ) ( \f{ \widehat{\bs{\si}}_\ell
} {\si}  )^2$ is a chi-square variable of $n_\ell - m$ degrees of
freedom. Moreover, $\bs{B}_{i, \ell} - \beta_i$ is independent of
$(n_\ell - m ) ( \f{ \widehat{\bs{\si}}_\ell }{\si} )^2$.  It
follows that {\small $( \bs{B}_{i, \ell} - \beta_i )  \li \{
\widehat{\bs{\si}}_\ell \; \sqrt{\left[
(\boldsymbol{X}_\ell^\intercal \boldsymbol{X}_\ell)^{-1}
\right]_{ii}}  \ri \}^{-1}$} possesses a Student $t$-distribution of
$n_\ell - m $ degrees of freedom.  Hence, by (\ref{tbound}), we have
\be \la{conclu1}
 \Pr \{ | \wh{\bs{\ba}}_i - \beta_i |
> \vep_i \} \leq 2 \sum_{\ell = 1}^\iy \ze \de_\ell = 2 (\tau + 1) \ze
\de \ee for $i = 1, \cd, m$.  By the definition of the stopping
rule, \bel \Pr \{ | \wh{\bs{\si}} - \sigma | > \vep \} & \leq &
\sum_{\ell = 1}^\iy \Pr \li \{ | \wh{\bs{\si}}_\ell - \sigma | >
\vep, \; \sqrt{ \frac{ n_\ell - m} { \chi_{n_\ell - m, \; \ze
\de_\ell }^2 } } \; \widehat{\bs{\si}}_\ell - \vep \leq
\widehat{\bs{\si}}_\ell \leq \sqrt{ \frac{ n_\ell - m } {
\chi_{n_\ell - m, \; 1 - \ze \de_\ell }^2 } }
\; \widehat{\bs{\si}}_\ell + \vep \ri \} \nonumber \\
& \leq &  \sum_{\ell = 1}^\iy \Pr \li \{ \wh{\bs{\si}}_\ell - \sigma
< - \vep, \; \sqrt{ \frac{ n_\ell - m} { \chi_{n_\ell - m, \; \ze
\de_\ell }^2 } } \; \widehat{\bs{\si}}_\ell - \vep \leq
\widehat{\bs{\si}}_\ell \ri \} \nonumber \\
&    & + \sum_{\ell = 1}^\iy \Pr \li \{ \wh{\bs{\si}}_\ell - \sigma
> \vep, \; \sqrt{ \frac{ n_\ell - m} { \chi_{n_\ell - m, \; 1 - \ze
\de_\ell }^2 } } \; \widehat{\bs{\si}}_\ell + \vep \geq
\widehat{\bs{\si}}_\ell \ri \} \nonumber \\
& \leq &  \sum_{\ell = 1}^\iy \Pr \li \{  \sqrt{ \frac{ n_\ell - m}
{ \chi_{n_\ell - m, \; \ze \de_\ell }^2 } } \;
\widehat{\bs{\si}}_\ell < \si \ri \} + \sum_{\ell = 1}^\iy \Pr \li
\{ \sqrt{ \frac{ n_\ell - m} { \chi_{n_\ell - m, \; 1 - \ze \de_\ell
}^2 } } \; \widehat{\bs{\si}}_\ell > \si \ri \}. \la{chibound} \eel
Recalling that $(n_\ell - m ) ( \f{ \widehat{\bs{\si}}_\ell } {\si}
)^2$ is a chi-square variable of $n_\ell - m$ degrees of freedom, we
have \be \la{conclu2}
 \Pr \li \{
\sqrt{ \frac{ n_\ell - m} { \chi_{n_\ell - m, \; \ze \de_\ell }^2 }
} \; \widehat{\bs{\si}}_\ell < \si \ri \} \leq \ze \de_\ell, \qqu
\Pr \li \{ \sqrt{ \frac{ n_\ell - m} { \chi_{n_\ell - m, \; 1 - \ze
\de_\ell }^2 } } \; \widehat{\bs{\si}}_\ell > \si \ri \} \leq \ze
\de_\ell \ee for all $\ell > 0$. Combining (\ref{chibound}) and
(\ref{conclu2}) yields \be \la{conclu3} \Pr \{ | \wh{\bs{\si}} -
\sigma | > \vep \} \leq 2 \sum_{\ell = 1}^\iy \ze \de_\ell = 2 (\tau
+ 1) \ze \de. \ee By virtue of (\ref{conclu1}) and (\ref{conclu3}),
we have \bee &
  & \Pr \{ | \wh{\bs{\si}} - \sigma | \leq \vep, \;\; | \wh{\bs{\ba}}_i
- \beta_i | \leq \vep_i \; \tx{for} \; i = 1, \cdots, m \} \\
& \geq &  1 - \sum_{i = 1}^m \Pr \{  | \wh{\bs{\ba}}_i - \beta_i | >
\vep_i  \}  - \Pr \{ | \wh{\bs{\si}} - \sigma | > \vep \}\\
& \geq & 1 - 2 m  (\tau + 1) \ze \de - 2 (\tau + 1) \ze \de = 1 - 2
(m + 1)  (\tau + 1) \ze \de.  \eee

The finite stopping property of the sampling scheme can be shown by
an argument similar to the proof of statement (I) of Theorem
\ref{Bino_Rev_noninverse_Chernoff}.  This completes the proof of the
theorem.

\subsection{Proof of Theorem \ref{Regression_Rev} } \la{App_Regression_Rev}

By the definition of the stopping rule, {\small \be \la{8966} \Pr \{
| \wh{\bs{\ba}}_i - \beta_i | > \vep_i | \beta_i | \}  \leq
\sum_{\ell = 1}^\iy \Pr \li \{ | \bs{B}_{i, \ell} - \beta_i | >
\vep_i | \beta_i |, \; t_{n_\ell - m, \; \ze \de_\ell} \;
\widehat{\bs{\si}}_\ell \; \sqrt{\left[
(\boldsymbol{X}_\ell^\intercal \boldsymbol{X}_\ell)^{-1}
\right]_{ii}} \leq \f{ \vep_i | \bs{B}_{i, \ell} | }{1 + \vep_i }
\ri \}  \qu \ee} for $i = 1, \cd, m$.  By identity
(\ref{RI_Indentity}), we have \bel &   & \Pr \li \{ | \bs{B}_{i,
\ell} - \beta_i | > \vep_i | \beta_i |, \; t_{n_\ell - m, \; \ze
\de_\ell} \; \widehat{\bs{\si}}_\ell \; \sqrt{\left[
(\boldsymbol{X}_\ell^\intercal \boldsymbol{X}_\ell)^{-1}
\right]_{ii}} \leq \f{ \vep_i }{1 + \vep_i }  | \bs{B}_{i, \ell} |
\ri \} \nonumber\\
& = & \Pr \li \{ \beta_i < \f{ \bs{B}_{i, \ell} } { 1 + \mrm{sgn}
(\bs{B}_{i, \ell}) \vep_i }, \; t_{n_\ell - m, \; \ze \de_\ell} \;
\widehat{\bs{\si}}_\ell \; \sqrt{\left[
(\boldsymbol{X}_\ell^\intercal \boldsymbol{X}_\ell)^{-1}
\right]_{ii}} \leq \f{ \vep_i }{1 + \vep_i }  | \bs{B}_{i, \ell} | \ri \} \nonumber\\
&   & + \Pr \li \{ \beta_i > \f{ \bs{B}_{i, \ell} } { 1 - \mrm{sgn}
(\bs{B}_{i, \ell}) \vep_i }, \; t_{n_\ell - m, \; \ze \de_\ell} \;
\widehat{\bs{\si}}_\ell \; \sqrt{\left[
(\boldsymbol{X}_\ell^\intercal \boldsymbol{X}_\ell)^{-1}
\right]_{ii}} \leq \f{ \vep_i }{1 + \vep_i }  | \bs{B}_{i, \ell} |
\ri \} \nonumber\\
& = & \Pr \li \{ \bs{B}_{i, \ell} - \beta_i > \f{ \vep_i |
\bs{B}_{i, \ell} |  } { 1 + \mrm{sgn} (\bs{B}_{i, \ell}) \vep_i },
\; t_{n_\ell - m, \; \ze \de_\ell} \; \widehat{\bs{\si}}_\ell \;
\sqrt{\left[ (\boldsymbol{X}_\ell^\intercal
\boldsymbol{X}_\ell)^{-1}
\right]_{ii}} \leq \f{ \vep_i }{1 + \vep_i }  | \bs{B}_{i, \ell} | \ri \} \nonumber\\
&   & + \Pr \li \{ \beta_i - \bs{B}_{i, \ell} > \f{ \vep_i |
\bs{B}_{i, \ell} | } { 1 - \mrm{sgn} (\bs{B}_{i, \ell}) \vep_i }, \;
t_{n_\ell - m, \; \ze \de_\ell} \; \widehat{\bs{\si}}_\ell \;
\sqrt{\left[ (\boldsymbol{X}_\ell^\intercal
\boldsymbol{X}_\ell)^{-1} \right]_{ii}} \leq \f{ \vep_i }{1 + \vep_i
}  | \bs{B}_{i, \ell} | \ri \} \nonumber \\
& \leq & \Pr \li \{ | \bs{B}_{i, \ell} - \beta_i | > \f{ \vep_i |
\bs{B}_{i, \ell} |  } { 1 + \vep_i } \geq t_{n_\ell - m, \; \ze
\de_\ell} \; \widehat{\bs{\si}}_\ell \; \sqrt{\left[
(\boldsymbol{X}_\ell^\intercal \boldsymbol{X}_\ell)^{-1}
\right]_{ii}} \ri \} \nonumber\\
& \leq  &  \Pr \li \{  \f{ | \bs{B}_{i, \ell} - \beta_i | } {
\widehat{\bs{\si}}_\ell \; \sqrt{\left[
(\boldsymbol{X}_\ell^\intercal \boldsymbol{X}_\ell)^{-1}
\right]_{ii}}  } > t_{n_\ell - m, \; \ze \de_\ell} \ri \} \nonumber \\
& = & 2 \ze \de_\ell \la{eq8833} \eel for $i = 1, \cd, m$, where the
last equality (\ref{eq8833}) follows from the fact that {\small $(
\bs{B}_{i, \ell} - \beta_i ) \li \{ \widehat{\bs{\si}}_\ell \;
\sqrt{\left[ (\boldsymbol{X}_\ell^\intercal
\boldsymbol{X}_\ell)^{-1} \right]_{ii}}  \ri \}^{-1}$} possesses a
Student $t$-distribution of $n_\ell - m $ degrees of freedom.
Combining (\ref{8966}) and (\ref{eq8833}) yields \be
\la{conclurev88} \Pr \{ | \wh{\bs{\ba}}_i - \beta_i | > \vep_i |
\beta_i | \}  \leq 2 \sum_{\ell = 1}^\iy \ze \de_\ell = 2 (\tau + 1)
\ze \de  \ee for $i = 1, \cd, m$.  By the definition of the stopping
rule, {\small \bel &   & \Pr \{ | \wh{\bs{\si}} - \sigma | > \vep
\si \} \leq  \sum_{\ell = 1}^\iy \Pr \li \{ | \wh{\bs{\si}}_\ell -
\sigma |
> \vep \si, \; \f{ \chi_{n_\ell - m, \; 1 - \ze \de_\ell }^2  } { (1
+ \vep)^2 } \leq n_\ell - m \leq \f{
\chi_{n_\ell - m, \; \ze \de_\ell }^2  } { (1 - \vep)^2 } \ri \} \nonumber \\
& \leq &  \sum_{\ell = 1}^\iy \li [ \Pr \li \{ \wh{\bs{\si}}_\ell  <
(1 - \vep ) \si \leq \si \sqrt{ \frac{ \chi_{n_\ell - m, \; \ze
\de_\ell }^2 } { n_\ell - m} }  \ri \}  + \Pr \li \{
\wh{\bs{\si}}_\ell
> (1 + \vep) \si \geq \si \sqrt{ \frac{ \chi_{n_\ell - m, \;
1 - \ze \de_\ell }^2 } { n_\ell - m} } \ri \} \ri ] \nonumber \\
& \leq &  \sum_{\ell = 1}^\iy \li [ \Pr \li \{  \sqrt{ \frac{ n_\ell
- m} { \chi_{n_\ell - m, \; \ze \de_\ell }^2 } } \;
\widehat{\bs{\si}}_\ell < \si \ri \} + \Pr \li \{ \sqrt{ \frac{
n_\ell - m} { \chi_{n_\ell - m, \; 1 - \ze
\de_\ell }^2 } } \; \widehat{\bs{\si}}_\ell > \si \ri \} \ri ] \nonumber\\
& =  & 2 (\tau + 1) \ze \de,  \qu \qu \la{chiboundrev} \eel} where
(\ref{chiboundrev}) follows from an argument similar to that of
(\ref{conclu2}).  Making use of (\ref{conclurev88}) and
(\ref{chiboundrev}), we have \bee &
  & \Pr \{ | \wh{\bs{\si}} - \sigma | \leq \vep \si, \;\; | \wh{\bs{\ba}}_i
- \beta_i | \leq \vep_i | \beta_i | \; \tx{for} \; i = 1, \cdots, m \} \\
& \geq &  1 - \sum_{i = 1}^m \Pr \{  | \wh{\bs{\ba}}_i - \beta_i | >
\vep_i | \beta_i | \}  - \Pr \{ | \wh{\bs{\si}} - \sigma | > \vep \si \}\\
& \geq & 1 - 2 m  (\tau + 1) \ze \de - 2 (\tau + 1) \ze \de = 1 - 2
(m + 1)  (\tau + 1) \ze \de.  \eee The finite stopping property of
the sampling scheme can be shown by an argument similar to the proof
of statement (I) of Theorem \ref{Bino_Rev_noninverse_Chernoff}. This
completes the proof of the theorem.

\sect{Proofs of Theorems for Estimation of Quantile}

\subsection{Proof of Theorem \ref{quantile_abs} } \la{App_quantile_abs}

By the definition of the stopping rule,  {\small \be  \la{good1abs}
\Pr \{ | \wh{\bs{\xi}}_{p} - \xi_p | > \vep \}  \leq \sum_{\ell =
1}^\iy \Pr \li \{ | \wh{\bs{\xi}}_{p, \ell} - \xi_p | > \vep, \;
X_{j_\ell: n_\ell} - \vep \leq \wh{\bs{\xi}}_{p, \ell} \leq
X_{i_\ell: n_\ell} + \vep \ri \}, \ee where \bel &   & \Pr \li \{ |
\wh{\bs{\xi}}_{p, \ell} - \xi_p | > \vep, \; X_{j_\ell:
n_\ell} - \vep \leq \wh{\bs{\xi}}_{p, \ell} \leq X_{i_\ell: n_\ell} + \vep \ri \} \nonumber\\
& \leq & \Pr \li \{  \xi_p < \wh{\bs{\xi}}_{p, \ell} - \vep \leq
X_{i_\ell: n_\ell} \ri \} +  \Pr \li \{  \xi_p >
\wh{\bs{\xi}}_{p, \ell} + \vep \geq X_{j_\ell: n_\ell} \ri \} \nonumber\\
& \leq &   \Pr \li \{  X_{i_\ell: n_\ell} > \xi_p \ri \}  + \Pr \li
\{  X_{j_\ell: n_\ell} < \xi_p \ri \} \la{good2abs} \eel} for all
$\ell > 0$.

Now, let $K_\ell$ denote the number of samples among $X_1, \cd,
X_{n_\ell}$ which are no greater than $\xi_p$.  Then, $\li \{
X_{i_\ell: n_\ell} > \xi_p \ri \} \subseteq \{ K_\ell < i_\ell \}$
and thus $\Pr \li \{ X_{i_\ell: n_\ell} > \xi_p \ri \} \leq \Pr \{
K_\ell < i_\ell \} = \sum_{k = 0}^{i_\ell - 1} \bi{n_\ell}{k} [F_X (
\xi_p ) ]^k [1 - F_X ( \xi_p ) ]^{n_\ell - k}$.  By the definition
of $\xi_p$, we have $F_X ( \xi_p ) \geq p$. Making use of the fact
that $\sum_{k = 0}^m \bi{n}{k} \se^k (1 - \se)^{n - k}$ is
monotonically decreasing with respect to $\se \in (0, 1)$, we have
that \be \la{good3} \Pr \li \{ X_{i_\ell: n_\ell} > \xi_p \ri \}
\leq \sum_{k = 0}^{i_\ell - 1} \bi{n_\ell}{k} p^k (1 - p)^{n_\ell -
k} \leq \ze \de_\ell, \ee where the last inequality follows from the
definition of $i_\ell$. On the other hand, let $K_\ell^*$ denote the
number of samples among $X_1, \cd, X_{n_\ell}$ which are smaller
than $\xi_p$.  Then, $ \{ X_{j_\ell: n_\ell} < \xi_p \} \subseteq \{
K_\ell^* \geq j_\ell \}$ and thus $\Pr \li \{ X_{j_\ell: n_\ell} <
\xi_p \ri \} \leq \Pr \{ K_\ell^* \geq j_\ell \} = \sum_{k =
j_\ell}^{n_\ell} \bi{n_\ell}{k} [F_X^{-} ( \xi_p ) ]^k [1 - F_X^{-}
( \xi_p ) ]^{n_\ell - k}$, where $F_X^{-} ( \xi_p ) = \Pr \{ X <
\xi_p \}$.  By the definition of $\xi_p$, we have $F_X^{-} ( \xi_p )
\leq p$. Making use of the fact that $\sum_{k = m}^n \bi{n}{k} \se^k
(1 - \se)^{n - k}$ is monotonically increasing with respect to $\se
\in (0, 1)$, we have that \be \la{good4}
 \Pr \li \{ X_{j_\ell:
n_\ell} < \xi_p \ri \} \leq  \sum_{k = j_\ell}^{n_\ell}
\bi{n_\ell}{k} p^k (1 - p)^{n_\ell - k} \leq \ze \de_\ell, \ee where
the last inequality follows from the definition of $j_\ell$.
Combining (\ref{good1abs}), (\ref{good2abs}), (\ref{good3}) and
(\ref{good4}) yields $\Pr \{ | \wh{\bs{\xi}}_{p} - \xi_p | > \vep \}
\leq 2 \sum_{\ell = 1}^\iy \ze \de_\ell = 2 (\tau + 1) \ze \de$.
The finite stopping property of the sampling scheme can be shown by
an argument similar to the proof of statement (I) of Theorem
\ref{Bino_Rev_noninverse_Chernoff}.

\subsection{Proof of Theorem \ref{quantile_rev} } \la{App_quantile_rev}

By the definition of the stopping rule,  {\small \be  \la{good1} \Pr
\{ | \wh{\bs{\xi}}_{p} - \xi_p | > \vep | \xi_p | \}  \leq
\sum_{\ell = 1}^\iy \Pr \{ | \wh{\bs{\xi}}_{p, \ell} - \xi_p | >
\vep | \xi_p |, \; [ 1 - \mrm{sgn} ( \wh{\bs{\xi}}_{p, \ell} ) \vep]
X_{j_\ell: n_\ell}  \leq \wh{\bs{\xi}}_{p, \ell} \leq [ 1 +
\mrm{sgn} ( \wh{\bs{\xi}}_{p, \ell} ) \vep] X_{i_\ell: n_\ell} \}.
\ee By identity (\ref{RI_Indentity}), we have \bel &   & \Pr \li \{
| \wh{\bs{\xi}}_{p, \ell} - \xi_p | > \vep | \xi_p |, \; [ 1 -
\mrm{sgn} ( \wh{\bs{\xi}}_{p, \ell} ) \vep] X_{j_\ell: n_\ell}  \leq
\wh{\bs{\xi}}_{p, \ell} \leq [ 1 + \mrm{sgn} ( \wh{\bs{\xi}}_{p,
\ell} ) \vep] X_{i_\ell: n_\ell} \ri \} \nonumber\\
 & \leq & \Pr \li \{  \xi_p < \f{ \wh{\bs{\xi}}_{p, \ell} } { 1 + \mrm{sgn}
(\wh{\bs{\xi}}_{p, \ell}) \vep }, \; [ 1 - \mrm{sgn} (
\wh{\bs{\xi}}_{p, \ell} ) \vep] X_{j_\ell: n_\ell}  \leq
\wh{\bs{\xi}}_{p, \ell} \leq  [ 1 +
\mrm{sgn} ( \wh{\bs{\xi}}_{p, \ell} ) \vep] X_{i_\ell: n_\ell} \ri \} \nonumber\\
&   & + \Pr \li \{  \xi_p > \f{ \wh{\bs{\xi}}_{p, \ell} } { 1 -
\mrm{sgn} (\wh{\bs{\xi}}_{p, \ell}) \vep }, \; [ 1 - \mrm{sgn} (
\wh{\bs{\xi}}_{p, \ell} ) \vep] X_{j_\ell: n_\ell}  \leq
\wh{\bs{\xi}}_{p, \ell} \leq  [ 1 + \mrm{sgn} ( \wh{\bs{\xi}}_{p,
\ell} ) \vep] X_{i_\ell: n_\ell} \ri \} \nonumber \\
& \leq & \Pr \li \{  \xi_p < \f{ \wh{\bs{\xi}}_{p, \ell} } { 1 +
\mrm{sgn} (\wh{\bs{\xi}}_{p, \ell}) \vep } \leq  X_{i_\ell: n_\ell}
\ri \} +  \Pr \li \{  \xi_p > \f{ \wh{\bs{\xi}}_{p, \ell} } { 1 -
\mrm{sgn} (\wh{\bs{\xi}}_{p, \ell})
\vep } \geq X_{j_\ell: n_\ell} \ri \} \nonumber\\
& \leq &   \Pr \li \{  X_{i_\ell: n_\ell} > \xi_p \ri \}  + \Pr \li
\{  X_{j_\ell: n_\ell} < \xi_p \ri \} \la{good2} \eel} for all $\ell
> 0$.
Combining (\ref{good3}), (\ref{good4}),  (\ref{good1}) and
(\ref{good2}) yields $\Pr \{ | \wh{\bs{\xi}}_{p} - \xi_p | > \vep |
\xi_p | \} \leq 2 \sum_{\ell = 1}^\iy \ze \de_\ell = 2 (\tau + 1)
\ze \de$.  The finite stopping property of the sampling scheme can
be shown by an argument similar to the proof of statement (I) of
Theorem \ref{Bino_Rev_noninverse_Chernoff}.

\subsection{Proof of Theorem \ref{quantile_mix} } \la{App_quantile_mix}

By the definition of the stopping rule and identity
(\ref{RI_Indentity}), we have \bee &  & \Pr \li \{ |
\wh{\bs{\xi}}_{p} - \xi_p | > \max ( \vep_a, \; \vep_r | \xi_p | ) \ri \}\\
&   \leq & \sum_{\ell = 1}^\iy \Pr \li \{ | \wh{\bs{\xi}}_{p, \ell}
-
\xi_p | > \max ( \vep_a, \; \vep_r | \xi_p | ), \ri.  \\
&   &  X_{j_\ell: n_\ell} - \max ( \vep_a, \; \mrm{sgn} ( \li.
\wh{\bs{\xi}}_{p, \ell} ) \vep_r X_{j_\ell: n_\ell} ) \leq
\wh{\bs{\xi}}_{p, \ell} \leq X_{i_\ell: n_\ell} + \max ( \vep_a, \;
\mrm{sgn} ( \wh{\bs{\xi}}_{p, \ell} ) \vep_r X_{i_\ell: n_\ell})
\ri \}\\
&   \leq & \sum_{\ell = 1}^\iy \Pr \li \{ \xi_p < \min \li (
\wh{\bs{\xi}}_{p, \ell} - \vep_a, \; \f{ \wh{\bs{\xi}}_{p, \ell} } {
1 + \mrm{sgn} (\wh{\bs{\xi}}_{p, \ell}) \vep_r } \ri ), \;
\wh{\bs{\xi}}_{p, \ell} \leq X_{i_\ell: n_\ell} + \max ( \vep_a, \;
\mrm{sgn} ( \wh{\bs{\xi}}_{p, \ell} ) \vep_r X_{i_\ell: n_\ell})
\ri \}\\
&   & +  \sum_{\ell = 1}^\iy \Pr \li \{ \xi_p > \max \li (
\wh{\bs{\xi}}_{p, \ell} + \vep_a, \; \f{ \wh{\bs{\xi}}_{p, \ell} } {
1 - \mrm{sgn} (\wh{\bs{\xi}}_{p, \ell}) \vep_r } \ri ), \;
X_{j_\ell: n_\ell} - \max ( \vep_a, \; \mrm{sgn} ( \wh{\bs{\xi}}_{p,
\ell} ) \vep_r X_{j_\ell: n_\ell} ) \leq \wh{\bs{\xi}}_{p, \ell} \ri \}\\
&   = & \sum_{\ell = 1}^\iy \Pr \li \{ \xi_p < \min \li (
\wh{\bs{\xi}}_{p, \ell} - \vep_a, \; \f{ \wh{\bs{\xi}}_{p, \ell} } {
1 + \mrm{sgn} (\wh{\bs{\xi}}_{p, \ell}) \vep_r } \ri )  \leq
X_{i_\ell: n_\ell}
\ri \}\\
&   & +  \sum_{\ell = 1}^\iy \Pr \li \{ \xi_p > \max \li (
\wh{\bs{\xi}}_{p, \ell} + \vep_a, \; \f{ \wh{\bs{\xi}}_{p, \ell} } {
1 - \mrm{sgn} (\wh{\bs{\xi}}_{p, \ell}) \vep_r } \ri ) \geq
X_{j_\ell: n_\ell} \ri \}\\
&  \leq & \sum_{\ell = 1}^\iy \Pr \li \{ X_{i_\ell: n_\ell}
> \xi_p \ri \} +  \sum_{\ell = 1}^\iy \Pr \li \{ X_{j_\ell:
n_\ell} < \xi_p \ri \} \leq 2 (\tau + 1) \ze \de,  \eee where the
last inequality follows from (\ref{good3}) and (\ref{good4}).  The
finite stopping property of the sampling scheme can be shown by an
argument similar to the proof of statement (I) of Theorem
\ref{Bino_Rev_noninverse_Chernoff}.

\sect{Adapted Branch and Bound Algorithms} \la{BBA}

In this Appendix, we shall discuss Branch and Bound Algorithms and propose a new method for improving their efficiency of global optimization.
Moreover, we will propose Adapted Branch and Bound Algorithms which can be applied to efficiently compare the maximum of a function with a
prescribed number.  The algorithms presented here have immediate applications for computing the minimum coverage probability of confidence
intervals and for quick determination whether the coverage probability exceeds a certain level.

Branch and bound (B\&B) is a general algorithm for finding optimal solutions of various optimization problems, especially in discrete and
combinatorial optimization. It consists of a systematic enumeration of all candidate solutions, where large subsets of fruitless candidates are
 safely discarded by using upper and lower estimated bounds of the quantity being optimized. The method was first proposed by A. H. Land and
A. G. Doig \cite{Land} in 1960 for discrete programming (see, \cite{BBB91, LW66, Moo91} and the references therein).

In many applications, it is desirable to compute the maximum of a function $\Psi (.)$ over a set $\varTheta$ contained in a hypercube
$\mcal{Q}_{\mrm{init}}$. By a ``hypercube'', we mean a multidimensional region of the form $\{(\se_1, \cd, \se_n): \udl{\se}_\ell \leq \se_\ell
\leq \ovl{\se}_\ell, \; 1 \leq \ell \leq n \}$. For a hypercube $\mcal{Q} \subseteq \mcal{Q}_{\mrm{init}}$, let $\max \Psi (\mcal{Q})$ denote
the maximum of the function $\Psi (.)$ over $\mcal{Q}$. Let $\Psi_{\mrm{lb}} (\mcal{Q})$ and $\Psi_{\mrm{ub}} (\mcal{Q})$ be respectively the
lower and upper bounds of $\Psi (\mcal{Q})$ such that the bounds converge as the ``size'' (i.e., ``diameter'') of $\mcal{Q}$ tends to $0$.  Let
$\ep > 0 $ be a pre-specified tolerance.  A typical B\&B algorithm is described below.
\bee &  & \underline{\tx{{\bf B\&B Algorithm (I):}}}\\
&  & \diamond \; \tx{Let} \; k \leftarrow 0 \;  \tx{and} \; \mscr{S}_0  \leftarrow  \{ \mcal{Q}_{\mrm{init}} \}.\\
&  & \diamond \; \tx{Let} \; L_0 \leftarrow \Psi_{\mrm{lb}} ( \mcal{Q}_{\mrm{init}} ) \;
\tx{and} \;  U_0 \leftarrow \Psi_{\mrm{ub}} ( \mcal{Q}_{\mrm{init}} ).\\
&  & \diamond \; \tx{While} \; U_k - L_k > \ep, \; \tx{do the following:}\\
&  & \indent \indent \star \; \tx{Eliminate any hypercube} \; \mcal{Q} \; \tx{from} \; \mscr{S}_k \; \tx{for which} \;  \Psi_{\mrm{ub}} ( \mcal{Q} ) < L_k.\\
&  & \indent \indent \star \; \tx{Pick} \; \mcal{Q}_k \in \mscr{S}_k \; \tx{for which} \;  \Psi_{\mrm{ub}} ( \mcal{Q}_k ) = U_k.\\
&  & \indent \indent \star \; \tx{Split} \; \mcal{Q}_k \; \tx{along one of its longest edges as} \; \mcal{Q}_{\mrm{I}} \; \tx{and} \;
\mcal{Q}_{\mrm{II}}.\\
&  & \indent \indent \star \; \tx{Form} \; \mscr{S}_{k+1} \; \tx{from} \; \mscr{S}_k \; \tx{by removing} \; \mcal{Q}_k \;\\
&  & \indent \indent \;\; \; \; \; \tx{and adding all} \;
\mcal{Q} \; \tx{from} \; \{ \mcal{Q}_{\mrm{I}}, \; \mcal{Q}_{\mrm{II}} \} \; \tx{such that $\mcal{Q} \cap \varTheta \neq \emptyset$}.\\
&  & \indent \indent  \star \; \tx{Let} \; L_{k+1} \leftarrow \max_{\mcal{Q} \in \mscr{S}_{k+1}} \Psi_{\mrm{lb}} (\mcal{Q}) \; \tx{and} \;
 U_{k+1} \leftarrow \max_{\mcal{Q} \in \mscr{S}_{k+1}} \Psi_{\mrm{ub}} (\mcal{Q}).\\
&  & \indent \indent \star \; \tx{Let} \; k \leftarrow k + 1.\\
&  & \diamond \; \tx{Return $U_k$ as an estimate for $\max \Psi ( \varTheta )$}.
 \eee
In some applications, it suffices to determine whether $\max \Psi ( \varTheta )$ is greater than a certain number $\ga$. For that
 purpose, we can reduce the computational complexity by revising the above B\&B algorithm  as follows.
\bee &  & \underline{\tx{{\bf Adapted B\&B Algorithm (I):}}}\\
&  & \diamond \; \tx{Let} \; k \leftarrow 0, \;  L_0 \leftarrow \Psi_{\mrm{lb}} ( \mcal{Q}_{\mrm{init}} ) \;
\tx{and} \;  U_0 \leftarrow \Psi_{\mrm{ub}} ( \mcal{Q}_{\mrm{init}} ).\\
&  & \diamond \; \tx{Let} \; \mscr{S}_0  \leftarrow  \{ \mcal{Q}_{\mrm{init}} \} \; \tx{if $U_0 > \ga$.  Otherwise, let $\mscr{S}_0$ be empty}.\\
&  & \diamond \; \tx{While $\mscr{S}_k$ is nonempty, $L_k < \ga$  and $U_k$ is greater than $\max \{ L_k + \ep, \;\ga \}$, do the following:}\\
&  & \indent \indent \star \; \tx{Pick} \; \mcal{Q}_k \in \mscr{S}_k \; \tx{for which} \;  \Psi_{\mrm{ub}} ( \mcal{Q}_k ) = U_k.\\
&  & \indent \indent \star \; \tx{Split} \; \mcal{Q}_k \; \tx{along one of its longest edges as} \; \mcal{Q}_{\mrm{I}} \; \tx{and} \;
\mcal{Q}_{\mrm{II}}.\\
&  & \indent \indent \star \; \tx{Form} \; \mscr{S}_{k+1} \; \tx{from} \; \mscr{S}_k \; \tx{by removing} \; \mcal{Q}_k\\
&  &  \indent \indent \indent \tx{and adding all
$\mcal{Q} \in \{  \mcal{Q}_{\mrm{I}}, \; \mcal{Q}_{\mrm{II}} \}$
such that $\Psi_{\mrm{ub}} (\mcal{Q}) > \ga$} \; \tx{and that $\mcal{Q} \cap \varTheta \neq \emptyset$}.\\
&  & \indent \indent  \star \; \tx{Let} \; L_{k+1} \leftarrow \max_{\mcal{Q} \in \mscr{S}_{k+1}} \Psi_{\mrm{lb}} (\mcal{Q}) \; \tx{and} \;
 U_{k+1} \leftarrow \max_{\mcal{Q} \in \mscr{S}_{k+1}} \Psi_{\mrm{ub}} (\mcal{Q}).\\
&  & \indent \indent \star \; \tx{Let} \; k \leftarrow k + 1.\\
&  & \diamond \;  \tx{If $\mscr{S}_k$ is empty and $L_k < \ga$, then declare that $\max \Psi ( \varTheta ) \leq \ga$}.\\
&  & \;\; \; \tx{Otherwise, declare that $\max \Psi ( \varTheta ) > \ga$}.
 \eee
A drawback of B\&B Algorithm (I) is that a large portion of computational effort may be wasted on branching operations which lead to no
decision. To overcome this issue, we propose the following new B\&B algorithm.
\bee &  & \underline{\tx{{\bf B\&B Algorithm (II):}}}\\
&  & \diamond \; \tx{Let} \; k \leftarrow 0 \;  \tx{and} \; \mscr{S}_0  \leftarrow  \{ \mcal{Q}_{\mrm{init}} \}.\\
&  & \diamond \; \tx{Let} \; L_0 \leftarrow \Psi_{\mrm{lb}} ( \mcal{Q}_{\mrm{init}} ) \;
\tx{and} \;  U_0 \leftarrow \Psi_{\mrm{ub}} ( \mcal{Q}_{\mrm{init}} ).\\
&  & \diamond \; \tx{While} \; U_k - L_k > \ep, \; \tx{do the following:}\\
&  & \indent  \star \; \tx{Eliminate any hypercube} \; \mcal{Q} \; \tx{from} \; \mscr{S}_k \; \tx{for which} \;  \Psi_{\mrm{ub}} ( \mcal{Q} ) < L_k.\\
&  & \indent \star \; \tx{Split each hypercube of $\mscr{S}_k$ along one of its longest edges as two smaller hypercubes}.\\
&  & \indent \star \; \tx{Let $\mscr{S}_{k+1} $ denote the set of all new hypercubes obtained from splitting hypercubes of $\mscr{S}_k$}.\\
&  & \indent  \star \; \tx{Eliminate any hypercube} \; \mcal{Q} \; \tx{from} \; \mscr{S}_{k+1} \; \tx{such that} \;  \mcal{Q} \cap \varTheta = \emptyset.\\
&  & \indent \star \; \tx{Let} \; L_{k+1} \leftarrow \max_{\mcal{Q} \in \mscr{S}_{k+1}} \Psi_{\mrm{lb}} (\mcal{Q}) \; \tx{and} \;
 U_{k+1} \leftarrow \max_{\mcal{Q} \in \mscr{S}_{k+1}} \Psi_{\mrm{ub}} (\mcal{Q}).\\
&  & \indent \star \; \tx{Let} \; k \leftarrow k + 1.\\
&  & \diamond \; \tx{Return $U_k$ as an estimate for $\max \Psi ( \varTheta )$}.
 \eee
 Similar to B\&B Algorithm (I), a drawback of Adapted B\&B Algorithm (I) is that a large portion of computational effort may be wasted on branching operations which lead to no
decision. To improve the efficiency, we propose in the following another Adapted B\&B algorithm for checking whether $\max \Psi (
\varTheta )$ is greater than a certain number $\ga$. \bee &  & \underline{\tx{{\bf Adapted B\&B Algorithm (II):}}}\\
&  & \diamond \;
\tx{Let} \; k \leftarrow 0, \;  L_0 \leftarrow \Psi_{\mrm{lb}} ( \mcal{Q}_{\mrm{init}} ) \;
\tx{and} \;  U_0 \leftarrow \Psi_{\mrm{ub}} ( \mcal{Q}_{\mrm{init}} ).\\
&  & \diamond \; \tx{Let} \; \mscr{S}_0  \leftarrow  \{ \mcal{Q}_{\mrm{init}} \} \; \tx{if $U_0 > \ga$.  Otherwise, let $\mscr{S}_0$ be empty}.\\
&  & \diamond \; \tx{While $\mscr{S}_k$ is nonempty, $L_k < \ga$ and $U_k$ is greater than $\max \{ L_k + \ep, \;\ga \}$, do the following:}\\
&  & \indent \indent \star \; \tx{Split each hypercube in $\mscr{S}_k$} \; \tx{along one of its longest edges as two new hypercubes}.\\
&  & \indent \indent \;\; \; \; \tx{Let $S_k$ denote the set of all new hypercubes obtained from this splitting procedure}.\\
&  & \indent \indent \star \; \tx{Eliminate any hypercube $\mcal{Q}$ from $S_k$
such that $\mcal{Q} \cap \varTheta = \emptyset$ or $\Psi_{\mrm{ub}} (\mcal{Q}) \leq \ga$}.\\
&  & \indent \indent \star \; \tx{Let} \; \mscr{S}_{k+1} \; \tx{be the set $S_k$ processed by the above elimination procedure}.\\
&  & \indent \indent  \star \; \tx{Let} \; L_{k+1} \leftarrow \max_{\mcal{Q} \in \mscr{S}_{k+1}} \Psi_{\mrm{lb}} (\mcal{Q}) \; \tx{and} \;
 U_{k+1} \leftarrow \max_{\mcal{Q} \in \mscr{S}_{k+1}} \Psi_{\mrm{ub}} (\mcal{Q}).\\
&  & \indent \indent \star \; \tx{Let} \; k \leftarrow k + 1.\\
&  & \diamond \;  \tx{If $\mscr{S}_k$ is empty and $L_k < \ga$, then declare that $\max \Psi ( \varTheta ) \leq \ga$}. \\
&  &   \; \; \; \tx{Otherwise, declare that $\max \Psi ( \varTheta ) > \ga$}.
 \eee

\bsk

\end{document}